\documentclass{amsproc}
\usepackage{breqn,float,amssymb,pdflscape,rotating,multirow,parskip}
\makeatletter
\@namedef{subjclassname@2020}{%
  \textup{2020} Mathematics Subject Classification}
\makeatother
\newtheorem{theorem}{Theorem}[section]

\newtheorem{proposition}[theorem]{Proposition}

\theoremstyle{definition}

\newtheorem{example}[theorem]{Example}

\theoremstyle{remark}

\DeclareMathOperator{\csch}{csch}

\DeclareMathOperator{\sech}{sech}

\numberwithin{equation}{section}
\allowdisplaybreaks
\begin{document}

\title{A collection of integrals, products and series}

%    Remove any unused author tags.

%    author one information
\author{Robert Reynolds}
\address[Robert Reynolds]{Department of Mathematics and Statistics, York University, Toronto, ON, Canada, M3J1P3}
\email[Corresponding author]{milver73@gmail.com}
\thanks{}

%%    author two information
%\author{ Allan Stauffer}
%\address[Allan Stauffer]{Department of Mathematics and Statistics, York University, Toronto, ON, Canada, M3J1P3}
%\email{stauffer@yorku.ca}
%\thanks{This research is supported by NSERC Canada under Grant 504070}



\date{}

\dedicatory{}

\begin{abstract}
This is a conspectus of definite integrals, products and series. These formulae involve special functions in the integrand and summand functions and closed form solutions. Some of the special cases are stated in terms of fundamental constants.
\end{abstract}

\maketitle
\section{Significance statement}
Research involving high energy and particle physics that is both computational and theoretical in nature use formulae involving integrals, series, and products. Researchers seeking to model processes involving elements like dark matter, lattice quantum fields, and quantum chromodynamics, even in conjunction with quantum computing, often peruse tables of formulae listing integrals, series, and products to assign equations to these processes. Upon finding equation(s) which fit the model being studied more advanced studies can be performed and a deeper understanding achieved.  This work contributes to this research by listing new formulae in easy to navigate tables, to aid researchers. 
\section{Introduction}
Tables of integrals, products and series are currently found in well known literature. Books by  Bierens de Haan, Fritz Oberhettinger, Eldon Hansen, Gradshteyn and Ryzhik to name a few have an exhaustive listing of these types of formulae. These books are listed in the reference section of this work for interested readers to peruse.\\\\
The aim of this work is to add to present literature by providing new definite integrals, products and series involving special functions and fundamental constants.  Definite integrals of special functions may have closed-form solutions or require numerical methods for evaluation and have particular uses in mathematical physics or some other branch of mathematics. Products involving special functions often arise in mathematical expressions and are encountered in various branches of mathematics, physics, and engineering. Different special functions have different properties, and their products may exhibit unique characteristics. Series involving special functions often involve expressing a function as a series expansion, where an infinite sum of terms contributes to the overall behaviour of the function. Special functions have different series representations, and these representations are useful for various purposes, including solving differential equations, understanding the behaviour of functions, and numerical computations. The plots in this work were rendered using Mathematica software by Wolfram. The formulae in this work were numerically verified using Mathematica by Wolfram.\\\\
\section{Preliminaries}
In this work, the methods in \cite{reyn4} and section (7) in \cite{plos} are applied to the stated theorems and their special cases are listed as examples that follow. The variables in the theorems are valid along the complex plane unless otherwise stated. The functions appearing in this work are; Hurwitz-Lerch zeta function $\Phi(z,s,a)$, given in section (25.14) in \cite{dlmf}, the incomplete gamma function $\Gamma(z,s)$, given in section (8) in \cite{dlmf}, logarithmic function $\log(z)$, given in section (4) in \cite{dlmf}, Laguerre polynomial $L_n^a(x)$, given in section (18) in \cite{dlmf},  Euler's polynomial $E_{n}(z)$, given in section (24) in \cite{dlmf}, Bessel function of the first kind $J_{n}(z)$,given in section (10) in \cite{dlmf},  Kummer confluent hypergeometric function $\, _1F_1\left(a,b,c\right)$, given in section (13) in \cite{dlmf}, the incomplete beta function $B_{x}(a,b)$ given in equation (8.17.7) in \cite{dlmf}, the Gudermannian function $\text{gd}(c)$ given in equation (19.10.2) in \cite{dlmf}.
\section{Main Theorems}
The theorems in this section have equations and method listed. Absence of a method indicates \cite{reyn4}. 
\section*{Theorems involving definite integrals}
In this section we tabulate definite integrals involving the logarithmic function, quotient gamma functions, logarithm of quotient trigonometric functions, product of logarithmic functions, logarithmic and rational functions, logarithmic and quotient gamma functions, product of trigonometric and logarithm of trigonometric functions, quotient exponential functions, product of arctangent and rational functions, logarithm of rational functions, product of quotient trigonometric functions and rational functions, generalized product of logarithmic functions, product of Hurwitz-Lerch zeta function and trigonometric functions, product of logarithmic function and arctangent function.
\begin{theorem}
From Eq. (3.241.4) in \cite{grad}.
\begin{multline}
\int_{0}^{\infty}\frac{x^m \log ^k(a x)}{\left(x^2+1\right)^3}dx\\
=\frac{1}{8} \pi ^{k-1} e^{\frac{1}{2} i \pi  (k+m)} \left(-2 i
   \pi  k (m-2) \Phi \left(-e^{i m \pi },1-k,\log (a)+\frac{1}{2}\right)\right. \\ \left.
   -(k-1) k \Phi \left(-e^{i m \pi },2-k,\log
   (a)+\frac{1}{2}\right)\right. \\ \left.
   +\pi ^2 (m-3) (m-1) \Phi \left(-e^{i m \pi },-k,\log (a)+\frac{1}{2}\right)\right)
\end{multline}
\end{theorem}
\begin{theorem}
Form Entry (23) in \cite{berndt1}
\begin{multline}
\int_0^{\infty } \frac{x^{-m} \Gamma (b+x) \log ^k\left(\frac{a}{x}\right)}{\Gamma (1+b+n+x)} \, dx\\
=-\sum
   _{p=0}^{\infty } \frac{e^{i m \pi } (b+p)^{-m} (2 i \pi )^{1+k} \Phi \left(e^{2 i m \pi },-k,\frac{\pi -i \log
   (a)+i \log (b+p)}{2 \pi }\right) (-n)_p}{p! \Gamma (1+n)}
\end{multline}
\end{theorem}
\begin{theorem}
From Eq. (1.9.51) in \cite{bateman}.
\begin{multline}
\int_{0}^{\infty}e^{-i m x} \left(e^{2 i m x} (\log (a)+i x)^{k-1}+\frac{m \left(e^{2 i m x} (\log (a)+i x)^k+(\log (a)-i x)^k\right)}{k}\right. \\ \left.
+(\log (a)-i
   x)^{k-1}\right) \log \left(\frac{\cos (\alpha )+\cosh (x)}{\cos (\beta )+\cosh (x)}\right)dx\\
=\frac{(2 \pi )^{k+1} }{k}\left(e^{(\pi -\alpha ) m} \Phi
   \left(e^{2 m \pi },-k,\frac{-\alpha +\log (a)+\pi }{2 \pi }\right)\right. \\ \left.
+e^{(\alpha +\pi ) m} \Phi \left(e^{2 m \pi },-k,\frac{\alpha +\log (a)+\pi }{2
   \pi }\right)\right. \\ \left.
   -e^{(\pi -\beta ) m} \left(\Phi \left(e^{2 m \pi },-k,\frac{-\beta +\log (a)+\pi }{2 \pi }\right)\right.\right. \\ \left.\left.
+e^{2 \beta  m} \Phi \left(e^{2 m
   \pi },-k,\frac{\beta +\log (a)+\pi }{2 \pi }\right)\right)\right)
\end{multline}
\end{theorem}
\begin{theorem}
From Eq. (2.1.5.7) in \cite{brychkov}.
\begin{multline}
\int_{0}^{\infty}\frac{x^{-1+m} \left(-1+x^u\right) \left(\left(-\frac{i \pi }{v}+\log (a x)\right)^k-e^{\frac{2 i m \pi}{v}} \left(\frac{i \pi }{v}+\log (a x)\right)^k\right)}{-1+x^v}dx\\=2^{2+k} e^{\frac{i \pi  (2 m+u)}{v}} \pi ^{1+k}
   \left(\frac{i}{v}\right)^{2+k} v \Phi \left(e^{\frac{2 i \pi  (m+u)}{v}},-k,\frac{\pi -i v \log (a)}{2 \pi
   }\right) \sin \left(\frac{\pi  u}{v}\right)
\end{multline}
\end{theorem}
\begin{theorem}
From Eq. (2.6.15.8) in \cite{prud1}.
\begin{multline}
\int_{0}^{\infty}\log \left(\frac{(x+1) \left(b^2+x\right)}{(b+x)^2}\right) \log ^k(a x)\frac{dx}{x}\\
=\frac{2^{k+3} (i \pi )^{k+2} \zeta \left(-k-1,\frac{-i \log (a)-i \log (b)+\pi }{2 \pi
   }\right)}{k+1}\\
   -\frac{(2 i \pi )^{k+2} \zeta \left(-k-1,\frac{-i \log (a)-2 i \log (b)+\pi }{2 \pi }\right)}{k+1}-\\
\frac{(2 i \pi )^{k+2} \zeta \left(-k-1,\frac{\pi -i \log (a)}{2 \pi
   }\right)}{k+1}
\end{multline}
\end{theorem}
\begin{theorem}
From Entry 9, in Table 292 on page 419 in \cite{bdh}.
\begin{multline}
\int_{0}^{\pi/4}\frac{\csc (2 x) \left(\cot ^m(x) \log (\tan (x)) \log ^k(a \cot (x))+\tan ^m(x) \log (\cot (x)) \log ^k(a \tan (x))\right)}{\log (\cot (x))
   \left(\cot ^p(x)-\tan ^p(x)\right)}dx\\
=\frac{i \pi ^{k+1} \left(\frac{i}{p}\right)^k e^{\frac{i \pi  m}{p}} \Phi \left(-e^{\frac{i m \pi
   }{p}},-k,1-\frac{i p \log (a)}{\pi }\right)}{2 p}-\frac{i \pi  \log ^k(a)}{4 p}
\end{multline}
\end{theorem}
\begin{theorem}
From Eq. (2.2.7.3) in \cite{prud1}.
\begin{multline}
\int_{-\alpha }^{\alpha } \frac{\left(\frac{\alpha +x}{\alpha
   -x}\right)^m \log ^k\left(\frac{a (\alpha +x)}{\alpha -x}\right)}{\alpha
   ^2+x^2} \, dx=\frac{\pi ^{k+1} e^{\frac{1}{2} i \pi  (k+m)} \Phi
   \left(-e^{i m \pi },-k,\frac{1}{2}-\frac{i \log (a)}{\pi }\right)}{\alpha
   }
\end{multline}
\end{theorem}
\begin{theorem}
From Eq. (3.414.31) in \cite{grad}.
\begin{multline}
\int _{0}^{\infty }\frac{e^{x (p-m)} (\log (a)-x)^k-e^{m x} (\log (a)+x)^k}{e^{p x}-1}dx\\
=-\frac{i \pi  \left(\log ^k(a)+2^{k+1} \pi
   ^k \left(\frac{i}{p}\right)^k e^{\frac{2 i \pi  m}{p}} \Phi \left(e^{\frac{2 i m \pi }{p}},-k,1-\frac{i p \log
   (a)}{2 \pi }\right)\right)}{p}
\end{multline}
\end{theorem}
\begin{theorem}
From Eq. (2.6.8.9) in \cite{prud1}.
\begin{multline}
\int_0^{\infty }\frac{x^{m-1} \log ^k(\alpha  x)}{a^2+2 a x \cos (\varphi )+x^2}dx\\
=(2 i)^k \pi ^{k+1} a^{m-2} e^{i (m (\pi -\varphi )+\varphi )} \csc (\varphi ) \left(e^{2 i (m-1) \varphi } \Phi \left(e^{2 i m \pi
   },-k,\frac{\varphi -i \log (a)-i \log (\alpha )+\pi }{2 \pi }\right)\right. \\ \left.
-\Phi \left(e^{2 i m \pi },-k,-\frac{\varphi +i \log (a)+i \log (\alpha )-\pi }{2 \pi }\right)\right)
\end{multline}
\end{theorem}
\begin{theorem}
Schr\"{o}der's integral and Gregory coefficients \cite{schroder}. From Eq. (3.251.11) in \cite{grad}.
\begin{multline}
\int_0^{\infty } \frac{x^t  \log ^k(a x)}{(b x+1)^{m}}  \, dx\\
=-k! \sum
   _{p=0}^{m-1} \left(\sum _{z=0}^p \frac{(2 i \pi
   )^{k-z+1} (-1)^{2 m+p+1} S_p^{(z)} (-t)_{m-p-1} \Phi \left(e^{2 i \pi 
   t},z-k,\frac{-i \log (a)+i \log (b)+\pi }{2 \pi }\right)}{e^{-i \pi  t} b^{t+1} p! (k-z)! \Gamma
   (m-p)}\right)
\end{multline}
\end{theorem}
\begin{theorem}
From Eq. (3.195) in \cite{grad}.
\begin{multline}
\int_{0}^{\infty}\frac{\left(\frac{x+1}{b+x}\right)^m \log ^k\left(\frac{a (x+1)}{b+x}\right)}{(x+1) (b+x)}dx\\
=\frac{a^{-m} (-m)^{-k} }{(b-1) m}(\Gamma (k+1,-m
   \log (a))-\Gamma (k+1,m (\log (b)-\log (a))))
\end{multline}
\end{theorem}
%
%\begin{theorem}
%From Eq. (4.293.10) in \cite{grad}.
%\begin{multline}
%\int_{0}^{\infty}x^{m-1} \tanh ^{-1}(b x) \log ^k(a x)dx\\
%%
%=\sum_{y=0}^{\infty}\frac{i \pi }{m} \left(-\frac{1}{b}\right)^m e^{-m \left(\log (a)+\log \left(-\frac{1}{b}\right)\right)}
%   \left((-1)^k m^{-k} \Gamma (k+1)+\left(\log (a)+\log \left(-\frac{1}{b}\right)+i (2 \pi  y+\pi )\right)^k \right. \\ \left.
%   \left(-m \left(\log (a)+\log
%   \left(-\frac{1}{b}\right)+i (2 \pi  y+\pi )\right)\right)^{-k} \right. \\ \left.
%   \left(-k \Gamma (k)+\Gamma \left(k+1,-m \left(i \pi  (2 y+1)+\log (a)+\log
%   \left(-\frac{1}{b}\right)\right)\right)\right)\right)\\
%   -\sum_{y=0}^{\infty}\frac{i \pi }{m} \left(\frac{1}{b}\right)^m e^{-m \left(\log (a)+\log
%   \left(\frac{1}{b}\right)\right)} \left((-1)^k m^{-k} \right. \\ \left.
%   \Gamma (k+1)+\left(\log (a)+\log \left(\frac{1}{b}\right)+i (2 \pi  y+\pi )\right)^k \left(-m
%   \left(\log (a)+\log \left(\frac{1}{b}\right)+i (2 \pi  y+\pi )\right)\right)^{-k}\right. \\ \left.
%    \left(-k \Gamma (k)+\Gamma \left(k+1,-m \left(i \pi  (2 y+1)+\log
%   (a)+\log \left(\frac{1}{b}\right)\right)\right)\right)\right)
%\end{multline}
%\end{theorem}
%%
%
\begin{theorem}
From Eq. (1.7.7.122) in \cite{obert}.
\begin{multline}
\int_{0}^{\infty}e^{-i m x} \left(e^{2 i m x} (\log (a)+i x)^{k-1}+\frac{m \left(e^{2 i m x} (\log (a)+i x)^k+(\log (a)-i x)^k\right)}{k}\right. \\ \left.
+(\log (a)-ix)^{k-1}\right) \tan ^{-1}(\sinh (\alpha ) \sech(b x))dx\\
=\frac{i \pi ^{k+1}}{k} \left(\frac{1}{b}\right)^k e^{\frac{(\pi -2 i \alpha ) m}{2 b}}
   \left(\Phi \left(-e^{\frac{m \pi }{b}},-k,\frac{-2 i \alpha +2 b \log (a)+\pi }{2 \pi }\right)\right. \\ \left.
-e^{\frac{2 i \alpha  m}{b}} \Phi \left(-e^{\frac{m \pi
   }{b}},-k,\frac{2 i \alpha +2 b \log (a)+\pi }{2 \pi }\right)\right)
\end{multline}
\end{theorem}
\begin{theorem}
From Eq. (2.1.16) in \cite{brychkov}.
\begin{multline}
\int_{0}^{\infty}\frac{x^{-m-\frac{1}{2}} \left(\beta +\sqrt{x (2 \beta +x)}+x\right)^{-m} \log ^k\left(\frac{a}{x \left(\beta +\sqrt{x (2 \beta
   +x)}+x\right)}\right)}{\sqrt{2 \beta +x}}dx\\
=-i \pi ^{k+1} 2^{2 k+m+1} e^{\frac{1}{2} i \pi  (k+4 m)} \beta ^{-2 m}\\
 \Phi \left(e^{4 i m \pi
   },-k,\frac{-i \log (2 a)+2 i \log (\beta )+2 \pi }{4 \pi }\right)
\end{multline}
\end{theorem}
\begin{theorem}
From Eq. (2.4.4.10(ii)) in \cite{prud1}.
\begin{multline}
\int_0^{\infty } \frac{ \cosh (c u) }{\cosh (b u)}\left((-u+\log (a))^k+e^{2 m
   u} (u+\log (a))^k\right) \, \frac{du}{e^{m u}}\\
=\frac{2^{-1-k}
   \left(\frac{i}{b}\right)^{-1+k} \pi ^{1+k}}{b^2} \left(2^{1+k} e^{\frac{i (b-c+m)
   \pi }{2 b}} \Phi \left(e^{\frac{i (b-c+m) \pi }{b}},-k,\frac{1}{2}-\frac{i b
   \log (a)}{\pi }\right)\right. \\ \left.
+e^{\frac{i (b+c+m) \pi }{2 b}} \left(\Phi
   \left(-e^{\frac{i (b+c+m) \pi }{2 b}},-k,1-\frac{2 i b \log (a)}{\pi
   }\right)+\Phi \left(e^{\frac{i (b+c+m) \pi }{2 b}},-k,1-\frac{2 i b \log
   (a)}{\pi }\right)\right)\right)
\end{multline}
\end{theorem}
\begin{theorem}
From Eq. (3.241.4) in \cite{grad}.
\begin{multline}\label{dilf}
\int_{0}^{\infty} \ldots \int_{0}^{\infty}\frac{x_1^{m-1} }{x_1^n+1}\left(\prod _{l=1}^j \frac{x_{l+1}^{2^{l-1}
   m+\frac{n}{2}-1}}{x_{l+1}^n+1}\right) \log ^k\left(a x_1 \left(\prod _{l=1}^j
   x_{l+1}^{2^{l-1}}\right)\right)dx_{1} \ldots dx_{j+1}\\
   =-i (2 \pi )^{j+k+1} n^{-j-1}
   \left(\frac{i 2^j}{n}\right)^k e^{\frac{i \pi  2^j m}{n}} \Phi
   \left(e^{\frac{i 2^{j+1} m \pi }{n}},-k,\frac{\pi -i 2^{-j} n \log (a)}{2 \pi
   }\right)
\end{multline}
\end{theorem}
\begin{theorem}
From Eq.  (2.4.6.11) in \cite{prud1}.
\begin{multline}\label{prud24611}
\int_0^{\infty } \frac{e^{-m x} \left((-x+\log (a))^k+e^{2 m x} (x+\log (a))^k\right)}{(\cos (t)+\cosh (c x))^2} \,dx\\
=\frac{\left(2 c^{-2-k} e^{\frac{i m (\pi -t)}{c}} (2 i \pi )^k\right) }{-1+\cos (2 t)}\left(c k \Phi \left(e^{\frac{2 i m \pi }{c}},1-k,\frac{\pi -t-i c \log (a)}{2 \pi }\right)\right. \\ \left.
+c e^{\frac{2 i m t}{c}} k \Phi \left(e^{\frac{2 i m \pi }{c}},1-k,\frac{\pi +t-i c \log (a)}{2 \pi }\right)\right. \\ \left.
+2 \pi  (i m+c \cot (t)) \Phi \left(e^{\frac{2 i m \pi }{c}},-k,\frac{\pi -t-i c \log (a)}{2 \pi }\right)\right. \\ \left.
+2 i e^{\frac{2 i m t}{c}} \pi  (m+i c \cot (t)) \Phi \left(e^{\frac{2 i m \pi
   }{c}},-k,\frac{\pi +t-i c \log (a)}{2 \pi }\right)\right)
\end{multline}
\end{theorem}
\begin{theorem}
From Eq.  (1.7.7.20) in \cite{obert}.
\begin{multline}\label{ober_17720}
\int_0^{\infty } \frac{e^{-m x} \left((-x+a)^k+e^{2 m x} (x+a)^k\right) \sinh (b x)}{\sinh (c \pi  x)}
   \, dx\\
=\left(\frac{i}{c}\right)^{1+k} \left(e^{\frac{i (-b+m)}{c}} \Phi \left(-e^{\frac{i
   (-b+m)}{c}},-k,1-i c a\right)-e^{\frac{i (b+m)}{c}} \Phi \left(-e^{\frac{i (b+m)}{c}},-k,1-i c
   a\right)\right)
\end{multline}
\end{theorem}
\begin{theorem}
From Eq. (1.7.7.1) in \cite{obert}
\begin{multline}\label{eq_sech}
\int_0^{\infty } \left(e^{-m x} (-x+\log (a))^k+e^{m x} (x+\log
   (a))^k\right) \sech(b x) \, dx\\
=\frac{2
   \left(\left(\frac{i}{b}\right)^k e^{\frac{i m \pi }{2 b}} \pi ^{1+k} \Phi
   \left(-e^{\frac{i m \pi }{b}},-k,\frac{1}{2}-\frac{i b \log (a)}{\pi
   }\right)\right)}{b}
\end{multline}
\end{theorem}
\begin{theorem}
From Eq. (1.7.7.2) in \cite{obert}.
\begin{multline}\label{eq_sech2}
\int_0^{\infty } e^{-m x} \left((-x+\log (a))^k+e^{2 m x} (x+\log (a))^k\right) \sech^2(b x) \,
   dx\\
=-\frac{2 \left(\frac{i}{b}\right)^k e^{\frac{i m \pi }{2 b}} \pi ^k \left(b k \Phi \left(e^{\frac{i m \pi
   }{b}},1-k,\frac{1}{2}-\frac{i b \log (a)}{\pi }\right)+i m \pi  \Phi \left(e^{\frac{i m \pi
   }{b}},-k,\frac{1}{2}-\frac{i b \log (a)}{\pi }\right)\right)}{b^2}
\end{multline}
\end{theorem}
\begin{theorem}
From Eq. (3.511.9) in \cite{grad}.
\begin{multline}\label{eq_csch2}
\int_0^{\infty } \csch^2(b x) \left(-2 \log ^k(a)+e^{-2 m x} (-2 x+\log (a))^k+e^{2 m x} (2 x+\log
   (a))^k\right) \, dx\\
=\frac{2 }{b^2 \log (a)}\left(\left(\frac{i}{b}\right)^k e^{\frac{2 i m \pi }{b}} (2 \pi )^k \left(b k \Phi
   \left(e^{\frac{2 i m \pi }{b}},1-k,1-\frac{i b \log (a)}{2 \pi }\right)\right.\right. \\ \left.\left.
+2 i m \pi  \Phi \left(e^{\frac{2 i m \pi
   }{b}},-k,1-\frac{i b \log (a)}{2 \pi }\right)\right) \log (a)+\log ^k(a) (i k \pi +(b+i m \pi ) \log
   (a))\right)
\end{multline}
\end{theorem}
\begin{theorem}
From Eq. (2.4.6.14) in \cite{prud1}.
\begin{multline}\label{eq_24614}
\int_0^{\infty } \frac{e^{-m x} \left((-x+\log (a))^k-e^{2 m x} (x+\log (a))^k\right) \sinh (b x)}{\cos
   (t)+\cosh (c x)} \, dx\\
=\frac{\left(i e^{-\frac{2 i (b (\pi +t)+m (2 \pi +t))}{c}} \pi  \csc \left(\frac{2 b \pi
   }{c}\right) \csc (t)\right) }{2 c}\\\left(\left(-\frac{i}{c}\right)^k \left(e^{\frac{3 i (b+m) (\pi +t)}{c}}-e^{\frac{i
   (3 m (\pi +t)+b (5 \pi +t))}{c}}\right) (2 \pi )^k \Phi \left(e^{\frac{2 i (b-m) \pi }{c}},-k,\frac{\pi -t+i c
   \log (a)}{2 \pi }\right)\right. \\ \left.
-\left(-\frac{i}{c}\right)^k e^{\frac{i (b (3 \pi +t)+m (\pi +3 t))}{c}}
   \left(-1+e^{\frac{2 i b (\pi +t)}{c}}\right) (2 \pi )^k \Phi \left(e^{\frac{2 i (b-m) \pi }{c}},-k,\frac{3 \pi
   -t+i c \log (a)}{2 \pi }\right)\right. \\ \left.
-\left(-\frac{i}{c}\right)^k \left(e^{\frac{i (b+m) (3 \pi +t)}{c}}-e^{\frac{i (5
   b \pi +3 m \pi +3 b t+m t)}{c}}\right) (2 \pi )^k \Phi \left(e^{\frac{2 i (b-m) \pi }{c}},-k,\frac{\pi +t+i c
   \log (a)}{2 \pi }\right)\right. \\ \left.
-\left(-\frac{i}{c}\right)^k \left(e^{\frac{i (3 b+m) (\pi +t)}{c}}-e^{\frac{i ((5 b+m)
   \pi +(b+m) t)}{c}}\right) (2 \pi )^k \Phi \left(e^{\frac{2 i (b-m) \pi }{c}},-k,\frac{3 \pi +t+i c \log (a)}{2
   \pi }\right)\right. \\ \left.
+\left(\frac{i}{c}\right)^k \left(e^{\frac{i (b+m) (5 \pi +t)}{c}}-e^{\frac{i (3 b (\pi +t)+m (5 \pi
   +t))}{c}}\right) (2 \pi )^k \Phi \left(e^{\frac{2 i (b+m) \pi }{c}},-k,\frac{\pi -t-i c \log (a)}{2 \pi
   }\right)\right. \\ \left.
-\left(\frac{i}{c}\right)^k \left(e^{\frac{i (b+m) (5 \pi +3 t)}{c}}-e^{\frac{i (3 b \pi +5 m \pi +b t+3
   m t)}{c}}\right) (2 \pi )^k \Phi \left(e^{\frac{2 i (b+m) \pi }{c}},-k,\frac{\pi +t-i c \log (a)}{2 \pi
   }\right)\right. \\ \left.
+\left(\frac{i}{c}\right)^k e^{\frac{i (3 b \pi +7 m \pi +3 m t)}{c}} \left(e^{\frac{3 i b
   t}{c}}-e^{\frac{i b (2 \pi +t)}{c}}\right) (2 \pi )^k \Phi \left(e^{\frac{2 i (b+m) \pi }{c}},-k,\frac{3 \pi +t-i
   c \log (a)}{2 \pi }\right)\right. \\ \left.
+\left(\frac{i}{c}\right)^k e^{\frac{i (b (3 \pi +t)+m (7 \pi +t))}{c}}
   \left(-1+e^{\frac{2 i b (\pi +t)}{c}}\right) (2 \pi )^k \Phi \left(e^{\frac{2 i (b+m) \pi }{c}},-k,-\frac{-3 \pi
   +t+i c \log (a)}{2 \pi }\right)\right. \\ \left.
-e^{\frac{3 i (b+m) (\pi +t)}{c}} \left(-\frac{i (\pi -t)}{c}+\log
   (a)\right)^k+e^{\frac{i (b+3 m) (\pi +t)}{c}} \left(-\frac{i (\pi -t)}{c}+\log (a)\right)^k\right. \\ \left.
-e^{\frac{i (b (\pi
   +t)+m (5 \pi +t))}{c}} \left(\frac{i (\pi -t)}{c}+\log (a)\right)^k+e^{\frac{i (3 b (\pi +t)+m (5 \pi +t))}{c}}
   \left(\frac{i (\pi -t)}{c}+\log (a)\right)^k\right. \\ \left.
+e^{\frac{i (b+m) (3 \pi +t)}{c}} \left(-\frac{i (\pi +t)}{c}+\log
   (a)\right)^k-e^{\frac{i (m (3 \pi +t)+b (\pi +3 t))}{c}} \left(-\frac{i (\pi +t)}{c}+\log (a)\right)^k\right. \\ \left.
+e^{\frac{i
   (b \pi +5 m \pi +3 m t)}{c}} \left(e^{\frac{3 i b t}{c}}-e^{\frac{i b (2 \pi +t)}{c}}\right) \left(\frac{i (\pi
   +t)}{c}+\log (a)\right)^k\right)
\end{multline}
\end{theorem}
\begin{theorem}
From Eq.  (3.984.5) in \cite{grad}.
\begin{multline}
\int_{0}^{\infty}e^{x (m-2 i \alpha )} x^{-1+s} \left(-e^{2 m x} \Phi \left(-e^{2 x (m-i \alpha )},1-s,\frac{1}{2}
   \left(3+\frac{a}{x}\right)\right)\right. \\ \left.
+\Phi \left(-e^{2 x (m-i \alpha )},1-s,\frac{1}{2}
   \left(\frac{a}{x}+1\right)\right)+e^{4 i x \alpha } \left(e^{2 m x} \Phi \left(-e^{2 x (m+i \alpha
   )},1-s,\frac{1}{2} \left(3+\frac{a}{x}\right)\right)\right.\right. \\ \left.\left.
-\Phi \left(-e^{2 x (m+i \alpha )},1-s,\frac{1}{2}
   \left(\frac{a}{x}+1\right)\right)\right)\right) \sec (x \alpha )dx\\
=\frac{\pi  }{2^s e^{a (m+i\alpha )}}\left(\frac{e^{2 i a \alpha } \Gamma
   (s,-a (m-i \alpha ))}{(-m+i \alpha )^s}-\frac{\Gamma (s,-a (m+i \alpha ))}{(-m-i \alpha )^s}\right)
\end{multline}
\end{theorem}
\begin{theorem}
From Eq. (2.2.12.10) in \cite{prud1}.
\begin{multline}
\int_0^1 x^{-m} \left(x^{-1+2 m} \alpha ^m (\alpha -x \alpha +\beta )^{-m} \log ^k\left(\frac{a x \alpha
   }{\alpha -x \alpha +\beta }\right)\right. \\ \left.
+\beta ^{1-m} (\alpha +\beta -x \beta )^{-1+m} \log ^k\left(\frac{a (\alpha
   +\beta -x \beta )}{x \beta }\right)\right) \, dx\\
=-e^{i m \pi } (2 i \pi )^{1+k} \Phi \left(e^{2 i m \pi
   },-k,\frac{\pi -i \log (a)}{2 \pi }\right)
\end{multline}
\end{theorem}
\begin{theorem}
From Eq. (1.2.4.90) in \cite{bateman1}.
\begin{multline}
\int_0^{\infty } e^{-i m x} \left((-i x+\log (a))^{-2+k}-e^{2 i m x} (i
   x+\log (a))^{-2+k}\right. \\ \left.
+\frac{2 k m}{(-1+k) k} \left((-i x+\log (a))^{-1+k}-e^{2 i m x} (i
   x+\log (a))^{-1+k}\right)\right. \\ \left.
+\left(1+m^2\right) \left((-i x+\log (a))^k-e^{2 i
   m x} (i x+\log (a))^k\right)\right) \log \left(1-e^{-2 x}\right)
   \sinh (x) \, dx\\
=\frac{2 a^{-i-m} }{i (-1+k)
   k}\left(-a^{2 i} (i-m)^{-1-k} \Gamma
   (1+k,-((-i+m) \log (a)))\right. \\ \left.
-(-i-m)^{-1-k} \Gamma (1+k,-((i+m) \log
   (a)))\right)\\
+\pi  \left(-2 e^{m \pi } \pi ^k \Phi \left(-e^{m \pi
   },-k,\frac{\pi +\log (a)}{\pi }\right)+\log ^k(a)\right)
\end{multline}
\end{theorem}
\begin{theorem}
From Eq. (2.4.14.9) in \cite{prud1}.
\begin{multline}
\int_{-\infty }^{\infty } \frac{e^{-i m x} (-i x+\log (a))^k}{\sinh (x)+\sinh (\alpha )} \, dx\\
=i e^{-i m \alpha } \pi  \left(-e^{m \pi } \pi ^k \left(\Phi \left(-e^{m \pi
   },-k,\frac{\pi -i \alpha +\log (a)}{\pi }\right)+\Phi \left(e^{m \pi },-k,\frac{\pi -i \alpha +\log (a)}{\pi }\right)\right.\right. \\ \left.\left.
+e^{2 i m \alpha } \left(\Phi \left(-e^{m \pi },-k,\frac{\pi +i
   \alpha +\log (a)}{\pi }\right)-\Phi \left(e^{m \pi },-k,\frac{\pi +i \alpha +\log (a)}{\pi }\right)\right)\right)\right. \\ \left.
+e^{2 i m \alpha } (i \alpha +\log (a))^k\right) \text{sech}(\alpha
   )
\end{multline}
\end{theorem}
\begin{theorem}
From Eq. (2.4.6.27) in \cite{prud1}.
\begin{multline}
\int_0^{\infty } \frac{e^{-m x} (-x+\log (a))^k+e^{m x} (x+\log (a))^k}{\cosh ^2(c x)+\sinh ^2(t)} \, dx\\
=\frac{2 \left(\frac{i}{c}\right)^{-1+k} e^{\frac{i m (\pi +2 i t)}{2 c}}}{c^2}
   \pi ^{1+k} \text{csch}(2 t) \left(e^{\frac{2 m t}{c}} \Phi \left(e^{\frac{i m \pi }{c}},-k,\frac{\pi -2 i t-2 i c \log (a)}{2 \pi }\right)\right. \\ \left.
-\Phi \left(e^{\frac{i m \pi
   }{c}},-k,\frac{\pi +2 i t-2 i c \log (a)}{2 \pi }\right)\right)
\end{multline}
\end{theorem}
\begin{theorem}
From Eq. (2.4.6.26) in \cite{prud1}
\begin{multline}
\int_0^{\infty } \frac{\cosh (c x) \left(e^{-m x} (-x+\log (a))^k+e^{m x} (x+\log (a))^k\right)}{ \left(1+\cosh ^2(c x)\right)} \, dx\\
=\frac{\left(\frac{i}{c}\right)^k e^{\frac{i
   m \left(\pi +2 i \log \left(1+\sqrt{2}\right)\right)}{2 c}} \pi ^{1+k} }{
   \sqrt{2} c}\left(\left(1+\sqrt{2}\right)^{\frac{2 m}{c}} \Phi \left(e^{\frac{i (c+m) \pi }{c}},-k,\frac{\pi -2 i \log
   \left(1+\sqrt{2}\right)-2 i c \log (a)}{2 \pi }\right)\right. \\ \left.
+\Phi \left(e^{\frac{i (c+m) \pi }{c}},-k,\frac{\pi +2 i \log \left(1+\sqrt{2}\right)-2 i c \log (a)}{2 \pi }\right)\right)
\end{multline}
\end{theorem}
\begin{theorem}
From Eq. (2.2.7.12) in \cite{prud1}.
\begin{multline}
\int_0^{\infty } \left(\frac{x^{-2 m} \log ^k\left(\frac{a}{x^2}\right)}{(-1+x)^2}+\frac{x^{2 m} \log ^k\left(a x^2\right)}{(-1+x)^2}-\frac{2
   \log ^k(a)}{(-1+x)^2}\right) \, dx\\
=2^{3+2 (-1+k)} k e^{4 i m \pi } (i \pi )^k \Phi \left(e^{4 i m \pi },1-k,1-\frac{i \log (a)}{4 \pi }\right)\\
+2^{3+2
   k} e^{4 i m \pi } m (i \pi )^{1+k} \Phi \left(e^{4 i m \pi },-k,1-\frac{i \log (a)}{4 \pi }\right)\\
+4 i \pi  k \log ^{-1+k}(a)+4 i m \pi  \log ^k(a)+2
   \log ^k(a)
\end{multline}
\end{theorem}
\begin{theorem}
\begin{multline}
\int_0^{\infty } \frac{\log ^k(a x) \log \left(x^t \left(1+\frac{1}{x}\right)_t\right)}{x} \, dx=\sum _{m=1}^t -\frac{(2 i \pi )^{2+k} \zeta
   \left(-1-k,\frac{\pi -i \log (a)+i \log (m)}{2 \pi }\right)}{1+k}
\end{multline}
\end{theorem}
\begin{theorem}
From  (4.4.4) in \cite{polyanin}. 
\begin{multline}
\int_0^{\infty } \frac{e^{-i m \pi } x^{-1+m} \log (x) \left((-i \pi +\log (a x))^k-e^{2 i m \pi } (i \pi
   +\log (a x))^k\right)}{x+\alpha } \, dx\\
=-2 i \pi  \alpha ^{-1+m} \left(e^{2 i m \pi } (2 i \pi )^{1+k} \Phi
   \left(e^{2 i m \pi },-k,\frac{2 \pi -i \log (a)-i \log (\alpha )}{2 \pi }\right)\right. \\ \left.
+(i \pi +\log (\alpha )) (\log
   (a)+\log (\alpha ))^k\right)
\end{multline}
\end{theorem}
\begin{theorem}
From Eq.(3.1) in \cite{berndti}.
\begin{multline}
\int_{0}^{\infty}\frac{x^{2 m+1} \left(\log ^k\left(a x^2\right)+i^{2 m} \log ^k\left(-a x^2\right)\right)}{\cos
   (x)+\cosh (x)}dx\\
=\left(\frac{i}{2}\right)^m \pi ^{2 (m+1)} \sum_{p=0}^{\infty} \frac{(-1)^p (1+2 p)^{1+2 m} \log
   ^k\left(\frac{1}{2} \pi ^2 i a (2 p+1)^2\right)}{\cosh \left(\frac{1}{2} (2 p+1) \pi \right)}
\end{multline}
\end{theorem}
\begin{theorem}
From Eq. (3.2.18) in \cite{bateman}.
\begin{multline}
\int_{-\infty }^{\infty } \frac{e^{x (2-i m-\lambda )} (-i x+\log (a))^k}{\left(1+e^x \alpha \right)
   \left(1+e^x \beta \right)} \, dx\\
=\frac{i e^{\pi  (m-i \lambda )} (2 \pi )^{1+k}}{\alpha  (\alpha -\beta ) \beta } \left(-\alpha ^{i m+\lambda }
   \beta  \Phi \left(e^{2 \pi  (m-i \lambda )},-k,\frac{\pi +\log (a)+i \log (\alpha )}{2 \pi }\right)\right. \\ \left.
+\alpha 
   \beta ^{i m+\lambda } \Phi \left(e^{2 \pi  (m-i \lambda )},-k,\frac{\pi +\log (a)+i \log (\beta )}{2 \pi
   }\right)\right)
\end{multline}
\end{theorem}
\begin{theorem}
From Eq. (2.6.40.7(a)) in \cite{prud1}. Extended Winckler form see pp. 359 in \cite{winckler}.
\begin{multline}
\int_0^{\infty } e^{-i m x} \left((-i x+\log (a))^{-1+k}+e^{2 i m x} (i x+\log (a))^{-1+k}\right. \\ \left.
+\frac{m
   \left((-i x+\log (a))^k+e^{2 i m x} (i x+\log (a))^k\right)}{k}\right) \log \left(1-e^{-x \alpha }\right) \,
   dx\\
=\frac{\frac{a^{-m} (-m)^{-k} \alpha  \Gamma (1+k,-m \log (a))}{m}+e^{\frac{2 m \pi }{\alpha }} (2 \pi
   )^{1+k} \alpha ^{-k} \Phi \left(e^{\frac{2 m \pi }{\alpha }},-k,1+\frac{\alpha  \log (a)}{2 \pi }\right)+\pi 
   \log ^k(a)}{k}
\end{multline}
\end{theorem}
\begin{theorem}
From Eq.  (2.6.40.7(b)) in \cite{prud1}. Extended Winckler form see pp. 359 in \cite{winckler}.
\begin{multline}
\int_0^{\infty } e^{-i m x} \left(k \left((-i x+\log (a))^{-1+k}+e^{2 i m x} (i x+\log (a))^{-1+k}\right)\right. \\ \left.
+m   \left((-i x+\log (a))^k+e^{2 i m x} (i x+\log (a))^k\right)\right) \log \left(1+e^{-x \alpha }\right) \,
   dx\\
=\frac{a^{-m} (-m)^{-k} \alpha  \Gamma (1+k,-m \log (a))}{m}+e^{\frac{m \pi }{\alpha }} (2 \pi )^{1+k}
   \left(\frac{1}{\alpha }\right)^k \Phi \left(e^{\frac{2 m \pi }{\alpha }},-k,\frac{\pi +\alpha  \log (a)}{2 \pi
   }\right)
\end{multline}
\end{theorem}
\begin{theorem}
From Eq. (2.4.4.12) in \cite{prud1}.
\begin{multline}
\int_0^{\infty } e^{-m x} \left(-(-x+\log (a))^k+e^{2 m x} (x+\log (a))^k\right) \sech(c x) \tanh (c
   x) \, dx\\
=\frac{2 \left(\frac{i}{c}\right)^k e^{\frac{i m \pi }{2 c}} \pi ^k \left(-i c k \Phi \left(-e^{\frac{i
   m \pi }{c}},1-k,\frac{1}{2}-\frac{i c \log (a)}{\pi }\right)+m \pi  \Phi \left(-e^{\frac{i m \pi
   }{c}},-k,\frac{1}{2}-\frac{i c \log (a)}{\pi }\right)\right)}{c^2}
\end{multline}
\end{theorem}
\begin{theorem}
From Eq. (2.5.16.23) in \cite{prud1}.
\begin{multline}
\int_0^{\pi } e^{-2 i x} \cot (x) \left(\Phi \left(e^{-2 (m+i x)},-k,a\right)-e^{4 i x} \Phi \left(e^{-2 m+2
   i x},-k,a\right)\right) \, dx\\=-2 i \pi  \Phi \left(e^{-2 m},-k,a\right)
\end{multline}
\end{theorem}
\begin{theorem}
From Eq. (1.7.7.107) in \cite{obert}
\begin{multline}
\int_0^{\infty } e^{-i m x} \left((-i x+\log (a))^{-1+k}+e^{2 i m x} (i x+\log (a))^{-1+k}\right. \\ \left.
+\frac{m \left((-i
   x+\log (a))^k+e^{2 i m x} (i x+\log (a))^k\right)}{k}\right) \log \left(\frac{1+\cosh (c x)}{\cos (b)+\cosh (c
   x)}\right) \, dx\\
=\frac{2 \left(\frac{1}{c}\right)^k e^{\frac{m (-b+\pi )}{c}} \pi }{k} \left(2^{1+k} e^{\frac{b
   m}{c}} \pi ^k \Phi \left(e^{\frac{2 m \pi }{c}},-k,\frac{\pi +c \log (a)}{2 \pi }\right)\right. \\ \left.
-(2 \pi )^k \left(\Phi
   \left(e^{\frac{2 m \pi }{c}},-k,\frac{-b+\pi +c \log (a)}{2 \pi }\right)+e^{\frac{2 b m}{c}} \Phi
   \left(e^{\frac{2 m \pi }{c}},-k,\frac{b+\pi +c \log (a)}{2 \pi }\right)\right)\right)
\end{multline}
\end{theorem}
\begin{theorem}
From Eq. (1.7.7.110) in \cite{obert}
\begin{multline}\label{eq_177110}
\int_0^{\infty } e^{-i m x} \left((-i x+\log (a))^{-1+k}+e^{2 i m x} (i x+\log (a))^{-1+k}\right. \\ \left.
+\frac{m \left((-i
   x+\log (a))^k+e^{2 i m x} (i x+\log (a))^k\right)}{k}\right) \log \left(1+\alpha ^2 \text{sech}^2(x)\right) \,
   dx\\
=\frac{\left(2 e^{\frac{m \pi }{2}} \pi ^{1+k} \left(\alpha +\sqrt{1+\alpha ^2}\right)^{-i m}\right) }{k}\left(\Phi
   \left(e^{m \pi },-k,\frac{\pi -2 i \sinh ^{-1}(\alpha )+2 \log (a)}{2 \pi }\right)\right. \\ \left.
+\left(\alpha +\sqrt{1+\alpha
   ^2}\right)^{2 i m} \Phi \left(e^{m \pi },-k,\frac{\pi +2 i \sinh ^{-1}(\alpha )+2 \log (a)}{2 \pi }\right)\right. \\ \left.
-2
   \left(\alpha +\sqrt{1+\alpha ^2}\right)^{i m} \Phi \left(e^{m \pi },-k,\frac{1}{2}+\frac{\log (a)}{\pi
   }\right)\right)
\end{multline}
\end{theorem}
\begin{theorem}
From Eq. (1.7.7.111) in \cite{obert}
\begin{multline}\label{eq_177111}
\int_0^{\infty } e^{-i m x} \left((-i x+\log (a))^{-1+k}+e^{2 i m x} (i x+\log (a))^{-1+k}\right. \\ \left.
+\frac{m \left((-i
   x+\log (a))^k+e^{2 i m x} (i x+\log (a))^k\right)}{k}\right) \log \left(1-(\alpha  \text{sech}(x))^2\right) \,
   dx\\
=\frac{\left(2 e^{m \cos ^{-1}(\alpha )} \pi ^{1+k}\right)}{k} \left(e^{2 m \sin ^{-1}(\alpha )} \Phi \left(e^{m \pi
   },-k,\frac{\pi -\cos ^{-1}(\alpha )+\log (a)}{\pi }\right)\right. \\ \left.
+\Phi \left(e^{m \pi },-k,\frac{\cos ^{-1}(\alpha )+\log
   (a)}{\pi }\right)-2 e^{m \sin ^{-1}(\alpha )} \Phi \left(e^{m \pi },-k,\frac{1}{2}+\frac{\log (a)}{\pi
   }\right)\right)
\end{multline}
\end{theorem}
\begin{theorem}
From Eq. (2.7.7.76) in \cite{obert}.
\begin{multline}\label{eq_27776}
\int_0^{\infty } e^{-i m x} \cot ^{-1}(\sinh (b x)) \left(k \left((-i x+\log (a))^{-1+k}-e^{2 i m x} (i x+\log
   (a))^{-1+k}\right)\right. \\ \left.
+m \left((-i x+\log (a))^k-e^{2 i m x} (i x+\log (a))^k\right)\right) \, dx\\
=-i \pi  \left(-2
   \left(\frac{1}{b}\right)^k e^{\frac{m \pi }{2 b}} \pi ^k \Phi \left(-e^{\frac{m \pi }{b}},-k,\frac{1}{2}+\frac{b
   \log (a)}{\pi }\right)+\log ^k(a)\right)
\end{multline}
\end{theorem}
\begin{theorem}
From Eq. (1.9.68) in \cite{bateman}.
\begin{multline}\label{eq_1968}
\int_0^{\infty } e^{-i m x} \cot ^{-1}\left(e^{2 x}\right) \left((-i
   x+\log (a))^{-2+k}+e^{2 i m x} (i x+\log (a))^{-2+k}\right. \\ \left.
+\frac{2 k m }{(-1+k) k}\left((-i
   x+\log (a))^{-1+k}+e^{2 i m x} (i x+\log
   (a))^{-1+k}\right)\right. \\ \left.
+\left(1+m^2\right) \left((-i x+\log (a))^k+e^{2 i m x} (i
   x+\log (a))^k\right)\right) \sinh (x) \, dx\\
=\frac{1}{2
   (-1+k) k}\left(-\sqrt{2}
   e^{\frac{m \pi }{4}} k \pi ^k \Phi \left(-e^{m \pi
   },1-k,\frac{1}{4}+\frac{\log (a)}{\pi }\right)\right. \\ \left.
+\sqrt{2} e^{\frac{3 m \pi
   }{4}} k \pi ^k \Phi \left(-e^{m \pi },1-k,\frac{3}{4}+\frac{\log (a)}{\pi
   }\right)\right. \\ \left.
+\pi  \left(\sqrt{2} e^{\frac{m \pi }{4}} \pi ^k \left(-\left((-1+m)
   \Phi \left(-e^{m \pi },-k,\frac{1}{4}+\frac{\log (a)}{\pi
   }\right)\right)\right.\right.\right. \\ \left.\left.\left.
+e^{\frac{m \pi }{2}} (1+m) \Phi \left(-e^{m \pi
   },-k,\frac{3}{4}+\frac{\log (a)}{\pi }\right)\right)-\log ^k(a)\right)\right)
\end{multline}
\end{theorem}
\begin{theorem}
From Eq. (4.294.10) in \cite{grad}.
\begin{multline}\label{eq_429410}
\int_0^{\infty } x^{-1+m} \tanh ^{-1}(b x) \log ^k(a x) \, dx\\
=\sum _{y=0}^{\infty } -i (-1)^k a^{-m} m^{-1-k}
   \pi  \left(\Gamma \left(1+k,-m \left(i (\pi +2 \pi  y)+\log \left(\frac{a}{b}\right)\right)\right)\right. \\ \left.
-\Gamma
   \left(1+k,m \left(-2 i \pi  (1+y)+\log \left(\frac{b}{a}\right)\right)\right)\right)
\end{multline}
\end{theorem}
\begin{theorem}
From Eq. (2.4.6.6) in \cite{prud1}.
\begin{multline}\label{eq_2466}
\int_0^{\infty } \frac{e^{-m x} (-x+\log (a))^k+e^{m x} (x+\log (a))^k}{b+\cosh (c x)} \, dx\\
=\frac{i^{1+k}
   c^{-1-k} e^{\frac{i m \left(\pi +i \cosh ^{-1}(b)\right)}{c}} (2 \pi )^{1+k} }{\sqrt{-1+b^2}}\left(\Phi \left(e^{\frac{2 i m \pi
   }{c}},-k,\frac{\pi +i \cosh ^{-1}(b)-i c \log (a)}{2 \pi }\right)\right. \\ \left.
-e^{\frac{2 m \cosh ^{-1}(b)}{c}} \Phi
   \left(e^{\frac{2 i m \pi }{c}},-k,-\frac{i \left(i \pi +\cosh ^{-1}(b)+c \log (a)\right)}{2 \pi
   }\right)\right)
\end{multline}
\end{theorem}
\begin{theorem}
From Eq. (2.2.7.11) in \cite{prud1}.
\begin{multline}\label{eq_22711}
\int_0^{\infty } \frac{x^{-1+\alpha } \left(b^m \log ^k(a b)-x^m \log ^k(a x)\right)}{b-x} \, dx\\
=-(2 i)^k
   b^{-1+m+\alpha } e^{i \pi  \alpha } \pi ^{1+k} \csc (\pi  \alpha ) \left(e^{2 i m \pi } \Phi \left(e^{2 i \pi  (m+\alpha
   )},-k,\frac{2 \pi -i \log (a)-i \log (b)}{2 \pi }\right)\right. \\ \left.
-\Phi \left(e^{2 i \pi  (m+\alpha )},-k,-\frac{i (\log (a)+\log
   (b))}{2 \pi }\right)\right)
\end{multline}
\end{theorem}
\begin{theorem}
From Eq. (2.2.7.13) in \cite{prud1}.
\begin{multline}\label{eq_22713}
\int_0^{\infty } \frac{x^{-1+m-\alpha +\beta } \left(-1+x^{2 \alpha }\right) \log ^k(a x)}{-1+x^{2 \beta }}
   \, dx\\
=\pi ^{1+k} \left(\frac{i}{\beta }\right)^{1+k} \left(e^{\frac{i \pi  (m-\alpha )}{\beta }} \Phi
   \left(-e^{\frac{i \pi  (m-\alpha )}{\beta }},-k,1-\frac{i \beta  \log (a)}{\pi }\right)\right. \\ \left.-e^{\frac{i \pi  (m+\alpha
   )}{\beta }} \Phi \left(-e^{\frac{i \pi  (m+\alpha )}{\beta }},-k,1-\frac{i \beta  \log (a)}{\pi
   }\right)\right)
\end{multline}
\end{theorem}
\subsubsection{Integrals of the form $\int_{0}^{\infty}x^{m-1}R(a,x)$ where $R(a,x)$ is a rational function}
\begin{theorem}
From Eq. (2.1.8.3) in \cite{brychkov}
\begin{multline}
\int_0^{\infty } \frac{x^{-1+m} \log ^k\left(\frac{\alpha  x}{\sqrt{a}+\sqrt{a+x}}\right)}{\sqrt{a+x} \left(\sqrt{a}+\sqrt{a+x}\right)^m} \, dx\\
=-2^{1+k+m} a^{\frac{1}{2} (-1+m)} e^{i m \pi
   } (i \pi )^{1+k} \Phi \left(e^{2 i m \pi },-k,\frac{1}{2}-\frac{i (\log (4)+\log (a)+2 \log (\alpha ))}{4 \pi }\right)
\end{multline}
\end{theorem}
\begin{theorem}
From Eq. (2.1.8.3) in \cite{brychkov}
\begin{multline}
\int_0^{\infty } \frac{\log ^k\left(\frac{\alpha }{-\sqrt{a}+\sqrt{a+x}}\right)}{x \sqrt{a+x}
   \left(-\sqrt{a}+\sqrt{a+x}\right)^m} \, dx\\
=2^{1+k-m} a^{-\frac{1}{2}-\frac{m}{2}} e^{i m \pi } (i \pi )^{1+k} \Phi
   \left(e^{2 i m \pi },-k,\frac{1}{2}+\frac{i (\log (4)+\log (a)-2 \log (\alpha ))}{4 \pi }\right)
\end{multline}
\end{theorem}
\begin{theorem}
From Eq. (2.1.8.7) in \cite{brychkov}
\begin{multline}
\int_0^{\infty } \frac{x^{-1+m} \log ^k\left(\frac{\alpha  x}{\sqrt{a+x+\sqrt{a (a+2 x)}}}\right)}{\sqrt{a+2 x} \left(a+x+\sqrt{a (a+2 x)}\right)^{m/2}} \, dx\\
=-2^{1+k+\frac{m}{2}}
   a^{\frac{1}{2} (-1+m)} e^{i m \pi } (i \pi )^{1+k} \Phi \left(e^{2 i m \pi },-k,\frac{1}{2}-\frac{i (\log (2)+\log (a)+2 \log (\alpha ))}{4 \pi }\right)
\end{multline}
\end{theorem}
\begin{theorem}
From Eq. (2.1.8.10) in \cite{brychkov}
\begin{multline}
\int_0^{\infty } \frac{x^{-1+m} \left(x+\sqrt{a^2+x^2}\right)^{m-1} \log ^k\left(\alpha  x \left(x+\sqrt{a^2+x^2}\right)\right)}{\sqrt{a^2+x^2}} \, dx\\
=-2^{1+k-m} a^{-2+2 m} e^{i m \pi } (i
   \pi )^{1+k} \Phi \left(e^{2 i m \pi },-k,\frac{1}{2}-\frac{i \left(2 \log (a)+\log \left(\frac{\alpha }{2}\right)\right)}{2 \pi }\right)
\end{multline}
\end{theorem}
\begin{theorem}
From Eq. (3.194.4(11)) in \cite{grad}. K\"{o}lbig form.
\begin{multline}\label{eq_3194411}
\int_0^{\infty } \frac{x^m \log ^k(a x)}{\left(1+\beta  x^q\right)^{n+1}} \, dx\\
=\frac{k! }{q^n}\sum _{t=0}^n \sum _{j=0}^n \frac{(-1)^{j+t} e^{\frac{i (1+m) \pi }{q}} (2 \pi
   )^{1+k-t} \left(-\frac{1}{q}\right)^t \left(\frac{i}{q}\right)^{-1+k-t} q^{-2+n} \beta ^{-\frac{1+m}{q}}  \left(-\frac{1+m-q}{q}\right)_{-j+n} S_j^{(t)}}{j! (-j+n)! (k-t)!}\\\Phi \left(e^{\frac{2 i (1+m) \pi }{q}},-k+t,\frac{\pi -i q \log (a)+i \log
   (\beta )}{2 \pi }\right)
\end{multline}
\end{theorem}
\begin{theorem}
From Eq. (2.1.9.16) in \cite{brychkov}
\begin{multline}\label{eq_21916}
\int_0^{\infty } \frac{x^{-1+m} \left(b+x+\sqrt{x (2 b+x)}\right)^{-\frac{1}{2}+m} \log ^k\left(a x
   \left(b+x+\sqrt{x (2 b+x)}\right)\right)}{\sqrt{2 b+x}} \, dx\\
-2^{\frac{3}{2}+2 k-m} b^{-1+2 m} e^{2 i m \pi } (i
   \pi )^{1+k} \Phi \left(e^{4 i m \pi },-k,\frac{2 \pi -i \log \left(\frac{a b^2}{2}\right)}{4 \pi }\right)
\end{multline}
\end{theorem}
\begin{theorem}
From Eq. (2.1.9.7)  in \cite{brychkov}
\begin{multline}
\int_0^{\infty } \frac{\sinh ^{m-1}(x) \log ^k\left(\frac{a \sinh (x)}{\sqrt{1+\cosh
   (x)}}\right)}{\sqrt{(1+\cosh (x))^m}} \, dx\\
=-i 2^{m/2} e^{\frac{1}{2} i (k+m) \pi } \pi ^{1+k} \Phi
   \left(e^{i m \pi },-k,-\frac{i (i \pi +\log (2)+2 \log (a))}{2 \pi }\right)
\end{multline}
\end{theorem}
\begin{theorem}
From Eq. (2.1.9.10)  in \cite{brychkov}
\begin{multline}
\int_0^{\infty } \frac{\log ^k\left(a x \left(x+\sqrt{b^2+x^2}\right)\right) \left(x^{-1+m}
   \left(x+\sqrt{b^2+x^2}\right)^{-1+m}\right)}{\sqrt{b^2+x^2}} \, dx\\
-2^{1+k-m} b^{-2+2 m} e^{i m \pi } (i \pi )^{1+k} \Phi
   \left(e^{2 i m \pi },-k,-\frac{i (i \pi +\log (a)-\log (2)+2 \log (b))}{2 \pi }\right)
\end{multline}
\end{theorem}
\begin{theorem}
From Eq. (2.1.9.10)  in \cite{brychkov}
\begin{multline}
\int_0^{\infty } \frac{\log ^k\left(\frac{a x}{-x+\sqrt{b^2+x^2}}\right) x^{-1+m}
   \left(-x+\sqrt{b^2+x^2}\right)^{1-m}}{\sqrt{b^2+x^2}} \, dx\\
=-2^{1+k-m} e^{i m \pi } (i \pi )^{1+k} \Phi \left(e^{2 i m \pi
   },-k,\frac{\pi +i \log (2)-i \log (a)}{2 \pi }\right)
\end{multline}
\end{theorem}
\begin{theorem}
From Eq. (1.2.2.38) in \cite{obert}
\begin{multline}\label{eq_12238}
\int_0^{\infty } x^{-1+m} \left(x+\sqrt{b^2+x^2}\right)^{-2+m} \log ^k\left(a x
   \left(x+\sqrt{b^2+x^2}\right)\right) \, dx\\
=2^{-1+k-m} b^{-2+2 m} e^{\frac{1}{2} i (k+2 m) \pi } \pi ^k \left(k
   \Phi \left(e^{2 i m \pi },1-k,-\frac{i \left(i \pi +\log \left(\frac{a}{2}\right)+\log \left(b^2\right)\right)}{2
   \pi }\right)\right. \\ \left.
+2 i (-2+m) \pi  \Phi \left(e^{2 i m \pi },-k,-\frac{i \left(i \pi +\log \left(\frac{a}{2}\right)+\log
   \left(b^2\right)\right)}{2 \pi }\right)\right)
\end{multline}
\end{theorem}
\begin{theorem}
From Eq.(2.2.4.1) in \cite{prud4}
\begin{multline}\label{eq_2241}
\int_0^{\infty } e^{b m x} \left(1-e^{-b x}\right)^m \log ^k\left(a e^{b x} \left(1-e^{-b x}\right)\right)
   \, dx=\frac{e^{i m \pi } (2 i \pi )^{1+k} \Phi \left(e^{2 i m \pi },-k,\frac{\pi -i \log (a)}{2 \pi
   }\right)}{b}
\end{multline}
\end{theorem}
\begin{theorem}
From Eq.(2.3.1.13) in \cite{prud4}
\begin{multline}\label{eq_23113}
\int_0^{\infty } \frac{\log ^k\left(\frac{a}{e^x (\cosh (b x)-1)^{1/b}}\right)}{\left(e^x (\cosh (b
   x)-1)^{1/b}\right)^m} \, dx\\
=-i i^k 2^{1+2 k+\frac{m}{b}} b^{-1-k} e^{\frac{2 i m \pi }{b}} \pi ^{1+k} \Phi
   \left(e^{\frac{4 i m \pi }{b}},-k,\frac{1}{2}-\frac{i (\log (2)+b \log (a))}{4 \pi }\right)
\end{multline}
\end{theorem}
\begin{theorem}
Fron Eq. (2.9.41) in \cite{bateman1}
\begin{multline}\label{eq_2941}
\int_0^{\infty } e^{-i m x} \left((-i x+\log (a))^{-2+k}-e^{2 i m x} (i x+\log (a))^{-2+k}\right. \\ \left.
+\frac{1}{(-1+k) k}\left(2 k m \left((-i x+\log(a))^{-1+k}-e^{2 i m x} (i x+\log (a))^{-1+k}\right)\right.\right. \\ \left.\left.
+\left(1+m^2\right) \left((-i x+\log (a))^k-e^{2 i m x} (i x+\log(a))^k\right)\right)\right) \log \left(1-e^{-2 x}\right) \sinh (x) \, dx\\
=-\frac{i}{(-1+k) k} \left(2 a^{-i-m} k \left(-a^{2 i}
   (i-m)^{-1-k} \Gamma (k,-((-i+m) \log (a)))\right.\right. \\ \left.\left.
-(-i-m)^{-1-k} \Gamma (k,-((i+m) \log (a)))\right)\right. \\ \left.
+\pi  \left(-2 e^{m \pi } \pi ^k \Phi \left(-e^{m \pi },-k,\frac{\pi +\log (a)}{\pi }\right)+\log ^k(a)\right)\right. \\ \left.
-2 i m \left(-(i-m)^{-1-k} (-((-i+m) \log (a)))^k+(-i-m)^{-1-k} (-((i+m) \log (a)))^k\right)\right)
\end{multline}
\end{theorem}
\begin{theorem}
From Eq. (2.2.4.8) in \cite{prud1}
\begin{multline}\label{eq_2248}
\int_0^z \frac{x^{m-1} }{\left(z^u-x^u\right)^{m/u}}
   \, \log ^k\left(\frac{a x}{\left(z^u-x^u\right)^{1/u}}\right)dx\\
   =-e^{\frac{i m \pi }{u}} \left(\frac{2 \pi  i}{u}\right)^{1+k} \Phi \left(e^{\frac{2 i m \pi
   }{u}},-k,\frac{\pi -i u \log (a)}{2 \pi }\right)
\end{multline}
\end{theorem}
\begin{theorem}
From Eq. (3.411.31) in \cite{grad}.
\begin{multline}\label{eq_341131}
\int_0^{\infty } \frac{-e^{(-m+p) x} (-x+\log (a))^k+e^{m x} (x+\log (a))^k}{-1+e^{p x}} \, dx\\
=e^{\frac{2 i m \pi }{p}} \left(\frac{i}{p}\right)^{1+k} (2 \pi )^{1+k} \Phi \left(e^{\frac{2 i m
   \pi }{p}},-k,1-\frac{i p \log (a)}{2 \pi }\right)+\frac{i \pi  \log ^k(a)}{p}
\end{multline}
\end{theorem}
\begin{theorem}
\begin{multline}
\int_0^{\infty } x^{-1+m} \log ^{-1+k}(a x) (k+m \log (a x)) \log \left(x^n \left(1+\frac{1}{x}\right)_n\right)
   \, dx\\
=-\sum _{s=1}^n e^{i m \pi } (2 i \pi )^{1+k} s^{-m} \Phi \left(e^{2 i m \pi },-k,\frac{\pi -i \log (a)+i \log
   (s)}{2 \pi }\right)
\end{multline}
where $Re(m)<0$
\end{theorem}
\begin{theorem}
\begin{multline}\label{eq_genroot_sec}
\int_0^1 x^n \left(-1+x^{-1-n}\right)^{-\frac{1}{2}+m} \log ^k\left(a \left(-1+x^{-1-n}\right)\right)
   \, dx\\
=\frac{e^{\frac{1}{2} i (k+2 m) \pi } }{1+n}\left(i k (2 \pi )^k \Phi \left(-e^{2 i m \pi },1-k,\frac{\pi
   -i \log (a)}{2 \pi }\right)\right. \\ \left.
+\pi  \left(-2^{1+k} m \pi ^k+(2 \pi )^k\right) \Phi \left(-e^{2 i m \pi
   },-k,\frac{\pi -i \log (a)}{2 \pi }\right)\right)
\end{multline}
\end{theorem}
\subsection{Theorems on multi-dimensional integrals}
\begin{theorem}
From Entry 23. in \cite{berndt2}
\begin{multline}\label{eq_23}
\int _0^{\infty }\int _0^{\infty }\frac{\log ^k\left(\frac{a y}{x^{1/r}}\right) x^{-\frac{m}{r}} y^{-1+m+\frac{r}{2}} \Gamma (b+x)}{\left(1+d y^r\right) \Gamma (1+b+n+x)}dydx\\
=\sum
   _{p=0}^{\infty } \frac{4^{1+k} d^{-\frac{1}{2}-\frac{m}{r}} e^{\frac{2 i m \pi }{r}} (b+p)^{-\frac{m}{r}} \pi ^{2+k} \left(\frac{i}{r}\right)^{-1+k} \Phi \left(e^{\frac{4 i m \pi
   }{r}},-k,\frac{2 \pi -i r \log (a)+i \log (d)+i \log (b+p)}{4 \pi }\right) (-n)_p}{r^2 p! \Gamma (1+n)}
\end{multline}
\end{theorem}
\begin{theorem}
From Eq. (3.26.12) in \cite{brychkov} and Eq. (3.326.2(10)) in \cite{grad}.
\begin{multline}\label{eq_32612}
\int _0^{\infty }\int _0^{\infty }\int _0^{\infty }\int _0^{\infty }\log ^k\left(\frac{a x y^2
   z^2}{t}\right) e^{-b \left(t+y^4+z^4\right)} t^{-m} x^{-1+m} y^{2+2 m} z^{2+2 m}\\ K\left(\frac{-x+\alpha }{2
   \alpha }\right) \theta (-x+\alpha )dtdzdydx\\
=\frac{i^{-1+k} 2^{-4+k-m} e^{i m \pi } \pi ^{\frac{5}{2}+k}
   \alpha ^m \Phi \left(e^{2 i m \pi },-k,-\frac{i \left(i \pi +\log \left(\frac{a}{2}\right)+\log (\alpha
   )\right)}{2 \pi }\right)}{b^{5/2}}
\end{multline}
\end{theorem}
\begin{theorem}
From Eq. (2.2.7.3) in \cite{prud1}
\begin{multline}
\int _{-d}^d\int _{-c}^c\frac{\log ^k\left(\frac{a (c+y) \sqrt{-1+\frac{2 d}{d+x}}}{c-y}\right) \left(-1+\frac{2 d}{d+x}\right)^{\frac{1}{2} (-2+m)}\label{eq_2273}
   \left(\frac{c+y}{c-y}\right)^m}{c^2+y^2}dydx\\
=\frac{i^{-1+k} d e^{i m \pi } (2 \pi )^{2+k} \Phi \left(e^{2 i m \pi },-k,\frac{\pi -i \log (a)}{2 \pi }\right)}{c}+\frac{(2 i)^{1+k}
   d e^{i m \pi } m \pi ^{2+k} \Phi \left(e^{2 i m \pi },-k,\frac{\pi -i \log (a)}{2 \pi }\right)}{c}\\
+\frac{(2 i)^k k d e^{i m \pi } \pi ^{1+k} \Phi \left(e^{2 i m \pi
   },1-k,\frac{\pi -i \log (a)}{2 \pi }\right)}{c}
\end{multline}
\end{theorem}
\begin{theorem}
From Eq. (3.251.1) and (3.241.1) in \cite{grad}.
\begin{multline}\label{eq_32511}
\int _0^1\int _0^{\infty }\frac{\log ^k\left(a \left(-1+x^{-1-n}\right) y^2\right) x^n \left(-1+x^{-1-n}\right)^m y^{2 m}}{1+y^2}dydx\\
=\frac{2^{-1+2 k} e^{2 i m \pi } (i \pi )^{1+k}
   \left(i k \Phi \left(e^{4 i m \pi },1-k,\frac{1}{2}-\frac{i \log (a)}{4 \pi }\right)-4 m \pi  \Phi \left(e^{4 i m \pi },-k,\frac{1}{2}-\frac{i \log (a)}{4 \pi
   }\right)\right)}{1+n}
\end{multline}
\end{theorem}
\begin{theorem}
From Eq. (3.1.3.9) and (8.334.3)  in \cite{prud1}
\begin{multline}\label{eq_3139}
\int_{\mathbb{R}^{4}_+}(r s)^{-\frac{m}{2}-1} (r+s)^{\frac{m+1}{2}} (x y)^{m/2} (x+y)^{\frac{1}{2} (-m-1)} e^{-p (r+x)-q (s+y)}\\
 \log ^k\left(\frac{a \sqrt{r+s}
   \sqrt{x y}}{\sqrt{r s} \sqrt{x+y}}\right)dxdydrds\\
   =\frac{2 i \pi ^{k+2} e^{\frac{1}{2} i \pi  (k+m)} \Phi \left(e^{i m \pi },-k,\frac{1}{2}-\frac{i \log
   (a)}{\pi }\right)}{p q}
\end{multline}
\end{theorem}
\begin{theorem}
From (3.241.4) in \cite{grad} and equation (3.9.1.3) in \cite{brychkov}. For all $k,a\in\mathbb{C},Re(v)>0,Re(p)>0,Re(m)<Re(u), -1/2<Re(m)<Re(v), Re(u)\leq 1$, then
\begin{multline}\label{eq_32414}
\int_{{\mathbb{R}^{4}_{+}}}y^m x^{u-1} t^{-m+\frac{u}{v}-1} z^{u \left(\frac{1}{v}-1\right)-m}
   D_u(t \alpha ) e^{-b \left(y+z^2\right)-\frac{1}{4} \alpha ^2 t^2}\\
   \left(p+q x^v\right)^{-m-1} \log ^k\left(\frac{a y}{t z \left(p+q
   x^v\right)}\right)dxdydzdt\\\\
   =\frac{1}{v}i^{k+1} \pi ^{k+\frac{3}{2}} p^{-m-1} e^{i \pi 
   \left(m-\frac{u}{v}\right)} \Gamma \left(\frac{u}{v}\right) b^{-\frac{(m+3)
   v+u (-v)+u}{2 v}} 2^{\frac{1}{2} \left(2 k+m-\frac{u}{v}+u\right)} \alpha
   ^{m-\frac{u}{v}}
    \left(\frac{p}{q}\right)^{u/v}\\
     \Phi \left(e^{2 i \pi 
   \left(m-\frac{u}{v}\right)},-k,-\frac{i (2 \log (a)-\log (b)-2 \log (p)+2
   \log (\alpha )+\log (2)+2 i \pi )}{4 \pi }\right)
\end{multline}
\end{theorem}
\begin{theorem}
From equation (3.18.5) in \cite{brychkov} and equation (3.326.2) in \cite{grad}. For all $k,a,b,c,\alpha\in\mathbb{C},1/2<Re(m)<1$,
\begin{dmath}\label{dilf}
\int_{0}^{\infty}\int_{0}^{\infty}\int_{0}^{\infty}x^{m-1} y^{2 m+4} z^{-m} \text{Ai}(x \alpha )^2 e^{-b y^6-c z} \log ^k\left(\frac{a x y^2}{z}\right)dxdydz\\\\
=-i \pi ^{k+\frac{1}{2}}
   3^{-\frac{m}{3}-\frac{11}{6}} b^{-\frac{m}{3}-\frac{5}{6}} c^{m-1} 2^{k-\frac{2 (m+1)}{3}} e^{\frac{1}{2} i \pi  (k+2 m)} \alpha ^{-m}\\
    \Phi
   \left(e^{2 i m \pi },-k,\frac{i (-3 \log (a)+\log (b)-3 \log (c)+3 \log (\alpha )+\log (12)-3 i \pi )}{6 \pi }\right)
\end{dmath}
\end{theorem}
\begin{theorem}
From Eq. (3.251.1) and (3.241.1) in \cite{grad}
\begin{multline}
\int _0^1\int _0^{\infty }\frac{\log ^k\left(a \left(-1+x^{-1-n}\right) y^2\right) x^n \left(-1+x^{-1-n}\right)^m y^{2 m}}{1+y^2}dydx\\
=\frac{2^{-1+2 k} e^{2 i m \pi } (i \pi )^{1+k}
   \left(i k \Phi \left(e^{4 i m \pi },1-k,\frac{1}{2}-\frac{i \log (a)}{4 \pi }\right)-4 m \pi  \Phi \left(e^{4 i m \pi },-k,\frac{1}{2}-\frac{i \log (a)}{4 \pi
   }\right)\right)}{1+n}
\end{multline}
\end{theorem}
\begin{theorem}
\begin{multline}\label{eq_beta_di}
\int _0^1\int _0^1\log ^k\left(a \left(-1+x^{-1-n}\right) \left(-1+y^{-1-n}\right)\right) x^n \left(-1+x^{-1-n}\right)^m y^n \left(-1+y^{-1-n}\right)^{\frac{1}{2}+m}dydx\\
=\frac{i
   4^{-1+k} e^{\frac{1}{2} i (k+4 m) \pi } \pi ^k }{(1+n)^2}\left(2 i k (1+4 m) \pi  \Phi \left(e^{4 i m \pi },1-k,\frac{1}{2}-\frac{i \log (a)}{4 \pi }\right)\right. \\ \left.
+(-1+k) k \Phi \left(e^{4 i m \pi
   },2-k,\frac{1}{2}-\frac{i \log (a)}{4 \pi }\right)-8 m (1+2 m) \pi ^2 \Phi \left(e^{4 i m \pi },-k,\frac{1}{2}-\frac{i \log (a)}{4 \pi }\right)\right)
\end{multline}
\end{theorem}
\begin{theorem}
From Eq. (3.197.10) in \cite{grad}.
\begin{multline}\label{eq_319710}
\int _0^1\int _0^1\frac{\log ^k\left(\frac{a x \left(-1+y^{-1-n}\right)}{1-x}\right)
   (1-x)^{-\frac{1}{2}-m} x^{-\frac{1}{2}+m} y^n \left(-1+y^{-1-n}\right)^m}{1+p x}dydx\\
=-\frac{i (4 i)^k e^{2 i m
   \pi } (1+p)^{-\frac{1}{2}-m} \pi ^{1+k}}{1+n} \left(-i k \Phi \left(e^{4 i m \pi },1-k,\frac{2 \pi -i \log (a)+i
   \log (1+p)}{4 \pi }\right)\right. \\ \left.
+4 m \pi  \Phi \left(e^{4 i m \pi },-k,\frac{2 \pi -i \log (a)+i \log (1+p)}{4 \pi
   }\right)\right)
\end{multline}
where $Im(n)\leq 0$.
\end{theorem}
\begin{theorem}
From Eq. (3.3.1.1) in \cite{brychkov} and Eq. (3.1.3.9) in \cite{prud1}
\begin{multline}\label{eq_3311_1}
\int _0^{\infty }\int _0^{\infty }\int _0^{\infty }\log ^k\left(\frac{a \sqrt{x y}}{\sqrt{z} \sqrt{x+y}}\right) e^{-p x-q y} (x y)^{m/2} (x+y)^{-\frac{1}{2}-\frac{m}{2}}
   z^{-1-\frac{m}{2}} \text{Ei}(-b z)dzdydx\\
=\sum _{y=0}^{\infty } \frac{4 i a^{-m} \pi ^{3/2} }{m
   \left(\sqrt{p}+\sqrt{q}\right) \sqrt{p q}}\left((-1)^k m^{-k} \Gamma (1+k)+\left(-k \Gamma (k)+\Gamma \left(1+k,-\frac{1}{2} m \left(i \pi
   +2 i \pi  y\right.\right.\right.\right. \\ \left.\left.\left.\left.
+2 \log (a)+\log (b)-2 \log \left(\sqrt{p}+\sqrt{q}\right)\right)\right)\right) \left(2 \log (a)+\log (b)\right.\right. \\ \left.\left.
+i \left(\pi +2 \pi  y+2 i \log
   \left(\sqrt{p}+\sqrt{q}\right)\right)\right)^k \left(-m \left(2 \log (a)+\log (b)+i \left(\pi +2 \pi  y+2 i \log \left(\sqrt{p}+\sqrt{q}\right)\right)\right)\right)^{-k}\right)
\end{multline}
\end{theorem}
\begin{theorem}
From Eq. (3.3.1.1) in \cite{brychkov} and Eq. (3.1.3.9) in \cite{prud1}
\begin{multline}\label{eq_3311_2}
\int _0^{\infty }\int _0^{\infty }\int _0^{\infty }e^{-p x-q y} (x y)^{m/2} (x+y)^{-\frac{1}{2}-\frac{m}{2}} z^{-1-\frac{m}{2}} \text{Ei}(-b z) \log ^{-1+k}\left(\frac{a \sqrt{x
   y}}{\sqrt{x+y} \sqrt{z}}\right)\\
 \left(k+m \log \left(\frac{a \sqrt{x y}}{\sqrt{x+y} \sqrt{z}}\right)\right)dzdydx\\
=\frac{4 i b^{m/2} e^{\frac{1}{2} i (k+m) \pi } \pi ^{\frac{3}{2}+k}
   \left(\sqrt{p}+\sqrt{q}\right)^{-1-m} \Phi \left(e^{i m \pi },-k,\frac{\pi -2 i \log (a)-i \log (b)+2 i \log \left(\sqrt{p}+\sqrt{q}\right)}{2 \pi }\right)}{\sqrt{p q}}
\end{multline}
\end{theorem}
\begin{theorem}
From Eq.  (3.368.2) and (3.387.1) in \cite{grad}.
\begin{multline}\label{eq_33682}
\int _0^{\frac{\pi }{4}}\int _0^{\frac{\pi }{4}}\int _0^{\frac{\pi }{4}}\log ^k\left(\frac{a \sin (2 x) \sin (2 y) \sin (z)}{\cos (2 x) (\cos (z)-\sin (z)) (\cos (y)+\sin (y))^2}\right)\\ \cos
   ^{-m}(2 x) \sec (x) \sec ^2(z) \sin ^{-\frac{1}{2}+m}(2 x) (\cos (y)+\sin (y))^{-2 m}\\
 \sin ^{-1+m}(2 y) (\cos (z)-\sin (z))^{\frac{1}{2}-m} \sin ^{-\frac{1}{2}+m}(z)dzdydx\\
=-\sqrt[4]{-1} 2^{2 k-m}
   e^{\frac{1}{2} i (k+3 m) \pi } \pi ^{2+k} \left(\Phi \left(e^{4 i m \pi },-k,\frac{3}{8}-\frac{i \log \left(\frac{a}{2}\right)}{4 \pi }\right)\right. \\ \left.
+i e^{i m \pi } \Phi \left(e^{4 i m \pi
   },-k,\frac{5}{8}-\frac{i \log \left(\frac{a}{2}\right)}{4 \pi }\right)\right)
\end{multline}
\end{theorem}
\begin{theorem}
\begin{multline}\label{eq_gr}
\int _0^{\frac{\pi }{2}}\int _0^1\int _0^{\frac{\pi }{2}}\int _0^1\log ^k\left(\frac{a \sin (x) \sqrt{\log (y)}}{\sin (z) \sqrt{\log (r)}}\right) \log ^{-\frac{m}{2}}(r) \log
   ^{\frac{m}{2}}(y) \sin ^m(x) \sin ^{-m}(z)drdzdydx\\
=\frac{1}{2} e^{\frac{1}{2} i (k+m) \pi } \pi ^{2+k} \Phi \left(-e^{i m \pi },-k,\frac{1}{2}-\frac{i \log (a)}{\pi }\right)
\end{multline}
\end{theorem}
\begin{theorem}
From Eq. (3.192.4) in \cite{grad} and Eq. (3.2.1.2) in \cite{brychkov}.
\begin{multline}\label{eq_31924_1}
\int _0^{\infty }\int _1^{\infty }\frac{\log ^k(a x (-1+y)) (-1+y)^{-\frac{1}{2}} \text{Li}_n(-b x)}{x y}dydx\\
=-\frac{i k! (-1)^n i^{k+n} 2^{2+2 (k+n)} \pi ^{2+k+n} \zeta \left(-k-n,\frac{2 \pi
   -i \log (a)+i \log (b)}{4 \pi }\right)}{(k+n)!}
\end{multline}
\end{theorem}
\begin{theorem}
From Eq. (3.18.5.1) in \cite{brychkov} and Eq. (3.326.2) in \cite{grad}.
\begin{multline}\label{eq_31851}
\int _0^{\infty }\int _0^{\infty }\int _0^{\infty }\int _0^{\infty }\log ^k\left(\frac{a x y^2}{z t^2}\right) e^{-d t^6-c y^6} t^{-2 (-3+m)} x^{-1+m} y^{2 (2+m)} z^{-m} \text{Ai}(b x)^2
   \text{Ai}(f z)^2dtdzdydx\\
=-\frac{1}{81} \left(i 2^{-3+k} b^{-m} c^{-\frac{5}{6}-\frac{m}{3}} d^{-\frac{7}{6}+\frac{m}{3}} e^{i m \pi } f^{-1+m} (i \pi )^k\right. \\ \left.
 \Phi \left(e^{2 i m \pi },-k,\frac{3 \pi -3
   i \log (a)+3 i \log (b)+i \log (c)-i \log (d)-3 i \log (f)}{6 \pi }\right)\right)
\end{multline}
\end{theorem}
\begin{theorem}
\begin{multline}\label{eq_gendi}
\int _0^1\int _0^1\frac{x^{-1-m} y^s \left(-1+y^{-1-s}\right)^{-\frac{1}{2}+\frac{m}{q}}}{\left(b+x^q\right) \left(1+b
   x^q\right)} \left(x^q
   \left(1+b x^q\right) \log ^k\left(\frac{a \left(-1+y^{-1-s}\right)^{1/q}}{x}\right)\right. \\ \left.
+x^{2 m} \left(b+x^q\right)
   \log ^k\left(a x \left(-1+y^{-1-s}\right)^{1/q}\right)\right)dydx\\
=\frac{\left((4 i)^k b^{-\frac{m}{q}} e^{\frac{2 i m \pi }{q}} \pi ^{1+k} q^{-2-k}\right)}{1+s}
   \left(k q \Phi \left(e^{\frac{4 i m \pi }{q}},1-k,\frac{1}{2}-\frac{i q \left(\log (a)+\log
   \left(b^{-1/q}\right)\right)}{4 \pi }\right)\right. \\ \left.
-2 i \pi  (-2 m+q) \Phi \left(e^{\frac{4 i m \pi
   }{q}},-k,\frac{1}{2}-\frac{i q \left(\log (a)+\log \left(b^{-1/q}\right)\right)}{4 \pi
   }\right)\right)
\end{multline}
\end{theorem}
\begin{theorem}
From Eq. (6.4.14) in \cite{bateman1}
\begin{multline}\label{eq_6414}
\int _0^1\int _0^1\int _0^1z^s \left(-1+z^{-1-s}\right)^m \log ^{-\frac{1}{2}+m}(x) \log ^{-\frac{1}{2}-m}(y)
   \log ^k\left(\frac{a \left(-1+z^{-1-s}\right) \log (x)}{\log (y)}\right)dzdydx\\
=\frac{(4 i)^k e^{2 i m \pi } \pi
   ^{1+k} \left(k \Phi \left(e^{4 i m \pi },1-k,\frac{1}{2}-\frac{i \log (a)}{4 \pi }\right)+4 i m \pi  \Phi
   \left(e^{4 i m \pi },-k,\frac{1}{2}-\frac{i \log (a)}{4 \pi }\right)\right)}{1+s}
\end{multline}
\end{theorem}
\subsection{Theorems involving integrals of series}
\begin{theorem}
\begin{multline}\label{eq_iis}
\int_0^1 \left(\sum _{s=1}^n s^m \left(-1+\frac{1}{x}\right)^m \log ^k\left(a s
   \left(-1+\frac{1}{x}\right)\right)\right) \, dx\\
=\sum _{c=1}^n \left(-k c^m e^{i m \pi } (2 i \pi )^k \Phi
   \left(e^{2 i m \pi },1-k,\frac{\pi -i \log (a)-i \log (c)}{2 \pi }\right)\right. \\ \left.
-m c^m e^{i m \pi } (2 i \pi
   )^{1+k} \Phi \left(e^{2 i m \pi },-k,\frac{\pi -i \log (a)-i \log (c)}{2 \pi }\right)\right)
\end{multline}
\end{theorem}
\section*{Theorems involving series}
\begin{theorem}
\begin{multline}
\sum _{j=0}^{\infty } \sum _{l=0}^{\infty } \frac{\left(-\frac{1}{a}\right)^{j+l} \Gamma (j-m+2) \Gamma (l+m+1)
   (1-k)_{j+l-1}}{(j+l+1) (j+l+2) \Gamma (j+1) \Gamma (l+1)}\\
=\frac{\pi  e^a (m-1) m a^{-k} \csc (\pi  m) \Gamma (k+1,a)}{2
   k}
\end{multline}
\end{theorem}
\begin{theorem}
From Eq.  (5.7.1.15) in \cite{prud2} and Eq. (6.3.3.1) in \cite{prud3}.
\begin{multline}
\sum_{j=0}^{\infty}\sum_{l=0}^{\infty}\sum_{q=0}^{\infty}\frac{E_q(\alpha ) 2^{-2 j-2 l+q} \left(-\frac{1}{a}\right)^{2 j+2 l+q} (1-k)_{2 j+2 l+q}}{\Gamma (l+1) \Gamma
   (q+1) \Gamma (2 j+l+2)}\\
=\frac{a^{1-k} \left(2^{2 k+1} \left(\zeta \left(-k,\frac{1}{4} (a+2 \alpha +1)\right)-\zeta
   \left(-k,\frac{1}{4} (a+2 \alpha +3)\right)\right)-(a+2 \alpha -1)^k\right)}{k}
\end{multline}
\end{theorem}
\begin{theorem}
From Eq. (1.4) in \cite{yu} and Eq. (8.440) in \cite{grad}.
\begin{multline}
\sum _{l=0}^{\infty } \sum _{n=0}^{\infty } \sum _{j=0}^{\lfloor n/2\rfloor} \frac{(\alpha +n) 2^{-\alpha -2 j-2 l} (-1)^{-j+l+n} \left(\frac{1}{t}\right)^{-2 j+n+1} \Gamma (-j+n+1)
   \binom{-\alpha }{n-j} a^{-\alpha +k-2 l-n}}{\alpha  \Gamma (j+1) \Gamma (l+1) \Gamma (l+n+\alpha +1) (k)_{-2
   l-n-\alpha +1}}\\
=\frac{2^{-\alpha } e^{a t} \Gamma (k) t^{\alpha -k-1} \Gamma (k-\alpha +1,a t)}{\Gamma (\alpha +1)
   \Gamma (k-\alpha +1)}
\end{multline}
\end{theorem}
\begin{theorem}
From Ex.3 on page 230 in \cite{bromwich}.
\begin{multline}
\sum_{n=-\infty}^{\infty}\frac{a^{-m+n} e^{i n \beta } \Gamma (1+k,(-m+n) \log (a))}{(n-m)^{k+1}}\\
=i e^{\frac{1}{2} i (k \pi +2 m \beta
   )} (2 \pi )^{1+k} \Phi \left(e^{2 i m \pi },-k,\frac{\beta -i \log (a)}{2 \pi }\right)
\end{multline}
\end{theorem}
\begin{theorem}
From Eq. (1) in \cite{chatterjea} and \cite{plos}.
\begin{equation}
\sum _{n=0}^{\infty } \frac{t^n L_n^a(x)}{(a+n+1)!}=t^{-a-1} \int_0^t e^s s^a (s x)^{-a/2} J_a\left(2 \sqrt{s
   x}\right) \, ds
\end{equation}
\end{theorem}
\begin{theorem}
From Eq.  (12.2.11) in \cite{buchholz} and \cite{plos}.
\begin{multline}
\sum _{l=0}^{\infty } \frac{(-s)^l \Gamma (l+2 v) L_l^u(z)}{(l+u+1)!}=\int_0^s \frac{r^u (r+1)^{-2 v} s^{-u-1}
   \Gamma (2 v) \, _1F_1\left(2 v;u+1;\frac{r z}{r+1}\right)}{\Gamma (u+1)} \, dr
\end{multline}
\end{theorem}
\begin{theorem}
From Eq.  (12.2.18a) in \cite{buchholz} and \cite{plos}.
\begin{multline}
\sum _{l=0}^{\infty } \frac{(a s)^l L_n^{l-n}(s)}{(l+1)!}\\
=-\frac{e^s \left(-(a-1)^n s^n (s-a s)^{-n} \Gamma
   (n+1,s-a s)+e^{i \pi  n} \Gamma (n+1,s)\right)}{a s n!}
\end{multline}
\end{theorem}
\begin{theorem}
From Eq. (1) in \cite{chatterjea}.
\begin{equation}
\sum _{n=0}^{\infty } \frac{(-t)^n (1-k)_{n-1} L_n^{\alpha }(x)}{(\alpha
   +n)!}=-\frac{\left(\frac{1}{t}+1\right)^k t^k \, _1F_1\left(-k;\alpha +1;\frac{t x}{t+1}\right)}{k \Gamma (\alpha
   +1)}
\end{equation}
\end{theorem}
\begin{theorem}
From Eq. (1) in \cite{chatterjea} and \cite{plos}.
\begin{multline}
\sum _{n=0}^{\infty } \frac{(-t)^n (1-k)_{n-1} L_n^{\alpha }(x)}{(\alpha +n+1)!}\\
=-t^{-\alpha -1} \int_0^t
   \frac{(-r-1)^k (-r)^{-k} r^{\alpha +k} \, _1F_1\left(-k;\alpha +1;\frac{r x}{r+1}\right)}{k \Gamma (\alpha +1)} \,
   dr
\end{multline}
\end{theorem}
\begin{theorem}
From Eq. (3) in \cite{chatterjea}.
\begin{multline}
\sum _{n=0}^{\infty } \frac{t^n C_n^{\left(p+\frac{1}{2}\right)}(x)}{(n+2 p+1)!}=\int_0^t \frac{2^p s^{2 p}
   t^{-2 p-1} \Gamma (p+1) e^{s x} \left(s \sqrt{1-x^2}\right)^{-p} J_p\left(s \sqrt{1-x^2}\right)}{\Gamma (2 p+1)} \,
   ds
\end{multline}
\end{theorem}
\begin{theorem}
From Eq. (5.2.13.29) in \cite{prud1} and \cite{plos}.
\begin{equation}
\sum _{k=0}^{\infty } \frac{\left((y-1) y^{-v}\right)^k \Gamma (u+k v)}{(k+1)! \Gamma (v
   k-k+u)}=\frac{y^u-y^v}{(y-1) (u-v)}
\end{equation}
\end{theorem}
\begin{theorem}
From Eq. (5.2.13.30) in \cite{prud1} and \cite{plos}.
\begin{multline}
\sum _{k=0}^{\infty } \frac{\left((y-1) y^{-v}\right)^k \Gamma (u+k v-1)}{(k+1)! \Gamma (v
   k-k+u)}\\
   =\frac{y^v}{(1-y) (u-v-1) (u-v)}+\frac{y^u \left(\frac{v}{-u y+v y+y}+\frac{v-1}{u-v}\right)}{(u-1)
   (1-y)}
\end{multline}
\end{theorem}
\begin{theorem}
From Eq. (5.13.22) in \cite{hansen} and \cite{plos}.
\begin{equation}
\sum _{k=0}^{\infty } \frac{(k+2)^k \left(e^{-y} y\right)^{k+1}}{(k+2)!}=\frac{1}{2} \left(-e^y y+2
   e^y-2\right)
\end{equation}
\end{theorem}
\begin{theorem}
From Ex. 28 on page 472 in \cite{bromwich} and \cite{dieckmann}.
\begin{multline}
\sum_{n=-\infty}^{\infty}e^{\alpha  n} \Phi \left(-e^{2 (m+n \alpha )},-k,\frac{1}{2} (\log (a)+1)\right)\\
=\sum_{n=-\infty}^{\infty}\frac{\pi  2^{-k-2} e^{m
   \left(-1-\frac{2 i \pi  n}{\alpha }\right)} \text{sech}\left(\frac{\pi ^2 n}{\alpha }\right) \left(e^{\frac{4 i
   \pi  m n}{\alpha }} \left(\log (a)+\frac{2 i \pi  n}{\alpha }\right)^k+\left(\log (a)-\frac{2 i \pi  n}{\alpha
   }\right)^k\right)}{\alpha }
\end{multline}
\end{theorem}
\begin{theorem}
From Eq. (4.2.1.7) in \cite{prud2}.
\begin{multline}
\sum_{p=1}^{n}\sum_{j=0}^{\infty}\frac{(-1)^{j+p} (2 p+v) (2 a)^{-2 j-2 p-v}}{\Gamma (1+j) \Gamma (1-2 j+k-2 p-v) \Gamma (1+j+2
   p+v)}\\
=\sum_{j=0}^{\infty}\frac{\left((-1)^j 2^{-2-2 j-v} a^{-2-2 j+k-v}\right) }{a^k
   \Gamma (1+j)}\\
   \left(-\frac{1}{\Gamma (-1-2 j+k-v) \Gamma
   (2+j+v)}\right. \\ \left.
   +\frac{\left(-\frac{1}{4}\right)^n a^{-2 n}}{\Gamma (-1-2 j+k-2 n-v) \Gamma (2+j+2 n+v)}\right)
\end{multline}
\end{theorem}
\begin{theorem}
From Eq. (4.2.1.16) in \cite{prud2}.
\begin{multline}
\sum _{j=0}^{\infty } \sum _{p=0}^n \frac{(-1)^{j+p} 2^{-2 j-m+2 p} \binom{m-n}{p} a^{-2 j+k-m+2 p}}{j! (n-p)!
   \Gamma (j+m-p+1) (-2 j+k-m+2 p)!}\\
=\frac{(-1)^n 2^{2 n-m} a^{k-m+2 n} \,
   _2F_1\left(-\frac{k}{2}+\frac{m}{2}-n,-\frac{k}{2}+\frac{m}{2}-n+\frac{1}{2};m-2 n+1;-\frac{1}{a^2}\right)}{\Gamma
   (n+1) \Gamma (m-2 n+1) \Gamma (k-m+2 n+1)}
\end{multline}
\end{theorem}
\begin{theorem}
From Eq. (2) in \cite{andrews}.
\begin{multline}
\sum_{p=0}^{\infty}(-1)^p \left(\frac{a^p \Gamma (1+k,(-m+p) \log (a))}{(p-m)^{k+1}}+\frac{a^{-p} \Gamma (1+k,-((m+p) \log
   (a)))}{(-m-p)^{k+1}}\right)\\
=a^m (-1)^m (2 i \pi )^{1+k} \Phi \left(e^{2 i m \pi },-k,\frac{\pi -i \log (a)}{2 \pi
   }\right)-\frac{\Gamma (1+k,-m \log (a))}{m (-m)^k}
\end{multline}
\end{theorem}
\begin{theorem}
From Eq. (3.31) in \cite{berndt}.
\begin{multline}
\sum_{p=1}^{\infty}\alpha  e^{2 \alpha  m p} \left(\coth \left(\alpha ^2 p\right)-1\right) \left((\log (a)+2 \alpha  p)^k-e^{-4
   \alpha  m p} (\log (a)-2 \alpha  p)^k\right)\\
+i\sum_{p=1}^{\infty} \beta  e^{-2 i \beta  m p} \left(\coth \left(\beta ^2
   p\right)-1\right) \left((\log (a)-2 i \beta  p)^k-e^{4 i \beta  m p} (\log (a)+2 i \beta 
   p)^k\right)\\
=-\frac{1}{\log (a)}\left
(2^{k+1} \log (a) \alpha ^{k+1} e^{2 \alpha  m} \Phi \left(e^{2 m \alpha },-k,\frac{\log
   (a)}{2 \alpha }+1\right)\right. \\ \left.
-2^{k+1} \log (a) (i \beta )^{k+1} e^{2 i \beta  m} \Phi \left(e^{2 i m \beta
   },-k,1-\frac{i \log (a)}{2 \beta }\right)\right. \\ \left.
+\log ^k(a) (2 k+\log (a) (\alpha -i \beta +2 m))\right)
\end{multline}
\end{theorem}
\begin{theorem}
From Eq. (6.1.75) in \cite{hansen}.
\begin{multline}
\sum _{p=0}^{\infty }\frac{\left((-1)^p e^{-2 i a p x}\right)}{\left(m^2+x^2 (1+2
   p)^2\right)^{k+1} (1+2 p)^2}\\
 \left(-e^{2 i a (1+2 p) x} (k-1) (-i m+x (1+2 p)) (i m+x+2 p x)^2 \right. \\ \left.(-m-i x (1+2 p))^k \Gamma (k-1,-a (m-i x (1+2 p)))\right. \\ \left.
-e^{2 i a (1+2 p) x} m
   (i m+2 x (1+2 p)) (-m-i (1+2 p) x)^{1+k} \Gamma (k,-a m+i a (1+2 p) x)\right. \\ \left.
-x^2 (1+2 p)^2 (i m+x (1+2 p)) (-m+i (1+2 p) x)^k \Gamma (k,-a (m+i x (1+2 p)))\right)\\
=\frac{e^{a (m+i x)} }{200 k}\left(-200 e^{\frac{m \pi }{2 x}} \pi ^{1+k} \left(\frac{1}{x}\right)^{-1+k} \Phi \left(-e^{\frac{m \pi }{x}},-k,\frac{1}{2}+\frac{a
   x}{\pi }\right)\right. \\ \left.
+\frac{a^{-1+k}}{m^3} \left(-200 i k m \left(C m^2+(-3 C+\pi ) x^2\right)\right.\right. \\ \left.\left.
+\frac{x}{m^2+25 x^2} \left(64 k m^2 x (-i m+5 x) \, _4F_3\left(\frac{5}{2},\frac{5}{2},3,\frac{5}{2}+\frac{i m}{2
   x};\frac{7}{2},\frac{7}{2},\frac{7}{2}+\frac{i m}{2 x};-1\right)\right.\right.\right. \\ \left.\left.\left.
-320 k m^2 x^2 \, _4F_3\left(\frac{5}{2},3,\frac{5}{2}-\frac{i m}{2 x}\frac{5}{2}+\frac{i m}{2
   x};\frac{7}{2},\frac{7}{2}-\frac{i m}{2 x},\frac{7}{2}+\frac{i m}{2 x};-1\right)\right.\right.\right. \\ \left.\left.\left.
+25 \left(m^2+25 x^2\right) \left(-\left(\left(2 a m^3+k (m+i x) (m+3 i x)\right) H_{-\frac{3}{4}-\frac{i
   m}{4 x}}\right)\right.\right.\right.\right. \\ \left.\left.\left.\left.
+\left(2 a m^3+k (m+i x) (m+3 i x)\right) H_{-\frac{i m+x}{4 x}}+\left(-2 a m^3+k (m-i x) (m-3 i x)\right) \right.\right.\right.\right. \\ \left.\left.\left.\left.
\left(-\psi ^{(0)}\left(\frac{3}{4}+\frac{i m}{4 x}\right)+\psi
   ^{(0)}\left(\frac{i m+x}{4 x}\right)\right)\right)\right)\right)\right)
\end{multline}
\end{theorem}
\begin{theorem}
From Eq. (6.3.4.3) in \cite{prud3}.
\begin{multline}\label{eq:thm}
\sum_{p=0}^{\infty}\frac{(-1)^p \left(\frac{i}{\pi  a}\right)^{2 p} (1-k)_{2 p} E_{1+2 p}(x)}{\Gamma (2+2
   p)}\\
   =-\frac{a \pi }{k}\left(\frac{2}{a \pi }\right)^k  \left(\zeta \left(-k,\frac{1}{2} (1+a \pi -x)\right)-\zeta
   \left(-k,\frac{1}{2} (2+a \pi -x)\right)\right. \\ \left.
   -\zeta \left(-k,\frac{1}{2} (a \pi +x)\right)+\zeta \left(-k,\frac{1}{2}
   (1+a \pi +x)\right)\right)
\end{multline}
\end{theorem}
\begin{theorem}
From Exercise. (81) in \cite{sharma}.
\begin{multline}
\sum_{p=1}^{n}(i (p+2))^k e^{2 i m (p+2)} \Phi \left(e^{2 i m (p+2)},-k,1-\frac{i \log (a)}{2 (p+2)}\right)\\
-(i (p+1))^k e^{2 i
   m (p+1)} \Phi \left(e^{2 i m (p+1)},-k,1-\frac{i \log (a)}{2 (p+1)}\right)\\
=(i (n+2))^k e^{2 i m (n+2)} \Phi
   \left(e^{2 i m (n+2)},-k,1-\frac{i \log (a)}{2 (n+2)}\right)\\
-(2 i)^k e^{4 i m} \Phi \left(e^{4 i m},-k,1-\frac{1}{4}
   i \log (a)\right)
\end{multline}
\end{theorem}
\begin{theorem}
\begin{multline}
\sum _{p=0}^{2n}  -1)^p \left(i 2^p\right)^k e^{i m 2^p} \Phi \left(e^{i 2^{p+1} m},-k,\frac{1}{2} \left(1-i 2^{-p} \log
   (a)\right)\right)\\
-\sum _{p=0}^{n}2 \left(i 4^p\right)^k e^{i m 4^p} \Phi \left(e^{i 2^{2 p+1} m},-k,\frac{1}{2} \left(1-i 4^{-p}
   \log (a)\right)\right)\\
=\left(i 4^n\right)^k e^{i m 2^{2 n+1}} \Phi \left(e^{i 2^{2 n+1} m},-k,1-i 2^{-2 n-1} \log
   (a)\right)\\
   -\left(\frac{i}{2}\right)^k e^{i m} \Phi \left(e^{i m},-k,1-i \log (a)\right)
\end{multline}
\end{theorem}
\begin{theorem}
From page 170 exercise 9 in \cite{hind}.
\begin{multline}\label{hind}
\sum _{p=0}^n \left(i 2^p\right)^k e^{i m 2^{p+1}} \left(\Phi \left(-e^{i 2^{p+1} m},-k,1-i 2^{-p-1} \log
   (a)\right)+\Phi \left(e^{i 2^{p+1} m},-k,1-i 2^{-p-1} \log (a)\right)\right)\\
=-2 \left(\left(i 2^{n+1}\right)^k
   e^{i m 2^{n+2}} \Phi \left(e^{i 2^{n+2} m},-k,1-i 2^{-n-2} \log (a)\right)\right. \\ \left.
-i^k e^{2 i m} \Phi \left(e^{2 i
   m},-k,1-\frac{1}{2} i \log (a)\right)\right)
\end{multline}
\end{theorem}
\begin{theorem}
From Eq. (2.4.6.13) in \cite{prud1} and Eq. (3.511.4) in \cite{grad}.
\begin{multline}
\sum _{p=0}^{n-1} (-1)^p e^{-\frac{i m p \pi }{n \alpha }} \left(e^{\frac{i m (1-2 n+2 p) \pi }{n \alpha }} \Phi
   \left(e^{-\frac{2 i m \pi }{\alpha }},-k,1-\frac{1+2 p}{4 n}+\frac{a \alpha }{2}\right)\right. \\ \left.
-\Phi \left(e^{-\frac{2 i m
   \pi }{\alpha }},-k,\frac{1+2 p+2 a n \alpha }{4 n}\right)-e^{\frac{2 i m \pi }{\alpha }} \Phi \left(e^{\frac{2 i m
   \pi }{\alpha }},-k,1-\frac{1+2 p}{4 n}+\frac{a \alpha }{2}\right)\right. \\ \left.
+e^{\frac{i m (1+2 p) \pi }{n \alpha }} \Phi
   \left(e^{\frac{2 i m \pi }{\alpha }},-k,\frac{1+2 p+2 a n \alpha }{4 n}\right)\right)\\
=\left(\frac{1}{2 n}\right)^k
   \left(-\Phi \left(-e^{-\frac{i m \pi }{n \alpha }},-k,\frac{1}{2}+a n \alpha \right)+e^{\frac{i m \pi }{n \alpha }}
   \Phi \left(-e^{\frac{i m \pi }{n \alpha }},-k,\frac{1}{2}+a n \alpha \right)\right)
\end{multline}
\end{theorem}
\begin{theorem}
From Eq.  (4.4.7.14) in \cite{prud1}.
\begin{multline}
\sum _{p=0}^{n-1} e^{i 2^{1-p} m} \left(\left(i 2^{-p}\right)^k \Phi \left(-e^{i 2^{1-p} m},-k,1-i 2^{-1+p} \log (a)\right)\right. \\ \left.
-\left(i 2^{1-p}\right)^k e^{i 2^{1-p} m} \Phi
   \left(-e^{i 2^{2-p} m},-k,1-i 2^{-2+p} \log (a)\right)\right)\\
=-(2 i)^k e^{4 i m} \Phi \left(-e^{4 i m},-k,1-\frac{1}{4} i \log (a)\right)\\+\left(i 2^{1-n}\right)^k e^{i
   2^{2-n} m} \Phi \left(-e^{i 2^{2-n} m},-k,1-i 2^{-2+n} \log (a)\right)
\end{multline}
\end{theorem}
%
%%\begin{theorem}
%From Eq. 
%\begin{multline}
%\
%\end{multline}
%\end{theorem}
%
\begin{theorem}
Multidimensional sum involving Bernoulli and Euler numbers. From Eqs. (50.5.7) and (51.3.2) in \cite{hansen}. 
\begin{multline}
\sum\limits_{x_{1},...,x_{j}\geq 0}\frac{\left(2^{\sum _{l=1}^{j-1} l x_{2+l}} \left(-1+2^{-1+x_1}\right)
   \pi ^{x_1+\sum _{l=1}^j x_{1+l}} B_{x_1} (a \pi )^{k-x_1-\sum _l^{1 j}
   x_{1+l}}\right) \prod _{l=1}^j \frac{E_{x_{1+l}}}{\Gamma
   \left(1+x_{1+l}\right)}}{\Gamma \left(1+x_1\right) (k)_{1-x_1-\sum _{l=1}^j
   x_{1+l}}}\\
=\frac{1}{2} \left(2^{1+j}\right)^k \pi ^k \zeta
   \left(1-k,\frac{1}{2} \left(1+2^{-j} a\right)\right)
\end{multline}
\end{theorem}
\begin{theorem}
From pg. 127 in \cite{bedford}.
\begin{multline}
\sum _{n=1}^{\infty } \left((-t)^n (1-k)_{n-1}-\sum _{j=l}^n p^{-j+n} q^j
   (-t)^j \binom{n}{j} (1-k)_{j-1}\right)\\
=-e^{1/t} t^k \Gamma
   \left(k,\frac{1}{t}\right)+\frac{e^{1/t} \left(\frac{q}{1-p}\right)^l t^k
   \Gamma (k) \Gamma \left(1+k-l,\frac{1}{t}\right)}{\Gamma (1+k-l)}
\end{multline}
\end{theorem}
\begin{theorem}
Euler-Hermite generating function. From Eq. (6.3.4) in \cite{prud3}.
\begin{multline}
\sum _{j=0}^{\infty } \sum _{q=0}^{\infty } \frac{(-1)^j a^{-2 j-q} E_{2 j+1}(\alpha ) H_q(\beta )}{\Gamma
   (2+2 j) \Gamma (-2 j+k-q) \Gamma (1+q)}\\
=-\sum _{n=0}^{\infty } \sum _{p=0}^{\infty } \frac{(-1)^p i^{1+k-n-2
   p} 2^{k-2 p} a^{1-k} \beta ^n}{n! (k-n-2 p)! p!} \left(-\zeta \left(-k+n+2 p,\frac{1}{2} (1-i a-\alpha )\right)\right. \\ \left.
+\zeta
   \left(-k+n+2 p,1-\frac{i a}{2}-\frac{\alpha }{2}\right)+\zeta \left(-k+n+2 p,\frac{1}{2} (-i a+\alpha
   )\right)\right. \\ \left.
-\zeta \left(-k+n+2 p,\frac{1}{2} (1-i a+\alpha )\right)\right)
\end{multline}
\end{theorem}
\begin{theorem}
From Eq. (\ref{ober_17720})
\begin{multline}
\sum _{p=0}^n \frac{1}{2} \left(i 2^{-p}\right)^{1+k} e^{-i 2^{-p} (b+m)} \left(-e^{i 2^{1-p} b} \Phi
   \left(-e^{i 2^{-p} (b-m)},-k,1-i 2^p a\right)\right. \\ \left.
-e^{i 2^{1-p} m} \Phi \left(-e^{i 2^{-p} (-b+m)},-k,1-i 2^p
   a\right)+\Phi \left(-e^{-i 2^{-p} (b+m)},-k,1-i 2^p a\right)\right. \\ \left.
+e^{i 2^{1-p} (b+m)} \Phi \left(-e^{i 2^{-p}
   (b+m)},-k,1-i 2^p a\right)\right)\\
=i 2^{-1-n} e^{-i 2^{-n} \left(1+2^{1+n}\right) (b+m)} \left(i^k 2^{1+k+n} e^{i
   \left(4 m+2^{-n} (b+m)\right)} \Phi \left(e^{-2 i (b-m)},-k,1-\frac{i a}{2}\right)\right. \\ \left.
+i^k 2^{1+k+n} e^{i \left(4b+2^{-n} (b+m)\right)} \Phi \left(e^{2 i (b-m)},-k,1-\frac{i a}{2}\right)\right. \\ \left.
-\left(i 2^{-n}\right)^k e^{2 i\left(b+m+2^{-n} m\right)} \Phi \left(e^{-i 2^{-n} (b-m)},-k,1-i 2^n a\right)\right. \\ \left.
-\left(i 2^{-n}\right)^k e^{2 i\left(b+2^{-n} b+m\right)} \Phi \left(e^{i 2^{-n} (b-m)},-k,1-i 2^n a\right)\right. \\ \left.
-i^k 2^{1+k+n} e^{i 2^{-n} (b+m)} \Phi
   \left(e^{-2 i (b+m)},-k,1-\frac{i a}{2}\right)\right. \\ \left.
-i^k 2^{1+k+n} e^{i 2^{-n} \left(1+2^{2+n}\right) (b+m)} \Phi
   \left(e^{2 i (b+m)},-k,1-\frac{i a}{2}\right)\right. \\ \left.
+\left(i 2^{-n}\right)^k e^{2 i (b+m)} \Phi \left(e^{-i 2^{-n}
   (b+m)},-k,1-i 2^n a\right)\right. \\ \left.+\left(i 2^{-n}\right)^k e^{i 2^{1-n} \left(1+2^n\right) (b+m)} \Phi \left(e^{i 2^{-n}
   (b+m)},-k,1-i 2^n a\right)\right)
\end{multline}
\end{theorem}
\begin{theorem}
From Eq. (22.1.3) in \cite{hansen} and (5.3.8.4) in \cite{prud1}.
\begin{multline}
\sum _{j=1}^{\infty } 2^{-j-k} \left(-2 \left(2^{-j}\right)^k e^{2^{-1-j} m} \left(-\Phi \left(-e^{2^{-1-j}
   m},-k,1+2^{1+j} \log (a)\right)\right.\right. \\ \left.\left.
+2^k e^{2^{-1-j} m} \left(2^k e^{2^{-j} m} \Phi \left(-e^{2^{1-j} m},-k,1+2^{-1+j}
   \log (a)\right)\right.\right.\right. \\ \left.\left.\left.
-\Phi \left(-e^{2^{-j} m},-k,1+2^j \log (a)\right)\right)\right)-2^k \log
   ^k(a)\right)\\
=-\left(\frac{5 a^{-m} (-m)^{-k} \Gamma (1+k,-m \log (a))}{m}+2^{1-k} e^{m/2} \Phi
   \left(e^{m/2},-k,1+2 \log (a)\right)\right. \\ \left.
-2 e^m \left(-\Phi \left(e^m,-k,1+\log (a)\right)+2^k e^m \Phi \left(e^{2
   m},-k,1+\frac{\log (a)}{2}\right)\right)+\log ^k(a)\right)
\end{multline}
\end{theorem}
\begin{theorem}
From Eq. (\ref{eq_sech})
\begin{multline}
\sum _{p=0}^n 2^p \left(i 2^{-p}\right)^{1+k} e^{-i 2^{-1-p} m \pi } \left(-\Phi
   \left(e^{-i 2^{-p} m \pi },-k,\frac{1}{2}-\frac{i 2^p \log (a)}{\pi }\right)\right. \\ \left.
+e^{i 2^{-p} m \pi } \Phi
   \left(e^{i 2^{-p} m \pi },-k,\frac{1}{2}-\frac{i 2^p \log (a)}{\pi }\right)\right)\\
=i e^{-\frac{1}{2} i
   \left(4+\left(2+2^{-n}\right) m\right) \pi }  \left(e^{i 2^{-1-n} \left(2^n k+m\right) \pi }
   \left(\Phi \left(e^{-i m \pi },-k,1-\frac{i \log (a)}{\pi }\right)\right.\right. \\ \left.\left.
-e^{2 i m \pi } \Phi \left(e^{i m \pi
   },-k,1-\frac{i \log (a)}{\pi }\right)\right)+\left(i 2^{-1-n}\right)^k e^{i m \pi } \left(-\Phi
   \left(e^{-i 2^{-1-n} m \pi },-k,1-\frac{i 2^{1+n} \log (a)}{\pi }\right)\right.\right. \\ \left.\left.
+e^{i 2^{-n} m \pi } \Phi
   \left(e^{i 2^{-1-n} m \pi },-k,1-\frac{i 2^{1+n} \log (a)}{\pi }\right)\right)\right)
\end{multline}
\end{theorem}
\begin{theorem}
From Eqs. (\ref{eq_sech2}) and  (\ref{eq_csch2}).
\begin{multline}
\sum _{p=0}^{n-1} \left(i 2^{-p}\right)^k e^{-i 2^{-1-p} m \pi } \pi ^k \left(-2^p k \Phi \left(e^{-i 2^{-p} m
   \pi },1-k,\frac{1}{2}-\frac{i 2^p \log (a)}{\pi }\right)\right. \\ \left.
+i m \pi  \Phi \left(e^{-i 2^{-p} m \pi
   },-k,\frac{1}{2}-\frac{i 2^p \log (a)}{\pi }\right)+e^{i 2^{-p} m \pi } \left(2^p k \Phi \left(e^{i 2^{-p} m \pi
   },1-k,\frac{1}{2}-\frac{i 2^p \log (a)}{\pi }\right)\right.\right. \\ \left.\left.
+i m \pi  \Phi \left(e^{i 2^{-p} m \pi
   },-k,\frac{1}{2}-\frac{i 2^p \log (a)}{\pi }\right)\right)\right)\\
=e^{-i m \pi } (i \pi )^k \left(k \Phi
   \left(e^{-i m \pi },1-k,1-\frac{i \log (a)}{\pi }\right)-i m \pi  \Phi \left(e^{-i m \pi },-k,1-\frac{i \log
   (a)}{\pi }\right)\right. \\ \left.
-e^{2 i m \pi } \left(k \Phi \left(e^{i m \pi },1-k,1-\frac{i \log (a)}{\pi }\right)+i m \pi 
   \Phi \left(e^{i m \pi },-k,1-\frac{i \log (a)}{\pi }\right)\right)\right)\\
+\left(i 2^{-n}\right)^k e^{-i 2^{-n} m
   \pi } \pi ^k \left(-2^n k \Phi \left(e^{-i 2^{-n} m \pi },1-k,1-\frac{i 2^n \log (a)}{\pi }\right)\right. \\ \left.
+i m \pi  \Phi
   \left(e^{-i 2^{-n} m \pi },-k,1-\frac{i 2^n \log (a)}{\pi }\right)+e^{i 2^{1-n} m \pi } \left(2^n k \Phi
   \left(e^{i 2^{-n} m \pi },1-k,1-\frac{i 2^n \log (a)}{\pi }\right)\right.\right. \\ \left.\left.
+i m \pi  \Phi \left(e^{i 2^{-n} m \pi
   },-k,1-\frac{i 2^n \log (a)}{\pi }\right)\right)\right)
\end{multline}
\end{theorem}
\begin{theorem}
From Eqs. (\ref{eq_24614}) and (41.12.1) in \cite{hansen}.
\begin{multline}
-\sum _{p=0}^{n-1} \left(-i 2^{-p}\right)^k \left(\zeta \left(-k,\frac{\pi +3 i 2^p \log (a)}{6 \pi
   }\right)+\zeta \left(-k,\frac{1}{6} \left(5+\frac{3 i 2^p \log (a)}{\pi
   }\right)\right)\right)\\
=e^{-\frac{1}{2} i k \pi } \left(\zeta \left(-k,\frac{1}{3}+\frac{i \log (a)}{2
   \pi }\right)+\zeta \left(-k,\frac{2}{3}+\frac{i \log (a)}{2 \pi }\right)\right)\\
-\left(-i
   2^{-n}\right)^k \left(\zeta \left(-k,\frac{1}{6} \left(2+\frac{3 i 2^n \log (a)}{\pi
   }\right)\right)+\zeta \left(-k,\frac{1}{6} \left(4+\frac{3 i 2^n \log (a)}{\pi
   }\right)\right)\right)
\end{multline}
\end{theorem}
\begin{theorem}
From Eqs. (\ref{eq_24614}) and (41.12.2) in \cite{hansen}.
\begin{multline}
\sum _{p=1}^n \left(-2^{1+p}\right)^k \left(\zeta \left(-k,\frac{1}{6} \left(1-\frac{3\times 2^{-p} \log
   (a)}{\pi }\right)\right)+\zeta \left(-k,\frac{1}{6} \left(5-\frac{3\times 2^{-p} \log (a)}{\pi
   }\right)\right)\right)\\
=(-2)^k \left(\zeta \left(-k,\frac{1}{3}-\frac{\log (a)}{2 \pi }\right)+\zeta
   \left(-k,\frac{2}{3}-\frac{\log (a)}{2 \pi }\right)\right)\\
-\left(-2^{1+n}\right)^k \left(\zeta
   \left(-k,\frac{1}{6} \left(2-\frac{3\times 2^{-n} \log (a)}{\pi }\right)\right)+\zeta \left(-k,\frac{1}{6}
   \left(4-\frac{3\times 2^{-n} \log (a)}{\pi }\right)\right)\right)
\end{multline}
\end{theorem}
\begin{theorem}
From Eq. (5.2.9.5) in \cite{prud1}.
\begin{equation}
\sum _{k=1}^{\infty } \frac{k^k (a \exp (-a))^k}{(k+2)!}=-\frac{-1+2 a+2 a^2-4 a e^a+e^{2 a}}{4 a^2}
\end{equation}
\end{theorem}
\begin{theorem}
\begin{equation}
\sum _{k=1}^{\infty } \frac{k^k (a \exp (-a))^k}{(k+2)!}=-\frac{-1+2 a+2 a^2-4 a e^a+e^{2 a}}{4 a^2}
\end{equation}
\end{theorem}
%
%\begin{theorem}
%From Eq. 
%\begin{multline}
%\
%\end{multline}
%\end{theorem}
%%
\begin{theorem}
From Eqs. (5.24.14) in \cite{hansen} and (6.3.4.2) in \cite{prud3}.
\begin{multline}
\sum _{n=0}^{\infty } \sum _{p=0}^{\infty } \frac{(1+n)^{-1+n} E_{2 p}(x) (-n+a)^{k-n-2 p} (-1)^{1+n+2 p}
   (1-k)_{-1+n+2 p}}{n! (2 p)!}\\
=\frac{2^k }{k}\left(\zeta \left(-k,\frac{1}{2} (2-x+a)\right)-\zeta \left(-k,\frac{1}{2}
   (3-x+a)\right)\right. \\ \left.
+\zeta \left(-k,\frac{1}{2} (1+x+a)\right)-\zeta \left(-k,\frac{1}{2}
   (2+x+a)\right)\right)
\end{multline}
\end{theorem}
\begin{theorem}
From Eq. (5.2.11.3) in \cite{prud1}.
\begin{equation}
\sum _{k=0}^{\infty } \frac{\Gamma (1+k+m) x^k}{\Gamma (2+k)}=\frac{\left(-1+(1-x)^{-m}\right) \Gamma
   (m)}{x}
\end{equation}
\end{theorem}
\begin{theorem}
\begin{equation}
\sum _{k=0}^{\infty } \frac{x^k}{(k+1)! \Gamma (v-k+1)}=\frac{-1+(1+x)^{1+v}}{x \Gamma (2+v)}
\end{equation}
\end{theorem}
\begin{theorem}
Bernoulli generating function. From Eq. (6.3.2) in \cite{prud3}.
\begin{multline}
\sum _{k=0}^{\infty } \frac{\left(e^{-u y} u\right)^k B_k(x+k
   y)}{(k+1)!}\\
=\frac{e^{u x} \left(-\left(e^u\right)^{-x+y} u
   B_{e^u}(x-y,0)+\Phi \left(e^u,2,x-y\right)\right)-e^{u y} \psi
   ^{(1)}(x-y)}{u}
\end{multline}
\end{theorem}
\begin{theorem}
Kummer confluent hypergeometric generating function. From Eq. (6.6.1.3) in \cite{prud3}.
\begin{equation}
\sum _{k=0}^{\infty } \frac{(b-a)_k t^k \, _1F_1(a;b+k;x)}{(k+1)!
   (b)_k}=\int_0^t \exp (s) \, _1F_1(a;b;x-s) \frac{\, ds}{t}
\end{equation}
\end{theorem}
\begin{theorem}
From Eq. (5.2.9.9) in \cite{prud1}.
\begin{multline}
\sum _{n=1}^{\infty } \frac{(a-n)^{k-n} (-n)^n (1-k)_{n-1}}{\Gamma (n+4)}\\
=-\frac{-18 a^k+\frac{-4 a^{3+k}+4
   (3+a)^{3+k}+3 \left(5 a^{2+k}-9 (2+a)^{2+k}\right) (3+k)+18 \left(-a^{1+k}+3 (1+a)^{1+k}\right) (2+k)
   (3+k)}{(1+k) (2+k) (3+k)}}{108 k}
\end{multline}
\end{theorem}
\begin{theorem}
StruveH and Euler's polynomials. From Eq. (6.4.3.3) in \cite{prud3}.
\begin{multline}
\sum _{j=1}^{\infty } \sum _{n=0}^{\infty } \sum _{p=0}^{\infty }
   \frac{(-1)^{j+n} 2^{-2 n+p} a^{-2 n-p} j^{2 n} E_p(\alpha )}{\Gamma
   \left(\frac{3}{2}+n\right) \Gamma (1+p) (k)_{-2 n-p-v} \Gamma (k) \Gamma
   \left(\frac{3}{2}+n+v\right)}\\
=\frac{2^{-1+2 k-2 v} a^{1-k+v} \left(-\zeta
   \left(1-k+v,\frac{1}{4} (a+2 \alpha )\right)+\zeta \left(1-k+v,\frac{1}{4}
   (2+a+2 \alpha )\right)\right)}{\sqrt{\pi } \Gamma (k-v) \Gamma
   \left(\frac{3}{2}+v\right)}
\end{multline}
\end{theorem}
\begin{theorem}
From Eq. (6.7.1.7) in \cite{prud3}.
\begin{multline}
\sum _{j=0}^{\infty } \sum _{n=0}^{\infty } \frac{a^{-j-n} \Gamma
   (j+n+\alpha ) \Gamma (n+\beta ) \Gamma (j+\gamma )}{\Gamma (1+j) (k)_{1-j-n}
   \Gamma (1+n) \Gamma (c+j+n) \Gamma (\beta ) \Gamma (\gamma )}\\
=\sum
   _{m=0}^{\infty } \frac{a^{-m} \Gamma (m+\alpha ) \Gamma (m+\beta +\gamma
   )}{(k)_{1-m} \Gamma (1+m) \Gamma (c+m) \Gamma (\beta +\gamma )}
\end{multline}
\end{theorem}
\begin{theorem}
From Eq. (6.4.3.20) in \cite{prud3}.
\begin{multline}
\sum _{j=0}^{\infty } \sum _{m=0}^{\infty } \frac{(-1)^j t^{2 j+m} \Gamma
   (c+j-\alpha ) \Gamma (j+m+\alpha ) \Gamma (2 j+m+\beta ) \Gamma (j+\beta
   -\gamma ) \Gamma (j+\gamma )}{\Gamma (1+j) \Gamma (c+j) \Gamma (1-2 j+k-m)
   \Gamma (1+m) \Gamma (c+2 j+m) \Gamma (2 j+\beta )}\\
=\sum _{p=0}^{\infty } \sum
   _{n=0}^{\infty } \frac{t^{n+p} \Gamma (n+\alpha ) \Gamma (p+\alpha ) \Gamma
   (p+\beta -\gamma ) \Gamma (n+\gamma ) \Gamma (c-\alpha )}{\Gamma (1+n) \Gamma
   (c+n) \Gamma (1+k-n-p) \Gamma (1+p) \Gamma (c+p) \Gamma (\alpha )}
\end{multline}
\end{theorem}
\begin{theorem}
From Eq. (6.9.2.1) in \cite{prud3}.
\begin{multline}
\sum _{j=0}^{\infty } \sum _{l=0}^{\infty } \frac{s^l t^j (-j)_l (-j+\alpha )_l (1-\beta )_j (\gamma )_l}{\Gamma
   (1+j) \Gamma (1+k-l) \Gamma (1+l) (1-\alpha )_j (-j+\beta )_l (\delta )_l}\\
=\sum _{n=0}^{\infty } \sum _{m=0}^{\infty
   } \frac{s^n (-1)^n t^{n+m} (1-\beta )_m (\gamma )_n}{\Gamma (1+m) \Gamma (1+k-n) \Gamma (1+n) (1-\alpha )_m (\delta
   )_n}
\end{multline}
\end{theorem}
\begin{theorem}
From Eq. (489) pp. 92. in \cite{jolley}. For $n$ even.
\begin{multline}
\sum _{r=1}^{\frac{n-2}{2}} \left(e^{\frac{2 i \pi  m r}{n}} \Phi \left(e^{2 i m \pi
   },-k,\frac{r}{n}+\frac{z}{2}\right)+(-1)^{2 m} e^{-\frac{2 i \pi  m r}{n}} \Phi \left(e^{2 i m \pi },-k,\frac{1}{2}
   \left(-\frac{2 r}{n}+z+2\right)\right)\right)\\
=n^{-k} \Phi \left(e^{\frac{2 i m \pi }{n}},-k,\frac{n z}{2}\right)-2^{-k} \Phi
   \left(e^{i m \pi },-k,z\right)
\end{multline}
\end{theorem}
\begin{theorem}
From Eq. (490) pp. 92. in \cite{jolley}. For $n$ odd.
\begin{multline}
\sum _{r=1}^{\frac{n-1}{2}} e^{-\frac{i \pi  m r}{n}} \left(\Phi
   \left(e^{i m \pi },-k,\frac{1}{2} \left(a-b-\frac{2
   r}{n}+1\right)\right)\right. \\ \left.
+e^{\frac{2 i \pi  m r}{n}} \Phi \left(e^{i m \pi
   },-k,\frac{1}{2} \left(a-b+\frac{2 r}{n}+1\right)\right)\right)\\
=(-1)^m \Phi
   \left(e^{i m \pi },-k,\frac{1}{2} (a-b+1)\right)-n^{-k} i^{m
   \left(\frac{1}{n}+1\right)} \Phi \left(e^{\frac{i m \pi }{n}},-k,\frac{1}{2}
   ((a-b) n+1)\right)
\end{multline}
\end{theorem}
\begin{theorem}
From Eq. (120) on pp. 71 in \cite{edwards}. For any positive integer $n$.
\begin{multline}
\sum _{r=0}^{n-1} e^{-\frac{i \pi  m r}{n}} \left((-1)^m \Phi \left(e^{i m \pi },-k,a-\frac{2 \pi  r+t}{2 n \pi
   }+1\right)+e^{\frac{i m (2 \pi  r+t)}{n}} \Phi \left(e^{i m \pi },-k,a+\frac{2 \pi  r+t}{2 n \pi
   }\right)\right)\\
=n^{-k} \left(e^{\frac{i \pi  m}{n}} \Phi \left(e^{\frac{i m \pi }{n}},-k,a n-\frac{t}{2 \pi
   }+1\right)+e^{\frac{i m t}{n}} \Phi \left(e^{\frac{i m \pi }{n}},-k,a n+\frac{t}{2 \pi }\right)\right)
\end{multline}
\end{theorem}
\begin{theorem}
From (41.2.34) in \cite{hansen} and Eq. (2.4.6.13) in \cite{prud1}.
\begin{multline}
\sum _{p=1}^{n-1} (-1)^p e^{-\frac{m p \pi }{n}} \left(-e^{\frac{2 m (-n+p) \pi }{n}} \Phi \left(e^{2 m \pi
   },-k,\frac{a n+p}{2 n}\right)\right. \\ \left.
+\Phi \left(e^{2 m \pi },-k,\frac{1}{2} \left(2+a-\frac{p}{n}\right)\right)\right)
   \sin \left(\frac{p \pi }{n}\right)\\
=i 2^{-1-k} e^{\frac{(-i+m-2 m n) \pi }{n}} n^{-k} \left(\Phi
   \left(-e^{\frac{(-i+m) \pi }{n}},-k,1+a n\right)-e^{\frac{2 i \pi }{n}} \Phi \left(-e^{\frac{(i+m) \pi }{n}},-k,1+a
   n\right)\right)
\end{multline}
\end{theorem}
\begin{theorem}
From Eq. (41.2.27) in \cite{hansen} and Eq. (2.4.6.13) in \cite{prud1}.
\begin{multline}
\sum _{p=0}^{n-1} (-1)^p e^{\frac{m p \pi }{n}} \cos \left(\frac{j (1+2 p) \pi }{2 n}\right) \left(\Phi
   \left(e^{2 m \pi },-k,\frac{1+2 a n+2 p}{4 n}\right)\right. \\ \left.
-e^{\frac{m (-1+2 n-2 p) \pi }{n}} \Phi \left(e^{2 m \pi
   },-k,1+\frac{a}{2}-\frac{1+2 p}{4 n}\right)\right)\\
=\frac{1}{2} e^{-\frac{m \pi }{2 n}} n^{-k} \left(e^{\frac{(-i
   j+m) \pi }{2 n}} \Phi \left(e^{\frac{2 (-i j+m) \pi }{n}},-k,\frac{1}{4} (1+2 a n)\right)\right. \\ \left.
-e^{\frac{3 (-i j+m) \pi
   }{2 n}} \Phi \left(e^{\frac{2 (-i j+m) \pi }{n}},-k,\frac{1}{4} (3+2 a n)\right)+e^{\frac{(i j+m) \pi }{2 n}}
   \left(\Phi \left(e^{\frac{2 (i j+m) \pi }{n}},-k,\frac{1}{4} (1+2 a n)\right)\right.\right. \\ \left.\left.
-e^{\frac{(i j+m) \pi }{n}} \Phi
   \left(e^{\frac{2 (i j+m) \pi }{n}},-k,\frac{1}{4} (3+2 a n)\right)\right)\right)
\end{multline}
\end{theorem}
\begin{theorem}
From Eq. (37.2.4) in \cite{hansen} and Eq. (4.3.10) in \cite{polyanin}.
\begin{multline}
\sum _{p=1}^n \frac{1}{2} \left(i 2^{-p}\right)^k e^{i 2^{-p} m \pi } \pi ^{1+k} \left(2 i \Phi \left(-e^{i 2^{-p} m \pi },-k,1-\frac{i 2^p \log (a)}{\pi }\right)\right. \\ \left.
+\Phi \left(-i e^{i
   2^{-p} m \pi },-k,1-\frac{i 2^p \log (a)}{\pi }\right)-\Phi \left(i e^{i 2^{-p} m \pi },-k,1-\frac{i 2^p \log (a)}{\pi }\right)\right)\\
=-i \pi  \left(e^{i m \pi } (i \pi )^k \Phi
   \left(-e^{i m \pi },-k,1-\frac{i \log (a)}{\pi }\right)\right. \\ \left.
-\left(i 2^{-n}\right)^k e^{i 2^{-n} m \pi } \pi ^k \Phi \left(-e^{i 2^{-n} m \pi },-k,1-\frac{i 2^n \log (a)}{\pi
   }\right)\right)
\end{multline}
\end{theorem}
\begin{theorem}
From Eq. (6.10.6) in \cite{hansen}.
\begin{multline}
\sum _{j=0}^{\infty } \frac{(b+j \alpha )^j (j \alpha +a)^{k-j} (1-k)_{j-1}}{\Gamma (j+1)}=-\frac{(-1)^{-k} e^{\frac{-a+b}{\alpha }} \alpha ^k \Gamma \left(1+k,\frac{b-a}{\alpha
   }\right)}{k}
\end{multline}
\end{theorem}
\begin{theorem}
\begin{multline}
\sum _{j=0}^{\infty } \frac{(b+j)^j s^j (a+j s)^{k-j} (1-k)_{j-1}}{\Gamma (j+1)}=-\frac{(-1)^{-k} e^{b-\frac{a}{s}} s^k \Gamma \left(1+k,b-\frac{a}{s}\right)}{k}
\end{multline}
\end{theorem}
\begin{theorem}
Euler-Neumann Generating Functions.
\begin{multline}
\sum _{n=0}^{\infty } \sum _{p=0}^{\infty } \sum _{q=0}^{\infty } \frac{(-1)^n 2^p (n \alpha +\beta )^{-1+n} E_p(x) (n \alpha +a)^{-1+k-n-p-2 q} (-1)^{2+n+p+2 q} (1-k)_{n+p+2 q}}{n!
   p! \Gamma (2+2 q)}\\
=\frac{1}{k \beta }\left(-(-1+a+2 x-\beta )^k+2^{1+2 k} \left(-\zeta \left(-k,1+\frac{1}{4} \log \left(e^{-1+a+2 x-\beta }\right)\right)\right.\right. \\ \left.\left.
+\zeta \left(-k,\frac{1}{4} \left(2+\log
   \left(e^{-1+a+2 x-\beta }\right)\right)\right)\right)\right)
\end{multline}
\end{theorem}
\begin{theorem}
From Eq. (5.2.11.17) in \cite{prud1}.
\begin{multline}
\sum _{p=0}^{\infty } \sum _{q=0}^{\infty }
   \frac{\left(-\frac{1}{a}\right)^{p+q} \left(\Gamma (1+p-v) \Gamma (q+v)
   (1-k)_{-1+p+q}\right)}{\Gamma (1+p) \Gamma (1+q)}=-\frac{e^a \pi  \csc (\pi 
   v) \Gamma (1+k,a)}{a^k k}
\end{multline}
\end{theorem}
\begin{theorem}
\begin{equation}
\sum _{k=0}^{\infty } \frac{(k+2)^k (y \exp (-y))^{k+1}}{(k+2)!}=-\frac{1}{2} e^y \left(-2+2 e^{-y}+y\right)
\end{equation}
\end{theorem}
\begin{theorem}
\begin{multline}
\sum _{k=0}^{\infty } \frac{a^k E_k(x)}{(k+1)!}=\frac{2 e^{a x} \, _2F_1\left(1,x;1+x;-e^a\right)+x \left(\psi ^{(0)}\left(\frac{x}{2}\right)-\psi
   ^{(0)}\left(\frac{1+x}{2}\right)\right)}{a x}
\end{multline}
\end{theorem}
\begin{theorem}
From Eq. (21) in \cite{brafman}.
\begin{multline}
\sum _{n=0}^{\infty } \frac{(1+a)_n \, _2F_1\left(-n,1+a+b+n;1+a;\frac{1-x}{2}\right) z^n}{(n+1)!}\\
=\int_0^z \frac{2^{a+b} \left(1-s+\sqrt{1+s^2-2 s x}\right)^{-a}
   \left(1+s+\sqrt{1+s^2-2 s x}\right)^{-b}}{\sqrt{1+s^2-2 s x} z} \, ds
\end{multline}
\end{theorem}
\begin{theorem}
From Eq. (5.2.9.6) in \cite{prud1}.
\begin{equation}
\sum _{j=1}^{\infty } \frac{(1-k)_{j-1} (-j)^{1+j} (t-j)^{-j+k}}{\Gamma
   (j+2)}=\frac{e^t \Gamma (2+k,t)-(1+t)^{1+k}}{k (1+k)}
\end{equation}
\end{theorem}
\begin{theorem}
From Eq. (17.3.14) in \cite{hansen}.
\begin{multline}
\sum _{p=0}^{\infty } (-1)^p a^{-2 i p \alpha } \cos (x+2 p x) \left(-(-i
   (m+\alpha +2 p \alpha ))^{-1-k} \Gamma (1+k,-i (m+\alpha +2 p \alpha ) \log
   (a))\right. \\ \left.
+a^{2 i (1+2 p) \alpha } (-i (m-(1+2 p) \alpha ))^{-1-k} \Gamma (1+k,-i
   (m-(1+2 p) \alpha ) \log (a))\right)\\
=-\frac{1}{2} i a^{i (m+\alpha )}
   e^{\frac{i m (\pi -2 x)}{2 \alpha }} \pi ^{1+k} \left(\frac{1}{\alpha
   }\right)^{1+k} \left(\Phi \left(-e^{\frac{i m \pi }{\alpha }},-k,\frac{\pi -2
   x+2 \alpha  \log (a)}{2 \pi }\right)\right. \\ \left.
+e^{\frac{2 i m x}{\alpha }} \Phi
   \left(-e^{\frac{i m \pi }{\alpha }},-k,\frac{\pi +2 x+2 \alpha  \log (a)}{2
   \pi }\right)\right)
\end{multline}
\end{theorem}
%%%
\begin{theorem}
From Exercise vii (f) (p.135) no. 20 in \cite{durell}.
\begin{multline}
-\sum _{p=1}^n  \left(i 3^p\right)^k e^{2 i 3^{-1+p} m} \left(\Phi \left(e^{2 i 3^p m},-k,\frac{1}{3}-\frac{1}{2} i 3^{-p} \log (a)\right)\right . \\ \left.
+e^{2 i 3^{-1+p} m} \Phi \left(e^{2
   i 3^p m},-k,\frac{2}{3}-\frac{1}{2} i 3^{-p} \log (a)\right)\right)\\
= \left(-i^k e^{2 i m} \Phi \left(e^{2 i m},-k,1-\frac{1}{2} i \log (a)\right)\right. \\ \left.
+\left(i 3^n\right)^k e^{2 i 3^n
   m} \Phi \left(e^{2 i 3^n m},-k,1-\frac{1}{2} i 3^{-n} \log (a)\right)\right)
\end{multline}
\end{theorem}
\begin{theorem}
\begin{multline}
\sum _{p=1}^n \left(i 3^{-p}\right)^{1+k} \left(\sqrt[3]{-1} e^{i 3^{-p} \pi  r} \Phi \left(-\sqrt[3]{-1} e^{i
   3^{-p} \pi  r},-k,1+3^p a\right)\right. \\ \left.
+e^{\pi  \left(-\frac{i}{3}+i 3^{-p} r\right)} \Phi \left((-1)^{2/3} e^{i 3^{-p}
   \pi  r},-k,1+3^p a\right)\right)\\
=-i e^{\frac{1}{2} i \pi  (k+2 r)} \Phi \left(e^{i \pi  r},-k,1+a\right)+\left(i
   3^{-n}\right)^{1+k} e^{i 3^{-n} \pi  r} \Phi \left(e^{i 3^{-n} \pi  r},-k,1+3^n a\right)
\end{multline}
\end{theorem}
\begin{theorem}
\begin{multline}
\sum _{p=1}^n \left(i 3^{-p}\right)^k e^{2 i 3^{-1-p} m} \left(\Phi \left(e^{2 i 3^{-p}
   m},-k,\frac{1}{3}-\frac{1}{2} i 3^p \log (a)\right)\right. \\ \left.
+e^{2 i 3^{-1-p} m} \Phi \left(e^{2 i 3^{-p}
   m},-k,\frac{2}{3}-\frac{1}{2} i 3^p \log (a)\right)\right)\\
=-\left(\left(\frac{i}{3}\right)^k e^{\frac{2 i m}{3}}
   \Phi \left(e^{\frac{2 i m}{3}},-k,1-\frac{3}{2} i \log (a)\right)\right. \\ \left.
-\left(i 3^{-1-n}\right)^k e^{2 i 3^{-1-n} m} \Phi
   \left(e^{2 i 3^{-1-n} m},-k,1-\frac{1}{2} i 3^{1+n} \log (a)\right)\right)
\end{multline}
\end{theorem}
\begin{theorem}
\begin{multline}
\sum _{p=1}^n 3^p \left(i 3^p\right)^k \left((-1)^{5/6} e^{i 3^p \pi  r} \Phi \left(-\sqrt[3]{-1} e^{i 3^p \pi 
   r},-k,1+3^{-p} a\right)\right. \\ \left.
+i e^{\frac{1}{3} i \pi  \left(-1+3^{1+p} r\right)} \Phi \left((-1)^{2/3} e^{i 3^p \pi 
   r},-k,1+3^{-p} a\right)\right)\\
=3 i \left((3 i)^k e^{3 i \pi  r} \Phi \left(e^{3 i \pi 
   r},-k,1+\frac{a}{3}\right)\right. \\ \left.
   -3^n \left(i 3^{1+n}\right)^k e^{i 3^{1+n} \pi  r} \Phi \left(e^{i 3^{1+n} \pi 
   r},-k,1+3^{-1-n} a\right)\right)
\end{multline}
\end{theorem}
\begin{theorem}
From Eq. (18) on pp. 292 in \cite{hall}.
\begin{multline}
\sum _{p=1}^n \left(i 3^p\right)^k e^{i 3^{-1+p} m} \left(\Phi \left(e^{2 i 3^p m},-k,\frac{1}{6}-\frac{1}{2} i
   3^{-p} \log (a)\right)\right. \\ \left.
+e^{4 i 3^{-1+p} m} \Phi \left(e^{2 i 3^p m},-k,\frac{5}{6}-\frac{1}{2} i 3^{-p} \log
   (a)\right)\right)\\
=-\left(-i^k e^{i m} \Phi \left(e^{2 i m},-k,\frac{1}{2}-\frac{1}{2} i \log (a)\right)\right. \\ \left.
+\left(i
   3^n\right)^k e^{i 3^n m} \Phi \left(e^{2 i 3^n m},-k,\frac{1}{2}-\frac{1}{2} i 3^{-n} \log
   (a)\right)\right)
\end{multline}
\end{theorem}
\begin{theorem}
From Eq. (18) on pp. 292 in \cite{hall}.
\begin{multline}
\sum _{p=1}^n \left(i 3^{-p}\right)^{1+k} e^{3^{-p} m} \left((-1)^{2/3} \Phi \left(\sqrt[3]{-1} e^{3^{-p}
   m},-k,1+3^p a\right)\right. \\ \left.
   +\Phi \left(-(-1)^{2/3} e^{3^{-p} m},-k,1+3^p a\right)\right)\\
=(-1)^{5/6} 3^{-n} \left(-i^k 3^n
   e^m \Phi \left(-e^m,-k,1+a\right)\right. \\ \left.
   +\left(i 3^{-n}\right)^k e^{3^{-n} m} \Phi \left(-e^{3^{-n} m},-k,1+3^n
   a\right)\right)
\end{multline}
\end{theorem}
\begin{theorem}
From Eq. (18) on pp. 292 in \cite{hall}.
\begin{multline}
\sum _{p=1}^n 3^p \left(a^k+\sqrt[3]{-1} \left(i 3^p\right)^k e^{i 3^p m} \left(\Phi \left(\sqrt[3]{-1} e^{i
   3^p m},-k,1-i 3^{-p} a\right)\right.\right. \\ \left.\left.
-\sqrt[3]{-1} \Phi \left(-(-1)^{2/3} e^{i 3^p m},-k,1-i 3^{-p}
   a\right)\right)\right)\\
=\frac{3}{2} \left(\left(-1+3^n\right) a^k+2 (3 i)^k e^{3 i m} \Phi \left(-e^{3 i
   m},-k,1-\frac{i a}{3}\right)\right. \\ \left.
-2\times 3^n \left(i 3^{1+n}\right)^k e^{i 3^{1+n} m} \Phi \left(-e^{i 3^{1+n} m},-k,1-i
   3^{-1-n} a\right)\right)
\end{multline}
\end{theorem}
\begin{theorem}
From Eq. (18) on pp. 292 in \cite{hall}.
\begin{multline}
-\sum _{p=1}^n \left(i 3^{-p}\right)^k e^{i 3^{-1-p} m} \left(\Phi \left(e^{2 i 3^{-p}
   m},-k,\frac{1}{6}+\frac{3^p a}{2}\right)\right. \\ \left.
   +e^{4 i 3^{-1-p} m} \Phi \left(e^{2 i 3^{-p} m},-k,\frac{5}{6}+\frac{3^p
   a}{2}\right)\right)\\
=\left(\frac{i}{3}\right)^k e^{\frac{i m}{3}} \Phi \left(e^{\frac{2 i
   m}{3}},-k,\frac{1}{2}+\frac{3 a}{2}\right)\\
   -\left(i 3^{-1-n}\right)^k e^{i 3^{-1-n} m} \Phi \left(e^{2 i 3^{-1-n}
   m},-k,\frac{1}{2} \left(1+3^{1+n} a\right)\right)
\end{multline}
\end{theorem}
\begin{theorem}
From Eq. (17.3.14) in \cite{hansen}.
\begin{multline}
\sum _{p=0}^{\infty } (-1)^p a^{-2 i p \alpha } \cos (x+2 p x) \left(-(-i
   (m+\alpha +2 p \alpha ))^{-1-k} \Gamma (1+k,-i (m+\alpha +2 p \alpha ) \log
   (a))\right. \\ \left.
+a^{2 i (1+2 p) \alpha } (-i (m-(1+2 p) \alpha ))^{-1-k} \Gamma (1+k,-i
   (m-(1+2 p) \alpha ) \log (a))\right)\\
=-\frac{1}{2} i a^{i (m+\alpha )}
   e^{\frac{i m (\pi -2 x)}{2 \alpha }} \pi ^{1+k} \left(\frac{1}{\alpha
   }\right)^{1+k} \left(\Phi \left(-e^{\frac{i m \pi }{\alpha }},-k,\frac{\pi -2
   x+2 \alpha  \log (a)}{2 \pi }\right)\right. \\ \left.
+e^{\frac{2 i m x}{\alpha }} \Phi
   \left(-e^{\frac{i m \pi }{\alpha }},-k,\frac{\pi +2 x+2 \alpha  \log (a)}{2
   \pi }\right)\right)
\end{multline}
\end{theorem}
\begin{theorem}
From Eq. (6.1.62) in \cite{hansen}
\begin{multline}
\sum _{p=0}^{\infty } (-1)^p a^{-m-\frac{1}{2} (1+2 p) \beta } \left(a^{\beta +2 p \beta } (-2 m+\beta +2 p \beta
   )^{-1-k} \Gamma \left(1+k,\frac{1}{2} (-2 m+\beta +2 p \beta ) \log (a)\right)\right. \\ \left.
-(-2 m-(1+2 p) \beta )^{-1-k} \Gamma
   \left(1+k,-\frac{1}{2} (2 m+\beta +2 p \beta ) \log (a)\right)\right)\\
=-i e^{\frac{i m \pi }{\beta }} \left(\frac{\pi 
   i}{\beta }\right)^{1+k} \Phi \left(-e^{\frac{2 i m \pi }{\beta }},-k,\frac{\pi -i \beta  \log (a)}{2 \pi
   }\right)
\end{multline}
\end{theorem}
\begin{theorem}
From Eq. (37.2.2) in \cite{hansen}.
\begin{multline}
\sum _{p=1}^n 2^k \left(i 2^{-p}\right)^{1+k} e^{i 2^{-1-p} m} \left(\Phi \left(e^{i 2^{1-p}
   m},-k,\frac{1}{4}-i 2^{-1+p} \log (a)\right)\right. \\ \left.
+e^{i 2^{-p} m} \Phi \left(e^{i 2^{1-p} m},-k,\frac{3}{4}-i 2^{-1+p}
   \log (a)\right)\right. \\ \left.
   -2 e^{i 2^{-1-p} m} \Phi \left(e^{i 2^{1-p} m},-k,\frac{1}{2} \left(1-i 2^p \log
   (a)\right)\right)\right)\\
=-i 2^{-n} \left(2^n e^{\frac{1}{2} i (m+k \pi )} \Phi \left(e^{i m},-k,\frac{1}{2}-i
   \log (a)\right)\right. \\ \left.
-\left(i 2^{-n}\right)^k e^{i 2^{-1-n} m} \Phi \left(e^{i 2^{-n} m},-k,\frac{1}{2}-i 2^n \log
   (a)\right)\right)
\end{multline}
\end{theorem}
\begin{theorem}
From Eq. (37.2.2) in \cite{hansen}.
\begin{multline}
\sum _{p=1}^n \left(i 2^p\right)^k e^{i 2^p m} \left(2 i \Phi \left(-e^{i 2^p m},-k,1+2^{-p} a\right)+\Phi
   \left(-i e^{i 2^p m},-k,1+2^{-p} a\right)\right. \\ \left.
-\Phi \left(i e^{i 2^p m},-k,1+2^{-p} a\right)\right)\\
=2 i \left((2 i)^k
   e^{2 i m} \Phi \left(-e^{2 i m},-k,1+\frac{a}{2}\right)\right. \\ \left.
   -\left(i 2^{1+n}\right)^k e^{i 2^{1+n} m} \Phi
   \left(-e^{i 2^{1+n} m},-k,1+2^{-1-n} a\right)\right)
\end{multline}
\end{theorem}
\begin{theorem}
From Eq. (7.10) pp. 119 in \cite{hobson}.
\begin{multline}\label{eq_3116}
\sum _{p=1}^{2 n} (-1)^p e^{\frac{i p \pi }{n}} \Phi \left(-e^{i \left(m+\frac{p \pi
   }{n}\right)},-k,a\right)\\
   =(-1)^{n+1}  e^{i m (-1+n)} (2n)^{1+k} \Phi \left(e^{2 i m n},-k,\frac{-1+a+n}{2
   n}\right)
\end{multline}
\end{theorem}
\begin{theorem}
From Eq. (\ref{eq_2466}) and Eq. (41.12.2) in \cite{hansen}.
\begin{multline}
\sum _{p=0}^n \left(2^p\right)^{1+k} \left(-\zeta \left(-k,3\times 2^{-p} a\right)+3^k \left(\zeta \left(-k,2^{-p}
   a\right)+\zeta \left(-k,\frac{1}{6}+2^{-p} a\right)+\zeta \left(-k,\frac{5}{6}+2^{-p} a\right)\right)\right)\\
=2^{-k}
   \left(-3^k \zeta (-k,2 a)+\zeta (-k,6 a)+2^{1+k+n} \left(2^n\right)^k \left(3^k \zeta \left(-k,2^{-n} a\right)-\zeta
   \left(-k,3\times 2^{-n} a\right)\right)\right)
\end{multline}
\end{theorem}
\begin{theorem}
Eq. (5.11.2.13) in \cite{prud2}.
\begin{equation}\label{eq_511213}
\sum _{k=0}^{\infty } \frac{(1+a k)^k \left(e^{-a z} z\right)^k
   L_k^{\alpha }\left(\frac{x}{1+a k}\right)}{(1+\alpha )_k}=\frac{e^z (x
   z)^{-\frac{\alpha }{2}} J_{\alpha }\left(2 \sqrt{x z}\right) \Gamma (1+\alpha
   )}{1-a z}
\end{equation}
\end{theorem}
\begin{theorem}
From Eq. (\ref{eq_511213}) and section (7) in \cite{plos}.  Finite Hankel transform distribution similar forms in \cite{pathak}
\begin{multline}
\sum _{k=0}^{\infty } \frac{(1+a k)^k \left(e^{-a z} z\right)^k L_k^{\alpha }\left(\frac{x}{1+a
   k}\right)}{(k+\alpha +1)!}\\
=e^{a z (1+\alpha )} z^{-1-\alpha } \int_0^z \frac{e^{-t (-1+a+a \alpha )}
   t^{1+\alpha } (t x)^{-\frac{\alpha }{2}} J_{\alpha }\left(2 \sqrt{t x}\right)}{t} \, dt
\end{multline}
\end{theorem}
\begin{theorem}
Eq. (5.11.2.14) in \cite{prud2}
\begin{multline}\label{eq_511214}
\sum _{k=0}^{\infty } \frac{(1+a k)^{-1+k} \left(e^{-a z} z\right)^k L_k^{\alpha }\left(\frac{x}{1+a
   k}\right)}{(k+\alpha )!}=\frac{e^z \, _1F_2\left(\frac{1}{a};1+\frac{1}{a},1+\alpha ;-x z\right)}{\Gamma (1+\alpha
   )}
\end{multline}
\end{theorem}
\begin{theorem}
From Eq. (\ref{eq_511214}) and section (7) in \cite{plos}
\begin{multline}
\sum _{k=0}^{\infty } \frac{(1+a k)^{-1+k} \left(e^{-a z} z\right)^k L_k^{\alpha }\left(\frac{x}{1+a
   k}\right)}{(k+\alpha +1)!}\\
=-\left(e^{a z (1+\alpha )} z^{-1-\alpha }\right) \int_0^z \frac{e^{-t (-1+a+a \alpha )}
   t^{\alpha } (-1+a t) \, _1F_2\left(\frac{1}{a};1+\frac{1}{a},1+\alpha ;-t x\right)}{\Gamma (1+\alpha )} \,
   dt
\end{multline}
\end{theorem}
\begin{theorem}
From Eq. (5.11.4.6) in \cite{prud2} and section (7) in \cite{plos}
\begin{multline}
\sum _{k=0}^{\infty } \frac{L_k^{k \alpha +\beta }(x) \left(\frac{v}{(1+v)^{\alpha
   +1}}\right)^k}{k+1}\\
=\frac{e^x (1+v)}{v} \left((1+v)^{\alpha } x^{\alpha -\beta } \Gamma (-\alpha +\beta
   ,x)\right. \\ \left.
-(1+v)^{\beta } ((1+v) x)^{\alpha -\beta } \Gamma (-\alpha +\beta ,(1+v) x)\right)
\end{multline}
\end{theorem}
\begin{theorem}
From Eq. (48.17.2) in \cite{hansen} and section (7) in \cite{plos}. Extended Toscano Laguerre-form.
\begin{multline}
\sum _{k=0}^{\infty } \frac{4^{-k} \left(1-v^2\right)^k L_k^{c+k}(x)}{k+1}=\frac{2 e^x \left(-2 E_{2-c}(x)+2^c
   \left(\frac{1}{1+v}\right)^{-1+c} E_{2-c}\left(\frac{2 x}{1+v}\right)\right)}{-1+v^2}
\end{multline}
\end{theorem}
\begin{theorem}
From Eq. (66.4.1) in \cite{hansen} and section (7) in \cite{plos}
\begin{multline}
\sum _{k=0}^{\infty } \frac{\left(v (1+v)^{-1-c}\right)^k \, _1F_1(-k;a+c k;x) (a)_{k+c k}}{(k+1)! (a)_{c
   k}}\\
=\frac{e^x \left((1+v)^{1+c} E_{2-a+c}(x)-(1+v)^a E_{2-a+c}((1+v) x)\right)}{v}
\end{multline}
\end{theorem}
\begin{theorem}
From Eq. (48.19.14) in \cite{hansen} and section (7) in \cite{plos}
\begin{multline}
\sum _{k=0}^{\infty } \frac{4^{-k} \left(-1+v^2\right)^k L_{k+m}^{-1-2 k-2 m}(x) (1+m)_k}{(k+1)!}=\int_1^v
   \frac{2 e^{\frac{1}{2} (1-t) x} t^{-2 m} L_m^{-1-2 m}(t x)}{-1+v^2} \, dt
\end{multline}
\end{theorem}
\subsubsection{The Gegenbauer polynomials}
\begin{theorem}
From Eq. (5.13.12) in \cite{prud2} and section (7) in \cite{plos}
\begin{multline}
\sum _{k=0}^{\infty } \frac{\left(b^2+(1+a k)^2\right)^{k/2} \left(e^{-a z} z\right)^k
   C_k^{(v)}\left(\frac{1+a k}{\sqrt{b^2+(1+a k)^2}}\right)}{(k+2 v)!}\\
=e^{2 a v z} z^{-2 v} \int_0^z \frac{e^{t-2 a t
   v} t^{2 v} \, _0F_1\left(;\frac{1}{2}+v;-\frac{1}{4} b^2 t^2\right)}{t \Gamma (2 v)} \, dt
\end{multline}
\end{theorem}
\begin{theorem}
From Eq. (5.13.13) in \cite{prud2} and section (7) in \cite{plos}
\begin{multline}
\sum _{k=0}^{\infty } \frac{(1+a k)^{-1+\frac{k}{2}} (1+b+a k)^{k/2} \left(e^{-a z} z\right)^k
   C_k^{(v)}\left(\frac{2+b+2 a k}{2 \sqrt{(1+a k) (1+b+a k)}}\right)}{(k+2 v)!}\\
=-e^{2 a v z} z^{-2 v} \int_0^z
   \frac{e^{t-2 a t v} t^{-1+2 v} (-1+a t) \, _2F_2\left(\frac{1}{a},v;1+\frac{1}{a},2 v;b t\right)}{\Gamma (2 v)} \,
   dt
\end{multline}
where $Im(z)\neq 0$
\end{theorem}
\begin{theorem}
From Eq. (5.13.14) in \cite{prud2} and section (7) in \cite{plos}
\begin{multline}
\sum _{k=0}^{\infty } \frac{\left(-4 b+(1+a k)^2\right)^{k/2} \left(e^{-a z} z\right)^k C_k^{(v)}\left(\frac{1+a
   k}{\sqrt{-4 b+(1+a k)^2}}\right)}{(1+a k) (k+2 v)!}\\
=-e^{2 a v z} z^{-2 v} \int_0^z \frac{e^{t-2 a t v} t^{-1+2 v}
   (-1+a t) \, _1F_2\left(\frac{1}{2 a};1+\frac{1}{2 a},\frac{1}{2}+v;b t^2\right)}{\Gamma (2 v)} \, dt
\end{multline}
where $Im(a)\neq 0,Im(z)\neq 0$
\end{theorem}
\begin{theorem}
From Eq. (5.24.14) in \cite{hansen} and section (7) in \cite{plos}
\begin{equation}\label{eq_52414}
\sum _{k=0}^{\infty } \frac{(k+1)^{k-1} (y \exp (-y))^k}{(k+1)!}=-\frac{1}{2} e^y (-2+y)
\end{equation}
\end{theorem}
\begin{theorem}
From Eq. (5.13.22) in \cite{hansen} and section (7) in \cite{plos}
\begin{equation}\label{eq_51322}
\sum _{k=0}^{\infty } \frac{(k+2)^k \left(e^{-y} y\right)^{k+1}}{(k+2)!}=\frac{1}{2} \left(-e^y y+2
   e^y-2\right)
\end{equation}
\end{theorem}
\begin{theorem}
From Eq. (5.13.20) in \cite{hansen} and Eq. (6.3.4.2(a)) in \cite{prud3}.
\begin{multline}
\sum _{n=0}^{\infty } \sum _{p=1}^{\infty } \frac{(-1)^{2 n+p} (a-p)^{k-2 n-p} p^{-1+p} E_{2 n}(x) (1-k)_{-1+2
   n+p}}{(2 n)! \Gamma (2+p)}\\
=-\frac{2^{-1+k} }{k (1+k)}\left(4 \zeta \left(-1-k,\frac{1}{2} (1+a-x)\right)-8 \zeta
   \left(-1-k,\frac{1}{2} (2+a-x)\right)\right. \\ \left.
+4 \zeta \left(-1-k,\frac{1}{2} (3+a-x)\right)+4 \zeta
   \left(-1-k,\frac{a+x}{2}\right)-8 \zeta \left(-1-k,\frac{1}{2} (1+a+x)\right)\right. \\ \left.
+4 \zeta \left(-1-k,\frac{1}{2}(2+a+x)\right)+(1+k) \left(k \left(\zeta \left(1-k,\frac{1}{2} (1+a-x)\right)
\right.\right.\right. \\ \left.\left.\left.
-\zeta \left(1-k,\frac{1}{2} (2+a-x)\right)+\zeta \left(1-k,\frac{a+x}{2}\right)-\zeta \left(1-k,\frac{1}{2} (1+a+x)\right)\right)\right.\right. \\ \left.\left.
+2 \left(\zeta\left(-k,\frac{1}{2} (1+a-x)\right)-\zeta \left(-k,\frac{1}{2} (2+a-x)\right)+\zeta
   \left(-k,\frac{a+x}{2}\right)\right.\right.\right. \\ \left.\left.\left.
-\zeta \left(-k,\frac{1}{2} (1+a+x)\right)\right)\right)\right)
\end{multline}
\end{theorem}
\begin{theorem}
From Eq. (5.13.20) in \cite{hansen} and section (7) in \cite{plos}
\begin{equation}\label{eq_51320}
\sum _{k=1}^{\infty } \frac{k^{k-1} (y \exp (-y))^k}{(k+1)!}=\frac{1-e^y+y+y^2}{y}
\end{equation}
\end{theorem}
\begin{theorem}
From Eq. (45.4.2) in \cite{hansen} and section (7) in \cite{plos}
\begin{multline}
\sum _{k=0}^{\infty } \frac{\left(v (1+v)^{-1-c}\right)^k P_k^{(a+c k,b-k-c k)}(x)}{k+n}\\
=\frac{(1+v)^{(1+c) n}}{2 n (1+n)}
   \left(2 (1+n) F_1\left(n;-a+n+c n,a+b;1+n;-v,\frac{1}{2} v (-1+x)\right)\right. \\ \left.
+n v (-1+x) F_1\left(1+n;-a+n+c
   n,1+a+b;2+n;-v,\frac{1}{2} v (-1+x)\right)\right)
\end{multline}
\end{theorem}
\begin{theorem}
From Eq. (45.4.3) in \cite{hansen} and section (7) in \cite{plos}
\begin{multline}
\sum _{k=0}^{\infty } \frac{\left((1-z)^{-1-c} z\right)^k P_k^{(a-k-c k,b+c k)}(x)}{k+n}\\
=\frac{(1-z)^{(1+c) n}}{2 n (1+n)}
   \left(2 (1+n) F_1\left(n;-b+n+c n,a+b;1+n;z,\frac{1}{2} (1+x) z\right)\right. \\ \left.
+n (1+x) z F_1\left(1+n;-b+n+c
   n,1+a+b;2+n;z,\frac{1}{2} (1+x) z\right)\right)
\end{multline}
\end{theorem}
\begin{theorem}
From Eq. (45.4.4) in \cite{hansen} and section (7) in \cite{plos}
\begin{multline}
\sum _{k=0}^{\infty } \frac{2^k \left(-\frac{w (1+w)^{-c}}{-1+x}\right)^k P_k^{(a-k,b-c
   k)}(x)}{k+n}\\
=\frac{(1+w)^{c n}}{n (1+n)} \left((1+n) F_1\left(n;a+b+c n,-a;1+n;-w,\frac{2 w}{-1+x}\right)\right. \\ \left.-n w
   F_1\left(1+n;1+a+b+c n,-a;2+n;-w,\frac{2 w}{-1+x}\right)\right)
\end{multline}
\end{theorem}
\begin{theorem}
From Eq. (45.4.5) in \cite{hansen} and section (7) in \cite{plos}
\begin{multline}
\sum _{k=0}^{\infty } \frac{2^k \left(\frac{(1-u)^{-c} u}{1+x}\right)^k P_k^{(a-c k,b-k)}(x)}{k+n}\\
=\frac{(1-u)^{c
   n}}{n (1+n)} \left((1+n) F_1\left(n;a+b+c n,-b;1+n;u,\frac{2 u}{1+x}\right)\right. \\ \left.
+n u F_1\left(1+n;1+a+b+c n,-b;2+n;u,\frac{2
   u}{1+x}\right)\right)
\end{multline}
\end{theorem}
\begin{theorem}
From Eq. (57.14.8) in \cite{hansen} and section (7) in \cite{plos}
\begin{multline}
\sum _{k=0}^{\infty } \frac{4^{-2 k} J_{a+2 k}(z) (z \sin (2 t))^{2 k}}{(k+1)! (1+a)_k}\\
=\csc ^2(2 t) \int_0^t
   2^{2+2 a} J_a(z \cos (s)) J_a(z \sin (s)) \cos (2 s) \Gamma (1+a) \sin (2 s) (z \sin (2 s))^{-a} \, ds
\end{multline}
\end{theorem}
\begin{theorem}
From Eq. (65.4.4) in \cite{hansen} and section (7) in \cite{plos}
\begin{multline}
\sum _{k=0}^{\infty } \frac{\left(v (1+v)^{-1-c}\right)^k \, _2F_1(a,-k;b+c k;x) (b)_{k+c k}}{(k+1)! (b)_{c
   k}}\\
=\frac{(1+v) x}{(-1+a) v (-1+x)^2} \left((1+v)^c \left(\frac{x}{-1+x}\right)^{-b+c} \,
   _2F_1\left(1-a,2-b+c;2-a;\frac{1}{1-x}\right)\right. \\ \left.
-(1+v)^b \left(\frac{(1+v) x}{-1+x}\right)^{-b+c} (1+v x)^{1-a} \,
   _2F_1\left(1-a,2-b+c;2-a;\frac{1+v x}{1-x}\right)\right)
\end{multline}
\end{theorem}
\begin{theorem}
From Eq. (65.4.2) in \cite{hansen} and section (7) in \cite{plos}
\begin{multline}
\sum _{k=0}^{\infty } \frac{4^{-k} \left(x^2 \left(1+\frac{1}{-1+x}+x\right)^k+(u-x) (u+x)
   \left(\frac{-u^2+x^2}{-1+x}\right)^k\right) }{(u-x) (u+x) \Gamma (2+k)
   (c)_{2 k}}\\\times
\, _2F_1(-a+c,-b+c+k;c+2 k;x) (a)_{2 k} (b)_k\\
=\int_0^u \frac{2 t (1-x)^{a+b-c} F_1\left(a;b,b;c;\frac{t+x}{2},\frac{1}{2} (-t+x)\right)}{u^2-x^2} \,
   dt
\end{multline}
\end{theorem}
%
%-------special cases
%
\section{Special cases in terms of special functions and fundamental constants}
The constants used in this section are; Stieltjes constant $\gamma_{n}$ given in Eq. (25.2.5) in \cite{dlmf}, Glaisher's constant $A$ given in Eq. (5.17.7) in \cite{dlmf}, Harmonic number function $H_{n}$ given in Eq.(25.11.33) in \cite{dlmf}, Catalan's constant $C$ given in Eq.(25.11.40) in \cite{dlmf}, Euler's constant $\gamma$ given in Eq. (5.4.11) in \cite{dlmf}, Polygamma Functions $\psi^{(n)}(z)$ given in section (5.15) in \cite{dlmf}, $q$-digamma function $\psi_{q}(z)$ given on pp. 91-101 in \cite{borwein}, the the $q$-Pochhammer symbol $(a;q)_n$ given on pp. 261-286 in \cite{berndt4}.
\subsection{Definite integrals}
\begin{example}
Malmsten integral see Exercise 35 on pp. 80 in \cite{blagouchine}.
\begin{multline}
\int_{1}^{\infty}\frac{x^2 \log (\log (x))
   \left(x^{m/n}+\left(\frac{1}{x}\right)^{m/n}\right)}{\left(x^2+1\right)^3}dx\\
   =\frac{i \pi }{32 \pi  n^2} \left(8 m n
   \text{gd}\left(\frac{i \pi  m}{2 n}\right)+\pi  (m-n) (m+n)\right. \\ \left.
    \left(H_{\frac{1}{4}
   \left(\frac{m}{n}-1\right)}-H_{\frac{m+n}{4 n}}-(\pi -2 i \log (\pi )) \sec \left(\frac{\pi  m}{2
   n}\right)\right)\right. \\ \left.
   +2 \pi  n (3 m-2 n)\right)-4 e^{\frac{i \pi  m}{2 n}} \left(\pi ^2 (n-m) (m+n)
   \Phi'\left(-e^{\frac{i \pi  m}{n}},0,\frac{1}{2}\right)\right. \\ \left.
   +n^2 \Phi \left(-e^{\frac{i
   m \pi }{n}},2,\frac{1}{2}\right)\right)
\end{multline}
\end{example}
\begin{example}
\begin{multline}
\int_0^{\infty } \frac{x^{-m-q} \left(x^m-x^q\right) \Gamma (b+x)}{\Gamma (1+b+n+x) \log
   \left(\frac{1}{x}\right)} \, dx\\
=\sum _{p=0}^{\infty } \frac{\left(e^{i m \pi } (b+p)^{-m} \Phi \left(e^{2 i m \pi
   },1,\frac{\pi +i \log (b+p)}{2 \pi }\right)-e^{i \pi  q} (b+p)^{-q} \Phi \left(e^{2 i \pi  q},1,\frac{\pi +i \log
   (b+p)}{2 \pi }\right)\right) (-n)_p}{p! \Gamma (1+n)}
\end{multline}
\end{example}
\begin{example}
\begin{multline}
\int_0^{\infty } \frac{\left(-1+\sqrt[4]{x}\right) \Gamma (1+x)}{\sqrt{x} \Gamma (5+x) \log
   \left(\frac{1}{x}\right)} \, dx\\
=\sum _{p=0}^{\infty } \frac{i \left(\Phi \left(-1,1,\frac{\pi +i \log (1+p)}{2 \pi
   }\right)+(-1)^{3/4} \sqrt[4]{1+p} \Phi \left(i,1,\frac{\pi +i \log (1+p)}{2 \pi }\right)\right) (-3)_p}{6
   \sqrt{1+p} p!}
\end{multline}
\end{example}
\begin{example}
\begin{multline}
\int_0^{\infty } \frac{\Gamma \left(\frac{1}{3}+x\right) \log ^k\left(\frac{1}{x}\right)}{\sqrt{x} \Gamma
   \left(\frac{10}{3}+x\right)} \, dx\\
=\sum _{p=0}^{\infty } \frac{(4 i)^k \pi ^{1+k} \left(\zeta \left(-k,\frac{\pi
   +i \log \left(\frac{1}{3}+p\right)}{4 \pi }\right)-\zeta \left(-k,\frac{3}{4}+\frac{i \log
   \left(\frac{1}{3}+p\right)}{4 \pi }\right)\right) (-2)_p}{\sqrt{\frac{1}{3}+p} p!}
\end{multline}
\end{example}
\begin{example}
\begin{multline}
\int_0^{\infty } \frac{\Gamma (b+x) \log ^k\left(\frac{a}{x}\right)}{\sqrt{x} \Gamma (1+b+n+x)} \, dx\\
=\sum
   _{p=0}^{\infty } \frac{2^{1+2 k} e^{\frac{i k \pi }{2}} \pi ^k \Gamma (-n+p) \left(-\zeta \left(-k,\frac{\pi -i
   \log (a)+i \log (b+p)}{4 \pi }\right)+\zeta \left(-k,\frac{3 \pi -i \log (a)+i \log (b+p)}{4 \pi }\right)\right)
   \sin (n \pi )}{\sqrt{b+p} \Gamma (1+p)}
\end{multline}
\end{example}
\begin{example}
\begin{multline}
\int_0^{\infty } \frac{x^{-m} \Gamma (b+x)}{\Gamma (1+b+n+x) \left(a^2+c^2+2 i c \log (x)-\log ^2(x)\right)}
   \, dx\\
=\sum _{p=0}^{\infty } \frac{e^{i m \pi } (b+p)^{-m} \left(\Phi \left(e^{2 i m \pi },1,\frac{\pi -i (-a+i
   c)+i \log (b+p)}{2 \pi }\right)-\Phi \left(e^{2 i m \pi },1,\frac{\pi -i (a+i c)+i \log (b+p)}{2 \pi
   }\right)\right) (-n)_p}{2 a p! \Gamma (1+n)}
\end{multline}
\end{example}
\begin{example}
\begin{multline}
\int_0^{\infty } \frac{\log \left(\log \left(\frac{1}{x}\right)\right)}{\sqrt{x} (1+x) (2+x) (3+x) (4+x)} \,
   dx\\
=\frac{\pi }{24 \sqrt{2}} \left(i \left(-6+\sqrt{2}+2 \sqrt{6}\right) \pi +\sqrt{2} \log (16)+\sqrt{6} \log (256)-2 \log
   (4096)\right. \\ \left.
+2 \left(-6+\sqrt{2}+2 \sqrt{6}\right) \log (\pi )+8 \sqrt{2} \log \left(\frac{\Gamma
   \left(\frac{3}{4}\right)}{\Gamma \left(\frac{1}{4}\right)}\right)\right. \\ \left.
+\log \left(\frac{\Gamma \left(\frac{\pi +i \log
   (2)}{4 \pi }\right)^{24} \left(\frac{\Gamma \left(\frac{\pi +2 i \log (2)}{4 \pi }\right)}{\Gamma
   \left(\frac{3}{4}+\frac{i \log (2)}{2 \pi }\right)}\right)^{4 \sqrt{2}} \left(\frac{\Gamma \left(\frac{\pi +i \log
   (3)}{4 \pi }\right)}{\Gamma \left(\frac{3}{4}+\frac{i \log (3)}{4 \pi }\right)}\right)^{-8 \sqrt{6}}}{\Gamma
   \left(\frac{3}{4}+\frac{i \log (2)}{4 \pi }\right)^{24}}\right)\right)
\end{multline}
\end{example}
\begin{example}
\begin{multline}
\int_0^{\infty } \frac{\Gamma (b+x) \log ^k\left(\frac{1}{x}\right)}{\sqrt{x} \Gamma (1+b+n+x)} \, dx\\
=-\sum
   _{p=0}^{\infty } \frac{i^{2+k} (2 \pi )^{1+k} \left(-2^k \zeta \left(-k,\frac{1}{2} \left(1+\frac{\pi +i \log
   (b+p)}{2 \pi }\right)\right)+2^k \zeta \left(-k,\frac{\pi +i \log (b+p)}{4 \pi }\right)\right) (-n)_p}{\sqrt{b+p}
   p! \Gamma (1+n)}
\end{multline}
\end{example}
\begin{example}
\begin{multline}
\int_0^{\infty } \frac{\log \left(\log \left(\frac{1}{x}\right)\right)}{\sqrt{x} (1+2 x) (3+2 x)} \,
   dx\\
=-\frac{i \left(-3+\sqrt{3}\right) \pi ^2}{12 \sqrt{2}}\\
-\frac{\pi  \left(\sqrt{3} \log (16)-6 \log (\pi )+2
   \sqrt{3} \log (\pi )-\log \left(\frac{4096 \left(\frac{\Gamma \left(\frac{\pi +i \log \left(\frac{3}{2}\right)}{4
   \pi }\right)}{\Gamma \left(\frac{3}{4}+\frac{i \log \left(\frac{3}{2}\right)}{4 \pi }\right)}\right)^{4 \sqrt{3}}
   \Gamma \left(\frac{3}{4}-\frac{i \log (2)}{4 \pi }\right)^{12}}{\Gamma \left(\frac{\pi -i \log (2)}{4 \pi
   }\right)^{12}}\right)\right)}{12 \sqrt{2}}
\end{multline}
\end{example}
\begin{example}
\begin{multline}
\int_0^{\infty } \frac{\Gamma (i+x)}{\sqrt{x} \Gamma ((3+i)+x) \log ^2\left(\frac{i}{x}\right)} \, dx\\
=\frac{80
   (-1)^{3/4} C+\frac{10 \psi ^{(1)}\left(\frac{5}{16}+\frac{i \log (2)}{8 \pi }\right)}{\sqrt{1+i}}-\frac{10 \zeta
   \left(2,\frac{7}{8}+\frac{i \log (1+i)}{4 \pi }\right)}{\sqrt{1+i}}-\frac{5 \left(\zeta
   \left(2,\frac{3}{8}+\frac{i \log (2+i)}{4 \pi }\right)-\zeta \left(2,\frac{7}{8}+\frac{i \log (2+i)}{4 \pi
   }\right)\right)}{\sqrt{2+i}}}{80 \pi }
\end{multline}
\end{example}
\begin{example}
\begin{multline}
\int_0^{\infty } \frac{\Gamma (1+x) \log \left(\log \left(\frac{1}{x}\right)\right)}{\sqrt{x} \Gamma (4+x)} \,
   dx\\
=\frac{1}{12} \pi  \left(-i \left(5+3 \sqrt{2}-5 \sqrt{3}\right) \pi +\sqrt{3} \log (16)-\sqrt{2} \log (4096)+2
   \left(3-3 \sqrt{2}+\sqrt{3}\right) \log (\pi )\right. \\ \left.
-12 \log \left(\Gamma \left(\frac{1}{4}\right)\right)+12 \log
   \left(\Gamma \left(\frac{3}{4}\right)\right)\right. \\ \left.
+\log \left(4096 \left(\frac{\Gamma \left(\frac{\pi +i \log (2)}{4 \pi
   }\right)}{\Gamma \left(\frac{3}{4}+\frac{i \log (2)}{4 \pi }\right)}\right)^{12 \sqrt{2}} \left(-\frac{\Gamma
   \left(\frac{\pi +i \log (3)}{4 \pi }\right)}{\Gamma \left(\frac{3}{4}+\frac{i \log (3)}{4 \pi }\right)}\right)^{-4
   \sqrt{3}}\right)\right)
\end{multline}
\end{example}
\begin{example}
\begin{multline}
\int_0^{\infty } \frac{\Gamma (1+x) \log \left(\frac{1}{x}\right) \log \left(\log
   \left(\frac{1}{x}\right)\right)}{\sqrt{x} \Gamma (4+x)} \, dx\\
=\frac{1}{12} \pi  \left(-24 i C+2 \left(-\sqrt{3}
   \log (3)+\sqrt{2} \log (8)\right) \log (4 \pi )\right. \\ \left.
-i \pi  \left(\sqrt{3} \log (3)-\sqrt{2} \log (8)-16 \left(3
   \sqrt{2} \zeta ^{(1,0)}\left(-1,\frac{\pi +i \log (2)}{4 \pi }\right)\right.\right.\right. \\ \left.\left.\left.
-3 \sqrt{2} \zeta
   ^{(1,0)}\left(-1,\frac{3}{4}+\frac{i \log (2)}{4 \pi }\right)+\sqrt{3} \left(-\zeta ^{(1,0)}\left(-1,\frac{\pi +i
   \log (3)}{4 \pi }\right)\right.\right.\right.\right. \\ \left.\left.\left.\left.
+\zeta ^{(1,0)}\left(-1,\frac{3}{4}+\frac{i \log (3)}{4 \pi
   }\right)\right)\right)\right)\right)
\end{multline}
\end{example}
\begin{example}
\begin{multline}
\int_0^{\infty } \frac{\Gamma (b+x) \log \left(\log \left(\frac{a}{x}\right)\right)}{\sqrt[4]{x} \Gamma
   (1+b+n+x)} \, dx\\
=\sum _{p=0}^{\infty } \frac{\sqrt{2} \pi  (-n)_p }{\sqrt[4]{b+p} p! \Gamma (1+n)}\left(\frac{i \pi }{2}+\log (2 \pi )-(1-i)
   \left(i \log \left(\frac{\Gamma \left(\frac{1}{4} \left(1+\frac{\pi -i \log (a)+i \log (b+p)}{2 \pi
   }\right)\right)}{2 \Gamma \left(\frac{1}{4} \left(3+\frac{\pi -i \log (a)+i \log (b+p)}{2 \pi
   }\right)\right)}\right)\right.\right. \\ \left.\left.
+\log \left(\frac{\Gamma \left(\frac{\pi -i \log (a)+i \log (b+p)}{8 \pi }\right)}{2 \Gamma
   \left(\frac{1}{4} \left(2+\frac{\pi -i \log (a)+i \log (b+p)}{2 \pi
   }\right)\right)}\right)\right)\right)
\end{multline}
\end{example}
\begin{example}
\begin{multline}
\int_0^{\infty } \frac{\left(-1+\sqrt[6]{x}\right) \Gamma (2+x)}{x^{2/3} \Gamma (6+x) \log
   \left(\frac{1}{x}\right)} \, dx\\
=\sum _{p=0}^{\infty } \frac{\left(-i \sqrt[6]{2+p} \Phi \left(-1,1,\frac{\pi +i
   \log (2+p)}{2 \pi }\right)+(-1)^{2/3} \Phi \left(-\sqrt[3]{-1},1,\frac{\pi +i \log (2+p)}{2 \pi }\right)\right)
   (-3)_p}{6 (2+p)^{2/3} p!}
\end{multline}
\end{example}
\begin{example}
\begin{multline}
\int_0^{\infty } \frac{\Gamma (3+x) \log \left(\log \left(\frac{\sqrt[3]{-1}}{x}\right)\right)}{\sqrt{x}
   \Gamma ((9+i)+x)} \, dx\\
=\sum _{p=0}^{\infty } \frac{i \Gamma ((-5-i)+p) }{2 \sqrt{3+p} \Gamma (1+p)}\left(i \pi +\log (16)+2 \log (\pi )-4
   \log (-8 \pi +3 i \log (3+p))\right. \\ \left.
+4 \log (-2 \pi +3 i \log (3+p))-4 \text{log$\Gamma $}\left(-\frac{2}{3}+\frac{i \log
   (3+p)}{4 \pi }\right)\right. \\ \left.
+4 \text{log$\Gamma $}\left(-\frac{1}{6}+\frac{i \log (3+p)}{4 \pi }\right)\right) \sinh (\pi)
\end{multline}
\end{example}
%%
%\begin{example}
%
%\end{example}
%%
%\begin{example}
%
%\end{example}
%
\begin{example}
\begin{multline}
\int_{0}^{\pi/4}\log (\tan (x))   \left(\frac{\sqrt{\tan (x)}}{\left(\log ^2(\tan (x))+\pi ^2\right)^2}+\frac{\sqrt{\cot (x)}}{\left(\log
   ^2(\cot (x))+\pi ^2\right)^2}\right)\frac{\sin (2 x)}{\sin (4 x)}dx\\
=\frac{1-2 C}{8 \pi ^2}
\end{multline}
\end{example}
\begin{example}
\begin{equation}
\int_{0}^{\pi/4}\frac{\sec (2 x) \log (\tan (x)) \log \left(\frac{\log (i \cot (x))}{\log (i \tan (x))}\right)}{\log (\cot
   (x))}dx=i \pi  \log \left(-\frac{2 \Gamma \left(\frac{1}{4}\right)}{\Gamma \left(-\frac{1}{4}\right)}\right)
\end{equation}
\end{example}
\begin{example}
\begin{equation}
\int_{-\alpha }^{\alpha } \frac{\log \left(\frac{\alpha +x}{\alpha
   -x}\right) \log \left(\log \left(\frac{\alpha +x}{\alpha
   -x}\right)\right)}{\alpha ^2+x^2} \, dx=-\frac{i \pi  C}{\alpha }
\end{equation}
\end{example}
\begin{example}
\begin{multline}
\int_{0}^{\infty}\frac{\log \left(\frac{(x+1) \left(b^2+x\right)}{(b+x)^2}\right) \left(\frac{1}{c+\log (x)}-\frac{1}{a+\log (x)}\right)}{x}dx\\
=-2 i \pi  \log
   \left(\frac{\Gamma \left(\frac{\pi -i a}{2 \pi }\right) \Gamma \left(-\frac{i (a+2 \log (b)+i \pi )}{2 \pi }\right) \Gamma \left(-\frac{i (c+\log
   (b)+i \pi )}{2 \pi }\right)^2}{\Gamma \left(\frac{\pi -i c}{2 \pi }\right) \Gamma \left(-\frac{i (a+\log (b)+i \pi )}{2 \pi }\right)^2 \Gamma
   \left(-\frac{i (c+2 \log (b)+i \pi )}{2 \pi }\right)}\right)
\end{multline}
\end{example}
\begin{example}
\begin{multline}
\int_{0}^{\infty}\frac{\log (x) \log (\log (x))}{(x-1) (x+i)}dx=\left(\frac{1}{32}+\frac{i}{32}\right) \pi  \left(-i \pi  \log \left(\frac{64 e^6 \pi
   ^6}{A^{72}}\right)-16 i C+3 \pi ^2\right)
\end{multline}
\end{example}
\begin{example}
\begin{multline}
\int_{0}^{\infty}\frac{\sqrt{\log (x)} \log (\log (x))}{(x-1) (x+i)}dx\\
=\pi ^{3/2} \left(2 i
   \left(\zeta'\left(-\frac{1}{2}\frac{3}{4}\right)-\zeta '\left(-\frac{1}{2}\right)\right)+(\pi -2 i \log (2 \pi )) \left(\zeta
   \left(-\frac{1}{2},\frac{3}{4}\right)-\zeta \left(-\frac{1}{2}\right)\right)\right)
\end{multline}
\end{example}
\begin{example}
Mellin transform involving the logarithm of quotient cosine and hyperbolic cosine functions.
\begin{multline}
\int_{0}^{\infty}x^{s-1} \log \left(\frac{\cos (\alpha )+\cosh (x)}{\cos (\beta )+\cosh (x)}\right)dx\\
=\frac{2^s \pi ^{s+1} \csc \left(\frac{\pi  s}{2}\right)
   \left(\zeta \left(-s,\frac{\pi -\alpha }{2 \pi }\right)+\zeta \left(-s,\frac{\alpha +\pi }{2 \pi }\right)-\zeta \left(-s,\frac{\pi -\beta }{2 \pi
   }\right)-\zeta \left(-s,\frac{\beta +\pi }{2 \pi }\right)\right)}{s}
\end{multline}
\end{example}
\begin{example}
\begin{equation}
\int_{0}^{\infty}\frac{\left(x^2-\pi ^2\right) \log \left(\frac{2}{\text{sech}(x)+2}\right)}{\left(x^2+\pi ^2\right)^2}dx=\log \left(\frac{2 e}{3
   \sqrt{3}}\right)
\end{equation}
\end{example}
\begin{example}
\begin{multline}
\int_{0}^{\infty}\log \left(\frac{\cos (\alpha )+\cosh (x)}{\cos (\beta )+\cosh (x)}\right)\frac{dx}{\pi ^2 a^2+x^2}=\frac{1}{a}\log \left(\frac{\Gamma \left(\frac{\pi  a-\beta +\pi }{2 \pi }\right) \Gamma
   \left(\frac{\pi  a+\beta +\pi }{2 \pi }\right)}{\Gamma \left(\frac{\pi  a-\alpha +\pi }{2 \pi }\right) \Gamma \left(\frac{\pi  a+\alpha +\pi }{2 \pi }\right)}\right)
\end{multline}
\end{example}
\begin{example}
\begin{multline}
\int_0^{\infty } \frac{(\pi  x-\sinh (\pi  x)) \text{csch}\left(\frac{\pi  x}{2}\right) \text{sech}\left(\frac{\pi 
   x}{2}\right) \log \left(\frac{\cos (\alpha )+\cosh (x)}{\cos (\beta )+\cosh (x)}\right)}{4 \pi 
   x^2}dx\\
=\sum_{a=1}^{\infty}\frac{(-1)^a }{a}\log \left(\frac{\Gamma \left(\frac{a-\beta +\pi }{2 \pi }\right) \Gamma \left(\frac{a+\beta
   +\pi }{2 \pi }\right)}{\Gamma \left(\frac{a-\alpha +\pi }{2 \pi }\right) \Gamma \left(\frac{a+\alpha +\pi }{2
   \pi }\right)}\right)
\end{multline}
\end{example}
\begin{example}
\begin{equation}
\int_0^{\infty } \frac{\sinh \left(\frac{1}{2} (m-n) x\right) \sinh \left(\frac{1}{2} (m+n-p) x\right)}{x
   \sinh \left(\frac{p x}{2}\right)} \, dx=\frac{1}{2} \log \left(\frac{\sin \left(\frac{n \pi }{p}\right)}{\sin
   \left(\frac{m \pi }{p}\right)}\right)
\end{equation}
where $m,n,p\in\mathbb{R_{+}}$.
\end{example}
\begin{example}
\begin{equation}
\int _{0}^{\infty }\frac{\text{csch}(x) \log \left(1-\frac{2 x}{i a+x}\right)}{2 \pi  i}dx=\log \left(\frac{(1+i) \sqrt{\pi }
   \Gamma \left(1+\frac{a}{2 \pi }\right)}{\sqrt{i a} \Gamma \left(\frac{a+\pi }{2 \pi }\right)}\right)
\end{equation}
\end{example}
\begin{example}
\begin{multline}
\int _{0}^{\infty }e^x \left((\pi -2 i x) \log \left(-x-\frac{i \pi }{2}\right)-(\pi +2 i x) \log \left(x-\frac{i \pi
   }{2}\right)\right) (\coth (x)-1)dx\\
=i \pi  \left(-4 C-\pi  \log \left(-\frac{i \pi }{2}\right)\right)
\end{multline}
\end{example}
\begin{example}
\begin{multline}
\int _{0}^{\infty }\text{csch}(x) \left(\frac{i \log \left(-\frac{1}{2} (i \pi )-x\right)}{\pi -2 i x}+\frac{\log
   \left(-\frac{1}{2} (i \pi )+x\right)}{i \pi -2 x}\right)dx\\
   =\log \left(\frac{e^{\frac{1}{4} (2 \gamma +i (-2+\pi ))
   \pi } \pi ^{1+\pi }}{2^{1+\frac{5 \pi }{2}} \Gamma \left(\frac{5}{4}\right)^{2 \pi }}\right)
\end{multline}
\end{example}
\begin{example}
\begin{multline}
\int _{0}^{\infty }\left((-p x-i \pi ) \log \left(-x-\frac{i \pi }{p}\right)+i (\pi +i p x) \log \left(x-\frac{i \pi
   }{p}\right)\right) \text{csch}\left(\frac{p x}{2}\right)dx\\
=\frac{\pi }{p} \left(-8 C-2 \pi  \log \left(\frac{\pi
   }{p}\right)+i \pi ^2\right)
\end{multline}
\end{example}
\begin{figure}[H]
\includegraphics[scale=0.5]{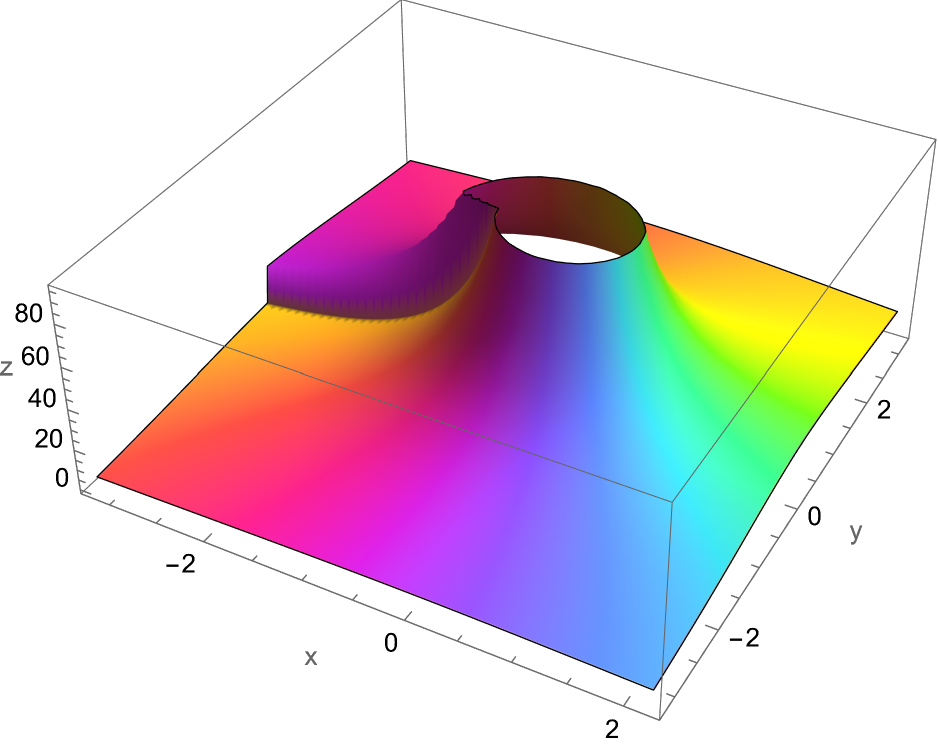}
\caption{Plot of  $f(p)=\frac{\pi  \left(-8 C-2 \pi  \log \left(\frac{\pi }{p}\right)+i \pi ^2\right)}{p}$, $p\in\mathbb{C}$.}
   \label{fig:fig2}
\end{figure}
\vspace{-6pt}
\begin{example}
\begin{multline}
\int_{0}^{\infty}x^{s-1} \tan ^{-1}(\sinh (\alpha ) \text{sech}(b x))dx\\
=\frac{i 2^{s-1} \pi ^{s+1} \left(\frac{1}{b}\right)^s \csc \left(\frac{\pi  s}{2}\right)}{s}
   \left(\zeta \left(-s,\frac{\pi -2 i \alpha }{4 \pi }\right)-\zeta \left(-s,\frac{2 i \alpha +\pi }{4 \pi }\right)\right. \\ \left.
-\zeta \left(-s,\frac{3}{4}-\frac{i\alpha }{2 \pi }\right)+\zeta \left(-s,\frac{i \alpha }{2 \pi }+\frac{3}{4}\right)\right)
\end{multline}
\end{example}
\begin{example}
\begin{multline}
\int_{0}^{\infty}\left((\log (a)-i x)^{k-1}+(\log (a)+i x)^{k-1}\right) \tan ^{-1}(\sinh (\alpha ) \text{sech}(b x))dx\\
=\frac{i 2^k \pi ^{k+1}
   \left(\frac{1}{b}\right)^k}{k} \left(\zeta \left(-k,\frac{-2 i \alpha +2 b \log (a)+\pi }{4 \pi }\right)-\zeta \left(-k,\frac{-2 i \alpha +2 b \log (a)+3\pi }{4 \pi }\right)\right. \\ \left.
-\zeta \left(-k,\frac{2 i \alpha +2 b \log (a)+\pi }{4 \pi }\right)+\zeta \left(-k,\frac{2 i \alpha +2 b \log (a)+3 \pi }{4 \pi }\right)\right)
\end{multline}
\end{example}
\begin{example}
\begin{multline}
\int_{0}^{\infty}\left(\frac{1}{(\log (a)+i x)^2}+\frac{1}{(\log (a)-i x)^2}\right) \tan
   ^{-1}(\sinh (\alpha ) \text{sech}(b x))dx\\
=\frac{1}{2} i b \left(\psi
   ^{(0)}\left(\frac{-2 i \alpha +2 b \log (a)+\pi }{4 \pi }\right)-\psi
   ^{(0)}\left(\frac{-2 i \alpha +2 b \log (a)+3 \pi }{4 \pi }\right)\right. \\ \left.
-\psi
   ^{(0)}\left(\frac{2 i \alpha +2 b \log (a)+\pi }{4 \pi }\right)+\psi
   ^{(0)}\left(\frac{2 i \alpha +2 b \log (a)+3 \pi }{4 \pi
   }\right)\right)
\end{multline}
\end{example}
\begin{example}
\begin{equation}
\int_{0}^{\infty}\frac{\left(12 x^2-1\right) \tanh ^{-1}(\text{sech}(\pi  x))}{\left(4
   x^2+1\right)^3}dx=\frac{1}{8} (\pi -2 \pi  C)
\end{equation}
\end{example}
\begin{example}
\begin{equation}
\int_{0}^{\infty}\frac{\left(x^4-6 x^2+1\right) \tanh ^{-1}(\text{sech}(\pi 
   x))}{\left(x^2+1\right)^4}dx=\frac{1}{12} \pi  (3 \zeta (3)-2)
\end{equation}
\end{example}
\begin{example}
\begin{multline}
\int_{0}^{\infty}\left(\frac{1}{(1+i x)^6}-\frac{1}{(x+i)^6}\right) \tanh
   ^{-1}(\text{sech}(\pi  x))dx=\frac{1}{40} \pi  (15 \zeta (5)-8)
\end{multline}
\end{example}
\begin{example}
\begin{multline}
\int_{0}^{\infty}\left(\frac{1}{(1+i x)^{3/2}}+\frac{1}{(1-i x)^{3/2}}\right) \tanh
   ^{-1}(\text{sech}(\pi  x))dx=-4 \left(\sqrt{2}-1\right) \pi  \zeta
   \left(\frac{1}{2}\right)-2 \pi
\end{multline}
\end{example}
\begin{example}
\begin{multline}
\int_{0}^{\infty}\left(\frac{1}{\sqrt{1+i x}}+\frac{1}{\sqrt{1-i x}}\right) \tanh
   ^{-1}(\text{sech}(\pi  x))dx=\left(1-2 \sqrt{2}\right) \zeta
   \left(\frac{3}{2}\right)+2 \pi
\end{multline}
\end{example}
\begin{example}
\begin{multline}
\int_{0}^{\infty}\left(\frac{\log (\log (a)-i x)}{\log (a)-i x}+\frac{\log (\log (a)+i x)}{\log (a)+i x}\right) \tan
   ^{-1}(\sinh (\alpha ) \text{sech}(b x))dx\\
=\frac{1}{2} i \pi  \left(\zeta''\left(0,\frac{2 b \log
   (a)-2 i \alpha +\pi }{4 \pi }\right)-\zeta''\left(0,\frac{2 b \log (a)-2 i \alpha +3 \pi }{4
   \pi }\right)\right. \\ \left.
-\zeta''\left(0,\frac{2 b \log (a)+2 i \alpha +\pi }{4 \pi
   }\right)+\zeta''\left(0,\frac{2 b \log (a)+2 i \alpha +3 \pi }{4 \pi }\right)\right. \\ \left.
+2 \log\left(\frac{2 \pi }{b}\right) \left(\log (-4 (-2 b \log (a)+2 i \alpha +\pi ))+\log (2 b \log (a)+2 i \alpha -3
   \pi )\right.\right. \\ \left.\left.
-\log (2 b \log (a)+2 i \alpha -\pi )-\log (8 b \log (a)-8 i \alpha -12 \pi )\right.\right. \\ \left.\left.
-\text{log$\Gamma
   $}\left(-\frac{-2 i \alpha -2 b \log (a)+\pi }{4 \pi }\right)+\text{log$\Gamma $}\left(-\frac{2 i \alpha -2 b \log
   (a)+\pi }{4 \pi }\right)\right.\right. \\ \left.\left.
-\text{log$\Gamma $}\left(-\frac{2 i \alpha -2 b \log (a)+3 \pi }{4 \pi
   }\right)+\text{log$\Gamma $}\left(\frac{2 i \alpha +2 b \log (a)-3 \pi }{4 \pi }\right)\right)\right)
\end{multline}
\end{example}
\begin{example}
\begin{multline}
\int_{0}^{\infty}\frac{\tan ^{-1}(\sinh (\alpha ) \text{sech}(b x))}{\log
   ^2(a)+x^2}dx\\
   =-\frac{i \pi }{2 \log (a)} \log \left(\frac{\Gamma \left(\frac{-2 i \alpha +2 b
   \log (a)+\pi }{4 \pi }\right) \Gamma \left(\frac{2 i \alpha +2 b \log (a)+3
   \pi }{4 \pi }\right)}{\Gamma \left(\frac{-2 i \alpha +2 b \log (a)+3 \pi }{4
   \pi }\right) \Gamma \left(\frac{2 i \alpha +2 b \log (a)+\pi }{4 \pi
   }\right)}\right)
\end{multline}
\end{example}
\begin{example}
\begin{equation}
\int_{0}^{\infty}\frac{\tanh ^{-1}(\text{sech}(\pi  x))}{x^2+1}dx=-\frac{1}{2} \pi  \log
   \left(\frac{2}{\pi }\right)
\end{equation}
\end{example}
\begin{example}
\begin{equation}
\int_{0}^{\infty}\frac{\tanh ^{-1}(\sin (\alpha ) \text{sech}(b x))}{a^2+x^2}dx=-\frac{\pi }{2 a}
   \log \left(\frac{\Gamma \left(\frac{2 a b-2 \alpha +3 \pi }{4 \pi }\right)
   \Gamma \left(\frac{2 a b+2 \alpha +\pi }{4 \pi }\right)}{\Gamma \left(\frac{2
   a b-2 \alpha +\pi }{4 \pi }\right) \Gamma \left(\frac{2 a b+2 \alpha +3 \pi
   }{4 \pi }\right)}\right)
\end{equation}
\end{example}
\begin{example}
\begin{multline}
\int_{0}^{\infty}\frac{\tanh ^{-1}(\sin (\alpha ) \text{sech}(b
   x))}{\left(a^2+x^2\right)^2}dx\\
=\frac{a b}{8 a^3} \left(-\psi ^{(0)}\left(\frac{2 a b-2
   \alpha +\pi }{4 \pi }\right)+\psi ^{(0)}\left(\frac{2 a b-2 \alpha +3 \pi }{4
   \pi }\right)+\psi ^{(0)}\left(\frac{2 a b+2 \alpha +\pi }{4 \pi }\right)\right. \\ \left.
-\psi
   ^{(0)}\left(\frac{2 a b+2 \alpha +3 \pi }{4 \pi }\right)\right)-2 \pi  \log
   \left(\frac{\Gamma \left(\frac{2 a b-2 \alpha +3 \pi }{4 \pi }\right) \Gamma
   \left(\frac{2 a b+2 \alpha +\pi }{4 \pi }\right)}{\Gamma \left(\frac{2 a b-2
   \alpha +\pi }{4 \pi }\right) \Gamma \left(\frac{2 a b+2 \alpha +3 \pi }{4 \pi
   }\right)}\right)
\end{multline}
\end{example}
\begin{example}
\begin{equation}
\int_{0}^{\infty}\frac{\tanh ^{-1}(\text{sech}(x))}{\left(x^2+\pi ^2\right)^2}dx=\frac{\log
   \left(\frac{2 \pi }{e}\right)}{4 \pi ^2}
\end{equation}
\end{example}
\begin{example}
\begin{equation}
\int_{0}^{\infty}\frac{\coth ^{-1}\left(\sqrt{2} \cosh (9 \pi 
   x)\right)}{x^2+4}dx=-\frac{1}{4} \pi  \log \left(\frac{\Gamma
   \left(\frac{75}{8}\right) \Gamma \left(\frac{77}{8}\right)}{\Gamma
   \left(\frac{73}{8}\right) \Gamma \left(\frac{79}{8}\right)}\right)
\end{equation}
\end{example}
\begin{example}
\begin{multline}
\int_{0}^{\infty}\text{sech}(t) \tanh ^{-1}(\sin (\alpha ) \text{sech}(a b \sinh
   (t)))dx\\
   =-\frac{1}{2} \pi  \log \left(\frac{\Gamma \left(\frac{2 a b-2 \alpha +3
   \pi }{4 \pi }\right) \Gamma \left(\frac{2 a b+2 \alpha +\pi }{4 \pi
   }\right)}{\Gamma \left(\frac{2 a b-2 \alpha +\pi }{4 \pi }\right) \Gamma
   \left(\frac{2 a b+2 \alpha +3 \pi }{4 \pi }\right)}\right)
\end{multline}
\end{example}
\begin{example}
\begin{multline}
\exp\left(\int_{0}^{\infty}\frac{2 i \log (a) \tan ^{-1}(\sinh (\alpha ) \text{sech}(b x))}{\pi 
   \left(\log ^2(a)+x^2\right)}dx\right)\\
   =\frac{\Gamma \left(\frac{-2 i \alpha +2 b \log
   (a)+\pi }{4 \pi }\right) \Gamma \left(\frac{2 i \alpha +2 b \log (a)+3 \pi
   }{4 \pi }\right)}{\Gamma \left(\frac{-2 i \alpha +2 b \log (a)+3 \pi }{4 \pi
   }\right) \Gamma \left(\frac{2 i \alpha +2 b \log (a)+\pi }{4 \pi
   }\right)}
\end{multline}
\end{example}
\begin{example}
\begin{multline}
\int_{0}^{\infty}\frac{2 (m x \sin (m x)+\cos (m x)-r x \sin (r x)-\cos (r x)) \tan
   ^{-1}(\sinh (\alpha ) \text{sech}(b x))}{x^2}dx\\
=i b \left(e^{\frac{(\pi -2 i
   \alpha ) m}{2 b}} \Phi \left(-e^{\frac{m \pi }{b}},1,\frac{1}{2}-\frac{i
   \alpha }{\pi }\right)-e^{\frac{(\pi +2 i \alpha ) m}{2 b}} \Phi
   \left(-e^{\frac{m \pi }{b}},1,\frac{i \alpha }{\pi
   }+\frac{1}{2}\right)\right. \\ \left.
-e^{\frac{(\pi -2 i \alpha ) r}{2 b}} \Phi
   \left(-e^{\frac{\pi  r}{b}},1,\frac{1}{2}-\frac{i \alpha }{\pi
   }\right)+e^{\frac{(\pi +2 i \alpha ) r}{2 b}} \Phi \left(-e^{\frac{\pi 
   r}{b}},1,\frac{i \alpha }{\pi }+\frac{1}{2}\right)\right)
\end{multline}
\end{example}
\begin{example}
\begin{multline}
\int_{0}^{\infty}\frac{x^{-m-\frac{1}{2}} \left(\beta +\sqrt{x (2 \beta +x)}+x\right)^{-m} \log \left(x \left(\beta +\sqrt{x (2 \beta
   +x)}+x\right)\right)}{\sqrt{2 \beta +x}}dx\\
=\pi  2^m \beta ^{-2 m} \csc (2 \pi  m) (2 \log (\beta )+2 \pi  \cot (2 \pi  m)-\log (2))
\end{multline}
\end{example}
\begin{example}
\begin{multline}
\int_{0}^{\infty}\frac{\log ^k\left(\frac{1}{2 x \left(x+\sqrt{x (x+2)}+1\right)}\right)}{x^{3/4} \sqrt{x+2} \sqrt[4]{x+\sqrt{x (x+2)}+1}}dx\\
=-i 2^{2 k+\frac{5}{4}} e^{\frac{1}{2} i \pi  (k+1)} \pi
   ^{k+1} \left(2^k \zeta \left(-k,\frac{1}{4}\right)-2^k \zeta \left(-k,\frac{3}{4}\right)\right)
\end{multline}
\end{example}
\begin{example}
\begin{multline}
\int_{0}^{\infty}\frac{\log \left(\log \left(\frac{1}{2 x \left(x+\sqrt{x (x+2)}+1\right)}\right)\right)}{x^{3/4} \sqrt{x+2} \sqrt[4]{x+\sqrt{x (x+2)}+1}}dx\\
=\frac{\pi  \left(-4 \text{log$\Gamma
   $}\left(-\frac{3}{4}\right)+4 \text{log$\Gamma $}\left(-\frac{1}{4}\right)+i \pi +\log \left(\frac{64 \pi ^2}{81}\right)\right)}{2^{3/4}}
\end{multline}
\end{example}
\begin{example}
\begin{multline}
\int_{0}^{\infty}\frac{\log \left(\frac{1}{2 x \left(x+\sqrt{x (x+2)}+1\right)}\right) \log \left(\log \left(\frac{1}{2 x \left(x+\sqrt{x (x+2)}+1\right)}\right)\right)}{x^{3/4} \sqrt{x+2}
   \sqrt[4]{x+\sqrt{x (x+2)}+1}}dx=-8 i \sqrt[4]{2} \pi  C
\end{multline}
\end{example}
\begin{example}
\begin{multline}
\int_{0}^{\infty}\frac{\log \left(\log \left(-\frac{1}{2 x \left(x+\sqrt{x (x+2)}+1\right)}\right)\right)}{x^{3/4} \sqrt{x+2} \sqrt[4]{x+\sqrt{x (x+2)}+1} \log \left(-\frac{1}{2 x
   \left(x+\sqrt{x (x+2)}+1\right)}\right)}dx\\
=\frac{(-\pi +2 i \log (8 \pi )) \left(\psi ^{(0)}\left(\frac{3}{8}\right)-\psi ^{(0)}\left(\frac{7}{8}\right)\right)-2 i \left(\gamma
   _1\left(\frac{3}{8}\right)-\gamma _1\left(\frac{7}{8}\right)\right)}{4\ 2^{3/4}}
\end{multline}
\end{example}
\begin{example}
\begin{multline}
\int_{0}^{\infty}\frac{dx}{x^{3/4} \sqrt{x+2} \sqrt[4]{x+\sqrt{x (x+2)}+1} \left(\log ^2\left(2 x \left(x+\sqrt{x (x+2)}+1\right)\right)+\pi ^2\right)}\\
=\frac{\pi +\log \left(3-2
   \sqrt{2}\right)}{2 \sqrt[4]{2} \pi }
\end{multline}
\end{example}
\begin{example}
\begin{multline}
\int_{0}^{\infty}\frac{dx}{2 x^{3/4} \sqrt{x+2} \sqrt[4]{x+\sqrt{x (x+2)}+1} \left(\log ^2\left(2 x \left(x+\sqrt{x (x+2)}+1\right)\right)+\frac{\pi ^2}{4}\right)}\\
=\frac{\sin \left(\frac{\pi
   }{8}\right) \left(\sqrt{2} \pi +2 \log \left(\tan \left(\frac{\pi }{16}\right)\right)\right)+2 \cos \left(\frac{\pi }{8}\right) \log \left(\cot \left(\frac{3 \pi
   }{16}\right)\right)}{2^{3/4} \pi }
\end{multline}
\end{example}
\begin{example}
\begin{multline}
\int_0^{\infty } \frac{\cosh (c u)}{\cosh (\pi
    u)} \log \left(a^2+u^2\right) \, du\\
=\frac{1}{2} e^{-\frac{1}{2} (i c)} \left(\left(-1+\frac{1}{1+e^{ic}}\right) \log (4)-2 \Phi'\left(-e^{-ic},0,\frac{1}{2}+a\right)\right. \\ \left.
-e^{i c}\left(\Phi'\left(-i e^{\frac{i c}{2}},0,1+2a\right)+\Phi'\left(e^{\frac{1}{2} i (c+\pi)},0,1+2 a\right)\right)\right)
\end{multline}
\end{example}
\begin{example}
\begin{multline}
\int_0^{\infty } \frac{\cosh (c u)}{\cosh (b
   u)} \log \left(a^2+u^2\right) \, du\\
=\frac{e^{-\frac{i c \pi }{2 b}} \pi  }{2 b}\left(2 \log \left(\frac{\pi}{\sqrt{2} b}\right) \left(1+i \tan \left(\frac{c \pi }{2 b}\right)\right)-2\Phi'\left(-e^{-\frac{i c \pi}{b}},0,\frac{1}{2}+\frac{a b}{\pi }\right)\right. \\ \left.
-e^{\frac{i c \pi }{b}}
   \left(\Phi'\left(-i e^{\frac{i c \pi }{2 b}},0,1+\frac{2 a b}{\pi }\right)+\Phi'\left(e^{\frac{i (b+c) \pi }{2 b}},0,1+\frac{2 a b}{\pi }\right)\right)\right)
\end{multline}
\end{example}
\begin{example}
\begin{multline}
\int_0^{\infty } \frac{ \cosh (c u) \sinh (m u)}{\cosh (b u)}\tan ^{-1}\left(\frac{a}{u}\right) \,
   du\\
=\frac{\pi }{16 b} \left(-H_{-\frac{b+c-m}{4 b}}+H_{-\frac{3 b+c-m}{4 b}}+H_{-\frac{b-c+m}{4 b}}-H_{-\frac{3 b-c+m}{4
   b}}+H_{\frac{-3 b+c+m}{4 b}}-H_{\frac{-b+c+m}{4 b}}+H_{-\frac{b+c+m}{4 b}}\right. \\ \left.
-H_{-\frac{3 b+c+m}{4 b}}-2 i \log (2)
   \left(\sec \left(\frac{(c-m) \pi }{2 b}\right)-\sec \left(\frac{(c+m) \pi }{2 b}\right)\right)\right. \\ \left.
-4 e^{\frac{i (b-c-m) \pi }{2 b}} \Phi'\left(e^{\frac{i (b-c-m) \pi }{b}},0,\frac{1}{2}+\frac{a
   b}{\pi }\right)-2 e^{\frac{i (b+c-m) \pi }{2 b}} \left(\Phi'\left(-e^{\frac{i (b+c-m)
   \pi }{2 b}},0,1+\frac{2 a b}{\pi }\right)\right.\right. \\ \left.\left.
+\Phi'\left(e^{\frac{i (b+c-m) \pi }{2b}},0,1+\frac{2 a b}{\pi }\right)\right)+4 e^{\frac{i (b-c+m) \pi }{2 b}}
   \Phi'\left(e^{\frac{i (b-c+m) \pi }{b}},0,\frac{1}{2}+\frac{a b}{\pi }\right)\right. \\ \left.
+2e^{\frac{i (b+c+m) \pi }{2 b}} \left(\Phi'\left(-e^{\frac{i (b+c+m) \pi }{ b}},0,1+\frac{2 a b}{\pi }\right)+\Phi'\left(e^{\frac{i (b+c+m) \pi }{2 b}},0,1+\frac{2 a b}{\pi }\right)\right)\right)
\end{multline}
\end{example}
\begin{example}
\begin{multline}
\int_0^{\infty } \frac{\cosh (c u) \sinh (m u)}{\cosh (b u) } \, \frac{du}{u}=\frac{1}{2} \log \left(\tan
   \left(\frac{(b-c+m) \pi }{4 b}\right) \tan \left(\frac{(b+c+m) \pi }{4 b}\right)\right)
\end{multline}
\end{example}
\begin{example}
\begin{multline}
\int_0^1 \left( \frac{\cosh (c u) \sinh (m u)}{\cosh (b u)}+\cosh \left(\frac{c}{u}\right) \text{sech}\left(\frac{b}{u}\right) \sinh \left(\frac{m}{u}\right)\right) \,
   \frac{du}{u}\\
=\frac{1}{2} \log \left(\tan \left(\frac{(b-c+m) \pi }{4 b}\right) \tan \left(\frac{(b+c+m) \pi }{4 b}\right)\right)
\end{multline}
\end{example}
\begin{example}
\begin{multline}
\int_0^{\infty } \frac{\cosh (c u) \sinh (m u)}{\cosh (b u) } \, \frac{du}{u}\\
=\frac{1}{4} \left(2 \tanh
   ^{-1}\left(e^{\frac{i (b-c-m) \pi }{2 b}}\right)-2 \tanh ^{-1}\left(e^{\frac{i (b-c+m) \pi }{2 b}}\right)+\log
   \left(\frac{b+c-m}{b}\right)\right. \\ \left.
+2 \log \left(\frac{-5 b-c+m}{b}\right)-\log \left(\frac{-3 b-c+m}{b}\right)-\log
   \left(\frac{b+c+m}{b}\right)\right. \\ \left.
+\log \left(-\frac{3 b+c+m}{b}\right)-2 \log \left(-\frac{5 b+c+m}{b}\right)+2
   \text{log$\Gamma $}\left(\frac{b+c-m}{4 b}\right)\right. \\ \left.
-\text{log$\Gamma $}\left(\frac{b+c-m}{2
   b}\right)-\text{log$\Gamma $}\left(-\frac{3 b+c-m}{2 b}\right)+2 \text{log$\Gamma $}\left(-\frac{5 b+c-m}{4
   b}\right)\right. \\ \left.
-2 \text{log$\Gamma $}\left(\frac{b+c+m}{4 b}\right)+\text{log$\Gamma $}\left(\frac{b+c+m}{2
   b}\right)+\text{log$\Gamma $}\left(-\frac{3 b+c+m}{2 b}\right)\right. \\ \left.
-2 \text{log$\Gamma $}\left(-\frac{5 b+c+m}{4 b}\right)\right)
\end{multline}
\end{example}
\begin{example}
\begin{multline}
\int_0^{\infty }u^{-1+s}  \frac{\cosh (c u)}{\cosh (b u)} \, du\\
=2^{-1-s} \left(\frac{i}{b}\right)^{1+s} b \pi
   ^s \left(i 2^s e^{-\frac{i c \pi }{2 b}} \Phi \left(-e^{-\frac{i c \pi
   }{b}},1-s,\frac{1}{2}\right)-\text{Li}_{1-s}\left(-i e^{\frac{i c \pi }{2
   b}}\right)+\text{Li}_{1-s}\left(e^{\frac{i (b+c) \pi }{2 b}}\right)\right)\\
 \left(i+\tan \left(\frac{1}{2} \pi 
   (-1+s)\right)\right)
\end{multline}
\end{example}
\begin{example}
\begin{multline}
\int_0^{\infty } \frac{u }{\left(a^2+u^2\right) }\frac{\cosh (c u) \sinh (m u)}{\cosh (b u)} \, du\\
=\frac{1}{4} \left(e^{\frac{i (b-c-m) \pi }{2 b}} \Phi \left(e^{\frac{i (b-c-m) \pi
   }{b}},1,\frac{1}{2}+\frac{a b}{\pi }\right)+e^{\frac{i (b+c-m) \pi }{2 b}} \left(\Phi \left(-e^{\frac{i (b+c-m) \pi }{2 b}},1,1+\frac{2 a b}{\pi }\right)\right.\right. \\ \left.\left.
+\Phi \left(e^{\frac{i (b+c-m)\pi }{2 b}},1,1+\frac{2 a b}{\pi }\right)\right)-e^{\frac{i (b-c+m) \pi }{2 b}} \Phi \left(e^{\frac{i (b-c+m) \pi }{b}},1,\frac{1}{2}+\frac{a b}{\pi }\right)\right. \\ \left.
-e^{\frac{i (b+c+m) \pi}{2 b}} \left(\Phi \left(-e^{\frac{i (b+c+m) \pi }{2 b}},1,1+\frac{2 a b}{\pi }\right)+\Phi \left(e^{\frac{i (b+c+m) \pi }{2 b}},1,1+\frac{2 a b}{\pi }\right)\right)\right)
\end{multline}
\end{example}
\begin{example}
\begin{multline}
\int_0^{\infty } \frac{u }{\left(a^4-u^4\right) }\frac{\cosh (c u) \sinh (m u)}{\cosh (b u)} \, du\\
=-\frac{1}{8 a^2}\left(e^{\frac{i (b-c-m) \pi }{2 b}} \Phi \left(e^{\frac{i (b-c-m) \pi}{b}},1,\frac{1}{2}+\frac{i a b}{\pi }\right)
-e^{\frac{i (b-c-m) \pi }{2 b}} \Phi \left(e^{\frac{i (b-c-m) \pi }{b}},1,\frac{1}{2}+\frac{a b}{\pi }\right)\right. \\ \left.
+e^{\frac{i (b+c-m) \pi }{2b}} \left(\Phi \left(-e^{\frac{i (b+c-m) \pi }{2 b}},1,1+\frac{2 i a b}{\pi }\right)+\Phi \left(e^{\frac{i (b+c-m) \pi }{2 b}},1,1+\frac{2 i a b}{\pi }\right)\right)\right. \\ \left.
-e^{\frac{i (b+c-m) \pi }{2 b}} \left(\Phi \left(-e^{\frac{i (b+c-m) \pi }{2 b}},1,1+\frac{2 a b}{\pi }\right)+\Phi \left(e^{\frac{i (b+c-m) \pi }{2 b}},1,1+\frac{2 a b}{\pi}\right)\right)\right. \\ \left.
-e^{\frac{i (b-c+m) \pi }{2 b}} \Phi \left(e^{\frac{i (b-c+m) \pi }{b}},1,\frac{1}{2}+\frac{i a b}{\pi }\right)+e^{\frac{i (b-c+m) \pi }{2 b}} \Phi \left(e^{\frac{i(b-c+m) \pi }{b}},1,\frac{1}{2}+\frac{a b}{\pi }\right)\right. \\ \left.
-e^{\frac{i (b+c+m) \pi }{2 b}} \left(\Phi \left(-e^{\frac{i (b+c+m) \pi }{2 b}},1,1+\frac{2 i a b}{\pi }\right)+\Phi\left(e^{\frac{i (b+c+m) \pi }{2 b}},1,1+\frac{2 i a b}{\pi }\right)\right)\right. \\ \left.
+e^{\frac{i (b+c+m) \pi }{2 b}} \left(\Phi \left(-e^{\frac{i (b+c+m) \pi }{2 b}},1,1+\frac{2 a b}{\pi
   }\right)+\Phi \left(e^{\frac{i (b+c+m) \pi }{2 b}},1,1+\frac{2 a b}{\pi }\right)\right)\right)
\end{multline}
\end{example}
\begin{example}
\begin{multline}
\int_0^{\infty } \frac{\tan ^{-1}\left(\frac{x}{a}\right) \sinh (m x)}{(\cos (t)+\cosh (c x))^2} \, dx\\
=-\frac{ie^{-\frac{i m (\pi -t)}{c}} \csc ^2(t) }{4 c^2}\left(-c \Phi \left(e^{-\frac{2 i m \pi }{c}},1,\frac{a c+\pi -t}{2 \pi }\right)-ce^{-\frac{2 i m t}{c}} \Phi \left(e^{-\frac{2 i m \pi }{c}},1,\frac{a c+\pi +t}{2 \pi }\right)\right. \\ \left.
+2 \pi  (-i m+c \cot (t)) \Phi'\left(e^{-\frac{2 i m \pi }{c}},0,\frac{a c+\pi -t}{2 \pi }\right)\right. \\ \left.
-2 i e^{-\frac{2 i mt}{c}} \pi  (m-i c \cot (t)) \Phi'\left(e^{-\frac{2 i m \pi }{c}},0,\frac{a c+\pi +t}{2 \pi
   }\right)+e^{\frac{2 i m (\pi -t)}{c}} \left(c \Phi \left(e^{\frac{2 i m \pi }{c}},1,\frac{a c+\pi -t}{2 \pi }\right)\right.\right. \\ \left.\left.
-2 i \pi(m-i c \cot (t)) \Phi'\left(e^{\frac{2 i m \pi }{c}},0,\frac{a c+\pi -t}{2 \pi
   }\right)+e^{\frac{2 i m t}{c}} \left(c \Phi \left(e^{\frac{2 i m \pi }{c}},1,\frac{a c+\pi +t}{2 \pi }\right)\right.\right.\right. \\ \left.\left.\left.
+2 \pi  (-i m+c\cot (t)) \Phi'\left(e^{\frac{2 i m \pi }{c}},0,\frac{a c+\pi +t}{2 \pi}\right)\right)\right)\right)
\end{multline}
\end{example}
\begin{example}
\begin{multline}
\int_0^{\infty } \frac{\cosh (m x) }{(\cos (t)+\cosh (c x))^2} \log \left(a^2+x^2\right)\, dx\\
=\frac{e^{-\frac{i m (\pi+t)}{c}} \csc ^2(t) }{2 c^2}\left(-c e^{\frac{2 i m t}{c}} \Phi \left(e^{-\frac{2 i m \pi }{c}},1,\frac{a c+\pi -t}{2 \pi }\right)-c\Phi \left(e^{-\frac{2 i m \pi }{c}},1,\frac{a c+\pi +t}{2 \pi }\right)\right. \\ \left.
-c e^{\frac{2 i m \pi }{c}} \Phi \left(e^{\frac{2 i m\pi }{c}},1,\frac{a c+\pi -t}{2 \pi }\right)+4 e^{\frac{i m (\pi +t)}{c}} \pi  \csc \left(\frac{m \pi }{c}\right) (\log(2)\right. \\ \left.
-\log (c)+\log (\pi )) \left(m \cos \left(\frac{m t}{c}\right)-c \cot (t) \sin \left(\frac{m t}{c}\right)\right)\right. \\ \left.
-2 i e^{\frac{2 i m t}{c}} m \pi  \Phi'\left(e^{-\frac{2 i m \pi }{c}},0,\frac{a c+\pi -t}{2 \pi }\right)+2 c e^{\frac{2 i m t}{c}} \pi  \cot (t) \Phi'\left(e^{-\frac{2 i m \pi }{c}},0,\frac{ac+\pi -t}{2 \pi }\right)\right. \\ \left.
-2 i m \pi  \Phi'\left(e^{-\frac{2 i m \pi }{c}},0,\frac{a c+\pi +t}{2\pi }\right)-2 c \pi  \cot (t) \Phi'\left(e^{-\frac{2 i m \pi }{c}},0,\frac{a c+\pi +t}{2 \pi}\right)\right. \\ \left.
+2 i e^{\frac{2 i m \pi }{c}} m \pi  \Phi'\left(e^{\frac{2 i m \pi }{c}},0,\frac{a c+\pi-t}{2 \pi }\right)+2 c e^{\frac{2 i m \pi }{c}} \pi  \cot (t) \Phi'\left(e^{\frac{2 i m \pi}{c}},0,\frac{a c+\pi -t}{2 \pi }\right)\right. \\ \left.
-e^{\frac{2 i m (\pi +t)}{c}} \left(c \Phi \left(e^{\frac{2 i m \pi }{c}},1,\frac{ac+\pi +t}{2 \pi }\right)+2 \pi  (-i m+c \cot (t)) \Phi'\left(e^{\frac{2 i m \pi }{c}},0,\frac{ac+\pi +t}{2 \pi }\right)\right)\right)
\end{multline}
\end{example}
\begin{example}
\begin{multline}
\int_0^{\infty } \frac{\sinh (m x)}{x (\cos (t)+\cosh (c x))^2} \, dx\\
=-\frac{i e^{-\frac{i m (\pi -t)}{c}} \csc ^2(t)}{8 c \pi }
   \left(2 \pi  (-i m+c \cot (t)) \Phi \left(e^{-\frac{2 i m \pi }{c}},1,\frac{\pi -t}{2 \pi }\right)\right. \\ \left.
-2 e^{-\frac{2 i m t}{c}}
   \pi  (i m+c \cot (t)) \Phi \left(e^{-\frac{2 i m \pi }{c}},1,\frac{\pi +t}{2 \pi }\right)-c \Phi \left(e^{-\frac{2 i m \pi
   }{c}},2,\frac{\pi -t}{2 \pi }\right)\right. \\ \left.
-c e^{-\frac{2 i m t}{c}} \Phi \left(e^{-\frac{2 i m \pi }{c}},2,\frac{\pi +t}{2 \pi
   }\right)+e^{\frac{2 i m (\pi -t)}{c}} \left(-2 \pi  (i m+c \cot (t)) \Phi \left(e^{\frac{2 i m \pi }{c}},1,\frac{\pi -t}{2
   \pi }\right)\right.\right. \\ \left.\left.
+2 e^{\frac{2 i m t}{c}} \pi  (-i m+c \cot (t)) \Phi \left(e^{\frac{2 i m \pi }{c}},1,\frac{\pi +t}{2 \pi
   }\right)+c \left(\Phi \left(e^{\frac{2 i m \pi }{c}},2,\frac{\pi -t}{2 \pi }\right)\right.\right.\right. \\ \left.\left.\left.
+e^{\frac{2 i m t}{c}} \Phi
   \left(e^{\frac{2 i m \pi }{c}},2,\frac{\pi +t}{2 \pi }\right)\right)\right)\right)
\end{multline}
\end{example}
\begin{example}
\begin{multline}
\int_0^{\infty } \frac{x^{-1+s} \sinh (m x)}{(\cos (t)+\cosh (c x))^2} \, dx\\
=-i 2^{-3+s} c^{-1-s} e^{\frac{i m (-\pi
   +t)}{c}} \pi ^{-1+s} \csc ^2(t) \sec \left(\frac{\pi  s}{2}\right) \left(2 \pi  (-i m+c \cot (t)) \Phi \left(e^{-\frac{2 i
   m \pi }{c}},1-s,\frac{\pi -t}{2 \pi }\right)\right. \\ \left.
-2 e^{-\frac{2 i m t}{c}} \pi  (i m+c \cot (t)) \Phi \left(e^{-\frac{2 i m \pi
   }{c}},1-s,\frac{\pi +t}{2 \pi }\right)+c (-1+s) \Phi \left(e^{-\frac{2 i m \pi }{c}},2-s,\frac{\pi -t}{2 \pi }\right)\right. \\ \left.
+ce^{-\frac{2 i m t}{c}} (-1+s) \Phi \left(e^{-\frac{2 i m \pi }{c}},2-s,\frac{\pi +t}{2 \pi }\right)-e^{\frac{2 i m (\pi-t)}{c}} \left(2 \pi  (i m+c \cot (t)) \Phi \left(e^{\frac{2 i m \pi }{c}},1-s,\frac{\pi -t}{2 \pi }\right)\right.\right. \\ \left.\left.
+2 i e^{\frac{2 i m t}{c}} \pi  (m+i c \cot (t)) \Phi \left(e^{\frac{2 i m \pi }{c}},1-s,\frac{\pi +t}{2 \pi }\right)+c (-1+s) \left(\Phi \left(e^{\frac{2 i m \pi }{c}},2-s,\frac{\pi -t}{2 \pi }\right)\right.\right.\right. \\ \left.\left.\left.
+e^{\frac{2 i m t}{c}} \Phi \left(e^{\frac{2 i m \pi}{c}},2-s,\frac{\pi +t}{2 \pi }\right)\right)\right)\right)
\end{multline}
\end{example}
\begin{example}
\begin{multline}
\int_0^{\infty } \frac{\left((a-x)^k-(a+x)^k\right) (\sinh (b x) \sinh (m x))}{\sinh (c \pi  x)} \, dx\\
=\frac{1}{2}
   \left(\frac{i}{c}\right)^{1+k} e^{-\frac{i (b+m)}{c}} \left(-e^{\frac{2 i b}{c}} \Phi \left(-e^{\frac{i (b-m)}{c}},-k,1-i a
   c\right)-e^{\frac{2 i m}{c}} \Phi \left(-e^{\frac{i (-b+m)}{c}},-k,1-i a c\right)\right. \\ \left.
+\Phi \left(-e^{-\frac{i (b+m)}{c}},-k,1-i
   a c\right)+e^{\frac{2 i (b+m)}{c}} \Phi \left(-e^{\frac{i (b+m)}{c}},-k,1-i a c\right)\right)
\end{multline}
\end{example}
\begin{example}
\begin{multline}
\int_0^{\infty } \frac{\left((a-x)^k+(a+x)^k\right) (\cosh (m x) \sinh (b x))}{\sinh (c \pi  x)} \, dx\\
=\frac{1}{2}
   \left(\frac{i}{c}\right)^{1+k} e^{-\frac{i (b+m)}{c}} \left(-e^{\frac{2 i b}{c}} \Phi \left(-e^{\frac{i (b-m)}{c}},-k,1-i a
   c\right)+\Phi \left(-e^{-\frac{i (b+m)}{c}},-k,1-i a c\right)\right. \\ \left.
+e^{\frac{2 i m}{c}} \left(\Phi \left(-e^{\frac{i
   (-b+m)}{c}},-k,1-i a c\right)-e^{\frac{2 i b}{c}} \Phi \left(-e^{\frac{i (b+m)}{c}},-k,1-i a c\right)\right)\right)
\end{multline}
\end{example}
\begin{example}
\begin{multline}
\int_0^{\infty } \frac{\tanh ^{-1}\left(\frac{x}{a}\right) (\sinh (b x) \sinh (m x))}{\sinh (c \pi  x)} \, dx\\
=-\frac{i
   e^{-\frac{i (b+m)}{c}}}{4 c} \left(e^{\frac{2 i b}{c}} \Phi'\left(-e^{\frac{i (b-m)}{c}},0,1-i a
   c\right)-\Phi'\left(-e^{-\frac{i (b+m)}{c}},0,1-i a c\right)\right. \\ \left.
+e^{\frac{2 i m}{c}}
   \left(\Phi'\left(-e^{\frac{i (-b+m)}{c}},0,1-i a c\right)-e^{\frac{2 i b}{c}}
   \Phi'\left(-e^{\frac{i (b+m)}{c}},0,1-i a c\right)\right)\right)
\end{multline}
\end{example}
\begin{example}
\begin{multline}
\int_0^{\infty } \frac{\cosh (m x) \log \left(a^2-x^2\right) \sinh (b x)}{\sinh (c \pi  x)} \, dx\\
=\frac{\log
   \left(\frac{i}{c}\right) \sin \left(\frac{b}{c}\right)}{c \cos \left(\frac{b}{c}\right)+c \cos
   \left(\frac{m}{c}\right)}+\frac{i e^{-\frac{i (b+m)}{c}}}{2 c} \left(e^{\frac{2 i b}{c}}
   \Phi'\left(-e^{\frac{i (b-m)}{c}},0,1-i a c\right)\right. \\ \left.
-e^{\frac{2 i m}{c}}
   \Phi'\left(-e^{\frac{i (-b+m)}{c}},0,1-i a
   c\right)-\Phi'\left(-e^{-\frac{i (b+m)}{c}},0,1-i a c\right)\right. \\ \left.
+e^{\frac{2 i (b+m)}{c}}
   \Phi'\left(-e^{\frac{i (b+m)}{c}},0,1-i a c\right)\right)
\end{multline}
\end{example}
\begin{example}
\begin{multline}
\int_0^{\infty } \left((a-x)^k-(a+x)^k\right) \left(\coth \left(\frac{\pi  x}{2}\right)-\coth \left(2^n\pi  x\right)\right) \sinh (b x) \sinh (m x) \, dx\\
=\sum _{p=0}^n \frac{1}{2} \left(i 2^{-p}\right)^{1+k} e^{-i 2^{-p} (b+m)} \left(-e^{i 2^{1-p} b} \Phi \left(-e^{i 2^{-p} (b-m)},-k,1-i 2^p a\right)\right. \\ \left.
-e^{i 2^{1-p}m} \Phi \left(-e^{i 2^{-p} (-b+m)},-k,1-i 2^p a\right)+\Phi \left(-e^{-i 2^{-p} (b+m)},-k,1-i 2^pa\right)\right. \\ \left.
+e^{i 2^{1-p} (b+m)} \Phi \left(-e^{i 2^{-p} (b+m)},-k,1-i 2^p a\right)\right)
\end{multline}
\end{example}
\begin{example}
\begin{multline}
\int_0^{\infty } \log \left(16 \pi ^2-x^2\right) \text{sech}(3 x) \sinh
   \left(\frac{x}{2}\right) \sinh (x) \, dx\\
=\frac{(-1)^{3/4} \pi  }{12 \sqrt{3}}\left(\left((22-10 i)-(12-4 i)
   \sqrt{3}\right) \log (2)-(10+2 i) \log (3)+(2-2 i) \sqrt{3} \log (577)\right. \\ \left.
-(2-2 i) \log (37505)+(4+4 i)
   \sqrt{3} \log \left(\frac{3}{\pi }\right)+(6+6 i) \log (\pi )-2
   \Phi'\left(-1,0,\frac{1}{4}-6 i\right)\right. \\ \left.
-2\Phi'\left(-1,0,\frac{1}{4}+6 i\right)+2 \left(i+\sqrt{3}\right)
   \Phi'\left(-1,0,\frac{3}{4}-6 i\right)+2 \left(i+\sqrt{3}\right)
   \Phi'\left(-1,0,\frac{3}{4}+6 i\right)\right. \\ \left.
+\left(-2-2 i \sqrt{3}\right)
   \Phi'\left(-1,0,\frac{5}{4}-6 i\right)+\left(-2-2 i \sqrt{3}\right)
   \Phi'\left(-1,0,\frac{5}{4}+6 i\right)\right. \\ \left.
+2 i\Phi'\left(-1,0,\frac{7}{4}-6 i\right)+2 i\Phi'\left(-1,0,\frac{7}{4}+6 i\right)\right. \\ \left.
+\left(-i+\sqrt{3}\right) \Phi'\left(\frac{1}{2}+\frac{i \sqrt{3}}{2},0,\frac{1}{4}-6 i\right)+\left(-i+\sqrt{3}\right) \Phi'\left(\frac{1}{2}+\frac{i \sqrt{3}}{2},0,\frac{1}{4}+6 i\right)\right. \\ \left.
+2 i \left(i+\sqrt{3}\right) \Phi'\left(\frac{1}{2}+\frac{i \sqrt{3}}{2},0,\frac{3}{4}-6 i\right)+2 i
   \left(i+\sqrt{3}\right) \Phi'\left(\frac{1}{2}+\frac{i\sqrt{3}}{2},0,\frac{3}{4}+6 i\right)\right. \\ \left.
-4 i \Phi'\left(\frac{1}{2}+\frac{i \sqrt{3}}{2},0,\frac{5}{4}-6 i\right)-4 i \Phi'\left(\frac{1}{2}+\frac{i \sqrt{3}}{2},0,\frac{5}{4}+6 i\right)\right. \\ \left.
+\left(1+i \sqrt{3}\right) \left(\Phi'\left(\frac{1}{2}+\frac{i \sqrt{3}}{2},0,\frac{7}{4}-6i\right)+\Phi'\left(\frac{1}{2}+\frac{i \sqrt{3}}{2},0,\frac{7}{4}+6 i\right)\right)\right)
\end{multline}
\end{example}
\begin{example}
\begin{multline}
\int_{0}^{\infty}x^{s-1} e^{x (m-i \alpha )} \left(\Phi \left(-e^{2 x (m-i \alpha )},1-s,\frac{1}{2}\right)-e^{2 i \alpha  x}
   \Phi \left(-e^{2 x (m+i \alpha )},1-s,\frac{1}{2}\right)\right)dx\\
=\pi  2^{-s-1} \Gamma (s) \left((-m+i \alpha
   )^{-s}-(-m-i \alpha )^{-s}\right)
\end{multline}
\end{example}
\begin{example}
\begin{equation}
\int_{0}^{\infty}x^{s-1} e^{m x} \Phi \left(-e^{2 m x},1-s,\frac{1}{2}\right)dx=\frac{1}{2} \pi  \left(-\frac{1}{2 m}\right)^s \Gamma (s)
\end{equation}
\end{example}
\begin{example}
\begin{multline}
\prod_{k=0}^{n-1}\int_{0}^{\infty}e^{m x} x^{\frac{k+s}{n}-1} \Phi \left(-e^{2 m x},1-\frac{k+s}{n},\frac{1}{2}\right)dx\\
=\pi ^{\frac{3 n}{2}-\frac{1}{2}}
   2^{-n-s} n^{\frac{1}{2}-s} \Gamma (s) \left(-\frac{1}{m}\right)^{\frac{1}{2} (n+2 s-1)}
\end{multline}
\end{example}
\begin{example}
\begin{multline}
\prod_{a=1}^{q}\int_{0}^{\infty}e^{-2 i x \alpha } x^{-1+s} \left(\Phi \left(-e^{-2 i x \alpha },1-s,\frac{1}{2} \left(1+\frac{a}{x}\right)\right)-\Phi \left(-e^{-2 i x \alpha },1-s,\frac{1}{2} \left(3+\frac{a}{x}\right)\right)\right. \\ \left.
+e^{4 i x\alpha } \left(-\Phi \left(-e^{2 i x \alpha },1-s,\frac{1}{2} \left(1+\frac{a}{x}\right)\right)+\Phi \left(-e^{2 i x \alpha },1-s,\frac{1}{2} \left(3+\frac{a}{x}\right)\right)\right)\right) \sec (x \alpha
   )dx\\
=\prod_{a=1}^{q}\frac{\pi  }{2^s}\left(\frac{e^{i a \alpha } \Gamma (s,i a \alpha )}{(i \alpha )^s}-\frac{e^{-i a \alpha } \Gamma (s,-i a \alpha )}{(-i \alpha )^s}\right)
\end{multline}
\end{example}
\begin{example}
\begin{multline}
\prod_{a=1}^{1}\int_{0}^{\infty}x^{s-1} e^{-2 i \alpha  x} \sec (\alpha  x) \left(\Phi \left(-e^{-2 i x \alpha },1-s,\frac{1}{2} \left(\frac{a}{x}+1\right)\right)\right. \\ \left.
-\Phi \left(-e^{-2 i x \alpha },1-s,\frac{1}{2}
   \left(\frac{a}{x}+3\right)\right)+e^{4 i \alpha  x} \left(\Phi \left(-e^{2 i x \alpha },1-s,\frac{1}{2} \left(\frac{a}{x}+3\right)\right)\right.\right. \\ \left.\left.
-\Phi \left(-e^{2 i x \alpha },1-s,\frac{1}{2}
   \left(\frac{a}{x}+1\right)\right)\right)\right)dx\\
=\frac{8}{3} i \pi ^2 \left(\pi ^2 \text{Ci}(\pi )-1\right)
\end{multline}
\end{example}
\begin{example}
\begin{multline}
\prod_{\alpha=1}^{2}\int_{0}^{\infty}\frac{e^{-2 i \alpha  x}}{x^7} \left(\Phi \left(-e^{-2 i x \alpha },7,\frac{1}{4} \left(2+\frac{\pi }{x}\right)\right)-\Phi \left(-e^{-2 i x \alpha },7,\frac{1}{4} \left(6+\frac{\pi }{x}\right)\right)\right. \\ \left.
+e^{4 i \alpha x} \left(\Phi \left(-e^{2 i x \alpha },7,\frac{1}{4} \left(6+\frac{\pi }{x}\right)\right)-\Phi \left(-e^{2 i x \alpha },7,\frac{1}{4} \left(2+\frac{\pi }{x}\right)\right)\right)\right) \sec (\alpha x)dx\\
=\frac{2048 \left(\pi ^5 \text{Ci}\left(\frac{\pi }{2}\right)-768+16 \pi ^2-2 \pi ^4\right) \left(-2 \pi ^5 \text{Si}(\pi )+48-4 \pi ^2+2 \pi ^4+\pi ^6\right)}{2025 \pi ^8}
\end{multline}
\end{example}
\begin{example}
\begin{multline}
\int_{0}^{\infty}\int_{0}^{\infty}u^{x-1} v^{y-1} e^{m (u+v)} \Phi \left(-e^{2 m u},1-x,\frac{1}{2}\right) \Phi \left(-e^{2 m
   v},1-y,\frac{1}{2}\right)dudv\\
=\pi ^2 2^{-x-y-2} \Gamma (x) \Gamma (y) \left(-\frac{1}{m}\right)^{x+y}
\end{multline}
\end{example}
\begin{example}
\begin{multline}
\sum_{k=0}^{n}\int_{0}^{\infty}\int_{0}^{\infty}\frac{2^{2 \alpha +n+2} u^{\alpha +k-1} \left(-\frac{1}{m}\right)^{-2 \alpha -n} e^{m (u+v)}
   v^{\alpha -k+n-1} \cos (x (n-2 k)) }{\pi ^2 k! (n-k)!}\\
\Phi \left(-e^{2 m u},-k-\alpha +1,\frac{1}{2}\right) \Phi
   \left(-e^{2 m v},k-n-\alpha +1,\frac{1}{2}\right)dudv\\
=\Gamma (\alpha )^2 C_n^{(\alpha
   )}(\cos (x))
\end{multline}
\end{example}
\begin{example}
\begin{multline}
\int_{0}^{\infty}\left(\text{sech}(x (m-i \alpha ))-\text{sech}(x (m+i \alpha ))\right. \\ \left.
-\text{sech}\left(\frac{1}{2} x (m-2 i \alpha)\right)+\text{sech}\left(\frac{1}{2} x (m+2 i \alpha )\right)\right)dx\\
=-\frac{3 i \pi  \alpha  m^2}{\left(\alpha
   ^2+m^2\right) \left(4 \alpha ^2+m^2\right)}
\end{multline}
\end{example}
\begin{example}
\begin{multline}
\int_0^1\left( \frac{\sqrt{\alpha }}{\sqrt{\alpha -x \alpha +\beta } \log \left(\frac{e^{i a \pi } x \alpha
   }{\alpha -x \alpha +\beta }\right)}+\frac{\sqrt{\beta }}{\sqrt{\alpha +\beta -x \beta } \log \left(\frac{e^{i a
   \pi } (\alpha +\beta -x \beta )}{x \beta }\right)}\right)\frac{dx}{\sqrt{x}} \, \\
=\frac{1}{2} i \left(\psi
   ^{(0)}\left(\frac{1+a}{4}\right)-\psi ^{(0)}\left(\frac{3+a}{4}\right)\right)
\end{multline}
\end{example}
\begin{example}
\begin{multline}
\int_0^1\left( \frac{\sqrt{\alpha } \log \left(-\frac{x \alpha }{\alpha -x \alpha +\beta }\right) \log
   \left(\log \left(-\frac{x \alpha }{\alpha -x \alpha +\beta }\right)\right)}{\sqrt{\alpha -x \alpha +\beta
   }}\right. \\ \left.
+\frac{\sqrt{\beta } \log \left(-\frac{\alpha +\beta -x \beta }{x \beta }\right) \log \left(\log
   \left(-\frac{\alpha +\beta -x \beta }{x \beta }\right)\right)}{\sqrt{\alpha +\beta -x \beta }}\right) \frac{ dx}{\sqrt{x}} \,\\
=i
   \pi ^2 \log \left(\frac{4 i \sqrt[3]{2} e \pi }{A^{12}}\right)
\end{multline}
\end{example}
\begin{example}
\begin{multline}
\int_0^1 \left( \frac{\sqrt{\alpha } \log ^2\left(-\frac{x \alpha }{\alpha -x \alpha +\beta }\right) \log
   \left(\log \left(-\frac{x \alpha }{\alpha -x \alpha +\beta }\right)\right)}{\sqrt{\alpha -x \alpha +\beta
   }}\right. \\ \left.+\frac{\sqrt{\beta } \log ^2\left(-\frac{\alpha +\beta -x \beta }{x \beta }\right) \log \left(\log
   \left(-\frac{\alpha +\beta -x \beta }{x \beta }\right)\right)}{\sqrt{\alpha +\beta -x \beta }}\right)\frac{dx}{\sqrt{x}} \, =14
   \pi  \zeta (3)
\end{multline}
\end{example}
\begin{example}
\begin{multline}
\int_0^1 \left( \frac{\sqrt{\alpha } \log \left(\frac{x \alpha }{\alpha -x \alpha +\beta }\right) \log
   \left(\log \left(\frac{x \alpha }{\alpha -x \alpha +\beta }\right)\right)}{\sqrt{\alpha -x \alpha +\beta
   }}\right. \\ \left.
+\frac{\sqrt{\beta } \log \left(\frac{\alpha +\beta -x \beta }{x \beta }\right) \log \left(\log
   \left(\frac{\alpha +\beta -x \beta }{x \beta }\right)\right)}{\sqrt{\alpha +\beta -x \beta }}\right)\frac{dx}{\sqrt{x}} \, =-4
   i C \pi 
\end{multline}
\end{example}
\begin{example}
\begin{multline}
\int_0^{\infty } \log \left(1-e^{-2 x}\right) \left(\left(1+m^2\right) x
   \cos (m x)+\left(2 m+\left(1+m^2\right) \log (a)\right) \sin (m x)\right)
   \sinh (x) \, dx\\
=\frac{1}{4} \left(\frac{8-8 m^2+8 \left(m+m^3\right) \log
   (a)}{\left(1+m^2\right)^2}-\pi  \text{sech}^2\left(\frac{m \pi }{2}\right)
   (\pi +\log (a) \sinh (m \pi ))\right)
\end{multline}
\end{example}
\begin{example}
\begin{multline}
\int_0^{\infty } \left(x \log (x)-\frac{1}{x}\right) \log \left(1-e^{-2
   x}\right) \sinh (x) \, dx=2 (1-\gamma )-\frac{1}{12} \pi ^2 \log
   \left(\frac{16 \pi ^3 \exp (3)}{A^{36}}\right)
\end{multline}
\end{example}
\begin{example}
\begin{multline}
\int_0^{\infty } \frac{x \left(\pi ^4-2 x^2+x^4+2 \pi ^2
   \left(3+x^2\right)\right) \log \left(1-e^{-2 x}\right) \sinh (x)}{\left(\pi
   ^2+x^2\right)^3} \, dx=\frac{1}{2}+2 \text{Ci}(\pi )-\log (2)
\end{multline}
\end{example}
\begin{example}
\begin{multline}
\int_0^{\infty } \log \left(1-e^{-2 x}\right) \left(\frac{2 x}{\pi
   ^2+x^2}-\left(2 \pi  \tan ^{-1}\left(\frac{x}{\pi }\right)+x \log \left(\pi
   ^2+x^2\right)\right)\right) \sinh (x) \, dx\\
=\frac{1}{6} \left(-\pi ^2 \log
   \left(\frac{16 e^3}{A^{36} \pi ^3}\right)-24 (1+\text{Ci}(\pi )+\log (\pi
   ))\right)
\end{multline}
\end{example}
\begin{example}
\begin{multline}
\int_0^{\infty } \left(\frac{2 a \pi  x}{\left(a^2 \pi
   ^2+x^2\right)^2}+\cot ^{-1}\left(\frac{a \pi }{x}\right)\right) \log
   \left(1-e^{-2 x}\right) \sinh (x) \, dx\\
=i e^{-i a \pi } \left(-\Gamma (0,-i
   a \pi )+e^{2 i a \pi } \Gamma (0,i a \pi )\right)+\pi  \log
   \left(\frac{(-1+a) \Gamma \left(\frac{1}{2} (-1+a)\right)}{\sqrt{2} \sqrt{a}
   \Gamma \left(\frac{a}{2}\right)}\right)
\end{multline}
\end{example}
\begin{example}
\begin{multline}
\int_0^{\infty } e^{x/2} \left(\frac{e^{-x} (-8+a \pi  (4 i-3 a \pi -6 i
   x)+x (-4+3 x))}{(a \pi +i x)^3}\right. \\ \left.
+\frac{-a \pi  (4+3 i a \pi +6 x)+i (-8+x
   (4+3 x))}{(i a \pi +x)^3}\right) \log \left(1-e^{-2 x}\right) \sinh (x) \,
   dx\\
=-\frac{4 i}{a}+8 i e^{\frac{1}{2} (-3) i a \pi } \left(e^{2 i a \pi }
   \Gamma \left(0,\frac{i a \pi }{2}\right)+\Gamma \left(0,\frac{1}{2} (-3) i a
   \pi \right)\right)-8 \Phi (-i,1,1+a)
\end{multline}
\end{example}
\begin{example}
Entry (2.4.14.11) in \cite{prud1}.
\begin{multline}
\int_{-\infty }^{\infty } \frac{e^{m x} x^k \text{cosh}(\alpha )}{\sinh (x)+\sinh (\alpha )} \, dx\\
=i e^{\frac{i k \pi }{2}-m \alpha } \pi  \left((i \alpha )^k-e^{i m \pi } \pi ^k \left(\Phi
   \left(-e^{i m \pi },-k,1+\frac{i \alpha }{\pi }\right)\right.\right. \\ \left.\left.
+e^{2 m \alpha } \left(\Phi \left(-e^{i m \pi },-k,1-\frac{i \alpha }{\pi }\right)+\Phi \left(e^{i m \pi },-k,1-\frac{i \alpha
   }{\pi }\right)\right)-\Phi \left(e^{i m \pi },-k,1+\frac{i \alpha }{\pi }\right)\right)\right) 
\end{multline}
\end{example}
\begin{example}
Entry (2.4.14.12) in \cite{prud1}. 
\begin{multline}
\int_{-\infty }^{\infty } \frac{e^{m x} x^k}{(\sinh (x)+\sinh (\alpha ))^2} \, dx\\
=\frac{e^{\frac{i k \pi }{2}-m \alpha } }{\alpha}\text{sech}^2(\alpha ) \left(e^{i m \pi +2 m \alpha } k
   \pi ^k \alpha  \Phi \left(-e^{i m \pi },1-k,1-\frac{i \alpha }{\pi }\right)\right. \\ \left.
-e^{i m \pi } k \pi ^k \alpha  \Phi \left(-e^{i m \pi },1-k,1+\frac{i \alpha }{\pi }\right)+e^{i m \pi +2 m
   \alpha } k \pi ^k \alpha  \Phi \left(e^{i m \pi },1-k,1-\frac{i \alpha }{\pi }\right)\right. \\ \left.
+e^{i m \pi } k \pi ^k \alpha  \Phi \left(e^{i m \pi },1-k,1+\frac{i \alpha }{\pi }\right)+i e^{i
   m \pi +2 m \alpha } \pi ^{1+k} \alpha  \Phi \left(-e^{i m \pi },-k,1-\frac{i \alpha }{\pi }\right) (m-\tanh (\alpha ))\right. \\ \left.
+i e^{i m \pi +2 m \alpha } \pi ^{1+k} \alpha  \Phi \left(e^{i m
   \pi },-k,1-\frac{i \alpha }{\pi }\right) (m-\tanh (\alpha ))\right. \\ \left.
-i e^{i m \pi } \pi ^{1+k} \alpha  \Phi \left(-e^{i m \pi },-k,1+\frac{i \alpha }{\pi }\right) (m+\tanh (\alpha ))\right. \\ \left.
+i e^{i m\pi } \pi ^{1+k} \alpha  \Phi \left(e^{i m \pi },-k,1+\frac{i \alpha }{\pi }\right) (m+\tanh (\alpha ))-i \pi  (i \alpha )^k (k-m \alpha -\alpha  \tanh (\alpha ))\right)
\end{multline}
\end{example}
\begin{example}
Entry 2.4.14.14 in \cite{prud1}.
\begin{equation}
\int_{-\infty }^{\infty } e^{-i m x} x^k \text{csch}(x) \, dx=2 i (-i)^{-k} \pi ^{1+k} \text{Li}_{-k}\left(-e^{m \pi }\right)
\end{equation}
\end{example}
\begin{example}
Frullani integral given in Eq. (2.4.14.15) in \cite{prud1} and pp. 406-407 in \cite{jeffreys}.
\begin{multline}
\int_{-\infty }^{\infty } \frac{-e^{m x}+e^{p x}}{x (\sinh (x)+\sinh (\alpha ))} \, dx\\
=\frac{1}{\alpha }\left(i e^{-m \alpha } \pi -i e^{-p \alpha } \pi +e^{m (i \pi +\alpha )} \alpha 
   \left(\Phi \left(-e^{i m \pi },1,1-\frac{i \alpha }{\pi }\right)\right.\right. \\ \left.\left.
+\Phi \left(e^{i m \pi },1,1-\frac{i \alpha }{\pi }\right)+e^{-2 m \alpha } \left(\Phi \left(-e^{i m \pi },1,1+\frac{i
   \alpha }{\pi }\right)-\Phi \left(e^{i m \pi },1,1+\frac{i \alpha }{\pi }\right)\right)\right)\right. \\ \left.
-e^{p (i \pi +\alpha )} \alpha  \left(\Phi \left(-e^{i p \pi },1,1-\frac{i \alpha }{\pi
   }\right)+\Phi \left(e^{i p \pi },1,1-\frac{i \alpha }{\pi }\right)+e^{-2 p \alpha } \left(\Phi \left(-e^{i p \pi },1,1+\frac{i \alpha }{\pi }\right)\right.\right.\right. \\ \left.\left.\left.
-\Phi \left(e^{i p \pi
   },1,1+\frac{i \alpha }{\pi }\right)\right)\right)\right) \text{sech}(\alpha )
\end{multline}
\end{example}
\begin{example}
Highly oscillatory integrand.
\begin{equation}
\int_{-\infty }^{\infty } e^{u x} \text{csch}^2(x) \sinh (m x) \, dx=\frac{\pi  (u \sin (m \pi )-m \sin (\pi  u))}{\cos (m \pi )-\cos (\pi  u)}
\end{equation}
\end{example}
\begin{example}
Hilbert transform see section (1.14) in \cite{dlmf}, in terms of the Hypergeometric function.
\begin{equation}
\int_{-\infty }^{\infty } \frac{\text{csch}(x) (x \cos (m x)+a \pi  \sin (m x))}{a^2 \pi ^2+x^2} \, dx=\frac{1}{a}-2 e^{-m \pi } \Phi \left(-e^{-m \pi },1,1+a\right)
\end{equation}
\end{example}
\begin{example}
\begin{multline}
\int_{-\infty }^{\infty } \frac{\pi  \text{csch}(x) \left(-2 a \pi  x \cos (m x)+\left(-a^2 \pi ^2+x^2\right) \sin (m x)\right)}{\left(a^2 \pi ^2+x^2\right)^2} \,
   dx\\
=-\frac{1}{a^2}+2 e^{-m \pi } \Phi \left(-e^{-m \pi },2,1+a\right)
\end{multline}
\end{example}
\begin{example}
\begin{multline}
\int_0^{\infty } \frac{x^{-1+s} \cosh (c x) \sinh (m x)}{3+\cosh (2 c x)} \, dx\\
=\frac{i \left(1+\sqrt{2}\right)^{-\frac{m}{c}} c^{-s} e^{-\frac{i m \pi }{2 c}} \pi ^s}{8 \sqrt{2}} \left(\Phi
   \left(-e^{-\frac{i m \pi }{c}},1-s,\frac{1}{2}-\frac{i \sinh ^{-1}(1)}{\pi }\right)\right. \\ \left.
+\left(1+\sqrt{2}\right)^{\frac{2 m}{c}} \Phi \left(-e^{-\frac{i m \pi
   }{c}},1-s,\frac{1}{2}+\frac{i \sinh ^{-1}(1)}{\pi }\right)\right. \\ \left.
-e^{\frac{i m \pi }{c}} \left(\left(1+\sqrt{2}\right)^{\frac{2 m}{c}} \Phi \left(-e^{\frac{i m \pi
   }{c}},1-s,\frac{1}{2}-\frac{i \sinh ^{-1}(1)}{\pi }\right)\right.\right. \\ \left.\left.
+\Phi \left(-e^{\frac{i m \pi }{c}},1-s,\frac{1}{2}+\frac{i \sinh ^{-1}(1)}{\pi }\right)\right)\right) \sec \left(\frac{\pi
    s}{2}\right)
\end{multline}
\end{example}
\begin{example}
\begin{equation}
\int_0^{\infty } \frac{\cosh (c x) \cosh (m x)}{3+\cosh (2 c x)} \, dx=\frac{\pi  \cosh \left(\frac{m \sinh ^{-1}(1)}{c}\right) \sec \left(\frac{m \pi }{2 c}\right)}{4 \sqrt{2}
   c}
\end{equation}
\end{example}
\begin{example}
\begin{multline}
\int_0^{\infty } \frac{x^{-1+s} \cosh (x)}{3+\cosh (2 x)} \, dx\\
=2^{-\frac{7}{2}+s} \pi ^s \csc \left(\frac{\pi  s}{2}\right) \left(\zeta \left(1-s,\frac{\pi -2 i \sinh
   ^{-1}(1)}{4 \pi }\right)+\zeta \left(1-s,\frac{\pi +2 i \sinh ^{-1}(1)}{4 \pi }\right)\right. \\ \left.
-\zeta \left(1-s,\frac{3}{4}-\frac{i \sinh ^{-1}(1)}{2 \pi }\right)-\zeta
   \left(1-s,\frac{3}{4}+\frac{i \sinh ^{-1}(1)}{2 \pi }\right)\right)
\end{multline}
\end{example}
\begin{example}
\begin{multline}
\int_0^{\infty } \frac{x \cosh (c x) \sinh (m x)}{3+\cosh (2 c x)} \, dx\\
=\frac{\pi  \sec \left(\frac{m \pi }{2 c}\right) \left(2 \sinh ^{-1}(1) \sinh \left(\frac{m \sinh
   ^{-1}(1)}{c}\right)+\pi  \cosh \left(\frac{m \sinh ^{-1}(1)}{c}\right) \tan \left(\frac{m \pi }{2 c}\right)\right)}{8 \sqrt{2} c^2}
\end{multline}
\end{example}
\begin{example}
\begin{multline}
\int_0^{\infty } \frac{x^{-2 m} \left(\log ^k\left(\frac{1}{x}\right)+x^{4 m} \left(-\log \left(\frac{1}{x}\right)\right)^k\right)}{(-1+x)^2} \, dx\\
=2^{1-k} \left((4 i \pi )^k \left(k
   \text{Li}_{1-k}\left(e^{4 i m \pi }\right)+4 i m \pi  \text{Li}_{-k}\left(e^{4 i m \pi }\right)\right)\right)
\end{multline}
\end{example}
\begin{example}
\begin{multline}
\int_0^{\infty } \frac{\log ^2(x) \left(\log \left(\log \left(\frac{1}{x}\right)\right)+x \log (\log (x))\right)}{(-1+x)^2 \sqrt{x}} \, dx\\
=\frac{2}{3} \pi ^2 \left(9+\log
   \left(\frac{16384}{A^{72}}\right)+6 \log (i \pi )+84 i \pi  \zeta '(-2)\right)
\end{multline}
\end{example}
\begin{example}
\begin{multline}
\int_0^{\infty } \frac{\log ^k\left(\frac{1}{x}\right)+x \log ^k(x)}{(-1+x)^2 \sqrt{x}} \, dx=2^{1+k} (i \pi )^k \left(\left(-1+2^k\right) k \zeta (1-k)+i \left(-1+2^{1+k}\right) \pi
    \zeta (-k)\right)
\end{multline}
\end{example}
\begin{example}
\begin{multline}
\int_0^{\infty } \frac{-2 a \log (a)+(a-2 \log (x)) \log (a-2 \log (x))+(a+2 \log (x)) \log (a+2 \log (x))}{(-1+x)^2} \, dx\\
=a (-2-i \pi -2 \log (4 \pi )+2 \log (a))\\
+2 \pi  \left(\pi
   -4 i \log (-i a)+2 i \log \left(8 \pi ^2 a\right)\right)-8 i \pi  \text{log$\Gamma $}\left(-\frac{a i}{4 \pi }\right)
\end{multline}
\end{example}
\begin{example}
\begin{multline}
\int_0^{\infty } \frac{\sin (m x) \left(-1+\tanh \left(\frac{x \alpha }{2}\right)\right)}{x^2+a^2} \, dx\\
=\frac{1}{2 a}\left(-i e^{-a m} \pi +e^{-\frac{m \pi }{\alpha }} \Phi \left(e^{-\frac{2 m \pi
   }{\alpha }},1,\frac{\pi +a \alpha }{2 \pi }\right)\right. \\ \left.
-e^{\frac{m \pi }{\alpha }} \Phi \left(e^{\frac{2 m \pi }{\alpha }},1,\frac{\pi +a \alpha }{2 \pi }\right)+2 \text{Chi}(a m) \sinh (a m)-2
   \cosh (a m) \text{Shi}(a m)\right)
\end{multline}
\end{example}
\begin{example}
\begin{multline}
\int_0^{\infty } x^{-1+s} \sin (m x) \left(-1+\tanh \left(\frac{x \alpha }{2}\right)\right) \, dx\\
=\frac{1}{2} \csc \left(\frac{\pi  s}{2}\right) \left(\left(-1+e^{-i \pi  s}\right) m^{-s}
   \Gamma (s)\right. \\ \left.
+e^{-\frac{m \pi }{\alpha }} (2 \pi )^s \left(\frac{1}{\alpha }\right)^s \left(\Phi \left(e^{-\frac{2 m \pi }{\alpha }},1-s,\frac{1}{2}\right)-e^{\frac{2 m \pi }{\alpha }} \Phi
   \left(e^{\frac{2 m \pi }{\alpha }},1-s,\frac{1}{2}\right)\right)\right)
\end{multline}
\end{example}
\begin{example}
Cosine transform.
\begin{multline}
\int_0^{\infty } \cos (m x) \left(-1+\tanh \left(\frac{x \alpha }{2}\right)\right) \, dx\\
=\frac{i}{m}+\frac{2 e^{-\frac{m \pi }{\alpha }}
   \left(\Phi'\left(e^{-\frac{2 m \pi }{\alpha }},0,\frac{1}{2}\right)+e^{\frac{2 m \pi }{\alpha }} \Phi'\left(e^{\frac{2 m \pi }{\alpha
   }},0,\frac{1}{2}\right)\right)}{\alpha }
\end{multline}
\end{example}
\begin{example}
\begin{multline}
\int_0^{\infty } \left((-i x+\log (a))^k+(i x+\log (a))^k\right) \sin (m x) \left(-1+\tanh \left(\frac{x \alpha }{2}\right)\right) \, dx\\
=-\frac{a^{-m} (-m)^{-k} \Gamma (1+k,-m \log
   (a))+a^m m^{-k} \Gamma (1+k,m \log (a))}{m}\\
-e^{-\frac{m \pi }{\alpha }} (2 \pi )^{1+k} \left(\frac{1}{\alpha }\right)^{1+k} \left(-\Phi \left(e^{-\frac{2 m \pi }{\alpha }},-k,\frac{\pi +\alpha
    \log (a)}{2 \pi }\right)\right. \\ \left.
+e^{\frac{2 m \pi }{\alpha }} \Phi \left(e^{\frac{2 m \pi }{\alpha }},-k,\frac{\pi +\alpha  \log (a)}{2 \pi }\right)\right)
\end{multline}
\end{example}
\begin{example}
\begin{multline}
\int_0^{\infty } \frac{(b-i \log (x)) \log \left(x^t \left(1+\frac{1}{x}\right)_t\right)}{x \left(a^2+b^2-2 i b \log (x)-\log ^2(x)\right)^2} \,
   dx\\
=\sum _{m=1}^t \frac{i \left(\psi ^{(0)}\left(\frac{\pi -i \log \left(e^{-a+i b}\right)+i \log (m)}{2 \pi }\right)-\psi ^{(0)}\left(\frac{\pi -i \log
   \left(e^{a+i b}\right)+i \log (m)}{2 \pi }\right)\right)}{4 a}
\end{multline}
\end{example}
\begin{example}
\begin{equation}
\int_0^{\infty } \frac{\left(4 \pi ^4-3 \log ^4(x)\right) \log (1+x)}{x \left(4 \pi ^4+\log ^4(x)\right)^2} \, dx=\frac{2-\pi  \coth
   \left(\frac{\pi }{2}\right)}{4 \pi ^2}
\end{equation}
\end{example}
\begin{example}
\begin{multline}
\int_0^{\infty } \log \left(x^2+a^2\right) \log (1+\cos (\alpha ) \text{sech}(x)) \, dx\\
=\frac{1}{4} \left(\pi ^2-4 \alpha ^2\right) (-1+\log (2)+\log (\pi ))+4 \pi ^2
   \left(-\zeta'\left(-1,\frac{\pi -\alpha +a}{2 \pi }\right)\right. \\ \left.
-\zeta'\left(-1,\frac{\pi +\alpha +a}{2 \pi
   }\right)+\zeta'\left(-1,\frac{\pi +2 a}{4 \pi }\right)+\zeta'\left(-1,\frac{3}{4}+\frac{a}{2 \pi }\right)\right)
\end{multline}
\end{example}
\begin{example}
\begin{multline}
\int_0^{\infty } x^{-1+s} \log (1+\cos (\alpha ) \text{sech}(x)) \, dx\\
=\frac{2^{-s} \pi ^{1+s} \csc \left(\frac{\pi  s}{2}\right) \left(4^s \zeta \left(-s,\frac{\pi -\alpha }{2 \pi
   }\right)+4^s \zeta \left(-s,\frac{\pi +\alpha }{2 \pi }\right)+\left(-1+2^s\right) \zeta (-s)\right)}{s}
\end{multline}
\end{example}
\begin{example}
\begin{multline}
\int_0^{\infty } x \log (1+\cos (\alpha ) \text{sech}(x)) \, dx\\
=-\frac{1}{16} (3 \zeta (3))+4 \pi ^2 \left(\zeta'\left(-2,\frac{\pi -\alpha }{2 \pi
   }\right)+\zeta'\left(-2,\frac{\pi +\alpha }{2 \pi }\right)\right)
\end{multline}
\end{example}
\begin{example}
Errata for equation  (4.4.4)  in \cite{polyanin}.
\begin{equation}
\int_0^{\infty } \frac{x^{-1+m} \log (x)}{x+\alpha } \, dx=\pi  \alpha ^{-1+m} \csc (m \pi ) (-\pi  \cot (m \pi
   )+\log (\alpha ))
\end{equation}
\end{example}
\begin{example}
\begin{multline}
\int_{0}^{\infty}\frac{i^{-2 m} x^{1-2 m} \left(\pi  i^{2 m+1} \left(x^{4 m}-1\right)+2 \left(1+i^{2 m}\right) \left(-1+i^{2 m}
   x^{4 m}\right) \log (x)\right)}{2 \log (x) (2 \log (x)+i \pi ) (\cos (x)+\cosh
   (x))}dx\\
=\sum_{p=0}^{\infty}\frac{\left(\frac{i}{2}\right)^{-m} \pi ^{2-2 m} (-1)^p (2 p+1)^{1-2 m} \left(-1+\left(\frac{i}{2}\right)^{2
   m} (2 \pi  p+\pi )^{4 m}\right) \text{sech}\left(\pi  \left(p+\frac{1}{2}\right)\right)}{\log \left(\frac{1}{2} i
   (2 \pi  p+\pi )^2\right)}
\end{multline}
\end{example}
\begin{example}
\begin{multline}
\int_{-\infty }^{\infty } \frac{e^{\frac{3 x}{2}}}{\left(-i+e^x\right) \left(1+\sqrt[4]{-1} e^x\right) x}
   \, dx\\
=-\frac{1}{2} \left(-1+\sqrt[4]{-1}\right) \pi +\sqrt{2} \log \left(\cot \left(\frac{3 \pi
   }{16}\right)\right)+\log \left(\cot \left(\frac{\pi }{8}\right) \cot \left(\frac{3 \pi }{16}\right) \tan
   \left(\frac{\pi }{16}\right)\right)\\
+(-1)^{3/4} \log \left(\tan \left(\frac{\pi }{8}\right)\right)
\end{multline}
\end{example}
\begin{example}
\begin{multline}
\int_{-\infty }^{\infty } \frac{e^{\frac{5 x}{4}} \log (-i x+\log (a))}{\left(1+e^x \alpha \right)
   \left(1+e^x \beta \right)} \, dx\\
=\frac{\sqrt{2} \pi }{\sqrt[4]{\alpha } (\alpha -\beta ) \sqrt[4]{\beta }} \left(-\sqrt[4]{\beta } \log (2 \pi )+\sqrt[4]{\alpha }
   \log (8 \pi )\right. \\ \left.
+(1+i) \sqrt[4]{\beta } \left(-i \log \left(\frac{\Gamma \left(\frac{\pi +\log (a)+i \log (\alpha
   )}{8 \pi }\right)}{2 \Gamma \left(\frac{5 \pi +\log (a)+i \log (\alpha )}{8 \pi }\right)}\right)+\log
   \left(\frac{\Gamma \left(\frac{3 \pi +\log (a)+i \log (\alpha )}{8 \pi }\right)}{2 \Gamma \left(\frac{7 \pi
   +\log (a)+i \log (\alpha )}{8 \pi }\right)}\right)\right)\right. \\ \left.
-(1-i) \sqrt[4]{\alpha } \left(\log
   \left(\frac{\Gamma \left(\frac{\pi +\log (a)+i \log (\beta )}{8 \pi }\right)}{\Gamma \left(\frac{5 \pi +\log
   (a)+i \log (\beta )}{8 \pi }\right)}\right)+i \log \left(\frac{\Gamma \left(\frac{3 \pi +\log (a)+i \log
   (\beta )}{8 \pi }\right)}{\Gamma \left(\frac{7 \pi +\log (a)+i \log (\beta )}{8 \pi
   }\right)}\right)\right)\right)
\end{multline}
\end{example}
\begin{example}
\begin{multline}
\int_{-\infty }^{\infty } \frac{e^{\frac{3 x}{2}} \log (-i x+\log (a))}{\left(1+e^x \alpha \right)
   \left(1+e^x \beta \right)} \, dx\\
=\frac{2 \pi }{\sqrt{\alpha } (\alpha -\beta )
   \sqrt{\beta }} \left(\frac{1}{2} \left(\sqrt{\alpha }-\sqrt{\beta }\right) \log
   (4 \pi )+\sqrt{\beta } (\log (-3 \pi +\log (a)+i \log (\alpha ))\right. \\ \left.
-\log (-\pi +\log (a)+i \log (\alpha
   )))+\sqrt{\alpha } (-\log (-3 \pi +\log (a)+i \log (\beta ))\right. \\ \left.
+\log (-\pi +\log (a)+i \log (\beta
   )))+\sqrt{\beta } \left(\text{log$\Gamma $}\left(\frac{-3 \pi +\log (a)+i \log (\alpha )}{4 \pi
   }\right)\right.\right. \\ \left.\left.
-\text{log$\Gamma $}\left(\frac{-\pi +\log (a)+i \log (\alpha )}{4 \pi }\right)\right)+\sqrt{\alpha }
   \left(-\text{log$\Gamma $}\left(\frac{-3 \pi +\log (a)+i \log (\beta )}{4 \pi }\right)\right.\right. \\ \left.\left.
+\text{log$\Gamma
   $}\left(\frac{-\pi +\log (a)+i \log (\beta )}{4 \pi }\right)\right)\right)
\end{multline}
\end{example}
\begin{example}
\begin{multline}
\int_{-\infty }^{\infty } \frac{e^{\frac{3 x}{2}} \log (-i x+\log (a))}{\left(1+e^x \alpha \right)
   \left(1+e^x \beta \right) (-i x+\log (a))} \, dx\\
=\frac{\log (4 \pi )}{2
   \sqrt{\alpha } (\alpha -\beta ) \sqrt{\beta }} \left(\sqrt{\beta } \left(\psi
   ^{(0)}\left(\frac{\pi +\log (a)+i \log (\alpha )}{4 \pi }\right)-\psi ^{(0)}\left(\frac{3 \pi +\log (a)+i \log
   (\alpha )}{4 \pi }\right)\right)\right. \\ \left.
+\sqrt{\alpha } \left(-\psi ^{(0)}\left(\frac{\pi +\log (a)+i \log (\beta )}{4
   \pi }\right)+\psi ^{(0)}\left(\frac{3 \pi +\log (a)+i \log (\beta )}{4 \pi }\right)\right)\right)\\
+\sqrt{\beta
   } \left(-\gamma _1\left(\frac{\pi +\log (a)+i \log (\alpha )}{4 \pi }\right)+\gamma _1\left(\frac{3 \pi +\log
   (a)+i \log (\alpha )}{4 \pi }\right)\right)\\
+\sqrt{\alpha } \left(\gamma _1\left(\frac{\pi +\log (a)+i \log
   (\beta )}{4 \pi }\right)-\gamma _1\left(\frac{3 \pi +\log (a)+i \log (\beta )}{4 \pi }\right)\right)
\end{multline}
\end{example}
\begin{example}
\begin{multline}
\int_{-\infty }^{\infty } \frac{-e^{x (2-i m-\lambda )}+e^{x (2-i r-\lambda )}}{x \left(1+e^x \alpha
   \right) \left(1+e^x \beta \right)} \, dx\\
=\frac{e^{-i \pi  \lambda }}{\alpha  (\alpha -\beta ) \beta } \left(e^{m \pi } \alpha ^{i m+\lambda }
   \beta  \Phi \left(e^{2 \pi  (m-i \lambda )},1,\frac{\pi +i \log (\alpha )}{2 \pi }\right)\right. \\ \left.
-e^{m \pi } \alpha 
   \beta ^{i m+\lambda } \Phi \left(e^{2 \pi  (m-i \lambda )},1,\frac{\pi +i \log (\beta )}{2 \pi }\right)+e^{\pi
    r} \left(-\alpha ^{i r+\lambda } \beta  \Phi \left(e^{2 \pi  (r-i \lambda )},1,\frac{\pi +i \log (\alpha )}{2
   \pi }\right)\right.\right. \\ \left.\left.
+\alpha  \beta ^{i r+\lambda } \Phi \left(e^{2 \pi  (r-i \lambda )},1,\frac{\pi +i \log (\beta
   )}{2 \pi }\right)\right)\right)
\end{multline}
\end{example}
\begin{example}
\begin{multline}
\int_{-\infty }^{\infty } \frac{e^{-x (-2+\lambda )} \cos (m x) (-i x+\log (a))^k}{\left(1+e^x \alpha
   \right) \left(1+e^x \beta \right)} \, dx\\
=-\frac{i 2^k e^{-\pi  (m+i \lambda )} \pi ^{1+k} \alpha ^{-1-i m} \beta
   ^{-1-i m}}{\alpha -\beta }\\ \left(\alpha ^{\lambda } \beta ^{1+i m} \left(e^{2 m \pi } \alpha ^{2 i m} \Phi \left(e^{2 \pi  (m-i
   \lambda )},-k,\frac{\pi +\log (a)+i \log (\alpha )}{2 \pi }\right)\right.\right. \\ \left.\left.
+\Phi \left(e^{-2 \pi  (m+i \lambda
   )},-k,\frac{\pi +\log (a)+i \log (\alpha )}{2 \pi }\right)\right)\right. \\ \left.
-\alpha ^{1+i m} \beta ^{\lambda } \left(e^{2
   m \pi } \beta ^{2 i m} \Phi \left(e^{2 \pi  (m-i \lambda )},-k,\frac{\pi +\log (a)+i \log (\beta )}{2 \pi
   }\right)\right.\right. \\ \left.\left.
+\Phi \left(e^{-2 \pi  (m+i \lambda )},-k,\frac{\pi +\log (a)+i \log (\beta )}{2 \pi
   }\right)\right)\right)
\end{multline}
\end{example}
\begin{example}
\begin{multline}
\int_0^{\frac{1}{2}} \left(\sum _{p=0}^{\infty } \sum _{q=0}^{\infty }
   \frac{\left(-\frac{1}{a}\right)^{p+q} v \Gamma (1+p-v) \Gamma (q+v)
   (1-k)_{-1+p+q}}{\Gamma (1+p) \Gamma (1+q)}\right) \, dv=-\frac{2 a C e^a
   E_{-k}(a)}{k \pi }
\end{multline}
\end{example}
\begin{example}
\begin{multline}
\int_0^{\infty } x^{-1+s} \log \left(1-e^{-x \alpha }\right) (m x \cos (m x)+s \sin (m x)) \,
   dx\\
=\frac{1}{2} \left(e^{\frac{i \pi  s}{2}} m^{-1-s} \alpha  \Gamma (1+s)-\left(\frac{2}{\alpha }\right)^s \pi
   ^{1+s} \left(\text{Li}_{-s}\left(e^{-\frac{2 m \pi }{\alpha }}\right)-\text{Li}_{-s}\left(e^{\frac{2 m \pi
   }{\alpha }}\right)\right) \sec \left(\frac{\pi  s}{2}\right)\right)
\end{multline}
\end{example}
\begin{example}
\begin{multline}
\int_0^{\infty } q^{-\frac{i x \alpha }{2 \pi }} x^{-1+s} \log \left(1-e^{-x \alpha }\right) \left(2 i \pi 
   s+x \alpha  \log (q)+q^{\frac{i x \alpha }{\pi }} (-2 i \pi  s+x \alpha  \log (q))\right) \, dx\\
=e^{\frac{i \pi 
   s}{2}} \pi  \alpha  (\alpha  \log (q))^{-1-s} \left(2^{2+s} \pi ^{1+s} \Gamma (1+s)+i e^{-\frac{1}{2} i \pi 
   (k+s)} (2 \pi )^{-k+s} \left(\frac{1}{\alpha }\right)^s \Gamma (1+k)\right. \\ \left. 
\left(\zeta \left(1+k,-\frac{i \log (q)}{2
   \pi }\right)-\zeta \left(1+k,\frac{i \log (q)}{2 \pi }\right)\right.\right. \\ \left.\left.
+e^{i k \pi } \left(\zeta \left(1+k,1-\frac{i \log
   (q)}{2 \pi }\right)-\zeta \left(1+k,1+\frac{i \log (q)}{2 \pi }\right)\right)\right) \log (q) (\alpha  \log (q))^s
   \sec \left(\frac{\pi  s}{2}\right)\right)
\end{multline}
\end{example}
\begin{example}
\begin{multline}
\int_0^{\infty } q^{-\frac{i x \alpha }{2 \pi }} x^{-1+s} \log \left(1-e^{-x \alpha }\right) \left(2 i \pi 
   s+x \alpha  \log (q)+q^{\frac{i x \alpha }{\pi }} (-2 i \pi  s+x \alpha  \log (q))\right) \, dx\\
=e^{\frac{i \pi 
   s}{2}} (2 \pi )^{2+s} \alpha  \Gamma (1+s) (\alpha  \log (q))^{-1-s}-i 2^{1+s} \pi ^{2+s} \left(\frac{1}{\alpha
   }\right)^s (\pi -i \log (q)) \sec \left(\frac{\pi  s}{2}\right)
\end{multline}
\end{example}
\begin{example}
\begin{multline}
\int_0^{\infty } \frac{\log \left(1-e^{-x \alpha }\right) \left(m \left(a^2+x^2\right) \cos (m x)-2 x \sin
   (m x)\right)}{\left(a^2+x^2\right)^2} \, dx\\
=\frac{\alpha }{4 a} \left(-e^{-a m} \Gamma (0,-a m)+e^{a m} \Gamma (0,a
   m)-e^{-\frac{2 m \pi }{\alpha }} \Phi \left(e^{-\frac{2 m \pi }{\alpha }},1,1+\frac{a \alpha }{2 \pi
   }\right)\right. \\ \left.
+e^{\frac{2 m \pi }{\alpha }} \Phi \left(e^{\frac{2 m \pi }{\alpha }},1,1+\frac{a \alpha }{2 \pi
   }\right)\right)
\end{multline}
\end{example}
\begin{example}
\begin{multline}
\int_0^{\infty } e^{-m x} x^{-1+k} \left((-1)^k (-k+m x)+e^{2 m x} (k+m x)\right) \coth ^{-1}\left(e^{x \alpha
   }\right) \, dx\\
=\frac{1}{2} \left(-\left(\frac{2 \pi  i}{\alpha }\right)^{1+k}\right) \alpha  \left(e^{\frac{i m
   \pi }{\alpha }} \Phi \left(e^{\frac{2 i m \pi }{\alpha }},-k,\frac{1}{2}\right)-\text{Li}_{-k}\left(e^{\frac{2 i m
   \pi }{\alpha }}\right)\right)
\end{multline}
\end{example}
\begin{example}
\begin{multline}
\int_0^{\infty } x^{-1+k} \coth ^{-1}\left(e^{x \alpha }\right) (m x \cos (m x)+k \sin (m x)) \, dx\\
=-2^{-2+k}
   \pi ^{1+k} \sec \left(\frac{k \pi }{2}\right) \left(e^{-\frac{m \pi }{\alpha }} \left(\frac{1}{\alpha }\right)^k
   \left(\Phi \left(e^{-\frac{2 m \pi }{\alpha }},-k,\frac{1}{2}\right)-e^{\frac{2 m \pi }{\alpha }} \Phi
   \left(e^{\frac{2 m \pi }{\alpha }},-k,\frac{1}{2}\right)\right)\right. \\ \left.
+\alpha ^{-k}
   \left(-\text{Li}_{-k}\left(e^{-\frac{2 m \pi }{\alpha }}\right)+\text{Li}_{-k}\left(e^{\frac{2 m \pi }{\alpha
   }}\right)\right)\right)
\end{multline}
\end{example}
\begin{example}
\begin{equation}
\int_0^{\infty } \frac{\coth ^{-1}\left(e^{x \alpha }\right)}{a^2+x^2} \, dx=\frac{\pi  \log
   \left(\frac{\sqrt{a} \sqrt{\alpha } \Gamma \left(\frac{a \alpha }{2 \pi }\right)}{\sqrt{2 \pi } \Gamma
   \left(\frac{\pi +a \alpha }{2 \pi }\right)}\right)}{2 a}
\end{equation}
\end{example}
\begin{example}
\begin{equation}
\int_0^{\infty } \frac{\coth ^{-1}\left(e^{b \sinh (t)}\right)}{\cosh (t)} \, dt=\frac{1}{2} \pi  \log
   \left(\frac{\sqrt{b} \Gamma \left(\frac{b}{2 \pi }\right)}{\sqrt{2 \pi } \Gamma \left(\frac{b+\pi }{2 \pi
   }\right)}\right)
\end{equation}
\end{example}
\begin{figure}[H]
\includegraphics[scale=0.5]{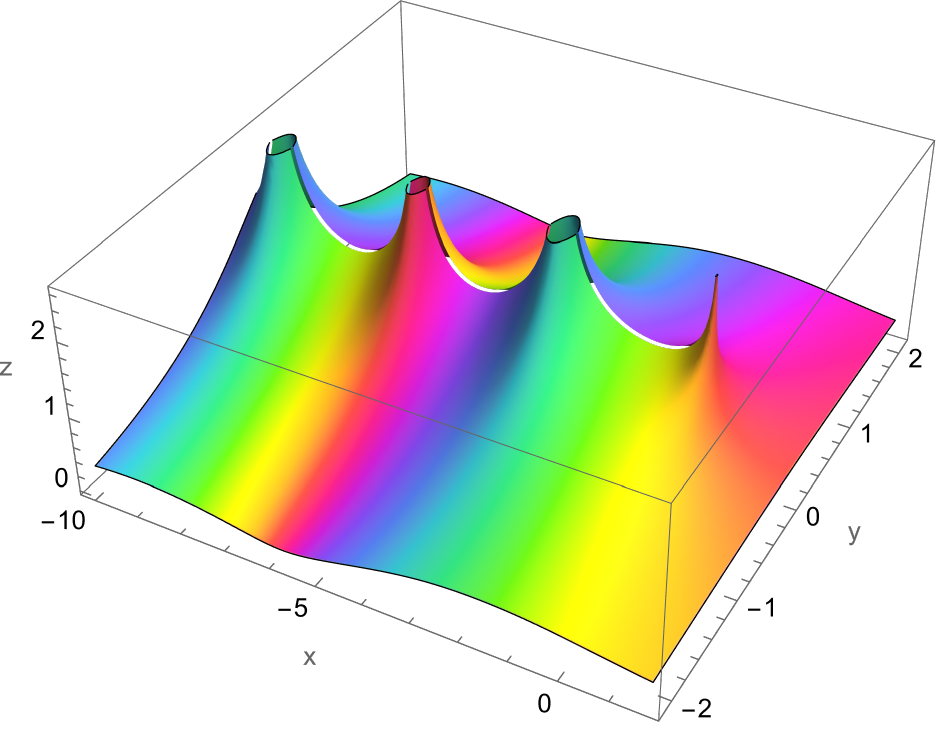}
\caption{Plot of  $f(b)=\log \left(\frac{\sqrt{b} \Gamma \left(\frac{b}{2 \pi }\right)}{\sqrt{2 \pi } \Gamma \left(\frac{b+\pi }{2 \pi
   }\right)}\right)$, $b\in\mathbb{C}$.}
   \label{fig:fig2}
\end{figure}
\vspace{-6pt}
\begin{example}
\begin{multline}
\int_0^{\infty } q^{-x} \coth ^{-1}\left(e^{2 \pi  x \alpha }\right) \left(-\frac{\alpha +(-i b+x \alpha )
   \log (q)}{(b+i x \alpha )^2}+\frac{q^{2 x} (-\alpha +(i b+x \alpha ) \log (q))}{(b-i x \alpha )^2}\right) \,
   dx\\
=\frac{1}{2} \pi  \left(\frac{1}{b}+2 q^{-\frac{i}{2 \alpha }} \Phi \left(q^{-\frac{i}{\alpha
   }},1,\frac{1}{2}-b\right)-2 q^{-\frac{i}{\alpha }} \Phi \left(q^{-\frac{i}{\alpha }},1,1-b\right)\right)
\end{multline}
\end{example}
\begin{example}
\begin{multline}
\int_0^{\infty } q^x \log \left(1-e^{-2 \pi  x \alpha }\right) \left(\frac{-1+(-a+x) \log
   (q)}{(a-x)^2}+\frac{q^{-2 x} (-1-(a+x) \log (q))}{(a+x)^2}\right) \, dx\\
=\frac{\pi  \left(i+2 a q^a \alpha  \Gamma
   (0,a \log (q))-a q^{i/\alpha } \alpha  \left(\Phi \left(q^{\frac{2 i}{\alpha }},1,\frac{1}{2}+\frac{i a \alpha
   }{2}\right)+q^{i/\alpha } \Phi \left(q^{\frac{2 i}{\alpha }},1,1+\frac{i a \alpha
   }{2}\right)\right)\right)}{a}
\end{multline}
\end{example}
\begin{example}
\begin{multline}
\int_0^{\infty } q^{-x} x^{k-1} \left(-k+x \log (q)+(-1)^{-k} q^{2 x} (k+x \log (q))\right) \log \left(1-q^{-4
   x}\right) \, dx\\
=i^{-k} \left(-4 \Gamma (1+k) (-i \log (q))^{-k}+2 \pi  \log ^{-k}(q) \left((2 \pi )^k \left(\zeta
   \left(-k,\frac{1}{4}\right)-\zeta \left(-k,\frac{3}{4}\right)\right)\right.\right. \\ \left.\left.
-i \left(-1+2^{1+k}\right) \pi ^k \zeta
   (-k)\right)\right)
\end{multline}
\end{example}
\begin{example}
\begin{multline}
\int_0^{\infty } q^x \left((-i x+\log (a))^{-1+k}+q^{-2 x} (i x+\log (a))^{-1+k}\right. \\ \left.
+\frac{i \left((-i x+\log
   (a))^k+q^{-2 x} (i x+\log (a))^k\right) \log (q)}{k}\right) \log \left(1-q^{-4 x}\right) \, dx\\
=\frac{1}{k}\left(\pi  \log
   ^k(a)-4 i a^{-i \log (q)} \Gamma (1+k,-i \log (a) \log (q)) (-i \log (q))^{-k}\right. \\ \left.
+(2 \pi )^{1+k} \left(-\zeta \left(-k,\frac{\pi +\log (a) \log (q)}{2 \pi }\right)+i \left(\zeta \left(-k,\frac{\pi +2 \log (a) \log (q)}{4 \pi   }\right)\right.\right.\right. \\ \left.\left.\left.
-\zeta \left(-k,\frac{3}{4}+\frac{\log (a) \log (q)}{2 \pi }\right)\right)+\zeta \left(-k,1+\frac{\log (a)   \log (q)}{2 \pi }\right)\right) \log ^{-k}(q)\right)
\end{multline}
\end{example}
\begin{example}
\begin{multline}
\int_0^{\infty } q^x \left(\frac{1}{(a-i x)^2}+\frac{q^{-2 x}}{(a+i x)^2}-i \left(\frac{1}{a-i x}+\frac{q^{-2
   x}}{a+i x}\right) \log (q)\right) \log \left(1-q^{-4 x}\right) \, dx\\
=-\frac{\pi }{a}+4 q^{-i a} \Gamma (0,-i a
   \log (q)) \log (q)-\log (q) \psi ^{(0)}\left(\frac{\pi +a \log (q)}{2 \pi }\right)\\
+i \log (q) \psi
   ^{(0)}\left(\frac{\pi +2 a \log (q)}{4 \pi }\right)-i \log (q) \psi ^{(0)}\left(\frac{3}{4}+\frac{a \log (q)}{2
   \pi }\right)+\log (q) \psi ^{(0)}\left(1+\frac{a \log (q)}{2 \pi }\right)
\end{multline}
\end{example}
\begin{example}
\begin{multline}
\int_0^{\infty } \frac{e^{-m x} \log \left(1-e^{-4 m x}\right) }{a^2+x^2}\left(a+e^{2 m x} (a+i x)-i x+i m
   \left(a^2+x^2\right)\right. \\ \left.
 \left(e^{2 m x} \log (a-i x)+\log (a+i x)\right)\right) \, dx\\
=\log
   \left(\left(\frac{1}{a}\right)^{4 i} a^{\pi } e^{-4 i e^{-i a m} \Gamma (0,-i a m)} m^{(1-i) \pi } (2 \pi
   )^{(-1+i) \pi } \left(\frac{\Gamma \left(\frac{a m+\pi }{2 \pi }\right) \left(\frac{\Gamma
   \left(\frac{3}{4}+\frac{a m}{2 \pi }\right)}{\Gamma \left(\frac{2 a m+\pi }{4 \pi }\right)}\right)^i}{\Gamma
   \left(1+\frac{a m}{2 \pi }\right)}\right)^{2 \pi }\right)
\end{multline}
\end{example}
\begin{example}
Mellin transform.
\begin{multline}
\int_0^{\infty } x^{-1+s} \text{sech}(c x) \sinh (m x) \tanh (c x) \, dx\\
=\frac{1}{2} i^s c^{-1-s} e^{-\frac{i m
   \pi }{2 c}} \pi ^{-1+s} \left(-i+\cot \left(\frac{\pi  s}{2}\right)\right) \left(m \pi  \Phi \left(-e^{-\frac{i m
   \pi }{c}},1-s,\frac{1}{2}\right)\right. \\ \left.
+i c (-1+s) \Phi \left(-e^{-\frac{i m \pi }{c}},2-s,\frac{1}{2}\right)+e^{\frac{i m
   \pi }{c}} \left(m \pi  \Phi \left(-e^{\frac{i m \pi }{c}},1-s,\frac{1}{2}\right)\right.\right. \\ \left.\left.
-i c (-1+s) \Phi \left(-e^{\frac{i
   m \pi }{c}},2-s,\frac{1}{2}\right)\right)\right)
\end{multline}
\end{example}
\begin{example}
\begin{equation}
\int_0^{\infty } \left(\frac{\tanh (m x)}{x}\right)^2 \, dx=\frac{14 m \zeta (3)}{\pi ^2}
\end{equation}
\end{example}
\begin{example}
\begin{multline}
\int_0^{\infty } \frac{-2 x+2^{1+2 n} x+3 \coth (2 x)-3\times 2^n \coth \left(2^{1+n} x\right)}{x^3} \, dx=\frac{12
   \left(-1+2^{3 n}\right) \zeta (3)}{\pi ^2}
\end{multline}
\end{example}
\begin{example}
\begin{equation}
\int_0^{\infty } \frac{3+x^2-3 x \coth (x)}{3 x^4} \, dx=\frac{\zeta (3)}{\pi ^2}
\end{equation}
\end{example}
\begin{example}
\begin{equation}
\int_0^{\infty } x^2 \text{sech}(2 m x) \sinh (m x) \tanh (2 m x) \, dx=\frac{\pi ^2 (8+3 \pi )}{64 \sqrt{2}
   m^3}
\end{equation}
\end{example}
\begin{example}
\begin{equation}
\int_0^{\infty } (-2+\cosh (2 m x)) \text{sech}^4(m x) \, dx=0
\end{equation}
\end{example}
\begin{example}
\begin{equation}
\int_0^{\infty } \left(\frac{1}{x^2}+\left(-2-\frac{1}{x^2}\right)
   \text{sech}^2(m x)+3 \text{sech}^4(m x)\right) \, dx=\frac{14 m \zeta
   (3)}{\pi ^2}
\end{equation}
\end{example}
\begin{example}
\begin{equation}
\prod _{p=1}^n \int_0^{\infty } \frac{2^p \left(2^p x-\tanh \left(2^p
   x\right)\right)}{x^3} \, dx=2^{\frac{3}{2} n (1+n)} \pi ^{-2 n} (7 \zeta
   (3))^n
\end{equation}
\end{example}
\begin{example}
\begin{equation}
2^{\frac{3}{2} n (1+n)} \prod _{p=1}^n \int_0^{\infty } \frac{2^{-p}
   \left(2^{-p} x-\tanh \left(2^{-p} x\right)\right)}{x^3} \, dx=\left(\frac{7
   \zeta (3)}{\pi ^2}\right)^n
\end{equation}
\end{example}
\begin{example}
\begin{multline}
\int_0^{\infty } \frac{\log (x) \tanh ^2(m x)}{x^2} \, dx=-\frac{m \left((7+16 \log (2)+14 \log (m)-14 \log
   (\pi )) \zeta (3)+14 \zeta '(3)\right)}{\pi ^2}
\end{multline}
\end{example}
\begin{example}
\begin{multline}
\int_0^{\infty } \log (x) \text{sech}(2 m x) \sinh (m x) \tanh (2 m x) \, dx=\frac{4 \coth
   ^{-1}\left(\sqrt{2}\right)+\sqrt{2} \pi  \log \left(\frac{2 \pi  \Gamma \left(\frac{5}{8}\right) \Gamma
   \left(\frac{7}{8}\right)}{m \Gamma \left(\frac{1}{8}\right) \Gamma \left(\frac{3}{8}\right)}\right)}{8 m}
\end{multline}
\end{example}
\begin{example}
\begin{multline}
\int_0^{\infty } \frac{\log ^2(x) \tanh ^2(m x)}{x^2} \, dx\\
=\frac{m }{2 \pi ^2}\left(\left(7 \pi ^2+32 \log (2) (1+\log
   (2))+4 \log \left(\frac{m}{\pi }\right) (7+16 \log (2)+7 \log (m)-7 \log (\pi ))\right)\right. \\ \left.
 \zeta (3)+4 \left(7+16 \log
   (2)+14 \log \left(\frac{m}{\pi }\right)\right) \zeta '(3)+28 \zeta ''(3)\right)
\end{multline}
\end{example}
\begin{example}
\begin{equation}
\int_0^{\infty } x^{-1+s} \tanh ^2(m x) \, dx=-4^{1-s} \left(-4+2^s\right) m^{-s} \Gamma (s) \zeta
   (-1+s)
\end{equation}
\end{example}
\begin{example}
\begin{multline}
\int_0^{\infty } x^{-1+s} \left(\frac{m x}{x}-\frac{\tanh (m x)}{x}\right) \, dx=4^{1-s} \left(-4+2^s\right)
   m^{1-s} \Gamma (-1+s) \zeta (-1+s)
\end{multline}
\end{example}
\begin{example}
\begin{multline}
\int_0^{\infty } x^{-1+s} \left(-2 x+2^{1+2 n} x+3 \coth (2 x)-3\times 2^n \coth \left(2^{1+n} x\right)\right) \,
   dx\\
=-3\times 2^{3-2 (1+s)} \left(2^{1+n}\right)^{-1-s} \left(2^{1+2 n+s}-\left(2^{1+n}\right)^{1+s}\right) \Gamma (s)
   \zeta (s)
\end{multline}
\end{example}
\begin{example}
\begin{multline}
\int_0^{\infty } \log \left(a^2-x^2\right) \text{sech}(c \pi  x) \sin (m x) \tanh (c \pi  x) \,
   dx\\
=-\frac{e^{\frac{m}{2 c}}}{c^2 \pi } \left(-c \Phi \left(-e^{m/c},1,\frac{1}{2}-i a c\right)-\frac{2 m \log
   \left(\frac{i}{c}\right)}{1+e^{m/c}}+e^{-\frac{m}{c}} \left(c \Phi \left(-e^{-\frac{m}{c}},1,\frac{1}{2}-i a
   c\right)\right.\right. \\ \left.\left.
+m \Phi'\left(-e^{-\frac{m}{c}},0,\frac{1}{2}-i a c\right)\right)+m
   \Phi'\left(-e^{m/c},0,\frac{1}{2}-i a c\right)\right)
\end{multline}
\end{example}
\begin{example}
\begin{multline}
\int_0^{\infty } \tanh ^{-1}\left(\frac{x}{a}\right) \cos (m x) \text{sech}(c \pi  x) \tanh (c \pi  x) \,
   dx\\
=-\frac{i e^{-\frac{m}{2 c}} }{2 c^2 \pi }\left(c \Phi \left(-e^{-\frac{m}{c}},1,\frac{1}{2}-i a c\right)+m
   \Phi'\left(-e^{-\frac{m}{c}},0,\frac{1}{2}-i a c\right)\right. \\ \left.
+e^{m/c} \left(c \Phi
   \left(-e^{m/c},1,\frac{1}{2}-i a c\right)-m \Phi'\left(-e^{m/c},0,\frac{1}{2}-i a
   c\right)\right)\right)
\end{multline}
\end{example}
\begin{example}
\begin{multline}
\int_0^{\infty } \frac{\cosh (m x) \text{sech}(c x) \tanh (c x)}{x} \, dx\\
=-\frac{i e^{-\frac{i m \pi }{2 c}}}{2 c \pi }
   \left(-2 e^{\frac{i m \pi }{2 c}} m \pi  \tan ^{-1}\left(e^{-\frac{i m \pi }{2 c}}\right)+i c \Phi
   \left(-e^{-\frac{i m \pi }{c}},2,\frac{1}{2}\right)\right. \\ \left.
+e^{\frac{i m \pi }{c}} \left(2 e^{-\frac{i m \pi }{2 c}} m \pi 
   \tan ^{-1}\left(e^{\frac{i m \pi }{2 c}}\right)+i c \Phi \left(-e^{\frac{i m \pi
   }{c}},2,\frac{1}{2}\right)\right)\right)
\end{multline}
\end{example}
\begin{example}
\begin{multline}
\int_0^{\infty } \frac{\text{sech}(2 m x) \sinh (m x) \tanh (2 m x)}{-a^2+x^2} \, dx\\
=\frac{(-1)^{3/4} }{4 a \pi }\left(\pi
    \Phi \left(-i,1,\frac{1}{2}-\frac{2 i a m}{\pi }\right)+2 i \Phi \left(-i,2,\frac{1}{2}-\frac{2 i a m}{\pi
   }\right)\right. \\ \left.
-i \pi  \Phi \left(i,1,\frac{1}{2}-\frac{2 i a m}{\pi }\right)-2 \Phi \left(i,2,\frac{1}{2}-\frac{2 i a
   m}{\pi }\right)\right)
\end{multline}
\end{example}
\begin{example}
\begin{multline}
\int_0^{\infty } \frac{q^{x/2} \left(-1+q^x\right)^2 \left(1+q^x\right)}{\left(1+q^{2 x}\right)^2
   \left(-a^2+x^2\right)} \, dx\\
=\frac{\sqrt[4]{-1}}{4 a \pi } \left(i \pi  \Phi \left(-i,1,\frac{1}{2}-\frac{i a \log (q)}{\pi
   }\right)-2 \Phi \left(-i,2,\frac{1}{2}-\frac{i a \log (q)}{\pi }\right)\right. \\ \left.
+\pi  \Phi \left(i,1,\frac{1}{2}-\frac{i a
   \log (q)}{\pi }\right)-2 i \Phi \left(i,2,\frac{1}{2}-\frac{i a \log (q)}{\pi }\right)\right)
\end{multline}
\end{example}
\begin{example}
\begin{multline}
\int_0^{\infty } \frac{\left(a^2 \pi ^2+x^2\right) \sec (c x) \sinh (m x) \tan (c x)}{\left(-a^2 \pi
   ^2+x^2\right)^2} \, dx\\
=\frac{i e^{-\frac{m \pi }{2 c}}}{2 \pi ^2} \left(m \pi  \Phi \left(-e^{-\frac{m \pi
   }{c}},2,\frac{1}{2}+a c\right)+2 c \Phi \left(-e^{-\frac{m \pi }{c}},3,\frac{1}{2}+a c\right)\right. \\ \left.
+e^{\frac{m \pi }{c}}
   \left(m \pi  \Phi \left(-e^{\frac{m \pi }{c}},2,\frac{1}{2}+a c\right)-2 c \Phi \left(-e^{\frac{m \pi
   }{c}},3,\frac{1}{2}+a c\right)\right)\right)
\end{multline}
\end{example}
\begin{example}
\begin{equation}
\int_0^{\infty } \frac{\left(a^2 \pi ^2+x^2\right) \tanh ^2(m x)}{\left(-a^2 \pi ^2+x^2\right)^2} \, dx=\frac{m
  }{\pi ^2} \psi ^{(2)}\left(\frac{1}{2}+i a m\right)
\end{equation}
\end{example}
\begin{example}
\begin{equation}
\int_0^{\infty } \frac{\left(-1+q^x\right)^2 \left(a^2 \pi ^2+x^2\right)}{\left(1+q^x\right)^2 \left(-a^2 \pi
   ^2+x^2\right)^2} \, dx=\frac{\log (q) \psi ^{(2)}\left(\frac{1}{2}+\frac{1}{2} i a \log (q)\right)}{2 \pi
   ^2}
\end{equation}
\end{example}
\begin{example}
\begin{equation}
\int_0^{\infty } \frac{\left(4 \pi ^2-x^2\right) \tanh ^2(m x)}{\left(4 \pi ^2+x^2\right)^2} \, dx=-\frac{m
   \psi ^{(2)}\left(\frac{1}{2}-2 m\right)}{\pi ^2}
\end{equation}
\end{example}
\begin{figure}[H]
\includegraphics[scale=0.5]{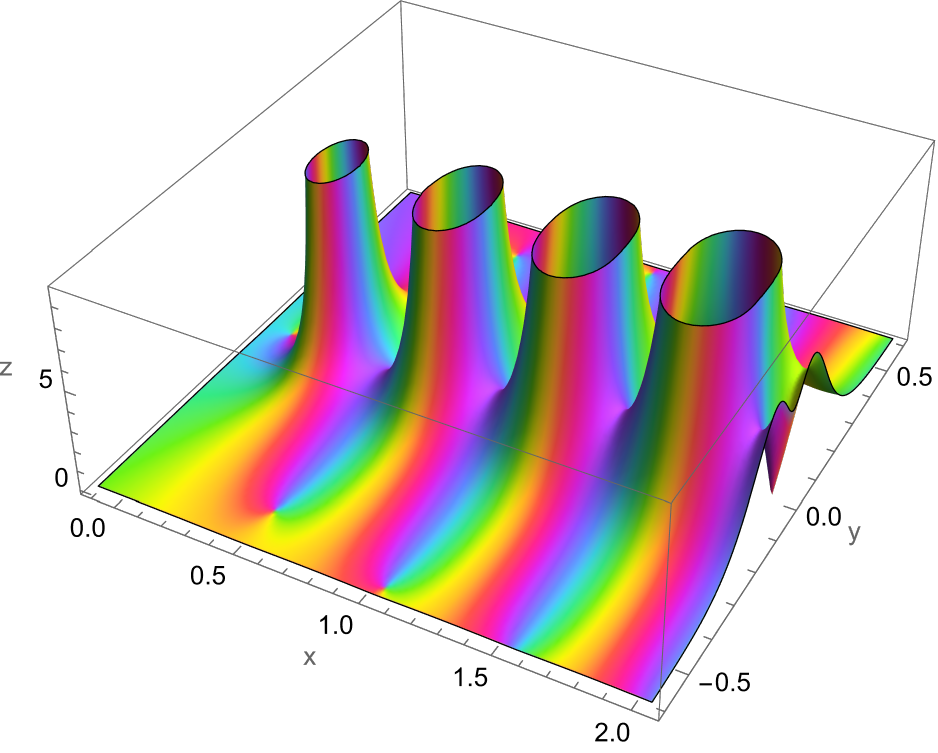}
\caption{Plot of  $f(m)=\frac{m \psi ^{(2)}\left(\frac{1}{2}-2 m\right)}{\pi ^2}$, $m\in\mathbb{C}$.}
   \label{fig:fig2}
\end{figure}
\vspace{-6pt}
\begin{example}
\begin{multline}
\int_0^{\infty } \frac{\left(a^2 \pi ^2+x^2\right) \text{sech}(2 m x) \sinh (m x) \tanh (2 m x)}{\left(-a^2 \pi
   ^2+x^2\right)^2} \, dx\\
=\frac{\sqrt[4]{-1} m}{2 \pi ^2} \left(\pi  \Phi \left(-i,2,\frac{1}{2}+2 i a m\right)+4 i \Phi
   \left(-i,3,\frac{1}{2}+2 i a m\right)-i \pi  \Phi \left(i,2,\frac{1}{2}+2 i a m\right)-\right. \\ \left.
4 \Phi
   \left(i,3,\frac{1}{2}+2 i a m\right)\right)
\end{multline}
\end{example}
\begin{example}
\begin{multline}
\int_0^{\infty } \frac{\left(-\pi ^2+x^2\right) \text{sech}(2 m x) \sinh (m x) \tanh (2 m x)}{\left(\pi
   ^2+x^2\right)^2} \, dx\\
=\frac{\sqrt[4]{-1} m }{2 \pi ^2}\left(\pi  \Phi \left(-i,2,\frac{1}{2}-2 m\right)+4 i \Phi
   \left(-i,3,\frac{1}{2}-2 m\right)\right. \\ \left.
-i \pi  \Phi \left(i,2,\frac{1}{2}-2 m\right)-4 \Phi \left(i,3,\frac{1}{2}-2
   m\right)\right)
\end{multline}
\end{example}
\begin{example}
\begin{multline}
\int_0^{\infty } \frac{\left(-\pi ^2+x^2\right) \text{sech}\left(\frac{x}{2}\right) \sinh
   \left(\frac{x}{4}\right) \tanh \left(\frac{x}{2}\right)}{\left(\pi ^2+x^2\right)^2} \, dx=\frac{-48 C \pi +5 \pi
   ^3-18 \zeta (3)}{192 \sqrt{2} \pi ^2}
\end{multline}
\end{example}
\begin{example}
\begin{equation}
\int_0^{\pi } \cot (x) \log \left(\frac{\sinh (m-i x) \sinh (r+i x)}{\sinh (r-i x) \sinh (m+i
   x)}\right) \, dx=2 \pi  i \log \left(\frac{e^{r-m} \sinh (m)}{\sinh (r)}\right)
\end{equation}
\end{example}
\begin{example}
\begin{multline}
\frac{i}{2\pi}\int_0^{\pi } \frac{  \cot (x) }{e^{2 i x} }\left(\Phi'\left(-e^{-2 i
   x},0,a\right)-e^{4 i x} \Phi'\left(-e^{2 i x},0,a\right)\right) \, dx=\log
   \left(\frac{\Gamma \left(\frac{a}{2}\right)}{\sqrt{2} \Gamma \left(\frac{1+a}{2}\right)}\right)
\end{multline}
\end{example}
\begin{figure}[H]
\includegraphics[scale=0.5]{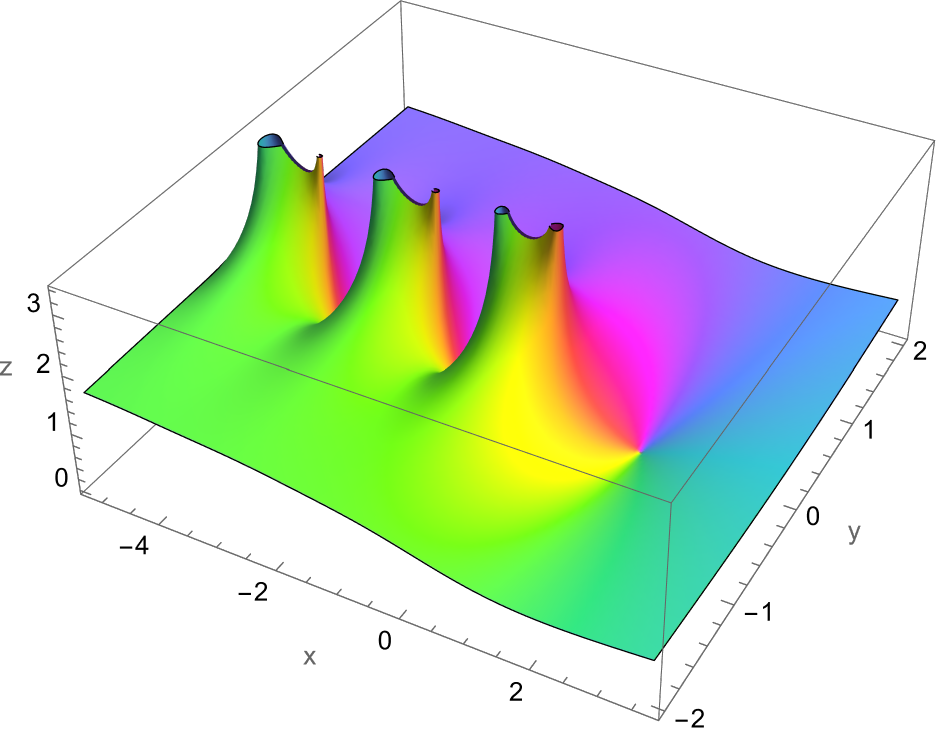}
\caption{Plot of  $f(a)=\log \left(\frac{\Gamma \left(\frac{a}{2}\right)}{\sqrt{2} \Gamma \left(\frac{a+1}{2}\right)}\right)$, $a\in\mathbb{C}$.}
   \label{fig:fig2}
\end{figure}
\vspace{-6pt}
\begin{example}
\begin{multline}
\int_0^{\pi } \log \left(\frac{\left(\frac{\sinh (m-i x)}{\sinh (m+i x)}\right)^{\frac{i \cot (x)}{2 \pi
   }}}{\cosh \left(\frac{x \cot (x)}{\pi }\right)+\sinh \left(\frac{x \cot (x)}{\pi }\right)}\right) \, dx=\log
   \left(\frac{1}{2} (1+\coth (m))\right)
\end{multline}
\end{example}
\begin{figure}[H]
\includegraphics[scale=0.5]{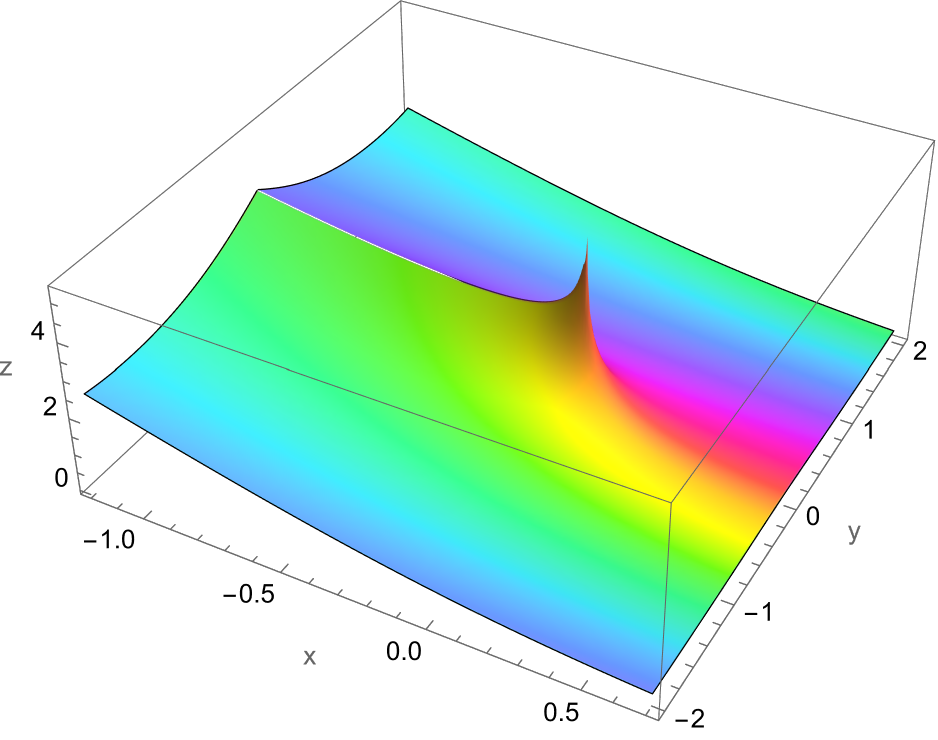}
\caption{Plot of  $f(m)=\log \left(\frac{1}{2} (\coth (m)+1)\right)$, $m\in\mathbb{C}$.}
   \label{fig:fig2}
\end{figure}
\vspace{-6pt}
\begin{example}
\begin{multline}
\int_0^{\pi } \log \left(\frac{\left(\frac{\left(-1+e^{2 i x} z\right) \left(e^{2 i
   x}-z^2\right)}{\left(e^{2 i x}-z\right) \left(-1+e^{2 i x} z^2\right)}\right)^{\frac{i e^{-2 i x} \cot (x)}{4
   \pi }} \left(\frac{\left(e^{2 i x}-z\right) \left(-1+e^{2 i x} z\right)}{z^2+e^{4 i x} z^2-e^{2 i x}
   \left(1+z^4\right)}\right)^{\frac{1}{2 \pi }}}{\left(\frac{\left(e^{2 i x}-z\right) \left(-1+e^{2 i x}
   z^2\right)}{\left(-1+e^{2 i x} z\right) \left(e^{2 i x}-z^2\right)}\right)^{\frac{i e^{2 i x} \cot (x)}{4 \pi }}
   \left(\frac{1+z^4-2 z^2 \cos (2 x)}{1+z^2-2 z \cos (2 x)}\right)^{\frac{\cos (2 x)}{2 \pi }}}\right) \, dx=\log
   \left(\frac{1}{1+z}\right)
\end{multline}
\end{example}
\begin{example}
\begin{equation}
\frac{i}{8\pi}\int_0^{\pi } \frac{  \cot (x) }{e^{2 i x} }\left(\Phi \left(-e^{-2 i x},2,\frac{1}{2}\right)-e^{4 i x} \Phi
   \left(-e^{2 i x},2,\frac{1}{2}\right)\right) \, dx=C
\end{equation}
\end{example}
\begin{example}
\begin{multline}
\int_0^{\pi } \frac{\cot (x)}{e^{2 i x} } \left(e^{4 i x} \Phi'\left(-e^{2 i x},1,\frac{1}{2}\right)-\Phi'\left(-e^{-2 i
   x},1,\frac{1}{2}\right)\right)
   \, dx=i \pi ^2 \left(\gamma +\log \left(\frac{8 \pi ^3}{\Gamma \left(\frac{1}{4}\right)^4}\right)\right)
\end{multline}
\end{example}
\begin{example}
\begin{multline}
\int_0^{\pi } e^{-2 i x} \cot (x) \left(-\Phi'\left(-e^{-2 i x},1,a\right)+e^{4 i
   x} \Phi'\left(-e^{2 i x},1,a\right)\right) \, dx=2 i \pi  \log \left(\frac{\Gamma
   \left(\frac{a}{2}\right)}{\sqrt{2} \Gamma \left(\frac{1+a}{2}\right)}\right)
\end{multline}
\end{example}
\begin{example}
\begin{equation}
\int_0^{\pi } \cot (x) \left(\text{Li}_{-k}\left(\frac{e^{-2 i x}}{q}\right)-\text{Li}_{-k}\left(\frac{e^{2
   i x}}{q}\right)\right) \, dx=-2 i \pi  \text{Li}_{-k}\left(\frac{1}{q}\right)
\end{equation}
\end{example}
\begin{example}
\begin{multline}
\int_0^{\pi } \cot (x) \left(-\zeta \left(1+k,\frac{2 x-i \log (q)}{2 \pi }\right)+\zeta \left(1+k,-\frac{2
   x+i \log (q)}{2 \pi }\right)\right. \\ \left.
+i^{2 k} \left(\zeta \left(1+k,\frac{2 \pi -2 x+i \log (q)}{2 \pi }\right)-\zeta
   \left(1+k,\frac{2 (\pi +x)+i \log (q)}{2 \pi }\right)\right)\right) \, dx\\
=2 i \pi  \left(\zeta
   \left(1+k,-\frac{i \log (q)}{2 \pi }\right)-i^{2 k} \zeta \left(1+k,1+\frac{i \log (q)}{2 \pi
   }\right)\right)
\end{multline}
\end{example}
\begin{example}
\begin{equation}
\frac{1}{2 i \pi }\int_0^{\pi } \cot (x) \log \left(\frac{e^{-2 i x}-q}{e^{2 i x}-q}\right) \, dx=\log
   \left(\frac{q}{q-1}\right)
\end{equation}
\end{example}
\begin{figure}[H]
\includegraphics[scale=0.5]{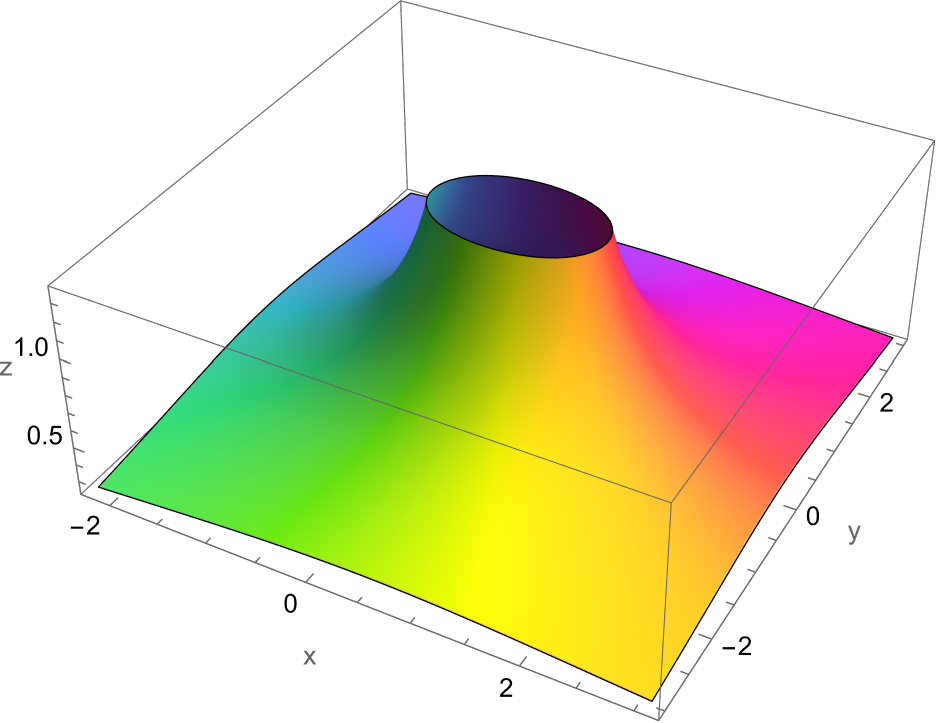}
\caption{Plot of  $f(q)=\log \left(\frac{q}{q-1}\right)$, $q\in\mathbb{C}$.}
   \label{fig:fig2}
\end{figure}
\vspace{-6pt}
\begin{example}
A Mellin transform.
\begin{multline}
\int_0^{\infty } x^{-1+k} \log \left(\frac{1+\cosh (c x)}{\cos (b)+\cosh (c x)}\right) (-k \cos (m x)+m x \sin
   (m x)) \, dx\\
=\frac{1}{2} c^{-k} e^{-\frac{m (b+\pi )}{c}} \pi  \csc \left(\frac{k \pi }{2}\right) \left(-2^{1+k}
   e^{\frac{b m}{c}} \pi ^k \left(\Phi \left(e^{-\frac{2 m \pi }{c}},-k,\frac{1}{2}\right)+e^{\frac{2 m \pi }{c}} \Phi
   \left(e^{\frac{2 m \pi }{c}},-k,\frac{1}{2}\right)\right)\right. \\ \left.
+(2 \pi )^k \left(e^{\frac{2 b m}{c}} \Phi
   \left(e^{-\frac{2 m \pi }{c}},-k,\frac{-b+\pi }{2 \pi }\right)+\Phi \left(e^{-\frac{2 m \pi }{c}},-k,\frac{b+\pi
   }{2 \pi }\right)\right.\right. \\ \left.\left.
+e^{\frac{2 m \pi }{c}} \left(\Phi \left(e^{\frac{2 m \pi }{c}},-k,\frac{-b+\pi }{2 \pi
   }\right)+e^{\frac{2 b m}{c}} \Phi \left(e^{\frac{2 m \pi }{c}},-k,\frac{b+\pi }{2 \pi
   }\right)\right)\right)\right)
\end{multline}
\end{example}
\begin{example}
\begin{multline}
\int_0^{\infty } x^{-1+s} \log \left(\frac{1+\cosh (c x)}{\cos (2 b \pi )+\cosh (c x)}\right) \,
   dx\\
=-\frac{c^{-s} \pi ^{1+s} \csc \left(\frac{\pi  s}{2}\right) \left(-2 \zeta (-s)+2^s \left(\zeta
   \left(-s,\frac{1}{2}-b\right)+\zeta \left(-s,\frac{1}{2}+b\right)+2 \zeta (-s)\right)\right)}{s}
\end{multline}
\end{example}
\begin{example}
\begin{multline}
\int_0^{\infty } x^{-1+k} \log \left(\frac{1+\cosh (2 m x)}{\cos (b)+\cosh (2 m x)}\right) (-k \cosh (m x)-m x
   \sinh (m x)) \, dx\\
=\frac{1}{2} e^{-\frac{1}{2} i (b+\pi )} m^{-k} \pi ^{1+k} \csc \left(\frac{k \pi }{2}\right)
   \left(-2^k \zeta \left(-k,\frac{-b+\pi }{4 \pi }\right)+2^k \zeta \left(-k,\frac{b+\pi }{4 \pi }\right)\right. \\ \left.
+2^k \zeta
   \left(-k,\frac{1}{2} \left(1+\frac{-b+\pi }{2 \pi }\right)\right)+e^{i b} \left(2^k \zeta \left(-k,\frac{-b+\pi }{4
   \pi }\right)-2^k \zeta \left(-k,\frac{1}{2} \left(1+\frac{-b+\pi }{2 \pi }\right)\right)\right)\right. \\ \left.
-2^k \zeta
   \left(-k,\frac{1}{2} \left(1+\frac{b+\pi }{2 \pi }\right)\right)-e^{i b} \left(2^k \zeta \left(-k,\frac{b+\pi }{4
   \pi }\right)-2^k \zeta \left(-k,\frac{1}{2} \left(1+\frac{b+\pi }{2 \pi }\right)\right)\right)\right)
\end{multline}
\end{example}
\begin{example}
Frullani integral see (3.434) in \cite{grad}.
\begin{multline}
\int_0^{\infty } \frac{\log \left(\frac{(\cosh (2 b \pi )+\cosh (c x)) (1+\cosh (d x))}{(1+\cosh (c x)) (\cosh
   (2 b \pi )+\cosh (d x))}\right)}{x} \, dx=2 \log \left(\frac{c}{d}\right) \log (\text{sech}(b \pi ))
\end{multline}
\end{example}
\begin{example}
\begin{multline}
\int_0^{\infty }\log \left(\frac{\left(b^2+2^x b^x\right) \left(1+2^x b^{2+x}\right)
   \left(1+\left(1+b^2\right)^x\right)^2}{\left(1+2^x b^x\right)^2 \left(b^2+\left(1+b^2\right)^x\right) \left(1+b^2
   \left(1+b^2\right)^x\right)}\right) \frac{ \, dx}{2 x}\\
=\log \left(\frac{2 b}{1+b^2}\right) \log \left(\frac{\log (2
   b)}{\log \left(1+b^2\right)}\right)
\end{multline}
\end{example}
\begin{figure}[H]
\includegraphics[scale=0.5]{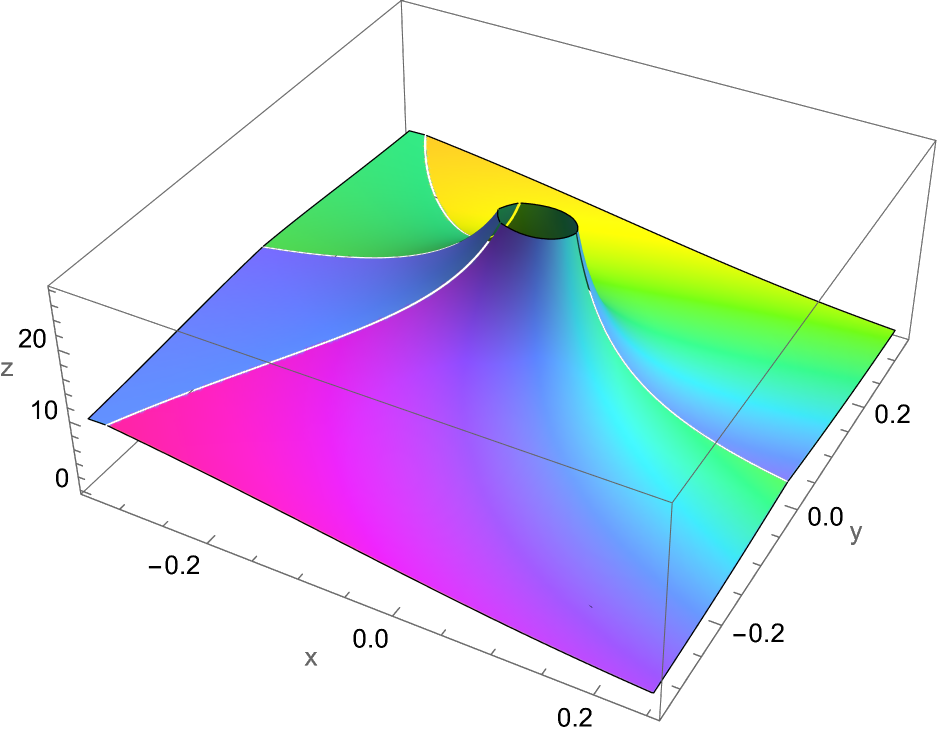}
\caption{Plot of  $f(b)=\log \left(\frac{2 b}{b^2+1}\right) \log \left(\frac{\log (2 b)}{\log \left(b^2+1\right)}\right)$, $b\in\mathbb{C}$.}
   \label{fig:fig2}
\end{figure}
\vspace{-6pt}
\begin{example}
\begin{multline}
\int_0^{\infty }\log \left(\frac{1+\cosh (x)}{\frac{1}{\sqrt{2}}+\cosh (x)}\right) \frac{ \left(\left(\pi
   ^2+x^2\right) \cosh \left(\frac{x}{2}\right) \log \left(\pi ^2+x^2\right)+4 x \sinh
   \left(\frac{x}{2}\right)\right)}{\pi ^2+x^2} \, dx\\
=\pi  \left(\log (256)+2 \sqrt{2+\sqrt{2}} \log \left(\frac{4 \tan
   \left(\frac{\pi }{16}\right)}{\pi }\right)\right)
\end{multline}
\end{example}
\begin{example}
\begin{multline}
\int_0^{\infty } \frac{\log \left(\frac{1+\cosh (x)}{\frac{1}{\sqrt{2}}+\cosh (x)}\right) \left(\left(\pi
   ^2+x^2\right) \cosh \left(\frac{x}{2}\right)-4 x \sinh \left(\frac{x}{2}\right)\right)}{\left(\pi ^2+x^2\right)^2}
   \, dx\\
=\frac{\log (256)+\sqrt{2+\sqrt{2}} \left(-16-\psi ^{(0)}\left(\frac{1}{16}\right)+\psi
   ^{(0)}\left(\frac{7}{16}\right)+\psi ^{(0)}\left(\frac{9}{16}\right)-\psi ^{(0)}\left(\frac{15}{16}\right)\right)}{4
   \pi }
\end{multline}
\end{example}
\begin{example}
\begin{multline}
\int_0^{\infty } \frac{\log (1+\text{sech}(x)) \left(\left(\pi ^2+x^2\right) \cosh \left(\frac{x}{2}\right)-4 x
   \sinh \left(\frac{x}{2}\right)\right)}{\left(\pi ^2+x^2\right)^2} \, dx=\frac{-2 \sqrt{2}+\log \left(12+8
   \sqrt{2}\right)}{\pi }
\end{multline}
\end{example}
\begin{example}
\begin{multline}
\int_0^{\infty } \frac{\log \left(\frac{2 (1+\cosh (x))}{\sqrt{2}+2 \cosh (x)}\right) \left(\left(\pi
   ^2+x^2\right) \cosh \left(\frac{x}{2}\right) \log \left(\pi ^2+x^2\right)+4 x \sinh
   \left(\frac{x}{2}\right)\right)}{\pi ^2+x^2} \, dx\\
=2 \pi  \left(\log (16)+\sqrt{2+\sqrt{2}} \log \left(\frac{4 \tan
   \left(\frac{\pi }{16}\right)}{\pi }\right)\right)
\end{multline}
\end{example}
\begin{example}
\begin{multline}
\int_0^{\infty } \frac{\log \left(1+\frac{1}{1+2 \cosh (x)}\right) \left(\left(4 \pi ^2+x^2\right) \cosh
   \left(\frac{x}{2}\right)-4 x \sinh \left(\frac{x}{2}\right)\right)}{\left(4 \pi ^2+x^2\right)^2} \, dx=\frac{16-9
   \sqrt{3}}{8 \pi }
\end{multline}
\end{example}
\begin{example}
\begin{multline}
\int_0^{\infty } \frac{\log \left(\frac{1+2 \cosh (x)}{2+2 \cosh (x)}\right) \left(\left(4 \pi ^2+x^2\right)
   \cosh \left(\frac{x}{2}\right) \log \left(4 \pi ^2+x^2\right)+4 x \sinh \left(\frac{x}{2}\right)\right)}{4 \pi
   ^2+x^2} \, dx\\
=\log \left(3^{-4 \left(2+\sqrt{3}\right) \pi } 4^{2 \left(2+\sqrt{3}\right) \pi } \pi ^{2
   \left(-2+\sqrt{3}\right) \pi } \left(\frac{\Gamma \left(-\frac{1}{4}\right)}{\Gamma
   \left(-\frac{3}{4}\right)}\right)^{8 \pi } \left(\frac{5 \Gamma \left(-\frac{5}{6}\right) \Gamma
   \left(-\frac{2}{3}\right)}{\Gamma \left(-\frac{1}{3}\right) \Gamma \left(-\frac{1}{6}\right)}\right)^{2 \sqrt{3} \pi
   }\right)
\end{multline}
\end{example}
\begin{example}
\begin{multline}
\int_0^{\infty } \frac{\log \left(\frac{1+\cosh (x)}{\cos (2 b \pi )+\cosh
   (x)}\right) }{a^2 \pi ^2+x^2} \left((a \pi +i x)^k (i a \pi +x) \left((a \pi +i x) \cosh
   \left(\frac{x}{2}\right)+2 i k \sinh \left(\frac{x}{2}\right)\right)\right. \\ \left.
+(a \pi -i
   x)^k (a \pi +i x) \left((i a \pi +x) \cosh \left(\frac{x}{2}\right)+2 k \sinh
   \left(\frac{x}{2}\right)\right)\right)\, dx\\
=-i (4 \pi
   )^{1+k} \left(-2 \zeta \left(-k,\frac{1+a}{4}\right)+2 \zeta
   \left(-k,\frac{3+a}{4}\right)\right. \\ \left.
+\cos (b \pi ) \left(\zeta \left(-k,\frac{1}{4}
   (1+a-2 b)\right)-\zeta \left(-k,\frac{1}{4} (3+a-2 b)\right)\right.\right. \\ \left.\left.
+\zeta
   \left(-k,\frac{1}{4} (1+a+2 b)\right)-\zeta \left(-k,\frac{1}{4} (3+a+2
   b)\right)\right)\right)
\end{multline}
\end{example}
\begin{example}
\begin{multline}
\int_0^{\infty } \frac{\log \left(\frac{1+\cosh (x)}{\cos (2 b \pi )+\cosh
   (x)}\right) \left(\left(a^2 \pi ^2+x^2\right) \cosh \left(\frac{x}{2}\right)-4
   x \sinh \left(\frac{x}{2}\right)\right)}{\left(a^2 \pi ^2+x^2\right)^2} \,
   dx\\
=\frac{1}{2 a \pi }\left(-2 \psi ^{(0)}\left(\frac{1+a}{4}\right)+2 \psi
   ^{(0)}\left(\frac{3+a}{4}\right)+\cos (b \pi ) \left(\psi
   ^{(0)}\left(\frac{1}{4} (1+a-2 b)\right)\right.\right. \\ \left.\left.
-\psi ^{(0)}\left(\frac{1}{4} (3+a-2
   b)\right)+\psi ^{(0)}\left(\frac{1}{4} (1+a+2 b)\right)-\psi
   ^{(0)}\left(\frac{1}{4} (3+a+2 b)\right)\right)\right)
\end{multline}
\end{example}
\begin{example}
\begin{multline}
\int_0^{\infty } \frac{\log \left(\frac{1+\cosh (x)}{\cos (2 b \pi )+\cosh
   (x)}\right) \left(\left(a^2 \pi ^2+x^2\right) \cosh \left(\frac{x}{2}\right)
   \log \left(a^2 \pi ^2+x^2\right)+4 x \sinh
   \left(\frac{x}{2}\right)\right)}{a^2 \pi ^2+x^2} \, dx\\
=8 \pi  \log
   \left(\frac{2 \sqrt{\pi } \Gamma \left(\frac{3+a}{4}\right)}{\Gamma
   \left(\frac{1+a}{4}\right)}\right)+4 \pi  \cos (b \pi ) \log
   \left(\frac{\Gamma \left(\frac{1}{4} (1+a-2 b)\right) \Gamma \left(\frac{1}{4}
   (1+a+2 b)\right)}{4 \pi  \Gamma \left(\frac{1}{4} (3+a-2 b)\right) \Gamma
   \left(\frac{1}{4} (3+a+2 b)\right)}\right)
\end{multline}
\end{example}
\subsubsection{Examples derived from Eq.s (\ref{eq_177110}) and (\ref{eq_177110}) }
\begin{example}
\begin{multline}
\int_0^{\infty } e^{-i m x} \tanh ^{-1}\left(\alpha ^2 \text{sech}^2(x)\right) \left((-i x+\log
   (a))^{-1+k}+e^{2 i m x} (i x+\log (a))^{-1+k}\right. \\ \left.
+\frac{m \left((-i x+\log (a))^k+e^{2 i m x} (i x+\log
   (a))^k\right)}{k}\right) \, dx\\
=-\frac{\pi ^{1+k} }{k}\left(e^{m \left(\cos ^{-1}(\alpha )+2 \sin ^{-1}(\alpha )\right)}
   \Phi \left(e^{m \pi },-k,\frac{\pi -\cos ^{-1}(\alpha )+\log (a)}{\pi }\right)\right. \\ \left.
+e^{m \cos ^{-1}(\alpha )} \Phi
   \left(e^{m \pi },-k,\frac{\cos ^{-1}(\alpha )+\log (a)}{\pi }\right)\right. \\ \left.
-e^{\frac{m \pi }{2}} \left(\alpha
   +\sqrt{1+\alpha ^2}\right)^{-i m} \left(\Phi \left(e^{m \pi },-k,\frac{\pi -2 i \sinh ^{-1}(\alpha )+2 \log (a)}{2
   \pi }\right)\right.\right. \\ \left.\left.
+\left(\alpha +\sqrt{1+\alpha ^2}\right)^{2 i m} \Phi \left(e^{m \pi },-k,\frac{\pi +2 i \sinh
   ^{-1}(\alpha )+2 \log (a)}{2 \pi }\right)\right)\right)
\end{multline}
\end{example}
\begin{example}
\begin{multline}
\int_0^{\infty } e^{-i m x} \left((-i x+\log (a))^{-1+k}+e^{2 i m x} (i x+\log (a))^{-1+k}\right. \\ \left.
+\frac{m \left((-i
   x+\log (a))^k+e^{2 i m x} (i x+\log (a))^k\right)}{k}\right) \log \left(1-\alpha ^4 \text{sech}^4(x)\right) \,
   dx\\
=\frac{2 e^{\frac{m \pi }{2}} \pi ^{1+k} }{k}\left(e^{m \sin ^{-1}(\alpha )} \Phi \left(e^{m \pi },-k,\frac{\pi -\cos
   ^{-1}(\alpha )+\log (a)}{\pi }\right)\right. \\ \left.
+e^{-m \sin ^{-1}(\alpha )} \Phi \left(e^{m \pi },-k,\frac{\cos ^{-1}(\alpha
   )+\log (a)}{\pi }\right)\right. \\ \left.
+\left(\alpha +\sqrt{1+\alpha ^2}\right)^{-i m} \Phi \left(e^{m \pi },-k,\frac{\pi -2 i
   \sinh ^{-1}(\alpha )+2 \log (a)}{2 \pi }\right)\right. \\ \left.
+\left(\alpha +\sqrt{1+\alpha ^2}\right)^{i m} \Phi \left(e^{m \pi
   },-k,\frac{\pi +2 i \sinh ^{-1}(\alpha )+2 \log (a)}{2 \pi }\right)\right. \\ \left.-4 \Phi \left(e^{m \pi
   },-k,\frac{1}{2}+\frac{\log (a)}{\pi }\right)\right)
\end{multline}
\end{example}
\begin{example}
\begin{multline}
\int_0^{\infty } x^{k-1} \tanh ^{-1}\left(\alpha ^2 \text{sech}^2(x)\right) \, dx\\
=-\frac{\pi ^{1+k} \csc
   \left(\frac{k \pi }{2}\right) }{2 k}\left(\zeta \left(-k,\frac{\cos ^{-1}(\alpha )}{\pi }\right)+\zeta
   \left(-k,1-\frac{\cos ^{-1}(\alpha )}{\pi }\right)\right. \\ \left.
-\zeta \left(-k,\frac{1}{2}-\frac{i \sinh ^{-1}(\alpha )}{\pi
   }\right)-\zeta \left(-k,\frac{1}{2}+\frac{i \sinh ^{-1}(\alpha )}{\pi }\right)\right)
\end{multline}
\end{example}
\begin{example}
\begin{multline}
\int_0^{\infty } x^{k-1} \log \left(1-\alpha ^4 \text{sech}^4(x)\right) \, dx\\
=\frac{\pi ^{1+k} \csc
   \left(\frac{k \pi }{2}\right)}{k} \left(\zeta \left(-k,\frac{\cos ^{-1}(\alpha )}{\pi }\right)+\zeta
   \left(-k,1-\frac{\cos ^{-1}(\alpha )}{\pi }\right)\right. \\ \left.
+\zeta \left(-k,\frac{1}{2}-\frac{i \sinh ^{-1}(\alpha )}{\pi
   }\right)+\zeta \left(-k,\frac{1}{2}+\frac{i \sinh ^{-1}(\alpha )}{\pi }\right)-4 \left(-1+2^{-k}\right) \zeta
   (-k)\right)
\end{multline}
\end{example}
\begin{example}
\begin{multline}
\int_0^{\infty } e^{-x} x^{-1+k} \left(-k+x+(-1)^{-k} e^{2 x} (k+x)\right) \tanh ^{-1}\left(\alpha ^2
   \text{sech}^2(x)\right) \, dx\\
=i^{1-k} 2^k \pi ^{1+k} \left(e^{i \cos ^{-1}(\alpha )} \left(\zeta
   \left(-k,\frac{\cos ^{-1}(\alpha )}{2 \pi }\right)-\zeta \left(-k,\frac{\pi +\cos ^{-1}(\alpha )}{2 \pi
   }\right)\right)\right. \\ \left.
+e^{-i \cos ^{-1}(\alpha )} \left(\zeta \left(-k,1-\frac{\cos ^{-1}(\alpha )}{2 \pi }\right)-\zeta
   \left(-k,\frac{\pi +2 \sin ^{-1}(\alpha )}{4 \pi }\right)\right)\right. \\ \left.
-i \left(\alpha +\sqrt{1+\alpha ^2}\right)
   \left(\zeta \left(-k,\frac{\pi -2 i \sinh ^{-1}(\alpha )}{4 \pi }\right)-\zeta \left(-k,\frac{3}{4}-\frac{i \sinh
   ^{-1}(\alpha )}{2 \pi }\right)\right.\right. \\ \left.\left.
+\frac{\zeta \left(-k,\frac{\pi +2 i \sinh ^{-1}(\alpha )}{4 \pi }\right)-\zeta
   \left(-k,\frac{3}{4}+\frac{i \sinh ^{-1}(\alpha )}{2 \pi }\right)}{\left(\alpha +\sqrt{1+\alpha
   ^2}\right)^2}\right)\right)
\end{multline}
\end{example}
\subsubsection{Examples derived from (\ref{eq_27776})}
\begin{example}
\begin{multline}
\int_0^{\infty } x^{-1+k} \cot ^{-1}(\sinh (b x)) \left(-k \cos \left(\frac{k \pi }{2}+m x\right)+k \cos
   \left(\frac{k \pi }{2}+p x\right)\right. \\ \left.
+m x \sin \left(\frac{k \pi }{2}+m x\right)-p x \sin \left(\frac{k \pi }{2}+p
   x\right)\right) \, dx\\
=-\left(\frac{1}{b}\right)^k \pi ^{1+k} \left(e^{\frac{m \pi }{2 b}} \Phi \left(-e^{\frac{m
   \pi }{b}},-k,\frac{1}{2}\right)-e^{\frac{p \pi }{2 b}} \Phi \left(-e^{\frac{p \pi
   }{b}},-k,\frac{1}{2}\right)\right)
\end{multline}
\end{example}
\begin{example}
\begin{multline}
\int_0^{\infty } \frac{\left(\cot ^{-1}(\sinh (b x))-\cot ^{-1}(\sinh (c x))\right) (m x \cos (m x)-p x \cos (p
   x)-\sin (m x)+\sin (p x))}{x^2} \, dx\\
=2 \left(b \tan ^{-1}\left(e^{\frac{m \pi }{2 b}}\right)-c \tan
   ^{-1}\left(e^{\frac{m \pi }{2 c}}\right)-b \tan ^{-1}\left(e^{\frac{p \pi }{2 b}}\right)+c \tan
   ^{-1}\left(e^{\frac{p \pi }{2 c}}\right)\right)
\end{multline}
\end{example}
\begin{example}
\begin{multline}
\int_0^{\infty } x^{-1+k} \cot ^{-1}(\sinh (b x)) \, dx=\frac{\pi   }{k}\left(\frac{2 \pi }{b}\right)^k\left(\zeta
   \left(-k,\frac{1}{4}\right)-\zeta \left(-k,\frac{3}{4}\right)\right) \sec \left(\frac{k \pi }{2}\right)
\end{multline}
\end{example}
\begin{example}
\begin{equation}
\int_0^{\infty } \frac{-\cot ^{-1}(\sinh (b x))+\cot ^{-1}(\sinh (c x))}{x} \, dx=\frac{1}{2} \pi  \log
   \left(\frac{b}{c}\right)
\end{equation}
\end{example}
\begin{example}
\begin{multline}
\int_0^{\infty } \cot ^{-1}(\sinh (b x)) \left(\frac{\log (a \pi -i x)}{a \pi -i x}-\frac{\log (a \pi +i x)}{a
   \pi +i x}\right) \, dx\\
=\frac{1}{8} i \pi  \left((2-4 a b) \log ^2\left(\frac{1}{b}\right)+4 (1-2 a b) \log
   \left(\frac{1}{b}\right) \log (2 \pi )-4 \log ^2(a \pi )\right. \\ \left.
+(2+4 a b) \log ^2\left(\frac{2 \pi }{b}\right)+\log
   \left(4 \pi ^2\right) \left(\log (2 \pi )-a b \log \left(4 \pi ^2\right)\right)\right. \\ \left.
+16 \log \left(\frac{2 \pi
   }{b}\right) \log \left(\frac{\Gamma \left(\frac{3}{4}+\frac{a b}{2}\right)}{\Gamma \left(\frac{1}{4}+\frac{a
   b}{2}\right)}\right)+8 \zeta''\left(0,\frac{1}{4}+\frac{a b}{2}\right)-8
   \zeta''\left(0,\frac{3}{4}+\frac{a b}{2}\right)\right)
\end{multline}
\end{example}
\begin{example}
\begin{multline}
\int_0^{\infty } e^{\frac{b x}{2}} \cot ^{-1}(\sinh (b x)) \left(\frac{1}{-i x+\log (a)}-\frac{e^{-b x}}{i
   x+\log (a)}\right. \\ \left.
+\frac{1}{2} i b \left(\log (-i x+\log (a))-e^{-b x} \log (i x+\log (a))\right)\right) \, dx\\
=i \pi 
   \left(\sqrt{2} \log \left(\frac{\pi }{b}\right)-2 \sqrt[4]{-1} \log \left(\frac{\Gamma \left(\frac{\pi +2 b \log
   (a)}{8 \pi }\right)}{2 \Gamma \left(\frac{5}{8}+\frac{b \log (a)}{4 \pi }\right)}\right)\right. \\ \left.
+2 (-1)^{3/4} \log
   \left(\frac{\Gamma \left(\frac{3}{8}+\frac{b \log (a)}{4 \pi }\right)}{2 \Gamma \left(\frac{7}{8}+\frac{b \log
   (a)}{4 \pi }\right)}\right)-\log (\log (a))\right)
\end{multline}
\end{example}
\begin{example}
\begin{multline}
\int_0^{\infty } e^{-i m x} \left(k \left(\left(\frac{\pi }{2 b}-i x\right)^{-1+k}-e^{2 i m x} \left(\frac{\pi
   }{2 b}+i x\right)^{-1+k}\right)\right. \\ \left.
+m \left(\left(\frac{\pi }{2 b}-i x\right)^k-e^{2 i m x} \left(\frac{\pi }{2 b}+i
   x\right)^k\right)\right) \cot ^{-1}(\sinh (b x)) \, dx\\
=-i \left(\frac{1}{b}\right)^k \pi  \left(\left(\frac{\pi
   }{2}\right)^k+\frac{2^{-k} e^{-\frac{m \pi }{2 b}} \Gamma (1+k) \left(i^{1+k} \zeta \left(1+k,\frac{b-i m}{2
   b}\right)+i^{-1-k} \zeta \left(1+k,\frac{b+i m}{2 b}\right)\right)}{\pi }\right)
\end{multline}
\end{example}
\begin{example}
\begin{multline}
\int_0^{\infty } e^{-i m x} \left(\frac{e^{2 i m x} (b-m \pi -i b m x)}{(\pi +i b x)^2}+\frac{m \pi +b (-1-i m
   x)}{(\pi -i b x)^2}\right) \cot ^{-1}\left(\sinh \left(\frac{b x}{2}\right)\right) \, dx\\
=-i+\frac{i e^{-\frac{m \pi
   }{b}} \left(m \pi +b \log \left(2 \cosh \left(\frac{m \pi }{b}\right)\right)\right)}{b}
\end{multline}
\end{example}
\begin{example}
\begin{equation}
\int_0^{\infty } \cot ^{-1}(\sinh (b x)) \log \left(\frac{\pi -2 i b x}{\pi +2 i b x}\right) \,
   dx=\frac{i \pi ^2}{2 b}\log \left(\frac{4 \sqrt[3]{2} e}{A^{12}}\right)
\end{equation}
\end{example}
\begin{example}
\begin{multline}
\int_0^{\infty } \tan ^{-1}\left(\frac{2 q^x}{-1+q^{2 x}}\right) \log \left(\frac{\pi -2 i x \log (q)}{\pi +2 i
   x \log (q)}\right) \, dx=\frac{\left(i \pi ^2\right) \log \left(\frac{4 \sqrt[3]{2} e}{A^{12}}\right)}{2 \log
   (q)}
\end{multline}
\end{example}
\begin{example}
\begin{multline}
\int_0^{\infty } \cot ^{-1}(\sinh (b x)) \left(-2 i b x+(\pi -2 i b x) \log \left(\frac{\pi }{2 b}-i
   x\right)-(\pi +2 i b x) \log \left(\frac{\pi }{2 b}+i x\right)\right) \, dx\\
=-\frac{i \pi  \left(\pi ^2 \log
   \left(\frac{\pi }{2 b}\right)+14 \zeta (3)\right)}{4 b}
\end{multline}
\end{example}
\begin{example}
\begin{equation}
\int_0^{\infty } \cot ^{-1}(\sinh (b x)) \, dx=\frac{2 C}{b}
\end{equation}
\end{example}
\begin{example}
\begin{equation}
\frac{\log (q)}{2} \int_0^{\infty } \tan ^{-1}\left(\frac{2 q^x}{q^{2 x}-1}\right) \, dx=C
\end{equation}
\end{example}
\subsubsection{Example derived using equation (\ref{eq_1968})}
\begin{example}
\begin{multline}
\int_0^{\infty } x \cot ^{-1}\left(e^{2 x}\right) \sinh (x) \,
   dx\\
=\frac{1}{\sqrt{2}}\left(\log \left(\frac{15}{7}\right)-\text{log$\Gamma
   $}\left(-\frac{7}{8}\right)+\text{log$\Gamma
   $}\left(-\frac{5}{8}\right)+\text{log$\Gamma
   $}\left(-\frac{3}{8}\right)-\text{log$\Gamma $}\left(-\frac{1}{8}\right)\right. \\ \left.
+2\pi \left(\zeta'\left(-1,\frac{1}{8}\right)+\zeta'\left(-1,\frac{3}{8}\right)-\zeta'\left(-1,\frac{5}{8}\right)-\zeta'\left(-1,\frac{7}{8}\right)\right)\right)
   \end{multline}
\end{example}
\begin{example}
\begin{equation}
\int_0^{\infty } \cot ^{-1}\left(e^{2 x}\right) \sinh (x) \,
   dx=\frac{1}{4} \left(-\pi +\sqrt{2} \pi \right)
\end{equation}
\end{example}
\begin{example}
\begin{multline}
\int_0^{\infty } \frac{\left(\pi ^6-x^4 \left(-6+x^2\right)+\pi ^4 \left(6+x^2\right)-\pi ^2 x^2
   \left(36+x^2\right)\right) \cot ^{-1}\left(e^{2 x}\right) \sinh (x)}{\left(\pi ^2+x^2\right)^4} \,
   dx\\
=\frac{13312+648 \pi ^3-6 \pi  \left(-320+9 \sqrt{2}+9 \psi ^{(1)}\left(\frac{1}{8}\right)+9 \psi
   ^{(1)}\left(\frac{3}{8}\right)\right)+27 \psi ^{(2)}\left(\frac{1}{8}\right)-27 \psi
   ^{(2)}\left(\frac{3}{8}\right)}{216 \sqrt{2} \pi ^2}
\end{multline}
\end{example}
\begin{example}
\begin{multline}
\int_0^{\infty }15 \left(\frac{1}{(\pi -i x)^6}+\frac{1}{(\pi +i x)^6}+\frac{\pi ^4-6 \pi ^2 x^2+x^4}{10
   \left(\pi ^2+x^2\right)^4}\right) \cot ^{-1}\left(e^{2 x}\right) \sinh (x) \frac{\, dx}{2 \pi } \\
=\frac{1}{20736 \sqrt{2} \pi ^5}\left(31719424+124416 \pi
   ^5-6 \pi  \left(-335872+648 \sqrt{2}\right.\right. \\ \left.\left.
+27 \psi ^{(3)}\left(\frac{1}{8}\right)+27 \psi
   ^{(3)}\left(\frac{3}{8}\right)\right)+81 \psi ^{(4)}\left(\frac{1}{8}\right)-81 \psi
   ^{(4)}\left(\frac{3}{8}\right)\right)
\end{multline}
\end{example}
\begin{example}
\begin{multline}
\int_0^{\infty } \cot ^{-1}\left(e^{2 x}\right) \left(-2+\frac{2 \left(\pi ^2-x^2\right)}{\left(\pi
   ^2+x^2\right)^2}-\log \left(\pi ^2+x^2\right)\right) \sinh (x) \, dx\\
=\frac{4 \sqrt{2}}{3}+\frac{\pi }{2}+\log
   \left(\pi ^{\pi /2} \tan ^2\left(\frac{\pi }{8}\right) \left(\frac{99 \Gamma \left(-\frac{11}{8}\right) \Gamma
   \left(-\frac{9}{8}\right)}{70 e \pi  \Gamma \left(-\frac{7}{8}\right) \Gamma
   \left(-\frac{5}{8}\right)}\right)^{\frac{\pi }{\sqrt{2}}}\right)
\end{multline}
\end{example}
\begin{example}
\begin{multline}
\int_0^{\infty } \frac{\left(\pi ^4+256 x^2 \left(-6+x^2\right)+32 \pi ^2 \left(1+x^2\right)\right) \cot
   ^{-1}\left(e^{2 x}\right) \sinh (x)}{\left(\pi ^2+16 x^2\right)^3} \, dx\\
=\frac{48 C+\pi  \left(5 \pi +6 \left(-4
   \sqrt{2}+\log (4)\right)\right)}{96 \sqrt{2} \pi ^2}
\end{multline}
\end{example}
\begin{example}
\begin{multline}
\int_0^{\infty } \cot ^{-1}\left(e^{2 x}\right) \left(\frac{32 \pi }{\pi ^2+16 x^2}+4 i x \log \left(\frac{\pi
   +4 i x}{\pi -4 i x}\right)+\pi  \left(-2+\log \left(\frac{\pi ^2}{16}+x^2\right)\right)\right) \sinh (x) \, dx\\
=\log
   \left(\left(\left(2^{1+\frac{2 \sqrt{2}}{3}} \sqrt{e} A^{-6 \sqrt{2}} \pi
   ^{-\frac{1}{2}+\frac{1}{\sqrt{2}}}\right)^{\pi } \exp \left(-2 \sqrt{2} C\right) \left(\frac{9 \Gamma
   \left(-\frac{3}{4}\right)^2}{\pi  \Gamma \left(-\frac{1}{4}\right)^2}\right)^{\sqrt{2}}\right)^{\pi
   }\right)
\end{multline}
\end{example}
\begin{example}
\begin{multline}
\int_0^{\infty } \left(\frac{96}{(\pi -4 i x)^4}+\frac{1}{(\pi -4 i x)^2}+\frac{96}{(\pi +4 i
   x)^4}+\frac{1}{(\pi +4 i x)^2}\right) \cot ^{-1}\left(e^{2 x}\right) \sinh (x) \, dx\\
=\frac{-16 \pi +\frac{48 C \pi
   +7 \pi ^3-18 \zeta (3)}{6 \sqrt{2}}}{32 \pi ^2}
\end{multline}
\end{example}
\begin{example}
\begin{multline}
\int_0^{\infty } \frac{\left(\pi ^4+256 x^2 \left(-6+x^2\right)+32 \pi ^2 \left(1+x^2\right)\right) \cot
   ^{-1}\left(e^{2 x}\right) \sinh (x)}{\left(\pi ^2+16 x^2\right)^3} \, dx\\
=\frac{-4 \pi +\frac{48 C+\pi  (5 \pi +\log
   (4096))}{6 \sqrt{2}}}{16 \pi ^2}
\end{multline}
\end{example}
\begin{example}
\begin{multline}
\int_0^{\infty } \frac{4 \pi  \left(\pi ^4+16 x^2 \left(-6+x^2\right)+8 \pi ^2 \left(1+x^2\right)\right) \cot
   ^{-1}\left(e^{2 x}\right) \sinh (x)}{\left(\pi ^2+4 x^2\right)^3} \, dx\\
=-\frac{1}{4 \sqrt{2} \pi }\left(\zeta \left(2,-\frac{1}{8}\right)+2
   \pi  \left(2 \sqrt{2}-\psi ^{(0)}\left(-\frac{1}{8}\right)-\psi ^{(0)}\left(\frac{1}{8}\right)+\psi
   ^{(0)}\left(\frac{3}{8}\right)\right.\right. \\ \left.\left.
+\psi ^{(0)}\left(\frac{5}{8}\right)\right)-\psi ^{(1)}\left(\frac{1}{8}\right)-\psi
   ^{(1)}\left(\frac{3}{8}\right)+\psi ^{(1)}\left(\frac{5}{8}\right)\right)
\end{multline}
\end{example}
\begin{example}
\begin{multline}
\int_0^{\infty } \frac{\left(288 \pi ^2+81 \pi ^4-1536 x^2+288 \pi ^2 x^2+256 x^4\right) \cot ^{-1}\left(e^{2
   x}\right) \sinh (x)}{\left(9 \pi ^2+16 x^2\right)^3} \, dx\\
=\frac{-48 \sqrt{2}+48 \sqrt{2} C-16 \pi +24 \sqrt{2} \pi
   -5 \sqrt{2} \pi ^2+4 \sqrt{2} \pi  \log (8)}{576 \pi ^2}
\end{multline}
\end{example}
\begin{example}
\begin{multline}
\int_0^{\infty } x^{-2+k} \cot ^{-1}\left(e^{2 x}\right) \left(-2 k m x \cos (m
   x)+\left(k-k^2+\left(1+m^2\right) x^2\right) \sin (m x)\right) \sinh (x) \, dx\\
=\frac{e^{-\frac{1}{4} (m \pi )} \pi
   ^k \csc \left(\frac{k \pi }{2}\right)}{4 \sqrt{2}} \left(-k \Phi \left(-e^{-m \pi },1-k,\frac{1}{4}\right)+(1+m) \pi  \Phi\left(-e^{-m \pi },-k,\frac{1}{4}\right)\right. \\ \left.
+e^{-\frac{1}{2} (m \pi )} \left(k \Phi \left(-e^{-m \pi
   },1-k,\frac{3}{4}\right)-(-1+m) \pi  \Phi \left(-e^{-m \pi },-k,\frac{3}{4}\right)\right)\right. \\ \left.
+e^{\frac{m \pi }{2}}
   \left(k \Phi \left(-e^{m \pi },1-k,\frac{1}{4}\right)+(-1+m) \pi  \Phi \left(-e^{m \pi
   },-k,\frac{1}{4}\right)\right.\right. \\ \left.
\left.
-e^{\frac{m \pi }{2}} \left(k \Phi \left(-e^{m \pi },1-k,\frac{3}{4}\right)+(1+m) \pi  \Phi
   \left(-e^{m \pi },-k,\frac{3}{4}\right)\right)\right)\right)
\end{multline}
\end{example}
\begin{example}
\begin{multline}
\int_0^{\infty } \frac{\cot ^{-1}\left(e^{2 x}\right) \left(x^2 \left(-6+x^2\right)+2 \left(1+x^2\right) \log
   ^2(a)+\log ^4(a)\right) \sinh (x)}{\left(x^2+\log ^2(a)\right)^3} \, dx\\
=\frac{1}{16 \pi  \log ^2(a)}\left(-4 \pi ^2+\sqrt{2} \log (a) \left(2
   \pi  \left(-\psi ^{(0)}\left(\frac{\pi +4 \log (a)}{8 \pi }\right)-\psi ^{(0)}\left(\frac{3}{8}+\frac{\log (a)}{2
   \pi }\right)\right.\right.\right. \\ \left.\left.\left.
+\psi ^{(0)}\left(\frac{5}{8}+\frac{\log (a)}{2 \pi }\right)+\psi ^{(0)}\left(\frac{7}{8}+\frac{\log
   (a)}{2 \pi }\right)\right)+\psi ^{(1)}\left(\frac{\pi +4 \log (a)}{8 \pi }\right)\right.\right. \\ \left.\left.
-\psi
   ^{(1)}\left(\frac{3}{8}+\frac{\log (a)}{2 \pi }\right)-\psi ^{(1)}\left(\frac{5}{8}+\frac{\log (a)}{2 \pi
   }\right)+\psi ^{(1)}\left(\frac{7}{8}+\frac{\log (a)}{2 \pi }\right)\right)\right)
\end{multline}
\end{example}
\begin{example}
\begin{multline}
\int_0^{\infty } e^{-i m x} \cot ^{-1}\left(e^{2 x}\right) \left(\frac{e^{2 i m x} (i-i a m \pi +m x)}{(-i a
   \pi +x)^3}+\frac{-i+i a m \pi +m x}{(i a \pi +x)^3}\right. \\ \left.
+\frac{e^{i m x} \left(1+m^2\right) (a \pi  \cos (m x)+x \sin (m
   x))}{a^2 \pi ^2+x^2}\right) \sinh (x) \, dx\\
=\frac{1}{4} \left(-\frac{1}{a}+\frac{\sqrt{2} e^{\frac{m \pi }{4}}}{\pi }
   \left((\pi -m \pi ) \Phi \left(-e^{m \pi },1,\frac{1}{4}+a\right)+\Phi \left(-e^{m \pi
   },2,\frac{1}{4}+a\right)\right.\right. \\ \left.
\left.
+e^{\frac{m \pi }{2}} \left((1+m) \pi  \Phi \left(-e^{m \pi },1,\frac{3}{4}+a\right)-\Phi
   \left(-e^{m \pi },2,\frac{3}{4}+a\right)\right)\right)\right)
\end{multline}
\end{example}
\begin{example}
\begin{multline}
\int_0^{\infty } e^x \left(\frac{e^{-2 x} (1-i \pi +x)}{(\pi +i x)^3}+\frac{i (-1+i \pi +x)}{(i \pi
   +x)^3}\right) \cot ^{-1}\left(e^{2 x}\right) \sinh (x) \, dx\\
=-\frac{19}{12}+\frac{4 i ((-8+10 i)+9 C)}{9 \pi }+\pi
\end{multline}
\end{example}
\begin{example}
\begin{multline}
\int_0^{\infty } e^x \left(\frac{8-4 i \pi -8 x}{(\pi -2 i x)^3}-\frac{4 i e^{-2 x} (\pi +2 i (1+x))}{(\pi +2 i
   x)^3}\right) \cot ^{-1}\left(e^{2 x}\right) \sinh (x) \, dx\\
=\frac{3}{2}-\frac{4 i ((-1-i)+C)}{\pi }
\end{multline}
\end{example}
\begin{example}
\begin{multline}
\int_0^{\infty } e^x \left(\frac{\left(\frac{\pi }{4}-i x\right)^{-2+k} \left(-1+k+\frac{i \pi }{2}+2
   x\right)}{-1+k}\right. \\ \left.
+e^{-2 x} \left(\frac{\pi }{4}+i x\right)^{-2+k} \left(1-\frac{2 \left(-\frac{1}{4} (i \pi
   )+x\right)}{-1+k}\right)\right) \cot ^{-1}\left(e^{2 x}\right) \sinh (x) \, dx\\
=-\frac{2^{-1-2 k} \pi ^k \left(\pi
   +2^{1+k} \left((1+i)-i 2^k\right) k \zeta (1-k)+2^{1+k} \left(-1+2^{1+k}\right) \pi  \zeta (-k)\right)}{(-1+k)
   k}
\end{multline}
\end{example}
\begin{example}
\begin{multline}
\int_0^{\infty } e^x \left(\frac{16 (4-i \pi -4 x)}{(\pi -4 i x)^3}-\frac{16 i e^{-2 x} (\pi +4 i (1+x))}{(\pi
   +4 i x)^3}\right) \cot ^{-1}\left(e^{2 x}\right) \sinh (x) \, dx\\
=-1+\left(\frac{1}{6}+\frac{i}{12}\right) \pi +\log
   (2)
\end{multline}
\end{example}
\begin{example}
\begin{multline}
\int_0^{\infty } e^x \cot ^{-1}\left(e^{2 x}\right) \left(-i (\pi -4 i x)^{-1+k}+(-1+k) (\pi -4 i x)^{-2+k}
   (-2+2 k+i \pi +4 x)\right. \\ \left.
\log \left(\frac{1}{4} (\pi -4 i x)\right)+e^{-2 x} (\pi +4 i x)^{-2+k} \left(-i \pi +4
   x+(-1+k) (-2+2 k+i \pi -4 x)\right.\right. \\ \left.\left.
 \log \left(\frac{1}{4} (\pi +4 i x)\right)\right)\right) \sinh (x) \, dx\\
=\frac{1}{16 k^2}\left(\pi
   ^{1+k} \left(-1+k \left(2+k \log \left(\frac{4}{\pi }\right)+\log \left(\frac{\pi }{4}\right)\right)\right)\right. \\ \left.
+2 \pi 
   \left((2 \pi )^k \left(1+(-1+k) k \log \left(\frac{\pi }{2}\right)\right)\right.\right. \\ \left.\left.
-2^{1+k} \pi ^k \left(k+2^k (1-2 k+(-1+k)
   k \log (\pi ))\right)\right) \zeta (-k)\right. \\ \left.
+4 k \cos \left(\frac{k \pi }{2}\right) \Gamma (1+k) \left((1+i) \left(1+k
   \log \left(\frac{2}{\pi }\right)+\log \left(\frac{\pi }{2}\right)\right)\right.\right. \\ \left.\left.
+i 2^k (-1+(-1+k) \log (\pi ))\right) \zeta
   (k)+2^{1+k} (-1+k) k \pi ^k \left(\left((1+i)-i 2^k\right) k \zeta '(1-k)\right.\right. \\ \left.\left.
+\left(-1+2^{1+k}\right) \pi  \zeta
   '(-k)\right)\right)
\end{multline}
\end{example}
\subsubsection{Examples derived from Eq. (429410)}
\begin{example}
\begin{multline}
\int_0^{\infty } \frac{x^{-1+m} \tan ^{-1}(b x)}{a i+\log (x)} \, dx\\
=\sum _{y=0}^{\infty } \left(e^{i
   a}\right)^{-m} \pi  \left(\Gamma \left(0,-i m \left(a+\pi +2 \pi  y-i \log
   \left(-\frac{i}{b}\right)\right)\right)\right. \\ \left.
-\Gamma \left(0,-i m \left(a+\pi +2 \pi  y-i \log
   \left(\frac{i}{b}\right)\right)\right)\right)
\end{multline}
\end{example}
\begin{example}
\begin{multline}
\int_0^{\infty } \frac{x^{-1+m} \tan ^{-1}(b x) \log (\log (x))}{\sqrt{\log (x)}} \, dx\\
=\sum _{y=0}^{\infty }
   \left(4 \pi  \, _2F_2\left(\frac{1}{2},\frac{1}{2};\frac{3}{2},\frac{3}{2};m \left(i (\pi +2 \pi  y)+\log
   \left(-\frac{i}{b}\right)\right)\right) \sqrt{i (\pi +2 \pi  y)+\log \left(-\frac{i}{b}\right)}\right. \\ \left.
-4 \pi  \,
   _2F_2\left(\frac{1}{2},\frac{1}{2};\frac{3}{2},\frac{3}{2};m \left(i (\pi +2 \pi  y)+\log
   \left(\frac{i}{b}\right)\right)\right) \sqrt{i \pi  (1+2 y)+\log \left(\frac{i}{b}\right)}\right. \\ \left.
+\frac{\pi ^{3/2}}{\sqrt{m}}
   \left(-\text{erfi}\left(\sqrt{m} \sqrt{i (\pi +2 \pi  y)+\log \left(-\frac{i}{b}\right)}\right) \log \left(i (\pi
   +2 \pi  y)+\log \left(-\frac{i}{b}\right)\right)\right.\right. \\ \left.\left.
+\text{erfi}\left(\sqrt{m} \sqrt{i (\pi +2 \pi  y)+\log
   \left(\frac{i}{b}\right)}\right) \log \left(i (\pi +2 \pi  y)+\log
   \left(\frac{i}{b}\right)\right)\right)\right)
\end{multline}
\end{example}
\begin{example}
\begin{multline}
\int_0^{\infty } \frac{\tan ^{-1}(x)}{x^{3/2} \sqrt[3]{\log (x)}} \, dx=\sum _{y=0}^{\infty } 2^{2/3} \pi 
   \left(-\Gamma \left(\frac{2}{3},\frac{1}{4} i \pi  (3+4 y)\right)+\Gamma \left(\frac{2}{3},\frac{1}{4} i (\pi +4
   \pi  y)\right)\right)
\end{multline}
\end{example}
\begin{example}
\begin{multline}
\int_0^{\infty } \frac{\tan ^{-1}(x)}{x^{5/3} \sqrt[4]{\log (x)}} \, dx=\sum _{y=0}^{\infty }
   \left(\frac{3}{2}\right)^{3/4} \pi  \left(-\Gamma \left(\frac{3}{4},\frac{1}{3} i \pi  (3+4 y)\right)+\Gamma
   \left(\frac{3}{4},\frac{1}{3} i (\pi +4 \pi  y)\right)\right)
\end{multline}
\end{example}
\subsubsection{Equations in this section are derived using (\ref{eq_2466})}
\begin{example}
\begin{multline}
\int_0^{\infty } \frac{\tan ^{-1}\left(\frac{x}{a}\right) \sinh (m x)}{b+\cosh (c x)} \,
   dx\\
=\frac{e^{-\frac{m \left(i \pi +\cosh ^{-1}(b)\right)}{c}} \pi }{2 \sqrt{-1+b^2} c}
   \left(\Phi'\left(e^{-\frac{2 i m \pi }{c}},0,\frac{a c+\pi -i \cosh ^{-1}(b)}{2
   \pi }\right)\right. \\ \left.
-\left(b+\sqrt{-1+b} \sqrt{1+b}\right)^{\frac{2 m}{c}}
   \left(\Phi'\left(e^{-\frac{2 i m \pi }{c}},0,\frac{a c+\pi +i \cosh ^{-1}(b)}{2
   \pi }\right)\right.\right. \\ \left.\left.
+e^{\frac{2 i m \pi }{c}} \Phi'\left(e^{\frac{2 i m \pi
   }{c}},0,\frac{a c+\pi -i \cosh ^{-1}(b)}{2 \pi }\right)\right)+e^{\frac{2 i m \pi }{c}}
   \Phi'\left(e^{\frac{2 i m \pi }{c}},0,\frac{a c+\pi +i \cosh ^{-1}(b)}{2 \pi
   }\right)\right)
\end{multline}
\end{example}
\begin{example}
\begin{multline}
\int_0^{\infty } \frac{\cosh (m x) \log \left(-x^2+\log ^2(a)\right)}{b+\cosh (c x)} \, dx\\
=\frac{i \pi 
   \left(-1-i \cot \left(\frac{m \pi }{c}\right)\right)}{2 \sqrt{-1+b^2} c} \left(\cos \left(\frac{m \left(\pi -i \cosh
   ^{-1}(b)\right)}{c}\right)-i \sin \left(\frac{m \left(\pi -i \cosh ^{-1}(b)\right)}{c}\right)\right)\\
 \left(2
   e^{\frac{2 i m \pi }{c}} \log \left(\frac{2 i \pi }{c}\right) \left(-1+\cosh \left(\frac{2 m \cosh
   ^{-1}(b)}{c}\right)+\sinh \left(\frac{2 m \cosh ^{-1}(b)}{c}\right)\right)\right. \\ \left.
-\left(-1+e^{\frac{2 i m \pi
   }{c}}\right) \left(\cosh \left(\frac{2 m \cosh ^{-1}(b)}{c}\right)+\sinh \left(\frac{2 m \cosh
   ^{-1}(b)}{c}\right)\right)\right. \\ \left.
 \Phi'\left(e^{-\frac{2 i m \pi }{c}},0,\frac{\pi +i \cosh
   ^{-1}(b)-i c \log (a)}{2 \pi }\right)\right. \\ \left.
+\left(-1+e^{\frac{2 i m \pi }{c}}\right)
   \left(\Phi'\left(e^{-\frac{2 i m \pi }{c}},0,-\frac{i \left(i \pi +\cosh ^{-1}(b)+c
   \log (a)\right)}{2 \pi }\right)\right.\right. \\ \left.\left.
+e^{\frac{2 i m \pi }{c}} \left(-\Phi'\left(e^{\frac{2
   i m \pi }{c}},0,\frac{\pi +i \cosh ^{-1}(b)-i c \log (a)}{2 \pi }\right)\right.\right.\right. \\ \left.\left.\left.+\left(\cosh \left(\frac{2 m \cosh
   ^{-1}(b)}{c}\right)+\sinh \left(\frac{2 m \cosh ^{-1}(b)}{c}\right)\right)\right.\right.\right. \\ \left.\left.\left.
   \Phi'\left(e^{\frac{2 i m \pi }{c}},0,-\frac{i \left(i \pi +\cosh ^{-1}(b)+c \log
   (a)\right)}{2 \pi }\right)\right)\right)\right)
\end{multline}
\end{example}
\begin{example}
\begin{multline}
\int_0^{\infty } \frac{x^{-1+s} \sinh (m x)}{b+\cosh (c x)} \, dx\\
=\frac{2^{-2+s} \left(b+\sqrt{-1+b}
   \sqrt{1+b}\right)^{-\frac{m}{c}} c^{-s} e^{-\frac{i m \pi }{c}} \pi ^s }{\sqrt{-1+b^2}}\left(\Phi \left(e^{-\frac{2 i m \pi
   }{c}},1-s,\frac{\pi -i \cosh ^{-1}(b)}{2 \pi }\right)\right. \\ \left.
 -\left(b+\sqrt{-1+b} \sqrt{1+b}\right)^{\frac{2 m}{c}}
   \left(\Phi \left(e^{-\frac{2 i m \pi }{c}},1-s,\frac{\pi +i \cosh ^{-1}(b)}{2 \pi }\right)\right.\right. \\ \left.\left.
+e^{\frac{2 i m \pi
   }{c}} \Phi \left(e^{\frac{2 i m \pi }{c}},1-s,\frac{\pi -i \cosh ^{-1}(b)}{2 \pi }\right)\right)+e^{\frac{2 i m
   \pi }{c}} \Phi \left(e^{\frac{2 i m \pi }{c}},1-s,\frac{\pi +i \cosh ^{-1}(b)}{2 \pi }\right)\right) \sec
   \left(\frac{\pi  s}{2}\right)
\end{multline}
\end{example}
\begin{example}
\begin{multline}
\int_0^{\infty } \frac{x^{-1+s} \cosh (m x)}{b+\cosh (c x)} \, dx\\
=\frac{i 2^{-2+s} \left(b+\sqrt{-1+b}
   \sqrt{1+b}\right)^{-\frac{m}{c}} c^{-s} e^{-\frac{i m \pi }{c}} \pi ^s \csc \left(\frac{\pi  s}{2}\right)}{\sqrt{-1+b^2}}\\
   \left(-\Phi \left(e^{-\frac{2 i m \pi }{c}},1-s,\frac{\pi -i \cosh ^{-1}(b)}{2 \pi }\right)+\left(b+\sqrt{-1+b}
   \sqrt{1+b}\right)^{\frac{2 m}{c}}\right. \\ \left.
 \left(\Phi \left(e^{-\frac{2 i m \pi }{c}},1-s,\frac{\pi +i \cosh ^{-1}(b)}{2
   \pi }\right)-e^{\frac{2 i m \pi }{c}} \Phi \left(e^{\frac{2 i m \pi }{c}},1-s,\frac{\pi -i \cosh ^{-1}(b)}{2 \pi
   }\right)\right)\right. \\ \left.
+e^{\frac{2 i m \pi }{c}} \Phi \left(e^{\frac{2 i m \pi }{c}},1-s,\frac{\pi +i \cosh ^{-1}(b)}{2
   \pi }\right)\right)
\end{multline}
\end{example}
\begin{example}
\begin{multline}
\int_0^{\infty } \frac{e^{m x} x^{-1+s}}{\cos (b \pi )+\cosh (c x)} \, dx\\
=2^{-1+s} c^{-s} e^{-\frac{i \pi 
   (2 (1+b) m+c s)}{2 c}} \pi ^s \csc (b \pi ) \csc (\pi  s) \left(e^{\frac{i \pi  (2 b m+c s)}{c}} \Phi
   \left(e^{-\frac{2 i m \pi }{c}},1-s,\frac{1}{2}-\frac{b}{2}\right)\right. \\ \left.
-e^{i \pi  s} \Phi \left(e^{-\frac{2 i m \pi
   }{c}},1-s,\frac{1+b}{2}\right)+e^{\frac{2 i m \pi }{c}} \Phi \left(e^{\frac{2 i m \pi
   }{c}},1-s,\frac{1}{2}-\frac{b}{2}\right)\right. \\ \left.
-e^{\frac{2 i (1+b) m \pi }{c}} \Phi \left(e^{\frac{2 i m \pi
   }{c}},1-s,\frac{1+b}{2}\right)\right)
\end{multline}
\end{example}
\begin{example}
\begin{multline}
\int_0^{\infty } \frac{\sinh (r x)-\sinh (m x)}{x (b+\cosh (c x))} \, dx\\
=\frac{1}{2 \sqrt{-1+b^2}}\left(e^{\frac{m \left(i \pi +\cosh
   ^{-1}(b)\right)}{c}} \Phi \left(e^{\frac{2 i m \pi }{c}},1,\frac{\pi -i \cosh ^{-1}(b)}{2 \pi }\right)\right. \\ \left.
-e^{\frac{i m
   \left(\pi +i \cosh ^{-1}(b)\right)}{c}} \Phi \left(e^{\frac{2 i m \pi }{c}},1,\frac{\pi +i \cosh ^{-1}(b)}{2 \pi
   }\right)-e^{\frac{r \left(i \pi +\cosh ^{-1}(b)\right)}{c}} \Phi \left(e^{\frac{2 i \pi  r}{c}},1,\frac{\pi -i \cosh
   ^{-1}(b)}{2 \pi }\right)\right. \\ \left.
+e^{\frac{i r \left(\pi +i \cosh ^{-1}(b)\right)}{c}} \Phi \left(e^{\frac{2 i \pi 
   r}{c}},1,\frac{\pi +i \cosh ^{-1}(b)}{2 \pi }\right)\right)
\end{multline}
\end{example}
\begin{example}
Eq. (2.4.7.8) in \cite{prud1}.
\begin{multline}
\int_0^{\infty } \frac{\cosh (b x)}{\left(x^2+z^2\right) (\cos (t)+\cosh (c x))} \, dx\\
=\frac{e^{-\frac{i b (\pi
   +t)}{c}} \csc (t) }{4 z}\left(e^{\frac{2 i b t}{c}} \Phi \left(e^{-\frac{2 i b \pi }{c}},1,\frac{\pi -t+c z}{2 \pi }\right)-\Phi
   \left(e^{-\frac{2 i b \pi }{c}},1,\frac{\pi +t+c z}{2 \pi }\right)\right. \\ \left.
+e^{\frac{2 i b \pi }{c}} \left(\Phi \left(e^{\frac{2 i b
   \pi }{c}},1,\frac{\pi -t+c z}{2 \pi }\right)-e^{\frac{2 i b t}{c}} \Phi \left(e^{\frac{2 i b \pi }{c}},1,\frac{\pi +t+c
   z}{2 \pi }\right)\right)\right)
\end{multline}
\end{example}
\subsubsection{Derived from equation (\ref{eq_22711})}
\begin{example}
\begin{multline}
\int_0^{\infty } \frac{\log ^k(a b)-\log ^k(a x)}{(b-x) \sqrt{x}} \, dx\\
=\frac{i 4^k e^{\frac{i k \pi }{2}} \pi ^{1+k}
   \left(-2 \zeta \left(-k,\frac{2 \pi -i \log (a)-i \log (b)}{4 \pi }\right)+\zeta \left(-k,\frac{4 \pi -i \log (a)-i \log
   (b)}{4 \pi }\right)+\zeta \left(-k,-\frac{i (\log (a)+\log (b))}{4 \pi }\right)\right)}{\sqrt{b}}
\end{multline}
\end{example}
\begin{example}
\begin{multline}
\int_0^{\infty } \frac{\log (\log (a x))}{\sqrt{x} (-b+x)} \, dx=-\frac{i \pi  \log }{\sqrt{b}}\left(\frac{\Gamma \left(-\frac{i
   \log (a b)}{4 \pi }\right) \Gamma \left(1-\frac{i \log (a b)}{4 \pi }\right) \log (a b)}{\Gamma \left(\frac{1}{2}-\frac{i
   \log (a b)}{4 \pi }\right)^2}\right)
\end{multline}
\end{example}
\begin{example}
\begin{multline}
\int_0^1 \frac{(1+x) \log \left(\log \left(\frac{1}{x}\right)\right)}{\sqrt{x} (b+x) (1+b x)} \, dx\\-\frac{\sqrt{b}
   \pi  \log \left(\frac{\left(i+\sqrt{b}\right) \Gamma \left(\frac{1}{2}-\frac{i \log (-b)}{4 \pi
   }\right)^2}{\left(-i+\sqrt{b}\right) \Gamma \left(-\frac{i \log (-b)}{4 \pi }\right) \Gamma \left(1-\frac{i \log (-b)}{4
   \pi }\right)}\right)}{b (1+b)}+\frac{\pi  \log (\log (-b))}{\sqrt{b} (1+b)}
\end{multline}
\end{example}
\begin{example}
\begin{equation}
\int_0^{\infty } \frac{x^{-1-m+\alpha } \left(-1+x^{2 m}\right)}{(1+x) \log (x)} \, dx=\log \left(\frac{\tan
   \left(\frac{1}{2} \pi  (\alpha +m)\right)}{\tan \left(\frac{1}{2} \pi  (\alpha -m)\right)}\right)
\end{equation}
\end{example}
\subsubsection{Derived using equation (\ref{eq_22713})}
\begin{example}
\begin{multline}
\int_0^{\infty } \frac{x^{-1-m-\alpha +\beta } \left(x^{2 m}-1\right) \left(x^{2 \alpha
   }-1\right)}{\left(x^{2 \beta }-1\right) \log (x)} \, dx=\log \left(\frac{1+\cos \left(\frac{\pi  (m-\alpha
   )}{\beta }\right)}{1+\cos \left(\frac{\pi  (m+\alpha )}{\beta }\right)}\right)
\end{multline}
\end{example}
\begin{example}
\begin{multline}
\int_0^{\infty } \frac{x^{-1+m-\alpha +\beta } \left(-1+x^{2 \alpha }\right)}{\left(-1+x^{2 \beta }\right)
   \left(a^2 \pi ^2+\log ^2(x)\right)} \, dx\\
=-\frac{i }{2 a \pi }\left(e^{\frac{i \pi  (m-\alpha )}{\beta }} \left(\Phi
   \left(-e^{\frac{i \pi  (m-\alpha )}{\beta }},1,1-a \beta \right)-\Phi \left(-e^{\frac{i \pi  (m-\alpha )}{\beta
   }},1,1+a \beta \right)\right)\right. \\ \left.
+e^{\frac{i \pi  (m+\alpha )}{\beta }} \left(-\Phi \left(-e^{\frac{i \pi  (m+\alpha
   )}{\beta }},1,1-a \beta \right)+\Phi \left(-e^{\frac{i \pi  (m+\alpha )}{\beta }},1,1+a \beta
   \right)\right)\right)
\end{multline}
\end{example}
\begin{example}
\begin{equation}
\int_0^{\infty } \frac{x^{-1+\alpha } \log (i a \pi +\log (x))}{1+x^{2 \alpha }} \, dx=-\frac{\pi  }{2 \alpha }\log
   \left(-\frac{i \alpha  \Gamma \left(\frac{1}{4}+\frac{a \alpha }{2}\right)^2}{2 \pi  \Gamma
   \left(\frac{3}{4}+\frac{a \alpha }{2}\right)^2}\right)
\end{equation}
\end{example}
\begin{example}
\begin{equation}
\int_0^{\infty } \frac{x^{-1-m-\alpha +\beta } \left(1-x^{2 m}\right) \left(-1+x^{2 \alpha
   }\right)}{\left(-1+x^{2 \beta }\right) \log (x)} \, dx=2 \log \left(\frac{\cos \left(\frac{\pi  (m+\alpha )}{2
   \beta }\right)}{\cos \left(\frac{\pi  (m-\alpha )}{2 \beta }\right)}\right)
\end{equation}
\end{example}
\begin{example}
\begin{multline}
\int_0^{\infty } \frac{x^{-1+m-\alpha +\beta } \left(-1+x^{2 \alpha }\right) \log ^k(x)}{-1+x^{2 \beta }} \,
   dx\\
=i^{1+k} \pi ^{1+k} \beta ^{-1-k} \left(-\text{Li}_{-k}\left(-e^{\frac{i \pi  (m-\alpha )}{\beta
   }}\right)+\text{Li}_{-k}\left(-e^{\frac{i \pi  (m+\alpha )}{\beta }}\right)\right)
\end{multline}
\end{example}
\subsubsection{Derived from equation (\ref{eq_3194411})}
\begin{example}
\begin{multline}
\int_0^{\infty } x^{-1+m} \left(1+x^q \beta ^{-q}\right)^{-n} \log ^k(x) \, dx\\
=\sum _{t=0}^{-1+n} \sum _{j=0}^{-1+n} \frac{(-1)^{j+t} e^{\frac{i m \pi }{q}} (2 \pi )^{1+k-t}
   \left(-\frac{1}{q}\right)^t \left(\frac{i}{q}\right)^{-1+k-t} \left(\beta ^{-q}\right)^{-\frac{m}{q}}  \left(-\frac{m-q}{q}\right)_{-1-j+n} (1+k-t)_t S_j^{(t)}}{q^2 j! (-1-j+n)!}\\\Phi \left(e^{\frac{2 i m \pi }{q}},-k+t,\frac{\pi +i \log \left(\beta
   ^{-q}\right)}{2 \pi }\right)
\end{multline}
\end{example}
\begin{example}
\begin{multline}
\int_0^{\infty } x^{-1+\frac{m}{2}} \left(1+x^q\right)^{-n} \log ^{2 k}(x) \, dx\\
=\sum _{t=0}^{-1+n} \sum _{j=0}^{-1+n} \frac{(-1)^{j+t} e^{\frac{i m \pi }{2 q}} (2 \pi )^{1+2 k-t}
   \left(-\frac{1}{q}\right)^t \left(\frac{i}{q}\right)^{-1+2 k-t} \left(-\frac{\frac{m}{2}-q}{q}\right)_{-1-j+n} (1+2 k-t)_t S_j^{(t)}}{q^2 j!
   (-1-j+n)!}\\
\Phi \left(e^{\frac{i m \pi }{q}},-2 k+t,\frac{1}{2}\right) 
\end{multline}
\end{example}
\subsubsection{The following were derived using equation (\ref{eq_21916})}
\begin{example}
\begin{multline}
\int_0^{\infty } \frac{\log \left(\log \left(a x \left(b+x+\sqrt{x (2 b+x)}\right)\right)\right)}{x^{3/4}
   \sqrt{2 b+x} \sqrt[4]{b+x+\sqrt{x (2 b+x)}}} \, dx\\
=\frac{\pi }{2^{3/4}
   \sqrt{b}} \left(i \pi +\log (64)+2 \log (\pi )+4 \log
   \left(\frac{-2 \pi +i \log (2)-i \log \left(a b^2\right)}{-6 \pi +i \log (2)-i \log \left(a b^2\right)}\right)\right. \\ \left.
-4
   \text{log$\Gamma $}\left(-\frac{3}{4}-\frac{i \log \left(\frac{a b^2}{2}\right)}{8 \pi }\right)+4
   \text{log$\Gamma $}\left(-\frac{1}{4}-\frac{i \log \left(\frac{a b^2}{2}\right)}{8 \pi }\right)\right)
\end{multline}
\end{example}
\begin{example}
\begin{multline}
\int_0^{\infty } \frac{\log \left(\log \left(\frac{2 x \left(b+x+\sqrt{x (2
   b+x)}\right)}{b^2}\right)\right)}{x^{3/4} \sqrt{2 b+x} \sqrt[4]{b+x+\sqrt{x (2 b+x)}}} \, dx\\
=\frac{\pi  \left(i
   \pi +2 \log \left(\frac{8 \pi }{9}\right)-4 \text{log$\Gamma $}\left(-\frac{3}{4}\right)+4 \text{log$\Gamma
   $}\left(-\frac{1}{4}\right)\right)}{2^{3/4} \sqrt{b}}
\end{multline}
\end{example}
\begin{example}
\begin{multline}
\int_0^{\infty } \frac{1}{x^{3/4} \sqrt{2 b+x} \sqrt[4]{b+x+\sqrt{x (2 b+x)}} \log ^2\left(\frac{2 x
   \left(b+x+\sqrt{x (2 b+x)}\right)}{b^2}\right)} \, dx=-\frac{C}{2^{3/4} \sqrt{b} \pi }
\end{multline}
\end{example}
\begin{example}
\begin{multline}
\int_0^{\infty } \frac{\log \left(\frac{2 x \left(b+x+\sqrt{x (2 b+x)}\right)}{b^2}\right) \log \left(\log
   \left(\frac{2 x \left(b+x+\sqrt{x (2 b+x)}\right)}{b^2}\right)\right)}{x^{3/4} \sqrt{2 b+x} \sqrt[4]{b+x+\sqrt{x
   (2 b+x)}}} \, dx=-\frac{8 i \sqrt[4]{2} C \pi }{\sqrt{b}}
\end{multline}
\end{example}
\begin{example}
\begin{multline}
\int_0^{\infty } \frac{x^{-1+m} \left(b+x+\sqrt{x (2 b+x)}\right)^{-\frac{1}{2}+m}}{\sqrt{2 b+x} \log
   \left(\frac{2 x \left(b+x+\sqrt{x (2 b+x)}\right)}{b^2}\right)} \, dx=-2^{-\frac{1}{2}-m} b^{-1+2 m} \log (i \cot
   (m \pi ))
\end{multline}
\end{example}
\begin{example}
\begin{multline}
-\frac{\int_0^{\infty } \frac{\log ^2\left(2 x \left(-1+x+\sqrt{(-2+x) x}\right)\right) \log \left(\log
   \left(2 x \left(-1+x+\sqrt{(-2+x) x}\right)\right)\right)}{\sqrt{-2+x} x^{3/4} \sqrt[4]{-1+x+\sqrt{(-2+x) x}}} \,
   dx}{56 i \sqrt[4]{2} \pi }=\zeta (3)
\end{multline}
\end{example}
\begin{example}
\begin{multline}
\int_0^{\infty } \frac{\log ^2\left(\frac{2 x \left(b+x+\sqrt{x (2 b+x)}\right)}{b^2}\right) \log \left(\log
   \left(\frac{2 x \left(b+x+\sqrt{x (2 b+x)}\right)}{b^2}\right)\right)}{x^{3/4} \sqrt{2 b+x} \sqrt[4]{b+x+\sqrt{x
   (2 b+x)}}} \, dx\\
=\frac{\sqrt[4]{2} \pi  \left(\left(4 \pi ^2+0^2\right) \log (4 i \pi )+32 \pi ^2
   \Phi'\left(-1,-2,\frac{1}{2}\right)\right)}{\sqrt{b}}
\end{multline}
\end{example}
\begin{example}
\begin{multline}
\int_0^{\infty } \frac{\log \left(2 x \left(-1+x+\sqrt{(-2+x) x}\right)\right) \log \left(\log \left(2 x
   \left(-1+x+\sqrt{(-2+x) x}\right)\right)\right)}{\sqrt{-2+x} x^{3/4} \sqrt[4]{-1+x+\sqrt{(-2+x) x}}} \,
   dx\\
=-\frac{i \pi  \left(-16 i \pi  \log \left(\frac{A^3}{\sqrt[3]{2} \sqrt[4]{e}}\right)+i \pi  \left(2 i \pi
   +\log \left(256 \pi ^4\right)\right)\right)}{2^{3/4}}
\end{multline}
\end{example}
\begin{example}
\begin{multline}
\int_0^{\infty } \frac{x^{-1+m} \left(-1+x+\sqrt{(-2+x) x}\right)^{-\frac{1}{2}+m} \log ^k\left(2 x
   \left(-1+x+\sqrt{(-2+x) x}\right)\right)}{\sqrt{-2+x}} \, dx\\
=i 2^{\frac{3}{2}+2 k-m} e^{\frac{i k \pi }{2}} \pi
   ^{1+k} \text{Li}_{-k}\left(e^{4 i m \pi }\right)
\end{multline}
\end{example}
\begin{example}
\begin{multline}
\int_0^{\infty } \frac{x^{-1+m} \left(-1+x+\sqrt{(-2+x) x}\right)^{-\frac{1}{2}+m} \log ^k\left(2 x
   \left(-1+x+\sqrt{(-2+x) x}\right)\right)}{\sqrt{-2+x}} \, dx\\
=2^{\frac{1}{2}+k-m} \Gamma (1+k) \left(\zeta (1+k,-2
   m)-e^{i k \pi } \zeta (1+k,1+2 m)\right)
\end{multline}
\end{example}
\subsubsection{These were derived using equation (\ref{eq_12238})}
\begin{example}
\begin{multline}
\int_0^{\infty } \frac{\log ^k\left(a x \left(x+\sqrt{b^2+x^2}\right)\right)}{\sqrt{x}
   \left(x+\sqrt{b^2+x^2}\right)^{3/2}} \, dx\\
=\frac{2^{-\frac{3}{2}+k} e^{\frac{1}{2} i (1+k) \pi } \pi ^k }{b}\left(k
   \left(2^{-1+k} \zeta \left(1-k,-\frac{i \left(i \pi +\log \left(\frac{a}{2}\right)+\log \left(b^2\right)\right)}{4
   \pi }\right)\right.\right. \\ \left.\left.
-2^{-1+k} \zeta \left(1-k,\frac{1}{2} \left(1-\frac{i \left(i \pi +\log \left(\frac{a}{2}\right)+\log
   \left(b^2\right)\right)}{2 \pi }\right)\right)\right)\right. \\ \left.
-3 i \pi  \left(2^k \zeta \left(-k,-\frac{i \left(i \pi +\log
   \left(\frac{a}{2}\right)+\log \left(b^2\right)\right)}{4 \pi }\right)\right.\right. \\ \left.\left.
-2^k \zeta \left(-k,\frac{1}{2}
   \left(1-\frac{i \left(i \pi +\log \left(\frac{a}{2}\right)+\log \left(b^2\right)\right)}{2 \pi
   }\right)\right)\right)\right)
\end{multline}
\end{example}
\begin{example}
\begin{multline}
\int_0^{\infty } \frac{\log ^k\left(a x \left(x+\sqrt{b^2+x^2}\right)\right)}{x^{3/4}
   \left(x+\sqrt{b^2+x^2}\right)^{7/4}} \, dx\\
=\frac{2^{-\frac{5}{4}+k} e^{\frac{1}{2} i \left(\frac{1}{2}+k\right)
   \pi } \pi ^k }{b^{3/2}}\left(k \Phi \left(i,1-k,-\frac{i \left(i \pi +\log \left(\frac{a}{2}\right)+\log
   \left(b^2\right)\right)}{2 \pi }\right)\right. \\ \left.
-\frac{7}{2} i \pi  \Phi \left(i,-k,-\frac{i \left(i \pi +\log
   \left(\frac{a}{2}\right)+\log \left(b^2\right)\right)}{2 \pi }\right)\right)
\end{multline}
\end{example}
\begin{example}
\begin{multline}
\int_0^{\infty } x^{-1+m} \left(x+\sqrt{1+x^2}\right)^{-2+m} \log ^k\left(-2 x
   \left(x+\sqrt{1+x^2}\right)\right) \, dx\\
=2^{-1+k-m} e^{\frac{1}{2} i (k+2 m) \pi } \pi ^k \left(e^{-2 i m \pi } k
   \text{Li}_{1-k}\left(e^{2 i m \pi }\right)+2 i e^{-2 i m \pi } (-2+m) \pi  \text{Li}_{-k}\left(e^{2 i m \pi
   }\right)\right)
\end{multline}
\end{example}
\begin{example}
\begin{multline}
\int_0^{\infty } \frac{1}{\sqrt{x} \left(x+\sqrt{1+x^2}\right)^{3/2} \log ^2\left(-2 x
   \left(x+\sqrt{1+x^2}\right)\right)} \, dx=-\frac{\pi ^3-6 i \zeta (3)}{32 \sqrt{2} \pi ^2}
\end{multline}
\end{example}
\begin{example}
\begin{multline}
\int_0^{\infty } \frac{1}{\sqrt{x} \left(x+\sqrt{1+x^2}\right)^{3/2} \log \left(-2 x
   \left(x+\sqrt{1+x^2}\right)\right)} \, dx=-\frac{\pi +36 i \log (2)}{48 \sqrt{2}}
\end{multline}
\end{example}
\begin{example}
\begin{multline}
\int_0^{\infty } \frac{\log \left(-2 x \left(x+\sqrt{1+x^2}\right)\right) \log \left(\log \left(-2 x
   \left(x+\sqrt{1+x^2}\right)\right)\right)}{\sqrt{x} \left(x+\sqrt{1+x^2}\right)^{3/2}} \, dx\\
-\frac{i
   \left(\frac{1}{2}-\frac{3 i \pi }{4}\right) \pi ^2}{2 \sqrt{2}}-\frac{\left(\frac{1}{2}-\frac{3 i \pi }{4}\right)
   \pi  \log (2)}{\sqrt{2}}-\frac{\left(\frac{1}{2}-\frac{3 i \pi }{4}\right) \pi  \log (\pi )}{\sqrt{2}}\\
-\frac{\pi 
   \left(\frac{1}{2}+3 i \pi  \log \left(\frac{A^3}{\sqrt[3]{2}
   \sqrt[4]{e}}\right)-\Phi'(-1,0,1)\right)}{\sqrt{2}}
\end{multline}
\end{example}
\begin{example}
\begin{multline}
\int_0^{\infty } \frac{\log ^2\left(-2 x \left(x+\sqrt{1+x^2}\right)\right) \log \left(\log \left(-2 x
   \left(x+\sqrt{1+x^2}\right)\right)\right)}{\sqrt{x} \left(x+\sqrt{1+x^2}\right)^{3/2}} \, dx\\
=\frac{\pi ^3}{2
   \sqrt{2}}-\frac{i \pi ^2 \log (2)}{\sqrt{2}}-\frac{i \pi ^2 \log (\pi )}{\sqrt{2}}-i \sqrt{2} \pi ^2
   \left(\frac{1}{4}-2 \log \left(\frac{A^3}{\sqrt[3]{2} \sqrt[4]{e}}\right)+\frac{21 i \zeta (3)}{4 \pi
   }\right)
\end{multline}
\end{example}
\begin{example}
\begin{multline}
\int_0^{\infty } \frac{\log \left(\log \left(2 x \left(x+\sqrt{1+x^2}\right)\right)\right)}{\sqrt{x}
   \left(x+\sqrt{1+x^2}\right)^{3/2} \log \left(2 x \left(x+\sqrt{1+x^2}\right)\right)} \, dx\\
=\frac{8 C (2-i \pi -2
   \log (2 \pi ))+3 \pi ^2 \left(2 i \gamma +\pi +i \log \left(\frac{16 \pi ^4}{\Gamma
   \left(\frac{1}{4}\right)^8}\right)\right)+4 \Phi'\left(-1,2,\frac{1}{2}\right)}{16
   \sqrt{2} \pi }
\end{multline}
\end{example}
\begin{example}
\begin{multline}
\int_0^{\infty } e^{(m-2) t} \cosh (t) \log ^k\left(a e^t \sinh (t)\right) \sinh ^{m-1}(t) \, dt\\
=2^{-1+k-m}
   e^{\frac{1}{2} i (k+2 m) \pi } \pi ^k \left(k \Phi \left(e^{2 i m \pi },1-k,\frac{\pi -i \log
   \left(\frac{a}{2}\right)}{2 \pi }\right)\right. \\ \left.+2 i (-2+m) \pi  \Phi \left(e^{2 i m \pi },-k,\frac{\pi -i \log
   \left(\frac{a}{2}\right)}{2 \pi }\right)\right)
\end{multline}
\end{example}
\begin{example}
\begin{multline}
\int_0^{\infty } e^{(-2+m) t} \cosh (t) \log ^k\left(-2 e^t \sinh (t)\right) \sinh ^{-1+m}(t) \, dt\\
=,2^{-1+k-m}
   e^{\frac{1}{2} i (k+2 m) \pi } \pi ^k \left(e^{-2 i m \pi } k \text{Li}_{1-k}\left(e^{2 i m \pi }\right)+2 i e^{-2
   i m \pi } (-2+m) \pi  \text{Li}_{-k}\left(e^{2 i m \pi }\right)\right)
\end{multline}
\end{example}
\begin{example}
\begin{multline}
\int_0^{\infty } \frac{e^{(-2+m) t} \cosh (t) \sinh ^{-1+m}(t)}{i a+t+\log (\sinh (t))} \, dt\\
=\frac{2^{-2-m}
   e^{\frac{1}{2} i (-1+2 m) \pi }}{\pi } \left(2 i (-2+m) \pi  \Phi \left(e^{2 i m \pi },1,\frac{\pi -i \log
   \left(\frac{e^{i a}}{2}\right)}{2 \pi }\right)\right. \\ \left.
   -\Phi \left(e^{2 i m \pi },2,\frac{\pi -i \log \left(\frac{e^{i
   a}}{2}\right)}{2 \pi }\right)\right)
\end{multline}
\end{example}
\begin{example}
\begin{multline}
\int_0^{\infty } \frac{e^{(-2+m) t} \cosh (t) \sinh ^{-1+m}(t)}{(i a-b+t+\log (\sinh (t))) (i a+b+t+\log
   (\sinh (t)))} \, dt\\
=-\frac{2^{-3-m} e^{i m \pi } }{b \pi }\left(2 (-2+m) \pi  \left(\Phi \left(e^{2 i m \pi },1,\frac{a-i
   b+\pi +i \log (2)}{2 \pi }\right)\right.\right. \\ \left.\left.
-\Phi \left(e^{2 i m \pi },1,\frac{a+i b+\pi +i \log (2)}{2 \pi }\right)\right)+i
   \left(\Phi \left(e^{2 i m \pi },2,\frac{a-i b+\pi +i \log (2)}{2 \pi }\right)\right.\right. \\ \left.\left.-\Phi \left(e^{2 i m \pi
   },2,\frac{a+i b+\pi +i \log (2)}{2 \pi }\right)\right)\right)
\end{multline}
\end{example}
\begin{example}
\begin{multline}
\int_0^{\infty } \frac{e^{-\frac{3 t}{2}} \cosh (t)}{(i a-b+t+\log (\sinh (t))) (i a+b+t+\log (\sinh (t)))
   \sqrt{\sinh (t)}} \, dt\\
=-\frac{6 i \pi }{32 \sqrt{2} b
   \pi } \left(\psi ^{(0)}\left(\frac{a-i b+\pi +i \log (2)}{4 \pi }\right)-\psi
   ^{(0)}\left(\frac{a+i b+\pi +i \log (2)}{4 \pi }\right)\right. \\ \left.
-\psi ^{(0)}\left(\frac{a-i b+3 \pi +i \log (2)}{4 \pi
   }\right)+\psi ^{(0)}\left(\frac{a+i b+3 \pi +i \log (2)}{4 \pi }\right)\right)\\
+\psi ^{(1)}\left(\frac{a-i b+\pi +i
   \log (2)}{4 \pi }\right)-\psi ^{(1)}\left(\frac{a+i b+\pi +i \log (2)}{4 \pi }\right)\\
-\psi ^{(1)}\left(\frac{a-i
   b+3 \pi +i \log (2)}{4 \pi }\right)+\psi ^{(1)}\left(\frac{a+i b+3 \pi +i \log (2)}{4 \pi }\right)
\end{multline}
\end{example}
\begin{example}
\begin{multline}
\int_0^{\infty } \frac{e^{-\frac{3 t}{2}} \cosh (t)}{\left(2 i \pi +t+\log \left(\frac{\sinh
   (t)}{2}\right)\right) (2 i \pi +t+\log (2 \sinh (t))) \sqrt{\sinh (t)}} \, dt\\
=\frac{16-16 C-6 i \pi  \left(-4+\pi
   +H_{\frac{\pi +2 i \log (2)}{4 \pi }}-H_{-\frac{1}{4}+\frac{i \log (2)}{2 \pi }}\right)-\psi
   ^{(1)}\left(\frac{3}{4}+\frac{i \log (2)}{2 \pi }\right)+\psi ^{(1)}\left(\frac{5}{4}+\frac{i \log (2)}{2 \pi
   }\right)}{32 \sqrt{2} \pi  \log (2)}
\end{multline}
\end{example}
\subsubsection{Following derived using equation (\ref{eq_23113})}
\begin{example}
\begin{multline}
\int_0^{\infty } \frac{e^{-x} \log \left(\log \left(\frac{a e^{-x}}{\sqrt[8]{-1+\cosh (8
   x)}}\right)\right)}{\sqrt[8]{-1+\cosh (8 x)}} \, dx\\
=\frac{\pi }{8\times 2^{3/8}} \left(i \pi +\log \left(\frac{\pi
   ^2}{4}\right)-(2-2 i) \left(\log \left(\frac{\Gamma \left(-\frac{i (2 i \pi +\log (2)+8 \log (a))}{16 \pi
   }\right)}{2 \Gamma \left(-\frac{i (10 i \pi +\log (2)+8 \log (a))}{16 \pi }\right)}\right)\right.\right. \\ \left.\left.+i \log
   \left(\frac{\Gamma \left(-\frac{i (6 i \pi +\log (2)+8 \log (a))}{16 \pi }\right)}{2 \Gamma \left(-\frac{i (14 i
   \pi +\log (2)+8 \log (a))}{16 \pi }\right)}\right)\right)\right)
\end{multline}
\end{example}
\subsubsection{This section is derived from equation (\ref{eq_2941})}
\begin{example}
\begin{multline}
\int_0^{\infty } x^{-1+s} \log \left(1-e^{-2 x}\right) \left(-2 m (1+s) x \cos (m x)+\left(-s (1+s)+\left(1+m^2\right)
   x^2\right) \sin (m x)\right)\\ \sinh (x) \, dx\\
=\frac{1}{2}
   \left(\left(-(-i-m)^{-2-s}-(i-m)^{-2-s}+(-i+m)^{-2-s}+(i+m)^{-2-s}\right) \Gamma (2+s)\right. \\ \left.
+\pi ^{2+s}
   \left(-\text{Li}_{-1-s}\left(-e^{-m \pi }\right)+\text{Li}_{-1-s}\left(-e^{m \pi }\right)\right)\right) \sec
   \left(\frac{1}{2} \pi  (1+s)\right)
\end{multline}
\end{example}
\begin{example}
\begin{multline}
\int_0^{\infty } \left(\tan ^{-1}\left(\frac{x}{a}\right)+\frac{2 x a}{\left(a^2+x^2\right)^2}\right) \log \left(1-e^{-2
   x}\right) \sinh (x) \, dx\\
=-i e^{-i a} \Gamma (0,-i a)+i e^{i a} \Gamma (0,i a)+\pi  \log \left(\frac{\sqrt{2 \pi } \Gamma
   \left(\frac{a+\pi }{2 \pi }\right)}{\sqrt{a} \Gamma \left(\frac{a}{2 \pi }\right)}\right)
\end{multline}
\end{example}
\begin{example}
\begin{multline}
\int_0^{\infty } e^{-i m x} \left(\frac{8}{(\pi -2 i x)^3}-\frac{8 e^{2 i m x}}{(\pi +2 i x)^3}+\frac{-4 m+\pi +m^2 \pi
   -2 i \left(1+m^2\right) x}{(\pi -2 i x)^2}\right. \\ \left.
-\frac{e^{2 i m x} \left(-4 m+\pi +m^2 \pi +2 i \left(1+m^2\right) x\right)}{(\pi
   +2 i x)^2}\right) \log \left(1-e^{-2 x}\right) \sinh (x) \, dx\\
=i+e^{-\frac{1}{2} (m \pi )} \left(-2 i \tan
   ^{-1}\left(e^{\frac{m \pi }{2}}\right)-\Gamma \left(0,-\frac{1}{2} (-i+m) \pi \right)+\Gamma \left(0,-\frac{1}{2} (i+m) \pi
   \right)\right)
\end{multline}
\end{example}
\begin{example}
\begin{multline}
\int_0^{\infty } x^{-1+s} \log \left(1-e^{-2 x}\right) \left(\left(-s (1+s)+\left(1+m^2\right) x^2\right) \cos (m x)+2 m
   (1+s) x \sin (m x)\right)\\ \sinh (x) \, dx\\
=\frac{1}{2}
   \left(\left(-(-i-m)^{-2-s}-(i-m)^{-2-s}-(-i+m)^{-2-s}-(i+m)^{-2-s}\right) \Gamma (2+s)\right. \\ \left.
+\pi ^{2+s}
   \left(\text{Li}_{-1-s}\left(-e^{-m \pi }\right)+\text{Li}_{-1-s}\left(-e^{m \pi }\right)\right)\right) \sec \left(\frac{\pi 
   s}{2}\right)
\end{multline}
\end{example}
\begin{example}
\begin{multline}
\int_0^{\infty } \left(\tan ^{-1}\left(\frac{x}{a}\right)+\frac{2 x a}{\left(a^2+x^2\right)^2}\right) \log \left(1-e^{-2
   x}\right) \sinh (x) \, dx\\
=-i e^{-i a} \Gamma (0,-i a)+i e^{i a} \Gamma (0,i a)+\pi  \log \left(\frac{\sqrt{2 \pi } \Gamma
   \left(\frac{a+\pi }{2 \pi }\right)}{\sqrt{a} \Gamma \left(\frac{a}{2 \pi }\right)}\right)
\end{multline}
\end{example}
\begin{example}
\begin{multline}
\int_0^{\infty } e^{-i m x} \left(\frac{8}{(\pi -2 i x)^3}-\frac{8 e^{2 i m x}}{(\pi +2 i x)^3}+\frac{-4 m+\pi +m^2 \pi
   -2 i \left(1+m^2\right) x}{(\pi -2 i x)^2}\right. \\ \left.
-\frac{e^{2 i m x} \left(-4 m+\pi +m^2 \pi +2 i \left(1+m^2\right) x\right)}{(\pi
   +2 i x)^2}\right) \log \left(1-e^{-2 x}\right) \sinh (x) \, dx\\
=i+e^{-\frac{1}{2} (m \pi )} \left(-2 i \tan
   ^{-1}\left(e^{\frac{m \pi }{2}}\right)-\Gamma \left(0,-\frac{1}{2} (-i+m) \pi \right)+\Gamma \left(0,-\frac{1}{2} (i+m) \pi
   \right)\right)
\end{multline}
\end{example}
\begin{example}
\begin{multline}
\int_0^{\infty } x^{-1+s} \log \left(1-e^{-2 x}\right) \left(\left(-s (1+s)+\left(1+m^2\right) x^2\right) \cos (m x)+2 m
   (1+s) x \sin (m x)\right)\\ \sinh (x) \, dx\\
=\frac{1}{2}
   \left(\left(-(-i-m)^{-2-s}-(i-m)^{-2-s}-(-i+m)^{-2-s}-(i+m)^{-2-s}\right) \Gamma (2+s)\right. \\ \left.
+\pi ^{2+s}
   \left(\text{Li}_{-1-s}\left(-e^{-m \pi }\right)+\text{Li}_{-1-s}\left(-e^{m \pi }\right)\right)\right) \sec \left(\frac{\pi 
   s}{2}\right)
\end{multline}
\end{example}
\begin{example}
\begin{multline}
\int_0^{\infty } e^{-\frac{x}{2}} \log \left(1-e^{-2 x}\right) \left(-9 \left(-1+e^x\right) \pi ^2+x \left(-16-9 x+e^x
   (-16+9 x)\right)\right. \\ \left.
+2 i \pi  \left(8+9 x+e^x (-8+9 x)\right)+2 e^x \left(8+3 \pi ^2+2 i \pi  (4-3 x)+(8-3 x) x\right) \log (\pi
   -i x)\right. \\ \left.
-2 \left(8+3 \pi ^2+2 i \pi  (4+3 x)-x (8+3 x)\right) \log (\pi +i x)\right) \sinh (x) \, dx\\
=\frac{8}{27} i \left(32
   \left(-27 \Gamma \left(0,\frac{i \pi }{2}\right)+\Gamma \left(0,-\frac{1}{2} (3 i \pi )\right)+26 i \log (\pi )\right)\right. \\ \left.
-480\pi  (-1+\log (\pi ))+9 \pi ^2 \left((-4 i+3 \pi ) (-3+2 \log (\pi ))+6 i \pi  Li'_{-2}(-i)\right)\right)
\end{multline}
\end{example}
\begin{example}
\begin{multline}
\int_0^{\infty } i e^{-i m x} \left((-1+k) k (-i x+\log (a))^{-2+k}-e^{2 i m x} (-1+k) k (i x+\log (a))^{-2+k}\right. \\ \left.
+2 k m
   \left((-i x+\log (a))^{-1+k}-e^{2 i m x} (i x+\log (a))^{-1+k}\right)\right. \\ \left.+\left(1+m^2\right) \left((-i x+\log (a))^k-e^{2 i m x}
   (i x+\log (a))^k\right)\right) \log \left(1-e^{-2 x}\right) \sinh (x) \, dx\\
=2 a^{-i-m} k \left(-a^{2 i} (i-m)^{-1-k} \Gamma(k,-((-i+m) \log (a)))\right. \\ \left.-(-i-m)^{-1-k} \Gamma (k,-((i+m) \log (a)))\right)+\pi  \left(-2 e^{m \pi } \pi ^k \Phi \left(-e^{m \pi},-k,\frac{\pi +\log (a)}{\pi }\right)+\log ^k(a)\right)\\
-2 i m \left(-(i-m)^{-1-k} (-((-i+m) \log (a)))^k+(-i-m)^{-1-k} (-((i+m) \log (a)))^k\right)
\end{multline}
\end{example}
\begin{example}
\begin{multline}
\int_0^{\infty } \log \left(1-e^{-2 x}\right) \left(\frac{1}{-i x+\log (a)}-\frac{1}{i x+\log (a)}+(-i x+\log (a)) \log
   (-i x+\log (a))\right. \\ \left.
-(i x+\log (a)) \log (i x+\log (a))\right) \sinh (x) \, dx\\
=-i \left(4+2 a^{-i} \left(a^i \log (-i \log
   (a))+G_{1,2}^{2,0}\left(-i \log (a)\left|
\begin{array}{c}
 1 \\
 0,0 \\
\end{array}
\right.\right)\right.\right. \\ \left.\left.
+a^{2 i} \left(a^{-i} \log (i \log (a))+G_{1,2}^{2,0}\left(i \log (a)\left|
\begin{array}{c}
 1 \\
 0,0 \\
\end{array}
\right.\right)\right)\right)\right. \\ \left.+\pi  \left(-2 \pi  \left(2 \zeta \left(-1,\frac{\pi +\log (a)}{2 \pi }\right)-2 \zeta\left(-1,\frac{1}{2} \left(1+\frac{\pi +\log (a)}{\pi }\right)\right)\right) \log (\pi )\right.\right. \\ \left.\left.
+\log (a) \log (\log (a))+2 \pi \Phi'\left(-1,-1,\frac{\pi +\log (a)}{\pi }\right)\right)\right)
\end{multline}
\end{example}
\begin{example}
\begin{multline}
\int_0^{\infty } \log \left(1-e^{-2 x}\right) \left(80 \left(-\frac{1}{(\pi -2 i x)^4}+\frac{1}{(\pi +2 i
   x)^4}\right)\right. \\ \left.
+\frac{4 \left(24+(\pi -2 i x)^2\right) \log \left(\frac{1}{2} (\pi -2 i x)\right)}{(\pi -2 i x)^4}\right. \\ \left.-\frac{4
   \left(24+(\pi +2 i x)^2\right) \log \left(\frac{1}{2} (\pi +2 i x)\right)}{(\pi +2 i x)^4}\right) \sinh (x) \, dx\\
=\frac{i}{2 \pi }
   \left(\pi  \text{Ci}\left(\frac{\pi }{2}\right) (4+\log (256)-8 \log (\pi ))\right. \\ \left.
+8 \pi  \left(G_{2,3}^{3,0}\left(-\frac{1}{2} (i
   \pi )|
\begin{array}{c}
 1,1 \\
 -2,0,0 \\
\end{array}
\right)+G_{2,3}^{3,0}\left(\frac{i \pi }{2}|
\begin{array}{c}
 1,1 \\
 -2,0,0 \\
\end{array}
\right)\right)\right. \\ \left.-4 \left(2+\log (4)-6 \log (\pi )+4 C \log (\pi
   )+\Phi'\left(-1,2,\frac{3}{2}\right)\right)\right)
\end{multline}
\end{example}
\begin{example}
\begin{multline}
\int_0^{\infty } \frac{\log \left(1-e^{-2 x}\right) \sinh (x)}{\left(\pi ^2+x^2\right)^3}\left(4 m x \left(\pi ^2+x^2\right) \cos (mx)+\left(\left(1+m^2\right) \pi ^4-6 x^2+\left(1+m^2\right) x^4\right.\right. \\ \left.\left.+2 \pi ^2 \left(1+\left(1+m^2\right) x^2\right)\right) \sin (m x)\right)  \, dx\\
=\frac{1}{2 \pi }\left(e^{-m \pi } \left(\Gamma (0,-((-i+m) \pi ))+\Gamma (0,-((i+m) \pi
   ))\right.\right. \\ \left.\left.
-e^{2 m \pi } \left(\Gamma (0,(-i+m) \pi )+\Gamma (0,(i+m) \pi )+\log \left(1+e^{-m \pi }\right)\right)+\log \left(1+e^{m
   \pi }\right)\right)\right)
\end{multline}
\end{example}
\subsubsection{Formaule derived using equation (\ref{eq_2248})}
\begin{example}
\begin{multline}
\int_0^z \frac{x^{-1+m} \log ^k\left(a x \left(-x^{2 m}+z^{2 m}\right)^{-\frac{1}{2 m}}\right)}{\sqrt{-x^{2
   m}+z^{2 m}}} \, dx\\
   =-i \left(\frac{i}{m}\right)^{1+k} \pi ^{1+k} \left(2^k \zeta \left(-k,\frac{\pi -2 i m \log
   (a)}{4 \pi }\right)-2^k \zeta \left(-k,\frac{1}{2} \left(1+\frac{\pi -2 i m \log (a)}{2 \pi
   }\right)\right)\right)
\end{multline}
\end{example}
\begin{example}
\begin{multline}
\int_0^z \frac{x^{-1+m} \log \left(\log \left(a x \left(-x^{2 m}+z^{2 m}\right)^{-\frac{1}{2
   m}}\right)\right)}{\sqrt{-x^{2 m}+z^{2 m}}} \, dx\\
   =\frac{\pi }{2 m} \left(\log \left(\frac{2 i \pi }{m}\right)-2 \log (-3
   \pi -2 i m \log (a))+2 \log (-\pi -2 i m \log (a))\right. \\ \left.
   +2 \text{log$\Gamma $}\left(-\frac{\pi +2 i m \log (a)}{4 \pi
   }\right)-2 \text{log$\Gamma $}\left(-\frac{3}{4}-\frac{i m \log (a)}{2 \pi }\right)\right)
\end{multline}
\end{example}
\begin{example}
\begin{multline}
\int_0^z \frac{x^{-1+m} \log ^k\left(a x \left(-x^m+z^m\right)^{-1/m}\right)}{-x^m+z^m} \,
   dx=\left(\frac{i}{m}\right)^{1+k} (2 \pi )^{1+k} \zeta \left(-k,\frac{\pi -i m \log (a)}{2 \pi }\right)
\end{multline}
\end{example}
\begin{example}
\begin{multline}
\int_0^z \frac{x^{-1+m} \log \left(\log \left(a x \left(-x^m+z^m\right)^{-1/m}\right)\right)}{-x^m+z^m} \,
   dx\\
   =\frac{i }{m}\left(i m \log (a) \log \left(\frac{2 i \pi }{m}\right)+\pi  (\log (8)+3 \log (\pi )\right. \\ \left.
   -2 \log (-\pi -i m
   \log (a)))-2 \pi  \text{log$\Gamma $}\left(-\frac{\pi +i m \log (a)}{2 \pi }\right)\right)
\end{multline}
\end{example}
\begin{example}
\begin{multline}
\int_0^z x^{-1+m} \left(-x^u+z^u\right)^{-\frac{m}{u}} \log ^k\left(e^{\frac{i \pi }{u}} x
   \left(-x^u+z^u\right)^{-1/u}\right) \, dx\\
   =-e^{-\frac{i m \pi }{u}} (2 \pi )^{1+k} \left(\frac{i}{u}\right)^{1+k}
   \text{Li}_{-k}\left(e^{\frac{2 i m \pi }{u}}\right)
\end{multline}
\end{example}
\begin{example}
\begin{multline}
\int_0^z x^{-1+m} \left(-x^u+z^u\right)^{-\frac{m}{u}} \log ^k\left(e^{\frac{i \pi }{u}} x
   \left(-x^u+z^u\right)^{-1/u}\right) \, dx\\
   =\frac{i^{-k} e^{-\frac{i m \pi }{u}} \left(\frac{i}{u}\right)^k \Gamma
   (1+k) \left(-\zeta \left(1+k,-\frac{m}{u}\right)+i^{2 k} \zeta \left(1+k,\frac{m+u}{u}\right)\right)}{u}
\end{multline}
\end{example}
\begin{example}
\begin{multline}
\int_0^z \frac{x^{-1+m} \left(-x^u+z^u\right)^{-\frac{m}{u}}}{i \pi +u \log \left(x
   \left(-x^u+z^u\right)^{-1/u}\right)} \, dx\\
   =\frac{e^{-\frac{i m \pi }{u}} }{2u^2}\left(i \pi  (2 m+u)+u \log \left(4 \pi
   ^2\right)\right. \\ \left.
   -2 u \left(\log \left(\frac{m}{u}\right)+\log \left(-\frac{m+u}{u}\right)\right)-2 u
   \left(\text{log$\Gamma $}\left(\frac{m}{u}\right)+\text{log$\Gamma $}\left(-\frac{m+u}{u}\right)\right)\right)
\end{multline}
\end{example}
\begin{example}
\begin{equation}
\int_0^z \frac{x^m \left(-x^u+z^u\right)^{-\frac{m}{u}}-x^r \left(-x^u+z^u\right)^{-\frac{r}{u}}}{x \log
   \left(x \left(-x^u+z^u\right)^{-1/u}\right)} \, dx=\log \left(\frac{\tan \left(\frac{m \pi }{2 u}\right)}{\tan
   \left(\frac{\pi  r}{2 u}\right)}\right)
\end{equation}
\end{example}
\begin{example}
\begin{multline}
\int_0^z \frac{x^{-1+m} \log \left(e^{\frac{i \pi }{2 m}} x \left(-x^{2 m}+z^{2 m}\right)^{-\frac{1}{2
   m}}\right) \log \left(\log \left(e^{\frac{i \pi }{2 m}} x \left(-x^{2 m}+z^{2 m}\right)^{-\frac{1}{2
   m}}\right)\right)}{\sqrt{-x^{2 m}+z^{2 m}}} \, dx\\
   =-\frac{i \pi  \left(4 \pi  \log \left(\frac{A^3}{\sqrt[3]{2}
   \sqrt[4]{e}}\right)-\pi  \log \left(\frac{i \pi }{m}\right)\right)}{4 m^2}
\end{multline}
\end{example}
\begin{example}
\begin{multline}
\int_0^z \frac{x^{-1+m} \log ^2\left(e^{\frac{i \pi }{2 m}} x \left(-x^{2 m}+z^{2 m}\right)^{-\frac{1}{2
   m}}\right) \log \left(\log \left(e^{\frac{i \pi }{2 m}} x \left(-x^{2 m}+z^{2 m}\right)^{-\frac{1}{2
   m}}\right)\right)}{\sqrt{-x^{2 m}+z^{2 m}}} \, dx\\
   =\frac{7 \pi  \zeta (3)}{4 m^3}
\end{multline}
\end{example}
\begin{example}
\begin{equation}
\int_0^z \frac{x^{-1+m}}{\sqrt{-x^{2 m}+z^{2 m}} \log ^2\left(x \left(-x^{2 m}+z^{2 m}\right)^{-\frac{1}{2
   m}}\right)} \, dx=-\frac{4 C m}{\pi }
\end{equation}
\end{example}
\begin{example}
\begin{multline}
\int_0^z \frac{x^{-1+m} \log \left(x \left(-x^{2 m}+z^{2 m}\right)^{-\frac{1}{2 m}}\right) \log \left(\log
   \left(x \left(-x^{2 m}+z^{2 m}\right)^{-\frac{1}{2 m}}\right)\right)}{\sqrt{-x^{2 m}+z^{2 m}}} \, dx=-\frac{i C
   \pi }{m^2}
\end{multline}
\end{example}
\subsubsection{From equation (\ref{eq_23})}
\begin{example}
\begin{multline}
\int _0^{\infty }\int _0^{\infty }\frac{y^5 \Gamma (4 i+x) \sqrt{\log \left(\frac{i y}{\sqrt[8]{x}}\right)}}{\sqrt[4]{x} \left(1+y^8\right) \Gamma ((4+4
   i)+x)}dydx\\
=\left(-\frac{7}{10200}+\frac{11 i}{10200}\right) \pi ^{5/2} \left((30+15 i) (1+4 i)^{3/4} \sqrt{2} \zeta \left(-\frac{1}{2},\frac{3}{4}+\frac{i \log (1+4 i)}{8 \pi
   }\right)\right. \\ \left.
-(30+15 i) (1+4 i)^{3/4} \sqrt{2} \zeta \left(-\frac{1}{2},\frac{5}{4}+\frac{i \log (1+4 i)}{8 \pi }\right)\right. \\ \left.
+(39-48 i) \sqrt[4]{-22-4 i} \zeta
   \left(-\frac{1}{2},\frac{3}{4}+\frac{i \log (2+4 i)}{8 \pi }\right)\right. \\ \left.
+(48+39 i) (1+2 i)^{3/4} \sqrt[4]{2} \zeta \left(-\frac{1}{2},\frac{5}{4}+\frac{i \log (2+4 i)}{8 \pi }\right)\right. \\ \left.
+(6+7 i)\left(\sqrt{4+22 i} \zeta \left(-\frac{1}{2},\frac{3}{4}+\frac{i \log (3+4 i)}{8 \pi }\right)\right.\right. \\ \left.\left.
-\sqrt{4+22 i} \zeta \left(-\frac{1}{2},\frac{5}{4}+\frac{i \log (3+4 i)}{8 \pi
   }\right)+(3+4 i) (-1)^{7/8} \left(\zeta \left(-\frac{1}{2},\frac{11}{16}+\frac{i \log (4)}{8 \pi }\right)\right.\right.\right. \\ \left.\left.\left.-\zeta \left(-\frac{1}{2},\frac{19}{16}+\frac{i \log (4)}{8 \pi}\right)\right)\right)\right)
\end{multline}
\end{example}
\begin{example}
\begin{multline}
\int _0^{\infty }\int _0^{\infty }\frac{y^2 \Gamma (3 i+x) \log \left(\log \left(\frac{i y}{\sqrt[4]{x}}\right)\right)}{\sqrt[4]{x} \left(1+3 y^4\right) \Gamma ((3+3
   i)+x)}dydx\\
=\left(\frac{2}{195}+\frac{i}{780}\right) \pi ^2 \left(-\frac{1}{2} \left((1+8 i) (-1)^{7/8}+(5+i) (1+3 i)^{3/4} \sqrt[4]{3}\right.\right. \\ \left.\left.
-(2+i) (2+3 i)^{3/4} \sqrt[4]{3}\right) (\pi -2 i
   \log (2 \pi ))+2 \left((8-i) (-1)^{7/8} \log \left(\frac{\Gamma \left(\frac{7}{16}+\frac{i \log (3)}{4 \pi }\right)}{\Gamma \left(\frac{15}{16}+\frac{i \log (3)}{4 \pi
   }\right)}\right)\right.\right. \\ \left.\left.
+\sqrt[4]{3} \left((1-5 i) (1+3 i)^{3/4} \log \left(\frac{\Gamma \left(\frac{8 \pi -2 \tan ^{-1}(3)+i \log (90)}{16 \pi }\right)}{\Gamma \left(\frac{16 \pi -2 \tan
   ^{-1}(3)+i \log (90)}{16 \pi }\right)}\right)\right.\right.\right. \\ \left.\left.\left.
-(1-2 i) (2+3 i)^{3/4} \log \left(\frac{\Gamma \left(\frac{7 \pi +\tan ^{-1}\left(\frac{12}{5}\right)+i \log (117)}{16 \pi }\right)}{\Gamma
   \left(\frac{15 \pi +\tan ^{-1}\left(\frac{12}{5}\right)+i \log (117)}{16 \pi }\right)}\right)\right)\right)\right)
\end{multline}
\end{example}
\begin{example}
\begin{multline}
\int _0^{\infty }\int _0^{\infty }\frac{y^2 \Gamma (1+x) \log \left(\log \left(-\frac{y}{\sqrt[4]{x}}\right)\right)}{\sqrt[4]{x} \left(1+y^4\right) \Gamma (4+x)}dydx\\
=\frac{1}{24}
   \pi ^2 \left(i \left(3-3\times 2^{3/4}+3^{3/4}\right) \pi -18 \log (2)+3^{3/4} \log (4)\right. \\ \left.
+2 \left(\left(3+3^{3/4}\right) \log (\pi )-3\times 2^{3/4} \log (2 \pi )\right.\right. \\ \left.\left.
+\log \left(\frac{\Gamma
   \left(\frac{1}{4}\right)^6 \left(\frac{\Gamma \left(\frac{3}{4}+\frac{i \log (2)}{8 \pi }\right)}{\Gamma \left(\frac{5}{4}+\frac{i \log (2)}{8 \pi }\right)}\right)^{6\times 2^{3/4}}
   \left(\frac{\Gamma \left(\frac{3}{4}+\frac{i \log (3)}{8 \pi }\right)}{\Gamma \left(\frac{5}{4}+\frac{i \log (3)}{8 \pi }\right)}\right)^{-2\times 3^{3/4}}}{\Gamma
   \left(\frac{3}{4}\right)^6}\right)\right)\right)
\end{multline}
\end{example}
\begin{example}
\begin{multline}
\int _0^{\infty }\int _0^{\infty }\frac{y^2 \Gamma (1+x) \log \left(\log \left(\frac{i y}{\sqrt[4]{x}}\right)\right)}{\sqrt[4]{x} \left(1+y^4\right) \Gamma (4+x)}dydx\\
=\frac{1}{24}
   \pi ^2 \left(i \left(3-3\times 2^{3/4}+3^{3/4}\right) \pi +3^{3/4} \log (4)-\left(-1+2^{3/4}\right) \log (64)\right. \\ \left.
+2 \left(-3\times 2^{3/4}+3^{3/4}\right) \log (\pi )+12\times 2^{3/4} \log
   \left(\frac{\Gamma \left(\frac{1}{2}+\frac{i \log (2)}{8 \pi }\right)}{\Gamma \left(1+\frac{i \log (2)}{8 \pi }\right)}\right)\right. \\ \left.
-4\times 3^{3/4} \log \left(\frac{\Gamma
   \left(\frac{1}{2}+\frac{i \log (3)}{8 \pi }\right)}{\Gamma \left(1+\frac{i \log (3)}{8 \pi }\right)}\right)\right)
\end{multline}
\end{example}
\subsubsection{Formulae derived using equation (\ref{eq_32612})}
\begin{example}
\begin{multline}
\int _0^{\infty }\int _0^{\infty }\int _0^{\infty }\int _0^{\infty }\frac{e^{-4 \left(t+y^4+z^4\right)}
   y^3 z^3 K\left(\frac{1-x}{2}\right) \theta (1-x)}{\sqrt{t} \sqrt{x} \log \left(-\frac{2 x y^2
   z^2}{t}\right)}dtdzdydx=-\frac{i \pi ^{3/2} \log (2)}{1024 \sqrt{2}}
\end{multline}
\end{example}
\begin{example}
\begin{multline}
\int _0^{\infty }\int _0^{\infty }\int _0^{\infty }\int _0^{\infty }\frac{e^{-\pi 
   \left(t+y^4+z^4\right)} y^3 z^3 K\left(\frac{1-x}{2}\right) \theta (1-x)}{\sqrt{t} \sqrt{x} \log
   ^2\left(-\frac{2 x y^2 z^2}{t}\right)}dtdzdydx=-\frac{1}{768 \sqrt{2}}
\end{multline}
\end{example}
\begin{example}
\begin{multline}
\int _0^{\infty }\int _0^{\infty }\int _0^{\infty }\int _0^{\infty }\frac{e^{-t-y^4-z^4} y^3 z^3
   K\left(\frac{1-x}{2}\right) \theta (1-x)}{\sqrt{t} \sqrt{x} \log ^3\left(-\frac{2 x y^2
   z^2}{t}\right)}dtdzdydx=\frac{3 i \zeta (3)}{512 \sqrt{2 \pi }}
\end{multline}
\end{example}
\begin{example}
\begin{multline}
\int _0^{\infty }\int _0^{\infty }\int _0^{\infty }\int _0^{\infty }\frac{e^{-t-y^4-z^4} y^3 z^3
   K\left(\frac{1-x}{2}\right) \theta (1-x) \log \left(-\frac{2 x y^2 z^2}{t}\right) \log \left(\log \left(-\frac{2
   x y^2 z^2}{t}\right)\right)}{\sqrt{t} \sqrt{x}}dtdzdydx\\
=-\frac{\pi ^{9/2}}{64 \sqrt{2}}+\frac{i \pi ^{7/2} \log
   (2)}{32 \sqrt{2}}-\frac{i \pi ^{7/2} \log \left(\frac{A^3}{\sqrt[3]{2} \sqrt[4]{e}}\right)}{8 \sqrt{2}}+\frac{i
   \pi ^{7/2} \log (\pi )}{32 \sqrt{2}}
\end{multline}
\end{example}
\begin{example}
\begin{multline}
\int _0^{\infty }\int _0^{\infty }\int _0^{\infty }\int _0^{\infty }\frac{e^{-t-y^4-z^4} y^3 z^3
   K\left(\frac{1-x}{2}\right) \theta (1-x) \log ^2\left(-\frac{2 x y^2 z^2}{t}\right) \log \left(\log
   \left(-\frac{2 x y^2 z^2}{t}\right)\right)}{\sqrt{t} \sqrt{x}}dtdzdydx\\
=\frac{7 \pi ^{5/2} \zeta (3)}{16
   \sqrt{2}}
\end{multline}
\end{example}
\begin{example}
\begin{multline}
\int _0^{\infty }\int _0^{\infty }\int _0^{\infty }\int _0^{\infty }\frac{e^{-t-y^4-z^4} y^3 z^3
   K\left(\frac{1-x}{2}\right) \theta (1-x) \log \left(\frac{2 x y^2 z^2}{t}\right) \log \left(\log \left(\frac{2 x
   y^2 z^2}{t}\right)\right)}{\sqrt{t} \sqrt{x}}dtdzdydx\\
=-\frac{i C \pi ^{5/2}}{8 \sqrt{2}}
\end{multline}
\end{example}
\subsubsection{Formulae derived using equation (\ref{eq_2273})}
\begin{example}
\begin{multline}
\int _{-1}^1\int _{-1}^1\frac{\left(\frac{1+y}{1-y}\right)^{3/4}}{\left(-1+\frac{2}{1+x}\right)^{5/8} \left(1+y^2\right) \log \left(\frac{\sqrt{-1+\frac{2}{1+x}}
   (1+y)}{-1+y}\right)}dydx\\
=-\frac{1}{96} \sqrt[4]{-1} (48 C+\pi  (29 i \pi +60 \log (2)))
\end{multline}
\end{example}
\begin{example}
\begin{multline}
\int _{-1}^1\int _{-1}^1\frac{\sqrt{\frac{y+1}{1-y}}}{\left(\frac{2}{x+1}-1\right)^{3/4} \left(y^2+1\right) \log \left(\frac{\sqrt{\frac{2}{x+1}-1}
   (y+1)}{y-1}\right)}dydx
=-\frac{1}{24} \pi  (\pi +36 i \log (2))
\end{multline}
\end{example}
\begin{example}
\begin{multline}
\int _{-1}^1\int _{-1}^1\frac{\sqrt{\frac{y+1}{1-y}} \log ^k\left(\frac{\sqrt{\frac{2}{x+1}-1} (y+1)}{y-1}\right)}{\left(\frac{2}{x+1}-1\right)^{3/4}
   \left(y^2+1\right)}dydx\\
=2^k e^{\frac{1}{2} i \pi  (k+1)} \pi ^{k+1} \left(\left(1-2^k\right) k \zeta (1-k)+3 i \pi  \left(2^{k+1}-1\right) \zeta
   (-k)\right)
\end{multline}
\end{example}
\begin{example}
\begin{multline}
\int _{-1}^1\int _{-1}^1\frac{\sqrt{\frac{1+y}{1-y}} \log \left(\log \left(\frac{\sqrt{-1+\frac{2}{1+x}} (1+y)}{-1+y}\right)\right)}{\left(-1+\frac{2}{1+x}\right)^{3/4}
   \left(1+y^2\right) \log \left(\frac{\sqrt{-1+\frac{2}{1+x}} (1+y)}{-1+y}\right)}dydx\\
=\frac{1}{48} \pi  \left(-i \pi ^2+2 \pi  (1+\gamma +19 \log (2)-12 \log (A))+36 i \log (2)
   \left(2 \gamma +\log \left(\frac{1}{8 \pi ^2}\right)\right)\right)
\end{multline}
\end{example}
\begin{example}
\begin{multline}
\int _{-1}^1\int _{-1}^1\frac{(1+x) \left(\sqrt[4]{-1+\frac{2}{1+x}} \sqrt{\frac{1+y}{1-y}}-\sqrt[3]{-1+\frac{2}{1+x}} \left(\frac{1+y}{1-y}\right)^{2/3}\right)}{(-1+x)
   \left(1+y^2\right) \log \left(\frac{\sqrt{-1+\frac{2}{1+x}} (1+y)}{1-y}\right)}dydx\\
=2 C+\frac{3 i \pi ^2}{4}-\frac{8}{3} \pi  \tanh ^{-1}\left((-1)^{2/3}\right)-\frac{1}{2}
   \sqrt[6]{-1} \Phi \left(-\sqrt[3]{-1},2,\frac{1}{2}\right)
\end{multline}
\end{example}
\begin{example}
\begin{multline}
\int _{-1}^1\int _{-1}^1\frac{\sqrt{\frac{1+y}{1-y}} \left(\log \left(-1+\frac{2}{1+x}\right)+2 \log \left(\frac{1+y}{1-y}\right)\right)}{2 \left(-1+\frac{2}{1+x}\right)^{3/4}
   \left(1+y^2\right) \log \left(\frac{\sqrt{-1+\frac{2}{1+x}} (1+y)}{-1+y}\right)}dydx\\
   =\frac{1}{24} \pi ^2 (36+i \pi -36 \log (2))
\end{multline}
\end{example}
\begin{example}
\begin{multline}
\int _{-1}^1\int _{-1}^1\frac{\sqrt{\frac{1+y}{1-y}}}{\left(-1+\frac{2}{1+x}\right)^{3/4} \left(1+y^2\right) \log ^2\left(\frac{\sqrt{-1+\frac{2}{1+x}}
   (1+y)}{-1+y}\right)}dydx=-\frac{\pi ^3-6 i \zeta (3)}{16 \pi }
\end{multline}
\end{example}
\begin{example}
\begin{multline}
\int _{-1}^1\int _{-1}^1-\frac{k \sqrt{\frac{1+y}{1-y}} \log ^{-1+k}\left(\frac{\sqrt{-1+\frac{2}{1+x}} (1+y)}{-1+y}\right)}{\left(-1+\frac{2}{1+x}\right)^{3/4}
   \left(1+y^2\right)}dydx\\
=e^{\frac{i k \pi }{2}} \pi  \Gamma (1+k) \left(\left(-2+2^k\right) \sin \left(\frac{k \pi }{2}\right) \zeta (-1+k)-3 i \left(-1+2^k\right) \cos
   \left(\frac{k \pi }{2}\right) \zeta (k)\right)
\end{multline}
\end{example}
\begin{example}
\begin{multline}
\int _{-1}^1\int _{-1}^1\frac{\sqrt{\frac{1+y}{1-y}}}{\left(-1+\frac{2}{1+x}\right)^{3/4} \left(1+y^2\right) \sqrt{\log \left(\frac{\sqrt{-1+\frac{2}{1+x}}
   (1+y)}{-1+y}\right)}}dydx\\
=\left(\frac{1}{2}+\frac{i}{2}\right) \pi ^{3/2} \left(-\left(-2+\sqrt{2}\right) \zeta \left(-\frac{1}{2}\right)+3 i \left(-1+\sqrt{2}\right) \zeta
   \left(\frac{1}{2}\right)\right)
\end{multline}
\end{example}
\begin{example}
\begin{multline}
\int _{-1}^1\int _{-1}^1\frac{\left(\frac{1+y}{1-y}\right)^{3/4}}{\left(-1+\frac{2}{1+x}\right)^{5/8} \left(1+y^2\right) \log ^2\left(\frac{\sqrt{-1+\frac{2}{1+x}}
   (1+y)}{-1+y}\right)}dydx\\
=\frac{\sqrt[4]{-1} \left(-240 C \pi +11 i \pi ^3+18 \zeta (3)\right)}{384 \pi }
\end{multline}
\end{example}
\subsubsection{Derived from equation (\ref{eq_341131})}
\begin{example}
\begin{equation}
\int_0^{\infty } \frac{\left(e^{m x}+(-1)^{1+k} e^{(-m+p) x}\right) x^k}{-1+e^{p x}} \, dx=\frac{i \pi  \left(2^{1+k} \left(\frac{i}{p}\right)^k \pi ^k \text{Li}_{-k}\left(e^{\frac{2 i m \pi
   }{p}}\right)\right)}{p}
\end{equation}
\end{example}
\begin{example}
\begin{equation}
\int_0^{\infty } \frac{x^2 \coth (x) \text{csch}(x)}{\pi ^2+x^2} \, dx=\frac{1}{12} \left(-6+\pi ^2\right)
\end{equation}
\end{example}
\begin{example}
\begin{equation}
\int_0^{\infty } \frac{x \text{csch}(x)}{\pi ^2+4 x^2} \, dx=\frac{1}{8} (-2+\pi )
\end{equation}
\end{example}
\begin{example}
\begin{equation}
\int_0^{\infty } \frac{x^2 \coth (x) \text{csch}(x)}{\pi ^2+4 x^2} \, dx=\frac{1}{4} (-1+2 C)
\end{equation}
\end{example}
\begin{example}
\begin{equation}
\int_0^{\infty } \tanh ^{-1}\left(\frac{x}{a}\right) \text{csch}(m x) \, dx=-\frac{i \pi  \log \left(-\frac{i a m \Gamma \left(-\frac{i a m}{2 \pi }\right)^2}{2 \pi  \Gamma \left(\frac{-i a
   m+\pi }{2 \pi }\right)^2}\right)}{2 m}
\end{equation}
\end{example}
\begin{example}
\begin{multline}
\int_0^{\infty } x^k \text{csch}(m x) \, dx=\left(-1+2^{1+k}\right) \left(\frac{i}{m}\right)^{2+k} m \pi ^{1+k} \left(-i+\cot \left(\frac{k \pi }{2}\right)\right) \zeta (-k)
\end{multline}
\end{example}
\begin{example}
\begin{equation}
\int_0^{\infty } \frac{x^2 \coth (x) \text{csch}(x)}{a^2+x^2} \, dx=\frac{\pi }{2 a}+\frac{a \left(\psi ^{(1)}\left(1+\frac{a}{2 \pi }\right)-\psi ^{(1)}\left(\frac{a+\pi }{2 \pi
   }\right)\right)}{4 \pi }
\end{equation}
\end{example}
\begin{example}
\begin{multline}
\int_0^{\infty } \left(-\left(\frac{i \pi }{m}-x\right)^k+\left(\frac{i \pi }{m}+x\right)^k\right) \text{csch}(m x) \, dx\\
=\frac{\left(\frac{i}{m}\right)^{-1+k} \pi ^{1+k}
   \left(1+\left(-2+2^{2+k}\right) \zeta (-k)\right)}{m^2}
\end{multline}
\end{example}
\begin{example}
\begin{equation}
\int_0^{\infty } \frac{x \text{csch}(m x)}{\pi ^2+m^2 x^2} \, dx=\frac{-1+\log (4)}{2 m^2}
\end{equation}
\end{example}
\begin{example}
\begin{equation}
\int_0^{\infty } \frac{x \text{csch}(m x)}{\left(\pi ^2+m^2 x^2\right)^2} \, dx=\frac{-6+\pi ^2}{24 m^2 \pi ^2}
\end{equation}
\end{example}
\begin{example}
\begin{multline}
\int_0^{\infty } \frac{-2 i \pi  x \text{csch}(x)+\pi ^2 \text{sech}(x)-x^2 \text{sech}(x)}{\left(\pi ^2+x^2\right)^2} \, dx=\frac{(48+6 i)-48 C-i \pi ^2}{12 \pi }
\end{multline}
\end{example}
\begin{example}
\begin{equation}
\int_0^{\infty } \frac{(\pi -x) (\pi +x) \text{sech}(x)}{\left(\pi ^2+x^2\right)^2} \, dx=\frac{4-4 C}{\pi }
\end{equation}
\end{example}
\begin{example}
\begin{multline}
\int_0^{\infty } \text{csch}(m x) \left(\frac{\log (-x+\log (a))}{x-\log (a)}+\frac{\log (x+\log (a))}{x+\log (a)}\right) \, dx\\
=\left(-H_{-\frac{i m \log (a)}{2 \pi }}+H_{-\frac{\pi +i m \log
   (a)}{2 \pi }}\right) \log \left(\frac{2 i \pi }{m}\right)+\frac{i \pi  \log (\log (a))}{m \log (a)}\\
   -\gamma _1\left(\frac{\pi -i m \log (a)}{2 \pi }\right)+\gamma _1\left(1-\frac{i m \log (a)}{2 \pi
   }\right)
\end{multline}
\end{example}
\begin{example}
\begin{multline}
\int_0^{\infty } \left(\frac{i \pi  x \text{csch}(x)}{\pi ^2+x^2}-\frac{2 \pi  \tan ^{-1}\left(\frac{x}{\pi }\right) \text{csch}(x)}{\pi ^2+x^2}+\frac{x \text{csch}(x) \log \left(\pi
   ^2+x^2\right)}{\pi ^2+x^2}\right) \, dx\\
   =\log (i \pi )-H_{\frac{1}{2}} \log (2 i \pi )-\gamma _1+\gamma _1\left(\frac{3}{2}\right)
\end{multline}
\end{example}
\begin{example}
\begin{equation}
\int_0^{\infty } \frac{x \text{csch}(x)}{\pi ^2+x^2} \, dx=-\frac{1}{2}+\log (2)
\end{equation}
\end{example}
\begin{example}
\begin{multline}
\int_0^{\infty } \frac{\text{csch}(x) \left(-2 \pi  \tan ^{-1}\left(\frac{x}{\pi }\right)+x \log \left(\pi ^2+x^2\right)\right)}{\pi ^2+x^2} \, dx\\
=\log (\pi )+(-2+\log (4)) \log (2 \pi )-\gamma
   _1+\gamma _1\left(\frac{3}{2}\right)
\end{multline}
\end{example}
\begin{example}
\begin{multline}
\int_0^{\infty } \frac{-e^{3 m x} \log (-x+\log (a))+e^{m x} \log (x+\log (a))}{-1+e^{4 m x}} \, dx\\
=\frac{\pi }{4m} \left((-1-i) \log \left(\frac{i}{m}\right)+(1+i) (\log (2)-\log (\pi ))+i \log
   (\log (a))\right. \\ \left.
   +2 \left(\log \left(\frac{\Gamma \left(\frac{1}{4} \left(1-\frac{2 i m \log (a)}{\pi }\right)\right)}{2 \Gamma \left(\frac{1}{4} \left(3-\frac{2 i m \log (a)}{\pi }\right)\right)}\right)+i
   \log \left(\frac{\Gamma \left(\frac{1}{4} \left(2-\frac{2 i m \log (a)}{\pi }\right)\right)}{2 \Gamma \left(\frac{1}{4} \left(4-\frac{2 i m \log (a)}{\pi }\right)\right)}\right)\right)\right)
\end{multline}
where $Re(a)<0$
\end{example}
\begin{example}
\begin{multline}
\int_0^{\infty } \frac{-e^{(-m+p) x} (2 a \pi -x)^k+e^{m x} (2 a \pi +x)^k}{-1+e^{p x}} \, dx\\
=\frac{i \pi  \left(a^k (2 \pi )^k+2^{1+k} e^{\frac{2 i m \pi }{p}} \left(\frac{i}{p}\right)^k \pi ^k
   \Phi \left(e^{\frac{2 i m \pi }{p}},-k,1-i a p\right)\right)}{p}
\end{multline}
\end{example}
\begin{example}
\begin{multline}
\int_0^{\infty } \frac{-e^{x/4} (2 i \pi -x) \log (2 i \pi -x)+e^{x/4} (2 i \pi +x) \log (2 i \pi +x)}{-1+e^{x/2}} \, dx\\
=2 i \pi  \left(-\pi ^2-4 i \pi  \log (2)\right. \\ \left.
-4 i \pi  \left(\frac{i \pi
   }{2}+\log (2)\right)-4 i \pi  \log (\pi )+2 i \pi  \log (2 \pi )+8 i \pi  \Phi'\left(-1,-1,\frac{3}{2}\right)\right)
\end{multline}
\end{example}
\begin{example}
\begin{multline}
\int_0^{\infty } \frac{e^{x/2} ((-i \pi +x) \log (i \pi -x)+(i \pi +x) \log (i \pi +x))}{-1+e^x} \, dx\\
=\frac{1}{2} \pi ^2 \left(i \pi +\log \left(16 \pi ^2\right)-8
   \Phi'\left(-1,-1,\frac{3}{2}\right)\right)
\end{multline}
\end{example}
\begin{example}
\begin{equation}
\int_0^{\infty } \frac{x (-1+x \coth (x)) \text{csch}(x)}{\left(\pi ^2+x^2\right)^3} \, dx=\frac{7}{960}-\frac{\pi ^2+18 \zeta (3)}{48 \pi ^4}
\end{equation}
\end{example}
\begin{example}
\begin{equation}
\int_0^{\infty } \frac{x (-1+x \coth (x)) \text{csch}(x)}{\left(\pi ^2+4 x^2\right)^3} \, dx=-\frac{512 C+32 \pi ^3-\psi ^{(3)}\left(\frac{1}{4}\right)+\psi
   ^{(3)}\left(\frac{3}{4}\right)}{4096 \pi ^4}
\end{equation}
\end{example}
\begin{example}
\begin{multline}
\int_0^{\infty } \frac{x \text{csch}(m x)}{\left(a^2+x^2\right)^3} \, dx\\
=\frac{3 \pi }{16 a^5 m}+\frac{m \left(2 \pi  \psi ^{(1)}\left(1+\frac{a m}{2 \pi }\right)-2 \pi  \psi
   ^{(1)}\left(\frac{a m+\pi }{2 \pi }\right)+a m \left(-\psi ^{(2)}\left(1+\frac{a m}{2 \pi }\right)+\psi ^{(2)}\left(\frac{a m+\pi }{2 \pi }\right)\right)\right)}{64 a^3 \pi ^2}
\end{multline}
\end{example}
\begin{example}
\begin{multline}
\int_0^{\infty } x \tan ^{-1}\left(\frac{2 x}{\pi }\right) \coth \left(\frac{2 x}{3}\right) \text{csch}\left(\frac{2 x}{3}\right) \, dx=\frac{1}{8} (-3) \pi  \log \left(\frac{3456
   e^{3-\frac{2 \pi }{\sqrt{3}}} \Gamma \left(-\frac{1}{3}\right)^6}{15625 \Gamma \left(-\frac{5}{6}\right)^6}\right)
\end{multline}
\end{example}
\begin{example}
\begin{equation}
\int_0^{\infty } \frac{x^2 \text{csch}^3\left(\frac{x}{2}\right) \sinh (x)}{\pi ^2+x^2} \, dx=-4+8 C
\end{equation}
\end{example}
\begin{example}
\begin{equation}
\int_0^{\infty } \frac{x^2 \text{csch}^3\left(\frac{x}{2}\right) \sinh (x)}{\pi ^2+4 x^2} \, dx=\frac{1}{16} \left(-32+\psi ^{(1)}\left(\frac{1}{8}\right)-\psi
   ^{(1)}\left(\frac{5}{8}\right)\right)
\end{equation}
\end{example}
\subsubsection{Formulae derived from equation (\ref{eq_32511})}
\begin{example}
\begin{multline}
\int _0^1\int _0^{\infty }\frac{x^n \left(-1+\sqrt{-1+x^{-1-n}} y\right)}{\sqrt[4]{-1+x^{-1-n}} \sqrt{y} \left(1+y^2\right) \log \left(\left(-1+x^{-1-n}\right)
   y^2\right)}dydx=\frac{C}{1+n}
\end{multline}
\end{example}
\begin{example}
\begin{multline}
\int _0^1\int _0^{\infty }\frac{x^n \sqrt[4]{-1+x^{-1-n}} \sqrt{y}}{\left(1+y^2\right) \log \left(\left(1-x^{-1-n}\right) y^2\right)}dydx\\
=\frac{2 i \pi  \left(\psi
   ^{(0)}\left(\frac{3}{8}\right)-\psi ^{(0)}\left(\frac{7}{8}\right)\right)+\psi ^{(1)}\left(\frac{3}{8}\right)-\psi ^{(1)}\left(\frac{7}{8}\right)}{32 (1+n)}
\end{multline}
\end{example}
\begin{example}
\begin{multline}
\int _0^1\int _0^{\infty }\frac{x^n \sqrt[4]{-1+x^{-1-n}} \sqrt{y} \log \left(\log \left(a \left(-1+x^{-1-n}\right) y^2\right)\right)}{1+y^2}dydx\\
=\frac{\pi }{8 (1+n)} \left(\pi  (i \pi +\log
   (64)+2 \log (\pi )-4 \log (-6 \pi -i \log (a))+4 \log (-2 \pi -i \log (a)))\right. \\ \left.
   -4 \pi  \text{log$\Gamma $}\left(-\frac{3}{4}-\frac{i \log (a)}{8 \pi }\right)+4 \pi  \text{log$\Gamma
   $}\left(-\frac{1}{4}-\frac{i \log (a)}{8 \pi }\right)+2 i \psi ^{(0)}\left(\frac{1}{4}-\frac{i \log (a)}{8 \pi }\right)\right. \\ \left.
   -2 i \psi ^{(0)}\left(\frac{3}{4}-\frac{i \log (a)}{8 \pi }\right)\right)
\end{multline}
\end{example}
\begin{example}
\begin{multline}
\int _0^1\int _0^{\infty }\frac{\sqrt[4]{-1+\frac{1}{x^2}} x \sqrt{y} \log \left(\log \left(\left(1-\frac{1}{x^2}\right) y^2\right)\right)}{1+y^2}dydx\\
=\frac{1}{16} \pi  \left(\pi 
   \left(i \pi +\log \left(\frac{64 \pi ^2 \Gamma \left(-\frac{1}{8}\right)^4}{625 \Gamma \left(-\frac{5}{8}\right)^4}\right)\right)+2 i \left(\psi ^{(0)}\left(\frac{3}{8}\right)-\psi
   ^{(0)}\left(\frac{7}{8}\right)\right)\right)
\end{multline}
\end{example}
\begin{example}
\begin{multline}
\int _0^1\int _0^{\infty }\frac{\sqrt[4]{-1+\frac{1}{x^2}} x \sqrt{y} \log \left(\left(-1+\frac{1}{x^2}\right) y^2\right) \log \left(\log \left(\left(-1+\frac{1}{x^2}\right)
   y^2\right)\right)}{1+y^2}dydx\\
   =\frac{1}{4} \pi ^2 \left(2-4 i C+\log \left(\exp (\pi  i) \left(\frac{8 \pi }{9}\right)^2 \left(\frac{\Gamma \left(-\frac{1}{4}\right)}{\Gamma
   \left(-\frac{3}{4}\right)}\right)^4\right)\right)
\end{multline}
\end{example}
\begin{example}
\begin{multline}
\int _0^1\int _0^{\infty }\frac{\sqrt[4]{-1+\frac{1}{x}} \sqrt{y}}{\left(1+y^2\right) \log \left(\frac{(-1+x) y^2}{x}\right)}dydx\\
=\frac{1}{32} \left(-2 i \sqrt{2} \pi  \left(\pi +2
   \log \left(\tan \left(\frac{\pi }{8}\right)\right)\right)+\psi ^{(1)}\left(\frac{3}{8}\right)-\psi ^{(1)}\left(\frac{7}{8}\right)\right)
\end{multline}
\end{example}
\begin{example}
\begin{multline}
\int _0^1\int _0^{\infty }\frac{\sqrt[4]{-1+\frac{1}{x}} \sqrt{y}}{\left(1+y^2\right) \sqrt{\log \left(\frac{(-1+x) y^2}{x}\right)}}dydx\\
=\left(\frac{1}{32}+\frac{i}{32}\right)
   \sqrt{\pi } \left(-4 i \pi  \left(\zeta \left(\frac{1}{2},\frac{3}{8}\right)-\zeta \left(\frac{1}{2},\frac{7}{8}\right)\right)+\zeta \left(\frac{3}{2},\frac{3}{8}\right)-\zeta
   \left(\frac{3}{2},\frac{7}{8}\right)\right)
\end{multline}
\end{example}
\begin{example}
\begin{multline}
\int _0^1\int _0^{\infty }\frac{\sqrt[4]{-1+\frac{1}{x}} \sqrt{y} \log ^i\left(\frac{(-1+x) y^2}{x}\right)}{1+y^2}dydx\\
=i^i 2^{-2+3 i} \pi ^{1+i} \left(2 \pi  \left(\zeta
   \left(-i,\frac{3}{8}\right)-\zeta \left(-i,\frac{7}{8}\right)\right)+\zeta \left(1-i,\frac{3}{8}\right)-\zeta \left(1-i,\frac{7}{8}\right)\right)
\end{multline}
\end{example}
\begin{example}
\begin{multline}
\int _0^1\int _0^{\infty }\frac{\sqrt[4]{-1+\frac{1}{x}} \sqrt{y} \left(\log \left(-1+\frac{1}{x}\right)+2 \log (y)\right)}{\left(1+y^2\right) \log \left(i \left(-1+\frac{1}{x}\right)
   y^2\right)}dydx\\
   =\frac{1}{64} \pi  \left(-i \left(\psi ^{(1)}\left(\frac{5}{16}\right)-\psi ^{(1)}\left(\frac{13}{16}\right)\right)+4 \pi  \left(4+2 \cos \left(\frac{\pi }{8}\right) \log
   \left(\tan \left(\frac{3 \pi }{16}\right)\right)\right.\right. \\ \left.\left.
   -\left(\sqrt{2} \pi +2 \log \left(\tan \left(\frac{\pi }{16}\right)\right)\right) \sin \left(\frac{\pi }{8}\right)\right)\right)
\end{multline}
\end{example}
\begin{example}
\begin{multline}
\int _0^1\int _0^{\infty }\frac{\left(1-\sqrt{-1+\frac{1}{x}} y\right) \log ^2\left(\left(-1+\frac{1}{x}\right) y^2\right) \log \left(\log \left(\left(-1+\frac{1}{x}\right)
   y^2\right)\right)}{\sqrt[4]{-1+\frac{1}{x}} \sqrt{y} \left(1+y^2\right)}dydx=32 i C \pi ^2
\end{multline}
\end{example}
\begin{example}
\begin{equation}
\int _0^1\int _0^{\infty }\frac{i \log ^3\left(\left(-1+\frac{1}{x}\right) y^2\right) \log \left(\log \left(\left(-1+\frac{1}{x}\right) y^2\right)\right)}{18 \pi ^2
   \left(1+y^2\right)}dydx=\zeta (3)
\end{equation}
\end{example}
\begin{example}
\begin{multline}
\int _0^1\int _0^{\infty }\frac{15 \log ^4\left(\left(-1+\frac{1}{x}\right) y^2\right) \log \left(\log \left(\left(-1+\frac{1}{x}\right) y^2\right)\right)}{2 \pi ^5
   \left(1+y^2\right)}dydx\\
   =7+14 i \pi +60 \log (2)+28 \log (\pi )-3360 \zeta '(-3)
\end{multline}
\end{example}
\subsubsection{Formulae derived using equation (\ref{eq_3139})}
\begin{example}
\begin{multline}\label{eq1}
\int_{\mathbb{R}^{4}_+}\frac{\sqrt[4]{r+s} e^{-3 r-2 s-3 x-2 y} \left(\pi ^2-4 \log ^2\left(\frac{\sqrt{r+s} \sqrt{x y}}{\sqrt{r s} \sqrt{x+y}}\right)\right)}{(r s)^{3/4}
   \sqrt[4]{x y} \sqrt[4]{x+y} \left(4 \log ^2\left(\frac{\sqrt{r+s} \sqrt{x y}}{\sqrt{r s} \sqrt{x+y}}\right)+\pi ^2\right)^2}dxdydrds
   =\frac{48 C+\pi ^2}{576
   \sqrt{2}}
\end{multline}
\textup{and}
\begin{multline}
\int_{\mathbb{R}^{4}_+}\frac{\sqrt[4]{r s} \sqrt[4]{r+s} (x y)^{3/4} e^{-3 r-2 s-3 x-2 y} \log \left(\frac{\sqrt{r+s} \sqrt{x y}}{\sqrt{r s} \sqrt{x+y}}\right)}{r s x y
   \sqrt[4]{x+y} \left(4 \log ^2\left(\frac{\sqrt{r+s} \sqrt{x y}}{\sqrt{r s} \sqrt{x+y}}\right)+\pi ^2\right)^2}dxdydrds\\
   =\frac{1}{16 \pi }\left(\frac{C}{3 \sqrt{2}}-\frac{\pi
   ^2}{144 \sqrt{2}}\right)
\end{multline}
\end{example}
\begin{example}
\begin{equation}\label{eq2}
\int_{\mathbb{R}^{4}_+}\frac{\sqrt{r s} \sqrt{x y} e^{-3 r-2 s-3 x-2 y}}{r s x y \left(\log ^2\left(\frac{\sqrt{r+s} \sqrt{x y}}{\sqrt{r s} \sqrt{x+y}}\right)+\pi
   ^2\right)}dxdydrds=\frac{4-\pi }{6}
\end{equation}
\textup{and}
\begin{equation}
\int_{\mathbb{R}^{4}_+}\frac{\sqrt{r s} \sqrt{x y} e^{-3 r-2 s-3 x-2 y} \log \left(\frac{\sqrt{r+s} \sqrt{x y}}{\sqrt{r s} \sqrt{x+y}}\right)}{r s x y \left(\log
   ^2\left(\frac{\sqrt{r+s} \sqrt{x y}}{\sqrt{r s} \sqrt{x+y}}\right)+\pi ^2\right)}dxdydrds=0
\end{equation}
\end{example}
\begin{example}
\begin{multline}\label{eq3}
\int_{\mathbb{R}^{4}_+}\frac{e^{-2 r-3 s-2 x-3 y} \left((r+s)^{3/8} \sqrt[4]{x y}-\sqrt[4]{r s} \sqrt[8]{r+s} \sqrt[4]{x+y}\right)}{(r s)^{7/8} (x y)^{3/8} (x+y)^{3/8}
   \log \left(\frac{\sqrt{r+s} \sqrt{x y}}{\sqrt{r s} \sqrt{x+y}}\right)}dxdydrds\\
   =\frac{2}{3} \pi  \tanh ^{-1}\left(\cos \left(\frac{\pi }{8}\right)-\sin
   \left(\frac{\pi }{8}\right)\right)
\end{multline}
\end{example}
\begin{example}
\begin{multline}\label{eq4}
\int_{\mathbb{R}^{4}_+}\frac{e^{-r-s-x-y} \left(\sqrt[6]{r+s} \sqrt[24]{x y}-\sqrt[24]{r s} \sqrt[8]{r+s} \sqrt[24]{x+y}\right)}{(r s)^{2/3} (x y)^{3/8} \sqrt[6]{x+y} \log
   \left(\frac{\sqrt{r+s} \sqrt{x y}}{\sqrt{r s} \sqrt{x+y}}\right)}dxdydrds\\
   =2 \pi  \log \left(\sqrt{3} \tan \left(\frac{3 \pi }{16}\right)\right)
\end{multline}
\end{example}
\begin{example}
\begin{equation}\label{eq5}
\int_{\mathbb{R}^{4}_+}\frac{e^{-r-2 (s+y)-x}}{\sqrt{r s} \sqrt{x y} \left(\log \left(\frac{\sqrt{r+s} \sqrt{x y}}{\sqrt{r s} \sqrt{x+y}}\right)+i \pi \right)^2}dxdydrds=4
   (C-1)
\end{equation}
\end{example}
\begin{example}
\begin{equation}\label{eq6}
\int_{\mathbb{R}^{4}_+}\frac{e^{-r-2 (s+y)-x}}{\sqrt{r s} \sqrt{x y} \left(\log \left(\frac{\sqrt{r+s} \sqrt{x y}}{\sqrt{r s} \sqrt{x+y}}\right)+i \pi \right)^3}dxdydrds=-\frac{i
   \left(\pi ^3-32\right)}{4 \pi }
\end{equation}
\end{example}
\begin{proposition}
For all $a,k,p,q\in\mathbb{C}$ then,
\begin{multline}\label{eq7}
\int_{\mathbb{R}^{4}_+}\frac{e^{-p (r+x)-q (s+y)} \log ^k\left(\frac{a \sqrt{r+s} \sqrt{x y}}{\sqrt{r s} \sqrt{x+y}}\right)}{\sqrt{r s} \sqrt{x y}}dxdydrds\\
=\frac{2 i
   e^{\frac{1}{2} i \pi  (k-1)} \pi ^{k+2} \left(2^k \zeta \left(-k,\frac{1}{2} \left(\frac{1}{2}-\frac{i \log (a)}{\pi }\right)\right)-2^k \zeta
   \left(-k,\frac{1}{2} \left(\frac{3}{2}-\frac{i \log (a)}{\pi }\right)\right)\right)}{p q}
\end{multline}
\end{proposition}
\begin{proposition}
For all $k\in\mathbb{C}$ then,
\begin{equation}\label{eq8}
\int_{\mathbb{R}^{4}_+}\frac{e^{-r-s-x-y} \log ^k\left(\frac{i \sqrt{r+s} \sqrt{x y}}{\sqrt{r s} \sqrt{x+y}}\right)}{\sqrt{r s} \sqrt{x y}}dxdydrds=-2 \left(2^{k+1}-1\right)
   e^{\frac{i \pi  k}{2}} \pi ^{k+2} \zeta (-k)
\end{equation}
\end{proposition}
\begin{example}
\begin{multline}\label{eq9}
\int_{\mathbb{R}^{4}_+}\frac{e^{-r-s-x-y} \sqrt{\log \left(\frac{i \sqrt{r+s} \sqrt{x y}}{\sqrt{r s} \sqrt{x+y}}\right)}}{\sqrt{r s} \sqrt{x y}}dxdydrds\\
=-2 \left(2
   \sqrt{2}-1\right) e^{\frac{i \pi }{4}} \pi ^{5/2} \zeta \left(-\frac{1}{2}\right)
\end{multline}
\end{example}
\begin{example}
\begin{equation}\label{eq10}
\int_{\mathbb{R}^{4}_+}\frac{e^{-3 (r+x)-4 (s+y)}}{\sqrt{r s} \sqrt{x y} \log \left(\frac{i \sqrt{r+s} \sqrt{x y}}{\sqrt{r s} \sqrt{x+y}}\right)}dxdydrds=-\frac{1}{6} i \pi  \log
   (2)
\end{equation}
\end{example}
\begin{example}
\begin{equation}\label{eq11}
\int_{\mathbb{R}^{4}_+}\frac{e^{-r-s-x-y} \log \left(\log \left(\frac{i \sqrt{r+s} \sqrt{x y}}{\sqrt{r s} \sqrt{x+y}}\right)\right)}{\sqrt{r s} \sqrt{x y}}dxdydrds=\frac{1}{2} \pi
   ^2 (\log (4)+i \pi )
\end{equation}
\end{example}
\begin{example}
\begin{multline}\label{eq11}
\int_{\mathbb{R}^{4}_+}\frac{e^{-r-s-x-y} \log \left(\log \left(\frac{i \sqrt{r+s} \sqrt{x y}}{\sqrt{r s} \sqrt{x+y}}\right)\right)}{\sqrt{r s} \sqrt{x y} \log
   \left(\frac{i \sqrt{r+s} \sqrt{x y}}{\sqrt{r s} \sqrt{x+y}}\right)}dxdydrds=\pi  \log (2) (2 i \gamma +\pi -i (\log (2)+2 \log (\pi )))
\end{multline}
\end{example}
\begin{example}
\begin{multline}\label{eq11}
\int_{\mathbb{R}^{4}_+}\frac{e^{-r-s-x-y} \log \left(\log \left(\frac{i \sqrt{r+s} \sqrt{x y}}{\sqrt{r s} \sqrt{x+y}}\right)\right)}{\sqrt{r s} \sqrt{x y} \log
   ^2\left(\frac{i \sqrt{r+s} \sqrt{x y}}{\sqrt{r s} \sqrt{x+y}}\right)}dxdydrds=\frac{1}{12} \pi ^2 (-24 \log (A)+2 \gamma -i \pi +\log (16))
\end{multline}
\end{example}
\begin{proposition}
For all $a,p,q\in\mathbb{C}, Re(a)>0$ then,
\begin{multline}\label{eq12}
\int_{\mathbb{R}^{4}_+}\frac{\log \left(\log \left(\frac{a \sqrt{r+s} \sqrt{x y}}{\sqrt{r s} \sqrt{x+y}}\right)\right) e^{-p (r+x)-q
   (s+y)}}{\sqrt{r s} \sqrt{x y}}dxdydrds\\
   =\frac{\pi ^2 \left(4 \log \left(\frac{\sqrt{2 \pi } \Gamma \left(\frac{3}{4}-\frac{i
   \log (a)}{2 \pi }\right)}{\Gamma \left(\frac{\pi -2 i \log (a)}{4 \pi }\right)}\right)+i \pi \right)}{2 p q}
\end{multline}
\end{proposition}
\subsubsection{Formulae derived using equation (\ref{eq_32414})}
\begin{example}
The Polylogarithm Function $Li_{k}(z)$,
\begin{multline}\label{eq:poly1}
\int_{{\mathbb{R}^{4}_{+}}}y^m x^{u-1} D_u(t) e^{-\frac{t^2}{4}-2 \left(y+z^2\right)} \left(q
   x^v+1\right)^{-m-1} t^{-m+\frac{u}{v}-1} z^{u \left(\frac{1}{v}-1\right)-m}\\
   \log ^k\left(-\frac{y}{t z \left(q x^v+1\right)}\right)dxdydzdt\\\\
   =\frac{1}{v}i^{k+1} \pi
   ^{k+\frac{3}{2}} e^{-i \pi  \left(m-\frac{u}{v}\right)}
   \left(\frac{1}{q}\right)^{u/v} \Gamma \left(\frac{u}{v}\right)
   2^{\frac{1}{2} \left(2 k+m-\frac{u}{v}+u\right)-\frac{(m+3) v+u (-v)+u}{2
   v}}
    \text{Li}_{-k}\left(e^{2 i \pi 
   \left(m-\frac{u}{v}\right)}\right)
\end{multline}
\end{example}
\begin{example}
Catalan's Constant $K$,
\begin{dmath}
\int_{{\mathbb{R}^{4}_{+}}}\frac{D_1(t) e^{-\frac{t^2}{4}-2 \left(y+z^2\right)}}{\sqrt[4]{t}
   \sqrt{x^4+1} \sqrt{y} \sqrt[4]{z} \log ^2\left(-\frac{y}{t
   \left(x^4+1\right) z}\right)}dxdydzdt
   =-\frac{i e^{\frac{3 i \pi }{4}}
   \left(-\frac{\pi ^2}{48}+i K\right) \Gamma \left(\frac{1}{4}\right)}{16\
   2^{3/4} \sqrt{\pi }}
\end{dmath}
\end{example}
\begin{example}
The Hurwitz zeta function $\zeta(s,v)$,
\begin{multline}\label{eq:hurwitz}
\int_{{\mathbb{R}^{4}_{+}}}\frac{D_{\frac{1}{2}}(t \alpha ) e^{-b \left(y+z^2\right)-\frac{1}{4}
   \alpha ^2 t^2} \log ^k\left(\frac{a y}{t z \left(p+q
   x^3\right)}\right)}{\sqrt{t} \sqrt{x} \sqrt[3]{y} \left(p+q
   x^3\right)^{2/3}}dxdydzdt\\
   =-\frac{1}{3 \sqrt{\alpha } b^{7/6} p^{2/3}}i i^{k+1} 2^k \pi ^{k+\frac{3}{2}} \Gamma
   \left(\frac{1}{6}\right) \sqrt[6]{\frac{p}{q}}\\
    \left(2^k \zeta
   \left(-k,-\frac{i (2 \log (a)-\log (b)-2 \log (p)+2 \log (\alpha )+\log (2)+2
   i \pi )}{8 \pi }\right)\right.\\
   \left.-2^k \zeta \left(-k,\frac{1}{2} \left(1-\frac{i (2
   \log (a)-\log (b)-2 \log (p)+2 \log (\alpha )+\log (2)+2 i \pi )}{4 \pi
   }\right)\right)\right)
\end{multline}
\end{example}
\begin{example}
The Riemann zeta function $\zeta(s)$,
\begin{multline}\label{eq:zeta}
\int_{{\mathbb{R}^{4}_{+}}}\frac{D_{\frac{1}{2}}(t) e^{-\frac{t^2}{4}-2 \left(y+z^2\right)} \log
   ^k\left(-\frac{y}{t z \left(q x^3+1\right)}\right)}{\sqrt{t} \sqrt{x}
   \sqrt[3]{y} \left(q x^3+1\right)^{2/3}}dxdydzdt\\\\
   =-\frac{1}{3} i i^{k+1}
   2^{k-\frac{7}{6}} \pi ^{k+\frac{3}{2}} \sqrt[6]{\frac{1}{q}} \Gamma
   \left(\frac{1}{6}\right) \left(2^k \left(2^{-k}-1\right) \zeta (-k)-2^k \zeta
   (-k)\right)
\end{multline}
\end{example}
\begin{example}
Euler's constant $\gamma$
\begin{multline}
\int_{{\mathbb{R}^{4}_{+}}}\frac{D_{\frac{1}{2}}(t) e^{-\frac{t^2}{4}-2 \left(y+z^2\right)} \log \left(\log \left(-\frac{y}{t x^3 z+t
   z}\right)\right)}{\sqrt{t} \sqrt{x} \left(x^3+1\right)^{2/3} \sqrt[3]{y} \log \left(-\frac{y}{t x^3 z+t
   z}\right)}dxdydzdt\\\\
   =\frac{\sqrt{\pi } \log (2) \left(2 i \gamma +\pi -i \log \left(8 \pi ^2\right)\right) \Gamma
   \left(\frac{1}{6}\right)}{24 \sqrt[6]{2}}
\end{multline}
\end{example}
\begin{example}
Glaisher's constant $A$
\begin{multline}
\int_{{\mathbb{R}^{4}_{+}}}\frac{D_{\frac{1}{2}}(t) e^{-\frac{t^2}{4}-2 \left(y+z^2\right)} \log \left(\log \left(-\frac{y}{t x^3 z+t
   z}\right)\right)}{\sqrt{t} \sqrt{x} \left(x^3+1\right)^{2/3} \sqrt[3]{y} \log ^2\left(-\frac{y}{t x^3 z+t
   z}\right)}dxdydzdt\\\\
   =\frac{\pi ^{3/2} \Gamma \left(\frac{1}{6}\right) (-24 \log (A)+2 \gamma -i \pi +\log (4))}{576
   \sqrt[6]{2}}
\end{multline}
\end{example}
\begin{example}
The Constant $\log(2)$,
\begin{dmath}
\int_{{\mathbb{R}^{4}_{+}}}\frac{D_{\frac{1}{2}}(t) e^{\frac{1}{4} \left(-t^2-8
   \left(y+z^2\right)\right)}}{\sqrt{t} \sqrt{x} \left(x^3+1\right)^{2/3}
   \sqrt[3]{y} \log \left(-\frac{y}{t x^3 z+t z}\right)}dxdydzdt
   =-\frac{i \sqrt{\pi }
   \log (2) \Gamma \left(\frac{1}{6}\right)}{12 \sqrt[6]{2}}
\end{dmath}
\end{example}
\begin{example}
The Inverse Hyperbolic Tangent $\tanh^{-1}(z)$,
\begin{dmath}
\int_{{\mathbb{R}^{4}_{+}}}\frac{y^m x^{u-1} t^{-m+\frac{u}{v}-1} z^{u \left(\frac{1}{v}-1\right)-m}
   D_u(t \alpha ) \left(p+q x^v\right)^{-m-1} e^{-\frac{1}{4} \alpha ^2
   \left(\frac{8 \left(y+z^2\right)}{p^2}+t^2\right)}}{\log \left(-\frac{y}{p t
   z+q t z x^v}\right)}dxdydzdt\\\\
   =\frac{\sqrt{\pi } p^{-m-1}
   2^{-\frac{u}{v}+u-\frac{3}{2}} \Gamma \left(\frac{u}{v}\right) \alpha
   ^{m-\frac{u}{v}} \tanh ^{-1}\left(e^{i \pi 
   \left(m-\frac{u}{v}\right)}\right) \left(\frac{p}{q}\right)^{u/v}
   \left(\frac{\alpha ^2}{p^2}\right)^{-\frac{(m+3) v+u (-v)+u}{2
   v}}}{v}
\end{dmath}
\end{example}
\begin{example}
\begin{dmath}
\int_{{\mathbb{R}^{4}_{+}}}\frac{\alpha  t^{3/4} e^{-b \left(y+z^2\right)-\frac{1}{2} \alpha ^2 t^2} \log \left(\log \left(\frac{a y}{p t
   z+q t x^4 z}\right)\right)}{\sqrt{y} \sqrt[4]{z} \sqrt{p+q x^4}}dxdydzdt\\\\
   =\frac{\pi ^{3/2} \Gamma \left(\frac{1}{4}\right)
   \sqrt[4]{\frac{p}{q}} }{8\ 2^{3/8} \alpha ^{3/4} b^{7/8}
   \sqrt{p}}\left((-2+2 i) \log \left(\frac{\Gamma \left(\frac{-2 i \log (a)+i \log
   \left(\frac{b}{2}\right)+2 i \log (p)-2 i \log (\alpha )+2 \pi }{16 \pi }\right)}{\Gamma \left(\frac{-2 i \log (a)+i
   \log \left(\frac{b}{2}\right)+2 i \log (p)-2 i \log (\alpha )+10 \pi }{16 \pi }\right)}\right)\right.\\
   \left.-(2+2 i) \log
   \left(\frac{\Gamma \left(\frac{-2 i \log (a)+i \log \left(\frac{b}{2}\right)+2 i \log (p)-2 i \log (\alpha )+6 \pi
   }{16 \pi }\right)}{\Gamma \left(\frac{-2 i \log (a)+i \log \left(\frac{b}{2}\right)+2 i \log (p)-2 i \log (\alpha
   )+14 \pi }{16 \pi }\right)}\right)+i \pi +2 \log (\pi )+\log (64)\right)
\end{dmath}
\end{example}
\begin{example}
\begin{dmath}
\int_{{\mathbb{R}^{4}_{+}}}\frac{\sqrt{x} D_{\frac{3}{2}}(t \alpha ) e^{-b \left(y+z^2\right)-\frac{1}{4} \alpha ^2 t^2} \log \left(\log
   \left(\frac{a y}{p t z+q t x^3 z}\right)\right)}{\sqrt[4]{t} \sqrt[4]{y} z^{3/4} \left(p+q
   x^3\right)^{3/4}}dxdydzdt\\\\
   =\frac{\pi ^2 \sqrt{\frac{p}{q}} }{6 \sqrt[8]{2} \alpha ^{3/4} b^{7/8}
   p^{3/4}}\left((-2+2 i) \log \left(\frac{\Gamma \left(\frac{-2 i \log (a)+i
   \log \left(\frac{b}{2}\right)+2 i \log (p)-2 i \log (\alpha )+2 \pi }{16 \pi }\right)}{\Gamma \left(\frac{-2 i \log
   (a)+i \log \left(\frac{b}{2}\right)+2 i \log (p)-2 i \log (\alpha )+10 \pi }{16 \pi }\right)}\right)\right.\\
   \left.-(2+2 i) \log
   \left(\frac{\Gamma \left(\frac{-2 i \log (a)+i \log \left(\frac{b}{2}\right)+2 i \log (p)-2 i \log (\alpha )+6 \pi
   }{16 \pi }\right)}{\Gamma \left(\frac{-2 i \log (a)+i \log \left(\frac{b}{2}\right)+2 i \log (p)-2 i \log (\alpha
   )+14 \pi }{16 \pi }\right)}\right)+i \pi +2 \log (\pi )+\log (64)\right)
\end{dmath}
\end{example}
\begin{example}
\begin{multline}
\int_{{\mathbb{R}^{4}_{+}}}\frac{x^{u-1} z^{\frac{3}{4}-u} D_u(t \alpha ) e^{-b \left(y+z^2\right)-\frac{1}{4} \alpha ^2 t^2} \log
   \left(\log \left(\frac{a y}{p t z+q t z x^{4 u}}\right)\right)}{\sqrt[4]{t} \sqrt{y} \sqrt{p+q x^{4 u}}}dxdydzdt\\
   =\frac{\pi
   ^{3/2} 2^{\frac{u}{2}-\frac{31}{8}} \Gamma \left(\frac{1}{4}\right) b^{\frac{u}{2}-\frac{11}{8}}
   \sqrt[4]{\frac{p}{q}} }{\alpha ^{3/4} \sqrt{p} u}\left((-2+2 i) \log \left(\frac{\Gamma \left(-\frac{i (2 \log (a)-\log (b)-2 \log (p)+2 \log
   (\alpha )+\log (2)+2 i \pi )}{16 \pi }\right)}{2 \Gamma \left(\frac{1}{2}-\frac{i (2 \log (a)-\log (b)-2 \log (p)+2
   \log (\alpha )+\log (2)+2 i \pi )}{16 \pi }\right)}\right)\right.\\
   \left.-(2+2 i) \log \left(\frac{\Gamma
   \left(\frac{1}{4}-\frac{i (2 \log (a)-\log (b)-2 \log (p)+2 \log (\alpha )+\log (2)+2 i \pi )}{16 \pi }\right)}{2
   \Gamma \left(\frac{3}{4}-\frac{i (2 \log (a)-\log (b)-2 \log (p)+2 \log (\alpha )+\log (2)+2 i \pi )}{16 \pi
   }\right)}\right)+i \pi +\log \left(4 \pi ^2\right)\right)
\end{multline}
\end{example}
\begin{example}
\begin{dmath}
\int_{{\mathbb{R}^{4}_{+}}}\frac{\sqrt{x} D_{\frac{3}{2}}(t) e^{-\frac{t^2}{4}-2 \left(y+z^2\right)} \log \left(\log \left(\frac{i y}{t
   x^6 z+t z}\right)\right)}{\sqrt[4]{t} \sqrt{x^6+1} \sqrt{y} z^{3/4}}dxdydzdt\\\\
   =\frac{\pi ^{3/2} \Gamma
   \left(\frac{1}{4}\right) \left(i \pi +\log \left(64 \pi ^2\right)-(2-2 i) \left(\log \left(\frac{\Gamma
   \left(\frac{3}{16}\right)}{\Gamma \left(\frac{11}{16}\right)}\right)+i \log \left(\frac{\Gamma
   \left(\frac{7}{16}\right)}{\Gamma \left(\frac{15}{16}\right)}\right)\right)\right)}{12\ 2^{3/4}}
\end{dmath}
\end{example}
\begin{example}
\begin{multline}
\int_{{\mathbb{R}^{4}_{+}}}\frac{x^{2/3} D_{\frac{5}{3}}(t) e^{-\frac{t^2}{4}-2 \left(y+z^2\right)} \log \left(\log \left(\frac{i y}{2 t
   x^{20/3} z+t z}\right)\right)}{\sqrt[4]{t} \sqrt{2 x^{20/3}+1} \sqrt{y} z^{11/12}}dxdydzdt\\\\
   =\frac{3 \pi ^{3/2} \Gamma
   \left(\frac{1}{4}\right) }{40\ 2^{5/6}}\left(i \pi +\log (64)+2 \log (\pi )\right.\\
   \left.-(2-2 i) \left(\log \left(\Gamma
   \left(\frac{3}{16}\right)\right)-\log \left(\Gamma \left(\frac{11}{16}\right)\right)+i \log \left(\frac{\Gamma
   \left(\frac{7}{16}\right)}{\Gamma \left(\frac{15}{16}\right)}\right)\right)\right)
\end{multline}
\end{example}
\begin{example}
\begin{multline}
\int_{{\mathbb{R}^{4}_{+}}}\frac{x^{\frac{\pi }{2}-1} z^{\frac{3}{4}-\frac{\pi }{2}} D_{\frac{\pi }{2}}(t) e^{-\frac{t^2}{4}-2
   \left(y+z^2\right)} \log \left(\log \left(-\frac{y}{2 t x^{2 \pi } z+t z}\right)\right)}{\sqrt[4]{t} \sqrt{2 x^{2
   \pi }+1} \sqrt{y}}dxdydzdt\\\\
   =2^{\frac{1}{2} (\pi -9)} \sqrt{\pi } \Gamma \left(\frac{1}{4}\right) \left(i \pi +\log
   (64)+(1-i) \log \left(\frac{\pi  \Gamma \left(\frac{3}{4}\right)^2}{\Gamma
   \left(\frac{1}{4}\right)^2}\right)\right)
\end{multline}
\end{example}
\begin{example}
\begin{dmath}
\int_{{\mathbb{R}^{4}_{+}}}\frac{x^{\frac{2 \pi }{5}-1} z^{\frac{3}{4}-\frac{2 \pi }{5}} D_{\frac{2 \pi }{5}}(t) e^{-\frac{t^2}{4}-2
   \left(y+z^2\right)} \log \left(\log \left(-\frac{y}{3 t x^{8 \pi /5} z+t z}\right)\right)}{\sqrt[4]{t} \sqrt{3 x^{8
   \pi /5}+1} \sqrt{y}}dxdydzdt\\\\
   =\frac{5\times 2^{\frac{2 \pi }{5}-\frac{25}{4}} \sqrt{\pi } \Gamma \left(\frac{1}{4}\right) \left(i
   \pi +\log (64)+(1-i) \log \left(\frac{\pi  \Gamma \left(\frac{3}{4}\right)^2}{\Gamma
   \left(\frac{1}{4}\right)^2}\right)\right)}{\sqrt[4]{3}}
\end{dmath}
\end{example}
\begin{example}
\begin{multline}
\int_{{\mathbb{R}^{4}_{+}}}\frac{x^{\frac{2 \pi }{5}-1} z^{\frac{3}{4}-\frac{2 \pi }{5}} D_{\frac{2 \pi }{5}}(t) e^{-\frac{t^2}{4}-2
   \left(y+z^2\right)} \log \left(\log \left(-\frac{y}{4 t x^{8 \pi /5} z+t z}\right)\right)}{\sqrt[4]{t} \sqrt{4 x^{8
   \pi /5}+1} \sqrt{y}}dxdydzdt\\\\
   =5\ 2^{\frac{2 \pi }{5}-\frac{27}{4}} \sqrt{\pi } \Gamma \left(\frac{1}{4}\right) \left(i \pi
   +\log (64)+(1-i) \log \left(\frac{\pi  \Gamma \left(\frac{3}{4}\right)^2}{\Gamma
   \left(\frac{1}{4}\right)^2}\right)\right)
\end{multline}
\end{example}
\begin{example}
\begin{multline}
\int_{{\mathbb{R}^{4}_{+}}}\frac{\sqrt[5]{x} D_{\frac{6}{5}}(t) e^{-\frac{t^2}{4}-2 \left(y+z^2\right)} \log \left(\log \left(\frac{i
   y}{\pi  t x^{24/5} z+t z}\right)\right)}{\sqrt[4]{t} \sqrt{\pi  x^{24/5}+1} \sqrt{y} z^{9/20}}dxdydzdt\\\\
   =\frac{5 \pi ^{5/4}
   \Gamma \left(\frac{1}{4}\right) }{96 \sqrt[20]{2}}\left(i \pi +\log (64)+2 \log (\pi )\right.\\
   \left.-(2-2 i) \left(\log \left(\Gamma
   \left(\frac{3}{16}\right)\right)-\log \left(\Gamma \left(\frac{11}{16}\right)\right)+i \log \left(\frac{\Gamma
   \left(\frac{7}{16}\right)}{\Gamma \left(\frac{15}{16}\right)}\right)\right)\right)
\end{multline}
\end{example}
\begin{example}
\begin{multline}
\int_{{\mathbb{R}^{4}_{+}}}\frac{D_{\frac{1}{2}}\left(\frac{t}{\sqrt{2}}\right) e^{-\frac{t^2}{8}-y-z^2} \log \left(\log \left(-\frac{y}{t
   x^4 z+t z}\right)\right)}{\sqrt{t} \sqrt{x} \left(x^4+1\right)^{5/8} y^{3/8} \log \left(-\frac{y}{t x^4 z+t
   z}\right)}dxdydzdt\\\\
   =\frac{\sqrt{\pi } \log (2) \left(2 i \gamma +\pi -i \log \left(8 \pi ^2\right)\right) \Gamma
   \left(\frac{9}{8}\right)}{2^{3/4}}
\end{multline}
\end{example}
\begin{example}
\begin{multline}
\int_{{\mathbb{R}^{4}_{+}}}\frac{D_{\frac{1}{2}}\left(\frac{t}{\sqrt{2}}\right) e^{-\frac{t^2}{8}-y-z^2} \log \left(\log \left(-\frac{y}{t
   x^4 z+t z}\right)\right)}{\sqrt{t} \sqrt{x} \left(x^4+1\right)^{5/8} y^{3/8} \log ^2\left(-\frac{y}{t x^4 z+t
   z}\right)}dxdydzdt\\\\
   =\frac{\pi ^{3/2} \Gamma \left(\frac{1}{8}\right) (-24 \log (A)+2 \gamma -i \pi +\log (4))}{192\
   2^{3/4}}
\end{multline}
\end{example}
\begin{example}
\begin{multline}
\int_{{\mathbb{R}^{4}_{+}}}\frac{D_{\frac{1}{4}}\left(\frac{t}{\sqrt{2}}\right) e^{-\frac{t^2}{8}-y-z^2}}{t^{3/4} x^{3/4}
   \left(x^4+1\right)^{13/16} y^{3/16} \log ^2\left(-\frac{y}{t x^4 z+t z}\right)}dxdydzdt\\\\
   =\frac{(-1)^{3/4} \left(\pi ^2+48 i
   C\right) \Gamma \left(\frac{1}{16}\right)}{384\ 2^{7/8} \sqrt{\pi }}
\end{multline}
\end{example}
\begin{example}
\begin{multline}
\int_{{\mathbb{R}^{4}_{+}}}t^{-\frac{1}{2}+\frac{i}{2}} x^{-\frac{1}{2}+\frac{i}{2}} \left(x^4+1\right)^{-\frac{5}{8}+\frac{3 i}{8}}
   y^{-\frac{3}{8}-\frac{3 i}{8}} D_{\frac{1}{2}+\frac{i}{2}}\left(\frac{t}{\sqrt{2}}\right)\\
    e^{-\frac{t^2}{8}-y-z^2}
   \log \left(-\frac{y}{t x^4 z+t z}\right)dxdydzdt\\\\
   =\frac{i 2^{-\frac{3}{4}+\frac{i}{4}} e^{\pi /2} \pi ^{5/2} \Gamma
   \left(\frac{1}{8}+\frac{i}{8}\right)}{\left(1+e^{\pi }\right)^2}
\end{multline}
\end{example}
\begin{example}
\begin{multline}
\int_{{\mathbb{R}^{4}_{+}}}\frac{t^{-\frac{3}{4}+\frac{i}{4}} x^{-\frac{3}{4}+\frac{i}{4}} \left(\sqrt{3}
   x^4+1\right)^{-\frac{13}{16}+\frac{3 i}{16}} y^{-\frac{3}{16}-\frac{3 i}{16}}
   D_{\frac{1}{4}+\frac{i}{4}}\left(\frac{t}{\sqrt{2}}\right) e^{-\frac{t^2}{8}-y-z^2}}{\log \left(-\frac{y}{\sqrt{3} t
   x^4 z+t z}\right)}dxdydzdt\\\\
   =-(-1)^{\frac{1}{4}+\frac{i}{4}} 2^{-\frac{23}{8}+\frac{i}{8}} 3^{-\frac{1}{32}-\frac{i}{32}}
   \sqrt{\pi } \log \left(1+i e^{\pi /2}\right) \Gamma \left(\frac{1}{16}+\frac{i}{16}\right)
\end{multline}
\end{example}
\begin{example}
Definite integral in terms of the Hurwitz zeta function $\zeta(s,v)$.
\begin{dmath}\label{eq2}
\int_{0}^{\infty}\int_{0}^{\infty}\int_{0}^{\infty}\frac{y^5 e^{-y^6-z} \text{Ai}(x \alpha )^2 \log ^k\left(-\frac{x y^2}{z}\right)}{\sqrt{x} \sqrt{z}}dxdydz\\
=\frac{1}{9 \sqrt{\alpha
   }}2^{2 k-1} e^{\frac{i \pi  k}{2}} \pi ^{k+\frac{1}{2}} \left(\zeta
   \left(-k,\frac{i (3 \log (\alpha )+\log (12)-6 i \pi )}{12 \pi }\right)-\zeta \left(-k,\frac{i (3 \log (\alpha )+\log (12)-12 i \pi )}{12 \pi }\right)\right)
\end{dmath}
\end{example}
\begin{example}
Definite integral in terms of the Hurwitz zeta function $\zeta(s,v)$.
\begin{multline}\label{eq3}
\int_{0}^{\infty}\int_{0}^{\infty}\int_{0}^{\infty}\frac{y^5 \text{Ai}(x)^2 e^{-\frac{y^6}{12}-z} \log ^k\left(\frac{i x y^2}{z}\right)}{\sqrt{x} \sqrt{z}}dxdydz\\
=\frac{1}{3} 2^{2 k+1} e^{\frac{i \pi  k}{2}} \pi ^{k+\frac{1}{2}}\left(\zeta \left(-k,\frac{3}{8}\right)-\zeta \left(-k,\frac{7}{8}\right)\right)
\end{multline}
\end{example}
\begin{example}
\begin{dmath}
\int_{0}^{\infty}\int_{0}^{\infty}\int_{0}^{\infty}\frac{y^5 \text{Ai}\left(\frac{x}{2^{2/3} \sqrt[3]{3}}\right)^2 e^{-y^6-z}}{\sqrt{x} \sqrt{z} \left(\log ^2\left(\frac{x y^2}{z}\right)+\pi
   ^2\right)}dxdydz=\frac{\log (2)}{6\ 2^{2/3} 3^{5/6} \pi ^{3/2}}
\end{dmath}
\textup{and}
\begin{dmath}
\int_{0}^{\infty}\int_{0}^{\infty}\int_{0}^{\infty}\frac{y^5 \text{Ai}\left(\frac{x}{2^{2/3} \sqrt[3]{3}}\right)^2 e^{-y^6-z} \log \left(\frac{x y^2}{z}\right)}{\sqrt{x} \sqrt{z} \left(\log ^2\left(\frac{x y^2}{z}\right)+\pi
   ^2\right)}dxdydz=0
\end{dmath}
\end{example}
\begin{example}
\begin{dmath}
\int_{0}^{\infty}\int_{0}^{\infty}\int_{0}^{\infty}\frac{y^5 \text{Ai}(x)^2 e^{-\frac{y^6}{12}-z}}{\sqrt{x} \sqrt{z} \left(4 \log ^2\left(\frac{x y^2}{z}\right)+\pi ^2\right)}dxdydz=\frac{\psi
   ^{(0)}\left(\frac{7}{8}\right)-\psi ^{(0)}\left(\frac{3}{8}\right)}{12 \pi ^{3/2}}
\end{dmath}
\textup{and}
\begin{dmath}
\int_{0}^{\infty}\int_{0}^{\infty}\int_{0}^{\infty}\frac{y^5 \text{Ai}(x)^2 e^{-\frac{y^6}{12}-z} \log \left(\frac{x y^2}{z}\right)}{\sqrt{x} \sqrt{z} \left(\log ^2\left(\frac{x
   y^2}{z}\right)+\frac{\pi ^2}{4}\right)}dxdydz=0
\end{dmath}
\end{example}
\subsubsection{Derived using equation (\ref{eq_beta_di})}
\begin{example}
\begin{multline}
\int _0^1\int _0^1\frac{x^n \left(-1+x^{-1-n}\right)^{-m} y^n \left(-1+y^{-1-n}\right)^{\frac{1}{2}-m} \left(1-\left(-1+x^{-1-n}\right)^{2 m} \left(-1+y^{-1-n}\right)^{2
   m}\right)}{\log \left(\left(-1+x^{-1-n}\right) \left(-1+y^{-1-n}\right)\right)}dydx\\
=\frac{e^{-2 i m \pi }}{8 (1+n)^2 \pi } \left(-8 e^{2 i m \pi } m \pi ^2 \left((-1+2 m) \coth ^{-1}\left(e^{2 i m \pi
   }\right)-(1+2 m) \tanh ^{-1}\left(e^{2 i m \pi }\right)\right)\right. \\ \left.
+\Phi \left(e^{-4 i m \pi },3,\frac{1}{2}\right)+i \left((-1+4 m) \pi  \Phi \left(e^{-4 i m \pi },2,\frac{1}{2}\right)\right.\right. \\ \left.\left.
+e^{4
   i m \pi } \left((\pi +4 m \pi ) \Phi \left(e^{4 i m \pi },2,\frac{1}{2}\right)+i \Phi \left(e^{4 i m \pi },3,\frac{1}{2}\right)\right)\right)\right)
\end{multline}
\end{example}
\begin{example}
\begin{multline}
\int _0^1\int _0^1\frac{x^n y^n \sqrt[4]{-1+y^{-1-n}}}{\sqrt[4]{-1+x^{-1-n}} \log \left(-\left(\left(-1+x^{-1-n}\right) \left(-1+y^{-1-n}\right)\right)\right)}dydx\\
=-\frac{i
   \left(\zeta \left(3,\frac{3}{8}\right)-\zeta \left(3,\frac{7}{8}\right)+2 \sqrt{2} \pi ^2 \left(\pi +2 \log \left(\tan \left(\frac{\pi }{8}\right)\right)\right)\right)}{64 (1+n)^2 \pi
   }
\end{multline}
\end{example}
\begin{example}
\begin{multline}
\int _0^1\int _0^1\frac{x^n y^n \sqrt[4]{-1+y^{-1-n}} \log \left(\log \left(a \left(-1+x^{-1-n}\right) \left(-1+y^{-1-n}\right)\right)\right)}{\sqrt[4]{-1+x^{-1-n}}}dydx\\
=\frac{1}{16 (1+n)^2}\left(\pi ^2
   (i \pi +\log (64)+2 \log (\pi )-4 \log (-6 \pi -i \log (a))+4 \log (-2 \pi -i \log (a)))\right. \\ \left.
+4 \pi ^2 \left(-\text{log$\Gamma $}\left(-\frac{3}{4}-\frac{i \log (a)}{8 \pi
   }\right)+\text{log$\Gamma $}\left(-\frac{1}{4}-\frac{i \log (a)}{8 \pi }\right)\right)-\psi ^{(1)}\left(\frac{1}{4}-\frac{i \log (a)}{8 \pi }\right)\right. \\ \left.
+\psi ^{(1)}\left(\frac{3}{4}-\frac{i
   \log (a)}{8 \pi }\right)\right)
\end{multline}
\end{example}
\begin{example}
\begin{multline}
\int _0^1\int _0^1x^n y^n \sqrt{-1+y^{-1-n}}\\
 \log \left(a \left(-1+x^{-1-n}\right) \left(-1+y^{-1-n}\right)\right) \log \left(\log \left(a \left(-1+x^{-1-n}\right)
   \left(-1+y^{-1-n}\right)\right)\right)dydx\\
=\frac{\pi }{2 (1+n)^2} \left(2+i \pi +\log (a)+\frac{1}{2} i \pi  \log (a)+2 \log (2) (2-i \pi +\log (a))+2 \log (\pi )\right. \\ \left.
-2 i \pi  \log (\pi )+\log (a) \log
   (\pi )+4 i \pi  \log \left(-\frac{1}{2}-\frac{i \log (a)}{4 \pi }\right)\right. \\ \left.+4 i \pi  \text{log$\Gamma $}\left(-\frac{1}{2}-\frac{i \log (a)}{4 \pi }\right)+2 \psi
   ^{(0)}\left(\frac{1}{2}-\frac{i \log (a)}{4 \pi }\right)\right)
\end{multline}
\end{example}
\begin{example}
\begin{multline}
\int _0^1\int _0^1x^n y^n \sqrt{-1+y^{-1-n}} \log ^2\left(a \left(-1+x^{-1-n}\right) \left(-1+y^{-1-n}\right)\right)\\
 \log \left(\log \left(a \left(-1+x^{-1-n}\right)
   \left(-1+y^{-1-n}\right)\right)\right)dydx\\
=\frac{4 \pi ^2}{(1+n)^2} \left(\frac{\left(4 \pi ^2+3 \log ^2(a)\right) (1+i \pi +\log (16)+2 \log (\pi ))}{48 \pi }\right. \\ \left.
+\frac{\log (a) (3+i \pi +\log (16)+2
   \log (\pi ))}{4 \pi }+2 i \left(-\frac{1}{2} \log (2 \pi )+\log \left(-\frac{1}{2}-\frac{i \log (a)}{4 \pi }\right)\right.\right. \\ \left.\left.+\text{log$\Gamma $}\left(-\frac{1}{2}-\frac{i \log (a)}{4 \pi
   }\right)\right)-4 \pi  \zeta'\left(-1,\frac{1}{2}-\frac{i \log (a)}{4 \pi }\right)\right)
\end{multline}
\end{example}
\subsubsection{Derived using equation (\ref{eq_genroot_sec})}
\begin{example}
\begin{multline}
\int_0^1 \frac{x^n \left(-1+x^{-1-n}\right)^{-\frac{1}{2}-m} \left(-1+\left(-1+x^{-1-n}\right)^{2
   m}\right)}{\log \left(-1+x^{-1-n}\right)} \, dx\\
   =\frac{1}{2 (1+n) \pi }\left(2 i \pi  \left((1+2 m) \cot ^{-1}\left(e^{i m \pi
   }\right)+(-1+2 m) \tan ^{-1}\left(e^{i m \pi }\right)\right)\right. \\ \left.
   +e^{-i m \pi } \Phi \left(-e^{-2 i m \pi
   },2,\frac{1}{2}\right)-e^{i m \pi } \Phi \left(-e^{2 i m \pi },2,\frac{1}{2}\right)\right)
\end{multline}
\end{example}
\begin{example}
\begin{multline}
\int_0^1 \frac{x^n \log ^k\left(1-x^{-1-n}\right)}{\sqrt{-1+x^{-1-n}}} \, dx\\
=\frac{e^{\frac{i k \pi }{2}}
   \left(-2^k \left(-1+2^{1+k}\right) \pi ^{1+k} \zeta (-k)-2 i \left(-1+2^k\right) \cos \left(\frac{k \pi
   }{2}\right) \Gamma (1+k) \zeta (k)\right)}{1+n}
\end{multline}
\end{example}
\begin{example}
\begin{equation}
\int_0^1 \frac{x^n}{\sqrt{-1+x^{-1-n}} \log ^2\left(1-x^{-1-n}\right)} \, dx=-\frac{\pi ^3-18 i \zeta (3)}{48
   (1+n) \pi ^2}
\end{equation}
\end{example}
\begin{example}
\begin{multline}
\int_0^1 \frac{x^n \log \left(1-x^{-1-n}\right) \log \left(\log
   \left(1-x^{-1-n}\right)\right)}{\sqrt{-1+x^{-1-n}}} \, dx\\
   =-\frac{\pi  \left(i \pi +\frac{\pi ^2}{2}\right)}{2
   (1+n)}\\+\frac{i \left(i \pi -2 \pi ^2 \log \left(\frac{A^3}{\sqrt[3]{2} \sqrt[4]{e}}\right)+i \pi  \log (2 \pi
   )+\frac{1}{2} \pi ^2 \log (2 \pi )-2 i \pi  \Phi'(-1,0,1)\right)}{1+n}
\end{multline}
\end{example}
\begin{example}
\begin{multline}
\int_0^1 \frac{x^n \log ^2\left(1-x^{-1-n}\right) \log \left(\log
   \left(1-x^{-1-n}\right)\right)}{\sqrt{-1+x^{-1-n}}} \, dx\\
   =\frac{\pi ^3}{1+n}-\frac{i \pi ^2-8 i \pi ^2 \log
   \left(\frac{A^3}{\sqrt[3]{2} \sqrt[4]{e}}\right)+2 i \pi ^2 \log (2 \pi )-7 \pi  \zeta (3)}{1+n}
\end{multline}
\end{example}
\begin{example}
\begin{multline}
\int_0^1 \frac{x^n \log \left(\log \left(1-x^{-1-n}\right)\right)}{\sqrt{-1+x^{-1-n}} \log
   \left(1-x^{-1-n}\right)} \, dx\\
=\frac{\pi  \left(12 i \log (2) \left(2 \gamma +\log \left(\frac{1}{8 \pi
   ^2}\right)\right)+\pi  (2-i \pi +\log (1024)-2 \log (\pi ))\right)-24 Li'_{2}(-1)}{48 (1+n) \pi
   }
\end{multline}
\end{example}
\begin{example}
\begin{multline}
\int_0^1 x^n \log ^k\left(a \left(-1+x^{-1-n}\right)\right) \, dx=-\frac{e^{\frac{i k \pi }{2}} k (2 \pi )^k
   \zeta \left(1-k,\frac{\pi -i \log (a)}{2 \pi }\right)}{1+n}
\end{multline}
\end{example}
\begin{example}
\begin{multline}
\int_0^1 x^n \log \left(a \left(-1+x^{-1-n}\right)\right) \log \left(\log \left(a
   \left(-1+x^{-1-n}\right)\right)\right) \, dx\\
=\frac{1}{1+n}\left(\log (a)+\frac{1}{2} i \pi  \log (a)+\log (a) (\log (2)+\log
   (\pi ))\right. \\ \left.
+2 i \pi  \left(-\frac{1}{2} \log (2 \pi )+\log \left(-1+\frac{\pi -i \log (a)}{2 \pi
   }\right)+\text{log$\Gamma $}\left(-1+\frac{\pi -i \log (a)}{2 \pi }\right)\right)\right)
\end{multline}
\end{example}
\begin{example}
\begin{equation}
\int_0^1 x^n \log \left(1-x^{-1-n}\right) \log \left(\log \left(1-x^{-1-n}\right)\right) \, dx=-\frac{\pi 
   (-2 i+\pi )}{2 (1+n)}
\end{equation}
\end{example}
\begin{figure}[H]
\includegraphics[scale=0.5]{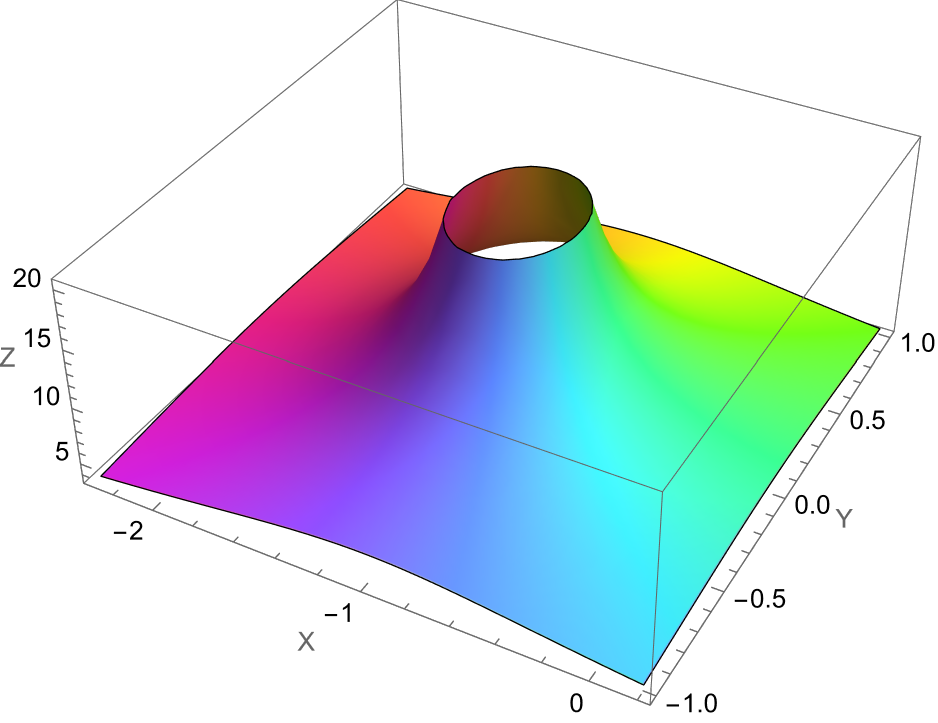}
\caption{Plot of  $f(n)=-\frac{\pi (-2 i+\pi )}{2 (1+n)}$, $n\in\mathbb{C}$.}
   \label{fig:fig2}
\end{figure}
\vspace{-6pt}
\begin{example}
\begin{equation}
\int_0^1 \log \left(\frac{x-1}{x}\right) \log \left(\log \left(\frac{x-1}{x}\right)\right) \, dx=\frac{1}{2}
   \pi  (2 i-\pi )
\end{equation}
\end{example}
\begin{figure}[H]
\includegraphics[scale=0.5]{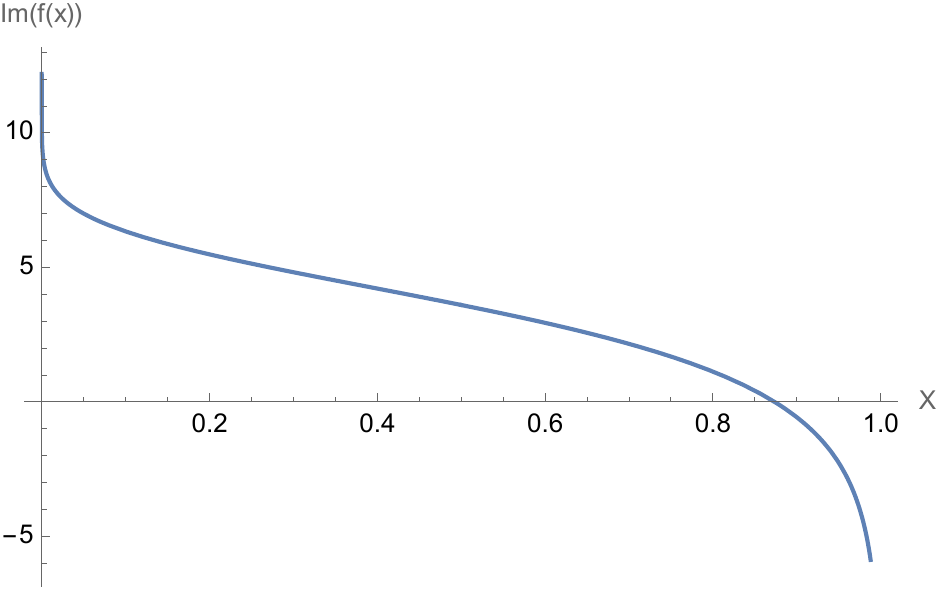}
\caption{Plot of  $f(x)=\log \left(\frac{x-1}{x}\right) \log \left(\log \left(\frac{x-1}{x}\right)\right)$, $x\in\mathbb{R}$.}
   \label{fig:fig2}
\end{figure}
\vspace{-6pt}
\begin{figure}[H]
\includegraphics[scale=0.5]{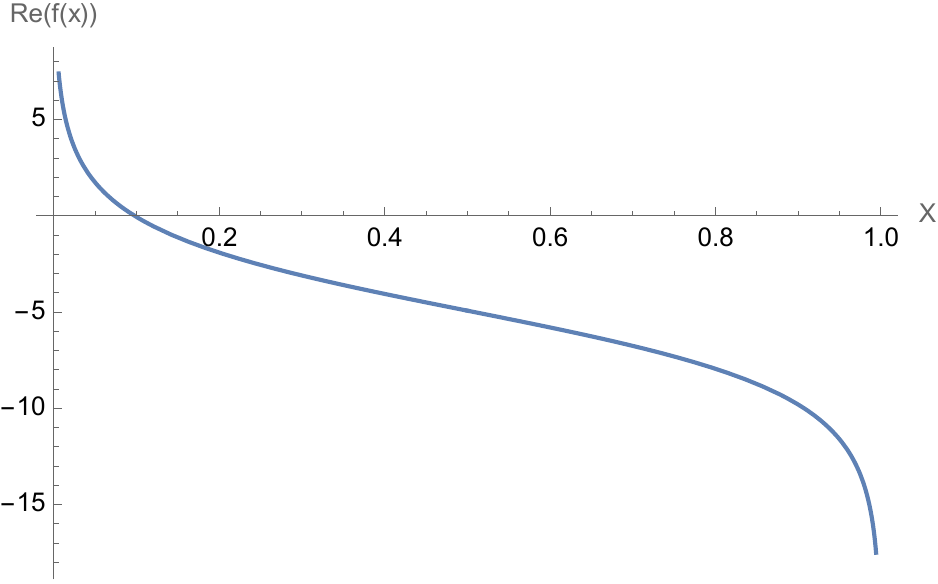}
\caption{Plot of  $f(x)=\log \left(\frac{x-1}{x}\right) \log \left(\log \left(\frac{x-1}{x}\right)\right)$, $x\in\mathbb{R}$.}
   \label{fig:fig2}
\end{figure}
\vspace{-6pt}
\begin{example}
\begin{multline}
\int_0^1 \frac{x^n}{a^2 \pi ^2-\log ^2\left(-1+x^{-1-n}\right)} \, dx=\frac{i \left(-\psi
   ^{(1)}\left(\frac{1}{2}-\frac{i a}{2}\right)+\psi ^{(1)}\left(\frac{1}{2}+\frac{i a}{2}\right)\right)}{4 a (1+n)
   \pi ^2}
\end{multline}
where $Im(a) \leq 0, Im(n)<0$
\end{example}
\begin{example}
\begin{equation}
\int_0^1 \frac{x^n \log \left(-1+x^{-1-n}\right)}{a^2-\log ^2\left(-1+x^{-1-n}\right)} \, dx = \begin{cases}
  \frac{i \pi }{2
   (1+n) (1+\cosh (a))}  ;& Im(n)>0 \\
 - \frac{i \pi }{2
   (1+n) (1+\cosh (a))} ;& Im(n)<0 
\end{cases}
\end{equation}
\end{example}
\begin{example}
\begin{multline}
\int_0^1 \frac{x^n \log \left(-1+x^{-1-n}\right)}{\left(a^2-\log ^2\left(-1+x^{-1-n}\right)\right)^2} \,
   dx= \begin{cases}
  \frac{i \pi  \sinh (a)}{4 a (1+n) (1+\cosh (a))^2}  ;& Im(n)>0 \\
  -\frac{i \pi  \sinh (a)}{4 a (1+n) (1+\cosh (a))^2}  ;& Im(n)<0
\end{cases}
\end{multline}
\end{example}
\begin{example}
\begin{equation}
\int_0^1 x^n \log ^k\left(1-x^{-1-n}\right) \, dx=-\frac{\left(1+e^{i k \pi }\right) }{1+n}\Gamma (1+k) \zeta
   (k)
\end{equation}
\end{example}
\begin{example}
\begin{equation}
\int_0^1 x^n \log \left(\log \left(1-x^{-1-n}\right)\right) \, dx=\frac{i \pi +\log \left(4 e^{-2 \gamma }
   \pi ^2\right)}{2 (1+n)}
\end{equation}
\end{example}
\begin{example}
\begin{equation}
\int_0^1 \frac{x^n}{\log \left(1-x^{-1-n}\right)} \, dx = \begin{cases}
  \frac{i \pi }{12 (1+n)} ;& Im(n)>0  \\
  -\frac{i \pi }{12 (1+n)} ;& Im(n)<0
\end{cases}
\end{equation}
\end{example}
\begin{figure}[H]
\includegraphics[scale=0.5]{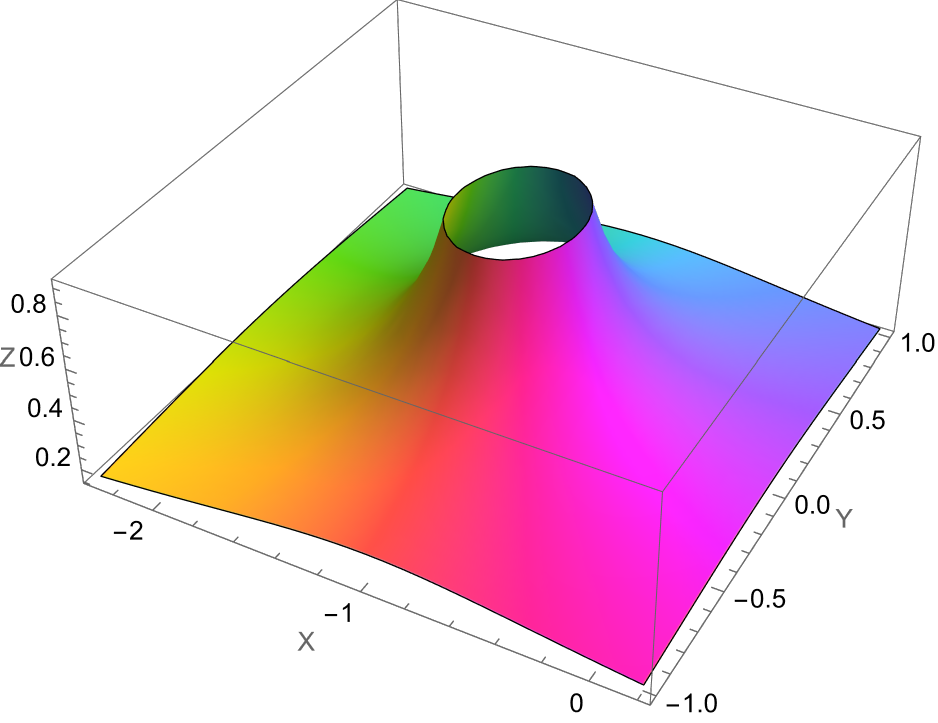}
\caption{Plot of  $f(n)=-\frac{i \pi }{12 (1+n)}$, $n\in\mathbb{C}$.}
   \label{fig:fig2}
\end{figure}
\vspace{-6pt}
\begin{example}
\begin{multline}
\int_0^1 \frac{x^n \log \left(\log \left(1-x^{-1-n}\right)\right)}{\log \left(1-x^{-1-n}\right)} \,
   dx = \begin{cases}
  \frac{\pi  \left(\pi +2 i \log \left(\frac{e^{1+\gamma }}{A^{12}}\right)\right)}{24 (1+n)} ;& Im(n)<0  \\
- \frac{\pi  \left(\pi +2 i \log \left(\frac{e^{1+\gamma }}{A^{12}}\right)\right)}{24 (1+n)} ;& Im(n)>0
\end{cases}
\end{multline}
\end{example}
\begin{example}
\begin{multline}
\int_0^1 \frac{x^n \log \left(x^{-1-n}-1\right)}{\left(a^2 \pi ^2-\log ^2\left(x^{-1-n}-1\right)\right)^2} \,
   dx = \begin{cases}
  \frac{i \text{csch}^3(a \pi ) \sinh ^4\left(\frac{a \pi }{2}\right)}{a (1+n)}  ;& Im(n)>0  \\
 -\frac{i \text{csch}^3(a \pi ) \sinh ^4\left(\frac{a \pi }{2}\right)}{a (1+n)}  ;& Im(n)<0  
\end{cases}
\end{multline}
\end{example}
\subsubsection{Derived using equation (\ref{eq_iis})}
\begin{example}
\begin{multline}
\int_0^1 \left(\sum _{s=1}^n \frac{s^{-m} \left(1-s^{2 m} \left(-1+\frac{1}{x}\right)^{2 m}\right)
   \left(-1+\frac{1}{x}\right)^{-m}}{\log \left(s \left(-1+\frac{1}{x}\right)\right)}\right) \, dx\\
=\sum
   _{c=1}^n \frac{c^{-m} e^{-i m \pi }}{2 \pi } \left(2 m \pi  \Phi \left(e^{-2 i m \pi },1,\frac{\pi -i \log (c)}{2
   \pi }\right)-i \Phi \left(e^{-2 i m \pi },2,\frac{\pi -i \log (c)}{2 \pi }\right)\right. \\ \left.
+c^{2 m} e^{2 i m \pi }
   \left(2 m \pi  \Phi \left(e^{2 i m \pi },1,\frac{\pi -i \log (c)}{2 \pi }\right)+i \Phi \left(e^{2 i m \pi
   },2,\frac{\pi -i \log (c)}{2 \pi }\right)\right)\right)
\end{multline}
\end{example}
\begin{example}
\begin{multline}
\int_0^1 \left(\sum _{s=1}^n \sqrt{s} \sqrt{-1+\frac{1}{x}} \log ^k\left(a s
   \left(-1+\frac{1}{x}\right)\right)\right) \, dx\\
=\sum _{c=1}^n 2^{-1+2 k} \sqrt{c} e^{\frac{i k \pi }{2}} \pi ^k
   \left(-i k \left(\zeta \left(1-k,\frac{\pi -i \log (a)-i \log (c)}{4 \pi }\right)\right.\right. \\ \left.\left.
-\zeta \left(1-k,\frac{3 \pi -i
   \log (a)-i \log (c)}{4 \pi }\right)\right)+2 \pi  \left(\zeta \left(-k,\frac{\pi -i \log (a)-i \log (c)}{4 \pi
   }\right)\right.\right. \\ \left.\left.
-\zeta \left(-k,\frac{3 \pi -i \log (a)-i \log (c)}{4 \pi }\right)\right)\right)
\end{multline}
\end{example}
\begin{example}
\begin{multline}
\int_0^1 \left(\sum _{s=1}^n \log ^k\left(a s \left(-1+\frac{1}{x}\right)\right)\right) \, dx=-\sum _{c=1}^n
   k (2 i \pi )^k \zeta \left(1-k,\frac{\pi -i \log (a)-i \log (c)}{2 \pi }\right)
\end{multline}
\end{example}
\begin{example}
\begin{multline}
\int_0^1 \left(\sum _{s=1}^n \log \left(a s \left(-1+\frac{1}{x}\right)\right) \log \left(\log \left(a s
   \left(-1+\frac{1}{x}\right)\right)\right)\right) \, dx\\
=\sum _{c=1}^n \frac{1}{2} \left(\log (a) \left(2+i \pi
   +\log \left(4 \pi ^2\right)\right)+\log (c) \left(2+i \pi +\log \left(4 \pi ^2\right)\right)\right. \\ \left.
-2 i \pi  \left(\log
   \left(8 \pi ^3\right)-2 \log (-\pi -i \log (a)-i \log (c))\right)+4 i \pi  \text{log$\Gamma $}\left(-\frac{\pi +i
   \log (a)+i \log (c)}{2 \pi }\right)\right)
\end{multline}
\end{example}
\begin{example}
\begin{multline}
\int_0^1 \left(\sum _{s=1}^n \frac{\log \left(\log \left(s \left(1-\frac{1}{x}\right)\right)\right)}{\sqrt{s}
   \sqrt{-1+\frac{1}{x}}}\right) \, dx\\
=\sum _{c=1}^n \left(\frac{i \Phi \left(-1,1,\frac{2 \pi -i \log (c)}{2 \pi
   }\right)}{\sqrt{c}}+\frac{\pi  \log (2 i \pi )}{2 \sqrt{c}}-\frac{\pi 
   \Phi'\left(-1,0,\frac{2 \pi -i \log (c)}{2 \pi }\right)}{\sqrt{c}}\right)
\end{multline}
\end{example}
\begin{example}
\begin{multline}
\int_0^1 \left(\sum _{s=1}^n \frac{\log \left(\frac{s (-1+x)}{x}\right) \log \left(\log \left(\frac{s
   (-1+x)}{x}\right)\right)}{\sqrt{s} \sqrt{-1+\frac{1}{x}}}\right) \, dx\\
=\sum _{c=1}^n \frac{\pi }{2
   \sqrt{c}} \left(-2+(-2+i
   \pi +\log (c)) \log (2 i \pi )-4 i \pi  \Phi'\left(-1,-1,1-\frac{i \log (c)}{2 \pi
   }\right)\right. \\ \left.
+4 \Phi'\left(-1,0,1-\frac{i \log (c)}{2 \pi }\right)\right)
\end{multline}
\end{example}
\begin{example}
\begin{multline}
\int_0^1 \left(\sum _{s=1}^n \frac{\log \left(\log \left(\frac{s (-1+x)}{x}\right)\right)}{\sqrt{s}
   \sqrt{-1+\frac{1}{x}} \log \left(\frac{s (-1+x)}{x}\right)}\right) \, dx\\
=\sum _{c=1}^n \frac{1}{8 \sqrt{c} \pi }\left(-4 i \pi  \Phi
   \left(-1,1,1-\frac{i \log (c)}{2 \pi }\right) \log (2 i \pi )\right. \\ \left.
-(-1+\log (2 i \pi )) \left(\psi
   ^{(1)}\left(\frac{1}{2}-\frac{i \log (c)}{4 \pi }\right)-\psi ^{(1)}\left(1-\frac{i \log (c)}{4 \pi
   }\right)\right)\right. \\ \left.
+4 i \pi  \Phi'\left(-1,1,1-\frac{i \log (c)}{2 \pi }\right)+4
   \Phi'\left(-1,2,1-\frac{i \log (c)}{2 \pi }\right)\right)
\end{multline}
\end{example}
\subsubsection{Derived using equation (\ref{eq_319710})}
\begin{example}
\begin{multline}
\int _0^1\int _0^1\frac{(-((-1+x) x))^{-\frac{1}{2}-m} y^n \left(-1+y^{-1-n}\right)^{-m} \left((1-x)^{2 m}-x^{2 m}
   \left(-1+y^{-1-n}\right)^{2 m}\right)}{(1+p x) \log \left(\frac{x \left(-1+y^{-1-n}\right)}{1-x}\right)}dydx\\
=\frac{e^{-2 i
   m \pi } (1+p)^{-\frac{1}{2}-m}}{4 (1+n)} \left((1+p)^{2 m} \left(4 m \pi  \Phi \left(e^{-4 i m \pi },1,\frac{1}{2}+\frac{i \log
   (1+p)}{4 \pi }\right)\right.\right. \\ \left.\left.
-i \Phi \left(e^{-4 i m \pi },2,\frac{1}{2}+\frac{i \log (1+p)}{4 \pi }\right)\right)+e^{4 i m \pi }
   \left(4 m \pi  \Phi \left(e^{4 i m \pi },1,\frac{1}{2}+\frac{i \log (1+p)}{4 \pi }\right)\right.\right. \\ \left.\left.
+i \Phi \left(e^{4 i m \pi
   },2,\frac{1}{2}+\frac{i \log (1+p)}{4 \pi }\right)\right)\right)
\end{multline}
\end{example}
\begin{example}
\begin{multline}
\int _0^1\int _0^1\frac{y^n \left(a+i b+\log \left(\frac{x
   \left(-1+y^{-1-n}\right)}{1-x}\right)\right)^k}{\sqrt{-((-1+x) x)} (1+p x)}dydx\\
=-\frac{(4 i)^k k \pi ^{1+k} \zeta
   \left(1-k,-\frac{i (a+i b+2 i \pi -\log (1+p))}{4 \pi }\right)}{(1+n) \sqrt{1+p}}
\end{multline}
\end{example}
\begin{example}
\begin{multline}
\int _0^1\int _0^1\frac{y^n \left(i b+a p x+\log \left(\frac{x
   \left(-1+y^{-1-n}\right)}{1-x}\right)\right)}{\sqrt{(1-x) x} \left(p^2 x^2-1\right) \left(a^2+b^2-2 i b \log \left(\frac{x
   \left(-1+y^{-1-n}\right)}{1-x}\right)-\log ^2\left(\frac{x \left(-1+y^{-1-n}\right)}{1-x}\right)\right)}dydx\\
=-\frac{i
   \left(\sqrt{1+p} \psi ^{(1)}\left(\frac{i a+b+2 \pi +i \log (1-p)}{4 \pi }\right)+\sqrt{1-p} \psi ^{(1)}\left(\frac{-i
   a+b+2 \pi +i \log (1+p)}{4 \pi }\right)\right)}{8 (1+n) \sqrt{1-p^2}}
\end{multline}
where $-1<Re(p)<1$
\end{example}
\begin{example}
\begin{multline}
\int _0^1\int _0^1\frac{2 i \pi +x \log (2)+2 \log \left(\frac{x (-1+y)}{(-1+x) y}\right)}{\sqrt{1-x} \sqrt{x}
   \left(-4+x^2\right) \left(\pi ^2+\log ^2(2)-2 i \pi  \log \left(\frac{x (-1+y)}{(-1+x) y}\right)-\log ^2\left(\frac{x
   (-1+y)}{(-1+x) y}\right)\right)}dydx\\
=-\frac{i \left(\sqrt{3} \left(-8 C+\pi ^2\right)+\psi ^{(1)}\left(\frac{3}{4}-\frac{i
   \log \left(\frac{4}{3}\right)}{4 \pi }\right)\right)}{8 \sqrt{6}}
\end{multline}
\end{example}
\begin{example}
\begin{multline}
\int _0^1\int _0^1\frac{2 i \pi +\frac{1}{2} x \log (2)+\log \left(\frac{x (-1+y)}{(-1+x) y}\right)}{\sqrt{-((-1+x) x)}
   \left(-4+x^2\right) \left(\log ^2(2)+\left(2 \pi -i \log \left(\frac{x (-1+y)}{(-1+x)
   y}\right)\right)^2\right)}dydx\\
=-\frac{i \left(\pi ^2+2 \sqrt{3} \psi ^{(1)}\left(1-\frac{i \log \left(\frac{4}{3}\right)}{4
   \pi }\right)\right)}{96 \sqrt{2}}
\end{multline}
\end{example}
\begin{example}
\begin{multline}
\int _0^1\int _0^1\frac{y^n \left(a+i b+\log \left(\frac{x \left(-1+y^{-1-n}\right)}{1-x}\right)\right) \log \left(a+i
   b+\log \left(\frac{x \left(-1+y^{-1-n}\right)}{1-x}\right)\right)}{\sqrt{-((-1+x) x)} (1+p x)}dydx\\
=\frac{\pi  }{2 (1+n) \sqrt{1+p}}\left(-b \pi
   +\log (1+p) (-i \pi -2 \log (4 e \pi ))+2 i b \log (4 e \pi )\right. \\ \left.
+a \left(i \pi +\log \left(16 e^2 \pi ^2\right)\right)-4 i \pi 
   \left(\log \left(32 \pi ^3\right)-2 \log (-i a+b-2 \pi +i \log (1+p))\right)\right. \\ \left.
+8 i \pi  \text{log$\Gamma $}\left(\frac{-i
   a+b-2 \pi +i \log (1+p)}{4 \pi }\right)\right)
\end{multline}
\end{example}
\begin{example}
\begin{multline}
\int _0^1\int _0^1\frac{\log \left(\log \left(-\frac{x (-1+y)}{(-1+x) y}\right)\right)}{\sqrt[4]{1-x} x^{3/4}
   \sqrt[4]{-1+\frac{1}{y}}}dydx\\
=\frac{1}{4} \pi  \left(-4 i \sqrt{2} \log \left(\cot \left(\frac{\pi }{8}\right)\right)+\pi 
   \left(2 i \sqrt{2}+i \pi +2 \log (4 \pi )-4
   \Phi'\left(-1,0,\frac{3}{4}\right)\right)\right)
\end{multline}
\end{example}
\begin{example}
\begin{multline}
\int _0^1\int _0^1\frac{\log \left(\log \left(\frac{i x (-1+y)}{(-1+x) y}\right)\right)}{\sqrt[4]{1-x} x^{3/4}
   \sqrt[4]{-1+\frac{1}{y}}}dydx\\
=\frac{1}{4} \pi  \left(8 i \left(\cos \left(\frac{\pi }{8}\right) \log \left(\cot
   \left(\frac{3 \pi }{16}\right)\right)+\log \left(\tan \left(\frac{\pi }{16}\right)\right) \sin \left(\frac{\pi
   }{8}\right)\right)\right. \\ \left.
+\pi  \left(i \pi +i \cot \left(\frac{3 \pi }{16}\right)+\log (16)+2 \log (\pi )+i \tan \left(\frac{3 \pi
   }{16}\right)-4 \Phi'\left(-1,0,\frac{5}{8}\right)\right)\right)
\end{multline}
\end{example}
\begin{example}
\begin{multline}
\int _0^1\int _0^1\frac{\log \left(\frac{2 x (-1+y)}{(-1+x) y}\right) \log \left(\log \left(\frac{2 x (-1+y)}{(-1+x)
   y}\right)\right)}{\sqrt[4]{1-x} x^{3/4} (1+x) \sqrt[4]{-1+\frac{1}{y}}}dydx\\
=\frac{\pi ^2 \left(-2-4 i C-i \pi -2 \log (4 \pi
   )+4 \Phi'\left(-1,0,\frac{1}{2}\right)\right)}{\sqrt[4]{2}}
\end{multline}
\end{example}
\begin{example}
\begin{multline}
\int _0^1\int _0^1\frac{\sqrt[4]{-1+\frac{1}{y}} \log \left(\log \left(-\frac{x (-1+y)}{(-1+x)
   y}\right)\right)}{(1-x)^{3/4} \sqrt[4]{x}}dydx\\
=\frac{1}{4} \pi  \left(4 i \sqrt{2} \log \left(\cot \left(\frac{\pi
   }{8}\right)\right)+\pi  \left(-2 i \sqrt{2}+i \pi +2 \log (4 \pi )-4
   \Phi'\left(-1,0,\frac{3}{4}\right)\right)\right)
\end{multline}
\end{example}
\begin{example}
A limiting case involving a nested logarithm function. 
\begin{multline}
\int _0^1\int _0^1\frac{y^n \log \left(\log \left(\frac{a x \left(-1+y^{-1-n}\right)}{1-x}\right)\right)}{\sqrt[4]{1-x}
   x^{3/4} (1+p x) \sqrt[4]{-1+y^{-1-n}}}dydx\\
=\frac{\pi}{4 (1+n) \sqrt[4]{1+p}}  \left(\pi  (i \pi +\log (64)+2 \log (\pi )-4 \log (-6 \pi -i \log
   (a)+i \log (1+p))\right. \\ \left.+4 \log (-2 \pi -i \log (a)+i \log (1+p)))+4 \pi  \text{log$\Gamma $}\left(-\frac{2 \pi +i \log (a)-i \log
   (1+p)}{8 \pi }\right)\right. \\ \left.
-4 \pi  \text{log$\Gamma $}\left(-\frac{6 \pi +i \log (a)-i \log (1+p)}{8 \pi }\right)-2 i \psi
   ^{(0)}\left(\frac{2 \pi -i \log (a)+i \log (1+p)}{8 \pi }\right)\right. \\ \left.
+2 i \psi ^{(0)}\left(\frac{6 \pi -i \log (a)+i \log
   (1+p)}{8 \pi }\right)\right)
\end{multline}
\end{example}
\begin{example}
\begin{multline}
\int _0^1\int _0^1\frac{\log \left(\log \left(-\frac{x (-1+y)}{(-1+x) y}\right)\right)}{\sqrt[4]{1-x} x^{3/4}
   \sqrt[4]{-1+\frac{1}{y}}}dydx\\
=\frac{1}{4} \pi  \left(\pi  \left(i \pi +2 \log \left(\frac{8 \pi }{25}\right)-4
   \text{log$\Gamma $}\left(-\frac{5}{8}\right)+4 \text{log$\Gamma $}\left(-\frac{1}{8}\right)\right)\right. \\ \left.-2 i \left(\psi
   ^{(0)}\left(\frac{3}{8}\right)-\psi ^{(0)}\left(\frac{7}{8}\right)\right)\right)
\end{multline}
\end{example}
\begin{example}
\begin{multline}
\int _0^1\int _0^1\frac{\log \left(\log \left(\frac{2 x (-1+y)}{(-1+x) y}\right)\right)}{\sqrt[4]{1-x} x^{3/4} (1+x)
   \sqrt[4]{-1+\frac{1}{y}}}dydx\\
=\frac{\pi ^2 \left(2 i+i \pi +\log \left(\frac{64 \pi ^2}{81}\right)-4 \text{log$\Gamma
   $}\left(-\frac{3}{4}\right)+4 \text{log$\Gamma $}\left(-\frac{1}{4}\right)\right)}{4 \sqrt[4]{2}}
\end{multline}
\end{example}
\subsubsection{Derived using equations (\ref{eq_3311_1}) and (\ref{eq_3311_2})}
\begin{example}
\begin{multline}
\int _0^{\infty }\int _0^{\infty }\int _0^{\infty }\frac{e^{-\frac{x}{4}-\frac{y}{4}} \text{Ei}(-z)}{\sqrt[4]{x y} \sqrt[4]{x+y} z^{3/4} \log \left(\frac{i \sqrt{x y}}{\sqrt{x+y}
   \sqrt{z}}\right)}dzdydx\\
=\sum _{y=0}^{\infty } -16 (-1)^{3/4} \pi ^{3/2} \Gamma \left(0,\frac{1}{2} i \pi  (1+y)\right)
\end{multline}
\end{example}
\begin{example}
\begin{multline}
\int _0^{\infty }\int _0^{\infty }\int _0^{\infty }\frac{e^{-\frac{x}{4}-\frac{y}{4}} \text{Ei}(-z)}{\sqrt[4]{x y} \sqrt[4]{x+y} z^{3/4} \log \left(-\frac{\sqrt{x y}}{\sqrt{x+y}
   \sqrt{z}}\right)}dzdydx=\sum _{y=0}^{\infty } 16 \pi ^{3/2} \Gamma \left(0,\frac{1}{4} i \pi  (3+2 y)\right)
\end{multline}
\end{example}
\begin{example}
\begin{multline}
\int _0^{\infty }\int _0^{\infty }\int _0^{\infty }\frac{e^{-\frac{x}{4}-\frac{y}{4}} \text{Ei}(-z)}{\sqrt[4]{x y} \sqrt[4]{x+y} z^{3/4} \log ^2\left(-\frac{\sqrt{x y}}{\sqrt{x+y}
   \sqrt{z}}\right)}dzdydx=\sum _{y=0}^{\infty } 8 \pi ^{3/2} \Gamma \left(-1,\frac{1}{4} i \pi  (3+2 y)\right)
\end{multline}
\end{example}
\begin{example}
\begin{multline}
\int _0^{\infty }\int _0^{\infty }\int _0^{\infty }\frac{e^{-\frac{x}{4}-\frac{y}{4}} \text{Ei}(-2 z)}{\sqrt[4]{x y} \sqrt[4]{x+y} z^{3/4} \log ^2\left(\frac{i \sqrt{x y}}{\sqrt{x+y}
   \sqrt{z}}\right)}dzdydx\\
=\sum _{y=0}^{\infty } -8 (-1)^{3/4} \pi ^{3/2} \Gamma \left(-1,\frac{1}{4} (2 i \pi  (1+y)+\log (2))\right)
\end{multline}
\end{example}
\begin{example}
\begin{multline}
\int _0^{\infty }\int _0^{\infty }\int _0^{\infty }\frac{e^{-\frac{x}{4}-\frac{y}{4}} \text{Ei}(-z) \log \left(\log \left(-\frac{\sqrt{x y}}{\sqrt{x+y}
   \sqrt{z}}\right)\right)}{\sqrt{x y} \sqrt{z} \log \left(-\frac{\sqrt{x y}}{\sqrt{x+y} \sqrt{z}}\right)}dzdydx\\
=\sum _{y=0}^{\infty } 16 i \pi ^{3/2} \left(\Gamma \left(0,\frac{1}{2} i \pi 
   (3+2 y)\right) \log \left(\frac{1}{2} i \pi  (3+2 y)\right)+G_{2,3}^{3,0}\left(\frac{1}{2} i \pi  (3+2 y)|
\begin{array}{c}
 1,1 \\
 0,0,0 \\
\end{array}
\right)\right)
\end{multline}
\end{example}
\begin{example}
\begin{equation}
\int _0^{\infty }\int _0^{\infty }\int _0^{\infty }\frac{e^{-\frac{x}{4}-\frac{y}{4}} \text{Ei}(-z) \left(-2-\log \left(\frac{i \sqrt{x y}}{\sqrt{x+y} \sqrt{z}}\right)\right)}{\sqrt{x
   y} \sqrt{z} \log ^3\left(\frac{i \sqrt{x y}}{\sqrt{x+y} \sqrt{z}}\right)}dzdydx=-\frac{1}{3} \left(4 \pi ^{3/2}\right)
\end{equation}
\end{example}
\begin{example}
\begin{multline}
\int _0^{\infty }\int _0^{\infty }\int _0^{\infty }\frac{e^{-\frac{x}{4}-\frac{y}{4}} \text{Ei}(-z) \left(3+\log \left(\frac{i \sqrt{x y}}{\sqrt{x+y} \sqrt{z}}\right)\right)}{\sqrt{x
   y} \sqrt{z} \log ^4\left(\frac{i \sqrt{x y}}{\sqrt{x+y} \sqrt{z}}\right)}dzdydx=-\frac{12 i \zeta (3)}{\pi ^{3/2}}
\end{multline}
\end{example}
\begin{example}
\begin{multline}
\int _0^{\infty }\int _0^{\infty }\int _0^{\infty }\frac{e^{-\frac{x}{4}-\frac{y}{4}} \text{Ei}(-z) \left(\frac{1}{2}-\log \left(\frac{i \sqrt{x y}}{\sqrt{x+y}
   \sqrt{z}}\right)\right)}{\sqrt{x y} \sqrt{z} \sqrt{\log \left(\frac{i \sqrt{x y}}{\sqrt{x+y} \sqrt{z}}\right)}}dzdydx\\
   =4 \sqrt[4]{-1} \left(-1+2 \sqrt{2}\right) \pi  \zeta
   \left(\frac{3}{2}\right)
\end{multline}
\end{example}
\begin{example}
\begin{multline}
\int _0^{\infty }\int _0^{\infty }\int _0^{\infty }\frac{e^{-\frac{x}{4}-\frac{y}{4}} \text{Ei}(-z) \left(1+\log \left(-\frac{\sqrt{x y}}{\sqrt{x+y} \sqrt{z}}\right)\right)}{\sqrt{x
   y} \sqrt{z} \log ^2\left(-\frac{\sqrt{x y}}{\sqrt{x+y} \sqrt{z}}\right)}dzdydx=-8 i (-4+\pi ) \sqrt{\pi }
\end{multline}
\end{example}
\subsubsection{Derived using equation (\ref{eq_33682})}
\begin{example}
\begin{multline}
\int _0^{\frac{\pi }{4}}\int _0^{\frac{\pi }{4}}\int _0^{\frac{\pi }{4}}\frac{\sec (x) \sec ^2(z)}{\sqrt{\cos (2 x)} \left(\pi ^2+\log ^2\left(\frac{2 \sin (2 y) \sin (z) \tan (2 x)}{(\cos
   (y)+\sin (y))^2 (\cos (z)-\sin (z))}\right)\right) (\cos (y)+\sin (y)) \sqrt{\sin (2 y)}}dzdydx\\
=\frac{-\psi ^{(0)}\left(\frac{5}{8}\right)+\psi ^{(0)}\left(\frac{7}{8}\right)}{4 \sqrt{2}}
\end{multline}
\end{example}
\begin{example}
\begin{multline}
\int _0^{\frac{\pi }{4}}\int _0^{\frac{\pi }{4}}\int _0^{\frac{\pi }{4}}\frac{\log \left(\frac{2 \sin (2 y) \sin (z) \tan (2 x)}{(\cos (y)+\sin (y))^2 (\cos (z)-\sin (z))}\right) \sec (x) \sec
   ^2(z)}{\sqrt{\cos (2 x)} \left(\pi ^2+\log ^2\left(\frac{2 \sin (2 y) \sin (z) \tan (2 x)}{(\cos (y)+\sin (y))^2 (\cos (z)-\sin (z))}\right)\right) (\cos (y)+\sin (y)) \sqrt{\sin (2
   y)}}dzdydx=0
\end{multline}
\end{example}
\subsubsection{Derived using equation (\ref{eq_gr})}
\begin{example}
\begin{equation}
\int _0^{\frac{\pi }{2}}\int _0^1\int _0^{\frac{\pi }{2}}\int _0^1\frac{1}{\pi ^2+\log ^2\left(\frac{\sqrt{\log (y)} \sin (x)}{\sqrt{\log (r)} \sin
   (z)}\right)}drdzdydx=1-\frac{\pi }{4}
\end{equation}
\end{example}
\begin{example}
\begin{multline}
\int _0^{\frac{\pi }{2}}\int _0^1\int _0^{\frac{\pi }{2}}\int _0^1\frac{-\pi ^2+\log ^2\left(\frac{\sqrt{\log (y)} \sin (x)}{\sqrt{\log (r)} \sin (z)}\right)}{\left(\pi ^2+\log
   ^2\left(\frac{\sqrt{\log (y)} \sin (x)}{\sqrt{\log (r)} \sin (z)}\right)\right)^2}drdzdydx=2 (-1+C)
\end{multline}
\end{example}
\subsubsection{Derived using equation (\ref{eq_31924_1})}
\begin{example}
\begin{equation}
\int _0^{\infty }\int _1^{\infty }\frac{1+(-4+x) x}{(1+x)^4 \sqrt{-1+y} y \log (x (1-y))}dydx=\frac{i \psi ^{(3)}\left(\frac{3}{4}\right)}{64 \pi ^2}
\end{equation}
\end{example}
\begin{example}
\begin{multline}
\int _0^{\infty }\int _1^{\infty }\frac{1}{(1+x)^2 \sqrt{-1+y} y \left(\pi ^2+\log ^2(x (-1+y))\right)}dydx=-\frac{2 C}{\pi }+\frac{\pi }{4}
\end{multline}
\end{example}
\begin{example}
\begin{equation}
\int _0^{\infty }\int _1^{\infty }\frac{1}{(1+x)^2 \sqrt{-1+y} y \log (x (1-y))}dydx=-\frac{1}{4} i \left(-8 C+\pi ^2\right)
\end{equation}
\end{example}
\subsubsection{Derived using equation (\ref{eq_31851})}
\begin{example}
\begin{equation}
\int _0^{\infty }\int _0^{\infty }\int _0^{\infty }\int _0^{\infty }\frac{e^{-t^6-y^6} t^5 y^5 \text{Ai}\left(\frac{x}{2}\right)^2 \text{Ai}\left(\frac{z}{3}\right)^2}{\sqrt{x} \sqrt{z}
   \left(\pi ^2+\log ^2\left(\frac{3 x y^2}{2 t^2 z}\right)\right)}dtdzdydx=\frac{\log (2)}{216 \sqrt{6} \pi ^2}
\end{equation}
\end{example}
\begin{example}
\begin{multline}
\int _0^{\infty }\int _0^{\infty }\int _0^{\infty }\int _0^{\infty }\frac{e^{-t^6-y^6} t^5 y^5 \text{Ai}\left(\frac{x}{2}\right)^2 \text{Ai}\left(\frac{z}{3}\right)^2}{\sqrt{x} \sqrt{z} \log
   ^2\left(-\frac{3 x y^2}{2 t^2 z}\right)}dtdzdydx=-\frac{1}{5184 \sqrt{6}}
\end{multline}
\end{example}
\begin{example}
\begin{multline}
\int _0^{\infty }\int _0^{\infty }\int _0^{\infty }\int _0^{\infty }\frac{e^{-t^6-y^6} t^5 y^5 \text{Ai}\left(\frac{x}{2}\right)^2 \text{Ai}\left(\frac{z}{3}\right)^2}{\sqrt{x} \sqrt{z} \log
   ^3\left(-\frac{3 x y^2}{2 t^2 z}\right)}dtdzdydx=\frac{i \zeta (3)}{1152 \sqrt{6} \pi ^3}
\end{multline}
\end{example}
\subsubsection{Derived using equations (\ref{eq_gendi}) and (\ref{eq_6414})}
\begin{example}
\begin{multline}
\int _0^1\int _0^1\frac{\log \left(-\frac{\sqrt{-1+\frac{1}{y}}}{x}\right)+x \log \left(-x
   \sqrt{-1+\frac{1}{y}}\right)}{\sqrt{x} \left(1+x^2\right) \sqrt[4]{-1+\frac{1}{y}} \log
   \left(-\frac{\sqrt{-1+\frac{1}{y}}}{x}\right) \log \left(-x \sqrt{-1+\frac{1}{y}}\right)}dydx\\
=-\frac{1}{48} \pi 
   (\pi +12 i \log (2))
\end{multline}
\end{example}
\begin{example}
\begin{multline}
\int _0^1\int _0^1\frac{\log \left(\frac{-1+y}{x y}\right) \log \left(\log \left(\frac{-1+y}{x
   y}\right)\right)+\log \left(\frac{x (-1+y)}{y}\right) \log \left(\log \left(\frac{x
   (-1+y)}{y}\right)\right)}{\sqrt{x} (1+x)}dydx\\
=i \pi ^2 \left(1+\frac{i \pi }{2}+\log (4)+\log (\pi )+4 \left(i
   \pi -\log (4)-\frac{1}{2} \log (2 \pi )+\text{log$\Gamma $}\left(-\frac{1}{4}\right)\right)\right)
\end{multline}
\end{example}
\begin{example}
\begin{equation}
\int _0^1\int _0^1\frac{\frac{1}{\log \left(\frac{1}{x+\frac{x}{-1+y}}\right)}+\frac{1}{\log
   \left(x-\frac{x}{y}\right)}}{\sqrt{x} (1+x)}dydx=\frac{\pi ^2-8 C}{4 i}
\end{equation}
\end{example}
\begin{example}
\begin{multline}
\int _0^1\int _0^1\frac{y^s \left(\frac{x \log \left(\log
   \left(-\frac{\sqrt{-1+y^{-1-s}}}{x}\right)\right)}{\log \left(-\frac{\sqrt{-1+y^{-1-s}}}{x}\right)}+\frac{\log
   \left(\log \left(-x \sqrt{-1+y^{-1-s}}\right)\right)}{\log \left(-x \sqrt{-1+y^{-1-s}}\right)}\right)}{\sqrt{x}
   \left(1+x^2\right) \sqrt[4]{-1+y^{-1-s}}}dydx\\
=\frac{-i \pi ^3+12 i \pi  \log (2) (2 \gamma -\log (8)-2 \log (\pi
   ))+\pi ^2 (2+\log (4096)-2 \log (\pi ))+12 \zeta '(2)}{96 (1+s)}
\end{multline}
\end{example}
\begin{example}
\begin{multline}
\int _0^1\int _0^1\int _0^1\frac{z^s \sqrt[4]{-1+z^{-1-s}} \log ^k\left(\frac{a \left(-1+z^{-1-s}\right) \log
   (x)}{\log (y)}\right)}{\sqrt[4]{\log (x)} \log ^{\frac{3}{4}}(y)}dzdydx\\
=\frac{i^k 2^{-1+3 k} \pi ^{1+k}}{1+s} \left(i k
   \left(\zeta \left(1-k,\frac{1}{4}-\frac{i \log (a)}{8 \pi }\right)-\zeta \left(1-k,\frac{3}{4}-\frac{i \log
   (a)}{8 \pi }\right)\right)\right. \\ \left.
+2 \pi  \left(-\zeta \left(-k,\frac{1}{4}-\frac{i \log (a)}{8 \pi }\right)+\zeta
   \left(-k,\frac{3}{4}-\frac{i \log (a)}{8 \pi }\right)\right)\right)
\end{multline}
\end{example}
\begin{example}
\begin{multline}
\int _0^1\int _0^1\int _0^1\frac{z^s \sqrt[4]{-1+z^{-1-s}} \log \left(\log \left(\frac{a
   \left(-1+z^{-1-s}\right) \log (x)}{\log (y)}\right)\right)}{\sqrt[4]{\log (x)} \log
   ^{\frac{3}{4}}(y)}dzdydx\\
=\frac{\pi }{4 (1+s)} \left(\pi  \left(-i \pi +\log \left(\frac{1}{64 \pi ^2}\right)+4 \log (-6 \pi
   -i \log (a))-4 \log (-2 \pi -i \log (a))\right)\right. \\ \left.
+4 \pi  \text{log$\Gamma $}\left(-\frac{3}{4}-\frac{i \log (a)}{8
   \pi }\right)-4 \pi  \text{log$\Gamma $}\left(-\frac{1}{4}-\frac{i \log (a)}{8 \pi }\right)\right. \\ \left.
-2 i \psi
   ^{(0)}\left(\frac{1}{4}-\frac{i \log (a)}{8 \pi }\right)+2 i \psi ^{(0)}\left(\frac{3}{4}-\frac{i \log (a)}{8 \pi
   }\right)\right)
\end{multline}
\end{example}
\begin{example}
\begin{multline}
\int _0^1\int _0^1\int _0^1\frac{z^s \sqrt[4]{-1+z^{-1-s}} \log \left(\frac{a \left(-1+z^{-1-s}\right) \log
   (x)}{\log (y)}\right) \log \left(\log \left(\frac{a \left(-1+z^{-1-s}\right) \log (x)}{\log
   (y)}\right)\right)}{\sqrt[4]{\log (x)} \log ^{\frac{3}{4}}(y)}dzdydx\\
=\frac{\pi ^2}{4
   (1+s)} \left(-8-4 i \pi +\log
   \left(\frac{1}{16777216 \pi ^8}\right)+\log (a) \left(-i \pi +\log \left(\frac{1}{64 \pi ^2}\right)\right)\right. \\ \left.
+16
   \log (-6 \pi -i \log (a))-16 \log (-2 \pi -i \log (a))+16 \text{log$\Gamma $}\left(-\frac{3}{4}-\frac{i \log
   (a)}{8 \pi }\right)\right. \\ \left.
-16 \text{log$\Gamma $}\left(-\frac{1}{4}-\frac{i \log (a)}{8 \pi }\right)+32 i \pi 
   \left(\zeta'\left(-1,\frac{1}{4}-\frac{i \log (a)}{8 \pi
   }\right)-\zeta'\left(-1,\frac{3}{4}-\frac{i \log (a)}{8 \pi }\right)\right)\right)
\end{multline}
\end{example}
\begin{example}
\begin{multline}
\int _0^1\int _0^1\int _0^1\frac{z^s \sqrt[8]{-1+z^{-1-s}} \log \left(\log
   \left(-\frac{\left(-1+z^{-1-s}\right) \log (x)}{\log (y)}\right)\right)}{\log ^{\frac{3}{8}}(x) \log
   ^{\frac{5}{8}}(y)}dzdydx\\
=\frac{\pi }{12 \sqrt{2} (1+s)} \left((16+16 i) \, _2F_1\left(\frac{3}{4},1;\frac{7}{4};i\right)-3 i \pi 
   \left(\pi -2 i \log (4 \pi )\right.\right. \\ \left.\left.
+(2+2 i) \left(\log \left(\Gamma \left(\frac{3}{16}\right)\right)-\log \left(2 \Gamma
   \left(\frac{11}{16}\right)\right)-i \log \left(\frac{2 \Gamma \left(\frac{15}{16}\right)}{\Gamma
   \left(\frac{7}{16}\right)}\right)\right)\right)\right)
\end{multline}
\end{example}
\begin{example}
\begin{multline}
\int _0^1\int _0^1\int _0^1\frac{z^s \left(-1+z^{-1-s}\right)^{-m} \log ^{-\frac{1}{2}-m}(x) \log
   ^{-\frac{1}{2}-m}(y) \left(-\left(-1+z^{-1-s}\right)^{2 m} \log ^{2 m}(x)+\log ^{2 m}(y)\right)}{\log
   \left(\frac{\left(-1+z^{-1-s}\right) \log (x)}{\log (y)}\right)}dzdydx\\
=\frac{1}{4 (1+s)}\left(-8 m \pi  \left(\coth
   ^{-1}\left(e^{2 i m \pi }\right)+\tanh ^{-1}\left(e^{2 i m \pi }\right)\right)\right. \\ \left.
+i e^{-2 i m \pi } \left(\Phi
   \left(e^{-4 i m \pi },2,\frac{1}{2}\right)-e^{4 i m \pi } \Phi \left(e^{4 i m \pi
   },2,\frac{1}{2}\right)\right)\right)
\end{multline}
\end{example}
\begin{example}
\begin{multline}
\int _0^1\int _0^1\int _0^1\frac{z^s \left(-\sqrt{-1+z^{-1-s}} \sqrt{\log (x)}+\sqrt{\log
   (y)}\right)}{\sqrt[4]{-1+z^{-1-s}} \log ^{\frac{3}{4}}(x) \log ^{\frac{3}{4}}(y) \log
   \left(\frac{\left(-1+z^{-1-s}\right) \log (x)}{\log (y)}\right)}dzdydx=\frac{2 C}{1+s}
\end{multline}
\end{example}
\begin{example}
\begin{multline}
\int _0^1\int _0^1\int _0^1\frac{z^s \left(\sqrt{-1+z^{-1-s}} \sqrt{\log (x)}-\sqrt{\log
   (y)}\right)}{\sqrt[4]{-1+z^{-1-s}} \log ^{\frac{3}{4}}(x) \log ^{\frac{3}{4}}(y) \log
   \left(\frac{\left(-1+z^{-1-s}\right) \log (x)}{\log (y)}\right)}dzdydx=-\frac{2 C}{1+s}
\end{multline}
\end{example}
\begin{example}
\begin{multline}
\int _0^1\int _0^1\int _0^1\frac{z^s \sqrt[4]{-1+z^{-1-s}} \log \left(\log
   \left(-\frac{\left(-1+z^{-1-s}\right) \log (x)}{\log (y)}\right)\right)}{\sqrt[4]{\log (x)} \log
   ^{\frac{3}{4}}(y)}dzdydx\\
=\frac{\pi  \left(-4 i \sqrt{2} \log \left(\cot \left(\frac{\pi }{8}\right)\right)+\pi 
   \left(2 i \sqrt{2}-i \pi -2 \log (4 \pi )+4
   \Phi'\left(-1,0,\frac{3}{4}\right)\right)\right)}{4 (1+s)}
\end{multline}
\end{example}
\begin{example}
\begin{multline}
\int _0^1\int _0^1\int _0^1\frac{z^s \log ^k\left(\frac{a \left(-1+z^{-1-s}\right) \log (x)}{\log
   (y)}\right)}{\sqrt{\log (x)} \sqrt{\log (y)}}dzdydx=\frac{(4 i)^k k \pi ^{1+k} \zeta
   \left(1-k,\frac{1}{2}-\frac{i \log (a)}{4 \pi }\right)}{1+s}
\end{multline}
\end{example}
\begin{example}
\begin{multline}
\int _0^1\int _0^1\int _0^1\frac{z^s}{\sqrt{\log (x)} \sqrt{\log (y)} \left(a+i b+\log
   \left(\frac{\left(-1+z^{-1-s}\right) \log (x)}{\log (y)}\right)\right)}dzdydx=\frac{i \psi ^{(1)}\left(\frac{-i
   a+b+2 \pi }{4 \pi }\right)}{4 (1+s)}
\end{multline}
\end{example}
\begin{example}
\begin{multline}
\int _0^1\int _0^1\int _0^1\frac{z^s}{\sqrt{\log (x)} \sqrt{\log (y)} \left(a^2+b^2-2 i b \log
   \left(\frac{\left(-1+z^{-1-s}\right) \log (x)}{\log (y)}\right)-\log ^2\left(\frac{\left(-1+z^{-1-s}\right) \log
   (x)}{\log (y)}\right)\right)}dzdydx\\
=-\frac{i \left(-\psi ^{(1)}\left(\frac{-i a+b+2 \pi }{4 \pi }\right)+\psi
   ^{(1)}\left(\frac{i a+b+2 \pi }{4 \pi }\right)\right)}{8 a (1+s)}
\end{multline}
\end{example}
\begin{example}
\begin{multline}
\int _0^1\int _0^1\int _0^1\frac{z^s}{\sqrt{\log (x)} \sqrt{\log (y)} \left(3 \pi ^2+4 i \pi  \log
   \left(\frac{\left(-1+z^{-1-s}\right) \log (x)}{\log (y)}\right)-\log ^2\left(\frac{\left(-1+z^{-1-s}\right) \log
   (x)}{\log (y)}\right)\right)}dzdydx\\
=\frac{2 (-1+C)}{\pi  (1+s)}
\end{multline}
\end{example}
\begin{example}
\begin{multline}
\int _0^1\int _0^1\int _0^1\frac{z^s \log \left(\frac{a \left(-1+z^{-1-s}\right) \log (x)}{\log (y)}\right)
   \log \left(\log \left(\frac{a \left(-1+z^{-1-s}\right) \log (x)}{\log (y)}\right)\right)}{\sqrt{\log (x)}
   \sqrt{\log (y)}}dzdydx\\
=\frac{\pi }{2 (1+s)} \left(\log (a) (-i \pi -2 (1+\log (4)+\log (\pi )))\right. \\ \left.
-8 i \pi  \left(-\frac{1}{2}
   \log (2 \pi )+\log \left(-\frac{1}{2}-\frac{i \log (a)}{4 \pi }\right)+\text{log$\Gamma
   $}\left(-\frac{1}{2}-\frac{i \log (a)}{4 \pi }\right)\right)\right)
\end{multline}
\end{example}
\begin{example}
\begin{multline}
\int _0^1\int _0^1\int _0^1\frac{z^s \left(\log \left(-a+i b+\log \left(\frac{\left(-1+z^{-1-s}\right) \log
   (x)}{\log (y)}\right)\right)+\log \left(a+i b+\log \left(\frac{\left(-1+z^{-1-s}\right) \log (x)}{\log
   (y)}\right)\right)\right)}{\sqrt{\log (x)} \sqrt{\log (y)}}dzdydx\\
=-\frac{\pi  \left(i \pi +2 \log (4 \pi )+\psi
   ^{(0)}\left(\frac{-i a+b+2 \pi }{4 \pi }\right)+\psi ^{(0)}\left(\frac{i a+b+2 \pi }{4 \pi
   }\right)\right)}{1+s}
\end{multline}
\end{example}
\begin{example}
\begin{multline}
\int _0^1\int _0^1\int _0^1\frac{z^s \log \left(\frac{-a+i b+\log \left(\frac{\left(-1+z^{-1-s}\right) \log
   (x)}{\log (y)}\right)}{a+i b+\log \left(\frac{\left(-1+z^{-1-s}\right) \log (x)}{\log
   (y)}\right)}\right)}{\sqrt{\log (x)} \sqrt{\log (y)}}dzdydx\\=\frac{\pi  \left(\psi ^{(0)}\left(\frac{-i a+b+2 \pi
   }{4 \pi }\right)-\psi ^{(0)}\left(\frac{i a+b+2 \pi }{4 \pi }\right)\right)}{1+s}
\end{multline}
\end{example}
\begin{example}

\end{example}
\subsection{Products}
\begin{example}
\begin{multline}
\prod _{p=1}^{\infty } \left(1-\frac{4 i p \pi }{2 i p \pi +a \alpha
   }\right)^{\frac{i \pi  }{2 \alpha ^2}\left(\coth \left(\frac{p \pi ^2}{\alpha
   ^2}\right)-1\right)} \left(\frac{2 a}{a+2 p \alpha
   }-1\right)^{\frac{1}{2} \left(1-\coth \left(p \alpha
   ^2\right)\right)}\\
=\frac{a^{\frac{i \left(\pi +i \alpha ^2\right)}{2 \alpha
   ^2}} e^{\frac{-4 \alpha +a \pi  (\pi -i a \alpha )}{4 a \alpha ^2}} \alpha
   ^{\frac{i \pi +2 a \alpha +\alpha ^2}{2 \alpha ^2}} \Gamma \left(1+\frac{a}{2
   \alpha }\right)}{\pi ^{\frac{a+\alpha }{2 \alpha }} \Gamma \left(1-\frac{i a
   \alpha }{2 \pi }\right)^{\frac{i \pi }{\alpha ^2}}}
\end{multline}
\end{example}
\begin{example}
\begin{multline}
\Gamma \left(\frac{7}{4}\right)=-\frac{\left(\pi ^{5/4} \Gamma
   \left(\frac{\pi -12 i}{\pi }\right)^{\frac{i \pi }{16}}\right) }{(-1)^{5/8}
   2^{\frac{1}{32} (112+3 i \pi )} 3^{\frac{1}{32} i (16 i+\pi )}
   e^{\frac{1}{192} \left(-8+3 \pi ^2\right)}}\\
   \prod\limits
   _{p=1}^{\infty } \left(-1+\frac{6}{3+4 p}\right)^{\frac{1}{2} (1-\coth (16
   p))} \left(1+\frac{2 p \pi }{12 i-p \pi }\right)^{\frac{1}{32} i \pi 
   \left(-1+\coth \left(\frac{p \pi ^2}{16}\right)\right)}
\end{multline}
\end{example}
\begin{example}
\begin{multline}
\prod_{p=1}^{\infty}\frac{
   (a+2 i p \beta )^{i e^{-\frac{p \pi  \beta }{2 \alpha }} \beta  \left(\coth \left(p \beta ^2\right)-1\right)}}{(a-2 i p \beta )^{i
   e^{\frac{p \pi  \beta }{2 \alpha }} \beta  \left(\coth \left(p \beta ^2\right)-1\right)}}\left(\frac{(a-2 p \alpha )^{(-1)^p}}{a+2 p \alpha }\right)^{i^p \alpha  \left(\coth \left(p \alpha ^2\right)-1\right)}\\
=2^{i \left((3+3 i) \alpha +\beta -\beta 
   \coth \left(\frac{\pi  \beta }{4 \alpha }\right)\right)} a^{\frac{i \pi }{2 \alpha }+\alpha -i \beta } e^{\frac{1}{2} \pi  \beta 
   \left(\coth \left(\frac{\pi  \beta }{4 \alpha }\right)-1\right)} \alpha ^{(-1+i) \alpha } \beta ^{-i \beta  \left(\coth
   \left(\frac{\pi  \beta }{4 \alpha }\right)-1\right)}\\
 \exp \left(\frac{2}{a}+2 i e^{-\frac{\pi  \beta }{2 \alpha }} \beta 
   \Phi'\left(e^{-\frac{\pi  \beta }{2 \alpha }},0,1-\frac{i a}{2 \beta }\right)\right) \left(\frac{\Gamma
   \left(\frac{1}{8} \left(4+\frac{a}{\alpha }\right)\right)}{\Gamma \left(1+\frac{a}{8 \alpha }\right) \left(\frac{\Gamma
   \left(\frac{1}{8} \left(2+\frac{a}{\alpha }\right)\right)}{\Gamma \left(\frac{1}{8} \left(6+\frac{a}{\alpha
   }\right)\right)}\right)^i}\right)^{2 \alpha }
\end{multline}
\end{example}
\begin{example}
\begin{multline}
\prod_{p=1}^{\infty}\frac{(6 i p+5 \pi )^{i e^{-\frac{3 p}{2}} \text{csch}(p)} }{(5 \pi -6 i p)^{i e^{p/2}
   (\coth (p)-1)}}\left(\frac{(5-6 p)^{(-1)^p} \left(\frac{3}{\pi
   }\right)^{1-(-1)^p}}{5+6 p}\right)^{i^p \pi  \left(-1+\coth \left(p \pi ^2\right)\right)}\\
=2^{i \left(1+(3+3 i) \pi -\coth \left(\frac{1}{4}\right)\right)} 3^{\frac{i}{2}-\frac{2
   i}{-1+\sqrt{e}}-(1+2 i) \pi -2 \psi _{\sqrt{e}}^{(0)}(1-i \pi )+2 \psi _{\sqrt{e}}^{(0)}(1+i \pi )}
   5^{-\frac{i}{2}+\pi } \pi ^{\frac{1}{2} i (-1+2 \pi )}\\
 \exp \left(\frac{6}{5 \pi }+\frac{1}{2} \pi  \left(-1+\coth
   \left(\frac{1}{4}\right)\right)+\frac{2 i \Phi'\left(\frac{1}{\sqrt{e}},0,1-\frac{5 i \pi
   }{6}\right)}{\sqrt{e}}\right)\\
    \left(\frac{\Gamma \left(\frac{23}{24}\right)}{\Gamma
   \left(\frac{11}{24}\right)}\right)^{2 i \pi } \left(\frac{\Gamma \left(\frac{17}{24}\right)}{\Gamma
   \left(\frac{29}{24}\right)}\right)^{2 \pi }
\end{multline}
\end{example}
\begin{example}
\begin{multline}
\prod _{\beta =1}^q \sum _{p=1}^{\infty } p \left(-\beta ^4+\pi ^2 \coth \left(\frac{\pi ^2 p}{\beta
   ^2}\right)+\beta ^4 \coth \left(\beta ^2 p\right)-\pi ^2\right)\\
=\frac{12^{-q} }{\pi ^3}\sin \left(\pi  \sqrt{3-i \sqrt{\pi
   ^2-9}}\right) \sin \left(\pi  \sqrt{3+i \sqrt{\pi ^2-9}}\right)\\ \Gamma \left(q-\sqrt{3-i \sqrt{-9+\pi ^2}}+1\right)
  \Gamma \left(q+\sqrt{3-i \sqrt{-9+\pi ^2}}+1\right)\\ \Gamma \left(q-\sqrt{3+i \sqrt{-9+\pi ^2}}+1\right) 
\Gamma
   \left(q+\sqrt{3+i \sqrt{-9+\pi ^2}}+1\right)
\end{multline}
\end{example}
\begin{example}
Triple product involving Euler's polynomial.
\begin{multline}
\prod_{j=0}^{\infty}\prod_{l=0}^{\infty}\prod_{q=0}^{\infty}\exp \left(\frac{E_q(\alpha ) 2^{-2 (j+l)+q-1} \left(-\frac{1}{x}\right)^{2 (j+l)+q} \Gamma (2 j+2 l+q+1)}{x
   \Gamma (l+1) \Gamma (q+1) \Gamma (2 j+l+2)}\right)\\
=\frac{2 \Gamma \left(\frac{1}{4} (x+2 \alpha +3)\right)}{\sqrt{2
   \alpha +x-1} \Gamma \left(\frac{1}{4} (x+2 \alpha +1)\right)}
\end{multline}
\end{example}
\begin{example}
Infinite product involving Euler's constant $\gamma$.
\begin{equation}
\prod _{p=1}^{\infty } \left(\frac{(2 p+1)^{\frac{1}{2 p+1}}}{(2 p-1)^{\frac{1}{2 p-1}}}\right)^{\frac{(-1)^p}{\pi }}=\frac{3 e^{-\frac{\gamma }{2}} }{\sqrt{2 \pi } }\frac{\Gamma
   \left(-\frac{3}{4}\right)}{\Gamma \left(-\frac{1}{4}\right)}
\end{equation}
\end{example}
\begin{example}
\begin{multline}
\prod _{k=0}^{n-1} \left(1-i \cot \left(\frac{k \pi }{n}+i x\right)-\coth \left(\frac{k \pi }{n}+i
   x\right)\right)\\
   =\frac{(-1)^n \left(-e^{2 i x};e^{\frac{2 \pi }{n}}\right){}_n \sin \left(\frac{1}{2} n (\pi +2 i
   x)\right)}{\left(e^{2 i x};e^{\frac{2 \pi }{n}}\right){}_n \sinh (n x)}
\end{multline}
\end{example}
\begin{example}
\begin{multline}
\prod_{p=1}^{\infty}\left(\frac{a+i p-1}{a-i p-1}\right)^{i (\coth (\pi  p)-1)} \left(\frac{a-p-1}{a+p-1}\right)^{\coth (\pi  p)-1}\\
=-\frac{(2 \pi )^{1-i} (-i (a-1))^{2 i} (a-1)^{-1-i} e^{\frac{1}{\pi  (a-1)}+\pi 
   \left(-\frac{1}{2}+i a\right)} \Gamma (-i (a-1))^{2 i}}{\Gamma (a-1)^2}
\end{multline}
\end{example}
\begin{example}
The exponential of the incomplete gamma function.
\begin{equation}
\prod _{p=1}^{\infty } \exp (\Gamma (0,i (p-i m) \pi )+\Gamma (0,-i (i
   m+p) \pi ))=\frac{1}{2} (\coth (\pi  m)+1) e^{-\Gamma (0,m \pi )}
\end{equation}
\end{example}
\begin{example}
\begin{multline}
\prod _{p=1}^{\infty }4^{-\frac{1}{2} i \pi  (2 a-1) (-1)^p} e^{\frac{m (-1)^{p+1}}{m^2-p^2}} ((m-p) (m+p))^{-\frac{1}{2} i \pi  (2 a-1) (-1)^p}\\
=e^{\frac{1}{2} \left(\frac{1}{m}-\pi  \csc (\pi 
   m)\right)} \left(\frac{\sin ^2\left(\frac{\pi  m}{2}\right) \csc (\pi  m)}{m}\right)^{\frac{1}{2} i \pi  (1-2 a)}
\end{multline}
\end{example}
\begin{example}
\begin{multline}
\prod _{p=1}^{\infty } e^{\frac{(-1)^{p+1}}{2 \left(\frac{1}{4}-p^2\right)}} 4^{-\frac{1}{2} i \pi  (2 a-1) (-1)^p} \left(\left(\frac{1}{2}-p\right)
   \left(p+\frac{1}{2}\right)\right)^{-\frac{1}{2} i \pi  (2 a-1) (-1)^p}=e^{\frac{2-\pi }{2}}
\end{multline}
\end{example}
\begin{example}
\begin{multline}
\prod _{k=1}^n 2 e^{\Gamma \left(0,\pi  \left(\frac{2
   k}{n}+\frac{x}{\pi }\right)\right)} \prod _{p=1}^{\infty } \exp
   \left(\Gamma \left(0,i \pi  \left(p-i \left(\frac{2
   k}{n}+\frac{x}{\pi }\right)\right)\right)+\Gamma \left(0,-i \pi 
   \left(p+i \left(\frac{2 k}{n}+\frac{x}{\pi
   }\right)\right)\right)\right)\\
   =\prod _{k=1}^n \left(\coth \left(\pi
    \left(\frac{2 k}{n}+\frac{x}{\pi }\right)\right)+1\right)
\end{multline}
\end{example}
\begin{example}
\begin{multline}
\prod _{p=0}^{\infty } \exp \left(\frac{(-1)^p e^{-i \pi  (2 p+1) z} \left(\Gamma (-1,-i (2 p+1) \pi  z)+e^{2 i \pi  (2 p+1) z} \Gamma (-1,i (2 p+1) \pi  z)\right)}{\pi  (2
   p+1)}\right)\\
=\frac{\sqrt{z} \Gamma \left(\frac{z}{2}+\frac{1}{4}\right)}{\sqrt{2} \Gamma \left(\frac{z}{2}+\frac{3}{4}\right)}
\end{multline}
\end{example}
\begin{example}
\begin{multline}
\prod _{p=0}^{\infty } \exp (\Gamma (-1,i p \pi -\log (q))+\Gamma (-1,-i
   \pi  (p+1)-\log (q)))=e^{\frac{2 q^2}{1-q^2}} \left(1-q^2\right)
\end{multline}
\end{example}
\begin{figure}[H]
\includegraphics[scale=0.5]{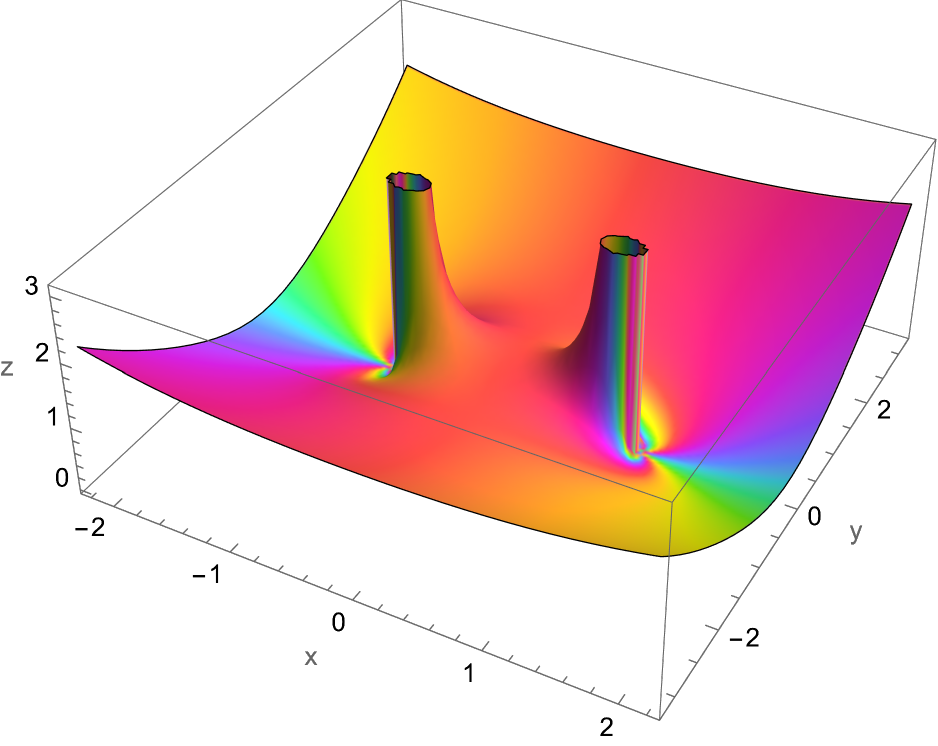}
\caption{Plot of  $f(q)=e^{\frac{2 q^2}{1-q^2}} \left(1-q^2\right)$, $q\in\mathbb{C}$.}
   \label{fig:fig2}
\end{figure}
\vspace{-6pt}
\begin{example}
\begin{multline}
\prod _{p=0}^{\infty } \exp \left((-1)^p e^{-2 i \pi  a (p+1)} \left(\frac{e^{2 i \pi  a (2 p+1)} \Gamma \left(0,\frac{1}{2} i a (4 p+1)
   \pi \right)}{4 p+1}+\frac{\Gamma \left(0,-\frac{1}{2} i a (4 p+3) \pi \right)}{4 p+3}\right)\right)\\
=\left(2^{1-i}
   a^{-\frac{1}{2}+\frac{i}{2}}\right)^{\frac{1}{2} \sqrt[4]{-1} \pi  i^{-a}} \left(\frac{\Gamma
   \left(\frac{a}{4}+\frac{5}{8}\right) \left(\frac{\Gamma \left(\frac{a}{4}+\frac{3}{8}\right)}{\Gamma
   \left(\frac{a}{4}+\frac{7}{8}\right)}\right)^i}{\Gamma \left(\frac{a}{4}+\frac{1}{8}\right)}\right)^{\frac{1}{2}
   \sqrt[4]{-1} \pi  i^{-a}}
\end{multline}
\end{example}
\begin{example}
\begin{multline}
\prod _{p=0}^{\infty } \exp \left(\frac{(-1)^p e^{-i a (2 p+3) x}}{\pi  (2 p+1)^2 (4 p (p+1)-3)} \left((2 p+3) (1-2 p)^2 e^{2 i a (2 p+1) x} \Gamma (-1,i a (2 p-1) x)\right.\right. \\ \left.\left.
-8 p (2 p+3) e^{2 i a (2 p+1) x} \Gamma (0,i a
   (2 p-1) x)-(2 p-1) (2 p+1)^2 \Gamma (0,-i a (2 p+3) x)\right)\right)\\
=e^{\frac{ C}{\pi  ia x}} \frac{\sqrt{2 \pi } }{\sqrt{a x}}\frac{\Gamma \left(\frac{2 a x+3 \pi }{4 \pi
   }\right)}{ \Gamma \left(\frac{2 a x+\pi }{4 \pi }\right)}
\end{multline}
\end{example}
\begin{example}
\begin{multline}
\prod _{p=0}^{\infty }\exp \left(\frac{(-1)^p  \left(\pi  (1-2 p)^2 z E_1(-i (2 p+3) \pi  z)+i (2 p+3) e^{2 i \pi  (2 p+1) z} E_2(i (2 p-1) \pi  z)\right)}{e^{i \pi  (2 p+3) z}(2 p+3) (\pi -2 \pi  p)^2
   z}\right)\\
=e^{-\frac{i (C-1)}{\pi ^2 z}}\frac{\sqrt{z}  }{\sqrt{2}}\frac{\Gamma \left(\frac{z}{2}+\frac{1}{4}\right)}{ \Gamma \left(\frac{z}{2}+\frac{3}{4}\right)}
\end{multline}
\end{example}
\begin{figure}[H]
\includegraphics[scale=0.5]{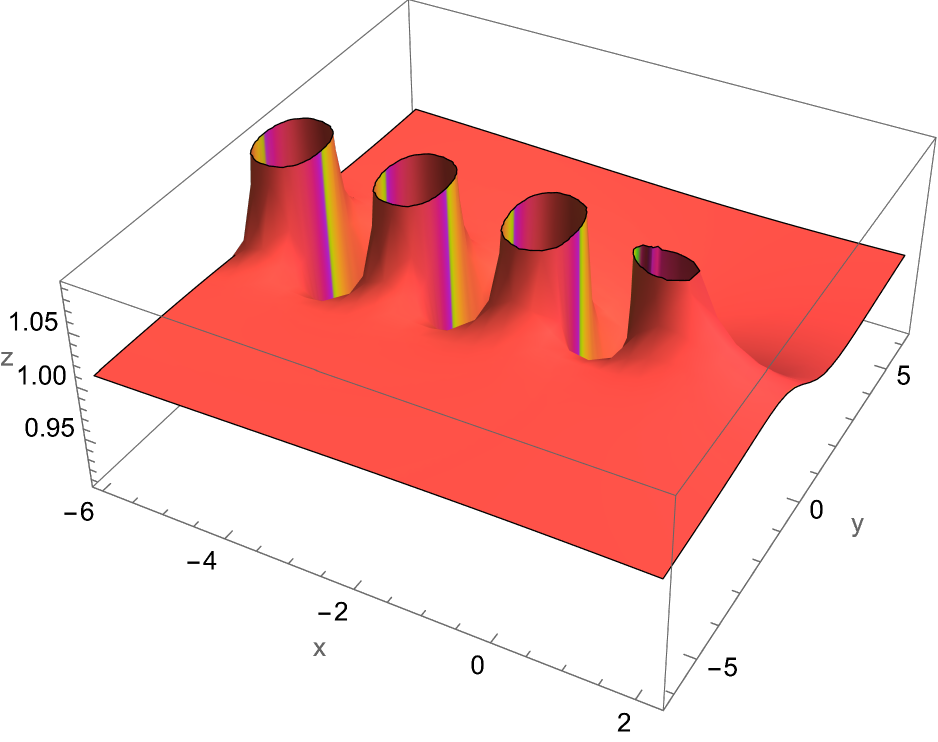}
\caption{Plot of  $f(z)=\frac{\sqrt{2} e^{-\frac{i C}{\pi ^2 z}} \Gamma \left(\frac{z}{2}+\frac{3}{4}\right)}{\sqrt{z} \Gamma
   \left(\frac{z}{2}+\frac{1}{4}\right)}$, $z\in\mathbb{C}$.}
   \label{fig:fig2}
\end{figure}
\vspace{-6pt}
\begin{example}
\begin{equation}
\prod _{p=0}^{\infty }\exp \left(\frac{12 }{\pi ^2}\left(\frac{E_1\left(\frac{1}{2} i (2 p-1) \pi \right)}{(1-2
   p)^2}+\frac{E_1\left(-\frac{1}{2} i (2 p+3) \pi \right)}{(2 p+3)^2}\right)\right)=\frac{128
   e^3}{A^{36}}
\end{equation}
\end{example}
\begin{example}
\begin{multline}
\prod _{p=0}^{\infty } \exp \left(\frac{2 i E_2\left(\frac{1}{2} i (2 p-1) \pi \right)}{\pi  (1-2 p)^2}-\frac{E_1\left(-\frac{1}{2} i (2
   p+3) \pi \right)}{2 p+3}\right)=2^{-i \pi } \pi ^{\frac{i \pi }{2}} e^{\frac{2 (C-1)}{\pi }}
\end{multline}
\end{example}
\begin{example}
\begin{multline}
\prod _{p=0}^{\infty }\exp \left(-\Gamma \left(0,\frac{i \pi  (i m+x+2 p x)}{2 x}\right)+\Gamma \left(0,\frac{i \pi  (i r+x+2 px)}{2 x}\right)\right. \\ \left.
-\Gamma \left(0,-\frac{\pi  (m+i (x+2 p x))}{2 x}\right)+\Gamma \left(0,-\frac{\pi  (r+i (x+2 px))}{2 x}\right)\right)\\
=e^{\frac{\pi  (m-r)}{2 x}} \frac{\cosh \left(\frac{m \pi }{2 x}\right)}{\cosh
   \left(\frac{\pi  r}{2 x}\right)}
\end{multline}
\end{example}
\begin{example}
\begin{multline}
\prod _{p=0}^{\infty }\exp \left(\Gamma \left(0,\frac{i \pi  \left(2 p x+x-\frac{i}{x}\right)}{2 x}\right)-\Gamma \left(0,\frac{i\pi  \left(2 p x+x+\frac{i}{x}\right)}{2 x}\right)\right. \\ \left.
+\Gamma \left(0,-\frac{\pi  \left(i (2 p
   x+x)-\frac{1}{x}\right)}{2 x}\right)-\Gamma \left(0,-\frac{\pi  \left(i (2 p x+x)+\frac{1}{x}\right)}{2
   x}\right)\right)\\
=e^{\frac{\pi }{x^2}}
\end{multline}
\end{example}
\begin{figure}[H]
\includegraphics[scale=0.5]{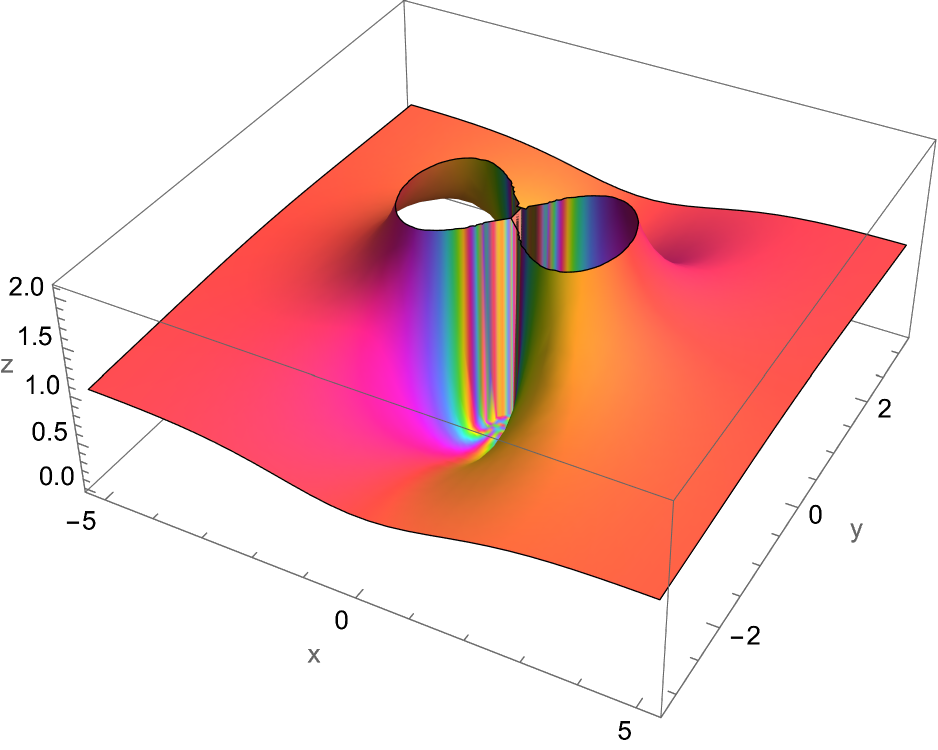}
\caption{Plot of  $f(x)=e^{\frac{\pi }{x^2}}$, $x\in\mathbb{C}$.}
   \label{fig:fig2}
\end{figure}
\vspace{-6pt}
\begin{example} 
\begin{multline}\label{eq:gamma}
\prod_{p=0}^{\infty}\exp \left(\frac{(-1)^p \left(\frac{i}{a \pi }\right)^{2 p} E_{2 p+1}(x)}{(a \pi ) (2
   p+1)}\right)=\frac{\Gamma \left(\frac{1}{2} (1+a \pi -x)\right) \Gamma \left(\frac{1}{2} (1+a \pi
   +x)\right)}{\Gamma \left(\frac{1}{2} (2+a \pi -x)\right) \Gamma \left(\frac{1}{2} (a \pi +x)\right)}
\end{multline}
\end{example}
\begin{example}
Asymptotic approximation where $Re(x)<0<Re(a)$. Similar forms and evaluations are given in the works of section (5.6) in \cite{dlmf}.
\begin{equation}
\prod_{p=0}^{\infty}\exp \left(\frac{(-1)^p \left(\frac{i}{a \pi }\right)^{2 p} E_{1+2 p}(2 x)}{a \pi  (1+2
   p)}\right) < \frac{x^{\frac{1}{2}-x}}{\sqrt{-x}}
\end{equation}
\end{example}
\begin{example}
\begin{multline}
\prod_{p=1}^{n}\left(\frac{\sin ^3\left(\frac{1}{2} (1+p) x\right)}{\sin ((1+p) x) \sin ^2\left(\frac{1}{4} (1+p)
   x\right)}\right)^{\frac{1}{p+1}} \left(\frac{\sin ((2+p) x) \sin ^2\left(\frac{1}{4} (2+p) x\right)}{\sin
   ^3\left(\frac{1}{2} (2+p) x\right)}\right)^{\frac{1}{p+2}}\\
=\frac{\csc \left(\frac{x}{2}\right) \sin
   ^{\frac{3}{2}}(x) \left(\frac{\tan \left(\frac{1}{4} (2+n) x\right)}{\tan \left(\frac{1}{2} (2+n)
   x\right)}\right)^{\frac{1}{2+n}}}{\sqrt{\sin (2 x)}}
\end{multline}
\end{example}
\begin{example}
\begin{multline}
\prod_{p=1}^{n}\left(\frac{\sin ((1+p) q)}{\sin (m (1+p))}\right)^{\frac{1}{1+p}} \left(\frac{\sin (m (2+p))}{\sin ((2+p)
   q)}\right)^{\frac{1}{2+p}}\\
=e^{-i (m-q)} \left(\frac{e^{i (2+n) (m-q)} \sin (m (2+n))}{\sin ((2+n)
   q)}\right)^{\frac{1}{2+n}} \sqrt{\frac{\cos (q) \sin (q)}{\sin (m) \cos (m)}}
\end{multline}
\end{example}
\begin{example}
\begin{multline}
\prod_{p=1}^{n}\frac{\sin ^{\frac{3}{1+p}}\left(\frac{1}{2} (1+p) x\right) \sin
   ^{\frac{2}{2+p}}\left(\frac{1}{4} (2+p) x\right) \sqrt[2+p]{\sin ((2+p)
   x)}}{\sin ^{\frac{3}{2+p}}\left(\frac{1}{2} (2+p) x\right) \sqrt[1+p]{\sin
   ((1+p) x)} \sin ^{\frac{2}{1+p}}\left(\frac{1}{4} (1+p) x\right)}\\
=\frac{\csc
   \left(\frac{x}{2}\right) \sin ^{\frac{3}{2}}(x) \left(\frac{\tan
   \left(\frac{1}{4} (2+n) x\right)}{\tan \left(\frac{1}{2} (2+n)
   x\right)}\right)^{\frac{1}{2+n}}}{\sqrt{\sin (2 x)}}
\end{multline}
\end{example}
\begin{example}
\begin{multline}
\prod _{p=1}^n  \sqrt{\frac{2+p}{1+p}}\left(\frac{(1+p)^{\frac{1}{1+p}}}{(2+p)^{\frac{1}{2+p}}}\right)^{2 a} \frac{
   \Gamma \left(\frac{2 a}{1+p}\right)}{\Gamma \left(\frac{2
   a}{2+p}\right)}=\sqrt{\frac{2+n}{2}}\left(\frac{2}{(2+n)^{\frac{2}{2+n}}}\right)^a \frac{ \Gamma (a)}{\Gamma
   \left(\frac{2 a}{2+n}\right)}
\end{multline}
\end{example}
\begin{example}
\begin{equation}
\frac{\prod\limits _{p=0}^n \Gamma \left(\frac{1}{2} \left(4^{-p} x+1\right)\right)^2}{\prod\limits _{p=0}^{2 n} \Gamma
   \left(\frac{1}{2} \left(2^{-p} x-1\right)\right)^{(-1)^p}}=\frac{\pi ^{n+\frac{1}{2}} (x-1)
   2^{\left(4^{-n}-2\right) x+2 n} \Gamma (x) \left(\frac{x}{4};\frac{1}{4}\right)_{\lfloor n\rfloor }}{\Gamma
   \left(2^{-2 n-1} x\right) \left(\frac{x}{2};\frac{1}{4}\right)_{\left\lfloor n-\frac{1}{2}\right\rfloor
   +1}}
\end{equation}
\end{example}
\begin{example}
\begin{multline}
\prod _{p=0}^{2 n} \left(\frac{\tan \left(2^{-1+p} r\right)}{\tan \left(2^{-1+p}
   m\right)}\right)^{\left(-\frac{1}{2}\right)^p} \prod _{p=0}^n \left(\frac{\tan \left(2^{-1+2 p}
   m\right)}{\tan \left(2^{-1+2 p} r\right)}\right)^{2^{1-2 p}}\\
=\left(\frac{\sin \left(\frac{m}{2}\right)}{\sin
   \left(\frac{r}{2}\right)}\right)^2 e^{i (m-r)} \left(\frac{e^{-i 4^n (m-r)} \sin \left(4^n r\right)}{\sin
   \left(4^n m\right)}\right)^{4^{-n}}
\end{multline}
\end{example}
\begin{example}
\begin{multline}
\prod _{p=0}^{2 n} \left(1-\frac{2}{q^{2^p}+1}\right)^{\left(-\frac{1}{2}\right)^p} \prod
   _{p=0}^n \left(\frac{2}{q^{4^p}-1}+1\right)^{2^{1-2 p}}=\frac{q^2 \left(1-q^{-2^{2
   n+1}}\right)^{4^{-n}}}{(q-1)^2}
\end{multline}
\end{example}
\begin{example}
\begin{multline}
\prod _{p=0}^n \pi ^{-(-1)^{4^p}} 2^{(-1)^{4^p} 4^{-p} \left(4 p z-3\ 4^p\right)} \left(2^{1-2 p}
   z-1\right)^{2 (-1)^{4^p}} \Gamma \left(2^{-2 p} z-\frac{1}{2}\right)^{2 (-1)^{4^p}}\\
 \prod _{p=0}^{2 n}
   \left(2^{-p} z-\frac{1}{2}\right)^{-(-1)^{p+2^p}} \Gamma \left(2^{-p}
   z-\frac{1}{2}\right)^{-(-1)^{p+2^p}}\\
=\frac{\pi  2^{-\frac{1}{9} 2^{1-2 n} \left(9 n+4^n-e^{2 i \pi  n} (3
   n+1)\right) z+n+1} \Gamma (z)}{\Gamma \left(z+\frac{1}{2}\right) \Gamma \left(4^{-n} z\right)}
\end{multline}
\end{example}
\begin{example}
\begin{equation}
\prod _{p=0}^{\infty } \left(\frac{2}{q^{2^p}-1}+1\right)^{\left(-\frac{1}{2}\right)^p}\prod
   _{p=0}^{\infty } \left(\frac{2}{q^{4^p}-1}+1\right)^{-2^{1-2 p}}=\left(1-\frac{1}{q}\right)^2
\end{equation}
\end{example}
\begin{figure}[H]
\includegraphics[scale=0.5]{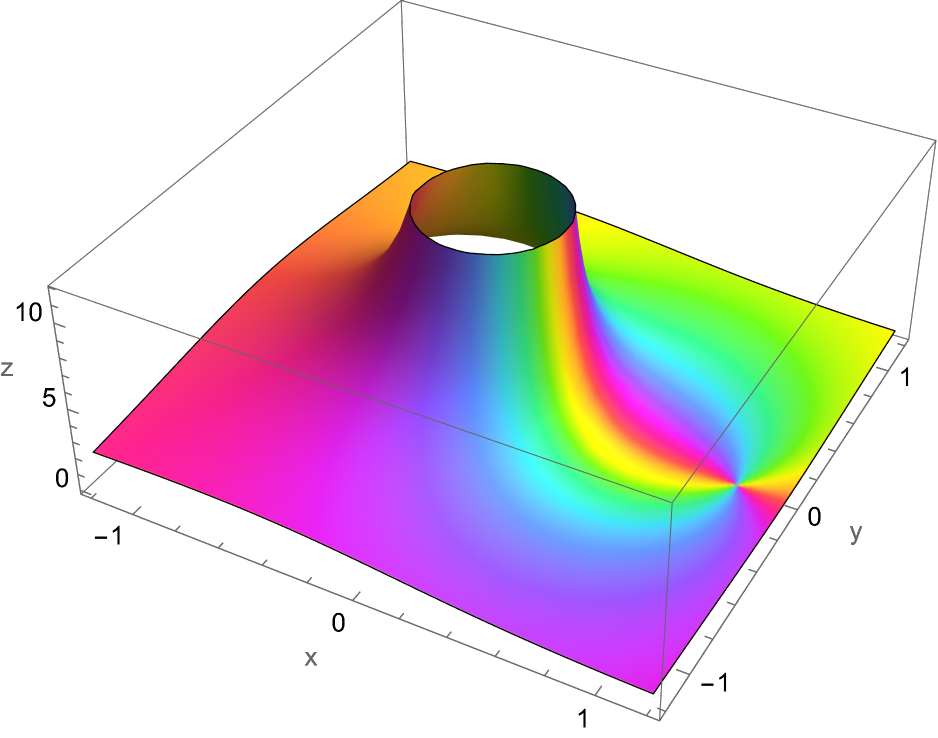}
\caption{Plot of  $f(q)=\left(1-\frac{1}{q}\right)^2$, $q\in\mathbb{C}$.}
   \label{fig:fig2}
\end{figure}
\vspace{-6pt}
\begin{example}
\begin{multline}
\prod _{p=0}^n e^{2 \left(\csc \left(2^{2 p-1} x\right)-\csc \left(4^p x\right)\right)}
   \left(1-\frac{1}{2} \sec ^2\left(4^{p-1} x\right)\right)^{a 2^{1-2 p}}\\
 \prod _{p=0}^{2 n} e^{(-1)^p
   \left(\csc \left(2^p x\right)-\csc \left(2^{p-1} x\right)\right)} \left(\tan \left(2^{p-1} x\right) \cot
   \left(2^{p-2} x\right)\right)^{a \left(-\frac{1}{2}\right)^p}\\
=\left(\frac{1}{1+e^{\frac{i x}{2}}}\right)^{2 a}
   \left(1+e^{i 4^n x}\right)^{a 4^{-n}} e^{\csc \left(\frac{x}{2}\right)-\csc \left(4^n x\right)}
\end{multline}
\end{example}
\begin{example}
\begin{multline}
\prod _{p=0}^n \exp \left(\frac{4 i q^{4^p} \left(\left(q^{4^p}-1\right)
   q^{4^p}+1\right)}{q^{4^{p+1}}-1}\right) \left(\frac{q^{2^{2 p+1}}+1}{\left(q^{4^p}+1\right)^2}\right)^{2^{1-2 p}
   z}\\
 \prod _{p=0}^{2 n} \exp \left(-\frac{2 i (-1)^p q^{2^p} \left(\left(q^{2^p}-1\right)
   q^{2^p}+1\right)}{q^{2^{p+2}}-1}\right) \left(\frac{2
   q^{2^p}}{q^{2^{p+1}}+1}+1\right)^{\left(-\frac{1}{2}\right)^p z}\\
=\left(\frac{1}{q+1}\right)^{2 z} \exp \left(-2 i
   \left(\frac{q^{2^{2 n+1}}}{q^{4^{n+1}}-1}-\frac{q}{q^2-1}\right)\right) \left(q^{2^{2 n+1}}+1\right)^{4^{-n}
   z}
\end{multline}
\end{example}
\begin{example}
\begin{multline}
\sum _{q=2}^n \left(\prod _{p=0}^{\infty }
   \left(\frac{2}{q^{2^p}-1}+1\right)^{\left(-\frac{1}{2}\right)^p} \prod _{p=0}^{\infty }
   \left(\frac{2}{q^{4^p}-1}+1\right)^{-2^{1-2 p}}\right)=-2 H_n+H_n^{(2)}+n
\end{multline}
\end{example}
\begin{example}
The exponential of the Polylogarithm function.
\begin{multline}
\prod_{p=0}^{2n}e^{\left(-\frac{1}{4}\right)^p
   \left(\text{Li}_2\left(-q^{2^p}\right)-\text{Li}_2\left(q^{2^p}\right)\right)}\prod_{p=0}^{n}e^{2^{1-4 p}
   \left(\text{Li}_2\left(q^{4^p}\right)-\text{Li}_2\left(-q^{4^p}\right)\right)}\\
   =e^{2 \text{Li}_2(q)-2^{-4 n-1}
   \text{Li}_2\left(q^{2^{2 n+1}}\right)}
\end{multline}
\end{example}
\begin{example}
\begin{equation}
\prod _{p=0}^{\infty } e^{\left(-\frac{1}{4}\right)^p
   \left(\text{Li}_2\left(-q^{2^p}\right)-\text{Li}_2\left(q^{2^p}\right)\right)} \prod _{p=0}^{\infty }
   e^{2^{1-4 p} \left(\text{Li}_2\left(q^{4^p}\right)-\text{Li}_2\left(-q^{4^p}\right)\right)}=e^{2
   \text{Li}_2(q)}
\end{equation}
\end{example}
\begin{figure}[H]
\includegraphics[scale=0.5]{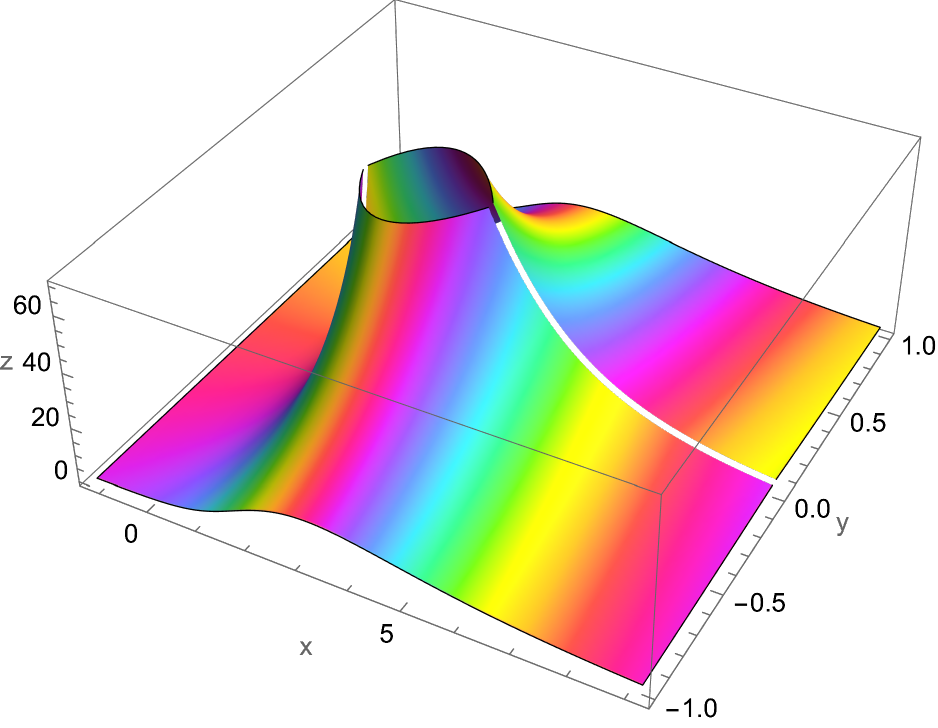}
\caption{Plot of  $f(q)=e^{2\text{Li}_2(q)}$, $q\in\mathbb{C}$.}
   \label{fig:fig2}
\end{figure}
\vspace{-6pt}
\begin{example}
\begin{multline}
\prod _{p=0}^{2 n} \exp \left(\left(-\frac{1}{2}\right)^{p+1} z^{2^{p+1}-1} \Phi
   \left(z^{2^{p+2}},1,\frac{1}{4} \left(2-2^{-p}\right)\right)\right)\\
 \prod _{p=0}^n \exp \left(4^{-p}
   z^{2^{2 p+1}-1} \Phi \left(z^{4^{p+1}},1,\frac{1}{4} \left(2-4^{-p}\right)\right)\right)\\
=\left(\frac{1+z }{1-z}\right)\exp
   \left(-2^{-2 n-1} z^{4^{n+1}-1} \Phi \left(z^{4^{n+1}},1,1-2^{-2 (n+1)}\right)\right)
\end{multline}
\end{example}
\begin{example}
\begin{equation}
\prod _{p=0}^n \left(\frac{\tan \left(2^p r\right)}{\tan \left(2^p
   m\right)}\right)^{2^{-p}}=\left(\frac{\sin (r)}{\sin (m)}\right)^2 \left(\frac{\sin \left(2^{1+n}
   m\right)}{\sin \left(2^{1+n} r\right)}\right)^{2^{-n}}
\end{equation}
\end{example}
\begin{example}
\begin{equation}
\prod _{p=0}^{\infty } \left(\frac{\tan \left(2^p r\right)}{\tan \left(2^p
   m\right)}\right)^{2^{-p-1}}=\frac{\sin (r)}{\sin (m)}
\end{equation}
\end{example}
\begin{example}
\begin{multline}
\prod _{p=0}^n \left(\frac{\left(m^{2^p}+1\right) \left(r^{2^p}-1\right)}{\left(m^{2^p}-1\right)
   \left(r^{2^p}+1\right)}\right)^{2^{-p}}=\frac{m }{r}\left(\frac{r-1}{m-1}\right)^2 \left(\frac{r^{2^n}
   \left(m^{2^{1+n}}-1\right)}{m^{2^n} \left(r^{2^{1+n}}-1\right)}\right)^{2^{-n}}
\end{multline}
\end{example}
\begin{example}
\begin{equation}
\prod _{p=0}^{\infty } \left(\frac{q^{2^p}+1}{q^{2^p}-1}\right)^{2^{-p}}=\frac{q^2}{(q-1)^2}
\end{equation}
\end{example}
\begin{figure}[H]
\includegraphics[scale=0.5]{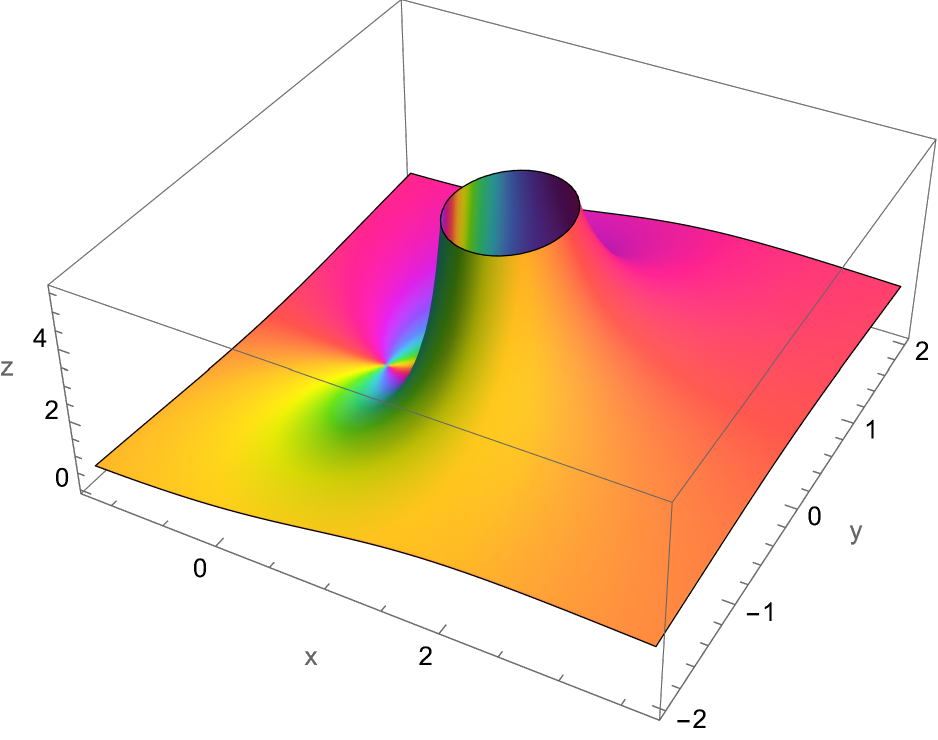}
\caption{Plot of  $f(q)=\frac{q^2}{(q-1)^2}$, $q\in\mathbb{C}$.}
   \label{fig:fig2}
\end{figure}
\vspace{-6pt}
\begin{example}
\begin{multline}
\prod _{p=0}^n e^{\csc \left(2^p z\right)-\csc \left(2^{1+p} z\right)} \left(\frac{\tan \left(2^{-1+p}
   z\right)}{\tan \left(2^p z\right)}\right)^{2^{-p} a}\\
=\exp \left(\int_{\frac{z}{2}}^z \frac{2^{2+n}+2 a \cot
   (x)-2 \csc (x) \left(a \cos \left(x-2^{2+n} x\right)+\csc (x) \sin ^2\left(2^{1+n} x\right)\right)}{-1+\cos
   \left(2^{2+n} x\right)} \, dx\right)
\end{multline}
\end{example}
\begin{example}
\begin{multline}
\prod _{p=0}^n \exp \left(2^p \left(\left(-1+4 \cos \left(2^p z\right)\right) \csc ^2\left(2^p
   z\right)+\sec ^2\left(2^p z\right)\right)\right) \left(\frac{\tan \left(2^p z\right)}{\tan \left(2^{-1+p}
   z\right)}\right)^{2^{-p} a^2}\\
=2^{2^{-n} \left(-1+2^{1+n}\right) a^2} \exp \left(2 \csc (z) (2 \cot (z)+\csc
   (z))+2^n \left(-3 \csc ^2\left(2^n z\right)+\sec ^2\left(2^n z\right)\right)\right)\\
 \cos ^{2
   a^2}\left(\frac{z}{2}\right) \sec ^{2^{-n} a^2}\left(2^n z\right)
\end{multline}
\end{example}
\begin{example}
\begin{multline}
\prod _{p=0}^n \left(\frac{2 z^{2^{p-1}}}{z^{2^p}+1}+1\right)^{a^2 2^{-p}} \exp \left(\frac{2^{p+3}
   z^{2^{p-1}} \left(z^{5\ 2^{p-1}}-z^{3\ 2^p}+z^{2^{p-1}}-3 z^{2^p}-3
   z^{2^{p+1}}-1\right)}{\left(z^{2^{p+1}}-1\right)^2}\right)\\
=\left(\sqrt[4]{z}+\frac{1}{\sqrt[4]{z}}\right)^{2
   a^2} \left(\frac{z^{2^{n-1}}}{z^{2^n}+1}\right)^{a^2 2^{-n}}\\
 \exp \left(-8 \left(-\frac{2^{n+1}
   \left(z^{2^n}+z^{2^{n+1}}+1\right)
   z^{2^n}}{\left(z^{2^{n+1}}-1\right)^2}+\frac{z^{3/2}}{(z-1)^2}+\frac{z}{(z-1)^2}+\frac{\sqrt{z}}{(z-1)^2}\right)
   \right)
\end{multline}
\end{example}
\begin{example}
\begin{multline}
\prod _{p=0}^{\infty } \left(\frac{2 z^{2^{p-1}}}{z^{2^p}+1}+1\right)^{a^2 2^{-p}} \exp \left(\frac{2^{p+3}
   z^{2^{p-1}} \left(z^{5\ 2^{p-1}}-z^{3\ 2^p}+z^{2^{p-1}}-3 z^{2^p}-3
   z^{2^{p+1}}-1\right)}{\left(z^{2^{p+1}}-1\right)^2}\right)\\
=e^{-\frac{8
   \left(z^{3/2}+z+\sqrt{z}\right)}{(z-1)^2}} \left(\sqrt[4]{z}+\frac{1}{\sqrt[4]{z}}\right)^{2 a^2}
   z^{-\frac{a^2}{2}}
\end{multline}
\end{example}
\begin{example}
\begin{multline}
\frac{2^{\frac{n}{2}+1} e^{-\frac{1}{4} \pi  \cot \left(\pi  2^{-n-2}\right)} \Gamma \left(2^{-n-3}
   a+1\right) }{\Gamma \left(\frac{1}{8} \left(2^{-n} a+4\right)\right)}\prod _{p=0}^n \left(i 2^p\right)^{\frac{e^{i \pi  2^{-n+p-1}}}{-1+e^{i \pi  2^{p-n}}}}\\
 \exp
   \left(\frac{1}{2} e^{i \pi  2^{-n+p-1}} \left(\Phi'\left(-e^{i \pi  2^{-n+p-1}},0,a
   2^{-p-1}+1\right)+\Phi'\left(e^{i \pi  2^{-n+p-1}},0,a
   2^{-p-1}+1\right)\right)\right)\\
=\exp \left(e^{i \pi 
   2^{-n-1}} \Phi'\left(e^{i \pi  2^{-n-1}},0,\frac{a}{2}+1\right)\right)
\end{multline}
\end{example}
\begin{example}
For all $j\in\mathbb Z_{\ne 0}$ then,
\begin{multline}
\prod _{p=0}^{\infty } \exp \left(\frac{2 i E_2\left(\frac{1}{2} i (-2 j+2 p+1) \pi \right)}{\pi  (2 j-2 p-1)}+\Gamma \left(-1,-\frac{1}{2} i (2 j+2 p+1) \pi \right)\right)=\frac{2}{e}
\end{multline}
\end{example}
\begin{example}
\begin{multline}
\prod _{p=0}^{n-1} \exp \left(\frac{i (-1)^p }{4 n e^{\frac{i m (1+2 p) \pi }{2 n \alpha }} }\left(-\Phi
   \left(e^{-\frac{2 i m \pi }{\alpha }},1,\frac{1+2 p}{4 n}\right)+e^{\frac{i m (1-2 n+2 p) \pi }{n \alpha }} \Phi
   \left(e^{-\frac{2 i m \pi }{\alpha }},1,1-\frac{1+2 p}{4 n}\right)\right.\right. \\ \left.\left.
+e^{\frac{i m (1+2 p) \pi }{n \alpha }} \Phi
   \left(e^{\frac{2 i m \pi }{\alpha }},1,\frac{1+2 p}{4 n}\right)-e^{\frac{2 i m \pi }{\alpha }} \Phi \left(e^{\frac{2
   i m \pi }{\alpha }},1,1-\frac{1+2 p}{4 n}\right)\right)\right)\\
=\cot \left(\frac{\pi  (m+n \alpha )}{4 n \alpha
   }\right)
\end{multline}
\end{example}
\begin{example}
\begin{multline}
\prod _{p=0}^{n-1} \exp \left(\frac{i \csc \left(\frac{m \pi
   }{n \alpha }\right)}{4 n e^{\frac{i \pi  (m+2 m p-n p \alpha )}{n \alpha }} } \left(\Phi \left(e^{-\frac{4 i m \pi }{\alpha }},1,\frac{1+p}{2 n}\right)-e^{\frac{2 i m (1-2 n+2
   p) \pi }{n \alpha }} \Phi \left(e^{-\frac{4 i m \pi }{\alpha }},1,1-\frac{p}{2 n}\right)\right.\right. \\ \left.\left.
-e^{\frac{2 i m (1+2 p) \pi
   }{n \alpha }} \Phi \left(e^{\frac{4 i m \pi }{\alpha }},1,\frac{1+p}{2 n}\right)+e^{\frac{4 i m \pi }{\alpha }} \Phi
   \left(e^{\frac{4 i m \pi }{\alpha }},1,1-\frac{p}{2 n}\right)\right)\right)\\
=\frac{1}{2} e^{\frac{m \pi  \cot
   \left(\frac{m \pi }{n \alpha }\right)}{n \alpha }} \sec \left(\frac{m \pi }{n \alpha }\right)
\end{multline}
\end{example}
\begin{example}
\begin{multline}
\prod _{p=0}^{n-1} \left(\frac{(1+2 p+n (4-2 a \alpha )) (1+2 p+2 n (-2+a \alpha ))}{(1+2 p-2 a n \alpha ) (1+2
   p+2 n (-4+a \alpha ))}\right)^{(-1)^p e^{-\frac{i p \pi }{2 n}} \left(-1+e^{\frac{i (1+2 p) \pi }{2 n}}\right)}\\
   \left(\frac{\Gamma \left(-\frac{1+4 n+2 p-2 a n \alpha }{8 n}\right)}{\Gamma \left(-\frac{1+2 p-2 a n \alpha }{8
   n}\right)}\right)^{(-1)^p e^{-\frac{i p \pi }{2 n}} \left(-1+e^{\frac{i (1+2 p) \pi }{2 n}}\right)}\\
   \left(\frac{\Gamma \left(\frac{1-4 n+2 p+2 a n \alpha }{8 n}\right)}{\Gamma \left(-1+\frac{1+2 p+2 a n \alpha }{8
   n}\right)}\right)^{(-1)^p e^{-\frac{i p \pi }{2 n}} \left(-1+e^{\frac{i (1+2 p) \pi }{2 n}}\right)}\\
=\exp
   \left(\Phi'\left(-e^{-\frac{i \pi }{2 n}},0,\frac{1}{2}+a n \alpha \right)-e^{\frac{i \pi
   }{2 n}} \Phi'\left(-e^{\frac{i \pi }{2 n}},0,\frac{1}{2}+a n \alpha \right)\right)
\end{multline}
\end{example}
\begin{example}
\begin{multline}
\prod _{p=0}^{n-1} \left(\frac{1}{\Gamma \left(\frac{1+2 p+2 a n \alpha }{4 n}\right) \Gamma \left(\frac{-1-2 p+2
   n (2+a \alpha )}{4 n}\right)}\right)^{(-1)^{\frac{(-1+n) p}{n}} \left(-1+(-1)^{\frac{1+2 p}{n}}\right)}\\
=\exp
   \left(\Phi'\left(-(-1)^{-1/n},0,\frac{1}{2}+a n \alpha \right)-(-1)^{1/n}
   \Phi'\left(-(-1)^{1/n},0,\frac{1}{2}+a n \alpha \right)\right)
\end{multline}
\end{example}
\begin{example}
\begin{equation}
\prod _{p=0}^{n-1} \frac{\Gamma \left(\frac{1}{2} \left(1+2^p x\right)\right)}{\Gamma \left(\frac{1}{4}
   \left(2+2^p x\right)\right)^2}=\frac{2^{-n+\frac{1}{2} \left(-1+2^n\right) x} \pi ^{-\frac{n}{2}} \Gamma
   \left(\frac{x}{4}\right) \Gamma \left(\frac{1}{4} \left(2+2^n x\right)\right)}{\Gamma \left(2^{-2+n}
   x\right) \Gamma \left(\frac{2+x}{4}\right)}
\end{equation}
\end{example}
\begin{example}
\begin{equation}
\prod _{p=0}^{n-1} \left(\frac{q^{2^{-1-p}}}{1+q^{2^{-p}}}\right)^{-2^{1+p}}
   \left(\frac{1}{q^{-2^{-p}}+q^{2^{-p}}}\right)^{2^p}=\frac{q
   \left(\frac{1}{q^{-2^{-n}}+q^{2^{-n}}}\right)^{-2^n}}{1+q^2}
\end{equation}
\end{example}
\begin{example}
\begin{multline}
\prod _{p=0}^{n-1} \frac{1}{8
   \pi ^3}\left(\frac{\pi  x}{n}-\cos ^{-1}\left(\cos \left(\frac{(1+2 p) \pi }{2
   n}\right)\right)\right) \left(\pi  \left(-2+\frac{x}{n}\right)+\cos ^{-1}\left(\cos \left(\frac{(1+2 p) \pi }{2
   n}\right)\right)\right)\\
 \Gamma \left(\frac{1}{4} \left(-2+\frac{2 x}{n}+\frac{i \left(\pi +2 \sin ^{-1}\left(\cos
   \left(\frac{(1+2 p) \pi }{2 n}\right)\right)\right) \sqrt{-1-\cos \left(\frac{(1+2 p) \pi }{2 n}\right)}}{\pi 
   \sqrt{1+\cos \left(\frac{(1+2 p) \pi }{2 n}\right)}}\right)\right) \\
\Gamma \left(\frac{1}{4} \left(-2+\frac{2
   x}{n}+\frac{i \left(\pi +2 \sin ^{-1}\left(\cos \left(\frac{(1+2 p) \pi }{2 n}\right)\right)\right) \sqrt{1+\cos
   \left(\frac{(1+2 p) \pi }{2 n}\right)}}{\pi  \sqrt{-1-\cos \left(\frac{(1+2 p) \pi }{2 n}\right)}}\right)\right)\\
=\frac{\Gamma \left(\frac{1}{2}+x\right)}{(2 n)^x \sqrt{2 \pi }}
\end{multline}
\end{example}
\begin{example}
\begin{multline}
\prod _{j=0}^{\infty } \prod _{n=0}^{\infty } \prod _{p=0}^{\infty } \exp
   \left(\frac{2^{1+2 p} (-t)^{j+2 p} t u^n \binom{-1+n+j (-1+\alpha )}{-j+n} E_{1+2
   p}(x) \Gamma (1+j+2 p)}{\Gamma (1+j) \Gamma (2+2 p)}\right)\\
=\frac{\Gamma
   \left(\frac{1}{4} \left(2+\frac{1}{t}+(1-u)^{-\alpha } u-2 x\right)\right) \Gamma
   \left(\frac{1}{4} \left(2+\frac{1}{t}+(1-u)^{-\alpha } u+2 x\right)\right)}{\Gamma
   \left(\frac{1}{4} \left(4+\frac{1}{t}+(1-u)^{-\alpha } u-2 x\right)\right) \Gamma
   \left(\frac{1}{4} \left(\frac{1}{t}+(1-u)^{-\alpha } u+2 x\right)\right)}
\end{multline}
\end{example}
\begin{example}
\begin{equation}
\prod _{p=0}^n \frac{1+\cos \left(2^{-p} (b+m)\right)}{1+\cos \left(2^{-p} (b-m)\right)}=\left(\frac{\sin (b+m)
   \sin \left(2^{-1-n} (b-m)\right)}{\sin \left(2^{-1-n} (b+m)\right) \sin (b-m)}\right)^2
\end{equation}
\end{example}
\begin{example}
\begin{equation}
\prod _{p=0}^n \frac{\cos \left(2^{-1-p} (b+m)\right)}{\cos \left(2^{-1-p} (b-m)\right)}=\frac{\sin
   \left(2^{-1-n} (b-m)\right) \sin (b+m)}{\sin (b-m) \sin \left(2^{-1-n} (b+m)\right)}
\end{equation}
\end{example}
\begin{example}
\begin{multline}
\prod _{p=0}^n \exp \left(-2^{-1-p} e^{-i 2^{1-p} m} \left(\Phi'\left(-e^{-i 2^{1-p}
   m},0,1+2^p a\right)\right.\right. \\ \left.\left.
   +e^{i 2^{2-p} m} \Phi'\left(-e^{i 2^{1-p} m},0,1+2^p a\right)\right)\right)\\
   \left(\frac{\Gamma \left(\frac{1}{2} \left(1+2^p a\right)\right)}{\Gamma \left(1+2^{-1+p}
   a\right)}\right)^{2^{-p}}\\
=\frac{\left(i^a 2^{2-2^{-n}-a (1+n)} \pi ^{1-2^{-1-n}} \Gamma \left(1+2^n a\right)^{2^{-n}}\right)
   }{\Gamma
   \left(\frac{2+a}{2}\right)^2}\\
\exp \left(-\frac{1}{2} i a \pi +e^{-4 i m} \Phi'\left(e^{-4 i m},0,1+\frac{a}{2}\right)+e^{4 i m} \Phi'\left(e^{4 i m},0,1+\frac{a}{2}\right)\right. \\ \left.
-2^{-1-n} e^{-i 2^{1-n} m} \left(\Phi'\left(e^{-i 2^{1-n} m},0,1+2^n a\right)+e^{i 2^{2-n} m}\Phi'\left(e^{i 2^{1-n} m},0,1+2^n a\right)\right)\right)
\end{multline}
\end{example}
\begin{example}
\begin{equation}
\prod _{p=0}^n \frac{\cos \left(2^{-p} b\right)+\cos \left(2^{-p} m\right)}{\cos \left(2^{-1-p} b\right)+\cos
   \left(2^{-1-p} m\right)}=\frac{\cos (b)+\cos (m)}{\cos \left(2^{-1-n} b\right)+\cos \left(2^{-1-n} m\right)}
\end{equation}
\end{example}
\begin{example}
See Eq. (1.432) in \cite{grad}. 
\begin{equation}
\prod _{p=0}^{\infty } \frac{\cos \left(2^{-p} b\right)+\cos \left(2^{-p} m\right)}{\cos \left(2^{-1-p} b\right)+\cos
   \left(2^{-1-p} m\right)}=\frac{1}{2} (\cos (b)+\cos (m))
\end{equation}
\end{example}
\begin{example}
\begin{equation}
\prod _{p=0}^{\infty } \frac{\cos \left(\frac{2^{-p}}{m}\right)+\cos
   \left(2^{-p} m\right)}{\cos \left(\frac{2^{-1-p}}{m}\right)+\cos
   \left(2^{-1-p} m\right)}=\frac{1}{2} \left(\cos \left(\frac{1}{m}\right)+\cos
   (m)\right)
\end{equation}
\end{example}
\begin{figure}[H]
\includegraphics[scale=0.5]{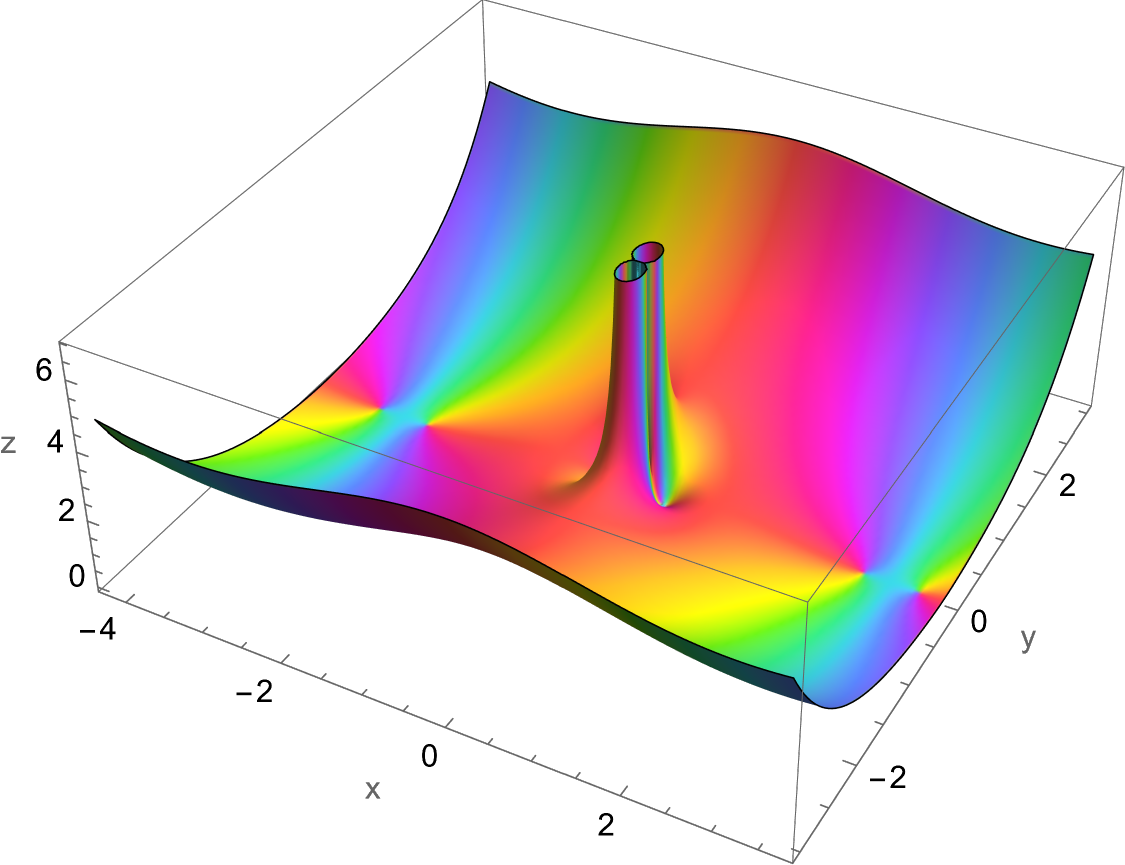}
\caption{Plot of  $f(m)=\frac{1}{2} \left(\cos \left(\frac{1}{m}\right)+\cos (m)\right)$, $m\in\mathbb{C}$.}
   \label{fig:fig2}
\end{figure}
\vspace{-6pt}
\begin{example}
\begin{equation}
\prod _{p=0}^{\infty } \frac{q^{-2^{-p}} \left(1+q^{2^{2-p}}+2 q^{2^{1-p}} \cosh \left(\frac{2^{-1-p}}{\log
   (q)}\right)\right)}{1+q^{2^{1-p}}+2 q^{2^{-p}} \cosh \left(\frac{2^{-2-p}}{\log (q)}\right)}=\frac{1+q^4+2 q^2 \cosh
   \left(\frac{1}{2 \log (q)}\right)}{4 q^2}
\end{equation}
\end{example}
\begin{example}
\begin{multline}
\prod _{j=1}^{\infty } \frac{\cosh \left(2^{-j} m\right) \cosh ^2\left(\frac{2^{-1-j} m}{b}\right) \cosh
   ^4\left(\frac{2^{-2-j} m}{b}\right)}{\cosh \left(\frac{2^{-j} m}{b}\right) \cosh ^2\left(2^{-1-j} m\right) \cosh
   ^4\left(2^{-2-j} m\right)}\\
   =\frac{b^5 \left(\sinh \left(\frac{m}{b}\right) \tanh ^2\left(\frac{m}{4 b}\right) \tanh
   ^2\left(\frac{m}{2 b}\right)\right)}{\tanh ^2\left(\frac{m}{4}\right) \tanh ^2\left(\frac{m}{2}\right) \sinh
   (m)}
\end{multline}
\end{example}
\begin{example}
\begin{equation}
\prod _{j=1}^{\infty } \frac{\cos ^4\left(2^{-4-j} \pi \right) \cos \left(2^{-1-j} \pi \right)}{\cos
   ^2\left(2^{-3-j} \pi \right) \cos ^3\left(2^{-2-j} \pi \right)}=16 \left(8+7 \sqrt{2}-4 \sqrt{10+7
   \sqrt{2}}\right)
\end{equation}
\end{example}
\begin{example}
\begin{multline}
\prod _{j=1}^{\infty } \frac{\cos ^4\left(\frac{1}{3} 2^{-3-j} \pi \right) \cos \left(\frac{2^{-j} \pi
   }{3}\right)}{\cos ^2\left(\frac{1}{3} 2^{-2-j} \pi \right) \cos ^3\left(\frac{1}{3} 2^{-1-j} \pi \right)}=32
   \left(24+15 \sqrt{3}-2 \sqrt{2} \left(9+5 \sqrt{3}\right)\right)
\end{multline}
\end{example}
\begin{example}
\begin{multline}
\prod _{j=1}^{\infty } \frac{\cos ^4\left(2^{-5-j} \pi \right) \cos \left(2^{-2-j} \pi \right)}{\cos
   ^2\left(2^{-4-j} \pi \right) \cos ^3\left(2^{-3-j} \pi \right)}=16 \left(2+\sqrt{2}\right)^{3/2}
   \left(-1+\frac{4}{2+\sqrt{2+\sqrt{2+\sqrt{2}}}}\right)
\end{multline}
\end{example}
\begin{example}
\begin{multline}
\prod _{j=1}^{\infty } \frac{\cos ^4\left(\frac{1}{5} 2^{-3-j} \pi \right) \cos \left(\frac{2^{-j} \pi
   }{5}\right)}{\cos ^2\left(\frac{1}{5} 2^{-2-j} \pi \right) \cos ^3\left(\frac{1}{5} 2^{-1-j} \pi
   \right)}\\
=\frac{-16 \sqrt{2} \left(5+3 \sqrt{5}\right)+32 \sqrt{5+\sqrt{5}} \left(-2+3
   \sqrt{2}-\sqrt{5}+\sqrt{10}\right)}{4+\sqrt{2 \left(4+\sqrt{2 \left(5+\sqrt{5}\right)}\right)}}
\end{multline}
\end{example}
\begin{example}
\begin{equation}
\prod _{j=1}^{\infty } \frac{\cos ^4\left(\frac{1}{3} 2^{-4-j} \pi \right) \cos \left(\frac{1}{3} 2^{-1-j}
   \pi \right)}{\cos ^2\left(\frac{1}{3} 2^{-3-j} \pi \right) \cos ^3\left(\frac{1}{3} 2^{-2-j} \pi \right)}=16
   \sqrt{2} \left(5+3 \sqrt{3}\right) \tan ^2\left(\frac{\pi }{48}\right)
\end{equation}
\end{example}
\begin{example}
\begin{multline}
\prod _{j=1}^{\infty } \frac{\left(1-2 \cos \left(2^{-2-j} \pi \right)\right)^4 \left(1-2 \cos
   \left(2^{-1-j} \pi \right)\right)^2}{-1+2 \cos \left(2^{-j} \pi \right)}=\frac{1}{243} \left(-17-12
   \sqrt{2}\right)
\end{multline}
\end{example}
\begin{example}
\begin{multline}
\prod _{j=1}^{\infty } \frac{\cos ^2\left(\frac{1}{5} 2^{1-j} \pi \right) \cos ^5\left(\frac{2^{-j} \pi }{5}\right)}{\cos
   ^4\left(\frac{1}{5} 2^{-2-j} \pi \right) \cos ^2\left(\frac{1}{5} 2^{-1-j} \pi \right) \cos \left(\frac{1}{5}
   2^{2-j} \pi \right)}=\frac{5 \left(95+42 \sqrt{5}+2 \sqrt{4450+1990 \sqrt{5}}\right)}{1024}
\end{multline}
\end{example}
\begin{example}
\begin{equation}
\prod _{j=1}^{\infty } \frac{\cos ^2\left(\frac{1}{5} 2^{1-j} \pi \right) \cos ^3\left(\frac{1}{5} 2^{2-j} \pi \right)}{\cos
   \left(\frac{1}{5} 2^{3-j} \pi \right) \cos ^4\left(\frac{2^{-j} \pi }{5}\right)}=\frac{1}{64}
   \left(-1-\sqrt{5}\right)
\end{equation}
\end{example}
\begin{example}
\begin{multline}
\prod _{j=1}^{\infty }\frac{\cos ^2\left(\frac{5}{7} 2^{-1-j} \pi \right) \cos ^3\left(\frac{5\ 2^{-j} \pi }{7}\right)}{\cos
   ^4\left(\frac{5}{7} 2^{-2-j} \pi \right) \cos \left(\frac{5}{7} 2^{1-j} \pi \right)}=\frac{\left(\cos
   \left(\frac{\pi }{28}\right)+\sin \left(\frac{3 \pi }{28}\right)\right)^2}{16 \sqrt{2} \left(\cos \left(\frac{\pi
   }{28}\right)-\sin \left(\frac{\pi }{28}\right)\right) \left(-1+\sin \left(\frac{\pi }{7}\right)\right)}
\end{multline}
\end{example}
\begin{example}
\begin{multline}
\prod _{j=1}^{\infty }\frac{\cos ^4\left(\frac{3}{7} 2^{-2-j} \pi \right) \cos \left(\frac{3}{7} 2^{1-j} \pi \right)}{\cos
   ^2\left(\frac{3}{7} 2^{-1-j} \pi \right) \cos ^3\left(\frac{3\ 2^{-j} \pi }{7}\right)}\\
   =\frac{128 \sqrt{2}
   \left(1+2 \cos \left(\frac{\pi }{14}\right)\right)^2 \sin ^3\left(\frac{\pi }{28}\right)}{-2+2 \cos
   \left(\frac{\pi }{14}\right)-3 \cos \left(\frac{\pi }{7}\right)+\cot \left(\frac{\pi }{28}\right)-3 \sin
   \left(\frac{\pi }{7}\right)}
\end{multline}
\end{example}
\begin{example}
\begin{multline}
\prod _{j=1}^{\infty }\frac{\cos ^2\left(\frac{3}{5} 2^{1-j} \pi \right) \cos ^5\left(\frac{3\ 2^{-j} \pi }{5}\right)}{\cos
   ^4\left(\frac{3}{5} 2^{-2-j} \pi \right) \cos ^2\left(\frac{3}{5} 2^{-1-j} \pi \right) \cos \left(\frac{3}{5}
   2^{2-j} \pi \right)}\\
   =\frac{5 \left(95-42 \sqrt{5}+2 \sqrt{4450-1990 \sqrt{5}}\right)}{1024}
\end{multline}
\end{example}
\begin{example}
\begin{equation}
\prod _{j=1}^{\infty }\frac{\cos ^4\left(\frac{1}{7} 2^{-2-j} \pi \right) \cos \left(\frac{1}{7} 2^{1-j} \pi \right)}{\cos
   ^2\left(\frac{1}{7} 2^{-1-j} \pi \right) \cos ^3\left(\frac{2^{-j} \pi }{7}\right)}=\frac{32 \left(1+\sin
   \left(\frac{3 \pi }{14}\right)\right) \tan ^2\left(\frac{\pi }{28}\right)}{\cos \left(\frac{\pi }{7}\right)-\sin
   \left(\frac{\pi }{14}\right)}
\end{equation}
\end{example}
\begin{example}
\begin{equation}
\prod _{j=1}^{\infty }\frac{\cos ^3\left(\frac{3}{7} 2^{1-j} \pi \right) \cos ^2\left(\frac{3\ 2^{-j} \pi }{7}\right)}{\cos
   ^4\left(\frac{3}{7} 2^{-1-j} \pi \right) \cos \left(\frac{3}{7} 2^{2-j} \pi \right)}=-\frac{\cos
   ^4\left(\frac{\pi }{14}\right)}{2 \left(1+2 \cos \left(\frac{\pi }{7}\right)\right)^3}
\end{equation}
\end{example}
\begin{example}
\begin{equation}
-i\prod _{j=1}^{\infty } \sqrt{\frac{\cos ^3\left(\frac{1}{3} 2^{1-j} \pi \right) \cos
   ^2\left(\frac{2^{-j} \pi }{3}\right)}{\cos ^4\left(\frac{1}{3} 2^{-1-j}
   \pi \right) \cos \left(\frac{1}{3} 2^{2-j} \pi \right)}}=\frac{3}{4
   \sqrt{2}}
\end{equation}
\end{example}
\begin{example}
\begin{multline}
-i\prod _{j=1}^{\infty } \frac{\cos \left(2^{-1-j} \pi \right) \cos ^4\left(2^{-4-j} \pi \right)}{\cos ^2\left(2^{-3-j} \pi
   \right) \cos ^3\left(2^{-2-j} \pi \right)}=16 \left(8+7 \sqrt{2}-4 \sqrt{10+7 \sqrt{2}}\right)
\end{multline}
\end{example}
\begin{example}
\begin{multline}
\prod _{j=1}^{\infty } \frac{e^{2^{-j} (-m+r)} \left(\cosh \left(2^{-j} m\right) \cosh
   ^2\left(2^{-1-j} r\right) \cosh ^4\left(2^{-2-j} r\right)\right)}{\cosh
   \left(2^{-j} r\right) \cosh ^2\left(2^{-1-j} m\right) \cosh
   ^4\left(2^{-2-j} m\right)}\\
=\frac{e^{-m+r} m^5 \left(\sinh (m) \sinh
   ^2\left(\frac{r}{2}\right) \sinh ^4\left(\frac{r}{4}\right)\right)}{r^5
   \left(\sinh (r) \sinh ^2\left(\frac{m}{2}\right) \sinh
   ^4\left(\frac{m}{4}\right)\right)}
\end{multline}
\end{example}
\begin{example}
\begin{multline}
\prod _{j=1}^{\infty } \cos ^4\left(2^{-4-j} \pi \right) \cos \left(2^{-1-j} \pi \right) \sec
   ^2\left(2^{-3-j} \pi \right) \sec ^3\left(2^{-2-j} \pi \right)\\
   =64 \sqrt{2} \csc ^2\left(\frac{\pi }{8}\right)
   \sin ^4\left(\frac{\pi }{16}\right)
\end{multline}
\end{example}
\begin{example}
\begin{multline}
\prod _{j=1}^{\infty } \cos ^4\left(2^{-4-j} \pi \right) \cos ^2\left(2^{-3-j} \pi \right) \cos
   \left(\frac{2^{-j} \pi }{3}\right)\\
    \sec ^4\left(\frac{1}{3} 2^{-2-j} \pi \right) \sec \left(2^{-2-j} \pi \right)
   \sec ^2\left(\frac{1}{3} 2^{-1-j} \pi \right)\\
=\frac{8192}{81} \sqrt{\frac{2}{3}} \left(1+\sqrt{3}\right)^4 \sin
   ^4\left(\frac{\pi }{16}\right) \sin ^2\left(\frac{\pi }{8}\right)
\end{multline}
\end{example}
\begin{example}
\begin{multline}
\prod _{j=1}^{\infty } \cos ^4\left(\frac{1}{3} 2^{-3-j} \pi \right) \cos \left(\frac{2^{-j} \pi }{3}\right)
   \sec ^2\left(\frac{1}{3} 2^{-2-j} \pi \right) \sec ^3\left(\frac{1}{3} 2^{-1-j} \pi \right)\\
=512 \left(3+2
   \sqrt{3}\right) \sin ^4\left(\frac{\pi }{24}\right)
\end{multline}
\end{example}
\begin{example}
\begin{multline}
\prod _{j=1}^{\infty } \cos ^4\left(\frac{1}{3} 2^{-3-j} \pi \right) \cos ^2\left(\frac{1}{3} 2^{-2-j} \pi
   \right) \cos \left(2^{-2-j} \pi \right)\\
    \sec ^4\left(2^{-4-j} \pi \right) \sec ^2\left(2^{-3-j} \pi \right) \sec
   \left(\frac{1}{3} 2^{-1-j} \pi \right)\\
=\frac{243}{64} \sqrt{\frac{7}{2}-2 \sqrt{3}} \csc ^4\left(\frac{\pi
   }{16}\right) \csc ^2\left(\frac{\pi }{8}\right) \sin ^4\left(\frac{\pi }{24}\right)
\end{multline}
\end{example}
\begin{example}
\begin{multline}
\prod _{j=1}^{\infty } \sqrt{\frac{\left(1-2 \cos \left(\frac{2^{-j} \pi }{3}\right)\right)^2}{\left(1-2
   \cos \left(\frac{1}{3} 2^{-2-j} \pi \right)\right)^8 \left(1-2 \cos \left(\frac{1}{3} 2^{-1-j} \pi
   \right)\right)^4}}\\
=\sqrt{59049 \left(7-4 \sqrt{3}\right) \csc ^8\left(\frac{\pi }{8}\right) \sin
   ^8\left(\frac{\pi }{24}\right)}
\end{multline}
\end{example}
\begin{example}
\begin{multline}
\prod _{j=1}^n \frac{e^{2^{2-j} (-m+r)} \left(\cosh \left(2^{2-j} m\right) \cosh ^2\left(2^{1-j} r\right)
   \cosh ^4\left(2^{-j} r\right)\right)}{\cosh \left(2^{2-j} r\right) \cosh ^2\left(2^{1-j} m\right) \cosh
   ^4\left(2^{-j} m\right)}\\
\frac{e^{4 \left(-1+2^{-n}\right) (m-r)} \left(\cosh (2 m) \cosh \left(2^{-n} m\right)
   \cosh (r) \cosh \left(2^{1-n} r\right) \sinh ^5\left(2^{-n} m\right) \sinh ^5(r)\right)}{\cosh (m) \cosh
   \left(2^{1-n} m\right) \cosh (2 r) \cosh \left(2^{-n} r\right) \sinh ^5(m) \sinh ^5\left(2^{-n}
   r\right)}
\end{multline}
\end{example}
\begin{example}
\begin{equation}
\prod _{p=0}^{n-1} \left(\frac{1-\cos \left(2^{-1-p} m\right)}{1+\cos \left(2^{-1-p}
   m\right)}\right)^{2^p}=\frac{2^{-1+2^n} \left(1-\cos \left(2^{-n} m\right)\right)^{2^n}}{1-\cos (m)}
\end{equation}
\end{example}
\begin{example}
\begin{multline}
\prod _{p=0}^{n-1} \exp \left(i 4^p e^{-i 2^{-1-p} m} \left(\Phi \left(e^{-i 2^{-p}
   m},2,\frac{1}{2}\right)-e^{i 2^{-p} m} \Phi \left(e^{i 2^{-p} m},2,\frac{1}{2}\right)\right)\right)\\
    \tan ^{2^{1+p}
   m}\left(2^{-2-p} m\right)\\
=\exp \left(-i \left(\text{Li}_2\left(e^{-i m}\right)-\text{Li}_2\left(e^{i m}\right)+4^n
   \left(-\text{Li}_2\left(e^{-i 2^{-n} m}\right)+\text{Li}_2\left(e^{i 2^{-n} m}\right)\right)\right)\right)\\
  \left(\frac{\left(2-2 \cos \left(2^{-n} m\right)\right)^{2^n}}{2-2 \cos (m)}\right)^m
\end{multline}
\end{example}
\begin{example}
\begin{equation}
\prod _{p=0}^{n-1} \Gamma \left(\frac{1}{6}+2^p x\right) \Gamma \left(\frac{5}{6}+2^p
   x\right)=\frac{2^{n-4 \left(-1+2^n\right) x} 27^{x-2^n x} \pi ^n \Gamma (x) \Gamma \left(3\times 2^n
   x\right)}{\Gamma (3 x) \Gamma \left(2^n x\right)}
\end{equation}
\end{example}
\begin{example}
\begin{multline}
\prod _{p=1}^n \frac{2 (-1)^{-2^{1-p} x} \left(2^{1+p}\right)^{-2^{1-p} x} \pi }{\Gamma
   \left(\frac{1}{6}+2^{-p} x\right) \Gamma \left(\frac{5}{6}+2^{-p} x\right)}\\
=\frac{(-1)^{2^{1-n} x} 2^{-2
   \left(1-2^{-n} (1+n)\right) x} 3^{\left(3-3\ 2^{-n}\right) x} e^{-2 i \pi  x} \Gamma (x) \Gamma \left(3\times2^{-n}
   x\right)}{\Gamma (3 x) \Gamma \left(2^{-n} x\right)}
\end{multline}
\end{example}
\begin{example}
\begin{equation}
\prod _{p=1}^{\infty } \frac{(-1)^{2^{1-p} x} \left(2^{1+p}\right)^{2^{1-p} x} \Gamma
   \left(\frac{1}{6}+2^{-p} x\right) \Gamma \left(\frac{5}{6}+2^{-p} x\right)}{2 \pi }=\frac{3^{1-3 x} 4^x e^{2 i
   \pi  x} \Gamma (3 x)}{\Gamma (x)}
\end{equation}
\end{example}
\begin{figure}[H]
\includegraphics[scale=0.5]{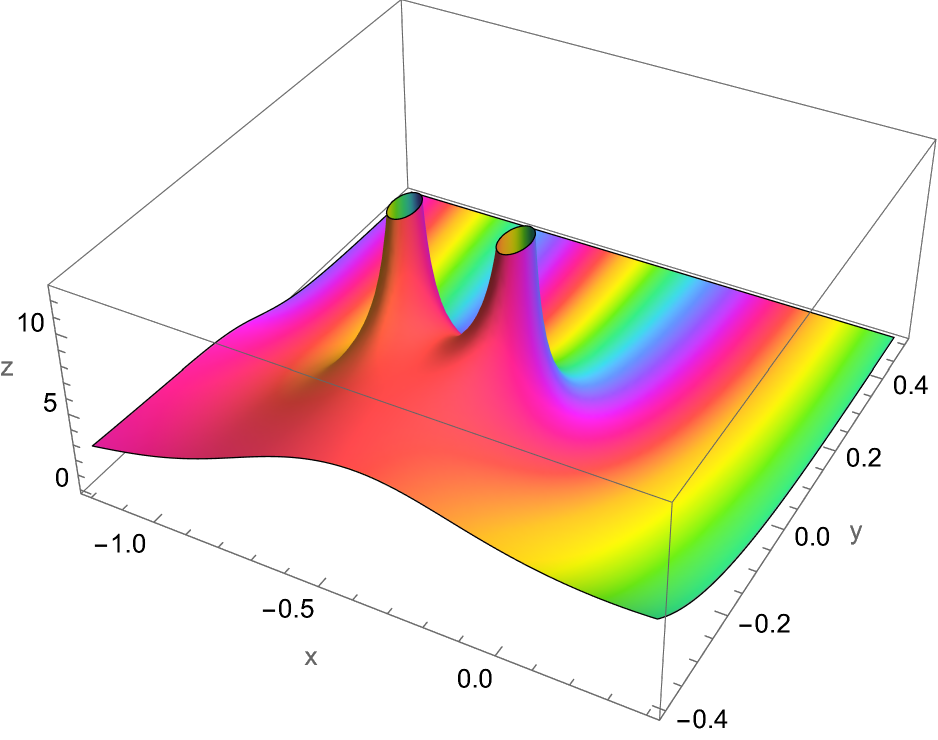}
\caption{Plot of  $f(x)=\frac{3^{1-3 x} 4^x e^{2 i \pi  x} \Gamma (3 x)}{\Gamma (x)}$, $x\in\mathbb{C}$.}
   \label{fig:fig2}
\end{figure}
\vspace{-6pt}
\begin{example}
Gauss's Multiplication Formula over even numbers.
\begin{equation}
\prod _{r=1}^{\frac{n-2}{2}} \Gamma \left(-\frac{r}{n}+z+1\right) \Gamma \left(\frac{r}{n}+z\right)=\frac{\pi
   ^{\frac{n}{2}-1} 2^{\frac{1}{2} (n+4 z-3)} n^{\frac{1}{2}-n z} \Gamma (n z)}{\Gamma (2 z)}
\end{equation}
\end{example}
\begin{example}
\begin{multline}
\prod _{r=1}^{\frac{n-2}{2}} \left(\frac{\Gamma \left(\frac{z n+n-r}{2 n}\right)}{\Gamma \left(\frac{1}{2}
   \left(-\frac{r}{n}+z+2\right)\right)}\right)^{e^{-\frac{i \pi  r}{n}}} \left(\frac{\Gamma \left(\frac{z n+n+r}{2
   n}\right)}{\Gamma \left(\frac{r+n z}{2 n}\right)}\right)^{e^{\frac{i \pi  r}{n}}}\\
=\frac{2^{\frac{1}{-1+e^{\frac{i \pi }{n}}}}
   n^{\frac{1}{-1+e^{\frac{i \pi }{n}}}} \left(\frac{\Gamma \left(\frac{z}{2}+\frac{1}{4}\right)}{\Gamma
   \left(\frac{z}{2}+\frac{3}{4}\right)}\right)^i \Gamma \left(\frac{z}{2}\right)
   }{\Gamma \left(\frac{z+1}{2}\right)}e^{-\Phi'\left(e^{\frac{i \pi }{n}},0,n z\right)}
\end{multline}
\end{example}
\begin{example}
Gauss's Multiplication Formula over odd numbers. 
\begin{multline}
\prod _{r=1}^{\frac{n-1}{2}} \Gamma \left(\frac{1}{2} \left(1-\frac{2 r}{n}+z\right)\right) \Gamma
   \left(\frac{1}{2} \left(1+\frac{2 r}{n}+z\right)\right)=\frac{n^{-\frac{1}{2} (n z)} (2 \pi )^{\frac{1}{2} (-1+n)}
   \Gamma \left(\frac{1}{2} (1+n z)\right)}{\Gamma \left(\frac{1+z}{2}\right)}
\end{multline}
\end{example}
\begin{example}
\begin{multline}
\prod _{r=1}^{\frac{n-1}{2}} \left(\frac{\Gamma \left(\frac{(3+a-b) n-2 r}{4 n}\right)}{\Gamma \left(\frac{n+a
   n-b n-2 r}{4 n}\right)}\right)^{e^{-\frac{i \pi  r}{n}}} \left(\frac{\Gamma \left(\frac{(3+a-b) n+2 r}{4
   n}\right)}{\Gamma \left(\frac{n+a n-b n+2 r}{4 n}\right)}\right)^{e^{\frac{i \pi  r}{n}}}\\
=\frac{1}{(2n)^{\frac{1}{2} \csc \left(\frac{\pi }{2 n}\right)} }\frac{\Gamma \left(\frac{1}{4} (1+a-b)\right) }{\Gamma \left(\frac{1}{4} (3+a-b)\right)}\exp \left(i
   e^{\frac{i \pi }{2 n}} \Phi'\left(e^{\frac{i \pi }{n}},0,\frac{1}{2} (1+a n-b
   n)\right)\right)
\end{multline}
\end{example}
\begin{example}
\begin{multline}
\prod _{r=1}^{\frac{n-1}{2}} \exp \left(e^{-\frac{i \pi  m r}{n}} \left(\Phi \left(e^{-i m \pi
   },1,\frac{r}{n}+\frac{1}{2}\right)-\Phi \left(e^{i m \pi },1,\frac{1}{2}-\frac{r}{n}\right)\right.\right. \\ \left.\left.
+e^{\frac{2 i \pi  m
   r}{n}} \left(\Phi \left(e^{-i m \pi },1,\frac{1}{2}-\frac{r}{n}\right)-\Phi \left(e^{i m \pi
   },1,\frac{r}{n}+\frac{1}{2}\right)\right)\right)\right)\\
=\left(-i \tan \left(\frac{\pi 
   m}{4}\right)\right)^{e^{\frac{i \pi  m}{2}}} \left(-i \cot \left(\frac{\pi  m}{4}\right)\right)^{e^{-\frac{1}{2} i
   \pi  m}} \left(i \tan \left(\frac{\pi  m}{4 n}\right)\right)^{e^{-\frac{1}{2} i \pi  m} n} \left(i \cot
   \left(\frac{\pi  m}{4 n}\right)\right)^{e^{\frac{i \pi  m}{2}} n}
\end{multline}
\end{example}
\begin{example}
Gauss's Multiplication Formula over integers.
\begin{multline}
\prod _{r=0}^{n-1} \Gamma \left(a-\frac{2 \pi  r+t}{2 n \pi }+1\right) \Gamma \left(a+\frac{2 \pi  r+t}{2 n \pi
   }\right)\\=(2 \pi )^{n-1} n^{-2 a n} \Gamma \left(a n-\frac{t}{2 \pi }+1\right) \Gamma \left(a n+\frac{t}{2 \pi
   }\right)
\end{multline}
\end{example}
\begin{example}
\begin{multline}
\prod _{r=0}^{n-1} \left(\frac{\Gamma \left(\frac{2 \pi  (n+a n+r)+t}{4 n \pi }\right)}{\Gamma \left(\frac{2
   \pi  (a n+r)+t}{4 n \pi }\right)}\right)^{e^{\frac{i (\pi  r+t)}{n}}} \left(\frac{\Gamma \left(-\frac{-2 (1+a) n
   \pi +2 \pi  r+t}{4 n \pi }\right)}{\Gamma \left(\frac{1}{4} \left(4+2 a-\frac{2 \pi  r+t}{n \pi
   }\right)\right)}\right)^{e^{-\frac{i \pi  r}{n}}}\\
=\exp \left(-e^{\frac{i \pi }{n}}
   \Phi'\left(e^{\frac{i \pi }{n}},0,1+a n-\frac{t}{2 \pi }\right)-e^{\frac{i t}{n}}
   \Phi'\left(e^{\frac{i \pi }{n}},0,a n+\frac{t}{2 \pi }\right)\right) (2 n)^{-i
   e^{\frac{i t}{2 n}} \cos \left(\frac{\pi -t}{2 n}\right) \csc \left(\frac{\pi }{2 n}\right)}
\end{multline}
\end{example}
\begin{example}
\begin{multline}
\prod _{r=0}^{n-1} \left(\frac{\Gamma \left(\frac{1}{4} \left(1+a-\frac{2 \pi  r+t}{2 n \pi
   }\right)\right)}{\Gamma \left(\frac{1}{4} \left(3+a-\frac{2 \pi  r+t}{2 n \pi }\right)\right)}\right)^{e^{\frac{\pi
    (n+r)}{2 n i}}} \left(\frac{\Gamma \left(\frac{1}{4} \left(2+a-\frac{2 \pi  r+t}{2 n \pi }\right)\right)}{\Gamma
   \left(\frac{1}{4} \left(4+a-\frac{2 \pi  r+t}{2 n \pi }\right)\right)}\right)^{e^{-\frac{i \pi  r}{2 n}}}\\
   \left(\frac{\Gamma \left(\frac{2 \pi  (a n+r)+t}{8 n \pi }\right)}{\Gamma \left(\frac{1}{4}
   \left(2+a+\frac{r+\frac{t}{2 \pi }}{n}\right)\right)}\right)^{-e^{\frac{i (\pi  r+t)}{2 n}}} \left(\frac{\Gamma
   \left(\frac{1}{4} \left(1+a+\frac{r+\frac{t}{2 \pi }}{n}\right)\right)}{\Gamma \left(\frac{1}{4}
   \left(3+a+\frac{r+\frac{t}{2 \pi }}{n}\right)\right)}\right)^{\exp \left(\frac{\pi -\frac{\pi  r+t}{n}}{2
   i}\right)}\\
=\exp \left(-e^{\frac{i \pi }{2 n}} \Phi'\left(e^{\frac{i \pi }{2 n}},0,1+a
   n-\frac{t}{2 \pi }\right)-e^{\frac{i t}{2 n}} \Phi'\left(e^{\frac{i \pi }{2 n}},0,a
   n+\frac{t}{2 \pi }\right)\right)\\
 (4 n)^{-i e^{\frac{i t}{4 n}} \cos \left(\frac{\pi -t}{4 n}\right) \csc
   \left(\frac{\pi }{4 n}\right)}
\end{multline}
\end{example}
\begin{example}
\begin{multline}
\prod _{p=1}^{n-1} \left(\frac{\Gamma \left(\frac{a n+p}{2 n}\right)}{\Gamma \left(\frac{1}{2}
   \left(2+a-\frac{p}{n}\right)\right)}\right)^{(-1)^p \sin \left(\frac{p \pi }{n}\right)}\\
=(2 n)^{-\frac{1}{2} \tan
   \left(\frac{\pi }{2 n}\right)} \exp \left(\frac{1}{2} i e^{-\frac{i \pi }{n}}
   \left(-\Phi'\left(-e^{-\frac{i \pi }{n}},0,1+a n\right)+e^{\frac{2 i \pi }{n}}
   \Phi'\left(-e^{\frac{i \pi }{n}},0,1+a n\right)\right)\right)
\end{multline}
\end{example}
\begin{example}
\begin{multline}
\prod _{p=1}^{n-1} \exp \left(i (-1)^p e^{-\frac{p \pi  (m+r)}{n}}
   \left(-e^{\frac{\pi  (2 m (-n+p)+p r)}{n}} \Phi \left(e^{2 m \pi
   },1,\frac{p}{2 n}\right)\right.\right. \\ \left.\left.
+e^{\frac{p \pi  r}{n}} \Phi \left(e^{2 m \pi
   },1,1-\frac{p}{2 n}\right)+e^{\frac{\pi  (m p-2 n r+2 p r)}{n}} \Phi
   \left(e^{2 \pi  r},1,\frac{p}{2 n}\right)\right.\right. \\ \left.\left.
-e^{\frac{m p \pi }{n}} \Phi
   \left(e^{2 \pi  r},1,1-\frac{p}{2 n}\right)\right) \sin \left(\frac{p \pi
   }{n}\right)\right)\\
=\left(\frac{e^{\frac{i \pi }{n}} \cosh \left(\frac{(i+m)
   \pi }{2 n}\right)}{\cosh \left(\frac{(-i+m) \pi }{2 n}\right)}\right)^{e^{-2 m
   \pi } n} \left(\frac{e^{-\frac{i \pi }{n}} \cosh \left(\frac{\pi  (-i+r)}{2
   n}\right)}{\cosh \left(\frac{\pi  (i+r)}{2 n}\right)}\right)^{e^{-2 \pi  r}
   n}
\end{multline}
\end{example}
\begin{example}
\begin{multline}
\prod _{p=1}^{n-1} (8 \pi )^{2 i e^{i p \pi } \sin ^2\left(\frac{p \pi
   }{n}\right)} \left(\left(a-\frac{p}{n}\right) \Gamma \left(\frac{1}{2}
   \left(a-\frac{p}{n}\right)\right)\right)^{2 (-1)^p e^{-\frac{i p \pi }{n}}
   \sin \left(\frac{p \pi }{n}\right)}\\
 \left(\left(-2+a+\frac{p}{n}\right) \Gamma
   \left(\frac{1}{2} \left(-2+a+\frac{p}{n}\right)\right)\right)^{i (-1)^p
   \left(-1+e^{\frac{2 i p \pi }{n}}\right)}\\
=2^{\frac{1}{2} \left(-i+\tan
   \left(\frac{\pi }{n}\right)\right)} \exp \left(-i e^{\frac{2 i \pi }{n}}
   \Phi'\left(-e^{\frac{2 i \pi }{n}},0,1+a
   n\right)\right) n^{\frac{1}{2} \tan \left(\frac{\pi }{n}\right)}
   \left(\frac{\Gamma \left(\frac{1}{2} (1+a n)\right)}{\Gamma \left(1+\frac{a
   n}{2}\right)}\right)^i
\end{multline}
\end{example}
\begin{example}
\begin{multline}
\prod _{p=1}^{n-1} \exp \left(2 (-1)^p e^{-\frac{i p \pi }{n}} \text{log$\Gamma $}\left(\frac{1}{2}
   \left(-2+a+\frac{p}{n}\right)\right) \sin \left(\frac{p \pi }{n}\right)\right) \left(-2+a+\frac{p}{n}\right)^{2
   (-1)^p e^{-\frac{i p \pi }{n}} \sin \left(\frac{p \pi }{n}\right)}\\
 \Gamma \left(\frac{1}{2}
   \left(2+a-\frac{p}{n}\right)\right)^{i (-1)^p \left(-1+e^{\frac{2 i p \pi }{n}}\right)}\\
=2^{\frac{1}{2}
   \left(-i+\left(-2+e^{i n \pi }\right) \tan \left(\frac{\pi }{n}\right)\right)} \exp \left(-i e^{-\frac{2 i \pi }{n}}
   \Phi'\left(-e^{-\frac{2 i \pi }{n}},0,1+a n\right)\right)\\
 n^{-\frac{1}{2} \tan
   \left(\frac{\pi }{n}\right)} \pi ^{\frac{1}{2} i \sin (n \pi ) \left(i+\tan \left(\frac{\pi }{n}\right)\right)}
   \left(\frac{\Gamma \left(\frac{1}{2} (1+a n)\right)}{\Gamma \left(1+\frac{a n}{2}\right)}\right)^i
\end{multline}
\end{example}
\begin{example}
\begin{multline}
\prod _{p=0}^{n-1} \left(\frac{\Gamma \left(\frac{-1+2 (2+a) n-2 p}{4 n}\right)}{\Gamma \left(\frac{1+2 a n+2
   p}{4 n}\right)}\right)^{(-1)^p \cos \left(\frac{j (1+2 p) \pi }{2 n}\right)}\\
=\exp \left(\frac{1}{2} e^{-\frac{3 i j
   \pi }{2 n}} \left(-e^{\frac{i j \pi }{n}} \Phi'\left(e^{-\frac{2 i j \pi
   }{n}},0,\frac{1}{4} (1+2 a n)\right)+\Phi'\left(e^{-\frac{2 i j \pi }{n}},0,\frac{1}{4}
   (3+2 a n)\right)\right.\right. \\ \left.\left.
-e^{\frac{2 i j \pi }{n}} \Phi'\left(e^{\frac{2 i j \pi
   }{n}},0,\frac{1}{4} (1+2 a n)\right)+e^{\frac{3 i j \pi }{n}} \Phi'\left(e^{\frac{2 i j
   \pi }{n}},0,\frac{1}{4} (3+2 a n)\right)\right)\right) n^{-\frac{1}{2} \sec \left(\frac{j \pi }{2
   n}\right)}
\end{multline}
\end{example}
\begin{example}
\begin{multline}
\prod _{p=0}^{n-1} 2^{-(-1)^p e^{\frac{i p \pi }{2 n}} \left(-1+e^{\frac{i (-1+2 n-2 p) \pi }{2 n}}\right) \cos
   \left(\frac{j (1+2 p) \pi }{2 n}\right)} \left(\frac{\Gamma \left(\frac{-1+2 (2+a) n-2 p}{8 n}\right)}{\Gamma
   \left(\frac{-1+2 (4+a) n-2 p}{8 n}\right)}\right)^{2 (-1)^p e^{\frac{i (-1+2 n-p) \pi }{2 n}} \cos \left(\frac{j (1+2
   p) \pi }{2 n}\right)}\\
 \left(\frac{\Gamma \left(\frac{1+2 (2+a) n+2 p}{8 n}\right)}{\Gamma \left(\frac{1+2 a n+2 p}{8
   n}\right)}\right)^{2 (-1)^p e^{\frac{i p \pi }{2 n}} \cos \left(\frac{j (1+2 p) \pi }{2 n}\right)}\\
=\exp
   \left(e^{\frac{i (1-3 j) \pi }{2 n}} \left(-e^{\frac{i (-1+2 j) \pi }{2 n}}
   \Phi'\left(e^{\frac{i (1-2 j) \pi }{n}},0,\frac{1}{4} (1+2 a
   n)\right)+\Phi'\left(e^{\frac{i (1-2 j) \pi }{n}},0,\frac{1}{4} (3+2 a
   n)\right)\right.\right. \\ \left.\left.
-e^{\frac{i (-1+4 j) \pi }{2 n}} \Phi'\left(e^{\frac{i (1+2 j) \pi
   }{n}},0,\frac{1}{4} (1+2 a n)\right)+e^{\frac{3 i j \pi }{n}} \Phi'\left(e^{\frac{i (1+2
   j) \pi }{n}},0,\frac{1}{4} (3+2 a n)\right)\right)\right)\\
 n^{-\frac{2 e^{-\frac{i \pi }{4 n}} \cos \left(\frac{\pi
   }{4 n}\right) \cos \left(\frac{j \pi }{2 n}\right)}{\cos \left(\frac{\pi }{2 n}\right)+\cos \left(\frac{j \pi
   }{n}\right)}}
\end{multline}
\end{example}
\begin{example}
\begin{equation}
\prod _{p=1}^n \frac{\Gamma \left(\frac{1}{2} \left(1+2^p x\right)\right)^2}{\Gamma \left(\frac{1}{2}+2^p x\right)}=\frac{2^{n-2 \left(-1+2^n\right) x} \pi ^{n/2} \Gamma \left(2^n
   x\right) \Gamma \left(\frac{1}{2}+x\right)}{\Gamma (x) \Gamma \left(\frac{1}{2}+2^n x\right)}
\end{equation}
\end{example}
\begin{example}
\begin{multline}
\prod _{p=1}^n \left(\frac{\cos \left(2^{-p} m\right) \left(1+\cos \left(2^{-p} r\right)\right) \sec \left(2^{-p} r\right)}{1+\cos \left(2^{-p} m\right)}\right)^{2^p}=\frac{(1+\cos
   (m)) \left(\frac{1+\cos \left(2^{-n} r\right)}{1+\cos \left(2^{-n} m\right)}\right)^{2^n}}{1+\cos (r)}
\end{multline}
\end{example}
\begin{example}
\begin{equation}
\prod _{p=1}^n \left(\frac{1+q^{2^{-p}}}{\sqrt{1+q^{2^{1-p}}}}\right)^{2^p}=\frac{\left(1+q^{2^{-n}}\right)^{2^n}}{1+q}
\end{equation}
\end{example}
\begin{example}
\begin{multline}
\prod _{p=1}^n \left(e^{i 2^{1-p} m \pi }\right)^{-2^{-1+p}} \left(32 \pi ^6 \left(1+i \tan \left(2^{-p} m \pi
   \right)\right)\right)^{2^p} \left(4 \pi ^3 \sec \left(2^{-1-p} m \pi \right)\right)^{-2^{p+1}}\\
=\frac{2^{2^{1+n}} \left(e^{i
   2^{-n} m \pi }\right)^{2^n} \sec ^{-2^{1+n}}\left(2^{-1-n} m \pi \right)}{\left(1+e^{i m \pi }\right)^2}
\end{multline}
\end{example}
\begin{example}
\begin{multline}
\prod _{p=1}^n \left(q^{2^{1-p}}\right)^{2^{-1+p}} \left(1+q^{2^{1-p}}\right)^{-2^p}
   \left(\frac{q^{2^{-1-p}}}{1+q^{2^{-p}}}\right)^{-2^{1+p}}=\frac{\left(q^{2^{-n}}\right)^{2^n}
   \left(\frac{q^{2^{-1-n}}}{1+q^{2^{-n}}}\right)^{-2^{1+n}}}{(1+q)^2}
\end{multline}
\end{example}
\begin{example}
\begin{multline}
\prod _{n=0}^{\infty } \prod _{p=0}^{\infty } \prod _{q=0}^{\infty } \exp \left(\frac{(-1)^{2 n+p+2 q} 2^p (a+n \alpha )^{-1-n-p-2 q} \beta  (n \alpha +\beta )^{-1+n} E_p(x) \Gamma
   (1+n+p+2 q)}{\Gamma (1+n) \Gamma (1+p) \Gamma (2+2 q)}\right)\\
=\frac{(-1+a+2 x-\beta ) \Gamma \left(\frac{1}{4} (-1+a+2 x-\beta )\right)^2}{4 \Gamma \left(\frac{1}{4} (1+a+2 x-\beta
   )\right)^2}
\end{multline}
\end{example}
\begin{example}
\begin{multline}
\prod_{p=0}^{\infty}\left(1+\frac{2 (m-r) (\alpha +2 p \alpha )}{(m+\alpha +2 p \alpha )
   (r-(1+2 p) \alpha )}\right)^{(-1)^p \cos ((1+2 p) x)}\\
=\exp \left(\frac{1}{2}
   i \left(e^{\frac{i m (\pi -2 x)}{2 \alpha }} \Phi \left(-e^{\frac{i m \pi
   }{\alpha }},1,\frac{1}{2}-\frac{x}{\pi }\right)+e^{\frac{i m (\pi +2 x)}{2
   \alpha }} \Phi \left(-e^{\frac{i m \pi }{\alpha }},1,\frac{1}{2}+\frac{x}{\pi
   }\right)\right.\right. \\ \left.\left.
-e^{\frac{i r (\pi -2 x)}{2 \alpha }} \Phi \left(-e^{\frac{i \pi 
   r}{\alpha }},1,\frac{1}{2}-\frac{x}{\pi }\right)-e^{\frac{i r (\pi +2 x)}{2
   \alpha }} \Phi \left(-e^{\frac{i \pi  r}{\alpha }},1,\frac{1}{2}+\frac{x}{\pi
   }\right)\right)\right)
\end{multline}
\end{example}
\begin{example}
\begin{multline}
\prod _{p=0}^{\infty } \exp \left(\frac{e^{i p \pi } \alpha  \left(2 m-i
   \left(m^2-(\alpha +2 p \alpha )^2\right) \log (a)\right) \sin (x+2 p
   x)}{(-m+\alpha +2 p \alpha )^2 (m+\alpha +2 p \alpha )^2}\right)\\
=\exp
   \left(\frac{\pi  \sec \left(\frac{m \pi }{2 \alpha }\right) \left(2 m x \cos
   \left(\frac{m x}{\alpha }\right)+\sin \left(\frac{m x}{\alpha }\right)
   \left(-2 \alpha +2 i m \alpha  \log (a)+m \pi  \tan \left(\frac{m \pi }{2
   \alpha }\right)\right)\right)}{8 m^2 \alpha }\right)
\end{multline}
\end{example}
\begin{example}
\begin{multline}
\prod _{p=0}^{\infty } \exp \left(\frac{2 (-1)^p e^{-i (1+2 p) \alpha }
   \cos ((1+2 p) x) \left(E_1(-i (1+2 p) \alpha )+e^{2 i (\alpha +2 p \alpha )}
   E_1(i (1+2 p) \alpha )\right)}{(1+2 p) \pi }\right)\\
=\frac{2 \pi  \Gamma
   \left(\frac{3 \pi -2 x+2 \alpha }{4 \pi }\right) \Gamma \left(\frac{3 \pi +2
   x+2 \alpha }{4 \pi }\right)}{\alpha  \Gamma \left(\frac{\pi -2 x+2 \alpha }{4
   \pi }\right) \Gamma \left(\frac{\pi +2 x+2 \alpha }{4 \pi }\right)}
\end{multline}
\end{example}
\begin{example}
\begin{multline}
\prod _{p=0}^{\infty } \left(\frac{m-(1+2 p) \alpha }{m+\alpha +2 p
   \alpha }\right)^{\frac{(-1)^p \cos (x+2 p x)}{(-m+\alpha +2 p \alpha
   )^2}}\\
=\exp \left(\frac{\pi }{8 \alpha ^2} \left((-2+2 \gamma -i \pi ) \sec \left(\frac{m
   \pi }{2 \alpha }\right) \left(-2 x \sin \left(\frac{m x}{\alpha }\right)+\pi 
   \cos \left(\frac{m x}{\alpha }\right) \tan \left(\frac{m \pi }{2 \alpha
   }\right)\right)\right.\right. \\ \left.
\left.
-4 i \pi  \left(e^{\frac{i m (\pi -2 x)}{2 \alpha }}
   \Phi'\left(-e^{\frac{i m \pi }{\alpha
   }},-1,\frac{1}{2}-\frac{x}{\pi }\right)+e^{\frac{i m (\pi +2 x)}{2 \alpha }}
   \Phi'\left(-e^{\frac{i m \pi }{\alpha
   }},-1,\frac{1}{2}+\frac{x}{\pi }\right)\right)\right)\right)\\
   \pi ^{\frac{\pi  \sec \left(\frac{m \pi }{2 \alpha }\right) \left(-2 x \sin
   \left(\frac{m x}{\alpha }\right)+\pi  \cos \left(\frac{m x}{\alpha }\right)
   \tan \left(\frac{m \pi }{2 \alpha }\right)\right)}{4 \alpha ^2}}
   \left(\frac{1}{\alpha }\right)^{\frac{\pi  \sec \left(\frac{m \pi }{2 \alpha
   }\right) \left(-2 x \sin \left(\frac{m x}{\alpha }\right)+\pi  \cos
   \left(\frac{m x}{\alpha }\right) \tan \left(\frac{m \pi }{2 \alpha
   }\right)\right)}{4 \alpha ^2}}
\end{multline}
\end{example}
\begin{example}
\begin{multline}
\prod _{p=0}^{\infty } \exp \left(\frac{4 i (-1)^{\frac{1}{4}+p}  e^{-\frac{1}{2} (i x)} \cos
   (x+2 p x) \left(\frac{a^{4 i m (1+2 p)} \Gamma (0,i m (1+4 p) \log (a))}{1+4 p}+\frac{\Gamma (0,-i m (3+4 p) \log
   (a))}{3+4 p}\right)}{a^{i m (3+4 p)}\pi }\right)\\
 (i m (1+4 p))^{-\frac{4 i (-1)^{\frac{1}{4}+p} e^{-\frac{1}{2} (i x)} \cos (x+2 p
   x)}{\pi +4 p \pi }} (-i m (3+4 p))^{-\frac{4 i (-1)^{\frac{1}{4}+p} e^{-\frac{1}{2} (i x)} \cos (x+2 p x)}{(3+4 p)
   \pi }}\\
 (i m (1+4 p) \log (a))^{\frac{4 i (-1)^{\frac{1}{4}+p} e^{-\frac{1}{2} (i x)} \cos (x+2 p x)}{\pi +4 p \pi
   }} (-i m (3+4 p) \log (a))^{\frac{4 i (-1)^{\frac{1}{4}+p} e^{-\frac{1}{2} (i x)} \cos (x+2 p x)}{(3+4 p) \pi
   }}\\
=\frac{ \left(\frac{\Gamma \left(\frac{\pi -2 x+4 m \log
   (a)}{8 \pi }\right)}{\Gamma \left(\frac{5 \pi -2 x+4 m \log (a)}{8 \pi }\right)}\right)^{e^{-i x}}
   \left(\frac{\Gamma \left(\frac{3 \pi -2 x+4 m \log (a)}{8 \pi }\right)}{\Gamma \left(\frac{7 \pi -2 x+4 m \log
   (a)}{8 \pi }\right)}\right)^{-i e^{-i x}} \Gamma \left(\frac{\pi +2 x+4 m \log (a)}{8 \pi }\right)
   \left(\frac{\Gamma \left(\frac{3 \pi +2 x+4 m \log (a)}{8 \pi }\right)}{\Gamma \left(\frac{7 \pi +2 x+4 m \log
   (a)}{8 \pi }\right)}\right)^{-i}}{\Gamma \left(\frac{5 \pi +2 x+4 m \log (a)}{8 \pi }\right)}\\
\left(\frac{1}{m}\right)^{\left(-\frac{1}{2}+\frac{i}{2}\right) \left(1+e^{-i x}\right)} (2 \pi
   )^{\left(-\frac{1}{2}+\frac{i}{2}\right) \left(1+e^{-i x}\right)}
\end{multline}
\end{example}
\begin{example}
\begin{equation}
\prod _{p=1}^n \frac{\Gamma \left(3^{-p} z\right)}{\Gamma \left(3^{1-p} z\right)}=\frac{\Gamma \left(3^{-n} z\right)}{\Gamma (z)}
\end{equation}
\end{example}
\begin{example}
\begin{equation}
\prod _{p=1}^n \left(1-\frac{3 m^{3^{-1+p}}}{1+m^{3^{-1+p}}+m^{2\times 3^{-1+p}}}\right)^{3^{-p}}=(1-m)
   \left(\frac{1}{1-m^{3^n}}\right)^{3^{-n}}
\end{equation}
\end{example}
\begin{example}
\begin{equation}
\prod _{p=1}^n \left(1+2 \cos \left(2\times 3^{-p} m\right)\right)=\frac{\sin (m)}{\sin \left(3^{-n}
   m\right)}
\end{equation}
\end{example}
\begin{example}
\begin{multline}
\prod _{p=1}^n \exp \left(3^{-p} e^{i 3^{-p} \pi } \left(\frac{\left(2+e^{-i 3^{-p} \pi }\right) \log \left(i
   3^{-p}\right)}{1+2 \cos \left(3^{-p} \pi \right)}-\sqrt[3]{-1} \Phi'\left(-\sqrt[3]{-1}
   e^{i 3^{-p} \pi },0,1+2\times 3^p z\right)\right.\right. \\ \left.\left.
+(-1)^{2/3} \Phi'\left((-1)^{2/3} e^{i 3^{-p} \pi
   },0,1+2\times 3^p z\right)\right)\right)\\
=\frac{(1+i) \left(i 3^{-n}\right)^{-3^{-n} \left(1+\frac{1}{-1+e^{i 3^{-n} \pi
   }}\right)} \exp \left(-3^{-n} e^{i 3^{-n} \pi } \Phi'\left(e^{i 3^{-n} \pi },0,1+2\times 3^n
   z\right)\right) z \Gamma (z)}{\Gamma \left(\frac{1}{2}+z\right)}
\end{multline}
\end{example}
\begin{example}
\begin{equation}
\prod _{p=1}^n \left(1+q^{3^{-p}}+q^{2\times 3^{-p}}\right)=\frac{q-1}{q^{3^{-n}}-1}
\end{equation}
\end{example}
\begin{example}
\begin{equation}
\prod _{p=1}^n \frac{\Gamma \left(3^p x\right)}{\Gamma \left(3^{-1+p}
   x\right)}=\frac{3^{3^n (1+n) x} \left(3^{-1-n}\right)^{3^n x} \Gamma
   \left(3^n x\right)}{\Gamma (x)}
\end{equation}
\end{example}
\begin{example}
\begin{multline}
\prod _{p=1}^n \left(\frac{\left(1+m^{3^{-p}}+m^{2\times 3^{-p}}\right)
   \left(-1+r^{3^{-p}}\right)^2}{\left(-1+m^{3^{-p}}\right)^2 \left(1+r^{3^{-p}}+r^{2\times
   3^{-p}}\right)}\right)^{3^p}=\frac{(-1+m)^3
   \left(\frac{-1+r^{3^{-n}}}{-1+m^{3^{-n}}}\right)^{3^{1+n}}}{(-1+r)^3}
\end{multline}
\end{example}
\begin{example}
\begin{equation}
\prod _{p=1}^n \left(1+q^{3^p}+q^{2\times 3^p}\right)=\frac{-1+q^{3^{1+n}}}{-1+q^3}
\end{equation}
\end{example}
\begin{example}
\begin{multline}
\prod _{p=1}^n \frac{1}{\left(1+q^{3^{-p}}+q^{2\times 3^{-p}}\right) \left(1+q^{3^p}+q^{2\times
   3^p}\right)}=\frac{\left(1+q+q^2\right) \left(-1+q^{3^{-n}}\right)}{-1+q^{3^{1+n}}}
\end{multline}
\end{example}
\begin{example}
\begin{multline}
\prod _{p=1}^n \exp \left(\frac{1}{2} 3^{-\frac{3}{2}+p} \left(\left(3 i+\sqrt{3}\right)
   \Phi'\left(-\frac{1}{2} i \left(-i+\sqrt{3}\right),0,1+3^{1-p} z\right)\right.\right. \\ \left.\left.
+\left(-3
   i+\sqrt{3}\right) \Phi'\left(\frac{1}{2} i \left(i+\sqrt{3}\right),0,1+3^{1-p}
   z\right)\right)\right)\\
=3^{\frac{3}{4}+\frac{1}{4} 3^n (-1+2 n)+z} \left(3^{1+n}\right)^{-\frac{3^n}{2}-z} (2 \pi
   )^{\frac{1}{2} \left(-1+3^n\right)} z \Gamma (z) \Gamma \left(1+3^{-n} z\right)^{-3^n}
\end{multline}
\end{example}
\begin{example}
\begin{equation}
\prod _{p=1}^n \frac{1+q^{3^{-p}}+q^{2\times 3^{-p}}}{1+q^{3^p}+q^{2\times 3^p}}=\frac{(-1+q)
   \left(-1+q^3\right)}{\left(-1+q^{3^{-n}}\right) \left(-1+q^{3^{1+n}}\right)}
\end{equation}
\end{example}
\begin{example}
Telescoping series involving the gamma function see Eq. (5.5.5) in \cite{dlmf}.
\begin{equation}
\prod _{p=1}^n \frac{\Gamma \left(\frac{1}{6}+\frac{3^{-p} x}{2}\right) \Gamma \left(\frac{5}{6}+\frac{3^{-p}
   x}{2}\right)}{2 \pi }=\frac{3^{\frac{3}{4} \left(-1+3^{-n}\right) x} \Gamma \left(\frac{1+x}{2}\right)}{\Gamma
   \left(\frac{1}{2} \left(1+3^{-n} x\right)\right)}
\end{equation}
\end{example}
\begin{example}
\begin{equation}
\prod _{p=1}^{\infty } \frac{\Gamma \left(\frac{1}{6}+\frac{3^{-p} x}{2}\right) \Gamma
   \left(\frac{5}{6}+\frac{3^{-p} x}{2}\right)}{2 \pi }=\frac{3^{-\frac{3 x}{4}} \Gamma
   \left(\frac{1+x}{2}\right)}{\sqrt{\pi }}
\end{equation}
\end{example}
\begin{example}
\begin{multline}
\prod _{p=1}^n \left(\frac{\left(1+q^{3^p}\right)^2 \left(1+q^{3^p}+q^{2\times
   3^p}\right)}{\left(-1+q^{3^p}\right)^2 \left(1-q^{3^p}+q^{2\times
   3^p}\right)}\right)^{3^{-p}}=-\frac{\left(1+q^3\right)
   }{-1+q^3}\left(\frac{1-q^{3^{1+n}}}{1+q^{3^{1+n}}}\right)^{3^{-n}}
\end{multline}
\end{example}
\begin{example}
\begin{equation}
\prod _{p=1}^{\infty } \left(\frac{\left(1+q^{3^p}\right)^2 \left(1+q^{3^p}+q^{2\times
   3^p}\right)}{\left(-1+q^{3^p}\right)^2 \left(1-q^{3^p}+q^{2\times
   3^p}\right)}\right)^{3^{-p}}=\frac{1+q^3}{1-q^3}
\end{equation}
\end{example}
\begin{example}
\begin{equation}
\prod _{p=1}^n \left(q^{2\times 3^{-p}}-q^{3^{-p}}+1\right)=\frac{q+1}{q^{3^{-n}}+1}
\end{equation}
\end{example}
\begin{example}
\begin{equation}
\prod _{p=1}^{\infty } \left(1-q^{3^{-p}}+q^{2\times 3^{-p}}\right)=\frac{1+q}{2}
\end{equation}
\end{example}
\begin{example}
\begin{equation}
\prod _{p=1}^{\infty } \prod _{k=1}^n \left(1-\left(-a
   q^k\right)^{3^{-p}}+\left(-a q^k\right)^{2\times 3^{-p}}\right)=\frac{2^{-n}
  }{1-a} (a;q)_{1+n}
\end{equation}
\end{example}
\begin{figure}[H]
\includegraphics[scale=0.5]{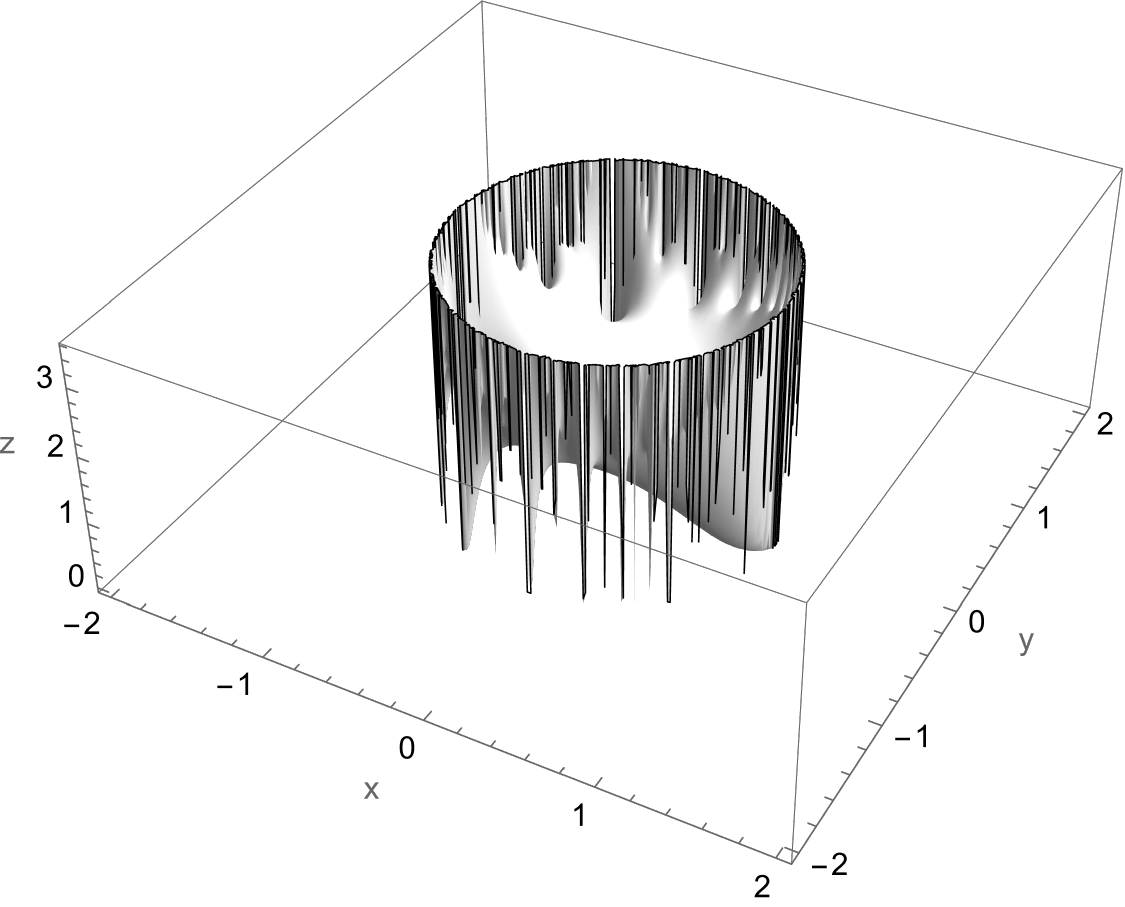}
\caption{Plot of  $f(z)=(z;z)_{\infty }$, $z\in\mathbb{C}$.}
   \label{fig:fig2}
\end{figure}
\vspace{-6pt}
\begin{example}
\begin{multline}
\prod _{p=1}^n \exp \left(\sqrt[3]{-1} 3^{-p}
   \left(-\Phi'\left(\sqrt[3]{-1},0,1+3^p a\right)+\sqrt[3]{-1}
   \Phi'\left(-(-1)^{2/3},0,1+3^p a\right)\right)\right)\\
=\frac{2^{\frac{1}{2}
   \left(1+3^{-n}\right)} 3^{\frac{3}{4} \left(-1+3^{-n}\right)} \Gamma \left(\frac{1+a}{2}\right)
   \left(\frac{\Gamma \left(1+\frac{3^n a}{2}\right)}{\Gamma \left(\frac{1}{2} \left(1+3^n
   a\right)\right)}\right)^{3^{-n}}}{a \Gamma \left(\frac{a}{2}\right)}
\end{multline}
\end{example}
\begin{example}
\begin{multline}
\prod _{p=1}^{\infty } \exp \left(\sqrt[3]{-1} 3^{-p}
   \left(-\Phi'\left(\sqrt[3]{-1},0,1+3^p a\right)+\sqrt[3]{-1}
   \Phi'\left(-(-1)^{2/3},0,1+3^p a\right)\right)\right)\\
=\frac{\Gamma
   \left(\frac{1+a}{2}\right)}{\sqrt{2} 3^{3/4} \Gamma \left(1+\frac{a}{2}\right)}
\end{multline}
\end{example}
\begin{example}
\begin{multline}
\prod _{p=1}^n \left(\frac{\left(1+q^{3^{-p}}\right)^2 \left(1+q^{3^{-p}}+q^{2\times
   3^{-p}}\right)}{\left(-1+q^{3^{-p}}\right)^2 \left(1-q^{3^{-p}}+q^{2\times 3^{-p}}\right)}\right)^{3^p}=-\frac{(-1+q)^3}{(1+q)^3}\left(-1-\frac{2}{-1+q^{3^{-n}}}\right)^{3^{1+n}}
\end{multline}
\end{example}
\begin{example}
\begin{equation}
\prod _{p=1}^n \frac{\Gamma \left(\frac{1}{6}+3^p x\right) \Gamma \left(\frac{5}{6}+3^p x\right)}{2 \pi
   }=\frac{3^{\frac{1}{2} (-9) \left(-1+3^n\right) x} \Gamma \left(\frac{1}{2}+3^{1+n} x\right)}{\Gamma
   \left(\frac{1}{2}+3 x\right)}
\end{equation}
\end{example}
\begin{example}
\begin{equation}
\sum _{p=1}^{\infty } \frac{9^p q^{3^{-1+p}} \left(1+q^{3^{-1+p}}
   \left(-4+q^{3^{-1+p}}\right)\right)}{\left(1+q^{3^{-1+p}} \left(-1+q^{3^{-1+p}}\right)\right)^2}=\frac{9
   q}{(1+q)^2}
\end{equation}
\end{example}
\begin{example}
\begin{equation}
\prod _{p=1}^n \left(1-q^{3^{-1+p}}+q^{2\times 3^{-1+p}}\right)=\frac{1+q^{3^n}}{1+q}
\end{equation}
\end{example}
\begin{example}
\begin{multline}
\prod _{p=1}^n \exp \left(\sqrt[3]{-1} 3^{-1+p}
   \left(-\Phi'\left(\sqrt[3]{-1},0,1+2\ 3^{1-p}
   x\right)+\sqrt[3]{-1} \Phi'\left(-(-1)^{2/3},0,1+2\
   3^{1-p} x\right)\right)\right)\\
=\frac{2^{\frac{1}{2}-\frac{3^n}{2}}
   3^{\frac{1}{4} \left(3+3^n (-3+2 n)\right)} \Gamma (1+x) \left(\frac{3^{n/2}
   \Gamma \left(\frac{1}{2}+3^{-n} x\right)}{x \Gamma \left(3^{-n}
   x\right)}\right)^{3^n}}{\Gamma \left(\frac{1}{2}+x\right)}
\end{multline}
\end{example}
\begin{example}
\begin{multline}
\prod_{p=0}^{\infty}\left(1+\frac{2 (m-r) (\alpha +2 p \alpha )}{(m+\alpha +2 p \alpha )
   (r-(1+2 p) \alpha )}\right)^{(-1)^p \cos ((1+2 p) x)}\\
=\exp \left(\frac{1}{2}
   i \left(e^{\frac{i m (\pi -2 x)}{2 \alpha }} \Phi \left(-e^{\frac{i m \pi
   }{\alpha }},1,\frac{1}{2}-\frac{x}{\pi }\right)+e^{\frac{i m (\pi +2 x)}{2
   \alpha }} \Phi \left(-e^{\frac{i m \pi }{\alpha }},1,\frac{1}{2}+\frac{x}{\pi
   }\right)\right.\right. \\ \left.\left.
-e^{\frac{i r (\pi -2 x)}{2 \alpha }} \Phi \left(-e^{\frac{i \pi 
   r}{\alpha }},1,\frac{1}{2}-\frac{x}{\pi }\right)-e^{\frac{i r (\pi +2 x)}{2
   \alpha }} \Phi \left(-e^{\frac{i \pi  r}{\alpha }},1,\frac{1}{2}+\frac{x}{\pi
   }\right)\right)\right)
\end{multline}
\end{example}
\begin{example}
\begin{multline}
\prod _{p=0}^{\infty } \exp \left(\frac{e^{i p \pi } \alpha  \left(2 m-i
   \left(m^2-(\alpha +2 p \alpha )^2\right) \log (a)\right) \sin (x+2 p
   x)}{(-m+\alpha +2 p \alpha )^2 (m+\alpha +2 p \alpha )^2}\right)\\
=\exp
   \left(\frac{\pi  \sec \left(\frac{m \pi }{2 \alpha }\right) \left(2 m x \cos
   \left(\frac{m x}{\alpha }\right)+\sin \left(\frac{m x}{\alpha }\right)
   \left(-2 \alpha +2 i m \alpha  \log (a)+m \pi  \tan \left(\frac{m \pi }{2
   \alpha }\right)\right)\right)}{8 m^2 \alpha }\right)
\end{multline}
\end{example}
\begin{example}
\begin{multline}
\prod _{p=0}^{\infty } \exp \left(\frac{2 (-1)^p e^{-i (1+2 p) \alpha }
   \cos ((1+2 p) x) \left(E_1(-i (1+2 p) \alpha )+e^{2 i (\alpha +2 p \alpha )}
   E_1(i (1+2 p) \alpha )\right)}{(1+2 p) \pi }\right)\\
=\frac{2 \pi  \Gamma
   \left(\frac{3 \pi -2 x+2 \alpha }{4 \pi }\right) \Gamma \left(\frac{3 \pi +2
   x+2 \alpha }{4 \pi }\right)}{\alpha  \Gamma \left(\frac{\pi -2 x+2 \alpha }{4
   \pi }\right) \Gamma \left(\frac{\pi +2 x+2 \alpha }{4 \pi }\right)}
\end{multline}
\end{example}
\begin{example}
\begin{multline}
\prod _{p=0}^{\infty } \left(\frac{m-(1+2 p) \alpha }{m+\alpha +2 p
   \alpha }\right)^{\frac{(-1)^p \cos (x+2 p x)}{(-m+\alpha +2 p \alpha
   )^2}}\\
=\exp \left(\frac{\pi }{8 \alpha ^2} \left((-2+2 \gamma -i \pi ) \sec \left(\frac{m
   \pi }{2 \alpha }\right) \left(-2 x \sin \left(\frac{m x}{\alpha }\right)+\pi 
   \cos \left(\frac{m x}{\alpha }\right) \tan \left(\frac{m \pi }{2 \alpha
   }\right)\right)\right.\right. \\ \left.
\left.
-4 i \pi  \left(e^{\frac{i m (\pi -2 x)}{2 \alpha }}
   \Phi'\left(-e^{\frac{i m \pi }{\alpha
   }},-1,\frac{1}{2}-\frac{x}{\pi }\right)+e^{\frac{i m (\pi +2 x)}{2 \alpha }}
   \Phi'\left(-e^{\frac{i m \pi }{\alpha
   }},-1,\frac{1}{2}+\frac{x}{\pi }\right)\right)\right)\right)\\
   \pi ^{\frac{\pi  \sec \left(\frac{m \pi }{2 \alpha }\right) \left(-2 x \sin
   \left(\frac{m x}{\alpha }\right)+\pi  \cos \left(\frac{m x}{\alpha }\right)
   \tan \left(\frac{m \pi }{2 \alpha }\right)\right)}{4 \alpha ^2}}
   \left(\frac{1}{\alpha }\right)^{\frac{\pi  \sec \left(\frac{m \pi }{2 \alpha
   }\right) \left(-2 x \sin \left(\frac{m x}{\alpha }\right)+\pi  \cos
   \left(\frac{m x}{\alpha }\right) \tan \left(\frac{m \pi }{2 \alpha
   }\right)\right)}{4 \alpha ^2}}
\end{multline}
\end{example}
\begin{example}
\begin{multline}
\prod _{p=0}^{\infty } \exp \left(\frac{4 i (-1)^{\frac{1}{4}+p}  e^{-\frac{1}{2} (i x)} \cos
   (x+2 p x) \left(\frac{a^{4 i m (1+2 p)} \Gamma (0,i m (1+4 p) \log (a))}{1+4 p}+\frac{\Gamma (0,-i m (3+4 p) \log
   (a))}{3+4 p}\right)}{a^{i m (3+4 p)}\pi }\right)\\
 (i m (1+4 p))^{-\frac{4 i (-1)^{\frac{1}{4}+p} e^{-\frac{1}{2} (i x)} \cos (x+2 p
   x)}{\pi +4 p \pi }} (-i m (3+4 p))^{-\frac{4 i (-1)^{\frac{1}{4}+p} e^{-\frac{1}{2} (i x)} \cos (x+2 p x)}{(3+4 p)
   \pi }}\\
 (i m (1+4 p) \log (a))^{\frac{4 i (-1)^{\frac{1}{4}+p} e^{-\frac{1}{2} (i x)} \cos (x+2 p x)}{\pi +4 p \pi
   }} (-i m (3+4 p) \log (a))^{\frac{4 i (-1)^{\frac{1}{4}+p} e^{-\frac{1}{2} (i x)} \cos (x+2 p x)}{(3+4 p) \pi
   }}\\
=\frac{ \left(\frac{\Gamma \left(\frac{\pi -2 x+4 m \log
   (a)}{8 \pi }\right)}{\Gamma \left(\frac{5 \pi -2 x+4 m \log (a)}{8 \pi }\right)}\right)^{e^{-i x}}
   \left(\frac{\Gamma \left(\frac{3 \pi -2 x+4 m \log (a)}{8 \pi }\right)}{\Gamma \left(\frac{7 \pi -2 x+4 m \log
   (a)}{8 \pi }\right)}\right)^{-i e^{-i x}} \Gamma \left(\frac{\pi +2 x+4 m \log (a)}{8 \pi }\right)
   \left(\frac{\Gamma \left(\frac{3 \pi +2 x+4 m \log (a)}{8 \pi }\right)}{\Gamma \left(\frac{7 \pi +2 x+4 m \log
   (a)}{8 \pi }\right)}\right)^{-i}}{\Gamma \left(\frac{5 \pi +2 x+4 m \log (a)}{8 \pi }\right)}\\
\left(\frac{1}{m}\right)^{\left(-\frac{1}{2}+\frac{i}{2}\right) \left(1+e^{-i x}\right)} (2 \pi
   )^{\left(-\frac{1}{2}+\frac{i}{2}\right) \left(1+e^{-i x}\right)}
\end{multline}
\end{example}
\begin{example}
\begin{equation}
\prod _{p=0}^{\infty } \left(1-\frac{2 (1+2 p) (m-r) \beta }{(m+(1+2 p) \beta ) (-r+(1+2 p) \beta
   )}\right)^{(-1)^p}=\frac{\tan \left(\frac{\pi  (r+\beta )}{4 \beta }\right)}{\tan \left(\frac{\pi  (m+\beta )}{4
   \beta }\right)}
\end{equation}
\end{example}
\begin{example}
\begin{equation}
\prod _{p=0}^{\infty } e^{-\frac{8 (-1)^p (1+2 p)}{3+16 p+16 p^2}} \left(\frac{(3+4 p)^2}{(1+4
   p)^2}\right)^{(-1)^p}=\left(3+2 \sqrt{2}\right) e^{-\frac{\pi }{\sqrt{2}}}
\end{equation}
\end{example}
\begin{example}
\begin{equation}
\prod _{p=0}^{\infty } \exp \left(-\frac{2 (-1)^p (1+2 p)}{2+9 p (1+p)}\right) \left(\frac{-(2+3 p)^2}{-(1+3
   p)^2}\right)^{\frac{2 (-1)^p}{3}}=3^{2/3} e^{-\frac{4 \pi }{9 \sqrt{3}}}
\end{equation}
\end{example}
\begin{example}
\begin{multline}
\prod _{p=0}^{\infty } \exp \left(-\frac{8 (-1)^p (1+2 p)}{(3+8 p) (5+8 p)}\right) \left(\frac{-(5+8 p)^2}{-(3+8
   p)^2}\right)^{\frac{(-1)^p}{2}}=\left(-1+\sqrt{2}+\sqrt{4-2 \sqrt{2}}\right) e^{-\frac{\pi }{4
   \sqrt{2+\sqrt{2}}}}
\end{multline}
\end{example}
\begin{example}
\begin{multline}
\prod _{p=0}^{\infty } \exp \left(-\frac{160 (-1)^p (1+2 p)}{-39+100 p+100 p^2}\right) \left(\frac{(13+10
   p)^2}{(3-10 p)^2}\right)^{(-1)^p}\\
\frac{1}{5} \left(5-\sqrt{5}+\sqrt{5 \left(5-2 \sqrt{5}\right)}\right)^2
   e^{\frac{8}{5} \left(-1+\sqrt{5}\right) \pi }
\end{multline}
\end{example}
\begin{example}
\begin{multline}
\prod _{p=0}^{\infty } e^{-\frac{12 (-1)^p (1+2 p)}{5+36 p+36 p^2}} \left(\frac{(5+6 p)^2}{(1+6
   p)^2}\right)^{\frac{(-1)^p}{2}}\left(2+\sqrt{3}\right) e^{-\frac{2 \pi }{3}}
\end{multline}
\end{example}
\begin{example}
\begin{multline}
\prod _{p=0}^{\infty } \exp \left(-\Gamma \left(0,-\frac{1}{2} i (\pi +2 p \pi )-\log (q)\right)-\Gamma
   \left(0,\frac{1}{2} i (\pi +2 p \pi +2 i \log (q))\right)\right)=1+q^2
\end{multline}
\end{example}
\begin{example}
\begin{multline}
\prod _{p=0}^{\infty } \exp \left(-\Gamma \left(0,\frac{i \pi  (-2 m+\beta +2 p \beta )}{2 \beta }\right)-\Gamma
   \left(0,-\frac{i \pi  (2 m+\beta +2 p \beta )}{2 \beta }\right)\right. \\ \left.
+\Gamma \left(0,\frac{i \pi  (-2 r+\beta +2 p \beta
   )}{2 \beta }\right)+\Gamma \left(0,-\frac{i \pi  (2 r+\beta +2 p \beta )}{2 \beta }\right)\right)\\
=e^{\frac{i \pi 
   (m-r)}{\beta }} \cos \left(\frac{m \pi }{\beta }\right) \sec \left(\frac{\pi  r}{\beta }\right)
\end{multline}
\end{example}
\begin{example}
\begin{multline}
\prod _{p=0}^{\infty } (-2 m+\beta +2 p \beta )^{\frac{(-1)^p}{-2 m+\beta +2 p \beta }}\\
 (-2 r+\beta +2 p \beta
   )^{-\frac{(-1)^p}{-2 r+\beta +2 p \beta }} (-2 m-(1+2 p) \beta )^{\frac{(-1)^p}{2 m+\beta +2 p \beta }} (-2 r-(1+2 p) \beta
   )^{-\frac{(-1)^p}{2 r+\beta +2 p \beta }}\\
=\exp \left(\frac{\pi }{\beta } \left(-\frac{1}{2} \left(\gamma +\log \left(\frac{i \pi }{\beta
   }\right)\right) \left(\sec \left(\frac{m \pi }{\beta }\right)-\sec \left(\frac{\pi  r}{\beta }\right)\right)\right.\right. \\ \left.\left.
+e^{\frac{i m \pi
   }{\beta }} \Phi'\left(-e^{\frac{2 i m \pi }{\beta }},0,\frac{1}{2}\right)-e^{\frac{i \pi  r}{\beta }}
   \Phi'\left(-e^{\frac{2 i \pi  r}{\beta }},0,\frac{1}{2}\right)\right)\right)
\end{multline}
\end{example}
\begin{example}
\begin{multline}
\prod _{p=0}^{\infty } \exp \left(\frac{(-1)^p e^{-i a \left(\frac{1}{2}+p\right) \pi  \beta } \left(E_1\left(-\frac{1}{2} i a
   (1+2 p) \pi  \beta \right)+e^{i a (1+2 p) \pi  \beta } E_1\left(\frac{1}{2} i a (1+2 p) \pi  \beta \right)\right)}{(1+2 p) \pi
   }\right)\\
=\frac{2 \Gamma \left(\frac{1}{4} (3+a \beta )\right)}{\sqrt{a} \sqrt{\beta } \Gamma \left(\frac{1}{4} (1+a \beta
   )\right)}
\end{multline}
\end{example}
\begin{example}
\begin{multline}
\prod _{p=0}^{\infty } \exp \left(\frac{(-1)^p e^{-i \left(\frac{1}{2}+p\right) \pi  x} \left(E_1\left(-\frac{1}{2} i (1+2 p)
   \pi  x\right)+e^{i (1+2 p) \pi  x} E_1\left(\frac{1}{2} i (1+2 p) \pi  x\right)\right)}{(1+2 p) \pi }\right)\\
=\frac{2 \Gamma
   \left(\frac{3+x}{4}\right)}{\sqrt{x} \Gamma \left(\frac{1+x}{4}\right)}
\end{multline}
\end{example}
\begin{example}
Application of Gautschi's Inequality using equation (5.6.4) in \cite{dlmf}.
\begin{multline}
\frac{2 \sqrt{x}}{\sqrt{1+4 x}}<\prod _{p=0}^{\infty } \exp \left(\frac{E_1\left(-\frac{1}{2} i (1+2 p) \pi  (1+4 x)\right)+e^{i
   (1+2 p) \pi  (1+4 x)} E_1\left(\frac{1}{2} i (1+2 p) \pi  (1+4 x)\right)}{(-1)^{-p} e^{i \left(\frac{1}{2}+p\right) \pi  (1+4 x)}
   (1+2 p) \pi }\right)\\
<\frac{2 \sqrt{1+x}}{\sqrt{1+4 x}}
\end{multline}
\end{example}
\begin{example}
\begin{multline}
\prod _{p=0}^{\infty } \exp \left(\frac{2 (-1)^p \sqrt{2} a^{-\left(\left(\frac{1}{4}+p\right) \beta
   \right)} \left(\frac{\Gamma \left(0,-\frac{1}{4} (1+4 p) \beta  \log (a)\right)}{1+4 p}+\frac{a^{\beta +2 p \beta } \Gamma
   \left(0,\frac{1}{4} (3+4 p) \beta  \log (a)\right)}{3+4 p}\right)}{\pi }\right)\\
=\frac{8 i \pi  \left(\frac{\Gamma \left(\frac{\pi -i
   \beta  \log (a)}{8 \pi }\right)}{\Gamma \left(\frac{5}{8}-\frac{i \beta  \log (a)}{8 \pi }\right)}\right)^{-1+i} \left(\frac{\Gamma
   \left(\frac{3}{8}-\frac{i \beta  \log (a)}{8 \pi }\right)}{\Gamma \left(\frac{7}{8}-\frac{i \beta  \log (a)}{8 \pi
   }\right)}\right)^{-1-i}}{\log (a) \beta }
\end{multline}
\end{example}
\begin{example}
\begin{multline}
\prod _{p=0}^{\infty } \exp \left(-\frac{(-1)^p  \left(e^{\frac{i \pi  x}{2}} (3+2 p) E_1\left(\frac{1}{2} i (1-2
   p) \pi  x\right)+e^{\frac{1}{2} i (3+4 p) \pi  x} (-1+2 p) E_1\left(\frac{1}{2} i (3+2 p) \pi  x\right)\right)}{e^{i p \pi  x}(-3+4 p (1+p)) \pi
   }\right)\\
=\frac{2 \Gamma \left(\frac{3+x}{4}\right)}{\sqrt{x} \Gamma \left(\frac{1+x}{4}\right)}
\end{multline}
\end{example}
\begin{example}
\begin{multline}
\prod _{p=0}^{\infty } \exp \left(\frac{(-1)^p e^{-i \left(\frac{3}{2}+p\right) \pi  x} \left(e^{i (1+2 p) \pi  x} (3+2 p)
   E_1\left(\frac{1}{2} i (-1+2 p) \pi  x\right)+(-1+2 p) E_1\left(-\frac{1}{2} i (3+2 p) \pi  x\right)\right)}{(-3+4 p (1+p)) \pi
   }\right)\\
=\frac{\sqrt{x} \Gamma \left(\frac{1+x}{4}\right)}{2 \Gamma \left(\frac{3+x}{4}\right)}
\end{multline}
\end{example}
\begin{example}
\begin{multline}
\prod _{p=0}^{\infty } \exp \left(\Gamma \left(0,\pi  \left(-\frac{i}{2}-i p+\frac{m}{\beta }\right)\right)+\Gamma \left(0,\pi 
   \left(\frac{i}{2}+i p+\frac{m}{\beta }\right)\right)\right)\\
=\frac{1}{2} \left(1+\tanh \left(\frac{m \pi }{\beta
   }\right)\right)
\end{multline}
\end{example}
\begin{example}
\begin{multline}
\prod _{j=1}^n \left(\prod _{p=0}^{\infty } \exp \left(-\Gamma \left(0,-\frac{1}{2} i (\pi +2 p \pi )-\log
   \left(q^j\right)\right)-\Gamma \left(0,\frac{1}{2} i \left(\pi +2 p \pi +2 i \log
   \left(q^j\right)\right)\right)\right)\right)\\
=\frac{1}{2} \left(-1;q^2\right){}_{1+n}
\end{multline}
\end{example}
\begin{example}
\begin{equation}
\prod _{p=1}^n \left(1+\frac{2 q^{2^{-p}}}{1+q^{2^{1-p}}}\right)=\frac{(-1+q)
   \left(1+q^{2^{-n}}\right)}{(1+q) \left(-1+q^{2^{-n}}\right)}
\end{equation}
\end{example}
\begin{example}
\begin{multline}
\prod _{p=1}^n \pi ^{-2^{-1-p}} \left(\frac{2^{\frac{1}{2} \left(-1+2^p z\right)} \Gamma \left(\frac{1}{4}
   \left(2+2^p z\right)\right)^2}{\Gamma \left(\frac{1}{2} \left(1+2^p
   z\right)\right)}\right)^{2^{-p}}\\
=2^{\frac{1}{2} \left(-1+2^{-n}+n z\right)} \pi ^{\frac{1}{2}
   \left(-1+2^{-n}\right)} \Gamma \left(\frac{1+z}{2}\right) \Gamma \left(\frac{1}{2} \left(1+2^n
   z\right)\right)^{-2^{-n}}
\end{multline}
\end{example}
\begin{example}
\begin{equation}
\prod _{p=1}^n \left(1+\frac{2 q^{2^{-1+p}}}{1+q^{2^p}}\right)^{2^{-p}}=(1+q)
   \left(1+q^{2^n}\right)^{-2^{-n}}
\end{equation}
\end{example}
\begin{example}
An infinite product two-ways.
\begin{equation}
\prod _{p=1}^{\infty } \left(1+\frac{2 q^{2^{-1+p}}}{1+q^{2^p}}\right)^{2^{-p}} = \begin{cases}
  1+\frac{1}{q} ; & q \in\mathbb{C} \\
  1+q; & q \in\mathbb{R}
\end{cases}
\end{equation}
\end{example}
\begin{figure}[H]
\includegraphics[scale=0.5]{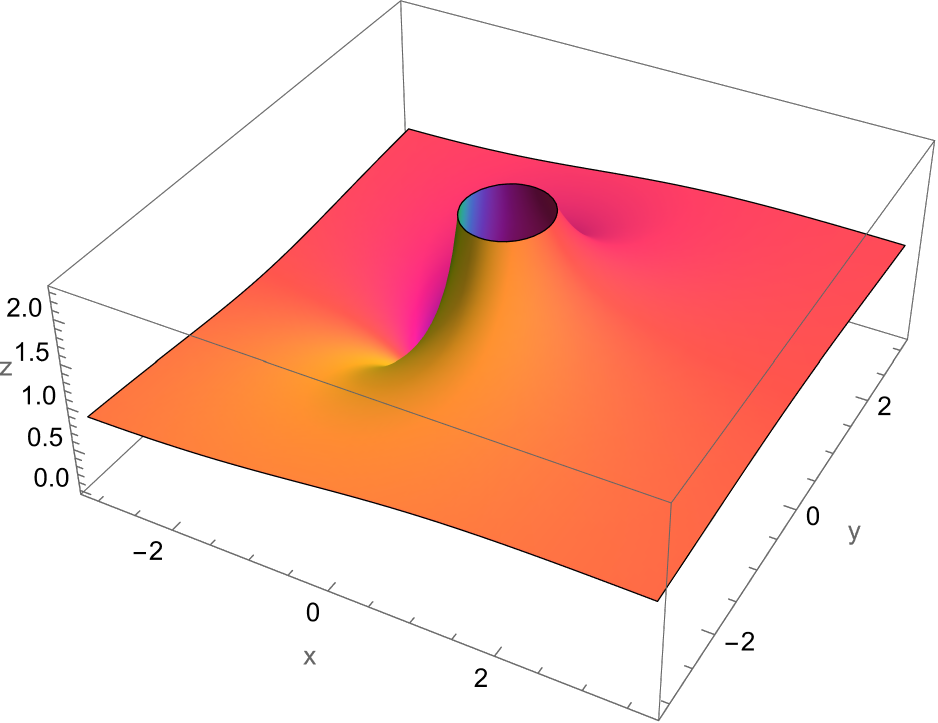}
\caption{Plot of  $f(q)=1+\frac{1}{q}$, $x\in\mathbb{C}$.}
   \label{fig:fig2}
\end{figure}
\vspace{-6pt}
\begin{example}
\begin{equation}
\prod _{j=1}^m \prod _{p=1}^{\infty } \left(1+\frac{2 q^{j 2^{-1+p}}}{1+q^{j
   2^p}}\right)^{2^{-p}}= \begin{cases}
  \frac{1}{2} \left(-1;\frac{1}{q}\right)_{1+m}  & q \in\mathbb{C} \\
  \frac{1}{2} (-1;q)_{1+m}  & q \in\mathbb{R}
\end{cases}
\end{equation}
\end{example}
\begin{example}
\begin{equation}
\prod _{p=1}^n \frac{\Gamma \left(\frac{1}{4} \left(2+2^{-p} z\right)\right)^2}{\sqrt{\pi } \Gamma
   \left(\frac{1}{2} \left(1+2^{-p} z\right)\right)}=\frac{2^{n+\frac{1}{2} \left(-1+2^{-n}\right) z} \Gamma
   \left(\frac{z}{4}\right) \Gamma \left(\frac{1}{4} \left(2+2^{-n} z\right)\right)}{\Gamma \left(2^{-2-n} z\right)
   \Gamma \left(\frac{2+z}{4}\right)}
\end{equation}
\end{example}
\subsubsection{Examples derived using equation (\ref{eq_3116})}
\begin{example}
\begin{equation}
\prod _{p=1}^{2 n} \exp \left((-1)^{-n+p} e^{\frac{i p \pi }{n}} q^{\frac{1+n}{2}} \Phi \left(-e^{\frac{i p
   \pi }{n}} \sqrt{q},1,1+n\right)\right)=1-q^n
\end{equation}
\end{example}
\begin{example}
\begin{multline}
\prod _{n=1}^{\infty } \left(\prod _{p=1}^{2 n} \exp \left((-1)^{-n+p} e^{\frac{i p \pi }{n}}
   q^{\frac{1+n}{2}} \Phi \left(-e^{\frac{i p \pi }{n}} \sqrt{q},1,1+n\right)\right)\right)=(q;q)_{\infty }
\end{multline}
\end{example}
\begin{figure}[H]
\includegraphics[scale=0.5]{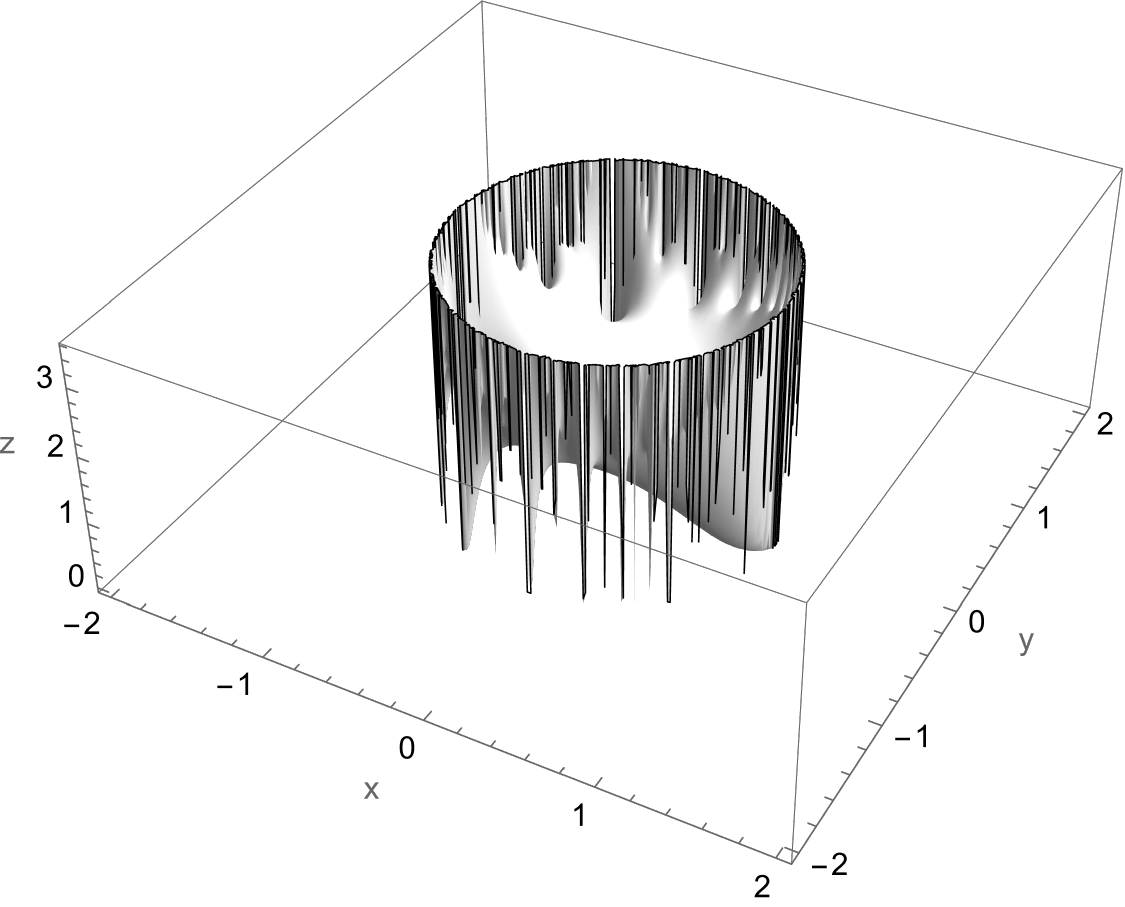}
\caption{Plot of  $f(z)=(z;z)_{\infty }$, $z\in\mathbb{C}$.}
   \label{fig:fig2}
\end{figure}
\vspace{-6pt}
\begin{example}
\begin{equation}
\prod _{p=1}^{2 n} \left(1+e^{\frac{i p \pi }{n}}
   q\right)^{(-1)^p}=\left(\frac{1-q^n}{1+q^n}\right)^{(-1)^n}
\end{equation}
\end{example}
\begin{example}
\begin{equation}
\prod _{n=1}^{\infty } \left(\prod _{p=1}^{2 n} \left(1+e^{\frac{i p \pi }{n}}
   q\right)^{(-1)^p}\right)=\frac{\left(-q;q^2\right){}_{\infty } \left(q^2;q^2\right){}_{\infty
   }}{\left(q;q^2\right){}_{\infty } \left(-q^2;q^2\right){}_{\infty }}
\end{equation}
\end{example}
\begin{example}
\begin{multline}
\prod _{p=1}^{2 n} \exp \left(\frac{(-1)^{-n+p} e^{-\frac{i (-1+n-2 p) \pi }{2 n}}
   \Phi'\left(-e^{\frac{i (1+2 p) \pi }{2 n}},0,a\right)}{2 n}\right)=\frac{2 \sqrt{n}
   \Gamma \left(\frac{-1+a+3 n}{4 n}\right)}{\Gamma \left(\frac{-1+a+n}{4 n}\right)}
\end{multline}
\end{example}
\begin{example}
\begin{equation}
\prod _{p=1}^{2 n} \left(\frac{e^{-i (m-r)} \cos \left(\frac{p \pi }{2 n}+r\right)}{\cos \left(m+\frac{p \pi
   }{2 n}\right)}\right)^{(-1)^p}=\left(\frac{\tan (n r)}{\tan (m n)}\right)^{(-1)^n}
\end{equation}
\end{example}
\begin{example}
\begin{equation}
\prod _{n=1}^{\infty } \left(\prod _{p=1}^{2 n} \left(\frac{1+e^{\frac{i p \pi }{n}} r}{1+e^{\frac{i p \pi
   }{n}} m}\right)^{(-1)^{-n+p}}\right)=\frac{(-1;m)_{\infty } (r;r)_{\infty }}{(-1;r)_{\infty } (m;m)_{\infty
   }}
\end{equation}
\end{example}
\begin{example}
A telescoping series similar to Eq. (5.5.6) in \cite{dlmf}.
\begin{multline}
\prod _{p=0}^n 3^{\frac{1}{2} \left(-2^p+3 x\right)} (2 \pi )^{-2^p} \left(\frac{\Gamma \left(2^{-1-p}
   x\right) \Gamma \left(\frac{1}{6}+2^{-1-p} x\right) \Gamma \left(\frac{5}{6}+2^{-1-p} x\right)}{\Gamma
   \left(3\times 2^{-1-p} x\right)}\right)^{2^p}\\
=\frac{3^{\frac{1}{2}-2^n} \left(2^{1+n}\right)^{-2 x} \Gamma (3 x)
   \left(\frac{\Gamma \left(2^{-1-n} x\right)}{\Gamma \left(3\times 2^{-1-n} x\right)}\right)^{2^{1+n}}}{\Gamma
   (x)}
\end{multline}
\end{example}
\subsection{Series}
\begin{example}
\begin{multline}
\alpha\sum _{p=1}^{\infty } \left(\left(\coth \left(\alpha ^2 p\right)-1\right) \left((\log (a)+2 \alpha 
   p)^k-(\log (a)-2 \alpha  p)^k\right)\right. \\ \left.
+\frac{i \beta  \left(\coth \left(\beta ^2 p\right)-1\right) \left((\log
   (a)-2 i \beta  p)^k-(\log (a)+2 i \beta  p)^k\right)}{\alpha }\right)\\
=-2^{k+1} \alpha ^{k+1} \zeta
   \left(-k,\frac{\log (a)}{2 \alpha }+1\right)+2^{k+1} (i \beta )^{k+1} \zeta \left(-k,1-\frac{i \log (a)}{2 \beta
   }\right)\\
-\left(\log ^{k-1}(a) (2 k+(\alpha -i \beta ) \log (a))\right)
\end{multline}
\end{example}
\begin{example}
\begin{multline}
\sum _{p=1}^{\infty } \frac{\left(\coth \left(\alpha ^2 p\right)-1\right) \cosh (2 \alpha  m p)-\left(\coth \left(\beta
   ^2 p\right)-1\right) \cos (2 \beta  m p)}{p}\\
=\log (\sinh (\alpha  m) \csc (\beta  m))-m^2
\end{multline}
\end{example}
\begin{example}
Ramanujan-Berndt type transformation in terms of the cotangent function.
\begin{multline}
\sum _{p=1}^{\infty } \left(\frac{\alpha ^2 p \coth \left(\alpha ^2
   p\right)}{a^2-\alpha ^2 p^2}+\frac{\beta ^2 p \coth \left(\beta ^2
   p\right)}{a^2+\beta ^2 p^2}\right)+\frac{1}{2 a^2}-\log \left(\frac{i \beta
   }{\alpha }\right)\\
=\begin{cases}
    \frac{1}{2} \pi  (\cot (a \beta )-i \coth (a \alpha
   )), & \text{if $Re(a)\geq 0$},\\
    -\frac{1}{2} \pi  (\cot (a \beta )-i \coth (a \alpha
   )), & \text{if $Re(a)<0$}.
  \end{cases}
\end{multline}
\end{example}
\begin{example}
\begin{multline}
\sum _{p=1}^{\infty } \left(\alpha  \left(\coth \left(\alpha ^2
   p\right)-1\right) \left(\frac{\log (a+2 \alpha  p)}{a+2 \alpha 
   p}-\frac{\log (a-2 \alpha  p)}{a-2 \alpha  p}\right)\right. \\ \left.
+i \beta  \left(\coth
   \left(\beta ^2 p\right)-1\right) \left(\frac{\log (a-2 i \beta  p)}{a-2 i
   \beta  p}-\frac{\log (a+2 i \beta  p)}{a+2 i \beta 
   p}\right)\right)\\
=\frac{1}{2} \left(\frac{-4+\log (a) (-2 a \alpha +2 i a
   \beta +4)}{\log ^2\left(e^a\right)}-2 \log (2 i \beta ) \psi
   ^{(0)}\left(1-\frac{i a}{2 \beta }\right)\right. \\ \left.
+\log (\alpha ) \log (4 \alpha )-2\log ^2(2 i \beta )+\log (i \beta ) \log (4 i \beta )
+2 \log^2(2)\right)\\
+\log (2 \alpha ) \psi ^{(0)}\left(\frac{a}{2 \alpha
   }+1\right)-\gamma _1\left(\frac{a}{2 \alpha }+1\right)+\gamma_1\left(1-\frac{i a}{2 \beta }\right) 
\end{multline}
\end{example}
\begin{example}
\begin{multline}
\sum _{p=1}^{\infty } \left(\left(\frac{e^{-\frac{p \pi  \beta }{\alpha
   }}}{i a-2 p \beta }-\frac{e^{\frac{p \pi  \beta }{\alpha }}}{i a+2 p \beta
   }\right) \beta  \left(-1+\coth \left(p \beta ^2\right)\right)\right. \\ \left.
   -\frac{4 (-1)^p p
   \alpha ^2 \left(-1+\coth \left(p \alpha ^2\right)\right)}{a^2-4 p^2 \alpha
   ^2}\right)\\
=\frac{1}{2} \left(\frac{-2 i a \pi +4 \alpha -2 a \alpha ^2+2 i a
   \alpha  \beta }{a^2 \alpha }\right. \\ \left.
   -H_{\frac{1}{4} \left(-2+\frac{a}{\alpha
   }\right)}+H_{\frac{a}{4 \alpha }}+2 e^{-\frac{\pi  \beta }{\alpha }} \Phi
   \left(e^{-\frac{\pi  \beta }{\alpha }},1,1-\frac{i a}{2 \beta
   }\right)\right)
\end{multline}
\end{example}
\begin{example}
\begin{equation}
\sum _{p=1}^{\infty } \left(\coth \left(\beta ^2 p\right)-1\right)
   \sinh \left(\frac{\pi  \beta  p}{\alpha }\right)=\frac{1}{2} \left(\coth
   \left(\frac{\pi  \beta }{2 \alpha }\right)-\frac{\pi }{\alpha  \beta
   }\right)
\end{equation}
\end{example}
\begin{example}
\begin{multline}
\sum_{p=1}^{\infty}p \left((-1)^p\left(\coth \left(\frac{\pi  p}{\beta }\right)-1\right)\frac{\pi   }{\beta }+\pi  \beta  \cosh (\pi\beta  p) (\coth (\pi  \beta  p)-1)\right)\\
=\frac{1}{4} \left(\frac{\pi }{\beta }-\pi  \beta 
   \text{csch}^2\left(\frac{\pi  \beta }{2}\right)-2\right)
\end{multline}
\end{example}
\begin{example}
\begin{equation}
\sum _{p=1}^{\infty } p^2 \left(\coth \left(\beta ^2 p\right)-1\right) \sinh \left(\frac{\pi  \beta  p}{\alpha }\right)=\frac{1}{8}
   \sinh \left(\frac{\pi  \beta }{\alpha }\right) \text{csch}^4\left(\frac{\pi  \beta }{2 \alpha }\right)
\end{equation}
\end{example}
\begin{example}
\begin{multline}
\sum _{p=1}^{\infty } \left(e^{-2 i m p} (3 \pi -2 i p)^{-s} (\coth (p)-1)-e^{2 i m p} (3 \pi +2 i p)^{-s}
   (\coth (p)-1)\right. \\ \left.
+i \pi  \left(e^{-2 \pi  m p} (3 \pi -2 \pi  p)^{-s} \left(\coth \left(\pi ^2 p\right)-1\right)-e^{2
   \pi  m p} (2 \pi  p+3 \pi )^{-s} \left(\coth \left(\pi ^2 p\right)-1\right)\right)\right)\\
=-i \left(e^{2 i m} (2
   i)^{1-s} \Phi \left(e^{2 i m},s,1-\frac{3 i \pi }{2}\right)-e^{2 \pi  m} (2 \pi )^{1-s} \Phi \left(e^{2 m \pi
   },s,\frac{5}{2}\right)\right. \\ \left.
+(3 \pi )^{-s-1} (2 s-3 \pi  (2 m+\pi -i))\right)
\end{multline}
\end{example}
\begin{example}
\begin{equation}
\sum _{p=1}^{\infty } p \left(\coth (p)+\pi ^2 \coth \left(\pi ^2 p\right)-\pi ^2-1\right)=\frac{1}{12}
   \left(\pi ^2-5\right)
\end{equation}
\end{example}
\begin{example}
\begin{multline}
\sum _{p=1}^{\infty } \left(\frac{4 i p \left(\coth \left(\pi ^2 p\right)-1\right)}{9-4 p^2}+\frac{1-\coth
   (p)}{3 \pi +2 i p}+\frac{\coth (p)-1}{3 \pi -2 i p}\right)\\
=\frac{-4 i+6 \pi -9 \pi ^3+6 i \pi ^2 \left(-7+3 \gamma
   +\log (64)-3 \log (\pi )+3 \psi ^{(0)}\left(1-\frac{3 i \pi }{2}\right)\right)}{18 \pi ^2}
\end{multline}
\end{example}
\begin{example}
\begin{multline}
\sum _{p=1}^{\infty } \left(\tan ^{-1}\left(\frac{2 p}{3 \pi }\right) (-1+\coth (p))+\pi  \tanh
   ^{-1}\left(\frac{2 p}{3}\right) \left(-1+\coth \left(p \pi ^2\right)\right)\right)\\
=\log
   \left(\sqrt{\frac{\left(\frac{3 \pi ^3}{16}\right)^{\pi } \left(\frac{\Gamma \left(\frac{3 i \pi }{2}\right)}{\Gamma
   \left(-\frac{1}{2} (3 i \pi )\right)}\right)^i}{e^{\frac{2}{3 \pi }+\frac{5 \pi }{2}}}}\right)
\end{multline}
\end{example}
\begin{example}
\begin{equation}
\sum _{p=1}^{\infty } p \left(\alpha ^2 \left(\coth \left(\alpha ^2 p\right)-1\right)+\beta ^2 \left(\coth
   \left(\beta ^2 p\right)-1\right)\right)=\frac{1}{12} \left(\alpha ^2+\beta ^2-6\right)
\end{equation}
\end{example}
\begin{example}
\begin{multline}
\log
   \left(\sqrt[4]{-\frac{i \alpha }{\beta }} e^{\frac{-a \alpha +i a \beta +2}{4 a^2}}\right)-\sum _{p=1}^{\infty } \left(\frac{\alpha ^2 p \left(\coth \left(\alpha ^2 p\right)-1\right)}{4 \alpha ^2
   p^2-a^2}-\frac{\beta ^2 p \left(\coth \left(\beta ^2 p\right)-1\right)}{a^2+4 \beta ^2 p^2}\right)\\
=\frac{1}{4}
   \left(H_{-\frac{i a}{2 \beta }}-H_{\frac{a}{2 \alpha }}\right)
\end{multline}
\end{example}
\begin{example}
\begin{multline}
\sum_{n=-\infty}^{\infty}\frac{n \sinh (2 m) e^{\alpha  n}}{\cosh (2 m)+\cosh (2 \alpha  n)}\\
=-\sum_{n=-\infty}^{\infty}\frac{\pi  \cosh (m)
   \text{sech}\left(\frac{\pi ^2 n}{\alpha }\right) \left(\pi  \tanh \left(\frac{\pi ^2 n}{\alpha }\right) \sin
   \left(\frac{2 \pi  m n}{\alpha }\right)-2 m \cos \left(\frac{2 \pi  m n}{\alpha }\right)\right)}{2 \alpha
   ^2}
\end{multline}
\end{example}
\begin{example}
\begin{multline}
\sum_{n=-\infty}^{\infty}\frac{\text{sech}^2(m+n \alpha )}{e^{n \alpha }}
=\sum_{n=-\infty}^{\infty}\frac{e^m \pi  \text{sech}\left(\frac{n \pi ^2}{\alpha
   }\right) \left(\alpha  \cos \left(\frac{2 m n \pi }{\alpha }\right)-2 n \pi  \sin \left(\frac{2 m n \pi }{\alpha
   }\right)\right)}{\alpha ^2}
\end{multline}
\end{example}
\begin{example}
\begin{equation}\label{eq:dieckmann}
\sum_{n=-\infty}^{\infty}e^{\alpha  (-n)} \tanh (m+\alpha  n)=\frac{\pi  e^m }{\alpha }\sum_{n=-\infty}^{\infty}\text{sech}\left(\frac{\pi ^2 n}{\alpha }\right) \cos
   \left(\frac{2 \pi  m n}{\alpha }\right)
\end{equation}
\end{example}
\begin{example}
A Dieckmann series in terms of the $q$-digamma function derived using \cite{dieckmann} and equation (\ref{eq:dieckmann}).
\begin{equation}
\sum_{n=-\infty}^{\infty}\frac{\tanh (n \pi  \alpha )}{e^{n \pi  \alpha }}=-\frac{1}{\alpha }+\frac{2 i \left(\psi _{e^{-\frac{\pi
   }{\alpha }}}^{(0)}\left(-\frac{1}{2} (i \alpha )\right)-\psi _{e^{-\frac{\pi }{\alpha }}}^{(0)}\left(\frac{i
   \alpha }{2}\right)\right)}{\pi }
\end{equation}
\end{example}
\begin{example}
\begin{equation}
\sum _{p=1}^{\infty } (-1)^p \left(\frac{\log (2 p-1)}{1-2 p}+\frac{\log (2 p+1)}{2 p+1}\right)=-\frac{1}{2} \pi  \log \left(\frac{2 e^{\gamma } \pi  \Gamma
   \left(-\frac{1}{4}\right)^2}{9 \Gamma \left(-\frac{3}{4}\right)^2}\right)
\end{equation}
\end{example}
\begin{example}
\begin{multline}
\sum _{p=1}^{\infty } \left(\frac{e^{i m p}}{(-1+a+p)^s}-\frac{e^{-i m
   p}}{(-1+a-p)^s}+i \left(\frac{e^{m p}}{(-1+a-i p)^s}-\frac{e^{-m p}}{(-1+a+i
   p)^s}\right)\right) \\
(-1+\coth (p \pi ))\\
=\frac{-i (-1+a) (m-(1+i) \pi )+s}{\pi 
   (-1+a)^{s+1}}+2 i e^{-m-\frac{i \pi  s}{2}} \Phi \left(e^{-m},s,(1+i)-i
   a\right)\\
-2 e^{i m} \Phi \left(e^{i m},s,a\right)
\end{multline}
\end{example}
\begin{example}
\begin{multline}
\left(1-\frac{2 p}{2 a+p-1}\right)^{\frac{1}{2} (-1)^p (\coth (p \pi )-1)} \left(\frac{(2 a+i p-1)^{e^{-2 p \pi
   }}}{2 a-i p-1}\right)^{\frac{1}{2} i \text{csch}(p \pi )}\\
=\frac{\exp \left(\frac{1}{2} \left(\frac{\pi }{e^{\pi
   }-1}+\frac{1}{(2 a-1) \pi }+2 i e^{-\pi } \Phi'\left(e^{-\pi },0,(1+i)-2 i
   a\right)\right)\right) \Gamma (a)}{\sqrt{a-\frac{1}{2}} \Gamma \left(a-\frac{1}{2}\right)}
\end{multline}
\end{example}
\begin{example}
\begin{multline}
\sum _{p=1}^{\infty } \left(\frac{e^{i m p}}{(-1+a+p)^s}-\frac{e^{-i m
   p}}{(-1+a-p)^s}+i \left(\frac{e^{m p}}{(-1+a-i p)^s}-\frac{e^{-m p}}{(-1+a+i
   p)^s}\right)\right) \\
(-1+\coth (p \pi ))\\
=\frac{-i (-1+a) (m-(1+i) \pi )+s}{\pi 
   (-1+a)^{s+1}}+2 i e^{-m-\frac{i \pi  s}{2}} \Phi \left(e^{-m},s,(1+i)-i
   a\right)\\
-2 e^{i m} \Phi \left(e^{i m},s,a\right)
\end{multline}
\end{example}
\begin{example}
\begin{multline}
\sum_{p=1}^{\infty}\coth (\pi  p)-1) \left(\log \left(\frac{a-p-1}{a+p-1}\right)+i \log \left(\frac{a+i p-1}{a-i p-1}\right)\right)\\
=\log
   \left(-\frac{(2 \pi )^{1-i} (-i (a-1))^{2 i}  \exp \left(2 i \text{log$\Gamma $}(-i (a-1))+\pi  \left(-\frac{1}{2}+i
   a\right)+\frac{1}{\pi  (a-1)}\right)}{(a-1)^{1+i}\Gamma (a-1)^2}\right)
\end{multline}
\end{example}
\
\begin{example}
\begin{multline}
\sum _{p=1}^{\infty } (\coth (\pi  p)-1) \left(\frac{\log (a-p-1)}{a-p-1}-\frac{\log (a+p-1)}{a+p-1}\right. \\ \left.
+i \left(\frac{\log (a-i p-1)}{-a+i
   p+1}+\frac{\log (a+i p-1)}{a+i p-1}\right)\right)\\
=-2 \gamma _1((1+i)-i a)+2 \gamma
   _1(a)\\
   +\frac{-\pi ^3 (a-1)^2+(-4+(4-4 i) \pi  (a-1)) \log (a-1)+4}{4 \pi  (a-1)^2}\\
+i \pi  \psi ^{(0)}((1+i)-i
   a)
\end{multline}
\end{example}
\begin{example}
\begin{multline}
\sum _{p=1}^{\infty } (\coth (\pi  p)-1) \left(\sin \left(\frac{\pi  p}{3}\right)+\sinh \left(\frac{\pi 
   p}{3}\right)\right)=\frac{1}{6} \left(-1-3 \sqrt{3}+3 \coth \left(\frac{\pi }{6}\right)\right)
\end{multline}
\end{example}
\begin{example}
\begin{multline}
\sum _{p=1}^{\infty } p (\coth (\pi  p)-1) \left(\cos \left(\frac{\pi  p}{3}\right)+\cosh \left(\frac{\pi 
   p}{3}\right)\right)=\frac{1}{4} \left(5-\frac{2}{\pi }-\coth ^2\left(\frac{\pi }{6}\right)\right)
\end{multline}
\end{example}
\begin{example}
Infinite series involving the incomplete beta function $B_{x}(a,b)$.
\begin{multline}
\sum _{p=1}^{\infty } p \left(\frac{(-1)^p}{p^2-a^2}-\frac{e^{-\pi p}}{a^2+p^2}\right) \coth (\pi 
   p)\\
=\frac{1}{2} \left(\frac{1}{\pi  a^2}+2 e^{-\pi } \Phi \left(e^{-\pi },1,1-i a\right)-\pi  \csc (\pi  a)\right. \\ \left.
-e^{-i
   \pi  a} B_{e^{-\pi }}(1-i a,0)+e^{i \pi  a} B_{e^{-\pi }}(i a+1,0)\right)
\end{multline}
\end{example}
\begin{example}
\begin{multline}
\sum _{p=1}^{\infty } \left(-\frac{4 e^{-i \pi  p}}{(9-2 p)^2}+\frac{e^{i \pi  p}}{\left(p+\frac{9}{2}\right)^2}-\frac{4 i
   e^{\pi  p}}{(2 p+9 i)^2}+\frac{4 i e^{\pi  (-p)}}{(-2 p+9 i)^2}\right) (\coth (\pi  p)-1)\\
\approx -8 C-2 i e^{-\pi } \Phi \left(e^{-\pi
   },2,1-\frac{9 i}{2}\right)+\frac{726052}{99225}+\frac{16}{729 \pi }
\end{multline}
\end{example}
\begin{example}
\begin{multline}
\sum _{p=1}^{\infty }\left( \Gamma (0,i (p-m) \pi )+\Gamma (0,-i (m+p) \pi )\right)\\
=\log \left(\frac{1}{2}
   (1+i \cot (\pi  m))\right)-\Gamma (0,-i m \pi )
\end{multline}
\end{example}
\begin{example}
\begin{multline}
\sum _{p=1}^{\infty } \left(\frac{\Gamma (1+k,i (p-m) \pi )}{(p-m)^{k+1}}+\frac{\Gamma (1+k,-i (m+p) \pi )}{(-m-p)^{k+1}}\right)\\
=\frac{(-1)^k \Gamma (1+k,-i m \pi )}{m^{k+1}}+\Gamma
   (1+k) \left(-(-1)^k \zeta (1+k,1+m)+\zeta (1+k,-m)\right)
\end{multline}
\end{example}
\begin{example}
\begin{multline}
\sum _{p=1}^{\infty } (-1)^p (E_2(i (1-2 p) \pi )+E_2(i (2 p+1) \pi ))=-4 i C+(1+4 i)+i \pi  \Gamma (0,i \pi )
\end{multline}
\end{example}
\begin{example}
\begin{equation}
\sum _{p=1}^{\infty } \frac{p \left(a^4+p^4\right) \coth (\pi 
   p)}{\left(a^4-p^4\right)^2}=\frac{\pi ^3 a^3 \csc ^2(\pi  a)-i \pi ^3 a^3
   \text{csch}^2(\pi  a)-2}{8 \pi  a^4}
\end{equation}
\end{example}
\begin{figure}[H]
\includegraphics[scale=0.5]{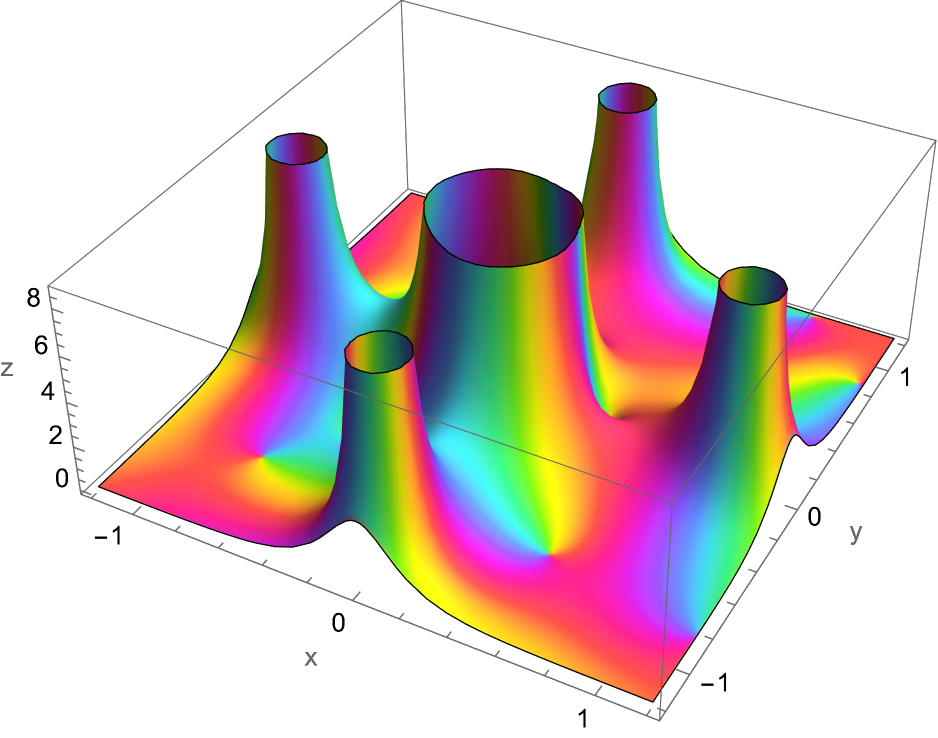}
\caption{Plot of  $f(a)=\frac{\pi ^3 a^3 \csc ^2(\pi  a)-i \pi ^3 a^3 \text{csch}^2(\pi  a)-2}{8 \pi  a^4}$, $a\in\mathbb{C}$.}
   \label{fig:fig2}
\end{figure}
\vspace{-6pt}
\begin{example}
\begin{multline}
\sum _{p=1}^{\infty } \frac{(-2 p+i) E_1\left(\frac{\pi }{2}-i p \pi
   \right)+(2 p+i) E_1\left(i \pi  p+\frac{\pi }{2}\right)}{4 p^2+1}\\
=i
   \left(-\pi  Li_{0}'\left(e^{-\pi
   }\right)+\text{Ei}\left(-\frac{\pi }{2}\right)+\frac{\pi  \log (2)}{e^{\pi
   }-1}\right)
\end{multline}
\end{example}
\begin{example}
\begin{multline}
\sum _{p=1}^{\infty } \left(\frac{E_1\left(\frac{1}{2} (2 i p+(1-i))
   \pi \right)}{(1+i) p-i}-\frac{E_1\left(\frac{1}{2} ((1-i)-2 i p) \pi
   \right)}{(1+i) p+i}\right)\\
=-(1+i) \pi 
   \left(Li_{0}'\left(-e^{-\pi }\right)+\frac{\log
   (2)}{1+e^{\pi }}\right)-i E_1\left(\left(\frac{1}{2}-\frac{i}{2}\right) \pi
   \right)
\end{multline}
\end{example}
\begin{example}
The generalized hypergeometric function $\, _3F_3(a,a,a;b,b,b;z)$.
\begin{multline}
\sum _{p=1}^{\infty } 2 \left(3 \pi  (-1-2 i p) \, _3F_3\left(1,1,1;2,2,2;\left(-i p-\frac{1}{2}\right) \pi \right)\right. \\ \left.
+3 \pi  (-1+2 i p) \, _3F_3\left(1,1,1;2,2,2;\left(i p-\frac{1}{2}\right) \pi
   \right)\right. \\ \left.
+6 \log \left(-\frac{i}{\pi }\right) \left(\text{Ei}\left(\left(-i
   p-\frac{1}{2}\right) \pi \right)+\text{Ei}\left(\left(i
   p-\frac{1}{2}\right) \pi \right)\right)\right. \\ \left.
+3 \left(\log ^2\left(\frac{1}{2}\pi  (1-2 i p)\right)+\log ^2\left(\frac{1}{2} \pi  (1+2 i p)\right)+2
   \gamma  \log \left(\frac{1}{2} \pi  (1-2 i p)\right)\right.\right. \\ \left.\left.
+2 \gamma  \log\left(\frac{1}{2} \pi  (1+2 i p)\right)\right)+\pi ^2+6 \gamma ^2\right)\\
=6
   \pi  \left(\, _3F_3\left(1,1,1;2,2,2;-\frac{\pi }{2}\right)+i
   \text{Ei}\left(-\frac{\pi }{2}\right)+\log (4)-i \log \left(e^{\pi
   }-1\right)+2 \log (\pi )\right)\\
-6 \left(2
   \text{PolyLog}^{(1,0)}\left(1,e^{-\pi }\right)+\log (\pi ) \left(-2
   \text{Ei}\left(-\frac{\pi }{2}\right)+2 \log \left(e^{\pi }-1\right)+\log
   (\pi )\right)\right. \\ \left.
+\log ^2(2)+\log (4) \left(\log \left(e^{\pi }-1\right)-\log
   (\pi )\right)\right)+(-1+6 i) \pi ^2\\
-6 \gamma  \left(\gamma +\log
   \left(\frac{\pi ^2}{4}\right)\right)
\end{multline}
\end{example}
\begin{example}
Sum product formula in terms of $q$-digamma function.
\begin{multline}
\sum _{k=1}^n 2 e^{\Gamma \left(0,\pi  \left(\frac{2 k}{n}+\frac{x}{\pi }\right)\right)} \prod
   _{p=1}^{\infty } \exp \left(\Gamma \left(0,i \pi  \left(p-i \left(\frac{2 k}{n}+\frac{x}{\pi
   }\right)\right)\right)+\Gamma \left(0,-i \pi  \left(p+i \left(\frac{2 k}{n}+\frac{x}{\pi
   }\right)\right)\right)\right)\\
=\frac{n }{2 \pi }\left(-\psi _{e^{\frac{2 \pi }{n}}}^{(0)}\left(\frac{n x}{2
   \pi }+1\right)+\psi _{e^{\frac{2 \pi }{n}}}^{(0)}\left(\frac{x n}{2 \pi }+n+1\right)-\psi
   _{e^{\frac{2 \pi }{n}}}^{(0)}\left(\frac{1}{2} n \left(\frac{x}{\pi }-i\right)+1\right)\right. \\ \left.
   +\psi
   _{e^{\frac{2 \pi }{n}}}^{(0)}\left(\frac{1}{2} n \left(\frac{x}{\pi
   }+(2-i)\right)+1\right)\right)
\end{multline}
\end{example}
\begin{example}
The digamma function $\psi ^{(0)}(z+1)$.
\begin{equation}
\psi ^{(0)}(z+1)=-\frac{\zeta '(1-z)}{\zeta (1-z)}-\frac{\zeta '(z)}{\zeta (z)}+\frac{1}{z}+\frac{1}{2} \pi  \tan \left(\frac{\pi  z}{2}\right)+\log (2 \pi )
\end{equation}
\end{example}
\begin{figure}[H]
\includegraphics[scale=0.7]{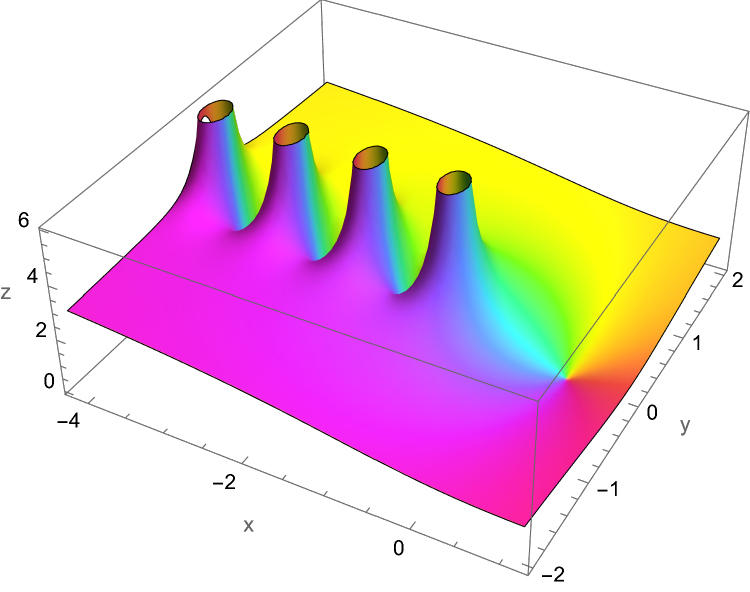}
\caption{Plot of  $f(z)=\psi ^{(0)}(z+1)$, $z\in\mathbb{C}$.}
   \label{fig:fig2}
\end{figure}
\vspace{-6pt}
\begin{example}
\begin{multline}
\sum _{p=0}^{\infty } \left(\Gamma \left(-1,\frac{i \pi  (i m-(2 p+1)
   x)}{2 x}\right)+\Gamma \left(-1,\frac{i \pi  (i m+(2 p+1) x)}{2
   x}\right)\right)\\
=\log \left(\frac{2 \cosh \left(\frac{\pi  m}{2 x}\right)
   \left(\sinh \left(\frac{\pi  m}{2 x}-\tanh \left(\frac{\pi  m}{2
   x}\right)\right)+\cosh \left(\frac{\pi  m}{2 x}-\tanh \left(\frac{\pi  m}{2
   x}\right)\right)\right)}{e}\right)
\end{multline}
\end{example}
\begin{example}
\begin{multline}
\sum _{p=0}^{\infty } (-1)^p \left(\frac{E_1\left(\frac{1}{2} i (4 p+1) \pi \right)}{4 p+1}+\frac{E_1\left(-\frac{1}{2} i (4 p+3) \pi \right)}{4 p+3}\right)\\
=\frac{1}{2} (-1)^{3/4}
   \pi  \left(\log \left(\Gamma \left(\frac{3}{8}\right)\right)-\log \left(2 \Gamma \left(\frac{7}{8}\right)\right)+i \log \left(\frac{2 \Gamma \left(\frac{9}{8}\right)}{\Gamma
   \left(\frac{5}{8}\right)}\right)\right)
\end{multline}
\end{example}
\begin{example}
\begin{multline}
\sum _{p=0}^{\infty } \left(\frac{E_1\left(\frac{1}{2} i (2 p-1) \pi \right)}{(2
   p-1)^3}-\frac{G_{2,3}^{3,0}\left(-\frac{1}{2} i (2 p+3) \pi \left |
\begin{array}{c}
 1,1 \\
 0,0,2 \\
\end{array}\right.
\right)}{(2 p+3)^3}\right)\\
=\frac{1}{16} i \left(-14 \pi  \zeta (3)+16+\pi ^3\right)
\end{multline}
\end{example}
\begin{example}. 
\begin{multline}
\sum_{p=0}^{\infty}\frac{(-1)^p \left(\frac{i}{a \pi }\right)^{2 p} E_{1+2 p}(x) (2)_{2 p}}{\Gamma (2+2 p)}\\
=-\frac{1}{2} a^2 \pi
   ^2 \left(\psi ^{(0)}\left(\frac{1}{2} (1+a \pi -x)\right)-\psi ^{(0)}\left(\frac{1}{2} (2+a \pi -x)\right)\right. \\ \left.
   -\psi
   ^{(0)}\left(\frac{1}{2} (a \pi +x)\right)+\psi ^{(0)}\left(\frac{1}{2} (1+a \pi +x)\right)\right)
\end{multline}
\end{example}
\begin{example}
\begin{multline}
\sum_{p=0}^{\infty}\frac{(-1)^p \left(\frac{i}{a \pi }\right)^{2 p} E_{1+2 p}(x) (2)_{2 p} (-1+\gamma +\psi ^{(0)}(2+2
   p))}{\Gamma (2+2 p)}\\
=-\frac{1}{2} a^2 \pi ^2 \left(-\psi ^{(0)}\left(\frac{1}{2} (1+a \pi -x)\right)-\log (2) \psi
   ^{(0)}\left(\frac{1}{2} (1+a \pi -x)\right)\right. \\ \left.
-\log \left(\frac{1}{a}\right) \psi ^{(0)}\left(\frac{1}{2} (1+a \pi
   -x)\right)+\log (\pi ) \psi ^{(0)}\left(\frac{1}{2} (1+a \pi -x)\right)\right. \\ \left.
+\psi ^{(0)}\left(\frac{1}{2} (2+a \pi
   -x)\right)+\log (2) \psi ^{(0)}\left(\frac{1}{2} (2+a \pi -x)\right)+\log \left(\frac{1}{a}\right) \psi
   ^{(0)}\left(\frac{1}{2} (2+a \pi -x)\right)\right. \\ \left.
-\log (\pi ) \psi ^{(0)}\left(\frac{1}{2} (2+a \pi -x)\right)+\psi
   ^{(0)}\left(\frac{1}{2} (a \pi +x)\right)+\log (2) \psi ^{(0)}\left(\frac{1}{2} (a \pi +x)\right)\right. \\ \left.
+\log
   \left(\frac{1}{a}\right) \psi ^{(0)}\left(\frac{1}{2} (a \pi +x)\right)-\log (\pi ) \psi ^{(0)}\left(\frac{1}{2}
   (a \pi +x)\right)-\psi ^{(0)}\left(\frac{1}{2} (1+a \pi +x)\right)\right. \\ \left.
-\log (2) \psi ^{(0)}\left(\frac{1}{2} (1+a \pi
   +x)\right)-\log \left(\frac{1}{a}\right) \psi ^{(0)}\left(\frac{1}{2} (1+a \pi +x)\right)\right. \\ \left.
+\log (\pi ) \psi
   ^{(0)}\left(\frac{1}{2} (1+a \pi +x)\right)+\gamma _1\left(\frac{1}{2} (1+a \pi -x)\right)\right. \\ \left.
-\gamma
   _1\left(\frac{1}{2} (2+a \pi -x)\right)-\gamma _1\left(\frac{1}{2} (a \pi +x)\right)+\gamma _1\left(\frac{1}{2}
   (1+a \pi +x)\right)\right)
\end{multline}
\end{example}
\begin{example}
\begin{equation}
\sum_{p=0}^{\infty}\frac{E_{2 p+1}\left(-\frac{1}{2}\right) (1-k)_{2 p}}{\Gamma (2 p+2)}=-\frac{2^{-k}
   \left(3^k-1\right)}{k}
\end{equation}
\end{example}
\begin{example}
\begin{multline}
\sum _{p=0}^{2 n} (-1)^{p+2^p} 4^{-p} \Phi \left((-1)^{2^{p+1}},2,\frac{1}{4}
   \left(2-2^{-p}\right)\right)\\
   -\sum _{p=0}^n (-1)^{4^p} 2^{1-4 p} \Phi \left((-1)^{2^{2 p+1}},2,\frac{1}{4}
   \left(2-4^{-p}\right)\right)\\
=16 C+(-1)^{2^{2 n+1}} 2^{-4 n} \Phi \left((-1)^{2^{2 n+1}},2,1-2^{-2
   (n+1)}\right)
\end{multline}
\end{example}
\begin{example}
\begin{multline}
\sum _{p=0}^n 2^{-p-1} \left(\text{log$\Gamma $}\left(i 2^p m\right)+\text{log$\Gamma $}\left(-i 2^p
   m-1\right)-\text{log$\Gamma $}\left(-i 2^p m-\frac{1}{2}\right)\right. \\ \left.
-\text{log$\Gamma $}\left(i 2^p
   m-\frac{1}{2}\right)+\log \left(\frac{m 2^p-i}{m^2 4^{p+1}+1}\right)\right)\\
=\log \left(\frac{i 2^{2^{1-n}-4}
   \pi ^{2^{-n-1}} \left(m 2^{n+1}-i\right)^{-2^{-n-1}} \text{csch}(\pi  m)}{m}\right)\\
-2^{-n-1}
   \left(\text{log$\Gamma $}\left(i 2^{n+1} m\right)+\text{log$\Gamma $}\left(-i 2^{n+1} m-1\right)\right)
\end{multline}
\end{example}
\begin{example}
\begin{multline}
-\sum _{p=0}^n 4^{-p} e^{i \pi  2^p} \left(\Phi \left(-e^{i 2^p \pi },2,1-2^{-p-1}\right)+\Phi \left(e^{i
   2^p \pi },2,1-2^{-p-1}\right)\right)\\
=8 C+2^{-2 n-1} e^{i \pi  2^{n+1}} \Phi \left(e^{i 2^{n+1} \pi
   },2,1-2^{-n-2}\right)
\end{multline}
\end{example}
\begin{example}
\begin{multline}
\sum _{p=0}^n \frac{4^p e^{i \pi  2^{-n+p-1}}}{\left(-1+e^{i \pi  2^{p-n}}\right)^3} \left(\left(-1+e^{i \pi  2^{p-n}}\right)^3
   \Phi'\left(-e^{i \pi  2^{-n+p-1}},-2,1\right)\right. \\ \left.
+\left(-1+e^{i \pi  2^{p-n}}\right)^3
   \Phi'\left(e^{i \pi  2^{-n+p-1}},-2,1\right)+2 \left(6 e^{i \pi  2^{p-n}}+e^{i \pi 
   2^{-n+p+1}}+1\right) \log \left(i 2^p\right)\right)\\
=2 e^{i \pi  2^{-n-1}}
   \Phi'\left(e^{i \pi  2^{-n-1}},-2,1\right)+\frac{7\ 2^{2 n+1} \zeta (3)}{\pi
   ^2}-\frac{1}{4} \pi  \cot \left(\pi  2^{-n-2}\right) \csc ^2\left(\pi  2^{-n-2}\right)
\end{multline}
\end{example}
\begin{example}
\begin{multline}
-\sum _{p=0}^n \frac{1}{2} i e^{-i \pi  2^{-n+p-1}} \csc ^2\left(\pi  2^{-n+p-1}\right) \left(\log \left(i
   2^p\right) \left(2^{p+1} \left(1+e^{i \pi  2^{p-n}}\right)\right.\right. \\ \left.\left.
+i \left(-1+e^{i \pi  2^{p-n}}\right) \log \left(e^{-i
   2^{n+1}}\right)\right)-2^p \left(-1+e^{i \pi  2^{p-n}}\right)^2 \right. \\ \left.
   \left(\Phi'\left(-e^{i \pi  2^{-n+p-1}},-1,1-i 2^{-p-1} \log \left(e^{-i
   2^{n+1}}\right)\right)\right.\right. \\ \left.\left.+\Phi'\left(e^{i \pi  2^{-n+p-1}},-1,1-i 2^{-p-1} \log
   \left(e^{-i 2^{n+1}}\right)\right)\right)\right)\\
=\frac{1}{2} \left(-i 2^{n+4}
   \Phi'\left(-1,-1,1-i 2^{-n-2} \log \left(e^{-i 2^{n+1}}\right)\right)\right. \\ \left.
-8 i e^{i \pi 
   2^{-n-1}} \Phi'\left(e^{i \pi  2^{-n-1}},-1,1-\frac{1}{2} i \log \left(e^{-i
   2^{n+1}}\right)\right)\right. \\ \left.
+i 2^{n+2} \log \left(i 2^{n+1}\right)+\pi  \csc ^2\left(\pi  2^{-n-2}\right)\right. \\ \left.-\log
   \left(e^{-i 2^{n+1}}\right) \left(\pi  \left(\cot \left(\pi  2^{-n-2}\right)+i\right)-2 \log \left(i
   2^{n+1}\right)\right)\right)
\end{multline}
\end{example}
\begin{example}
\begin{multline}
\lim_{n\to \infty}\sum_{j=1}^{n}\frac{1}{j}\left( \text{sech}\left(\frac{b j}{n}\right) \cosh \left(\frac{c j}{n}\right) \sinh \left(\frac{j m}{n}\right)+\text{sech}\left(\frac{b n}{j}\right) \cosh \left(\frac{c
   n}{j}\right) \sinh \left(\frac{m n}{j}\right)\right)\\
=\frac{1}{2} \log \left(\tan \left(\frac{\pi  (b-c+m)}{4 b}\right) \tan \left(\frac{\pi  (b+c+m)}{4 b}\right)\right)
\end{multline}
\end{example}
\begin{example}
\begin{multline}
\sum _{p=0}^{n-1} \frac{(-1)^p}{\sin ^2\left(\frac{m \pi }{\alpha }\right)}  \left((1+2 p) \sin \left(\frac{m (-1+4 n-2 p)
   \pi }{2 n \alpha }\right)+(1-4 n+2 p) \sin \left(\frac{m (1+2 p) \pi }{2 n \alpha }\right)\right)\\
=-\sec \left(\frac{m
   \pi }{2 n \alpha }\right) \tan \left(\frac{m \pi }{2 n \alpha }\right)
\end{multline}
\end{example}
\begin{example}
\begin{multline}
\sum _{p=0}^{n-1} e^{i p \pi  \left(1-\frac{m}{n \alpha }\right)} \left(e^{\frac{2 i m \pi }{\alpha }} (-1+4 n-2
   p)+e^{\frac{4 i m \pi }{\alpha }} (1+2 p)\right. \\ \left.
-e^{\frac{i m (1+2 p) \pi }{n \alpha }} (1+2 p)+e^{\frac{i m (1+2 n+2 p) \pi
   }{n \alpha }} (1-4 n+2 p)\right) \csc ^3\left(\frac{m \pi }{\alpha }\right)\\
=-\frac{2 i e^{\frac{i m (1+4 n) \pi }{2 n
   \alpha }} \tan \left(\frac{m \pi }{2 n \alpha }\right)}{\sin \left(\frac{m \pi }{\alpha }\right) \cos \left(\frac{m
   \pi }{2 n \alpha }\right)}
\end{multline}
\end{example}
\begin{example}
\begin{multline}
\sum _{p=0}^{n-1} 2^{-p} \left(2 \sec ^2\left(2^{1-p} m\right)-\sec ^2\left(2^{-p} m\right)-2 \sec
   ^2\left(2^{1-p} r\right)+\sec ^2\left(2^{-p} r\right)\right)\\
=2 \left(\sec ^2(2 m)-\sec ^2(2 r)+2^{-n}
   \left(-\sec ^2\left(2^{1-n} m\right)+\sec ^2\left(2^{1-n} r\right)\right)\right)
\end{multline}
\end{example}
\begin{example}
\begin{multline}
\sum _{p=0}^{n-1} 2^{-p}
   \left(-\frac{m^{2^{-p}}}{\left(1+m^{2^{-p}}\right)^2}+\frac{2}{\left(m^{-2^{-p}}+m^{2^{-p}}\right)^2}+\frac{
   r^{2^{-p}}}{\left(1+r^{2^{-p}}\right)^2}-\frac{2}{\left(r^{-2^{-p}}+r^{2^{-p}}\right)^2}\right)\\
=2
   \left(\frac{m^2}{\left(1+m^2\right)^2}-\frac{r^2}{\left(1+r^2\right)^2}+2^{-n}
   \left(-\frac{1}{\left(m^{-2^{-n}}+m^{2^{-n}}\right)^2}+\frac{1}{\left(r^{-2^{-n}}+r^{2^{-n}}\right)^2}\right
   )\right)
\end{multline}
\end{example}
\begin{example}
\begin{multline}
\sum _{p=0}^{\infty } 2^{-p}
   \left(-\frac{m^{2^{-p}}}{\left(1+m^{2^{-p}}\right)^2}+\frac{2}{\left(m^{-2^{-p}}+m^{2^{-p}}\right)^2}+\frac{
   r^{2^{-p}}}{\left(1+r^{2^{-p}}\right)^2}-\frac{2}{\left(r^{-2^{-p}}+r^{2^{-p}}\right)^2}\right)\\
=\frac{2
   m^2}{\left(1+m^2\right)^2}-\frac{2 r^2}{\left(1+r^2\right)^2}
\end{multline}
\end{example}
\begin{example}
\begin{multline}
\sum _{p=0}^{n-1} 2^{-p} \left(-\frac{2}{\left((i m)^{-2^{-p}}+(i m)^{2^{-p}}\right)^2}+\frac{(im)^{2^{-p}}}{\left(1+(i   m)^{2^{-p}}\right)^2}-\frac{m^{2^{-p}}}{\left(1+m^{2^{-p}}\right)^2}+\frac{2}{\left(m^{-2^{-p}}+m^{2^{-p}}\right)^2}\right)\\
=\frac{4 \left(m^2+m^6\right)}{\left(-1+m^4\right)^2}
\end{multline}
\end{example}
\begin{figure}[H]
\includegraphics[scale=0.7]{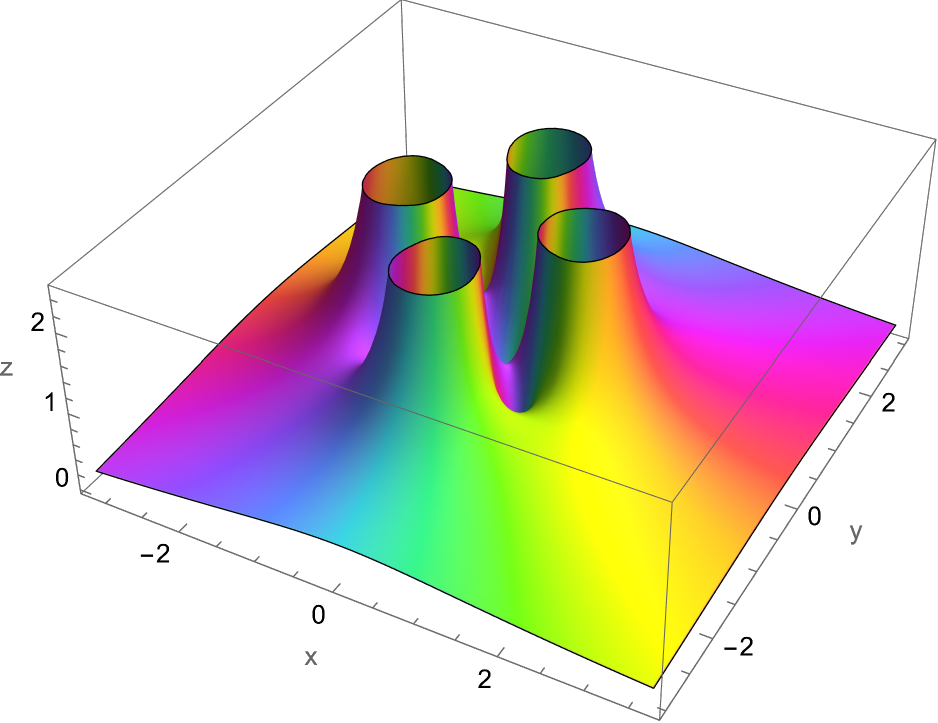}
\caption{Plot of  $f(m)=\frac{4 \left(m^6+m^2\right)}{\left(m^4-1\right)^2}$, $m\in\mathbb{C}$.}
   \label{fig:fig2}
\end{figure}
\vspace{-6pt}
\begin{example}
From Eq. (41.2.24) in \cite{hansen} and Eq. (\ref{prud24611}).
\begin{equation}
\sum _{k=0}^{n-1} (-1)^k \sinh \left(a \cosh ^{-1}\left(-\cos \left(\frac{\pi +2 k \pi }{2
   n}\right)\right)\right)=\frac{i \sin (a \pi )}{2 \cos \left(\frac{a \pi }{2 n}\right)}
\end{equation}
\end{example}
\begin{example}
\begin{multline}
\sum _{p=0}^{n-1} (-1)^p \exp \left(\frac{m \left(i (1-2 n) \pi -2 n \cosh ^{-1}\left(-\cos \left(\frac{\pi +2
   p \pi }{2 n}\right)\right)\right)}{2 c n}\right) \sqrt{-1+\cos \left(\frac{\pi +2 p \pi }{n}\right)}\\
 \csc \left(\frac{\pi +2 p \pi }{2 n}\right) \left(e^{\frac{2 m \cosh ^{-1}\left(-\cos \left(\frac{\pi +2 p \pi }{2
   n}\right)\right)}{c}} \Phi \left(e^{-\frac{2 i m \pi }{c}},-k,\frac{\pi +i \cosh ^{-1}\left(-\cos \left(\frac{\pi
   +2 p \pi }{2 n}\right)\right)-i c \log (a)}{2 \pi }\right)\right.\\ \left.
-\Phi \left(e^{-\frac{2 i m \pi }{c}},-k,-\frac{i \left(i
   \pi +\cosh ^{-1}\left(-\cos \left(\frac{\pi +2 p \pi }{2 n}\right)\right)+c \log (a)\right)}{2 \pi
   }\right)\right.\\ \left.
-e^{\frac{2 i m \pi }{c}} \Phi \left(e^{\frac{2 i m \pi }{c}},-k,\frac{\pi +i \cosh ^{-1}\left(-\cos
   \left(\frac{\pi +2 p \pi }{2 n}\right)\right)-i c \log (a)}{2 \pi }\right)\right.\\ \left.
+\exp \left(\frac{2 m \left(i \pi +\cosh
   ^{-1}\left(-\cos \left(\frac{\pi +2 p \pi }{2 n}\right)\right)\right)}{c}\right) \right.\\ \left.
\Phi \left(e^{\frac{2 i m \pi
   }{c}},-k,-\frac{i \left(i \pi +\cosh ^{-1}\left(-\cos \left(\frac{\pi +2 p \pi }{2 n}\right)\right)+c \log
   (a)\right)}{2 \pi }\right)\right)\\
=i 2^{\frac{1}{2}-k} n^{-k} \left(\Phi \left(-e^{-\frac{i m \pi }{c
   n}},-k,\frac{1}{2}-\frac{i c n \log (a)}{\pi }\right)\right.\\ \left.
-e^{\frac{i m \pi }{c n}} \Phi \left(-e^{\frac{i m \pi }{c
   n}},-k,\frac{1}{2}-\frac{i c n \log (a)}{\pi }\right)\right)
\end{multline}
\end{example}
\begin{example}
\begin{multline}
\sum _{p=0}^{n-1} \sqrt{-1+\cos \left(\frac{\pi +2 p \pi }{n}\right)} \csc \left(\frac{\pi +2 p \pi }{2
   n}\right) \left(\zeta \left(-k,\frac{\pi +b \pi -i \cosh ^{-1}\left(-\cos \left(\frac{\pi +2 p \pi }{2
   n}\right)\right)}{2 \pi }\right)\right. \\ \left.
+\zeta \left(-k,\frac{\pi +b \pi +i \cosh ^{-1}\left(-\cos \left(\frac{\pi +2 p \pi
   }{2 n}\right)\right)}{2 \pi }\right)\right)\\
=i 2^{\frac{1}{2}-k} n^{-k} \zeta \left(-k,\frac{1}{2}+b n\right)
\end{multline}
\end{example}
\begin{example}
\begin{multline}
\sum _{p=0}^{n-1} \frac{2^{-p} \left(2+4 \cos \left(2^{1-p} m\right)-2 \cos \left(2^{2-p} m\right)+2^p
   \log (a) \sin \left(2^{2-p} m\right)\right)}{\left(\cos \left(2^{-p} m\right)+\cos \left(3\ 2^{-p}
   m\right)\right)^2}\\
=\frac{2 \left(4-2^{3-n} \cos ^2(2 m)+4 \cos \left(2^{2-n} m\right)+\log (a) \left(\sin (4
   m)-\sin \left(2^{2-n} m\right)+\sin \left(4 \left(1-2^{-n}\right) m\right)\right)\right)}{2+2 \cos (4 m)+2
   \cos \left(2^{2-n} m\right)+\cos \left(4 \left(1-2^{-n}\right) m\right)+\cos \left(4 \left(1+2^{-n}\right)
   m\right)}
\end{multline}
\end{example}
\begin{example}
\begin{multline}
\sum _{p=0}^{n-1} \left(2^p \log \left(\frac{1-\tan ^2\left(2^{-p} x\right)}{1-\tan ^2\left(2^{-1-p}
   x\right)}\right)+\frac{2 \left(\sec \left(2^{-p} x\right) \tan \left(2^{-1-p} x\right)-\sec \left(2^{1-p}
   x\right) \tan \left(2^{-p} x\right)\right)}{\log (a)}\right)\\
=-\int_{\frac{x}{2}}^x \frac{4 \left(4-2^{3-n}
   \cos ^2(2 m)+4 \cos \left(2^{2-n} m\right)+\log (a) \left(\sin (4 m)-\sin \left(2^{2-n} m\right)+\sin
   \left(4 \left(1-2^{-n}\right) m\right)\right)\right)}{\left(2+2 \cos (4 m)+2 \cos \left(2^{2-n}
   m\right)+\cos \left(4 \left(1-2^{-n}\right) m\right)+\cos \left(4 \left(1+2^{-n}\right) m\right)\right) \log
   (a)} \, dm
\end{multline}
\end{example}
\begin{example}
\begin{multline}
\sum _{p=0}^{n-1} \frac{4^{-p} q^{2^{-p}} \left(-1+q^{2^{-p}} \left(5+q^{2^{-p}} \left(9+q^{2^{-p}}
   \left(11+q^{2^{-p}} \left(-11+q^{2^{-p}} \left(-9+q^{2^{-p}}
   \left(-5+q^{2^{-p}}\right)\right)\right)\right)\right)\right)\right)}{\left(1+q^{2^{1-p}}\right)^3
   \left(1+q^{2^{-p}}\right)^3}\\
=4
   \left(\frac{q^2}{\left(1+q^2\right)^3}-\frac{q^4}{\left(1+q^2\right)^3}+\frac{4^{-n} q^{2^{1-n}}
   \left(-1+q^{2^{1-n}}\right)}{\left(1+q^{2^{1-n}}\right)^3}\right)
\end{multline}
\end{example}
\begin{example}
\begin{multline}
\sum _{j=0}^{\infty } \sum _{n=0}^{\infty } \sum _{p=0}^{\infty } \frac{4^p
   (-t)^{j+2 p} u^n \binom{-1+n+j (-1+\alpha )}{-j+n} E_{1+2 p}(x) \Gamma (1+j+2
   p)}{\Gamma (1+j) \Gamma (2+2 p)}\\
=\frac{\log \left(\frac{\Gamma \left(\frac{1}{4}
   \left(2+\frac{1}{t}+(1-u)^{-\alpha } u-2 x\right)\right) \Gamma \left(\frac{1}{4}
   \left(2+\frac{1}{t}+(1-u)^{-\alpha } u+2 x\right)\right)}{\Gamma \left(\frac{1}{4}
   \left(4+\frac{1}{t}+(1-u)^{-\alpha } u-2 x\right)\right) \Gamma \left(\frac{1}{4}
   \left(\frac{1}{t}+(1-u)^{-\alpha } u+2 x\right)\right)}\right)}{2 t}
\end{multline}
\end{example}
\begin{example}
\begin{multline}
\sum _{p=0}^n \frac{4^{-p} b^{2^{-p}} \left(-1+b^{2^{1-p}}\right) m^{2^{-p}}
   \left(-1+m^{2^{1-p}}\right)}{\left(b^{2^{-p}}+m^{2^{-p}}\right)^2 \left(1+b^{2^{-p}} m^{2^{-p}}\right)^2}\\
=2^{-1-2 n} (b
   m)^{-2^{-n}} \left(\frac{(b m)^{2^{1-n}}}{\left(-1+b^{2^{-n}} m^{2^{-n}}\right)^3}+b^{2^{-n}} m^{2^{-n}} \left(\frac{2
   b^{2^{1-n}}}{\left(b^{2^{-n}}-m^{2^{-n}}\right)^2}+\frac{2}{-1+b^{-2^{-n}} m^{2^{-n}}}\right.\right. \\ \left.\left.
-\frac{1}{\left(-1+b^{2^{-n}}
   m^{2^{-n}}\right)^3}-\frac{3}{\left(-1+b^{2^{-n}} m^{2^{-n}}\right)^2}-\frac{2}{-1+b^{2^{-n}} m^{2^{-n}}}\right.\right. \\ \left.\left.
+2^{3+2 n} b^2 m^2
   \left(\frac{m^2}{\left(b^2-m^2\right)^3}-\frac{1}{\left(-1+b^2 m^2\right)^3}\right)+2^{3+2 n} b^4 m^2
   \left(-\frac{1}{\left(b^2-m^2\right)^3}\right.\right.\right. \\ \left.\left.\left.
+\frac{m^2}{\left(-1+b^2 m^2\right)^3}\right)\right)\right)
\end{multline}
\end{example}
\begin{example}
\begin{multline}
\sum _{p=0}^{\infty } \frac{4^{-p} b^{2^{-p}} \left(-1+b^{2^{1-p}}\right) m^{2^{-p}}
   \left(-1+m^{2^{1-p}}\right)}{\left(b^{2^{-p}}+m^{2^{-p}}\right)^2 \left(1+b^{2^{-p}} m^{2^{-p}}\right)^2}\\
=-\frac{4 b^2
   \left(-1+b^4\right) m^2 \left(-1+m^4\right)}{\left(-b^2+m^2+b^4 m^2-b^2 m^4\right)^2}+\frac{1}{\log
   ^2\left(\frac{b}{m}\right)}-\frac{1}{\log ^2(b m)}
\end{multline}
\end{example}
\begin{example}
\begin{multline}
\sum _{p=0}^{\infty } \frac{4^{-p} \left(-1+\left(\frac{1}{q}\right)^{2^{1-p}}\right) \left(-1+(i q)^{2^{1-p}}\right)
   \left(\frac{1}{q}\right)^{2^{-p}} (i q)^{2^{-p}}}{\left(\left(\frac{1}{q}\right)^{2^{-p}}+(i q)^{2^{-p}}\right)^2
   \left(1+\left(\frac{1}{q}\right)^{2^{-p}} (i q)^{2^{-p}}\right)^2}\\
=\frac{4}{\pi
   ^2}-\frac{\left(-1+q^4\right)^2}{\left(1+q^4\right)^2}+\frac{1}{\log ^2\left(-\frac{i}{q^2}\right)}
\end{multline}
\end{example}
\begin{example}
\begin{multline}
-\sum _{p=0}^{\infty } \frac{16^{-p} b^{2^{-1-p}} \left(-1+b^{2^{-p}}\right) m^{2^{-1-p}} \left(-1+m^{2^{-p}}\right)}{\left(b^{2^{-1-p}}+m^{2^{-1-p}}\right)^4 \left(1+b^{2^{-1-p}} m^{2^{-1-p}}\right)^4} \\
\left(m^{2^{-p}}+b^{2^{1-p}} m^{2^{-p}}-4 b^{2^{-1-p}}
   m^{2^{-1-p}} \left(1+m^{2^{-p}}\right)-4 b^{3\ 2^{-1-p}} m^{2^{-1-p}} \left(1+m^{2^{-p}}\right)\right. \\ \left.
+b^{2^{-p}} \left(1+m^{2^{-p}}
   \left(-20+m^{2^{-p}}\right)\right)\right)\\
=16 \left(4 b^2 m^2 \left(\frac{1}{(b-m)^4}-\frac{1}{(-1+b
   m)^4}\right)+b \left(\frac{m^3}{(b-m)^4}-\frac{m}{(-1+b m)^4}\right)\right. \\ \left.
+b^3 \left(\frac{m}{(b-m)^4}-\frac{m^3}{(-1+b m)^4}\right)-\frac{6}{\log ^4\left(\frac{b}{m}\right)}+\frac{6}{\log ^4(b
   m)}\right)
\end{multline}
\end{example}
\begin{example}
\begin{multline}
\sum _{j=1}^{\infty } 4^{-j} \left(\text{sech}^2\left(2^{-2-j} x\right)+2 \text{sech}^2\left(2^{-1-j}
   x\right)-4 \text{sech}^2\left(2^{-j} x\right)\right)\\
=\frac{1}{4} \left(\frac{80}{x^2}-\frac{8}{1+\cosh (x)}-5
   \text{csch}^2\left(\frac{x}{4}\right)+\text{sech}^2\left(\frac{x}{4}\right)\right)
\end{multline}
\end{example}
\begin{example}
\begin{multline}
\sum _{p=0}^n \frac{4^{-p} m^{2^{-1-p}} \left(4^{1+p} a^2
   \left(-1+m^{2^{-p}}\right)^2+\left(1+m^{2^{-p}} \left(6+m^{2^{-p}}\right)\right) \pi
   ^2\right)}{\left(-1+m^{2^{-p}}\right)^3}\\
=4 a^2
   \left(\frac{1}{1-m}+\frac{1}{-1+m^{2^{-1-n}}}\right)-\frac{4 m (1+m) \pi ^2}{(-1+m)^3}+\frac{4^{-n}
   m^{2^{-1-n}} \left(1+m^{2^{-1-n}}\right) \pi ^2}{\left(-1+m^{2^{-1-n}}\right)^3}
\end{multline}
\end{example}
\begin{example}
\begin{multline}
\sum _{p=0}^n \frac{16^{-p} m^{2^{-1-p}} }{\left(-1+m^{2^{-p}}\right)^5}\left(1+16^{1+p} a^4 \left(-1+m^{2^{-p}}\right)^4\right. \\ \left.
+3\times 2^{3+2p} a^2 \left(-1+m^{2^{-p}}\right)^2 \left(1+m^{2^{-p}} \left(6+m^{2^{-p}}\right)\right)\right. \\ \left.
+m^{2^{-p}}\left(76+m^{2^{-p}} \left(230+m^{2^{-p}}
   \left(76+m^{2^{-p}}\right)\right)\right)\right)\\
=\frac{2^{-4 n}
   \left(1+m^{2^{-1-n}}\right)}{\left(-1+m^{2^{-1-n}}\right)^5} \left(m^{2^{-1-n}}+3\times 2^{3+2 n} a^2 m^{2^{-1-n}}
   \left(-1+m^{2^{-1-n}}\right)^2\right. \\ \left.
+2^{3+4 n} a^4 \left(-1+m^{2^{-1-n}}\right)^4+m^{2^{-n}}
   \left(10+m^{2^{-1-n}}\right)\right)\\
-\frac{8 (1+m) \left(a^4 (-1+m)^4+12
   a^2 (-1+m)^2 m+2 m (1+m (10+m))\right)}{(-1+m)^5}
\end{multline}
\end{example}
\begin{example}
\begin{multline}
\sum _{p=0}^{n-1} \frac{2^{-p} m^{2^{-1-p}} \left(2^{1+p}+\log (m)+m^{2^{-p}} \left(-2^{1+p}+\log
   (m)\right)\right)}{\left(-1+m^{2^{-p}}\right)^2}\\
=2\left(\frac{1}{-1+m}+\frac{1}{1-m^{2^{-n}}}+\left(-\frac{m}{(-1+m)^2}+\frac{2^{-n}
   m^{2^{-n}}}{\left(-1+m^{2^{-n}}\right)^2}\right) \log (m)\right)
\end{multline}
\end{example}
\begin{example}
\begin{multline}
\sum _{p=0}^{\infty } \frac{2^{-p} m^{2^{-1-p}} \left(2^{1+p}+\log (m)+m^{2^{-p}} \left(-2^{1+p}+\log
   (m)\right)\right)}{\left(-1+m^{2^{-p}}\right)^2}=\frac{-1+m^2-2 m \log (m)}{(-1+m)^2}
\end{multline}
\end{example}
\begin{example}
\begin{multline}
\sum _{p=0}^{\infty } \frac{8^{-p} m^{2^{-1-p}}}{\left(-1+m^{2^{-p}}\right)^4} \left(-3\times 4^{1+p} \left(-1+m^{2^{-p}}\right)^2 a^2
   \left(2^{1+p} \left(-1+m^{2^{-p}}\right)-\left(1+m^{2^{-p}}\right) \log (m)\right)\right. \\ \left.
+\pi ^2 \left(3\times 2^{1+p}
   \left(-1+m^{2^{-p}}\right) \left(1+m^{2^{-p}} \left(6+m^{2^{-p}}\right)\right)\right.\right. \\ \left.\left.
-\left(1+m^{2^{-p}}\right)
   \left(1+m^{2^{-p}} \left(22+m^{2^{-p}}\right)\right) \log
   (m)\right)\right)\\
=\frac{12 \left(-1+m^2\right)}{(-1+m)^4} \left(a^2 (-1+m)^2-2 m \pi
   ^2\right)+8 m \left(-3 a^2 (-1+m)^2\right. \\ \left.
   +\left(1+4 m+m^2\right) \pi ^2\right) \log (m)
\end{multline}
\end{example}
\begin{example}
\begin{multline}
\sum _{p=0}^{n-1} 4^p \left(\psi ^{(1)}\left(\frac{1}{6}-2^{p-2}\right)+\psi
   ^{(1)}\left(\frac{5}{6}-2^{p-2}\right)\right)\\
=4^n \left(9 \zeta \left(2,-3 2^{n-2}\right)-\zeta
   \left(2,-2^{n-2}\right)\right)-8 \left(10 C+\pi ^2\right)
\end{multline}
\end{example}
\begin{example}
\begin{equation}
\sum _{k=1}^{\infty } \frac{k^k \left(a
   e^{-a}\right)^k}{(k+3)!}=-\frac{4-15 a+18 a^2+18 a^3-54 a^2 e^a+27 a e^{2
   a}-4 e^{3 a}}{108 a^3}
\end{equation}
\end{example}
\begin{example}
\begin{equation}
\sum _{k=1}^{\infty } \frac{k^{k+2} (a \exp (-a))^k}{(k+1)!}=\frac{1-2 a-e^a+3 a e^a-3 a^2 e^a+a^3 e^a}{(-1+a)^3
   a}
\end{equation}
\end{example}
\begin{example}
\begin{equation}
\sum _{k=1}^{\infty } \frac{k^{k+2} (a \exp (-a))^k}{k!}=\frac{a (2 a+1)}{(1-a)^5}
\end{equation}
\end{example}
\begin{example}
\begin{multline}
\sum _{j=0}^{\infty } \sum _{n=0}^{\infty } \sum _{l=0}^{\infty } \sum
   _{p=0}^{\infty } \frac{t^{j+l+n+p}}{(1+j+n) (1+l+p) (1-k)_{1-l-p}
   (k)_{1-j-n}}\\
=\sum _{m=0}^{\infty } \sum _{q=0}^{\infty } \frac{\pi  t^{m+q}
   \csc (k \pi )}{\Gamma (1+k-m) \Gamma (2-k-q)}
\end{multline}
\end{example}
\begin{example}
\begin{multline}
\sum _{j=0}^{\infty } \sum _{m=0}^{\infty } \frac{(-1)^j t^{2 j+m} \Gamma
   (c+j-\alpha ) \Gamma (j+m+\alpha ) \Gamma (2 j+m+\beta ) \Gamma (j+\beta
   -\gamma ) \Gamma (j+\gamma )}{(1+2 j+m) \Gamma (1+j) \Gamma (c+j) \Gamma (1-2
   j+k-m) \Gamma (1+m) \Gamma (c+2 j+m) \Gamma (2 j+\beta )}\\
=\sum _{p=0}^{\infty }
   \sum _{n=0}^{\infty } \frac{t^{n+p} \Gamma (c-\alpha ) \Gamma (n+\alpha )
   \Gamma (p+\alpha ) \Gamma (p+\beta -\gamma ) \Gamma (n+\gamma )}{(1+n+p) \Gamma
   (1+n) \Gamma (c+n) \Gamma (1+k-n-p) \Gamma (1+p) \Gamma (c+p) \Gamma (\alpha
   )}
\end{multline}
\end{example}
\begin{example}
\begin{multline}
\sum _{k=0}^{\infty } \frac{(1-b)_k t^k \, _3F_2(-k,a-k,c;b-k,d;x)}{(k+1)! (1-a)_k}=\int_0^t \frac{\,
   _1F_1(1-b;1-a;s) \, _1F_1(c;d;-x s)}{t} \, ds
\end{multline}
\end{example}
\begin{example}
\begin{multline}
\sum _{p=1}^{n-1} e^{i p \pi } \csc \left(\frac{\pi }{n}\right) \text{csch}^2(m \pi ) \sin \left(\frac{p \pi
   }{n}\right) \left(a n \cosh \left(\frac{m p \pi }{n}\right)\right. \\ \left.
-a n \cosh \left(m \left(2-\frac{p}{n}\right) \pi
   \right)+(-2 n+p) \sinh \left(\frac{m p \pi }{n}\right)+p \sinh \left(m \left(2-\frac{p}{n}\right) \pi
   \right)\right)\\
=\frac{a n \left(\cos \left(\frac{\pi }{n}\right)+\cosh \left(\frac{m \pi }{n}\right)\right)-\sinh
   \left(\frac{m \pi }{n}\right)}{\left(\cos \left(\frac{\pi }{n}\right)+\cosh \left(\frac{m \pi
   }{n}\right)\right)^2}
\end{multline}
\end{example}
\begin{example}
\begin{multline}
\sum _{p=1}^n \frac{(-1)^p e^{-\frac{i p \pi }{n}} \left(-1+e^{\frac{2 i p \pi }{n}}\right) m^{-1-\frac{p}{2
   n}} }{(-1+m)^2}\left(m ((-2+a (-1+m)) n+p-m p)\right. \\ \left.
+m^{p/n} (2 m n+a (n-m n)+p-m p)\right)\\
=-\frac{2 i e^{\frac{2 i \pi
   }{n}} m^{-1+\frac{1}{2 n}} \left(1+a n+m^{1/n} (-1+a n)+2 a m^{\frac{1}{2 n}} n \cos \left(\frac{\pi
   }{n}\right)\right) \sin \left(\frac{\pi }{n}\right)}{\left(e^{\frac{i \pi }{n}}+m^{\frac{1}{2 n}}\right)^2
   \left(1+e^{\frac{i \pi }{n}} m^{\frac{1}{2 n}}\right)^2}
\end{multline}
\end{example}
\begin{example}
\begin{multline}
\sum _{p=0}^{n-1} (-1)^p \cos \left(\frac{j (1+2 p) \pi }{2 n}\right) \left(2 a n \cosh \left(\frac{m (-1+4 n-2
   p) \pi }{2 n}\right)\right. \\ \left.
-2 a n \cosh \left(\frac{m (1+2 p) \pi }{2 n}\right)-(1+2 p) \sinh \left(\frac{m (-1+4 n-2 p) \pi
   }{2 n}\right)\right. \\ \left.
+(-1+4 n-2 p) \sinh \left(\frac{m (1+2 p) \pi }{2 n}\right)\right)\\
=\frac{2 \cos \left(\frac{j \pi }{2
   n}\right) \sinh ^2(m \pi )}{\left(\cos \left(\frac{j \pi
   }{n}\right)+\cosh \left(\frac{m \pi }{n}\right)\right)^2} \left(2 a n \cosh \left(\frac{m \pi }{2 n}\right) \left(\cos \left(\frac{j \pi
   }{n}\right)+\cosh \left(\frac{m \pi }{n}\right)\right)\right. \\ \left.
+\left(-2+\cos \left(\frac{j \pi }{n}\right)-\cosh
   \left(\frac{m \pi }{n}\right)\right) \sinh \left(\frac{m \pi }{2 n}\right)\right)
\end{multline}
\end{example}
\begin{example}
\begin{multline}
\sum _{p=1}^n \frac{2^{-p} q^{2^{-p}} \left(-\pi  \left(1+q^{2^{-p}} \left(-2+q^{2^{-p}} \left(-2+q^{2^{-p}} \left(-2+q^{2^{-p}}\right)\right)\right)\right)-i 2^p
   \left(-1+q^{2^{2-p}}\right) a\right)}{\left(1+q^{2^{1-p}}+q^{2^{-p}}+q^{3\ 2^{-p}}\right)^2}\\
=\pi  \left(\frac{q}{(1+q)^2}-\frac{2^{-n} q^{2^{-n}}}{\left(1+q^{2^{-n}}\right)^2}\right)-i
   \left(-\frac{1}{1+q}+\frac{1}{1+q^{2^{-n}}}\right) a
\end{multline}
\end{example}
\begin{example}
\begin{multline}
\sum _{p=1}^{\infty } \frac{2^{-p} q^{2^{-p}} \left(2^p a \left(-1+q^{2^{2-p}}\right)-\pi  \left(1+q^{2^{-p}} \left(-2+q^{2^{-p}} \left(-2+q^{2^{-p}}
   \left(-2+q^{2^{-p}}\right)\right)\right)\right)\right)}{\left(1+q^{2^{1-p}}+q^{2^{-p}}+q^{3\ 2^{-p}}\right)^2}\\
\frac{2 \pi  q+a \left(-1+q^2\right)}{2 (1+q)^2}
\end{multline}
\end{example}
\begin{example}
\begin{equation}
\sum _{k=0}^{\infty } \frac{(b+a k)^k (-1)^k e^{a y k} y^k}{(k+1)!}=\frac{e^{-a y-b y} \left(e^{a y}-e^{b y}\right)}{(a-b) y}
\end{equation}
\end{example}
\begin{example}
\begin{multline}
\sum _{k=0}^{\infty } \frac{(b+a k)^k (-1)^k e^{a y k} y^k}{(k+2)!}\\
=\frac{e^{(a-b) y-2 (a y+\log (y))} \left(-e^{b y} y+\frac{e^{a y} \left(a-b+2 a^2 y-a b y\right)}{(-2
   a+b)^2}\right)}{a-b}+\frac{e^{-2 (a y+\log (y))} (-1)}{(-2 a+b)^2}
\end{multline}
\end{example}
\begin{example}
\begin{multline}
\sum _{n=0}^{\infty } \sum _{p=0}^{\infty } \sum _{q=0}^{\infty } \frac{(-1)^{2 n+p+2 q} 2^p (a+n \alpha )^{-1-n-p-2 q} (n \alpha +\beta )^{-1+n} E_p(x) \Gamma (1+n+p+2 q)}{\Gamma
   (1+n) \Gamma (1+p) \Gamma (2+2 q)}\\
=\frac{\log \left(\frac{4 \left(\frac{\Gamma \left(\frac{1}{4} (-1+a+2 x-\beta )\right) (-1+a+2 x-\beta )}{\Gamma \left(\frac{1}{4} (-2-1+a+2 x-\beta
   )\right) (-2-1+a+2 x-\beta )}\right)^2}{-1+a+2 x-\beta }\right)}{\beta }
\end{multline}
\end{example}
\begin{example}
\begin{equation}
\sum _{k=1}^{\infty } \frac{\left(a e^{-a}\right)^k k^{-2+k}}{k!}=a-\frac{a^2}{2}
\end{equation}
\end{example}
\begin{example}
\begin{equation}
\sum _{k=1}^{\infty } \frac{\left(a e^{-a}\right)^k k^{-2+k}}{(k+1)!}=\frac{-2-2 a-a^3}{2
   a}+\frac{e^a}{a}
\end{equation}
\end{example}
\begin{example}
\begin{equation}
\sum _{k=0}^{\infty } \frac{a^k E_k}{(k+1)!}=\frac{2 \tan
   ^{-1}\left(\tanh \left(\frac{a}{2}\right)\right)}{a}
\end{equation}
\end{example}
\begin{example}
\begin{equation}
\sum _{n=0}^{\infty } \frac{\zeta (-n+s) z^n}{(n+1)!}=-\frac{(-z)^s
   \Gamma (-s)+\text{Li}_{1+s}(1)-\text{Li}_{1+s}\left(e^z\right)}{z}
   \end{equation}
\end{example}
\begin{example}
\begin{equation}
\sum _{n=0}^{\infty } \frac{\zeta (-n+s) z^n}{(n+2)!}=\frac{(-z)^s z
   \Gamma (-1-s)-z
   \text{Li}_{1+s}(1)-\text{Li}_{2+s}(1)+\text{Li}_{2+s}\left(e^z\right)}{z^2}
\end{equation}
\end{example}
\begin{example}
\begin{equation}
\sum _{r=0}^{\infty } \frac{z^r \zeta (-r+s,v)}{r!}=s (-z)^{-1+s} \Gamma
   (-s)+e^{v z} \Phi \left(e^z,s,v\right)
\end{equation}
\end{example}
\begin{example}
\begin{equation}
\sum _{r=0}^{\infty } \frac{z^r \zeta (-r+s,v)}{(r+1)!}=\int_0^z
   \frac{e^{v t} \Phi \left(e^t,s,v\right)+s \Gamma (-s) (-t)^{-1+s}}{z} \,
   dt
\end{equation}
\end{example}
\begin{example}
\begin{equation}
\sum _{k=1}^{\infty } \frac{x^{2 k} E_{2 k}}{(2 k+1)!}=\frac{-x+2 \tan
   ^{-1}\left(\tanh \left(\frac{x}{2}\right)\right)}{x}
\end{equation}
\end{example}
\begin{example}
\begin{equation}
\sum _{k=0}^{\infty } \frac{p^k \cos (k x)}{(k+1)!}=\frac{-\cos (x)+e^{p
   \cos (x)} \cos (x-p \sin (x))}{p}
\end{equation}
\end{example}
\begin{example}
\begin{multline}
\sum _{n=0}^{\infty } \frac{s^n (c)_n \,
   _2F_1(-n,a;c;z)}{(n+1)!}\\
   =\frac{\,
   _2F_1\left(1,2-c;2-a;\frac{1}{z}\right)-(1-s)^{1+a-c} (1+s (-1+z))^{1-a} \,
   _2F_1\left(1,2-c;2-a;\frac{1+s (-1+z)}{z}\right)}{(-1+a) s z}
\end{multline}
\end{example}
\begin{example}
\begin{multline}
\sum _{k=0}^{\infty } \frac{\, _2F_1(a-k,b;c;x) \left((c-a)_k
   y^k\right)}{(k+1)!}\\
   =-\frac{\left(\frac{-1+x}{x}\right)^{-a-b-c} \,
   _2F_1\left(1-b,-a-b-c;2-b;\frac{1}{x}\right)}{((-1+b) (-1+x))
   y}\\
   -\frac{(1-y)^{a+b+c} \left(-\frac{(-1+x) (-1+y)}{x}\right)^{-a-b-c}
   (1+(-1+x) y)^{1-b} \, _2F_1\left(1-b,-a-b-c;2-b;\frac{1+(-1+x)
   y}{x}\right)}{(-1+b) (-1+x) y}
\end{multline}
\end{example}
\begin{example}
\begin{multline}
\sum _{k=0}^{\infty } \frac{y^k \, _2F_1(a,-k;b+c k;x) (b)_{k+c k}}{k!
   (b)_{c k}}=\frac{\left(1+(1+v)^{1+c} y\right)^b \left(1+(1+v)^{1+c} x
   y\right)^{-a}}{1-c (1+v)^{1+c} y}
\end{multline}
\end{example}
\begin{example}
\begin{multline}
\sum _{k=0}^{\infty } \frac{y^k \, _2F_1(a,-k;b+c k;x) (b)_{k+c
   k}}{(k+1)! (b)_{c k}}\\
   =\frac{1}{y}\left(\frac{1}{(a-b) c}(1+v)^{-1-c} \left(1+(1+v)^{1+c} y\right)^b
   \left(1+(1+v)^{1+c} x y\right)^{-a}\right. \\ \left.
 \left(1+\frac{1+c}{-1+c (1+v)^{1+c}
   y}\right)^{-b} \left(1+\frac{c+x}{x \left(-1+c (1+v)^{1+c} y\right)}\right)^a\right. \\ \left.
   F_1\left(a-b;-b,a;1+a-b;-\frac{1+c}{-1+c (1+v)^{1+c} y},-\frac{c+x}{x
   \left(-1+c (1+v)^{1+c} y\right)}\right)\right. \\ \left.
-\frac{(-c)^{-b}
   (1+v)^{-1-c} \left(1-\frac{c+x}{x}\right)^a
   F_1\left(a-b;-b,a;1+a-b;1+c,\frac{c+x}{x}\right)}{(a-b) c}\right)
\end{multline}
\end{example}
\begin{example}
\begin{multline}
\sum _{k=0}^{\infty } \frac{\left((a)_k (b)_k\right) t^k \,
   _2F_1(a+k,b+k;c+k;x)}{(k+1)! (c)_k}\\
   =\frac{(-1+c) (-\,
   _2F_1(-1+a,-1+b;-1+c;x)+\, _2F_1(-1+a,-1+b;-1+c;t+x))}{(-1+a) (-1+b)
   t}
\end{multline}
\end{example}
\begin{example}
\begin{multline}
\sum _{k=0}^{\infty } \frac{(c-a)_k t^k \,
   _2F_1(a-k,b;c;x)}{(k+1)!}\\
   =\frac{(1-x)^{-b} \left(\,
   _2F_1\left(b,-1-a+c;c;\frac{x}{-1+x}\right)-(1-t)^{1+a-c} \,
   _2F_1\left(b,-1-a+c;c;-\frac{x}{(-1+t) (-1+x)}\right)\right)}{(1+a-c)
   t}
\end{multline}
\end{example}
\begin{example}
\begin{multline}
\sum _{n=0}^{\infty } \frac{(1+a)_n \, _2F_1\left(-n,1+a+b+n;1+a;\frac{1-x}{2}\right) z^n}{(n+1)!}\\
=\int_0^z \frac{2^{a+b} \left(1-s+\sqrt{1+s^2-2 s x}\right)^{-a}
   \left(1+s+\sqrt{1+s^2-2 s x}\right)^{-b}}{\sqrt{1+s^2-2 s x} z} \, ds
\end{multline}
\end{example}
\begin{example}
\begin{equation}
\sum _{j=1}^{\infty } \frac{(-j)^j (-j+t)^{-j}}{1+j}=-1-e^t \Gamma (0,t)+(1+t) \log
   \left(\frac{1+t}{t}\right)
\end{equation}
\end{example}
\begin{example}
\begin{equation}
\sum _{j=1}^{\infty } \frac{(-j)^j j (-j+t)^{-1-j}}{1+j}=e^t \Gamma (0,t)+\log
   \left(\frac{t}{1+t}\right)
\end{equation}
\end{example}
\begin{example}
\begin{equation}
\sum _{j=1}^{\infty } \frac{(-j)^j (t-j)^{k-j} (1-k)_{j-1}}{\Gamma (j+1)}=-e^t \Gamma (k,t)
\end{equation}
\end{example}
\begin{example}
\begin{multline}
\sum _{p=0}^{\infty } (-1)^p \left(\frac{1}{(m+(1+2 p) \alpha
   )^{k+1}}-\frac{1}{(m-(1+2 p) \alpha )^{k+1}}\right) \cos ((1+2 p)
   x)\\
=-\frac{(-i)^k e^{\frac{i m (\pi -2 x)}{2 \alpha }} \left(\frac{\pi
   }{\alpha }\right)^{1+k} \left(\Phi \left(-e^{\frac{i m \pi }{\alpha
   }},-k,\frac{1}{2}-\frac{x}{\pi }\right)+e^{\frac{2 i m x}{\alpha }} \Phi
   \left(-e^{\frac{i m \pi }{\alpha }},-k,\frac{1}{2}+\frac{x}{\pi
   }\right)\right)}{2 \Gamma (1+k)}
\end{multline}
\end{example}
\begin{example}
\begin{multline}
\sum _{p=1}^n \frac{9^{-p} q^{3^{-p}} \left(-1-q^{3^{-p}} \left(4+q^{3^{-p}}\right)-3^p a \left(1+2
   q^{3^{-p}}\right) \left(1+q^{3^{-p}}+q^{2\ 3^{-p}}\right)\right)}{\left(1+q^{3^{-p}}+q^{2\
   3^{-p}}\right)^2}\\
=\frac{q (1+a-a q)}{(-1+q)^2}+\frac{9^{-n} q^{3^{-n}} \left(-1+3^n a
   \left(-1+q^{3^{-n}}\right)\right)}{\left(-1+q^{3^{-n}}\right)^2}
\end{multline}
\end{example}
\begin{example}
\begin{multline}
\sum _{p=1}^{\infty } \frac{9^{-p} q^{3^{-p}} \left(-1-q^{3^{-p}} \left(4+q^{3^{-p}}\right)-3^p a \left(1+2
   q^{3^{-p}}\right) \left(1+q^{3^{-p}}+q^{2\times 3^{-p}}\right)\right)}{\left(1+q^{3^{-p}}+q^{2\times
   3^{-p}}\right)^2}\\
=\frac{q (1+a-a q)}{(-1+q)^2}+\frac{-1+a \log (q)}{\log ^2(q)}
\end{multline}
\end{example}
\begin{figure}[H]
\includegraphics[scale=0.7]{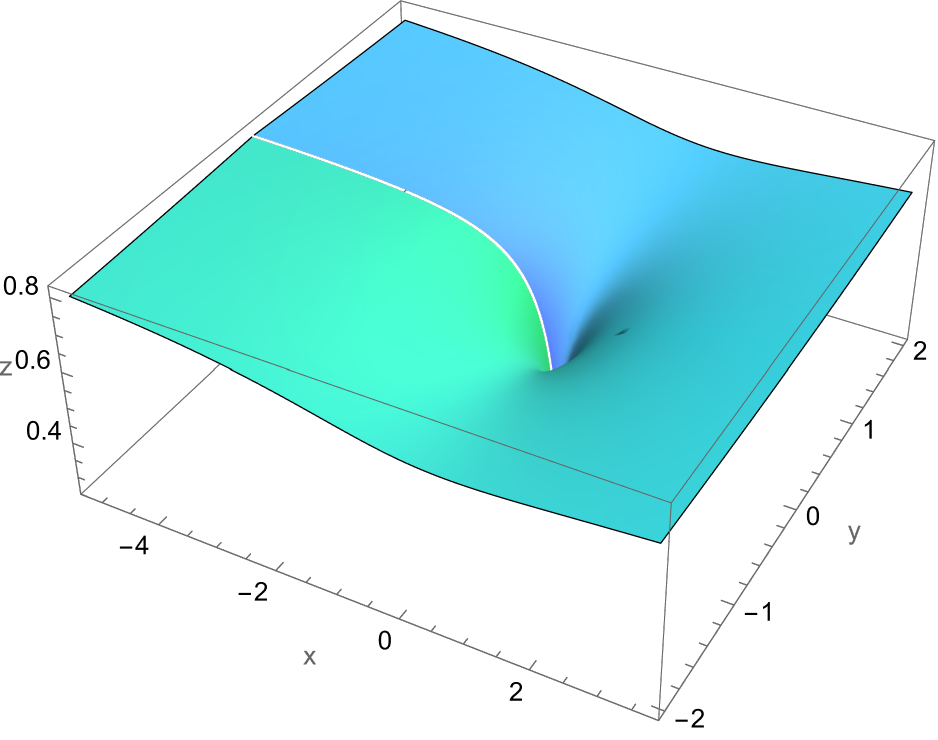}
\caption{Plot of  $f(q)=\frac{1}{(q-1)^2}+\frac{\log (q)-1}{\log ^2(q)}-1$, $q\in\mathbb{C}$.}
   \label{fig:fig2}
\end{figure}
\vspace{-6pt}
\begin{example}
\begin{multline}
\sum _{p=1}^n 3^{-3 p} q^{3^{-p}} \left(\frac{\sqrt[3]{-1} \left(1+3^p a\right)^2+\left(-1+2\times 3^p a
   \left(1+3^p a\right)\right) q^{3^{-p}}-(-1)^{2/3} 9^p a^2 q^{2\times 3^{-p}}}{\left(-1+(-1)^{2/3}
   q^{3^{-p}}\right)^3}\right. \\ \left.
+\frac{2+\left(-1+2\times 3^p a\right) \left(1+\sqrt[3]{-1} q^{3^{-p}}\right)+9^p
   \left(a+\sqrt[3]{-1} a q^{3^{-p}}\right)^2}{\left(1+\sqrt[3]{-1} q^{3^{-p}}\right)^3}\right)\\
=(-1)^{2/3}
   \left(-\frac{q (1+a (-2+a (-1+q)) (-1+q)+q)}{(-1+q)^3}\right. \\ \left.
+\frac{3^{-3 n} q^{3^{-n}} \left(1+q^{3^{-n}}+3^n a
   \left(-1+q^{3^{-n}}\right) \left(-2+3^n a
   \left(-1+q^{3^{-n}}\right)\right)\right)}{\left(-1+q^{3^{-n}}\right)^3}\right)
\end{multline}
\end{example}
\begin{example}
\begin{multline}
\sum _{p=1}^{\infty } 3^{-3 p} q^{3^{-p}} \left(\frac{\sqrt[3]{-1} \left(1+3^p a\right)^2+\left(-1+2\times 3^p a
   \left(1+3^p a\right)\right) q^{3^{-p}}-(-1)^{2/3} 9^p a^2 q^{2\times 3^{-p}}}{\left(-1+(-1)^{2/3}
   q^{3^{-p}}\right)^3}\right. \\ \left.
+\frac{2+\left(-1+2\times 3^p a\right) \left(1+\sqrt[3]{-1} q^{3^{-p}}\right)+9^p
   \left(a+\sqrt[3]{-1} a q^{3^{-p}}\right)^2}{\left(1+\sqrt[3]{-1} q^{3^{-p}}\right)^3}\right)\\
=(-1)^{2/3}
   \left(-\frac{q (1+a (-2+a (-1+q)) (-1+q)+q)}{(-1+q)^3}+\frac{2+a \log (q) (-2+a \log (q))}{\log
   ^3(q)}\right)
\end{multline}
\end{example}
\begin{figure}[H]
\includegraphics[scale=0.7]{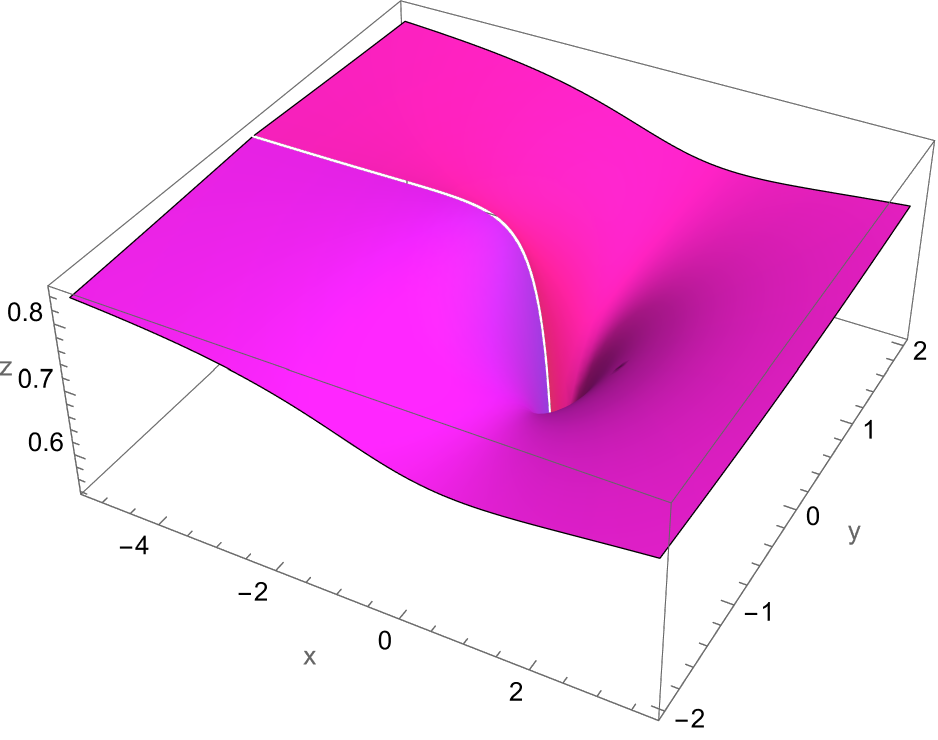}
\caption{Plot of  $f(q)=(-1)^{2/3} \left(\frac{(\log (q)-2) \log (q)+2}{\log ^3(q)}-\frac{q ((q-3) (q-1)+q+1)}{(q-1)^3}\right)$, $q\in\mathbb{C}$.}
   \label{fig:fig2}
\end{figure}
\vspace{-6pt}
\begin{example}
\begin{multline}
\sum _{p=1}^n \frac{3^p r^{3^p} \left(a \left(1+2 r^{3^p}\right) \left(1+r^{3^p}+r^{2\times 3^p}\right)+3^p
   \left(1+r^{3^p} \left(4+r^{3^p}\right)\right)\right)}{\left(1+r^{3^p}+r^{2\times 3^p}\right)^2}\\
=3 \left(\frac{r^3
   \left(3+a-a r^3\right)}{\left(-1+r^3\right)^2}-\frac{3^n r^{3^{1+n}} \left(3^{1+n}+a-a
   r^{3^{1+n}}\right)}{\left(-1+r^{3^{1+n}}\right)^2}\right)
\end{multline}
\end{example}
\begin{example}
\begin{multline}
\sum _{p=1}^{\infty } \frac{3^p r^{3^p} \left(a \left(1+2 r^{3^p}\right) \left(1+r^{3^p}+r^{2\times 3^p}\right)+3^p
   \left(1+r^{3^p} \left(4+r^{3^p}\right)\right)\right)}{\left(1+r^{3^p}+r^{2\times 3^p}\right)^2}=\frac{3 r^3 \left(3+a-a
   r^3\right)}{\left(-1+r^3\right)^2}
\end{multline}
\end{example}
\begin{figure}[H]
\includegraphics[scale=0.7]{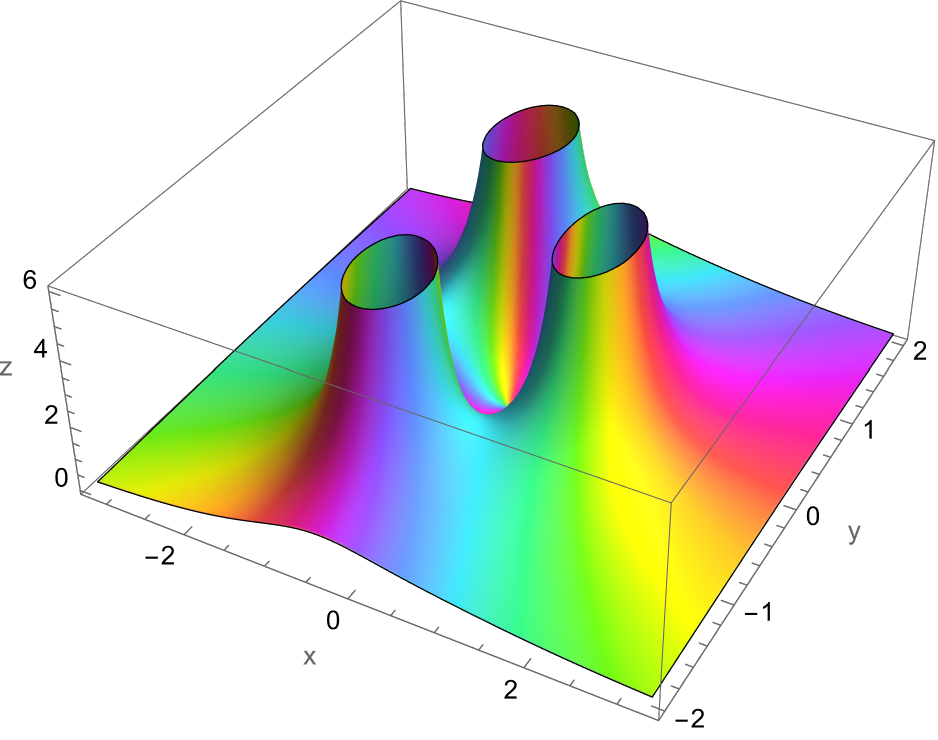}
\caption{Plot of  $f(r)=\frac{9 r^3}{\left(r^3-1\right)^2}$, $r\in\mathbb{C}$.}
   \label{fig:fig2}
\end{figure}
\vspace{-6pt}
\begin{figure}[H]
\includegraphics[scale=0.7]{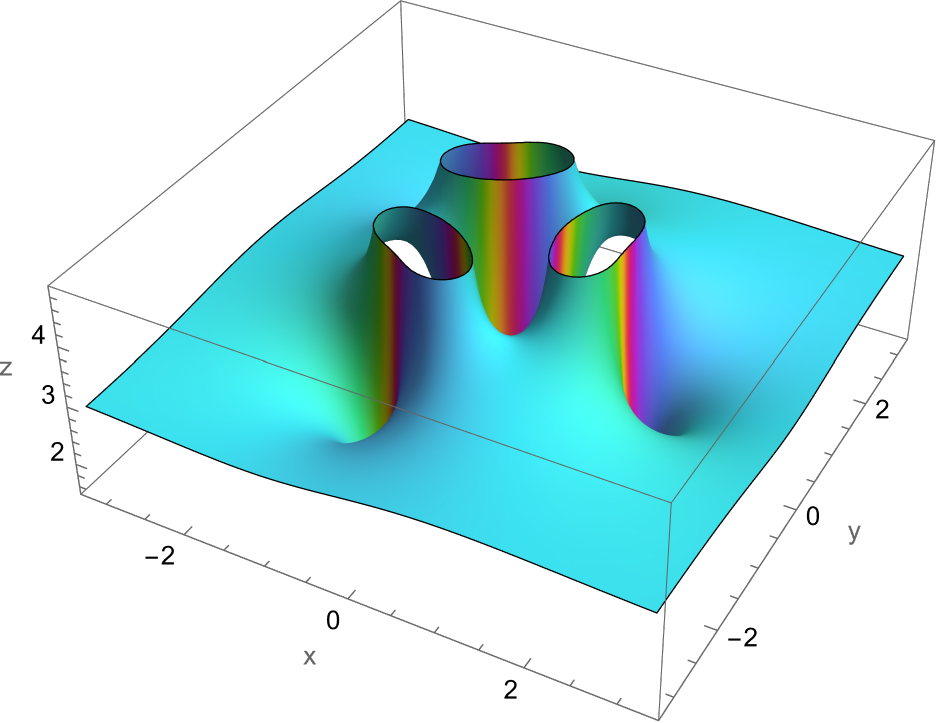}
\caption{Plot of  $f(r)=\frac{3 r^3 \left(4-r^3\right)}{\left(r^3-1\right)^2}$, $r\in\mathbb{C}$.}
   \label{fig:fig2}
\end{figure}
\vspace{-6pt}
\begin{example}
\begin{multline}
\sum _{p=1}^n 3^p r^{3^p} \left(\frac{(-1)^{2/3} \left(3^p+a\right)^2+\left(9^p-2 a \left(3^p+a\right)\right)
   r^{3^p}-\sqrt[3]{-1} a^2 r^{2\times 3^p}}{\left(1+\sqrt[3]{-1} r^{3^p}\right)^3}\right. \\ \left.
+\frac{-\left(3^p+a\right)^2+(-1)^{2/3}
   \left(-9^p+2 a \left(3^p+a\right)\right) r^{3^p}+\sqrt[3]{-1} a^2 r^{2\times 3^p}}{\left(-1+(-1)^{2/3}
   r^{3^p}\right)^3}\right)\\
=-3 \sqrt[3]{-1} \left(\frac{(3+a)^2 r^3+(9-2 a (3+a)) r^6+a^2
   r^9}{\left(-1+r^3\right)^3}\right. \\ \left.
+\frac{3^n r^{3^{1+n}} \left(-9^{1+n} \left(1+r^{3^{1+n}}\right)+a
   \left(-1+r^{3^{1+n}}\right) \left(2\times 3^{1+n}+a-a
   r^{3^{1+n}}\right)\right)}{\left(-1+r^{3^{1+n}}\right)^3}\right)
\end{multline}
\end{example}
\begin{example}
\begin{multline}
\sum _{j=1}^n \sum _{p=1}^n \frac{3^{-2 j+p} q^{3^{-j}+3^p}}{\left(1+q^{3^{-j}}+q^{2\ 3^{-j}}\right)^2 \left(1+q^{3^p}+q^{2\
   3^p}\right)^2} \left(1+q^{3^{-j}} \left(4+q^{3^{-j}}\right)\right. \\ \left.
+3^j a
   \left(1+2 q^{3^{-j}}\right) \left(1+q^{3^{-j}}+q^{2\ 3^{-j}}\right)\right) \left(a \left(1+2 q^{3^p}\right)
   \left(1+q^{3^p}+q^{2\ 3^p}\right)+3^p \left(1+q^{3^p}
   \left(4+q^{3^p}\right)\right)\right)\\
=-3 \left(\frac{q^3 \left(3+a-a q^3\right)}{\left(-1+q^3\right)^2}-\frac{3^n q^{3^{1+n}}
   \left(3^{1+n}+a-a q^{3^{1+n}}\right)}{\left(-1+q^{3^{1+n}}\right)^2}\right)\\
 \left(\frac{q (1+a-a
   q)}{(-1+q)^2}+\frac{9^{-n} q^{3^{-n}} \left(-1+3^n a
   \left(-1+q^{3^{-n}}\right)\right)}{\left(-1+q^{3^{-n}}\right)^2}\right)
\end{multline}
\end{example}
\begin{example}
\begin{multline}
\sum _{j=1}^{\infty } \sum _{p=1}^{\infty } \frac{3^{-2 j+p} q^{3^{-j}+3^p}}{\left(1+q^{3^{-j}}+q^{2\ 3^{-j}}\right)^2 \left(1+q^{3^p}+q^{2\
   3^p}\right)^2} \left(1+q^{3^{-j}} \left(4+q^{3^{-j}}\right)\right. \\ \left.
+3^j a
   \left(1+2 q^{3^{-j}}\right) \left(1+q^{3^{-j}}+q^{2\ 3^{-j}}\right)\right) \left(a \left(1+2 q^{3^p}\right)
   \left(1+q^{3^p}+q^{2\ 3^p}\right)+3^p \left(1+q^{3^p}
   \left(4+q^{3^p}\right)\right)\right)\\
=-\frac{3 q^3 \left(3+a-a q^3\right) \left(\frac{q (1+a-a q)}{(-1+q)^2}+\frac{-1+a \log (q)}{\log
   ^2(q)}\right)}{\left(-1+q^3\right)^2}
\end{multline}
\end{example}
\begin{figure}[H]
\includegraphics[scale=0.7]{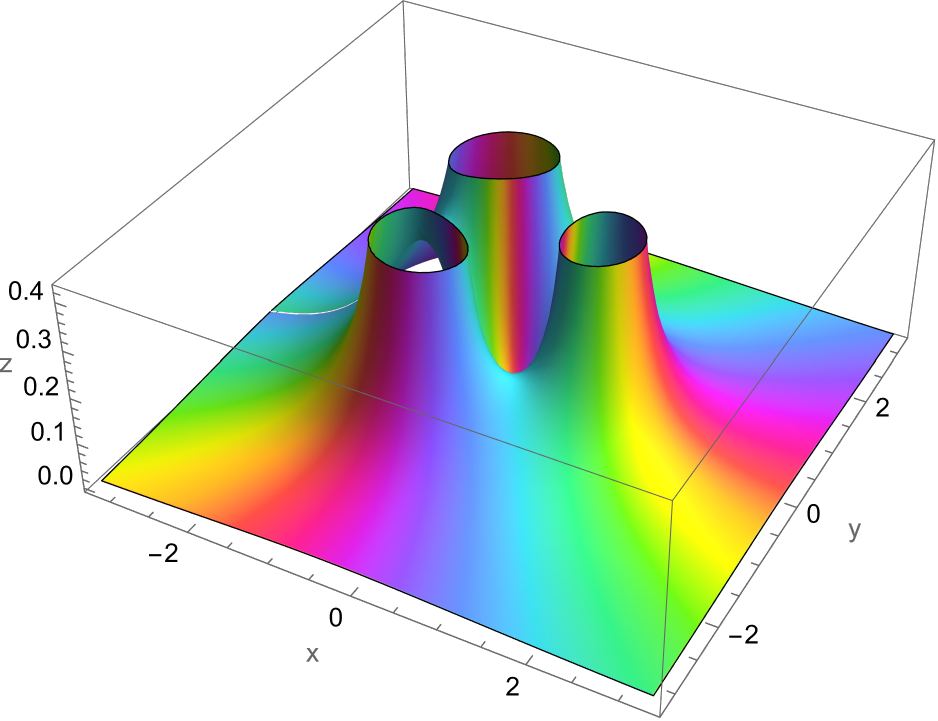}
\caption{Plot of  $f(q)=\frac{9 q^3 \left(\frac{1}{\log ^2(q)}-\frac{q}{(q-1)^2}\right)}{\left(q^3-1\right)^2}$, $q\in\mathbb{C}$.}
   \label{fig:fig2}
\end{figure}
\vspace{-6pt}
\begin{example}
\begin{multline}
\sum _{p=1}^n \frac{q^{\frac{3^{-1+p}}{2}} \left(-3 a \left(-1+q^{2\times 3^{-1+p}}
   \left(-1+q^{3^{-1+p}}+q^{3^p}\right)\right)+i 3^p \left(1+q^{2\times 3^{-1+p}} \left(5+5
   q^{3^{-1+p}}+q^{3^p}\right)\right)\right)}{\left(-1+q^{3^p}\right)^2}\\
=3 i \left(\frac{\sqrt{q} (1+i a
   (-1+q)+q)}{(-1+q)^2}\right. \\ \left.
-\frac{q^{\frac{3^n}{2}} \left(i a \left(-1+q^{3^n}\right)+3^n
   \left(1+q^{3^n}\right)\right)}{\left(-1+q^{3^n}\right)^2}\right)
\end{multline}
\end{example}
\begin{example}
\begin{equation}
\sum _{p=1}^n \frac{3 q^{\frac{3^{-1+p}}{2}} }{-1+q^{3^p}}\left(1+q^{2\times 3^{-1+p}}\right)=\frac{3
   \sqrt{q}}{-1+q}-\frac{3 q^{\frac{3^n}{2}}}{-1+q^{3^n}}
\end{equation}
\end{example}
\begin{example}
\begin{equation}
\sum _{p=1}^{\infty } \frac{q^{\frac{3^{-1+p}}{2}} }{-1+q^{3^p}}\left(1+q^{2\times
   3^{-1+p}}\right)=\frac{\sqrt{q}}{-1+q}
\end{equation}
\end{example}
\begin{example}
\begin{multline}
\sum _{p=1}^n \frac{9^{-p} \left(1+i \sqrt{3}\right) q^{3^{-p}} \left(-1-q^{3^{-p}}
   \left(-4+q^{3^{-p}}\right)+3^p a \left(-1+2 q^{3^{-p}}\right) \left(1+q^{3^{-p}}
   \left(-1+q^{3^{-p}}\right)\right)\right)}{\left(1+q^{3^{-p}} \left(-1+q^{3^{-p}}\right)\right)^2}\\
=2 \sqrt[3]{-1}
   \left(\frac{q (1+a+a q)}{(1+q)^2}-\frac{9^{-n} q^{3^{-n}} \left(1+3^n a
   \left(1+q^{3^{-n}}\right)\right)}{\left(1+q^{3^{-n}}\right)^2}\right)
\end{multline}
\end{example}
\begin{example}
\begin{multline}
\sum _{p=1}^{\infty } \frac{9^{-p} \left(1+i \sqrt{3}\right) q^{3^{-p}} \left(-1-q^{3^{-p}}
   \left(-4+q^{3^{-p}}\right)+3^p a \left(-1+2 q^{3^{-p}}\right) \left(1+q^{3^{-p}}
   \left(-1+q^{3^{-p}}\right)\right)\right)}{\left(1+q^{3^{-p}} \left(-1+q^{3^{-p}}\right)\right)^2}\\
=\frac{\left(1+i
   \sqrt{3}\right) q (1+a+a q)}{(1+q)^2}
\end{multline}
\end{example}
\begin{example}
\begin{multline}
\sum _{p=1}^n 3^{-3 p} q^{3^{-p}} \left(\frac{2}{\left(1+(-1)^{2/3} q^{3^{-p}}\right)^3}+\frac{-1+2\times 3^p
   a}{\left(1+(-1)^{2/3} q^{3^{-p}}\right)^2}+\frac{9^p a^2}{1+(-1)^{2/3} q^{3^{-p}}}\right. \\ \left.
+\frac{-(-1)^{2/3} \left(1+3^p
   a\right)^2+\left(1-2\times 3^p a \left(1+3^p a\right)\right) q^{3^{-p}}+\sqrt[3]{-1} 9^p a^2 q^{2\times
   3^{-p}}}{\left(-1+\sqrt[3]{-1} q^{3^{-p}}\right)^3}\right)\\
=\sqrt[3]{-1} \left(-\frac{q (1-q+a (1+q) (2+a+a
   q))}{(1+q)^3}\right. \\ \left.
+\frac{3^{-3 n} q^{3^{-n}} \left(1-q^{3^{-n}}\
+3^n a \left(1+q^{3^{-n}}\right) \left(2+3^n a
   \left(1+q^{3^{-n}}\right)\right)\right)}{\left(1+q^{3^{-n}}\right)^3}\right)
\end{multline}
\end{example}
\begin{example}
\begin{multline}
\sum _{p=1}^{\infty } 3^{-3 p} q^{3^{-p}} \left(\frac{2}{\left(1+(-1)^{2/3} q^{3^{-p}}\right)^3}+\frac{-1+2\times
   3^p a}{\left(1+(-1)^{2/3} q^{3^{-p}}\right)^2}+\frac{9^p a^2}{1+(-1)^{2/3} q^{3^{-p}}}\right. \\ \left.
+\frac{-(-1)^{2/3}
   \left(1+3^p a\right)^2+\left(1-2\times 3^p a \left(1+3^p a\right)\right) q^{3^{-p}}+\sqrt[3]{-1} 9^p a^2 q^{2\times
   3^{-p}}}{\left(-1+\sqrt[3]{-1} q^{3^{-p}}\right)^3}\right)\\
=-\frac{\sqrt[3]{-1} q (1-q+a (1+q) (2+a+a
   q))}{(1+q)^3}
\end{multline}
\end{example}
\begin{example}
\begin{equation}
\sum _{p=1}^{\infty } \frac{9^p q^{3^{-1+p}} \left(1+q^{3^{-1+p}}
   \left(-4+q^{3^{-1+p}}\right)\right)}{\left(1+q^{3^{-1+p}} \left(-1+q^{3^{-1+p}}\right)\right)^2}=\frac{9
   q}{(1+q)^2}
\end{equation}
\end{example}
\begin{example}
\begin{multline}
\sum _{p=0}^{\infty } (-1)^p \left(\frac{1}{(m+(1+2 p) \alpha
   )^{k+1}}-\frac{1}{(m-(1+2 p) \alpha )^{k+1}}\right) \cos ((1+2 p)
   x)\\
=-\frac{(-i)^k e^{\frac{i m (\pi -2 x)}{2 \alpha }} \left(\frac{\pi
   }{\alpha }\right)^{1+k} \left(\Phi \left(-e^{\frac{i m \pi }{\alpha
   }},-k,\frac{1}{2}-\frac{x}{\pi }\right)+e^{\frac{2 i m x}{\alpha }} \Phi
   \left(-e^{\frac{i m \pi }{\alpha }},-k,\frac{1}{2}+\frac{x}{\pi
   }\right)\right)}{2 \Gamma (1+k)}
\end{multline}
\end{example}
\begin{example}
\begin{multline}
\sum _{p=0}^{\infty } (-1)^p \left(\frac{1}{(-2 m+(1+2 p) \beta )^{k+1}}+\frac{(-1)^k}{(2 m+(1+2 p) \beta
   )^{k+1}}\right)\\
=\frac{e^{\frac{i m \pi }{\beta }} \pi ^{1+k} \left(\frac{i}{\beta }\right)^k \Phi \left(-e^{\frac{2 i
   m \pi }{\beta }},-k,\frac{1}{2}\right)}{\beta  \Gamma (1+k)}
\end{multline}
\end{example}
\begin{example}
\begin{multline}
\sum _{p=0}^{\infty } \frac{(-1)^p ((-1+2 p) \log (-3-2 p)+(3+2 p) \log (1-2 p))}{-3+4 p (1+p)}\\
=\frac{1}{2} \pi 
   \left(2 i+\gamma -i \pi +\log \left(\frac{2 \pi  \Gamma \left(-\frac{1}{4}\right)^2}{9 \Gamma
   \left(-\frac{3}{4}\right)^2}\right)\right)
\end{multline}
\end{example}
\begin{example}
\begin{multline}
\sum _{p=0}^{\infty } (-1)^p (1+2 p) \left(\frac{16 x}{4 x^2-(\beta +2 p \beta )^2}-\frac{2 \log \left(\frac{(-2
   x+\beta +2 p \beta )^2}{(2 x+\beta +2 p \beta )^2}\right)}{\beta +2 p \beta }\right)\\
=\frac{2 i \log (a)
   \left(\text{Li}_2\left(-i e^{-\frac{i \pi  x}{\beta }}\right)-\text{Li}_2\left(i e^{-\frac{i \pi  x}{\beta
   }}\right)-\text{Li}_2\left(-i e^{\frac{i \pi  x}{\beta }}\right)+\text{Li}_2\left(i e^{\frac{i \pi  x}{\beta
   }}\right)\right)}{\pi }\\
+\frac{2 i \pi  x \beta  \log (a)+4 \beta  \log \left(-i+\frac{2}{-i+e^{\frac{i \pi  x}{\beta
   }}}\right)-4 \pi  x \sec \left(\frac{\pi  x}{\beta }\right)}{\beta ^2}
\end{multline}
\end{example}
\begin{example}
\begin{multline}
\sum _{p=0}^{\infty } \left(\frac{\Gamma \left(k,\frac{i \pi  (-2 m+\beta +2 p \beta )}{2 \beta }\right)}{(-2
   m+\beta +2 p \beta )^k}+\frac{\Gamma \left(k,-\frac{i \pi  (2 m+\beta +2 p \beta )}{2 \beta }\right)}{(-2 m-(1+2 p)
   \beta )^k}\right)=\left(\frac{\pi  i}{\beta }\right)^k \text{Li}_{1-k}\left(-e^{\frac{2 i m \pi }{\beta
   }}\right)
\end{multline}
\end{example}
\begin{example}
\begin{multline}
\sum _{p=0}^{\infty } \left(\frac{\Gamma \left(k,\frac{1}{2} i (-1+2 p) \pi \right)}{(1-2 p)^k}+\frac{\Gamma
   \left(k,-\frac{1}{2} i (3+2 p) \pi \right)}{(3+2 p)^k}\right)=\left(-1+2^k\right) (-i \pi )^k \zeta (1-k)
\end{multline}
\end{example}
\begin{example}
\begin{multline}
\sum _{p=0}^{\infty } \left(\left(\Gamma \left(0,\frac{1}{2} i (-1+2 p) \pi \right)+\Gamma \left(0,-\frac{1}{2} i
   (3+2 p) \pi \right)\right) \log \left(-\frac{1}{2} (i \pi )\right)\right. \\ \left.
+G_{2,3}^{3,0}\left(\frac{1}{2} i (-1+2 p) \pi \left |
\begin{array}{c}
 1,1 \\
 0,0,0 \\
\end{array}
\right.\right)+G_{2,3}^{3,0}\left(-\frac{1}{2} i (3+2 p) \pi \left |
\begin{array}{c}
 1,1 \\
 0,0,0 \\
\end{array}
\right)\right.\right)\\
=\frac{1}{2} \log (2) \left(2 \gamma +i \pi +\log \left(\frac{1}{2 \pi ^2}\right)\right)
\end{multline}
\end{example}
\begin{example}
\begin{multline}
\sum _{p=0}^{\infty } \left(E_{-k}\left(\frac{1}{2} i (-1+2 p) \pi \right)+E_{-k}\left(-\frac{1}{2} i (3+2 p) \pi
   \right)\right)=2^{1+k} \left(-1+2^{1+k}\right) \zeta (-k)
\end{multline}
\end{example}
\begin{example}
\begin{multline}
\sum _{p=1}^n
   \left(\frac{2^p}{\left(1+q^{2^{-1+p}}\right)^2}+\frac{-2^p+a}{1+q^{2^{-1+p}}}-\frac{2^{1+p}}{\left(1+q^{2^p}\right)^2}+\frac{2^{1+p}-a}{1+q^{2^p}}\right)\\
=\frac{a+(-2+a) q}{(1+q)^2}+\frac{2^{1+n}
   q^{2^n}}{\left(1+q^{2^n}\right)^2}-\frac{a}{1+q^{2^n}}
\end{multline}
\end{example}
\begin{example}
\begin{multline}
\sum _{p=1}^{\infty }
   \left(\frac{2^p}{\left(1+q^{2^{-1+p}}\right)^2}+\frac{-2^p+a}{1+q^{2^{-1+p}}}-\frac{2^{1+p}}{\left(1+q^{2^p}\right)^2}+\frac{2^{1+p}-a}{1+q^{2^p}}\right)=-\frac{q (2+a+a q)}{(1+q)^2}
\end{multline}
\end{example}
\begin{figure}[H]
\includegraphics[scale=0.7]{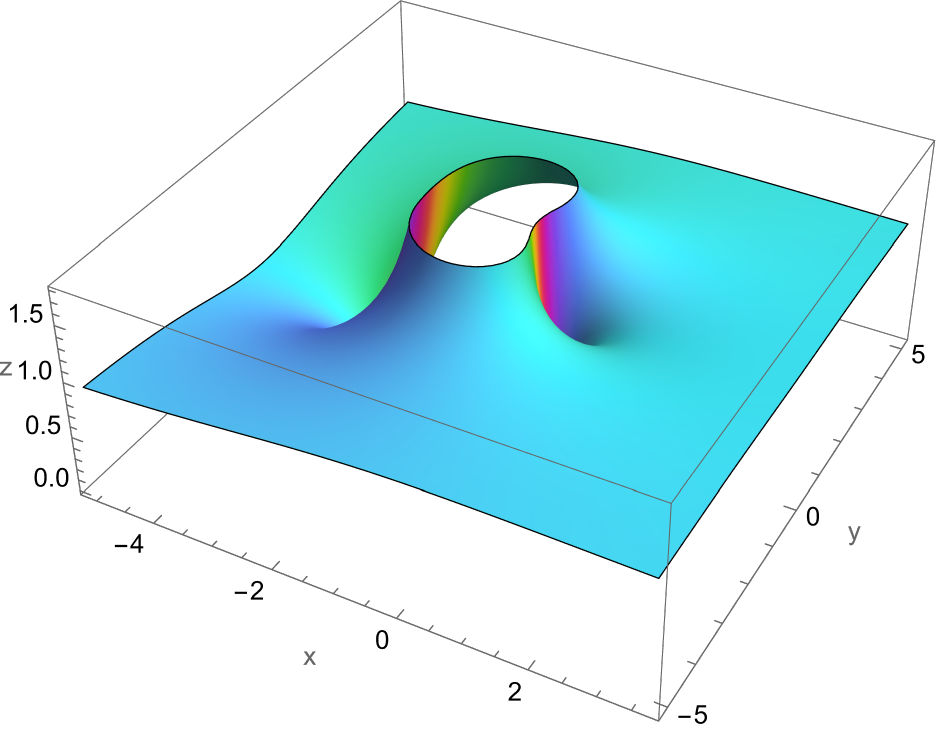}
\caption{Plot of  $f(q)=-\frac{q (q+3)}{(q+1)^2}$, $q\in\mathbb{C}$.}
   \label{fig:fig2}
\end{figure}
\vspace{-6pt}
\begin{example}
\begin{multline}
\sum _{p=1}^n 4^{-p} \left(\frac{1}{\left(1+q^{2^{-1-p}}\right)^2}+\frac{-1+2^{1+p}
   a}{1+q^{2^{-1-p}}}-\frac{2}{\left(1+q^{2^{-p}}\right)^2}-\frac{2 \left(-1+2^p
   a\right)}{1+q^{2^{-p}}}\right)\\
=\frac{\sqrt{q} (1-2 a (-1+q)+q)}{(-1+q)^2}+\frac{4^{-n} q^{2^{-1-n}}
   \left(-1-2^{1+n} a+\left(-1+2^{1+n} a\right) q^{2^{-n}}\right)}{\left(-1+q^{2^{-n}}\right)^2}
\end{multline}
\end{example}
\begin{example}
\begin{multline}
\sum _{p=1}^{\infty } 4^{-p} \left(\frac{1}{\left(1+q^{2^{-1-p}}\right)^2}+\frac{-1+2^{1+p}
   a}{1+q^{2^{-1-p}}}-\frac{2}{\left(1+q^{2^{-p}}\right)^2}-\frac{2 \left(-1+2^p
   a\right)}{1+q^{2^{-p}}}\right)\\
=\frac{\sqrt{q} (1-2 a (-1+q)+q)}{(-1+q)^2}+\frac{2 (-1+a \log (q))}{\log
   ^2(q)}
\end{multline}
\end{example}
\begin{example}
\begin{multline}
-\frac{\left(\beta ^{-q}\right)^{-\frac{m}{q}} \Gamma \left(n-\frac{m}{q}\right) \Gamma \left(\frac{m}{q}\right) \left(\log \left(\beta ^{-q}\right)+\psi ^{(0)}\left(n-\frac{m}{q}\right)-\psi
   ^{(0)}\left(\frac{m}{q}\right)\right)}{q^2 \Gamma (n)}\\
=\sum _{t=0}^{-1+n} \sum _{j=0}^{-1+n} \frac{(-1)^{j+t} e^{\frac{i m \pi }{q}} (2 \pi )^{2-t} \left(-\frac{1}{q}\right)^t
   \left(\frac{i}{q}\right)^{-t} \left(\beta ^{-q}\right)^{-\frac{m}{q}}  \left(-\frac{m-q}{q}\right)_{-1-j+n}
   (2-t)_t S_j^{(t)}}{q^2 j! (-1-j+n)!}\\
\Phi \left(e^{\frac{2 i m \pi }{q}},-1+t,\frac{\pi +i \log \left(\beta ^{-q}\right)}{2 \pi }\right)
\end{multline}
\end{example}
\begin{example}
\begin{multline}
\frac{\left(\beta ^{-q}\right)^{-\frac{m}{q}} \Gamma \left(n-\frac{m}{q}\right) \Gamma \left(\frac{m}{q}\right)}{q \Gamma (n)} \left(\log ^2\left(\beta ^{-q}\right)+\psi ^{(0)}\left(n-\frac{m}{q}\right)^2\right. \\ \left.
+2 \psi
   ^{(0)}\left(n-\frac{m}{q}\right) \left(\log \left(\beta ^{-q}\right)-\psi ^{(0)}\left(\frac{m}{q}\right)\right)-2 \log \left(\beta ^{-q}\right) \psi ^{(0)}\left(\frac{m}{q}\right)+\psi
   ^{(0)}\left(\frac{m}{q}\right)^2\right. \\ \left.
+\psi ^{(1)}\left(n-\frac{m}{q}\right)+\psi ^{(1)}\left(\frac{m}{q}\right)\right)\\
=\sum _{t=0}^{-1+n} \sum _{j=0}^{-1+n} \frac{(-1)^{j+t} e^{\frac{i m \pi
   }{q}} (2 \pi )^{3-t} \left(-\frac{1}{q}\right)^t \left(\frac{i}{q}\right)^{1-t} \left(\beta ^{-q}\right)^{-\frac{m}{q}}  \left(-\frac{m-q}{q}\right)_{-1-j+n} (3-t)_t S_j^{(t)}}{j! (-1-j+n)!}\\
\Phi \left(e^{\frac{2 i m \pi }{q}},-2+t,\frac{\pi +i \log \left(\beta ^{-q}\right)}{2 \pi
   }\right)
\end{multline}
\end{example}
\subsubsection{Formulae derived from equation (\ref{eq_52414}) and section (7) in \cite{plos}}
\begin{example}
\begin{equation}
\sum _{k=0}^{\infty } \frac{(k+1)^{k-1} (y \exp (-y))^k}{(k+2)!}=-\frac{e^y \left(2-2 e^y+2 y+y^3\right)}{2
   y^2}
\end{equation}
\end{example}
\begin{example}
\begin{multline}
\sum _{k=0}^{\infty } \frac{(k+1)^{k-1} (y \exp (-y))^k}{(k+3)!}=-\frac{e^y \left(-1+e^{2 y}-16 e^y y+2 y (1+y)
   \left(7+2 y^2\right)\right)}{16 y^3}
\end{multline}
\end{example}
\begin{example}
\begin{multline}
\sum _{k=0}^{\infty } \frac{(k+1)^{k-1} (y \exp (-y))^k}{(k+4)!}\\
=\frac{e^y \left(-16+16 e^{3 y}-243 e^{2 y}
   y+1944 e^y y^2-3 y (-65+6 y (85+y (85+6 y (5+3 y))))\right)}{3888 y^4}
\end{multline}
\end{example}
\begin{example}
\begin{multline}
\sum _{k=0}^{\infty } \frac{(k+1)^{k-1} (y \exp (-y))^k}{(k+5)!}\\
=\frac{e^y }{995328
   y^5}\left(243-243 e^{4 y}+4096 e^{3 y}
   y-31104 e^{2 y} y^2+165888 e^y y^3\right. \\ \left.
   -4 y (781+6 y (-865+12 y (415+y (415+12 y (13+6 y)))))\right)
\end{multline}
\end{example}
\begin{example}
\begin{multline}
\sum _{k=0}^{\infty } \frac{(k+1)^{k-1} (y \exp (-y))^k}{(k+6)!}\\
=\frac{e^y }{77760000000 y^6}\left(-995328+995328 e^{5
   y}-18984375 e^{4 y} y+160000000 e^{3 y} y^2\right. \\ \left.
   -810000000 e^{2 y} y^3+3240000000 e^y y^4-5 y (-2801547\right. \\ \left.
   +20 y (965041+30
   y (-153713+60 y (12019+y (12019+60 y (77+30 y))))))\right)
\end{multline}
\end{example}
\begin{example}
\begin{multline}
\sum _{k=0}^{\infty } \frac{(k+1)^{k-1} (y \exp (-y))^k}{(k+7)!}\\
=\frac{e^y }{4199040000000 y^7}\left(2500000-2500000 e^{6
   y}+53747712 e^{5 y} y\right. \\ \left.
   -512578125 e^{4 y} y^2+2880000000 e^{3 y} y^3-10935000000 e^{2 y} y^4+34992000000 e^y y^5\right. \\ \left.
   -3 y
   (12915904+15 y (-6418657+20 y (1568371+90 y (-66801+20 y (13489\right. \\ \left.
   +y (13489+180 y (29+10
   y)))))))\right)
\end{multline}
\end{example}
\begin{example}
\begin{multline}
\sum _{k=0}^{\infty } \frac{(k+1)^{k-1} (y \exp (-y))^k}{(k+8)!}\\
=\frac{e^y}{169446409937280000000 y^8} \left(-4199040000000+4199040000000
   e^{7 y}\right. \\ \left.
   +71490737500000 y-100884017500000 e^{6 y} y-582029398598592 y^2\right. \\ \left.+1084457023598592 e^{5 y}
   y^2+3048374148022665 y^3-6894792071015625 e^{4 y} y^3\right. \\ \left.
   -11819395880919900 y^4+29054597040000000 e^{3 y}
   y^4+38514655699821000 y^5\right. \\ \left.
   -88253338509000000 e^{2 y} y^5-138426709885380000 y^6+235342236024000000 e^y
   y^6\right. \\ \left.
   -138426709885380000 y^7-53552365952400000 y^8-16810159716000000 y^9\right)
\end{multline}
\end{example}
\begin{example}
\begin{multline}
\sum _{k=0}^{\infty } \frac{(k+1)^{k-1} (y \exp (-y))^k}{(k+9)!}\\
=-\frac{e^y}{22209679843299164160000000 y^9}
   \left(-20684376213046875+20684376213046875 e^{8 y}\right. \\ \left.
   +384901561175625000 y-550376570880000000 e^{7 y}
   y-3420799013537500000 y^2\right. \\ \left.
   +6611534970880000000 e^{6 y} y^2+19430599723471550208 y^3-47380650332371550208 e^{5 y}
   y^3\right. \\ \left.
   -80099863968742248960 y^4+225928546583040000000 e^{4 y} y^4+261589888296724377600 y^5\right. \\ \left.
   -761648828645376000000 e^{3
   y} y^5-750777617069512704000 y^6+1927923597508608000000 e^{2 y} y^6\right. \\ \left.
   +2455119676203909120000
   y^7-4406682508591104000000 e^y y^7+2455119676203909120000 y^8\right. \\ \left.+946256377960857600000 y^9+275417656786944000000
   y^{10}\right)
\end{multline}
\end{example}
\subsubsection{Formulae derived from equation (\ref{eq_51322}) and section (7) in \cite{plos}}
\begin{example}
\begin{equation}
\sum _{k=0}^{\infty } \frac{(k+2)^k (y \exp (-y))^{k+1}}{(k+3)!}=-\frac{2 e^y-2 e^{2 y}+2 e^y y+y^2+e^y
   y^3}{2 y^2}
\end{equation}
\end{example}
\begin{example}
\begin{multline}
\sum _{k=0}^{\infty } \frac{(k+2)^k (y \exp (-y))^{k+1}}{(k+4)!}=-\frac{3 e^{3 y}-48 e^{2 y} y+8 y^3+3 e^y
   \left(-1+2 y (1+y) \left(7+2 y^2\right)\right)}{48 y^3}
\end{multline}
\end{example}
\begin{example}
\begin{multline}
\sum _{k=0}^{\infty } \frac{(k+2)^k (y \exp (-y))^{k+1}}{(k+5)!}\\
=-\frac{1}{24}-\frac{5
   e^y}{36}-\frac{e^y}{243 y^4}+\frac{e^{4 y}}{243 y^4}+\frac{65 e^y}{1296 y^3}-\frac{e^{3 y}}{16 y^3}-\frac{85
   e^y}{216 y^2}+\frac{e^{2 y}}{2 y^2}-\frac{85 e^y}{216 y}-\frac{e^y y}{12}
\end{multline}
\end{example}
\begin{example}
\begin{multline}
\sum _{k=0}^{\infty } \frac{(k+2)^k (y \exp (-y))^{k+1}}{(k+6)!}\\
=-\frac{1}{120}-\frac{13
   e^y}{288}+\frac{e^y}{4096 y^5}-\frac{e^{5 y}}{4096 y^5}-\frac{781 e^y}{248832 y^4}+\frac{e^{4 y}}{243
   y^4}+\frac{865 e^y}{41472 y^3}-\frac{e^{3 y}}{32 y^3}-\frac{415 e^y}{3456 y^2}\\
   +\frac{e^{2 y}}{6 y^2}-\frac{415
   e^y}{3456 y}-\frac{e^y y}{48}
\end{multline}
\end{example}
\subsubsection{Formulae derived using equation (\ref{eq_51320}) and section (7) in \cite{plos}}
\begin{example}
\begin{equation}
\sum _{k=1}^{\infty } \frac{k^{k-1} (y \exp (-y))^k}{(k+2)!}=\frac{-1+e^{2 y}+6 y-8 e^y y+6 y^2+4 y^3}{8 y^2}
\end{equation}
\end{example}
\begin{example}
\begin{multline}
\sum _{k=1}^{\infty } \frac{k^{k-1} (y \exp (-y))^k}{(k+3)!}=\frac{8-8 e^{3 y}-57 y+81 e^{2 y} y+198 y^2-324 e^y y^2+198
   y^3+108 y^4}{648 y^3}
\end{multline}
\end{example}
\begin{example}
\begin{multline}
\sum _{k=1}^{\infty } \frac{k^{k-1} (y \exp (-y))^k}{(k+4)!}\\
=\frac{-81+81 e^{4 y}+700 y-1024 e^{3 y} y-2760 y^2+5184 e^{2
   y} y^2+7200 y^3-13824 e^y y^3+7200 y^4+3456 y^5}{82944 y^4}
\end{multline}
\end{example}
\begin{example}
\begin{multline}
\sum _{k=1}^{\infty } \frac{k^{k-1} (y \exp (-y))^k}{(k+5)!}\\
=\frac{1}{1296000000
   y^5}\left(82944-82944 e^{5 y}-850905 y+1265625 e^{4 y} y+3974300
   y^2\right. \\ \left.
   -8000000 e^{3 y} y^2-11397000 y^3+27000000 e^{2 y} y^3+24660000 y^4-54000000 e^y y^4\right. \\ \left.
   +24660000 y^5+10800000 y^6\right)
\end{multline}
\end{example}
\section{Table of arctangent definite integrals}
\begin{center}
\setlength{\tabcolsep}{9pt} % Default value: 6pt
\renewcommand{\arraystretch}{2.0} % Default value: 1
\begin{tabular}{ | l | c | r }
  \hline			
  $f(x)$ & $\int_{0}^{\infty}f(x)dx$  \\
  \hline
  $\cos (m x) \tan ^{-1}(\sinh (\alpha ) \text{sech}(b x))$ & $\frac{\pi  \text{sech}\left(\frac{\pi  m}{2 b}\right) \sin \left(\frac{\alpha 
   m}{b}\right)}{2 m}$  \\
  $\frac{\left(x^2-1\right) \tanh ^{-1}(\text{sech}(\pi  x))}{\left(x^2+1\right)^2}$ & $-\frac{1}{2} \pi  (\log (4)-1)$  \\
   $\frac{\left(12 x^2-1\right) \tanh ^{-1}(\text{sech}(\pi  x))}{\left(4 x^2+1\right)^3}$ & $\frac{1}{8} (\pi -2 \pi  C)$  \\
    $\frac{\left(x^4-6 x^2+1\right) \tanh ^{-1}(\text{sech}(\pi  x))}{\left(x^2+1\right)^4}$ & $\frac{1}{12} \pi  (3 \zeta (3)-2)$  \\
     $\left(\frac{1}{(1+i x)^6}-\frac{1}{(x+i)^6}\right) \tanh ^{-1}(\text{sech}(\pi  x))$ & $\frac{1}{40} \pi  (15 \zeta (5)-8)$  \\
      $\left(\frac{1}{(1+i x)^{3/2}}+\frac{1}{(1-i x)^{3/2}}\right) \tanh ^{-1}(\text{sech}(\pi  x))$ & $-4 \left(\sqrt{2}-1\right) \pi  \zeta
   \left(\frac{1}{2}\right)-2 \pi$  \\
       $\left(\frac{1}{\sqrt{1+i x}}+\frac{1}{\sqrt{1-i x}}\right) \tanh ^{-1}(\text{sech}(\pi  x))$ & $\left(1-2 \sqrt{2}\right) \zeta
   \left(\frac{3}{2}\right)+2 \pi$  \\
       $\frac{\tan ^{-1}(\sinh (\alpha ) \text{sech}(b x))}{\log ^2(a)+x^2}$ & $-\frac{i \pi  \log \left(\frac{\Gamma
   \left(\frac{-2 i \alpha +2 b \log (a)+\pi }{4 \pi }\right) \Gamma \left(\frac{2 i \alpha +2 b \log (a)+3 \pi }{4
   \pi }\right)}{\Gamma \left(\frac{-2 i \alpha +2 b \log (a)+3 \pi }{4 \pi }\right) \Gamma \left(\frac{2 i \alpha +2
   b \log (a)+\pi }{4 \pi }\right)}\right)}{2 \log (a)}$  \\
       $\frac{\tanh ^{-1}(\text{sech}(\pi  x))}{x^2+1}$ & $-\frac{1}{2} \pi  \log
   \left(\frac{2}{\pi }\right)$  \\
       $\frac{\tanh ^{-1}(\sin (\alpha ) \text{sech}(b x))}{a^2+x^2}$ & $-\frac{\pi 
   \log \left(\frac{\Gamma \left(\frac{2 a b-2 \alpha +3 \pi }{4 \pi }\right)
   \Gamma \left(\frac{2 a b+2 \alpha +\pi }{4 \pi }\right)}{\Gamma \left(\frac{2
   a b-2 \alpha +\pi }{4 \pi }\right) \Gamma \left(\frac{2 a b+2 \alpha +3 \pi
   }{4 \pi }\right)}\right)}{2 a}$  \\
       $\frac{\tanh ^{-1}(\text{sech}(x))}{\left(x^2+\pi
   ^2\right)^2}$ & $\frac{\log \left(\frac{2 \pi }{e}\right)}{4 \pi
   ^2}$  \\
       $\frac{\coth ^{-1}\left(\sqrt{2} \cosh (9 \pi 
   x)\right)}{x^2+4}$ & $-\frac{1}{4} \pi  \log \left(\frac{\Gamma
   \left(\frac{75}{8}\right) \Gamma \left(\frac{77}{8}\right)}{\Gamma
   \left(\frac{73}{8}\right) \Gamma \left(\frac{79}{8}\right)}\right)$  \\
       $\text{sech}(t) \tanh ^{-1}(\sin (\alpha ) \text{sech}(a b \sinh
   (t)))$ & $-\frac{1}{2} \pi  \log \left(\frac{\Gamma \left(\frac{2 a b-2 \alpha +3
   \pi }{4 \pi }\right) \Gamma \left(\frac{2 a b+2 \alpha +\pi }{4 \pi
   }\right)}{\Gamma \left(\frac{2 a b-2 \alpha +\pi }{4 \pi }\right) \Gamma
   \left(\frac{2 a b+2 \alpha +3 \pi }{4 \pi }\right)}\right)$  \\
        $\text{sech}(t) \tanh ^{-1}\left(\text{sech}\left(\frac{1}{4} \pi  \sinh
   (t)\right)\right)$ & $-\frac{1}{2} \pi  \log \left(\frac{\Gamma
   \left(\frac{5}{8}\right)^2}{\Gamma \left(\frac{1}{8}\right) \Gamma
   \left(\frac{9}{8}\right)}\right)$  \\
         $$ & $$  \\
  \hline  
\end{tabular}
\end{center}
\begin{center}
\setlength{\tabcolsep}{10pt} % Default value: 6pt
\renewcommand{\arraystretch}{3.0} % Default value: 1
\begin{tabular}{ |c|c|c|c| } 
\hline
$f(x_{1} \ldots x_{j+1})$ & $\int_{0}^{\infty} \ldots \int_{0}^{\infty} f(x_{1} \ldots x_{j+1})dx_{1} \ldots dx_{j+1}$  \\
\hline
$\frac{x_1^{m-1} }{x_1^n+1}\left(\prod _{l=1}^j \frac{x_{l+1}^{2^{l-1}
   m+\frac{n}{2}-1}}{x_{l+1}^n+1}\right)$ & $2^j \pi ^{j+1} n^{-j-1} \csc
   \left(\frac{\pi  2^j m}{n}\right)$  \\
$\frac{x_1^p \left(\prod _{l=1}^j \frac{x_{l+1}^{2^{l-1}
   p+\frac{n}{2}-1}}{x_{l+1}^n+1}\right)-x_1^m \left(\prod _{l=1}^j
   \frac{x_{l+1}^{2^{l-1} m+\frac{n}{2}-1}}{x_{l+1}^n+1}\right)}{x_1
   \left(x_1^n+1\right) \log \left(x_1 \left(\prod _{l=1}^j
   x_{l+1}^{2^{l-1}}\right)\right)}$ & $2 \pi ^j n^{-j} \left(\tanh
   ^{-1}\left(e^{\frac{i \pi  2^j m}{n}}\right)-\tanh ^{-1}\left(e^{\frac{i \pi 
   2^j p}{n}}\right)\right)$  \\
$\frac{\sqrt[6]{x_1} \left(\prod _{l=1}^j
   \frac{x_{l+1}^{2^{j-1}+\frac{2^l}{3}-1}}{x_{l+1}^{2^j}+1}\right)-\prod
   _{l=1}^j \frac{x_{l+1}^{2^{j-1}+2^{l-2}-1}}{x_{l+1}^{2^j}+1}}{\sqrt{x_1}
   \left(x_1^{2^j}+1\right) \log \left(x_1 \left(\prod _{l=1}^j
   x_{l+1}^{2^{l-1}}\right)\right)}$ & $\frac{1}{2} \left(2^j\right)^{-j} \pi ^j
   \log (3)$  \\ 
$\frac{\prod _{l=1}^j \frac{x_{l+1}^{2^{j-1}+2^{l-2}-1}}{x_{l+1}^{2^j}+1}
   \log ^k\left(x_1 \left(-\prod _{l=1}^j
   x_{l+1}^{2^{l-1}}\right)\right)}{\sqrt{x_1}
   \left(x_1^{2^j}+1\right)}$ & $\left(2^j\right)^{-j-1} i^k \left(1-2^{k+1}\right)
   (2 \pi )^{j+k+1} \zeta (-k)$  \\ 
$\frac{\prod _{l=1}^j
   \frac{x_{l+1}^{2^{j-1}+2^{l-2}-1}}{x_{l+1}^{2^j}+1}}{\sqrt{x_1}
   \left(x_1^{2^j}+1\right) \log \left(x_1 \left(-\prod _{l=1}^j
   x_{l+1}^{2^{l-1}}\right)\right)}$ & $-i \left(2^j\right)^{-j} \pi ^j \log
   (2)$  \\ 
$\frac{\prod _{l=1}^j
   \frac{x_{l+1}^{2^{j-1}+2^{l-2}-1}}{x_{l+1}^{2^j}+1}}{\sqrt{x_1}
   \left(x_1^{2^j}+1\right) \log ^3\left(x_1 \left(-\prod _{l=1}^j
   x_{l+1}^{2^{l-1}}\right)\right)}$ & $\frac{3}{16} i \left(2^j\right)^{-j} \pi
   ^{j-2} \zeta (3)$  \\ 
$\frac{\prod _{l=1}^j
   \frac{x_{l+1}^{2^{j-1}+2^{l-2}-1}}{x_{l+1}^{2^j}+1}}{\sqrt{x_1}
   \left(x_1^{2^j}+1\right) \log ^5\left(x_1 \left(-\prod _{l=1}^j
   x_{l+1}^{2^{l-1}}\right)\right)}$ & $-\frac{15}{256} i \left(2^j\right)^{-j} \pi
   ^{j-4} \zeta (5)$  \\ 
$\frac{\prod _{l=1}^j \frac{x_{l+1}^{2^{j-1}+2^{l-2}-1}}{x_{l+1}^{2^j}+1}
   \log \left(\log \left(x_1 \left(-\prod _{l=1}^j
   x_{l+1}^{2^{l-1}}\right)\right)\right)}{\sqrt{x_1} \left(x_1^{2^j}+1\right)
   \log \left(x_1 \left(-\prod _{l=1}^j
   x_{l+1}^{2^{l-1}}\right)\right)}$ & $\frac{1}{2} \left(2^j\right)^{-j} \pi ^j
   \log (2) \left(2 i \gamma +\pi -i \log \left(8 \pi ^2\right)\right)$  \\ 
$\frac{\prod _{l=1}^j \frac{x_{l+1}^{2^{j-1}+2^{l-2}-1}}{x_{l+1}^{2^j}+1}
   \log \left(\log \left(x_1 \left(-\prod _{l=1}^j
   x_{l+1}^{2^{l-1}}\right)\right)\right)}{\sqrt{x_1} \left(x_1^{2^j}+1\right)
   \log ^2\left(x_1 \left(-\prod _{l=1}^j
   x_{l+1}^{2^{l-1}}\right)\right)}$ & $\frac{1}{48} \left(2^j\right)^{-j} \pi
   ^{j+1} (-24 \log (A)+2 \gamma -i \pi +\log (4))$  \\ 
%$$ & $$  \\ 
 %  $$ & $$  \\
$\frac{\prod _{l=1}^j \frac{x_{l+1}^{2^{j-1}+2^{l-3}-1}}{x_{l+1}^{2^j}+1}
   \log \left(\log \left(a x_1 \left(\prod _{l=1}^j
   x_{l+1}^{2^{l-1}}\right)\right)\right)}{x_1^{3/4}
   \left(x_1^{2^j}+1\right)}$ & $\multirow{3}{13em}{$\frac{\left(2^j\right)^{-j} \pi ^{j+1} (i \pi +\log (64)+2 \log (\pi
   ))}{\sqrt{2}}+(-1+i) \sqrt{2} \left(2^j\right)^{-j} \pi ^{j+1} \log
   \left(\frac{\Gamma \left(\frac{\pi -i \log (a)}{8 \pi
   }\right)}{\Gamma \left(\frac{5}{8}-\frac{i \log (a)}{8 \pi
   }\right)}\right)+(-1-i) \sqrt{2} \left(2^j\right)^{-j} \pi ^{j+1} \log
   \left(\frac{\Gamma \left(\frac{3}{8}-\frac{i \log (a)}{8 \pi
   }\right)}{\Gamma \left(\frac{7}{8}-\frac{i \log (a)}{8 \pi
   }\right)}\right)$}$  \\ 
         $$ & $$  \\
   %            $$ & $$  \\
   %                  $$ & $$  \\
%$$ & $\multirow{4}{6em}{$$}$  \\ 
%$$ & $$  \\ 
    %$$ & $$  \\
        $$ & $$  \\
%         $$ & $$  \\
%          $$ & $$  \\
%$$ & $$  \\ 
\hline
\end{tabular}
\end{center}
%
%
%\section{Conclusion}
%%
%In this paper, we have presented a novel method for deriving a new Arctangent integral transform along with some interesting definite integrals similar to those published by Oberhettinger, using contour integration. The results presented were numerically verified for both real and imaginary and complex values of the parameters in the integrals using Mathematica by Wolfram.
%%

%

\begin{thebibliography}{999}
%
\bibitem{reyn4} Reynolds, R.; Stauffer, A.
{A Method for Evaluating Definite Integrals in Terms of Special Functions with Examples}.  \emph{Int. Math. Forum} (\textbf{2020}), \emph{15}, 235--244, doi:10.12988/imf.2020.91272 
%
\bibitem{plos}Reynolds R, Stauffer A (\textbf{2023}) 
 \emph{Chebyshev series: Derivation and evaluation}. PLoS ONE
18(3): e0282703. 
https://doi.org/10.1371/journal.pone.0282703
 %
 \bibitem{dlmf} Olver, F.W.J.; Lozier, D.W.; Boisvert, R.F.; Clark, C.W. (Eds.)
 \emph{NIST Digital Library of Mathematical Functions}; U.S. Department of Commerce, National Institute of Standards and Technology: Washington, DC, USA; Cambridge University Press: Cambridge, UK, 2010; With 1 CD-ROM (Windows, Macintosh and UNIX). MR 2723248 (\textbf{2012a}:33001).
%
\bibitem{grad} Gradshteyn, I.S.; Ryzhik, I.M.
\emph{Tables of Integrals, Series and Products}, 6th ed.; Academic Press: Cambridge, MA, USA, 
 (\textbf{2000}).
%
\bibitem{berndt1}Bruce C. Berndt, 
 \emph{Ramanujan's Notebooks Part II}, (\textbf{1989}), Springer-Verlag New York Entry (23)
 %
  \bibitem{bateman}Bateman, Harry, and Bateman Manuscript Project. (\textbf{1954}). 
 \emph{Tables of Integral Transforms} (version Published). Volume 1, New York: McGraw-Hill Book Company. Eq. (1.9.51)
 %
 \bibitem{brychkov}Brychkov, Y., Marichev, O., \& Savischenko, N. (\textbf{2018}). 
 \emph{Handbook of Mellin Transforms} (1st ed.). Chapman and Hall/CRC. 
 https://doi.org/10.1201/9780429434259
Eq. (2.1.5.7)
 %
 \bibitem{prud1}A. P. Prudnikov, Yu. A. Brychkov, and O. I. Marichev (\textbf{1986b}) 
 \emph{Integrals and Series: Special Functions, Vol. 2.}, Gordon \& Breach Science Publishers, New York.
 %
 \bibitem{bdh} D. Bierens de Haan (\textbf{1939}) 
 \emph{Nouvelles Tables d’Intégrales Définies}. G.E. Stechert \& Co., New York.
 %
 \bibitem{schroder}Ernst Schr\"{o}der, Zeitschrift fur Mathematik und Physik, vol. 25, pp. 106–117 (\textbf{1880}).
%
 \bibitem{obert}F. Oberhettinger (\textbf{1990}) Tables of Fourier Transforms and Fourier Transforms of Distributions. Springer-Verlag, Berlin. 
 %
\bibitem{bateman1}Bateman, Harry, and Bateman Manuscript Project.( \textbf{1954}). Tables of Integral Transforms (version Published). Volume I, New York: McGraw-Hill Book Company.
 %
 \bibitem{polyanin}Andrei D. Polyanin, Alexander V. Manzhirov, 
 \emph{Handbooks of mathematical equations, Handbook of integral equations}, Chapman \& Hall CRC (\textbf{2008}). Eq. (4.3.10)
 %
 \bibitem{berndti}Berndt, B.C. 
 \emph{Integrals associated with Ramanujan and elliptic functions}. Ramanujan J 41, 369-389 (2016). https://doi.org/10.1007/s11139-016-9793-1
 %
 \bibitem{winckler}von Anton Winckler, 
 \emph{Über die Eigenschaften einiger bestimmten Integrale}, ,(vorgelegt in der Sitzung vom 3. Jänner \textbf{1861}), K. K. Court and State Dr., 1861, pp.359
 %
  \bibitem{prud2}A. P. Prudnikov, Yu. A. Brychkov, and O. I. Marichev (\textbf{1986b}) 
 \emph{Integrals and Series: Special Functions, Vol. 2.}, Gordon \& Breach Science Publishers, New York. 
 %
 \bibitem{prud3}A. P. Prudnikov, Yu. A. Brychkov, and O. I. Marichev (\textbf{1990}) 
 \emph{Integrals and Series: More Special Functions, Vol. 3.},  Gordon and Breach Science Publishers, New York.
 %
 \bibitem{prud4}A. P. Prudnikov, Yu. A. Brychkov, and O. I. Marichev (1992a) 
 \emph{Integrals and Series: Direct Laplace Transforms}, Vol. 4. Gordon and Breach Science Publishers, New York.
 %
 \bibitem{prud5}A. P. Prudnikov, Yu. A. Brychkov, and O. I. Marichev (1992b) 
 \emph{Integrals and Series: Inverse Laplace Transforms}, Vol. 5. Gordon and Breach Science Publishers, New York. 
 %
  \bibitem{yu}Yu. A. Brychkov \& P. C. Sofotasios (\textbf{2023}) 
 \emph{On some properties of the Neumann polynomials}, Integral Transforms and Special Functions, 34:4, 316-333, 
 DOI:10.1080/10652469.2022.2118738. Eq. (1.4) 
 %
  \bibitem{bromwich}Bromwich, T.J., Watson, G.N.,  
 \emph{An Introduction to the Theory of Infinite Series}. Nature 78, 242 (\textbf{1908}). 
 https://doi.org/10.1038/078242a0
 %
  \bibitem{chatterjea}S. K. Chatterjea, 
 \emph{On a generating function of Laguerre polynomials}., Bollettino dell'Unione Matematica Italiana, Serie 3, Vol. 17, (\textbf{1962}), n.2, p. 179-182.
 %
  \bibitem{buchholz}H. Buchholz (\textbf{1969}) 
 \emph{The Confluent Hypergeometric Function with Special Emphasis on Its Applications}. Springer-Verlag, New York. 
 %
  \bibitem{hansen}E. R. Hansen (\textbf{1975}) 
 \emph{A Table of Series and Products}. Prentice-Hall, Englewood Cliffs, N.J..
 %
  \bibitem{dieckmann}Dieckmann, A. A Collection of Infinite Products and Series. Available online: http://www-elsa.physik.uni-bonn.de/~{}dieckman/InfProd/InfProd.html (accessed on 2 November \textbf{2023}).
 %
  \bibitem{andrews}G. E. Andrews, R. Askey, and R. Roy (\textbf{1999}) 
 \emph{Special Functions. Encyclopedia of Mathematics and its Applications}, Vol. 71, Cambridge University Press, Cambridge. Eq. (2) on page 626
 %
  \bibitem{berndt}Bruce C. Berndt 
 \emph{Modular transformations and generalizations of several formulae of Ramanujan}," Rocky Mountain Journal of Mathematics, Rocky Mountain J. Math. 7(1), 147-190, (Winter \textbf{1977}) Eq. (3.31)
 %
 %
  \bibitem{sharma}A. R. Vasishtha, Dr. S. K. Sharma, A. K. Vasishtha, 
\emph{Krishna's Series Trigonometry and Algebra: For the Degree Part First Students of C.C.S. University, Meerut and all other Indian Universities and for various Competitive Examination like I.A.S., P.C.S., etc}., Krishna Prakashan Media. 6th Edition, (\textbf{2009}).
%
\bibitem{hind}Hind, John. 
 \emph{The Elements of Plane and Spherical Trigonometry}. United Kingdom: Deighton, Bell \& Company, (\textbf{1855}). pp. 170 exercise 9.
 %
   \bibitem{bedford}Bedford, T.,Quigley, J., Walls, L., Alkali, B., Daneshkhah, A., Hardman, G. , Advances in Mathematical Modeling for Reliability. Iran: IOS Press, (\textbf{2008}).
 %
   \bibitem{jolley}Jolley, L. B. W. (Leonard Benjamin William), 
 \emph{Summation of series, b}, New York, Dover Publications, (\textbf{1886})
 %
  \bibitem{edwards}Edwards, Joseph, 1854-1931.
 \emph{ Differential Calculus With Applications And Numerous Examples: an Elementary Treatise}.,  London,etc: Macmillan and co., (\textbf{1886}). pp.71.
 %
  \bibitem{brafman}Brafman, Fred. 
 \emph{"Generating Functions of Jacobi and Related Polynomials."}, Proceedings of the American Mathematical Society 2, no. 6 (\textbf{1951}): 942[49. https://doi.org/10.2307/2031712.
 %
  \bibitem{durell}C. V. Durell and A. Robson, Advanced Trigonometry, Courier Corporation, 352 pages, (\textbf{1882}).
 %
 \bibitem{hall}Hall, Henry Sinclair., Knight, Samuel Ratcliffe. 
 \emph{Elementary Trigonometry}. India: Macmillan and Company, (\textbf{1893}.) Eq. (18) pp. 292
 %
 \bibitem{hobson}Hobson, Ernest William (1856-1933), 
 \emph{A treatise on plane trigonometry}, (\textbf{1918}), Cambridge, England, University Press
 %
 \bibitem{borwein}Borwein, J. M. and Borwein, P. B.
 \emph{ Evaluation of Sums of Reciprocals of Fibonacci Sequences}." §3.7 in Pi \& the AGM: A Study in Analytic Number Theory and Computational Complexity. New York: Wiley, pp. 91-101, (\textbf{1987}).
 %
 \bibitem{berndt2}Bruce C. Berndt, 
 \emph{Ramanujan's Notebooks  Part II}, (\textbf{1989}), Springer-Verlag New York Entry (23)
 %
   \bibitem{berndt4}Berndt, B. C. 
 \emph{"q-Series."}, Ch. 27 in Ramanujan's Notebooks, Part IV. New York:Springer-Verlag, pp. 261-286, (\textbf{1994}).
 %
 \bibitem{bateman1}Bateman, Harry, and Bateman Manuscript Project. (\textbf{1954}). 
 \emph{Tables of Integral Transforms }(version Published). Volume 1., New York: McGraw-Hill Book Company.
 %
 \bibitem{bateman2}Bateman, Harry, and Bateman Manuscript Project. (\textbf{1954}). 
 \emph{Tables of Integral Transforms }(version Published). Volume 2., New York: McGraw-Hill Book Company.
 %
 \bibitem{pathak} R. S. Pathak. O. P. Singh. 
 \emph{"Finite Hankel transforms of distributions.}." Pacific J. Math. 99 (2) 439 - 458, (\textbf{1982}). 
 %
  \bibitem{blagouchine}Blagouchine, I.V. 
 \emph{Rediscovery of Malmsten’s integrals, their evaluation by contour integration methods and some related results}. Ramanujan J 35, 21–110 (\textbf{2014}). 
 https://doi.org/10.1007/s11139-013-9528-5
 %
  \bibitem{jeffreys}Jeffreys, H. and Jeffreys, B. S. 
 \emph{Frullani's Integrals. §12.16 in Methods of Mathematical Physics}, 3rd ed. Cambridge, England: Cambridge University Press, pp. 406-407, (\textbf{1988}).
 %
 %
\end{thebibliography}
\end{document}